w. b. vasantha kandasamy

# SMARANDACHE FUZZY ALGEBRA

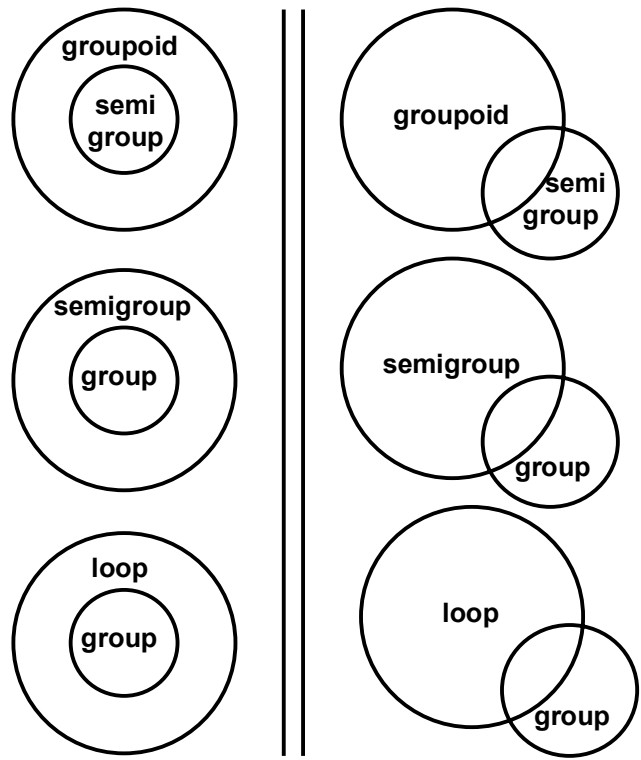

2003

# SMARANDACHE FUZZY ALGEBRA

W. B. Vasantha Kandasamy
Department of Mathematics
Indian Institute of Technology Madras
Chennai – 600 036, India
e-mail: vasantha@iitm.ac.in
web: http://mat.iitm.ac.in/~wbv

*The picture on the cover is a simple graphic illustration depicting the classical algebraic structures with single binary operations and their Smarandache analogues. The pictures on the left, composed of concentric circles, depicts the traditional conception of algebraic structures, and the pictures of the right, with their liberal intersections, describe Smarandache algebraic structures. In fact, Smarandache Algebra, like its predecessor, Fuzzy Algebra, arose from the need to define structures which were more compatible with the real world where the grey areas mattered. Lofti A Zadeh, the father of fuzzy sets, remarked that: "So, this whole thing started because of my perception at that time, that the world of classical mathematics – was a little too much of a black and white world, that the principle of the 'excluded middle' meant that every proposition must be either true or false. There was no allowance for the fact that classes do not have sharply defined boundaries." So, here is this book, which is an amalgamation of alternatives.*

# CONTENTS

**Preface** 5

# PART ONE

# PART TWO

# PREFACE

In 1965, Lofti A. Zadeh introduced the notion of a fuzzy subset of a set as a method for representing uncertainty. It provoked, at first (and as expected), a strong negative reaction from some influential scientists and mathematicians—many of whom turned openly hostile. However, despite the controversy, the subject also attracted the attention of other mathematicians and in the following years, the field grew enormously, finding applications in areas as diverse as washing machines to handwriting recognition. In its trajectory of stupendous growth, it has also come to include the theory of fuzzy algebra and for the past five decades, several researchers have been working on concepts like fuzzy semigroup, fuzzy groups, fuzzy rings, fuzzy ideals, fuzzy semirings, fuzzy near-rings and so on.

In this book, we study the subject of Smarandache Fuzzy Algebra. Originally, the revolutionary theory of Smarandache notions was born as a paradoxist movement that challenged the status quo of existing mathematics. The genesis of Smarandache Notions, a field founded by Florentine Smarandache, is alike to that of Fuzzy Theory: both the fields imperatively questioned the dogmas of classical mathematics.

Despite the fact that Fuzzy Algebra has been studied for over fifty years, there are only two books on fuzzy algebra. But both the books do not cover topics related to fuzzy semirings, fuzzy near-rings etc. so we have in this book, two parts: In Part 1 we have recalled all the definitions and properties of fuzzy algebra. In Part II we give Smarandache fuzzy algebraic notions. This is the first book in fuzzy algebra which covers the notions of fuzzy semirings and fuzzy near-rings though there are several papers on these two concepts.

This book has seven chapters, which are divided into two parts. Part I contains the first chapter, and Part II encloses the remaining six chapters. In the first chapter, which is subdivided into twelve sections, we deal with eleven distinct fuzzy algebraic concepts and in the concluding section list the miscellaneous properties of fuzzy algebra. The eleven fuzzy algebraic concepts which we analyze are fuzzy sets, fuzzy subgroups, fuzzy sub-bigroups, fuzzy rings, fuzzy birings, fuzzy fields, fuzzy semirings, fuzzy near-rings, fuzzy vector spaces, fuzzy semigroups and fuzzy half-groupoids. The results used in these sections are extensive and we have succeeded in presenting new concepts defined by several researchers. In the second chapter we introduce the notion of Smarandache fuzzy semigroups and its properties and also study Smarandache fuzzy bisemigroups. In the third chapter, we define the notion of Smarandache fuzzy half-groupoids and their generalizations (Smarandache fuzzy groupoids and bigroupoids, Smarandache fuzzy loops and biloops).

Chapter four deals with Smarandache fuzzy rings and Smarandache non-associative fuzzy rings. This chapter includes Smarandache fuzzy vector spaces and Smarandache birings. The study of Smarandache fuzzy semirings and its generalizations comprises the fifth chapter. Likewise, in the sixth chapter we analyze Smarandache fuzzy near-rings and its generalizations. In these six chapters, we have succeeded in introducing around 664 concepts related to Smarandache fuzzy algebra. The reader is expected to be well-versed with a strong background in Algebra, Fuzzy Algebra and Smarandache algebraic notions.

The final chapter in this book deals with the applications of Smarandache Fuzzy algebraic structures. I do not claim that I have exhausted all the possibilities of applications, all that I have done here is to put forth those concepts that clearly have relevant applications. When I informed my interest in writing this book, Dr. Minh Perez of the American Research Press, editor of the Smarandache Notions Journal, a close research associate and inspiration-provider par excellence, insisted, rather subtly, that I try to find applications for these Smarandache notions. I was worried a little bit about finding the right kind of applications to suit this book, and then I happened to come across an perceptive interview with the Father of Fuzzy Sets, Lofti. A. Zadeh. Emphasizing about the long time it takes for a new subject to secure its place in the spotlight, he says, "Now: Probabilistic computing. It is interesting that within Artificial Intelligence it is only within the past several years that it has become sort of accepted. Previous to that it was not accepted. There was an article in the New York Times about Bayesian things. It says this technology is 276 years old. Another example that comes to mind is holography. Garbor came up with his first paper in 1946; I saw the paper. No applications until the laser was invented! It's only after laser was invented that holography became useful. And then he got the Nobel Prize. Sometimes it has to await certain things. …. *So, sometimes it's a matter of some application that all of the sudden brings something to light. Sometimes it needs that kind of thing.*" Somewhere between those lines, I could find the hope that I had longed for. It made me attest to the fact that research is generally a legacy, and that our effort will subsequently stand up to speak for itself.

Since I am generalizing now, and speaking of hope and resurrection and the legacy of effort, and also about movements that challenge the dogmas and the irrationality of tradition, I am also aware of how all of this resonates with the social aspects of our life.

Thinking about society, about revolution and revolt, and about the crusades against domination and dogma, I dedicate this book to Periyar (Literally meaning, The Great Man), the icon of rationalism. He single-handedly led the non-brahmins of South India, to a cultural, political and social awakening, freeing them from the cruel bonds of slavery that traditional brahminism foisted upon them. He was the first political leader in India to fight for the concepts of *Self-Respect* and *Social Justice*; and in terms of social reform, he stands unparalleled. His writings and speeches, which I read with the rigour that is expected of serious research, are now a permanent part of my personal faith. Periyar's ideology and political praxis have influenced me overwhelmingly, and his thought drives me to dissent and to dare.

# PART ONE

# Fuzzy Algebra

**Chapter One**

# SOME RESULTS ON FUZZY ALGEBRA

This chapter has twelve sections. First section we introduce the concept of fuzzy sets. As there are very few books on fuzzy algebra we have tried our level best to introduce all the possible definitions of fuzzy groups, fuzzy rings, fuzzy vector spaces, fuzzy near rings. Section two is devoted to the definition of fuzzy groups and some of its basic properties. Section three solely deals with the study and introduction of fuzzy sub-bigroup of a group. Fuzzy rings and its properties are introduced in section four. Section five introduces the notions of fuzzy birings. Study of fuzzy fields is carried out in section six. Study of fuzzy semirings and their generalizations are given in section seven. Section eight gives the properties of fuzzy near-rings and its properties. We describe the notions of fuzzy vector spaces and fuzzy bivector spaces in section nine. A brief study of fuzzy semigroups is carried out in the tenth section. The generalization of fuzzy half groupoids and its generalizations are given in section eleven. The final section, which is quite radical in nature gives the miscellaneous properties in fuzzy algebraic structures.

## 1.1 Fuzzy Subsets

In 1965 Zadeh [144] mathematically formulated the fuzzy subset concept. He defined fuzzy subset of a non-empty set as a collection of objects with grade of membership in a continuum, with each object being assigned a value between 0 and 1 by a membership function. Fuzzy set theory was guided by the assumption that classical sets were not natural, appropriate or useful notions in describing the real life problems, because every object encountered in this real physical world carries some degree of fuzziness. Further the concept of grade of membership is not a probabilistic concept.

**DEFINITION 1.1.1:** *Let X be a non-empty set. A fuzzy set (subset) $\mu$ of the set X is a function $\mu: X \to [0, 1]$.*

**DEFINITION 1.1.2:** *Let $\mu$ be a fuzzy subset of a set X. For $t \in [0, 1]$, the set $X_\mu^t = \{x \in X \mid \mu(x) \geq t\}$ is called a t-level subset of the fuzzy subset $\mu$.*

**DEFINITION 1.1.3:** *A fuzzy set of a set X is called a fuzzy point if and only if it takes the value 0 for all $y \in X$ except one, say, $x \in X$. If its value at x is t, $(0 < t \leq 1)$ then we denote this fuzzy point by $x_t$.*

**DEFINITION 1.1.4:** *The complement of a fuzzy set $\mu$ of a set X is denoted by $\mu^c$ and defined as $\mu^c(x) = 1 - \mu(x)$ for every $x \in X$.*

We mainly give definitions, which pertain to algebraic operations, or to be more precise we are not interested in discussing concepts topologically or analytically like continuity, connected, increasing function or decreasing function. Just we proceed on to define when are two functions disjoint and the concept of min max functions.

**DEFINITION 1.1.5:** *Two fuzzy subsets $\mu$ and $\lambda$ of a set X are said to be disjoint if there exists no $x \in X$ such that $\mu(x) = \lambda(x)$.*

**DEFINITION 1.1.6:** *The union of two fuzzy sets $\lambda$ and $\mu$ of a set X, denoted by $\lambda \cup \mu$ is a fuzzy subset of the set X defined as $(\lambda \cup \mu)(x) = \max\{\lambda(x), \mu(x)\}$ for every $x \in X$. The intersection of two fuzzy (subsets) sets $\lambda$ and $\mu$ of a set X, written as $\lambda \cap \mu$, is a fuzzy subset of X defined as $(\lambda \cap \mu)(x) = \min\{\lambda(x), \mu(x)\}$ for every $x \in X$.*

**DEFINITION 1.1.7:** *Let $\lambda$ and $\mu$ be two fuzzy subsets of a set X. Then $\lambda$ is said to be contained in $\mu$, written as $\lambda \subseteq \mu$ if $\lambda(x) \leq \mu(x)$ for every $x \in X$. If $\lambda(x) = \mu(x)$ for every $x \in X$ then we say $\lambda$ and $\mu$ are equal and write $\lambda = \mu$.*

**DEFINITION 1.1.8:** *A fuzzy subset $\mu$ of a set X is said to normal if*

$$\sup_{x \in X} \mu(x) = 1.$$

*A fuzzy subset $\mu$ of a set X is said to be normalized if there exist $x \in X$ such that $\mu(x) = 1$.*

**DEFINITION 1.1.9:** *Let $f : X \to Y$ be a function. For a fuzzy set $\mu$ in Y, we define $(f^{-1}(\mu))(x) = \mu(f(x))$ for every $x \in X$.*

*For a fuzzy set $\lambda$ in X, $f(\lambda)$ is defined by*

$$(f(\lambda))(y) = \begin{cases} \sup \lambda(x) & \text{if } f(z) = y, z \in X \\ 0 & \text{if there is no such } x \end{cases}$$

*where $y \in Y$.*

**DEFINITION 1.1.10:** *Let X be any set. A fuzzy subset $\mu$ in the set X has the sup property if for any subset A of the set X there exists $x_0 \in A$ such that $\mu(x_0) = \sup\{\mu(x) \mid x \in A\}$.*

**DEFINITION 1.1.11:** *Let $\lambda$ and $\mu$ be fuzzy subsets of the sets X and Y respectively. The cartesian product of $\lambda$ and $\mu$ is defined as $\lambda \times \mu: X \times Y \to [0, 1]$ such that $(\lambda \times \mu)(x, y) = \min\{\lambda(x), \mu(y)\}$ for every $(x, y) \in X \times Y$. A fuzzy binary relation $R_\lambda$ on a set X is defined as a fuzzy subset of $X \times X$.*

*The composition of two fuzzy relations $R_\lambda$ and $R_\mu$ is defined by $(R_\lambda \circ R_\mu)(x, y) = \sup_{t \in X} \{\min R_\lambda(x, t), R_\mu(t, y)\}$, for every $x, y \in X$.*

**DEFINITION 1.1.12:** *Let $R_\lambda$ be a fuzzy binary relation on a set X. A fuzzy subset $\mu$ of the set X is said to be a pre class of $R_\lambda$ if $\min\{\mu(x), \mu(y)\} \leq R_\lambda(x, y)$ for every $x, y \in X$.*

*A fuzzy binary relation $R_\lambda$ on a set X is said to be a similarity relation on the set X if it is reflexive, symmetric and transitive that is, for every x, y, z ∈ X.*

$$R_\lambda (x, x) = 1$$
$$R_\lambda (x, y) = R_\lambda (y, x)$$
$$\min \{ R_\lambda (x, y), R_\lambda (y, z)\} \leq R_\lambda (x, z).$$

*Let $\mu$ be a fuzzy subset of a set X. If $\mu (x) = 0$ for every x ∈ X then we call $\mu$ as empty fuzzy set and denote it by $\phi_X$. If $\mu (x) = 1$ for every x ∈ X then we call $\mu$ as whole fuzzy set and denote it by $1_X$.*

**DEFINITION 1.1.13:** *A fuzzy binary relation S on X is said to be a similarity relation on X if it is reflexive, symmetric and transitive i.e.*

$$S (x, x) = 1.$$
$$S (x, y) = S (y, x).$$
$$S (x, y) \wedge S (y, z) \leq S (x, z) \text{ for all } x, y, z \text{ in } X.$$

For more about fuzzy sets please refer [17, 26, 59, 144].

## 1.2 Groups and fuzzy subgroups

Rosenfield [112] introduced the notion of fuzzy group and showed that many group theory results can be extended in an elementary manner to develop the theory of fuzzy group. The underlying logic of the theory of fuzzy group is to provide a strict fuzzy algebraic structure where level subset of a fuzzy group of a group G is a subgroup of the group. [14, 15] reduced fuzzy subgroup of a group using the general t-norm. However, [112] used the t-norm 'min' in his definition of fuzzy subgroup of a group. Fuzzy groups are further investigated by [32, 33] who mainly studied about the level subgroups of a fuzzy subgroup. [109] analyzed this level subgroups of a fuzzy subgroup in more detail and investigated whether the family of level subgroups of a fuzzy subgroup, determine the fuzzy subgroup uniquely or not. The concepts of fuzzy normal subgroup and fuzzy coset were introduced by [98]. For more about fuzzy groups please refer [2, 5, 14, 16, 30, 32, 55, 73, 83, 85, 86, 89, 93, 109, 112, 136, 137, 138, 139].

**DEFINITION 1.2.1:** *Let G be a group. A fuzzy subset $\mu$ of a group G is called a fuzzy subgroup of the group G if*

   i. $\mu(xy) \geq \min \{ \mu (x), \mu (y)\}$ *for every x, y ∈ G and*
   ii. $\mu(x^{-1}) = \mu (x)$ *for every x ∈ G.*

**DEFINITION 1.2.2:** *Let G be a group. A fuzzy subgroup A of G is called normal if $A(x) = A(y^{-1} x y)$ for all x, y ∈ G.*

**DEFINITION 1.2.3:** *Let A be a fuzzy subset of S. For t ∈ [0, 1] the set $A_t = \{ s ∈ S / A(x) \geq t\}$ is called a level subset of the fuzzy subset A.*

In consequence of the level subset we have the following theorem:

**THEOREM 1.2.1:** *Let G be a group and A be a fuzzy subgroup of G. Then the level subsets $A_t$, for $t \in [0, 1]$, $t \leq A(e)$ is a subgroup of G, where e is the identity of G.*

*Proof*: Direct, refer [16].

**THEOREM 1.2.2:** *A fuzzy subset $\mu$ of a group G is a fuzzy subgroup of the group G if and only if $\mu(xy^{-1}) \geq \min\{\mu(x), \mu(y)\}$ for every $x, y \in G$.*

*Proof*: Left for the reader as it is a direct consequence of the definition.

**THEOREM 1.2.3:** *Let $\mu$ be a fuzzy subset of a group G. Then $\mu$ is a fuzzy subgroup of G if and only if $G_\mu^t$ is a subgroup (called level subgroup) of the group G for every $t \in [0, \mu(e)]$, where e is the identity element of the group G.*

*Proof*: Left as an exercise for the reader to prove.

**DEFINITION 1.2.4:** *A fuzzy subgroup $\mu$ of a group G is called improper if $\mu$ is constant on the group G, otherwise $\mu$ is termed as proper.*

**DEFINITION 1.2.5:** *We can define a fuzzy subgroup $\mu$ of a group G to be fuzzy normal subgroup of a group G if $\mu(xy) = \mu(yx)$ for every $x, y \in G$. This is just an equivalent formation of the normal fuzzy subgroup. Let $\mu$ be a fuzzy normal subgroup of a group G. For $t \in [0, 1]$, the set $\mu_t = \{(x, y) \in G \times G / \mu(xy^{-1}) \geq t\}$ is called the t-level relation of $\mu$. For the fuzzy normal subgroup $\mu$ of G and for $t \in [0, 1]$, $\mu_t$ is a congruence relation on the group G.*

In view of all these the reader is expected to prove the following theorem:

**THEOREM 1.2.4:** *Let $\mu$ be a fuzzy subgroup of a group G and $x \in G$. Then $\mu(xy) = \mu(y)$ for every $y \in G$ if and only if $\mu(x) = \mu(e)$.*

**DEFINITION 1.2.6:** *Let $\mu$ be a fuzzy subgroup of a group G. For any $a \in G$, $a\mu$ defined by $(a\mu)x = \mu(a^{-1}x)$ for every $x \in G$ is called the fuzzy coset of the group G determined by a and $\mu$.*

The reader is expected to prove the following.

**THEOREM 1.2.5:** *Let $\mu$ be a fuzzy subgroup of a group G. Then $xG_\mu^t = G_{x\mu}^t$ for every $x \in G$ and $t \in [0, 1]$.*

We now define the order of the fuzzy subgroup $\mu$ of a group G.

**DEFINITION 1.2.7:** *Let $\mu$ be a fuzzy subgroup of a group G and $H = \{x \in G / \mu(x) = \mu(e)\}$ then $o(\mu)$, order of $\mu$ is defined as $o(\mu) = o(H)$.*

**THEOREM 1.2.6:** *Any subgroup H of a group G can be realised as a level subgroup of some fuzzy subgroup of G.*

The proof is left as an exercise to the reader. Some of the characterization about standard groups in fuzzy terms are given. The proof of all these theorems are left for the reader to refer and obtain them on their own.

**THEOREM 1.2.7:** *G is a Dedekind group if and only if every fuzzy subgroup of G is normal.*

By a Dedekind group we mean a group, which is abelian or Hamiltonian. (A group G is Hamiltonian if every subgroup of G is normal)

**THEOREM 1.2.8:** *Let G be a cyclic group of prime order. Then there exists a fuzzy subgroup A of G such that $A(e) = t_o$ and $A(x) = t_1$ for all $x \neq e$ in G and $t_o > t_1$.*

**THEOREM 1.2.9:** *Let G be a finite group of order n and A be a fuzzy subgroup of G. Let $Im(A) = \{t_i / A(x) = t_i$ for some $x \in G\}$. Then $\{A_{t_i}\}$ are the only level subgroups of A.*

Now we give more properties about fuzzy subgroups of a cyclic group.

**THEOREM [16]:** *Let G be a group of prime power order. Then G is cyclic if and only if there exists a fuzzy subgroup A of G such that for $x, y \in G$,*

  i.   *if $A(x) = A(y)$ then $\langle x \rangle = \langle y \rangle$.*
  ii.  *if $A(x) > A(y)$ then $\langle x \rangle \subset \langle y \rangle$.*

**THEOREM [16]:** *Let G be a group of square free order. Let A be a normal fuzzy subgroup of G. Then for $x, y \in G$,*

  i.   *if $o(x) \mid o(y)$ then $A(y) \leq A(x)$.*
  ii.  *if $o(x) = o(y)$ then $A(y) = A(x)$.*

**THEOREM [16]:** *Let G be a group of order $p_1, p_2, \ldots, p_r$ where the $p_i$'s are primes but not necessarily distinct. Then G is solvable if and only if there exists a fuzzy subgroup A of G such that $A_{t_0}, A_{t_1}, \ldots, A_{t_r}$ are the only level subgroups of A, $Im(A) = \{t_0, t_1, \ldots, t_r\}$, $t_0 > t_1 > \ldots > t_r$ and the level subgroups form a composition chain.*

**THEOREM [16]:** *Suppose that G is a finite group and that G has a composition chain $\langle e \rangle = A_0 \subset A_1 \subset \ldots \subset A_r = G$ where $A_i / A_{i-1}$ is cyclic of prime order, $i = 1, 2, \ldots, r$. Then there exists a composition chain of level subgroups of some fuzzy subgroup A of G and this composition chain is equivalent to $\langle e \rangle = A_0 \subset A_1 \subset \ldots \subset A_r = G$.*

The proof of these results can be had from [16].

**DEFINITION [98]:** *Let λ and μ be two fuzzy subgroups of a group G. Then λ and μ are said to be conjugate fuzzy subgroups of G if for some g ∈ G, λ(x) = μ ($g^{-1}xg$) for every x ∈ G.*

**THEOREM [139]:** *If λ and μ are conjugate fuzzy subgroups of the group G then o(λ) = o(μ).*

*Proof:* Refer [139] for proof.

Mukherjee and Bhattacharya [98] introduced fuzzy right coset and fuzzy left coset of a group G. Here we introduce the notion of fuzzy middle coset of a group G mainly to prove that o(α μ $α^{-1}$) = o (μ) for any fuzzy subgroup μ of the group G and α ∈ G.

**DEFINITION 1.2.8:** *Let μ be a fuzzy subgroup of a group G. Then for any a, b ∈ G a fuzzy middle coset a μ b of the group G is defined by (a μ b) (x) = μ ($a^{-1}$ x $b^{-1}$) for every x ∈ G.*

The following example from [139] is interesting which explains the notion of fuzzy middle coset.

*Example 1.2.1:* Let G = {1, –1, i, –i} be the group, with respect to the usual multiplication.

Define μ: G → [0, 1] by

$$\mu(x) = \begin{cases} 1 & \text{if } x = 1 \\ 0.5 & \text{if } x = -1 \\ 0 & \text{if } x = i, -i. \end{cases}$$

Clearly μ is a fuzzy subgroup of the group G. A fuzzy middle coset a μ b is calculated and given by

$$(a\mu b)(x) = \begin{cases} 0 & \text{if } x = 1, -1 \\ 0.5 & \text{if } x = -i \\ 1 & \text{if } x = i \end{cases}$$

for all a = –1 and b = – i.

*Example 1.2.2:* Consider the infinite group Z = {0, 1, –1, 2, –2, …} with respect to usual addition. Clearly 2Z is a proper subgroup of Z.

Define μ: Z → [0, 1] by

$$\mu(x) = \begin{cases} 0.9 & \text{if } x \in 2Z \\ 0.8 & \text{if } x \in 2Z+1. \end{cases}$$

It is easy to verify that μ is a fuzzy subgroup of the group Z. For any a ∈ 2Z and b ∈ 2Z + 1 the fuzzy middle coset a μ b is given by

$$(a\mu b)(x) = \begin{cases} 0.8 & \text{if } x \in 2Z \\ 0.9 & \text{if } x \in 2Z+1. \end{cases}$$

Hence it can be verified that this fuzzy middle coset aμb in not a fuzzy subgroup of Z.

We have the following theorem.

**THEOREM 1.2.10:** *If μ is a fuzzy subgroup of a group G then for any a ∈ G the fuzzy middle coset aμ $a^{-1}$ of the group G is also a fuzzy subgroup of the group G.*

*Proof:* Refer [137].

**THEOREM 1.2.11:** *Let μ be any fuzzy subgroup of a group G and a μ $a^{-1}$ be a fuzzy middle coset of the group G then o (a μ $a^{-1}$) = o(μ) for any a ∈ G.*

*Proof:* Let μ be a fuzzy subgroup of a group G and a ∈ G. By Theorem 1.2.10 the fuzzy middle coset aμ$a^{-1}$ is a fuzzy subgroup of the group G. Further by the definition of a fuzzy middle coset of the group G we have (a μ $a^{-1}$) (x) = μ ($a^{-1}$xa) for every x ∈ G. Hence for any a ∈ G, μ and aμ$a^{-1}$ are conjugate fuzzy subgroups of the group G as there exists a ∈ G such that (aμ$a^{-1}$)(x) = μ($a^{-1}$xa) for every x ∈ G. By using earlier theorems which states o(aμ$a^{-1}$) = o(μ) for any a ∈ G.

For the sake of simplicity and better understanding we give the following example.

*Example 1.2.3:* Let G = $S_3$ the symmetric group of degree 3 and $p_1$, $p_2$, $p_3$ ∈ [0, 1] such that $p_1 \geq p_2 \geq p_3$.

Define μ : G → [0, 1] by

$$\mu(x) = \begin{cases} p_1 & \text{if } x = e \\ p_2 & \text{if } x = (12) \\ p_3 & \text{otherwise.} \end{cases}$$

Clearly μ is a fuzzy subgroup of a group G and o(μ) = number of elements of the set {x ∈ G | μ (x) = μ(e)} = number of elements of the set {e} = 1. Now we can evaluate a μ $a^{-1}$ for every a ∈ G as follows:

For a = e we have a μ $a^{-1}$ = μ. Hence o (a μ $a^{-1}$) = o (μ) = 1.

For a = (12) we have

$$(a\mu a^{-1})(x) = \begin{cases} p_1 & \text{if } x = e \\ p_2 & \text{if } x = (12) \\ p_3 & \text{otherwise.} \end{cases}$$

Hence $o(a\mu a^{-1}) = 1$. For the values of $a = (13)$ and $(132)$ we have $a\mu a^{-1}$ to be equal which is given by

$$(a\mu a^{-1})(x) = \begin{cases} p_1 & \text{if } x = e \\ p_2 & \text{if } x = (23) \\ p_3 & \text{otherwise.} \end{cases}$$

Hence $o(a\mu a^{-1}) = 1$ for $a = (13)$ and $(132)$. Now for $a = (23)$ and $a = (123)$ we have $a\mu a^{-1}$ to be equal which is given by

$$(a\mu a^{-1})(x) = \begin{cases} p_1 & \text{if } x = e \\ p_2 & \text{if } x = (13) \\ p_3 & \text{otherwise.} \end{cases}$$

Thus $o(a\mu a^{-1}) = 1$. Hence $o(a\mu a^{-1}) = o(\mu) = 1$ for any $a \in G$.

From this example we see the functions $\mu$ and $a\mu a^{-1}$ are not equal for some $a \in G$. Thus it is interesting to note that if $\mu$ is fuzzy subgroup of an abelian group G then the functions $\mu$ and $a\mu a^{-1}$ are equal for any $a \in G$. However it is important and interesting to note that the converse of the statement is not true. That is if $a\mu a^{-1} = \mu$ for any $a \in G$ can hold good even if G is not abelian. This is evident from the following example.

*Example 1.2.4:* Let $G = S_3$ be the symmetric group of degree 3 and $p_1, p_2, p_3 \in [0, 1]$ be such that $p_1 \geq p_2 \geq p_3$.

Define $\mu: G \to [0, 1]$ by

$$\mu(x) = \begin{cases} p_1 & \text{if } x = e \\ p_2 & \text{if } x = (123) \text{ and } x = (132) \\ p_3 & \text{otherwise.} \end{cases}$$

Clearly $\mu$ is a fuzzy subgroup of G. For any $a \in G$ the fuzzy subgroup $(a\mu a^{-1})$ is given by

$$(a\mu a^{-1})(x) = \begin{cases} p_1 & \text{if } x = e \\ p_2 & \text{if } x = (123) \text{ or } x = (132) \\ p_3 & \text{otherwise.} \end{cases}$$

Thus we have $(a\mu a^{-1})(x) = \mu(x)$ for every $x \in G$. Hence $a\mu a^{-1} = \mu$ for any $a \in G$. Thus the functions $a\mu a^{-1}$ and $\mu$ are identical but G is not an abelian group. It is worthwhile to note that in general o(aμ) is not defined since aμ is not a fuzzy subgroup of the group G. The reader is advised to construct an example to prove the above claim.

**THEOREM 1.2.12:** *Let μ be a fuzzy subgroup of a finite group G then o (μ) | o(G).*

*Proof*: Let μ be a fuzzy subgroup of a finite group with e as its identity element. Clearly H = $\{x \in G \,|\, \mu(x) = \mu(e)\}$ is a subgroup of the group G for H is a t- level subset of the group G where t = μ (e). By Lagranges Theorem o(H) | o(G). Hence by the definition of the order of the fuzzy subgroup of the group G we have o (μ)|o(G).

The following theorem is left as an exercise for the reader to prove.

**THEOREM 1.2.13:** *Let λ and μ be any two improper fuzzy subgroups of a group G. Then λ and μ are conjugate fuzzy subgroups of the group G if and only if λ = μ.*

**DEFINITION 1.2.9:** *Let λ and μ be two fuzzy subsets of a group G. We say that λ and μ are conjugate fuzzy subsets of the group G if for some g ∈ G we have λ(x) = μ (g$^{-1}$xg) for every x ∈ G.*

We now give a relation about conjugate fuzzy subsets of a group G.

**THEOREM 1.2.14:** *Let λ and μ be two fuzzy subsets of an abelian group G. Then λ and μ are conjugate fuzzy subsets of the group G if and only if λ = μ.*

*Proof:* Let λ and μ be conjugate fuzzy subsets of group G then for some g ∈ G we have

$$\begin{aligned}\lambda(x) &= \mu(g^{-1}xg) \text{ for every } x \in G \\ &= \mu(g^{-1}gx) \text{ for every } x \in G \\ &= \mu(x) \text{ for every } x \in G.\end{aligned}$$

Hence λ = μ.

Conversely if λ = μ then for the identity element e of group G, we have λ(x) = μ(e$^{-1}$xe) for every x ∈ G. Hence λ and μ are conjugate fuzzy subsets of the group G.

The reader is requested to prove the following theorem as a matter of routine.

**THEOREM 1.2.15**: *Let λ be a fuzzy subgroup of a group G and μ be a fuzzy subset of the group G. If λ and μ are conjugate fuzzy subsets of the group G then μ is a fuzzy subgroup of the group G.*

The reader is requested verify if λ , μ : S$_3$ → [0, 1] as

$$\lambda(x) = \begin{cases} 0.5 & \text{if } x = e \\ 0.4 & \text{if } x = (123) \ \& \ x = (132) \\ 0.3 & \text{otherwise} \end{cases}$$

and

$$\mu(x) = \begin{cases} 0.6 & \text{if } x = e \\ 0.5 & \text{if } x = (23) \\ 0.3 & \text{otherwise} \end{cases}$$

where e is the identity element of $S_3$, to prove $\lambda$ and $\mu$ are not conjugate fuzzy subsets of the group $S_3$.

Now we proceed on to recall the notions of conjugate fuzzy relations of a group and the generalized conjugate fuzzy relations on a group.

**DEFINITION 1.2.10:** *Let $R_\lambda$ and $R_\mu$ be any two fuzzy relations on a group G. Then $R_\lambda$ and $R_\mu$ are said to be conjugate fuzzy relations on a group G if there exists $(g_1, g_2) \in G \times G$ such that $R_\lambda (x, y) = R_\mu = (g_1^{-1} x g_1, g_2^{-1} y g_2)$ for every $(x, y) \in G \times G$.*

**DEFINITION 1.2.11:** *Let $R_\lambda$ and $R_\mu$ be any two fuzzy relation on a group G. Then $R_\lambda$ and $R_\mu$ are said to be generalized conjugate fuzzy relations on the group G if there exists $g \in G$ such that $R_\lambda (x, y) = R_\mu (g^{-1} x g, g^{-1} y g)$ for every $(x, y) \in G \times G$.*

**THEOREM 1.2.16:** *Let $R_\lambda$ and $R_\mu$ be any two fuzzy relations on a group G. If $R_\lambda$ and $R_\mu$ are generalized conjugate fuzzy relations on the group G then $R_\lambda$ and $R_\mu$ are conjugate fuzzy relations on the group G.*

*Proof*: Let $R_\lambda$ and $R_\mu$ be generalized conjugate fuzzy relations on the group G. Then there exists $g \in G$ such that $R_\lambda (x, y) = R_\mu (g^{-1} x g, g^{-1} y g)$ for every $(x, y) \in G \times G$. Now choose $g_1 = g_2 = g$. Then for $(g_1, g_2) \in G \times G$ we have $R_\lambda (x, y) = R_\mu \left( g_1^{-1} x g_1, g_2^{-1} y g_2 \right)$ for every $(x, y) \in G \times G$. Thus $R_\lambda$ and $R_\mu$ are conjugate fuzzy relations on the group G.

The reader can prove that the converse of the above theorem in general is not true.

**THEOREM 1.2.17:** *Let $\mu$ be a fuzzy normal subgroup of a group G. Then for any $g \in G$ we have $\mu (gxg^{-1}) = \mu (g^{-1} xg)$ for every $x \in G$.*

*Proof*: Straightforward and hence left for the reader to prove.

**THEOREM 1.2.18:** *Let $\lambda$ and $\mu$ be conjugate fuzzy subgroups of a group G. Then*

  i.  *$\lambda \times \mu$ and $\mu \times \lambda$ are conjugate fuzzy relations on the group G and*
  ii. *$\lambda \times \mu$ and $\mu \times \lambda$ are generalized conjugate fuzzy relations on the group G only when at least one of $\lambda$ or $\mu$ is a fuzzy normal subgroup of G.*

*Proof*: The proof can be obtained as a matter of routine. The interested reader can refer [139].

Now we obtain a condition for a fuzzy relation to be a similarity relation on G.

**THEOREM 1.2.19:** *Let $R_\lambda$ be a similarity relation on a group G and $R_\mu$ be a fuzzy relation on the group G. If $R_\lambda$ and $R_\mu$ are generalized conjugate fuzzy relations on the group G then $R_\mu$ is a similarity relation on the group G.*

*Proof*: Refer [139].

Now we define some properties on fuzzy symmetric groups.

**DEFINITION [55]:** *Let $S_n$ denote the symmetric group on {1, 2, …, n}. Then we have the following:*

   i.  *Let $F(S_n)$ denote the set of all fuzzy subgroups of $S_n$.*
   ii. *Let $f \in F(S_n)$ then $Im\,f = \{f(x) \mid x \in S_n\}$.*
   iii. *Let $f, g \in F(S_n)$. If $|Im(f)| < |Im(g)|$ then we write $f < g$. By this rule we define max $F(S_n)$.*
   iv. *Let f be a fuzzy subgroup of $S_n$. If $f = max\,F(S_n)$ then we say that f is a fuzzy symmetric subgroup of $S_n$.*

**THEOREM 1.2.20:** *Let f be a fuzzy symmetric subgroup of the symmetric group $S_3$ then $o(Im\,f) = 3$.*

*Proof:* Please refer [139].

Here we introduce a new concept called co fuzzy symmetric group which is a generalization of the fuzzy symmetric group.

**DEFINITION [139]:** *Let $G(S_n) = \{ g \mid g$ is a fuzzy subgroup of $S_n$ and $g(C(\Pi))$ is a constant for every $\Pi \in S_n\}$ where $C(\Pi)$ is the conjugacy class of $S_n$ containing $\Pi$, which denotes the set of all $y \in S_n$ such that $y = x\Pi x^{-1}$ for $x \in S_n$. If $g = max\,G(S_n)$ then we call g as co-fuzzy symmetric subgroup of $S_n$.*

For better understanding of the definition we illustrate it by the following example.

*Example 1.2.5:* Let $G = S_3$ be the symmetric group of degree 3.

Define g: G → [0 1] as follows:

$$g(x) = \begin{cases} 1 & \text{if } x = e \\ 0.5 & \text{if } x = (123), (132) \\ 0 & \text{otherwise} \end{cases}$$

where e is the identity element of $S_3$. It can be easily verified that all level subsets of g are {e} {e, (123), (132)} and $S_3$. All these level subsets are subgroups of $S_3$, hence g is a fuzzy subgroup of $S_3$. Further g (C(Π)) is constant for every $Π \in S_3$ and o (Im (g)) ≥ o (Im g (μ)) for every subgroup μ of the symmetric group $S_3$. Hence g is a co-fuzzy symmetric subgroup of $S_3$.

Now we proceed on to prove the following theorem using results of [55].

**THEOREM 1.2.21:**

i. *If g is a co-fuzzy symmetric subgroup of the symmetric group $S_3$ then o(Im(g)) = 3.*
ii. *If g is a co-fuzzy symmetric subgroup of $S_4$ then o (Im (g)) = 4 and*
iii. *If g is a co-fuzzy symmetric subgroup of $S_n$ (n ≥ 5) then o (Im (g)) = 3.*

*Proof*: The proof follows verbatim from [55] when the definition of fuzzy symmetric group is replaced by the co-fuzzy symmetric group.

**THEOREM 1.2.22:** *Every co fuzzy symmetric subgroup of a symmetric group $S_n$ is a fuzzy symmetric subgroup of the symmetric group $S_n$.*

*Proof:* Follows from the very definitions of fuzzy symmetric subgroup and co fuzzy symmetric subgroup.

**THEOREM 1.2.23:** *Every fuzzy symmetric subgroup of a symmetric group $S_n$ need not in general be a co-fuzzy symmetric subgroup of $S_n$.*

*Proof:* By an example. Choose $p_1, p_2, p_3 \in [0, 1]$ such that $1 \geq p_1 \geq p_2 \geq p_3 \geq 0$.

Define f : $S_3 \to [0\ 1]$ by

$$f(x) = \begin{cases} p_1 & \text{if } x = e \\ p_2 & \text{if } x = (12) \\ p_3 & \text{otherwise} \end{cases}$$

It can be easily checked that f is a fuzzy subgroup of $S_3$ as all the level subsets of f are subgroups of $S_3$. Further o(Im (f)) = 3 ≥ o(Im (μ)) for every fuzzy subgroup μ of the symmetric group $S_3$. Hence f is a fuzzy symmetric subgroup of $S_3$ but f(12) ≠ f(13) in this example. By the definition of co fuzzy symmetric subgroup it is clear that f is not a co fuzzy symmetric subgroup of $S_3$. Hence the claim.

Now we proceed on to recall yet a new notion called pseudo fuzzy cosets and pseudo fuzzy double cosets of a fuzzy subset or a fuzzy subgroup. [98] has defined fuzzy coset as follows:

**DEFINITION [98]:** *Let μ be a fuzzy subgroup of a group G. For any a ∈ G, a μ defined by (a μ) (x) = μ($a^{-1}$x) for every x ∈ G is called a fuzzy coset of μ.*

One of the major marked difference between the cosets in fuzzy subgroup and a group is "any two fuzzy cosets of a fuzzy subgroup μ of a group G are either identical or disjoint" is not true.

This is established by the following example:

*Example 1.2.6:* Let G = { ±1, ± i } be the group with respect to multiplication.

Define μ: G → [0, 1] as follows:

$$\mu(x) = \begin{cases} \frac{1}{2} & \text{if } x = -1 \\ 1 & \text{if } x = 1 \\ \frac{1}{4} & \text{if } x = i, -i \end{cases}$$

The fuzzy cosets iμ and − iμ of μ are calculated as follows:

$$i\mu(x) = \begin{cases} \frac{1}{4} & \text{if } x = 1, -1 \\ 1 & \text{if } x = i \\ \frac{1}{2} & \text{if } x = -i \end{cases}$$

and

$$(-i\mu)(x) = \begin{cases} \frac{1}{4} & \text{if } x = 1, -1 \\ 1 & \text{if } x = -i \\ \frac{1}{2} & \text{if } x = i \end{cases}$$

It is easy to see that these fuzzy cosets iμ and –iμ are neither identical nor disjoint. For (iμ)(i) ≠ (–iμ) (i) implies iμ and –iμ are not identical and (iμ)(1) = (–iμ)(1) implies iμ and –iμ are not disjoint. Hence the claim.

Now we proceed on to recall the notion of pseudo fuzzy coset.

**DEFINITION 1.2.12:** *Let μ be a fuzzy subgroup of a group G and a ∈ G. Then the pseudo fuzzy coset $(a\mu)^P$ is defined by $((a\mu)^P)(x) = p(a) \mu(x)$ for every x ∈ G and for some p ∈ P.*

*Example 1.2.7:* Let G = {1, ω, ω²} be a group with respect to multiplication, where ω denotes the cube root of unity. Define μ: G → [0, 1] by

$$\mu(x) = \begin{cases} 0.6 & \text{if } x = 1 \\ 0.4 & \text{if } x = \omega, \omega^2 \end{cases}$$

It is easily checked the pseudo fuzzy coset $(a\mu)^P$ for $p(x) = 0.2$ for every $x \in G$ to be equal to 0.12 if $x = 1$ and 0.08 if $x = \omega, \omega^2$.

We define positive fuzzy subgroup.

**DEFINITION 1.2.13:** *A fuzzy subgroup $\mu$ of a group G is said to be a positive fuzzy subgroup of G if $\mu$ is a positive fuzzy subset of the group G.*

**THEOREM 1.2.24:** *Let $\mu$ be a positive fuzzy subgroup of a group G then any two pseudo fuzzy cosets of $\mu$ are either identical or disjoint.*

*Proof:* Refer [137]. As the proof is lengthy and as the main motivation of the book is to introduce Smarandache fuzzy concepts we expect the reader to be well versed in fuzzy algebra, we request the reader to supply the proof.

Now we prove the following interesting theorem.

**THEOREM 1.2.25:** *Let $\mu$ be a fuzzy subgroup of a group G then the pseudo fuzzy coset $(a\mu)^P$ is a fuzzy subgroup of the group G for every $a \in G$.*

*Proof:* Let $\mu$ be a fuzzy subgroup of a group G. For every x, y in G we have

$$\begin{aligned}
(a\mu)^P(xy^{-1}) &= p(a)\, \mu(xy^{-1}) \\
&\geq p(a) \min\{\mu(x), \mu(y)\} \\
&= \min\{p(a)\mu(x), p(a), \mu(y)\} \\
&= \min\{(a\mu)^P(x), (a\mu)^P(y)\}.
\end{aligned}$$

That is $(a\mu)^P(xy^{-1}) \geq \min\{(a\mu)^P(x), (a\mu)^P(y)\}$ for every x, y $\in$ G. This proves that $(a\mu)^P$ is a fuzzy subgroup of the group G. We illustrate this by the following example:

*Example 1.2.8:* Let G be the Klein four group. Then G = {e, a, b, ab} where $a^2 = e = b^2$, ab = ba and e the identity element of G.

Define $\mu: G \to [0, 1]$ as follows

$$\mu(x) = \begin{cases} \dfrac{1}{2} & \text{if } x = a \\ 1 & \text{if } x = e \\ \dfrac{1}{4} & \text{if } x = b, ab \end{cases}$$

Take the positive fuzzy subset p as follows:

$$p(x) = \begin{cases} 1 & \text{if } x = e \\ \dfrac{1}{2} & \text{if } x = a \\ \dfrac{1}{3} & \text{if } x = b \\ \dfrac{1}{4} & \text{if } x = ab \end{cases}$$

Now we calculate the pseudo fuzzy cosets of $\mu$. For the identity element e of the group G we have $(e\mu)^P = \mu$.

$$(a\mu)^P (x) = \begin{cases} \dfrac{1}{2} & \text{if } x = e \\ \dfrac{1}{4} & \text{if } x = a \\ \dfrac{1}{8} & \text{if } x = b, ab \end{cases}$$

$$(b\mu)^P (x) = \begin{cases} \dfrac{1}{3} & \text{if } x = e \\ \dfrac{1}{6} & \text{if } x = a \\ \dfrac{1}{12} & \text{if } x = b, ab \end{cases}$$

and

$$((ab)\mu)^P)(x) = \begin{cases} \dfrac{1}{4} & \text{if } x = e \\ \dfrac{1}{8} & \text{if } x = a \\ \dfrac{1}{16} & \text{if } x = b, ab \end{cases}$$

It is easy to check that all the above pseudo fuzzy cosets of $\mu$ are fuzzy subgroups of G. As there is no book on fuzzy algebraic theory dealing with all these concepts we have felt it essential to give proofs and examples atleast in few cases.

**THEOREM 1.2.26:** *Let $\mu$ be a fuzzy subgroup of a finite group G and $t \in [0, 1]$ then $o(G^t_{(a\mu)P}) \leq o(G^t_\mu) = o(aG^t_\mu)$ for any $a \in G$.*

*Proof:* The proof is left as an exercise for the reader to prove.

**THEOREM 1.2.27:** *A fuzzy subgroup $\mu$ of a group G is normalized if and only if $\mu(e) = 1$, where e is the identity element of the group G.*

*Proof:* If μ is normalized then there exists x ∈ G such that μ(x) = 1, but by properties of a fuzzy subgroup μ of the group G, μ(x) ≤ μ(e) for every x ∈ G. Since μ(x) = 1 and μ(e) ≥ μ(x) we have μ(e) ≥1. But μ(e) ≤ 1. Hence μ(e) = 1. Conversely if μ(e) = 1 then by the very definition of normalized fuzzy subset μ is normalized.

The proof of the following theorem is left as an exercise for the reader, which can be proved as a matter of routine. The only notion which we use in the theorem is the notion of pre class of a fuzzy binary relation $R_\mu$. Let μ be a fuzzy subgroup of a group G. Now we know that a fuzzy subset μ of a set X is said to be a pre class of a fuzzy binary relation $R_\mu$ on the set X if min {μ(x), μ(y)} ≤ $R_\mu$(x, y) for every x, y ∈ X.

**THEOREM 1.2.28:** *Let μ be a fuzzy subgroup of a group G and $R_\mu : G \times G \to [0\ 1]$ be given by $R_\mu(x, y) = \mu(xy^{-1})$ for every x, y ∈ G. Then*

  i. *$R_\mu$ is a similarity relation on the group G only when μ is normalized and*
  ii. *μ is a pre class of $R_\mu$ and in general the pseudo fuzzy coset $(a\mu)^P$ is a pre class of $R_\mu$ for any a ∈ G.*

**DEFINITION 1.2.14:** *Let μ be a fuzzy subset of a non-empty set X and a ∈ X. We define the pseudo fuzzy coset $(a\mu)^P$ for some p ∈ P by $(a\mu)^P(x) = p(a)\mu(x)$ for every x ∈ X.*

*Example 1.2.9:* Let X = {1, 2, 3, …, n} and μ: X → [0, 1] is defined by $\mu(x) = \frac{1}{x}$ for every x ∈ X. Then the pseudo fuzzy coset $(\alpha\mu)^P$: X → [0, 1] is computed in the following manner by taking $p(x) = \frac{1}{2x}$ for every x ∈ X; $(\alpha\mu)^P(x) = \frac{1}{2x^2}$ for every x ∈ X.

**THEOREM 1.2.29:** *Let λ and μ be any two fuzzy subsets of a set X. Then for a ∈ X $(a\mu)^P \subset (a\lambda)^P$ if and only if $\mu \subseteq \lambda$.*

*Proof:* Left as an exercise for the reader.

Now we proceed on to define the fuzzy partition of a fuzzy subset.

**DEFINITION 1.2.15:** *Let μ be a fuzzy subset of a set X. Then Σ = {λ: λ is a fuzzy subset of a set X and $\lambda \subseteq \mu$} is said to be a fuzzy partition of μ if*

  i. $\bigcup_{\lambda \in \Sigma} \lambda = \mu$ *and*
  ii. *any two members of Σ are either identical or disjoint*

However we illustrate by an example.

***Example 1.2.10:*** Let $X = N$ be the set of all natural numbers and $\mu$ be defined by $\mu(x) = \frac{1}{x}$ for every $x \in X$. Now consider the collection of fuzzy subsets of $X$ which is given by $\{\mu_i\}_{i=1}^{\infty}$ where $\mu_i$'s are such that $\mu_I(x) = (1 - \frac{1}{i})\frac{1}{x}$ for every $x \in X$.

For $x \in X$ we have

$$\left(\bigcup_{i=1}^{\infty} \mu_i\right)(x) = \sup\left\{\left(1 - \frac{1}{i}\right)\frac{1}{x}\right\} = \frac{1}{x}\left(\text{as } \frac{1}{i} \to 0, i \to \infty\right) = \mu(x).$$

Hence

$$\left(\bigcup_{i=1}^{\infty} \mu_i\right)(x) = \mu(x)$$

for every $x \in X$. That is

$$\left(\bigcup_{i=1}^{\infty} \mu_i\right) = \mu.$$

If $i \neq j$ then it is easy to verify that

$$\left(1 - \frac{1}{i}\right)\frac{1}{x} \neq \left(1 - \frac{1}{j}\right)\frac{1}{x}$$

for every $x \in X$. This proves $\mu_i(x) \neq \mu_j(x)$ for every $x \in X$. Hence $\mu_I$ and $\mu_j$ are disjoint.

Hence $\{\mu_i\}_{i=1}^{\infty}$ is a fuzzy partition of $\mu$.

The following theorem is left for the reader, however the proof can be found in [89].

**THEOREM 1.2.30:** *Let $\mu$ be a positive fuzzy subset of a set $X$ then*

  i. *any two pseudo fuzzy cosets of $\mu$ are either identical or disjoint.*
  ii. $\bigcup_{p \in P}((a\mu)^P) = \mu$,
  iii. $\bigcup_{a \in X}((a\mu)^P) \subseteq \bigcup_{p \in P}((a\mu)^P)$ *and the equality holds good if and only if p is normal,*
  iv. *The collection $\{(a\mu)^P \mid a \in X\}$ is a fuzzy partition of $\mu$ if and only if p is normal.*

The following theorem is yet another piece of result on pre class

**THEOREM 1.2.31:** *Let $\mu$ be a fuzzy subgroup of a group $G$ and $R_\mu : G \times G \to [0, 1]$ be given by $R_\mu (x, y) = \mu (xy^{-1})$ for every $x, y \in G$. If $\lambda$ is a fuzzy subset of the group $G$ such that $\lambda \subseteq \mu$ then $(a\lambda)^P$ is pre class of $R_\mu$ for any $a \in G$.*

*Proof*: Let $\mu$ be a fuzzy subgroup of a group $G$, $a \in G$ and $\lambda$ be a fuzzy subset of the group $G$ such that $\lambda \subseteq \mu$.

For $x, y \in G$ we have

$$
\begin{aligned}
\min \{((a\lambda)^P)(x), ((a\lambda)^P)y\} &= \min \{ p(a) \lambda(x), p(a) \lambda(y)\} \\
&\leq \min \{p(a) \mu(x), p(a) \mu(y)\} \text{ (since } \lambda \subseteq \mu.) \\
&= p(a) \min \{\mu(x), \mu(y)\} \\
&\leq 1 \bullet \mu (xy^{-1}) \text{ (since } p(a) \leq 1) \\
&= R_\mu (x, y).
\end{aligned}
$$

That is $\min \{((a\lambda)^P)(x), ((a\lambda)^P)y\} \leq R_\mu (x, y)$ for every $x, y \in G$. Hence $(a \lambda)^P$ is a pre class of $R_\mu$ for any $a \in G$.

Now we proceed on to define the notion of pseudo fuzzy double cosets.

**DEFINITION 1.2.16:** *Let $\mu$ and $\lambda$ be any two fuzzy subsets of a set $X$ and $p \in P$. The pseudo fuzzy double coset $(\mu x \lambda)^p$ is defined by $(\mu x \lambda)^P = (x\mu)^P \cap (x\lambda)^P$ for $x \in X$.*

We illustrate this concept by the following example:

*Example 1.2.11:* Let $X = \{1, 2, 3\}$ be a set. Take $\lambda$ and $\mu$ to be any two fuzzy subsets of $X$ given by $\lambda(1) = 0.2, \lambda(2) = 0.8, \lambda(3) = 0.4$. $\mu(1) = 0.5$ $\mu(2) = 0.6$ and $\mu(3) = 0.7$. Then for a positive fuzzy subset $p$ such that $p(1) = p(2) = p(3) = 0.1$, we calculate the pseudo fuzzy double coset $(\mu \times \lambda)^P$ and this is given below.

$$(\mu x \lambda)^P (y) = \begin{cases} 0.02 & \text{if } y = 1 \\ 0.06 & \text{if } y = 2 \\ 0.04 & \text{if } y = 3 \end{cases}$$

The following theorem is left to the reader as the proof can be obtained by a routine calculation.

**THEOREM 1.2.32:** *Let $\lambda$ and $\mu$ be any two positive fuzzy subsets of a set $X$ and $p \in P$. The set of all pseudo fuzzy double cosets $\{(\mu x \lambda)^P | x \in X\}$ is a fuzzy partition of $(\mu \cap \lambda)$ if and only if $p$ is normal.*

It can be easily verified that the intersection of any two similarity relations on a set $X$ is a similarity relation on the set $X$ and on the contrary the union of similarity relations and composition of similarity relations need not in general be similarity relations.

The following theorem is left as an exercise for the reader to prove.

**THEOREM 1.2.33:** *Let $\lambda$ and $\mu$ be any two fuzzy subgroups of a group G and $R_{\mu \cap \lambda}: G \times G \to [0,1]$ be given by $R_{\mu \cap \lambda}(x, y) = (\mu \cap \lambda)(xy^{-1})$ for every $x, y \in G$. Then*

   i. *$R_{\mu \cap \lambda}$ is a similarity relation on the group G only when both $\mu$ and $\lambda$ are normalized.*
   ii. *$(\mu x \lambda)^P$ is a pre class of $R_{\mu \cap \lambda}$ for any $x \in G$ where $p \in P$.*

Consequent of this theorem one can easily prove the following theorem:

**THEOREM 1.2.34:** *Let $\mu$ and $\lambda$ be any two fuzzy subgroups of a group G and $R_{\mu \cap \lambda}: G \times G \to [0, 1]$ be given by $R_{\mu \cap \lambda}(x, y) = (\mu \cap \lambda)(xy^{-1})$ for every $x, y \in G$. If $\eta$ is any fuzzy subset of the group G such that $\eta \subseteq \mu \cap \lambda$ then $\eta$ is a pre class of $R_{\mu \cap \lambda}$.*

*We will show by the following example that $R_{\mu \cap \lambda}$ is not a similarity relation on the group G.*

***Example 1.2.12:*** Let $G = \{1, \omega, \omega^2\}$ be the group with respect to the usual multiplication, where $\omega$ denotes the cube root of unity.

Define $\lambda, \mu: G \to [0, 1]$ by

$$\mu(x) = \begin{cases} 1 & \text{if } x = 1 \\ 0.6 & \text{if } x = \omega \\ 0.5 & \text{if } x = \omega^2 \end{cases}$$

and

$$\mu(x) = \begin{cases} 0.5 & \text{if } x = 1 \\ 0.4 & \text{if } x = \omega \\ 0.3 & \text{if } x = \omega^2 \end{cases}.$$

It can be found that for every $x \in G$. $R_{\mu \cap \lambda}(x, x) = (\mu \cap \lambda)(xx^{-1}) = (\mu \cap \lambda)(1) = 0.5$. Hence $R_{\mu \cap \lambda}$ is not reflexive and hence $R_{\mu \cap \lambda}$ is not a similarity relation on the group G.

In this section we study the concept of fuzzy subgroup using the definition of [70].

**DEFINITION [70]:** *Let G be a group and e denote the identity element of the group G. A fuzzy subset $\mu$ of the group G is called a fuzzy subgroup of group G if*

   i. *$\mu(xy^{-1}) \geq \min \{\mu(x), \mu(y)\}$ for every $x, y \in G$ and*
   ii. *$\mu(e) = 1$.*

Using the results of [70] we give some of the classical results.

**DEFINITION [70]:** *Let $\mu$ be a fuzzy normal subgroup of a group G and $\mu_t$ be a t-level congruence relation of $\mu$ on G. Let A be a non-empty subset of the group G. The congruence class of $\mu_t$ containing the element x of the group G is denoted by $[x]_\mu$.*

*The set $\underline{\mu}_t(A) = \{x \in G \mid [x]_\mu \subseteq A\}$ and $\overline{\mu}_t(A) = \{x \in G \mid [x]_\mu \cap A \neq \phi\}$ are called respectively the lower and upper approximations of the set A with respect to $\mu_t$.*

We give some simple proofs to the results of [70] using the notions of t-level relation and the coset.

**THEOREM (EXISTENCE THEOREM):** *Let $\mu$ be a fuzzy subgroup of a group G. The congruence class $[x]_\mu$ of $\mu_t$ containing the element x of the group G exist only when $\mu$ is a fuzzy normal subgroup of the group G.*

*Proof*: Let $\mu$ be a fuzzy subgroup of a group G. [70] has proved that if $\mu$ is a fuzzy normal subgroup of a group G then the t-level relation $\mu_t$ of $\mu$ is a congruence relation on the group G and hence the congruence class $[x]_\mu$ of $\mu_t$ containing the element x of the group G exist.

Now we prove that for the existence of the congruence class $[x]_\mu$ we must have the fuzzy subgroup $\mu$ of the group G to be fuzzy normal subgroup of group G. That is if $\mu$ is not a fuzzy normal subgroup of the group G then the congruence class $[x]_\mu$ of $\mu_t$ containing the element x of the group G does not exist.

To prove this, consider the alternating group $A_4$.

Choose $p_1, p_2, p_3 \in [0, 1]$ such that $1 > p_1 > p_2 > p_3 \geq 0$.

Define $\mu : A_4 \to [0, 1]$ by

$$\mu(x) = \begin{cases} 1 & \text{if } x = e \\ p_1 & \text{if } x = (1\,2)(3\,4) \\ p_2 & \text{if } x = (1\,4)(2\,3), (1\,3)(2\,4) \\ p_3 & \text{otherwise} \end{cases}$$

where e is the identity element of $A_4$.

The t-level subsets of $\mu$ are given by $\{e\}$, $\{e, (1\,2)(3\,4)\}$, $\{e, (1\,2)(3\,4), (1\,3)(2\,4), (1\,4)(2\,3)\}$ and $A_4$. All these t-level subsets are subgroups of the alternating group $A_4$. Hence $\mu$ is a fuzzy subgroup of the alternating group $A_4$. For x = (123) and y = (143), $\mu(xy) = \mu((123)(143)) = \mu((1\,2)(3\,4)) = p_1$ and $\mu(yx) = \mu((1\,4\,3)(1\,2\,3)) = \mu((1\,4)(2\,3)) = p_2$. As $p_1 > p_2$, $\mu(xy) \neq \mu(yx)$ for x = (1 2 3) and y = (1 4 3). Hence $\mu$ is not a fuzzy normal subgroup of $A_4$. Let x = (1 4)(2 3) and y = (1 3)(2 4) then for t = $p_1$, $\mu(xy^{-1}) = p_1 = t$.

Thus by the definition of t-level relation of μ we have (x, y) ∈ μ_t. Further we note that for a = (1 2 3), μ((ax)(ay)^{-1}) = p_2 < p_1. So by the definition of t-level relation of μ we have (ax, ay) ∉ μ_t for t = p_1 and α = (1 2 3). Hence it follows that μ_t is not a congruence relation on the alternating group A_4. So by the definition of congruence class, [x]_μ does not exist. That is if μ is not a fuzzy normal subgroup of the group G then the congruence class [x]_μ of μ_t containing the element x of the group G does not exist.

The following theorem is left as an exercise for the reader to prove.

**THEOREM 1.2.35:** *Let μ be a fuzzy normal subgroup of a group G and t ∈ [0, 1]. Then for every x ∈ G, $[x]_\mu = xG_\mu^t$ and $G_\mu^t$ is a normal subgroup of the group G.*

This theorem is however illustrated by the following example:

*Example 1.2.13:* Consider the Klein four group G = {a, b | $a^2 = b^2 = (ab)^2 = e$} where e is the identity element of G.

Define μ : G → [0, 1] by

$$\mu(x) = \begin{cases} 1 & \text{if } x = e \\ 0.6 & \text{if } x = a \\ 0.2 & \text{if } x = b, ab \end{cases}$$

Clearly all the t-level subsets of μ are normal subgroups of the group G. So μ is a fuzzy normal subgroup of the group G. For t = 0.5 we calculate the following

$[x]_\mu$ = {e, a} for all values of x = e and x = a
$[x]_\mu$ = {b, ab} for all values of x = b and x = ab

$xG_\mu^t$ = {e, a} for values of x = e and x = a and

$xG_\mu^t$ = {b, a b} for the values of x = b and x = ab. This verifies that $[x]_\mu = xG_\mu^t$ for every x ∈ G.

The following result is also left for the reader to prove.

**THEOREM 1.2.36:** *Let μ be a fuzzy normal subgroup of a group G, t ∈ [0, 1] and A be a non-empty subset of the group G. Then*

i. $\underline{\mu}_t(A) = \underline{G}_\mu^t(A)$

ii. $\hat{\mu}_t(A) = G_\mu^{\hat{t}}(A)$

*where $G_\mu^t$ is a normal subgroup of the group G.*

*Proof:* The proof of the following theorem is left for the reader as an exercise.

**THEOREM 1.2.37:** *Let $\mu$ and $\lambda$ be fuzzy normal subgroups of a group G and $t \in [0, 1]$. Let A and B be non-empty subsets of the group G. Then*

i. $\underline{\mu}_t(A) \subseteq A \subseteq \hat{\mu}_t(A)$
ii. $\hat{\mu}_t(A \cup B) = \hat{\mu}_t(A) \cup \hat{\mu}_t(B)$
iii. $A \cap B \neq \phi \Rightarrow \underline{\mu}_t(A \cap B) = \underline{\mu}_t(A) \cap \underline{\mu}_t(B)$
iv. $A \subseteq B \Rightarrow \underline{\mu}_t(A) \, \underline{\mu}_t(B)$
v. $A \subseteq B \Rightarrow \hat{\mu}_t(A) \cup \hat{\mu}_t(B)$
vi. $\underline{\mu}_t(A \cup B) \supseteq \underline{\mu}_t(A) \cup \underline{\mu}_t(B)$
vii. $A \cap B \neq \phi$, $\hat{\mu}_t(A \cap B) \subset \hat{\mu}_t(A) \cap \hat{\mu}_t(B)$
viii. $\mu_t \subset \lambda_t \Rightarrow \hat{\mu}_t(A) \subseteq \hat{\lambda}_t(A)$

Now we prove the following theorem:

**THEOREM 1.2.38:** *Let $\mu$ be a fuzzy normal subgroup of a group G and $t \in [0, 1]$. If A and B are non-empty subsets of the group G then $\hat{\mu}_t(A)\hat{\mu}_t(B) = \hat{\mu}_t(AB)$.*

*Proof:* Let $\mu$ be a fuzzy normal subgroup of a group G and $t \in [0, 1]$. Let A and B be any two non-empty subsets of the group G, then AB = $\{ab \mid a \in A \text{ and } b \in B\}$ is a non-empty subset of the group G. We have

$$\hat{\mu}_t(AB) = \hat{G}^t_\mu(AB)$$
$$= \hat{G}^t_\mu(A)\hat{G}^t_\mu(B)$$
$$= \hat{\mu}_t(A)\hat{\mu}_t(B)$$

Hence $\hat{\mu}_t(A)\hat{\mu}_t(B) = \hat{\mu}_t(AB)$.

**THEOREM 1.2.39:** *Let $\mu$ be a fuzzy normal subgroup of a group G and $t \in [0, 1]$. If A and B are non-empty subsets of the group G then $\underline{\mu}_t(A) \, \underline{\mu}_t(B) \subseteq \underline{\mu}_t(AB)$.*

*Proof:* Let $\mu$ be a fuzzy normal subgroup of a group G, $t \in [0, 1]$ and A and B by any two non-empty subsets of the group G. Then AB is non-empty as A and B are non-empty.

Consider
$$\underline{\mu}_t(A) \, \underline{\mu}_t(B) = \underline{G}^t_\mu(A) \, \underline{G}^t_\mu(B) \subseteq \underline{G}^t_\mu(AB) = \underline{\mu}_t(AB) .$$

Hence $\underline{\mu}_t(A) \, \underline{\mu}_t(B) \subseteq \underline{\mu}_t(AB)$.

The following theorem is left for the reader to prove.

**THEOREM 1.2.40:** *Let $\mu$ and $\lambda$ be fuzzy normal subgroups of a group G and $t \in [0, 1]$. If A is a non-empty subset of the group G then*

i. $(\mu \hat{\cap} \lambda)_t(A) = \hat{\mu}_t(A) \cap \hat{\lambda}_t(A)$

ii. $(\mu \cap \lambda)_t(A) = \underline{\mu}_t(A) \cap \underline{\lambda}_t(A)$.

**THEOREM 1.2.41:** *Let $\mu$ be a fuzzy normal subgroup of a group G and $t \in [0, 1]$. If A is a subgroup of the group G then $\hat{\mu}_t(A)$ is a subgroup of the group G.*

*Proof:* Let $\mu$ be a fuzzy normal subgroup of a group G and $t \in [0, 1]$. Then $G_\mu^t$ is a normal subgroup of a group G. A is a $\hat{G}_\mu^t$ rough subgroup of the group G. By the definition of rough subgroup, we have $\hat{G}_\mu^t(A)$ to be a subgroup of the group G. If $\mu$ is a fuzzy normal subgroup of a group G, $t \in [0, 1]$ and A is a non-empty subset of the group G then $\hat{\mu}_t(A) = \hat{G}_\mu^t(A)$ we have $\hat{\mu}_t(A)$ to be a subgroup of the group G.

Now we just recall some fuzzy relation and also the condition for the composition of two fuzzy subgroups to be a fuzzy subgroup.

**DEFINITION 1.2.17:** *Let $\mu$ be a fuzzy relation on S and let $\sigma$ be a fuzzy subset of S. Then $\mu$ is called a fuzzy relation on $\sigma$ if $\mu(x, y) \leq \min(\sigma(x), \sigma(y))$ for all $x, y \in S$*

*For any two fuzzy subsets $\sigma$ and $\mu$ of S; the cartesian product of $\mu$ and $\sigma$ is defined by $(\mu \times \sigma)(x, y) = \min(\mu(x), \sigma(y))$ for all $x, y \in S$.*

*Let $\sigma$ be a fuzzy subset of S. Then the strongest fuzzy relation on $\sigma$ is $\mu_\sigma$ defined by $\mu_\sigma(x, y) = (\sigma \times \sigma)(x, y) = \min(\sigma(x), \sigma(y))$ for all $x, y, \in S$.*

The following theorem can be easily verified.

**THEOREM 1.2.42:** *Let $\mu$ and $\sigma$ be fuzzy subsets of S. Then*

i. $\mu \times \sigma$ *is a fuzzy relation on S.*

ii. $(\mu \times \sigma)_t = \mu_t \times \sigma_t$ *for all $t \in [0, 1]$.*

The natural question would be when we have the strongest fuzzy relation can we ever have the weakest fuzzy subset of S; the answer is yes and it is defined as follows:

*If $\mu$ is a fuzzy relation on S, then the weakest fuzzy subset of S on which $\mu$ is a fuzzy relation is $\sigma_\mu$, defined by*

$$\sigma_\mu(x) = \sup_{y \in S} \{\max(\mu(x, y), \mu(y, x))\}$$

*for all $x \in S$. We define for any two fuzzy subset $\mu$ and $\sigma$ of G. $\mu \circ \sigma$ as*

$$(\mu \circ \sigma)(x) = \sup_{x=yz} \{\min (\mu(y), \sigma(z))\},$$

for all $x \in G$.

**DEFINITION [86]:** *A system of fuzzy singletons $\{(x_1)_{t_1}, \cdots, (x_k)_{t_k}\}$ where $0 < t_i < A(x_i)$ for $i = 1, 2, \ldots, k$ is said to be linearly independent in A if and only if $n_1(x_1)_{t_1} + \cdots + n_k(x_k)_{t_k} = 0_t$ implies $n_1 x_1 = \ldots = n_k x_k = 0$, where $n_i \in Z$, $i = 1, 2, \ldots, k$ and $t \in (0, 1]$. A system of fuzzy singletons is called dependent if it is not independent. An arbitrary system $\xi$ of fuzzy singleton is independent in A if and only if every finite sub-system of $\xi$ is independent.*

*We let $\xi$ denote a system of fuzzy singletons such that for all $x_t \in \xi$, $0 < t \leq A(x)$. $\xi^* = \{x \mid x_t \in \xi\}$ and $\xi_t = A_t \cap \xi^*$ for all $t \in (0, A(0)]$.*

**THEOREM 1.2.43:** *$\xi$ is independent in A if and only if the fuzzy subgroup of G generated by $\xi$ in A is a fuzzy direct sum of fuzzy subgroup of G whose support is cycle i.e. for*

$$\xi = \left\{(x_i)_{t_i} \mid 0 < t_i \leq A(x_i), i \in I\right\} \text{ holds } \langle \xi \rangle = \bigoplus_{i \in I} \langle (x_i)_{t_i} \rangle.$$

*Proof:* Left for the reader to prove as an exercise.

Now in the next section we introduce the concept of fuzzy bigroup which is very new and an interesting one.

## 1.3 Fuzzy sub-bigroup of a group

In this section we define fuzzy sub-bigroup of a bigroup [89, 135]. To define the notion of fuzzy sub-bigroup of a bigroup we define a new notion called the fuzzy union of any two fuzzy subsets of two distinct sets.

**DEFINITION [89, 135]:** *Let $\mu_1$ be a fuzzy subset of a set $X_1$ and $\mu_2$ be a fuzzy subset of a set $X_2$, then the fuzzy union of the fuzzy sets $\mu_1$ and $\mu_2$ is defined as a function.*

$\mu_1 \cup \mu_2 : X_1 \cup X_2 \to [0, 1]$ given by

$$(\mu_1 \cup \mu_2)(x) = \begin{cases} \max(\mu_1(x), \mu_2(x)) & \text{if } x \in X_1 \cap X_2 \\ \mu_1(x) & \text{if } x \in X_1 \text{ \& } x \notin X_2 \\ \mu_2(x) & \text{if } x \in X_2 \text{ \& } x \notin X_1 \end{cases}$$

We illustrate this definition by the following example:

***Example 1.3.1:*** Let $X_1 = \{1, 2, 3, 4, 5\}$ and $X_2 = \{2, 4, 6, 8, 10\}$ be two sets.

Define $\mu_1: X_1 \to [0, 1]$ by

$$\mu_1(x) = \begin{cases} 1 & \text{if } x = 1, 2 \\ 0.6 & \text{if } x = 3 \\ 0.2 & \text{if } x = 4, 5 \end{cases}$$

and define $\mu_2: X_2 \to [0, 1]$ by

$$\mu_2(x) = \begin{cases} 1 & \text{if } x = 2, 4 \\ 0.6 & \text{if } x = 6 \\ 0.2 & \text{if } x = 8, 10 \end{cases}$$

It is easy to calculate $\mu_1 \cup \mu_2$ and it is given as follows:

$$(\mu_1 \cup \mu_2)(x) = \begin{cases} 1 & \text{if } x = 1, 2, 4 \\ 0.6 & \text{if } x = 3, 6 \\ 0.2 & \text{if } x = 5, 8, 10 \end{cases}$$

Now we proceed on to define fuzzy sub-bigroup of a bigroup.

**DEFINITION 1.3.1:** *Let $G = (G_1 \cup G_2, +, \bullet)$ be a bigroup. Then $\mu: G \to [0, 1]$ is said to be a fuzzy sub-bigroup of the bigroup $G$ if there exists two fuzzy subsets $\mu_1$ (of $G_1$) and $\mu_2$ (of $G_2$) such that*

   i. *$(\mu_1, +)$ is a fuzzy subgroup of $(G_1, +)$*
   ii. *$(\mu_2, \bullet)$ is a fuzzy subgroup of $(G_2, \bullet)$ and*
   iii. *$\mu = \mu_1 \cup \mu_2$.*

We illustrate this by the following example

*Example 1.3.2:* Consider the bigroup $G = \{\pm i, \pm 0, \pm 1, \pm 2, \pm 3, \ldots\}$ under the binary operation '+' and '•' where $G_1 = \{0, \pm 1, \pm 2, \ldots\}$ and $G_2 = \{\pm 1, \pm i\}$.

Define $\mu : G \to [0, 1]$ by

$$\mu(x) = \begin{cases} \dfrac{1}{3} & \text{if } x = i, -i \\ 1 & \text{if } x \in \{0, \pm 2, \pm 4, \ldots\} \\ \dfrac{1}{2} & \text{if } x \in \{\pm 1, \pm 3, \ldots\} \end{cases}$$

It is easy to verify that $\mu$ is a fuzzy sub-bigroup of the bigroup G, for we can find
$\mu_1: G_1 \to [0, 1]$ by

$$\mu_1(x) = \begin{cases} 1 & \text{if } x \in \{0, \pm 2, \pm 4, \ldots\} \\ \dfrac{1}{2} & \text{if } x \in \{\pm 1, \pm 3, \ldots\} \end{cases}$$

and $\mu_2: G_2 \to [0, 1]$ given by

$$\mu_2(x) = \begin{cases} \dfrac{1}{2} & \text{if } x = 1, -1 \\ \dfrac{1}{3} & \text{if } x = i, -i \end{cases}$$

That is, there exists two fuzzy subgroups $\mu_1$ of $G_1$ and $\mu_2$ of $G_2$ such that $\mu = \mu_1 \cup \mu_2$.

Now we prove the following theorem.

**THEOREM 1.3.1:** *Every t-level subset of a fuzzy sub-bigroup $\mu$ of a bigroup G need not in general be a sub-bigroup of the bigroup G.*

*Proof:* The proof is by an example. Take $G = \{-1, 0, 1\}$ to be a bigroup under the binary operations '+' and '•' where $G_1 = \{0\}$ and $G_2 = \{-1, 1\}$ are groups respectively with respect to usual addition and usual multiplication.

Define $\mu: G \to [0, 1]$ by

$$\mu(x) = \begin{cases} \dfrac{1}{2} & \text{if } x = -1, 1 \\ \dfrac{1}{4} & \text{if } x = 0 \end{cases}$$

Then clearly $(\mu, +, \bullet)$ is a fuzzy sub-bigroup of the bigroup $(G, +, \bullet)$. Now consider the level subset $G_\mu^{\frac{1}{2}}$ of the fuzzy sub-bigroup $\mu$

$$G_\mu^{\frac{1}{2}} = \left\{ x \in G \mid \mu(x) \geq \dfrac{1}{2} \right\} = \{-1, 1\}.$$

It is easy to verify that $\{-1, 1\}$ is not a sub-bigroup of the bigroup $(G, +, \bullet)$. Hence the t-level subset

$$G_\mu^t \left( \text{for } t = \dfrac{1}{2} \right)$$

of the fuzzy sub-bigroup $\mu$ is not a sub-bigroup of the bigroup $(G, +, \bullet)$.

We define fuzzy sub-bigroup of a bigroup, to define this concept we introduce the notion of bilevel subset of a fuzzy sub-bigroup.

**DEFINITION 1.3.2:** *Let $(G = G_1 \cup G_2, +, \bullet)$ be a bigroup and $\mu = (\mu_1 \cup \mu_2)$ be a fuzzy sub-bigroup of the bigroup G. The bilevel subset of the fuzzy sub-bigroup $\mu$ of the bigroup G is defined as $G_\mu^t = G_{1\mu_1}^t \cup G_{2\mu_2}^t$ for every $t \in [0, \min\{\mu_1(e_1), \mu_2(e_2)\}]$, where $e_1$ denotes the identity element of the group $(G_1, +)$ and $e_2$ denotes the identity element of the group $(G_2, \bullet)$.*

**Remark:** The condition $t \in [0, \min\{\mu_1(e_1), \mu_2(e_2)\}]$ is essential for the bilevel subset to be a sub-bigroup for if $t \notin [0, \min\{\mu_1(e_1), \mu_2(e_2)\}]$ then the bilevel subset need not in general be a sub-bigroup of the bigroup G, which is evident from the following example:

*Example 1.3.3:* Take $\mu$ as in example 1.3.2 then the bi-level subset

$$G_\mu^t \left(\text{for } t = \frac{3}{4}\right)$$

of the fuzzy sub-bigroup $\mu$ is given by

$$G_\mu^t = \{0, \pm 2, \pm 4, \ldots\}$$

which is not a sub-bigroup of the bigroup G. Therefore the bilevel subset

$$G_\mu^t \left(\text{for } t = \frac{3}{4}\right)$$

is not a sub-bigroup of the bigroup G.

**THEOREM 1.3.2:** *Every bilevel subset of a fuzzy sub-bigroup $\mu$ of a bigroup G is a sub-bigroup of the bigroup G.*

*Proof:* Let $\mu \,(= \mu_1 \cup \mu_2)$ be the fuzzy subgroup of a bigroup $(G = G_1 \cup G_2, +, \bullet)$. Consider the bilevel subset $G_\mu^t$ of the fuzzy sub-bigroup $\mu$ for every $t \in [0, \min\{\mu_1(e_1), \mu_2(e_2)\}]$ where $e_1$ and $e_2$ denote the identity elements of the groups $G_1$ and $G_2$ respectively. Then $G_\mu^t = G_{2\mu_2}^t \cup G_{1\mu_1}^t$ where $G_{1\mu_1}^t$ and $G_{2\mu_2}^t$ are subgroups of $G_1$ and $G_2$ respectively (since $G_{1\mu_1}^t$ is a t-level subset of the group $G_1$ and $G_{2\mu_2}^t$ is a t-level subset of $G_2$).

Hence by the definition of sub-bigroup $G_\mu^t$ is a sub-bigroup of the bigroup $(G, +, \bullet)$. However to make the theorem explicit we illustrate by the following example.

*Example 1.3.4:* $G = \{0, \pm 1, \pm i\}$ is a bigroup with respect to addition modulo 2 and multiplication. Clearly $G_1 = \{0, 1\}$ and $G_2 = \{\pm 1, \pm i\}$ are group with respect to addition modulo 2 and multiplication respectively.

Define $\mu: G \to [0, 1]$ by

$$\mu(x) = \begin{cases} 1 & \text{if } x = 0 \\ 0.5 & \text{if } x = \pm 1 \\ 0.3 & \text{if } x = \pm i. \end{cases}$$

It is easy to verify that μ is a fuzzy sub-bigroup of the bigroup G as there exist two fuzzy subgroups $\mu_1: G \to [0, 1]$ and $\mu_2: G \to [0, 1]$ such that $\mu = \mu_1 \cup \mu_2$ where

$$\mu_1(x) = \begin{cases} 1 & \text{if } x = 0 \\ 0.4 & \text{if } x = 1 \end{cases}$$

and

$$\mu_2(x) = \begin{cases} 0.5 & \text{if } x = \pm 1 \\ 0.3 & \text{if } x = \pm i. \end{cases}$$

Now we calculate the bilevel subset $G_\mu^t$ for t = 0.5,

$$\begin{aligned} G_\mu^t = G_{1\mu_1}^t \cup G_{2\mu_2}^t &= \{x \in G_1 \mid \mu_1(x) \geq t\} \cup \{x \in G_2 \mid \mu_2(x) \geq t\} \\ &= \{0\} \cup \{\pm 1\} \\ &= \{0, \pm 1\}. \end{aligned}$$

That is $G_\mu^t = \{0, \pm 1\}$. It is easily verified that $G_\mu^t$ is a sub-bigroup of the bigroup G.

Now we proceed on to define fuzzy bigroup of a group.

**DEFINITION 1.3.3:** *A fuzzy subset μ of a group G is said to be a fuzzy sub-bigroup of the group G if there exists two fuzzy subgroups $\mu_1$ and $\mu_2$ of μ ($\mu_1 \neq \mu$ and $\mu_2 \neq \mu$) such that $\mu = \mu_1 \cup \mu_2$. Here by the term fuzzy subgroup λ of μ we mean that λ is a fuzzy subgroup of the group G and $\lambda \subseteq \mu$ (where μ is also a fuzzy subgroup of G).*

We illustrate the definition by the following example:

*Example 1.3.5*: Consider the additive group of integers. G = {0, ± 1, ± 2, …}.

Define μ: G → [0, 1] by

$$\mu(x) = \begin{cases} 1 & \text{if } x \in \{0, \pm 2, \pm 4, \ldots\} \\ 0.5 & \text{if } x \in \{\pm 1, \pm 3, \pm 5, \ldots\}. \end{cases}$$

It can be verified that μ is a fuzzy sub-bigroup of the group G, as there exists two fuzzy subgroups $\mu_1$ and $\mu_2$ of μ ($\mu_1 \neq \mu$ and $\mu_2 \neq \mu$) such that $\mu = \mu_1 \cup \mu_2$ where $\mu_1$ and $\mu_2$ are as given below.

$$\mu_1(x) = \begin{cases} 1 & \text{if } x \in \{0, \pm 2, \pm 4, \ldots\} \\ 0.25 & \text{if } x \in \{\pm 1, \pm 3, \ldots\} \end{cases}$$

and

$$\mu_2(x) = \begin{cases} 0.75 & \text{if } x \in \{0, \pm 2, \pm 4, \ldots\} \\ 0.5 & \text{if } x \in \{\pm 1, \pm 3, \ldots\}. \end{cases}$$

The following theorem relates the fuzzy sub-bigroup and the level subset.

**THEOREM 1.3.3:** *Let $\mu = \mu_1 \cup \mu_2$ be a fuzzy sub-bigroup of a group G, where $\mu_1$ and $\mu_2$ are fuzzy subgroups of the group G. For $t \in [0, \min\{\mu_1(e), \mu_2(e)\}]$, the level subset $G_\mu^t$ of $\mu$ can be represented as the union of two subgroups of the group G. That is $G_\mu^t = G_{1\mu_1}^t \cup G_{2\mu_2}^t$.*

*Proof*: Let $\mu$ be a fuzzy sub-bigroup of a group G and $t \in [0, \min\{\mu_1(e), \mu_2(e)\}]$. This implies that there exists fuzzy subgroups $\mu_1$ and $\mu_2$ of the group $G_1$ such that $\mu = \mu_1 \cup \mu_2$. Let $G_\mu^t$ be the level subset of $\mu$, then we have

$$\begin{aligned} x \in G_\mu^t &\Leftrightarrow \mu(x) \geq t \\ &\Leftrightarrow \max\{\mu_1(x), \mu_2(x)\} \geq t \\ &\Leftrightarrow \mu_1(x) \geq t \text{ or } \mu_2(x) \geq t \\ &\Leftrightarrow x \in G_{1\mu_1}^t \text{ or } x \in G_{2\mu_2}^t \end{aligned}$$

if and only if $x \in G_{1\mu_1}^t \cup G_{2\mu_2}^t$. Hence $G_\mu^t = G_{1\mu_1}^t \cup G_{2\mu_2}^t$.

***Example 1.3.6:*** Just in the example 1.3.5. For $t = 0.5$, $G_\mu^t = G$, $G_{1\mu_1}^t = 2G$ and $G_{2\mu_2}^t = G$. Hence $G_\mu^t = G = G_{1\mu_1}^t \cup G_{2\mu_2}^t$.

Now in the following theorem we give the condition for two fuzzy subgroups $\mu_1$ and $\mu_2$ of a fuzzy subgroup $\mu$ ($\mu_1 \neq \mu$ and $\mu_2 \neq \mu$) to be such that $\mu = \mu_1 \cup \mu_2$.

**THEOREM 1.3.4:** *Let $\mu$ be a fuzzy subgroup of a group G with $3 \leq o(Im(\mu)) < \infty$ then there exists two fuzzy subgroups $\mu_1$ and $\mu_2$ of $\mu$ ($\mu_1 \neq \mu$ and $\mu_2 \neq \mu$) such that $\mu = \mu_1 \cup \mu_2$.*

*Proof*: Let $\mu$ be a fuzzy subgroup of a group G. Suppose $Im(\mu) = \{a_1, a_2, \ldots, a_n\}$ where $3 \leq n < \infty$ and $a_1 > a_2 > \ldots > a_n$. Choose $b_1, b_2 \in [0, 1]$ be such that $a_1 > b_1 > a_2 > b_2 > a_3 > b_3 > \ldots > a_n$ and define $\mu_1, \mu_2: G \to [0, 1]$ by

$$\mu_1(x) = \begin{cases} a_1 & \text{if } x \in \mu_{a_1} \\ b_2 & \text{if } x \in \mu_{a_2} \setminus \mu_{a_1} \\ \mu(x) & \text{otherwise} \end{cases}$$

and

$$\mu_2(x) = \begin{cases} b_1 & \text{if } x \in \mu_{a_1} \\ a_2 & \text{if } x \in \mu_{a_2} \setminus \mu_{a_1} \\ \mu(x) & \text{otherwise.} \end{cases}$$

Thus it can be easily verified that both $\mu_1$ and $\mu_2$ are fuzzy subgroups of $\mu$. Further $\mu_1 \neq \mu$, $\mu_2 \neq \mu$ and $\mu = \mu_1 \cup \mu_2$.

Clearly the condition $3 \leq o(\text{Im}(\mu)) < \infty$ cannot be dropped in the above theorem. This is explained by the following example.

*Example 1.3.7:* Consider the group $G = \{1, -1, i, -i\}$ under the usual multiplication. Define $\mu : G \to [0, 1]$ by

$$\mu(x) = \begin{cases} 0 & \text{if } x = i, -i \\ 1 & \text{if } x = 1, -1 \end{cases}$$

Then it is easy to verify that $\mu$ is a fuzzy subgroup of the group G as all of its level subsets are subgroups of G. Further $o(\text{Im}(\mu)) = 2$. If $\mu_K$ is a fuzzy subgroup of $\mu$ such that $\mu_K \subseteq \mu$ ($\mu_K \neq \mu$) then $\mu_K$ takes the following form:

$$\mu_K(x) = \begin{cases} 0 & \text{if } x = i, -i \\ \alpha_i & \text{if } x = 1, -1 \end{cases}$$

with $0 \leq \alpha_K < 1$ for every K in the index set I. It is easy to verify that $\mu_j \cup \mu_K \neq \mu$ for any $j, K \in I$. Thus there does not exist two fuzzy subgroups $\mu_1$ and $\mu_2$ of $\mu$ ($\mu_1 \neq \mu$ and $\mu_2 \neq \mu$) such that $\mu = \mu_1 \cup \mu_2$.

Now we prove a very interesting theorem.

**THEOREM 1.3.5:** *Every fuzzy sub-bigroup of a group G is a fuzzy subgroup of the group G but not conversely.*

*Proof*: It follows from the definition of the fuzzy sub-bigroup of a group G that every sub-bigroup of a group G is a fuzzy subgroup of the group G.

However the converse of this theorem is not true. It is easy to see from example 1.3.7 that $\mu$ is a fuzzy subgroup of the group G and there does not exist two fuzzy subgroups $\mu_1$ and $\mu_2$ of $\mu$ ($\mu_1 \neq \mu$ and $\mu_2 \neq \mu$) such that $\mu = \mu_1 \cup \mu_2$. That is $\mu$ is not a fuzzy sub-bigroup of the group G. Now we obtain a necessary and sufficient condition for a fuzzy subgroup to be a fuzzy sub-bigroup of G.

**THEOREM 1.3.6:** *Let $\mu$ be a fuzzy subset of a group G with $3 \leq o(\text{Im}(\mu)) < \infty$. Then $\mu$ is a fuzzy subgroup of the group G if and only if $\mu$ is a fuzzy sub-bigroup of G.*

*Proof*: Let μ be a fuzzy subgroup of the group G with $3 \leq o(\text{Im}(\mu)) < \infty$ then there exists two fuzzy subgroups $\mu_1$ and $\mu_2$ of μ ($\mu_1 \neq \mu$ and $\mu_2 \neq \mu$) such that $\mu = \mu_1 \cup \mu_2$. Hence μ is a fuzzy sub-bigroup of the group G. Conversely, let μ be a fuzzy sub-bigroup of a group G, we know every fuzzy sub-bigroup of a group G is a fuzzy subgroup of the group G.

We shall illustrate this theorem by example.

*Example 1.3.8:* Define $\mu: G \to [0, 1]$ where $G = \{1, -1, i, -i\}$ by

$$\mu(x) = \begin{cases} 1 & \text{if } x = 1 \\ 0.9 & \text{if } x = -1 \\ 0.8 & \text{if } x = \pm i \end{cases}$$

It is easy to prove that μ is a fuzzy subgroup of the group G and $o(\text{Im}(\mu)) = 3$. Further, it can be easily verified that there exists two fuzzy subgroups $\mu_1$ and $\mu_2$ of μ ($\mu_1 \neq \mu$ and $\mu_2 \neq \mu$) such that $\mu = \mu_1\mu_2$ where $\mu_1, \mu_2 : G \to [0, 1]$ are defined by

$$\mu_1(x) = \begin{cases} 0.9 & \text{if } x = 1 \\ 0.8 & \text{if } x = -1, \pm i \end{cases}$$

and

$$\mu_2(x) = \begin{cases} 1 & \text{if } x = 1 \\ 0.9 & \text{if } x = -1 \\ 0.7 & \text{if } x = \pm i \end{cases}$$

The specialty about this section is that we have given examples to illustrate the theorem; this is mainly done to make one understand the concepts. As there are no books on bigroups and to the best of my knowledge the concept on fuzzy bigroups appeared in the year 2002 [135].

## 1.4 Fuzzy Rings and its properties

In this section we recall the concept of fuzzy rings and some of its basic properties. In 1971 [112, 145] introduced fuzzy sets in the realm of group theory and formulated the concept of a fuzzy subgroup of a group. Since then many researchers are engaged in extending the concept / results of abstract algebra to the broader frame work in fuzzy setting. However not all results on groups and rings can be fuzzified.

In 1982 [73] defined and studied fuzzy subrings as well as fuzzy ideals. Subsequently among [148, 149] fuzzified certain standard concepts on rings and ideals. Here we just recall some of the results on fuzzy rings and leave it for the reader to prove or get the proof of the results by referring the papers in the references. Now we recall definitions as given by [37].

**DEFINITION [145]:** *Let $\mu$ be any fuzzy subset of a set S and let $t \in [0, 1]$. The set $\{s \in S \mid \mu(x) \geq t\}$ is called a level subset of $\mu$ and is symbolized by $\mu_t$. Clearly $\mu_t \subseteq \mu_s$ whenever $t > s$.*

**DEFINITION [140]:** *Let '•' be a binary composition in a set S and $\mu$, $\sigma$ be any two fuzzy subsets of S. The product $\mu\sigma$ of $\mu$ and $\sigma$ is defined as follows:*

$$(\mu\sigma)(x) = \begin{cases} \sup_{x=y \bullet z}(\min(\mu(y)\sigma(z))) \text{ where } y, z \in S \\ 0 \quad \text{if } x \text{ is not expressible as } x = y \bullet z \text{ for all } y, z \in S. \end{cases}$$

**DEFINITION [140]:** *A fuzzy subset $\mu$ of a ring R is called a fuzzy subring of R if for all $x, y \in R$ the following requirements are met*

  i. $\mu(x - y) \geq \min(\mu(x), \mu(y))$ and
  ii. $\mu(xy) \geq \min(\mu(x), \mu(y))$

*Now if the condition (ii) is replaced by $\mu(xy) \geq \max(\mu(x), \mu(y))$ then $\mu$ is called a fuzzy ideal of R.*

**THEOREM 1.4.1:** *Let $\mu$ be any fuzzy subring / fuzzy ideal of a ring R. If, for some $x, y \in R$, $\mu(x) < \mu(y)$, then $\mu(x - y) = \mu(x) = \mu(y - x)$.*

*Proof:* Direct by the very definition of fuzzy subring / fuzzy ideal of a ring R.

Now we proceed on to define the notion of level subring/ level ideal of $\mu$, a fuzzy subring or fuzzy ideal of the ring R.

**DEFINITION [36]:** *Let $\mu$ be any fuzzy subring / fuzzy ideal of a ring R and let $0 \leq t \leq \mu(0)$. The subring / ideal $\mu_t$ is called a level subring / level ideal of $\mu$.*

**THEOREM 1.4.2:** *A fuzzy subset $\mu$ of a ring R is a fuzzy ideal of R if and only if the level subsets $\mu_t$, $t \in Im(\mu)$ are ideals of R.*

*Proof:* Left for the reader as it can be proved using the definition.

**THEOREM 1.4.3:** *If $\mu$ is any fuzzy ideal of a ring R, then two level ideals $\mu_{t_1}$ and $\mu_{t_2}$ (with $t_1 < t_2$) are equal if and only if there is no x in R such that $t_1 \leq \mu(x) \leq t_2$.*

*Proof:* Left an as exercise for the reader to prove as a matter of routine.

This theorem gives an insight that the level ideals of a fuzzy ideal need not be distinct.

**THEOREM 1.4.4:** *The level ideals of a fuzzy ideal $\mu$ form a chain. That is if Im $\mu = \{t_0, t_1, \ldots, t_n\}$ with $t_0 > \ldots > t_n$, then the chain of level ideals of $\mu$ is given by $\mu_{t_0} \subset \mu_{t_1} \subset \ldots \subset \mu_{t_n} = R$.*

*Proof*: Straightforward, hence omitted.

Now we proceed on to recall the definition of fuzzy prime ideal as given by [99, 100].

**DEFINITION [99]:** *A non-constant fuzzy ideal $\mu$ of a ring R is called fuzzy prime if for any fuzzy ideals $\mu_1$ and $\mu_2$ of R the condition $\mu_1 \mu_2 \subseteq \mu$ implies that either $\mu_1 \subseteq \mu$ or $\mu_2 \subseteq \mu$.*

It is left for the reader to prove the following theorem.

**THEOREM 1.4.5:** *The level ideal $\mu_t$, where $t = \mu(0)$ is a prime ideal of the ring R.*

**DEFINITION [99]:** *A fuzzy ideal $\mu$ of a ring R, not necessarily non-constant is called fuzzy prime if for any fuzzy ideals $\mu_1$ and $\mu_2$ of R the condition $\mu_1 \mu_2 \subseteq \mu$ implies that either $\mu_1 \subseteq \mu$ or $\mu_2 \subseteq \mu$.*

Consequent of this definition the following result is left for the reader to prove.

**THEOREM 1.4.6:** *Any constant fuzzy ideal $\mu$ of a ring R is fuzzy prime.*

Now we give the characterization theorem.

**THEOREM 1.4.7:** *If $\mu$ is any non-constant fuzzy ideal of a ring R, then $\mu$ is fuzzy prime if and only if $l \in $ Im $\mu$ : the ideal $\mu_t$, $t = \mu(0)$ is prime and the chain of level ideals of $\mu$ consists of $\mu_t \subset R$.*

*Proof*: Left for the reader to prove.

**DEFINITION 1.4.1:** *An ideal I of a ring R is said to be irreducible if I cannot be expressed as $I_1 \cap I_2$ where $I_1$ and $I_2$ are any two ideals of R properly containing I, otherwise I is termed reducible.*

The following results which will be used are given as theorems, the proof are for the reader to prove.

**THEOREM 1.4.8:** *Any prime ideal of ring R is irreducible.*

**THEOREM 1.4.9:** *In a commutative ring with unity, any ideal, which is both semiprime and irreducible, is prime.*

**THEOREM 1.4.10:** *Every ideal in a Noetherian ring is a finite intersection of irreducible ideals.*

**THEOREM 1.4.11:** *Every irreducible ideal in a Noetherian ring is primary.*

**DEFINITION [37]:** *A fuzzy ideal μ of a ring R is called fuzzy irreducible if it is not a finite intersection of two fuzzy ideals of R properly containing μ: otherwise μ is termed fuzzy reducible.*

**THEOREM [109]:** *If μ is any fuzzy prime ideal of a ring R, then μ is fuzzy irreducible.*

*Proof:* Refer [109].

However as a hint for the interested reader we mention the proof is given by the method of contradiction. We give the theorem from [109] the proof is omitted as it can be had from [109].

**THEOREM [109]:** *If μ is any non-constant fuzzy irreducible ideal of a ring R, then the following are true.*

   i.   $1 \in Im\ \mu$.
   ii.  *There exists $\alpha \in [0, 1]$ such that $\mu(x) = \alpha$ for all $x \in R \setminus \{x \in R\ /\ \mu(x)=1\}$.*
   iii. *The ideal $\{x \in R\ /\ \mu(x)=1\}$ is irreducible.*

It is important to recall some of the basic properties of fuzzy subrings and fuzzy ideals of a ring R.

The proofs of all these results are left for the reader to prove.

**THEOREM 1.4.12:** *The intersection of any family of fuzzy subrings (fuzzy ideals) of a ring R is again a fuzzy subring (fuzzy ideal) of R.*

**THEOREM 1.4.13:** *Let μ be any fuzzy subring and θ be any fuzzy ideal of a ring R. Then $\mu \cap \theta$ is a fuzzy ideal of the subring $\{x \in R\ /\ \mu(x) = \mu(0)\}$.*

**THEOREM 1.4.14:** *Let μ be any fuzzy subset of a field F. Then μ is a fuzzy ideal of F if and only if $\mu(x) = \mu(y) \leq \mu(0)$ for all $x, y \in F \setminus \{0\}$.*

**THEOREM 1.4.15:** *Let R be a ring. Then R is a field if and only if $\mu(x) = \mu(y) \leq \mu(0)$ where μ is any non-constant fuzzy ideal of R and $x, y \in R \setminus \{0\}$.*

**THEOREM 1.4.16:** *Let $I_0 \subset I_1 \subset \ldots \subset I_n = R$ be any chain of ideal of a ring R. Let $t_0, t_1, \ldots, t_n$ be some numbers lying in the interval $[0, 1]$ such that $t_0 > t_1 > \ldots > t_n$. Then the fuzzy subset μ of R defined by*

$$\mu(x) = \begin{cases} t_o & \text{if } x \in I_o \\ t_i & \text{if } x \in I_i \setminus I_{i-1}, i = 1,2,\ldots,n \end{cases}$$

*is a fuzzy ideal of R with $F_\mu = \{I_i \mid i = 0, 1, 2, \ldots, n\}$.*

Now we recall the definition of the product of any two fuzzy ideals of a ring R as given by [140].

**DEFINITION [140]:** *Let $\mu$ and $\theta$ be any fuzzy ideals of a ring R. The product $\mu \circ \theta$ of $\mu$ and $\theta$ is defined by*

$$(\mu \circ \theta)(x) = \sup_{\substack{x = \sum y_i z_i \\ i < \infty}} \left( \min_i \left( \min \left( \mu(y_i), \theta(z_i) \right) \right) \right),$$

*where $x, y_i, z_i \in R$.*

It is left for the reader to verify that $\mu \circ \theta$ is the smallest fuzzy ideal of R containing $\mu\theta$.

**Notation:** At times we will make use of this notation also. Let A ($\mu$) be any subset (fuzzy subset) of a ring R. The ideal (fuzzy subring / fuzzy ideal) generated by A($\mu$) is denoted by $\langle A \rangle$ ($\langle \mu \rangle$).

The following theorem which can be proved by a routine computation is left as an exercise for the reader to prove.

**THEOREM [36]:** *Let $\mu$ be a fuzzy subset of a ring R with card Im $\mu < \infty$. Define subrings $R_i$ of R by*

$$R_0 = \langle \{x \in R \mid \mu(x) = \sup_{z \in R} \mu(z)\} \rangle \text{ and}$$

$$R_i = \langle \{R_{i-1} \cup \{x \in R \mid \mu(x) = \sup_{z \in R - R_{i-1}} \mu(z)\} \rangle, 1 \leq i \leq k$$

*where k is such that $R_k = R$. Then $k <$ card Im $\mu$. Also the fuzzy subset $\mu*$ of R defined by*

$$\mu*(x) = \begin{cases} \sup_{z \in R} \mu(z) & \text{if } x \in R_0 \\ \sup_{z \in R \setminus R_{i-1}} \mu(z) & \text{if } x \in R_i \setminus R_{i-1}, \ 1 \leq i \leq k \end{cases}$$

*is a fuzzy subring generated by $\mu$ in R.*

To define the concept of fuzzy coset of a fuzzy ideal the following theorem will help; the proof of which is left for the reader as an exercise.

**THEOREM 1.4.17:**

i. Let $\mu$ be any fuzzy ideal of a ring R and let $t = \mu(0)$. Then the fuzzy subset $\mu*$ of $R/\mu_t$ defined by $\mu*(x + \mu_t) = \mu(x)$ for all $x \in R$, is a fuzzy ideal of $R/\mu_t$.

ii. *If A is an ideal of R and θ is a fuzzy ideal of R/A such that θ (x + A) = θ(A) only when x ∈ A, then there exists a fuzzy ideal μ of R such that $\mu_t$ = A where t = μ(0) and θ = μ*.*

Now we recall the definition of fuzzy coset.

**DEFINITION [36]:** *Let μ be any fuzzy ideal of a ring R and let x ∈ R. The fuzzy subset $\mu_x$* of R defined by $\mu_x$*(r) = μ (r – x) for all r ∈ R is termed as the fuzzy coset determined by x and μ.*

For more about these concepts please refer [109]. The following theorem is given, the proof of which is a matter of routine hence left for the reader as an exercise.

**THEOREM 1.4.18:** *Let μ be any fuzzy ideal of a ring R. Then $R_\mu$, the set of all fuzzy cosets of μ in R is a ring under the binary compositions.*

$$\mu_x* + \mu_y* = \mu*_{x+y} \text{ and}$$
$$\mu_x* \mu_y* = \mu*_{xy} \text{ for all } x, y \in R.$$

Now we proceed on to recall the definition of fuzzy quotient ideal of ring R.

**DEFINITION [109]:** *If μ is any fuzzy ideal of a ring R, then the fuzzy ideal μ' of $R_\mu$ defined by μ'($\mu_x^*$) = μ(x) for all x ∈ R is called the fuzzy quotient ideal determined by μ.*

The proof can be had from [109].

**THEOREM [109]:** *If μ is any fuzzy ideal of a ring R, then the map f : R → $R_\mu$ defined by f(x) = $\mu_x^*$ for all x ∈ R is a homomorphism with kernel $\mu_t$, where t = μ (0).*

**THEOREM [109]:** *If μ is any fuzzy ideal of a ring R, then each fuzzy ideal of $R_\mu$ corresponds in a natural way to a fuzzy ideal of R.*

*Proof*: Let μ' be any fuzzy ideal of $R_\mu$. It is entirely straightforward matter to show that the fuzzy subset θ of R defined by θ (x) = μ' ($\mu_x^*$) for all x ∈ R, is a fuzzy ideal of R.

Now we proceed on to define fuzzy semiprime ideal.

**DEFINITION [109]:** *A fuzzy ideal μ of a ring R is called fuzzy semiprime if for any fuzzy ideal θ of R, the condition $\theta^m \subseteq \mu$ (m ∈ $Z^+$) implies θ ⊂ μ.*

The following results are direct and can be obtained as a matter of routine and hence left for the reader as an exercise.

**THEOREM 1.4.19:** *Let μ be any fuzzy subset of a ring R. Then μ(x) = t if and only if x ∈ $\mu_t$ and x ∉ $\mu_s$ for all s > t.*

**THEOREM 1.4.20:** *A fuzzy ideal $\mu$ of a ring R is fuzzy semiprime if and only if $\mu_t$, $t \in$ Im $\mu$, is a semiprime ideal of R.*

**THEOREM 1.4.21:** *An ideal A of a ring R is semiprime if and only if $\psi_A$ is a fuzzy semiprime ideal of R.*

**THEOREM 1.4.22:** *If $\mu$ is any fuzzy semiprime ideal of a ring R, then $R_\mu$, the ring of fuzzy cosets of $\mu$ in R is free from non-zero nilpotent elements.*

**THEOREM 1.4.23:** *Let $\mu$ be any fuzzy ideal of a ring R such that Im $\mu = \{t, s\}$ with $t > s$. If the ring $R_\mu$ has no non-zero nilpotent elements, then the fuzzy ideal $\mu$ is fuzzy semiprime.*

**THEOREM [99]:** *A ring R is regular if and only if every fuzzy ideal of R is idempotent.*

*Proof*: Refer [99].

**THEOREM 1.4.24:** *A ring R is regular if and only if every fuzzy ideal of R is fuzzy semiprime.*

*Proof:* Straightforward hence left for the reader to prove.

Now we recall some results on fuzzy subrings from [40].

Now for the time being we assume that R is a commutative ring with unit and M will denote a maximal ideal of R. If R is quasi local, then M is the unique maximal ideal of R and we write (R, M), whenever we say a subring of R we assume the subring contains the identity of R. We let g denote the natural homomorphism of R onto R/M. If S is a subset of R, we let $\langle S \rangle$ denote the ideal of R generated by S.

Let X and Y be fuzzy subsets of R i.e. functions from R into the closed interval [0, 1]. Then X+Y is the fuzzy subset of R defined by $\forall\ z \in R$, $(X+Y)(z) = \sup\{\min(X(x), Y(y)) \mid z = x + y\}$. We say that $X \supseteq Y$ if $X(x) \geq Y(x)$ for all $x \in R$. Let $\{X_i \mid i \in I\}$ be a collection of fuzzy subset of R. Define the fuzzy subset

$$\left(\bigcap_{i \in I} X_i\right)$$

of R by, for all $z \in R$,

$$\left(\bigcap_{i \in I} X_i\right)(x) = \inf\{X_i(z) \mid i \in I\},$$

If U is the subset of X, we let $\delta_U$ denote the characteristic function of U in R. We let Im (X) denote the image of X. We say that X is finite valued if Im (X) is finite. We let $X^* = \{x \in R \mid X(x) > 0\}$ the support of X. For all $t \in [0, 1]$ we let $X_t = \{x \in R \mid X(x) \geq t\}$ a level subset of X in R.

Let A be a fuzzy subring of R and S a subring of R. Define $A|_s$ of R by $(A|_s)(x) = A(x)$ if $x \in S$ and $(A|_s)(x) = 0$ otherwise. It follows that $A|_s$ is a fuzzy subring of R. If B is a fuzzy subring of R such that $B^*$ is a field containing the identity of R and $B(x) = B(x^{-1})$ for all units x in R, then B is called the fuzzy subfield in R. R is said to have a coefficient field with respect to M if there exist a field F containing the identity of R such that $R = F + M$. Suppose that R has a coefficient field F with respect to M. If $A = A|_F + A|_M$ where $A|_F$ is a subfield in R then $A|_F$ is called the fuzzy coefficient field of A with respect to $A|_M$.

**THEOREM [40]:** *Suppose that M is nil and R has a coefficient field. Let A be a fuzzy subring of R such that for all t, s $\in$ Im(A) with $t > s$, $(A_s \cap F)^P \cap (A_t \cap F) = (A_t \cap F)^P$ then $A = A|_F + A|_M$.*

*Proof*: Refer [40].

Using the above theorem one can easily prove the following theorem, which is once again left for the reader as an exercise.

**THEOREM [40]:** *Suppose that M is nil, R contains a field K and R/M is separable algebraic over K. Let A be a fuzzy subring of R such that $A \supseteq \delta_K$. Then $A = A|_F + A|_M$ and $A|_F$ is a fuzzy subfield in R where F is the coefficient field containing K.*

*Proof*: Left for the reader to refer [40].

**THEOREM [40]:** *Let A be a fuzzy subring of R. Suppose that Im (A) = $\{t_0, t_1, …, t_n\}$ where $t_0 > t_1 > … > t_n$. Suppose further that (R, M) and $(A_t, M \cap A_t)$ are complete local rings for all $t \in$ Im (A) such that $R^P \subset A_{t_0}$. Then R has a coefficient field F such that $A = A|_F + A|_M$ where $A|_F$ is a fuzzy subfield of R.*

*Proof:* Refer [40].

Suppose F/K be a field extension and let C and D be fuzzy subfields of F such that D, $C \supseteq \delta_K$. If C and D are linearly disjoint over $\delta_K$ then we write $CD = C \otimes D$ where CD is the composite of C and D. The following results give the structure of the fuzzy field extension $A|_F$ over $A|_K$ where $K \subset A_{t_0}$. The proof of these theorems are left for the reader as an exercise, we only state the theorem; interested reader can refer [40].

**THEOREM [40]:** *Let F/K be a field extension. Suppose that A is a fuzzy subfield of F such that $A \supseteq \delta_K$ and that A is finite valued say Im(A) − {0} = $\{t_0, t_1, …, t_n\}$ where $t_0 > t_1 > t_2 > … > t_n$. If there exists intermediate field $H_i$ of $A_t/K$ such that $A_{t_i} = A_{t_{i-1}} \otimes_K H_i$ for $i = 1, …, n$ then there exists fuzzy subfield $A_i$ of F, $i = 0, 1, …, n$, such that $A = A_0 \otimes A_1 \otimes … \otimes A_n$ (over $\delta_K$), $K \subset A_{t_0} = A_0^*$ and $A_i^* = H_i$, $i = 1, …, n$.*

**THEOREM [40]:** *Let F/K be a field extension with characteristic $p > 0$. Suppose that $F^P \subseteq K$. Let A be a fuzzy subfield of F such that $A^* = F$ and $A \supseteq \delta_K$. Suppose that A is finite valued say Im (A) − {0} = $\{t_0, t_1, …, t_n\}$ where $t_0 > t_1 > t_2 > … > t_n$. Then there exists fuzzy subfield $A_i$ of $F_i$, $i = 0, 1, 2, …, n$, such that $A = A_0 \otimes A_1 \otimes … \otimes A_n$.*

*Further more for all i = 0, 1, 2, …, n there exists set $T_i$ of fuzzy subfields of F such that for all $L \in T_i$, $L \supseteq \delta_K$, $L^* / K$ is simple and*

$$A_i = \bigotimes_{L \in T_i} L \,(\text{over } \delta_K).$$

Now we proceed on to define fuzzy polynomial subring and fuzzy power series subring.

**DEFINITION [40]:** *Let S be a commutative ring with identity. Let $R = S[x_1, …, x_n]$ be a polynomial ring in the indeterminates $x_1$, $x_2$, …, $x_n$ over S. Let A be a fuzzy subring of R. If for all*

$$r = \sum_{i_n=0}^{m_n} \ldots \sum_{i_1=0}^{m_1} C_{i_1 i_2 \ldots i_n} \, x_1^{i_1} \ldots x_n^{i_n} \in R$$

*where $C_{i_1 i_2 \ldots i_n} \in S$. $A(r) = \min \{(\min\{A(C_{i_1 i_2 \ldots i_n})\mid i_j = 0, 1, …, m_j, j = 1, 2, …, n\}$, $\min \{t_j \mid j = 1, 2, …, n\}\}$, then A is called a fuzzy polynomial subring of R where $t_j = A(x_j)$ if $x_j$ appears non-trivially in the above representation of r and $t_j = 1$ otherwise, $j = 1, 2, …, n$.*

**DEFINITION [40]:** *Let S be a commutative ring with identity. Let $R = S[[x_1, …, x_n]]$ be a power series ring in the indeterminates $x_1$, $x_2$, …, $x_n$ over S. Let A be a fuzzy subring of R.*

*If for all*

$$r = \sum_{i_n=0}^{\infty} \ldots \sum_{i_1=0}^{\infty} C_{i_1 \ldots i_n}, x_1^{i_1} \ldots x_n^{i_n} \in R \text{ where } C_{i_1 \ldots i_n} \in S.$$

*$A(r) = \min \{\inf (A (C_{i_1 \ldots i_n})\mid i_j = 0, 1, …, \infty_i; j = 1, 2, …, n\}$, $\min (t_j) \mid j = 1, 2, …, n\}\}$*

*then A is called a fuzzy power series subring of R where $t_j$ is defined as in the above equation.*

The reader is expected to prove the following theorems.

**THEOREM [40]:** *Let $R = S[x_1, x_2, …, x_n]$ be a polynomial ring over S and $M = \langle x_1, …, x_n \rangle$. Let A be a fuzzy subring of R. Then A is a fuzzy polynomial subring of R if and only if $A = A|_S + A|_M$ and for all $r \in R$, $r \neq 0$, and for all $x_j$, $j = 1, 2, …, n$. $A(x_j r) = \min \{A(x_j), A(r)\}$.*

*Proof:* Left for the reader as an exercise.

The following theorem is interesting and left for the reader to prove.

**THEOREM [40]:** *Let $R = S[[x_1, \ldots, x_n]]$ be a power series ring in the in determinants $x_1, x_2, \ldots, x_n$ over S and $M = \langle x_1, \ldots, x_n \rangle$. Let A be a fuzzy subring of R.*

  i. *If A is a fuzzy power series subring of R, then $A = A|_S + A|_M$ and for all $r \in R$, $r \neq 0$ and for all $x_j$, $j = 1, 2, \ldots, n$, $A(x, r) = \min\{A(x_j), A(r)\}$.*

  ii. *Suppose that $A|_M$ is a fuzzy ideal in R. If $A = A|_S + A|_M$ and for all $r \in R$, $r \neq 0$ and for all $x_j$, $j = 1, 2, \ldots, n$. $A(x, r) = \min\{A(x_j), A(r)\}$ then A is a fuzzy power series subring of R.*

Now we proceed on to recall some more notions on fuzzy polynomial subrings.

**THEOREM [40]:** *Suppose that $R = S[x] (S[[x]])$ is a polynomial (power series) ring over S and that A is a fuzzy polynomial (power series) subring of R. If $A/M$ is a fuzzy ideal of R, then A is constant on $M - \langle 0 \rangle$ where $M = \langle x \rangle$.*

*Proof*: The reader is requested to refer [40].

We will also using the notation for $A^*$ and $A_*$ obtain some interesting results about fuzzy subrings.

**DEFINITION 1.4.2:** *Let A be a fuzzy subset of R. Define*

$$A^* = \{x \in R \mid A(x) > 0\} \text{ and}$$
$$A_* = \{x \in R \mid A(x) = A(0)\}.$$

The following theorem is Straightforward hence left for the reader to prove.

**THEOREM 1.4.25:** *Let A be a fuzzy subring (fuzzy ideal) of R. Then $A_*$ is a subring (ideal) of R. If L has the finite intersection property then $A^*$ is a subring (ideal) of R.*

**THEOREM 1.4.26:** *Let $\{A_\alpha \mid \alpha \in \Omega\}$ be a collection of fuzzy subrings of R. Suppose that L has the finite intersection property. Then*

$$\sum_{\alpha \in \Omega} A_\alpha^* = \left(\sum_{\alpha \in \Omega} A_\alpha\right)^*.$$

*Proof:* Follows from the fact that if

$$x \in \left(\sum_{\alpha \in \Omega} A_\alpha\right)^*$$

if and only if

$$\left(\sum_{\alpha \in \Omega} A_\alpha\right)(x) > 0$$

if and only if

$$\sup\left\{\inf\left\{A_\alpha(x_\alpha) \mid \alpha \in \Omega\right\} \mid x = \sum_{\alpha \in \Omega} x_\alpha\right\} > 0$$

if and only if

$$x = \sum_{\alpha \in \Omega} x_\alpha,$$

for some $x_\alpha \in A_\alpha^*$ if and only if

$$x \in \sum_{\alpha \in \Omega} A_\alpha^*.$$

The following result is also Straightforward left for the reader as exercise.

**THEOREM 1.4.27:** *Let $\{A_\alpha \mid \alpha \in \Omega\}$ be a collection of fuzzy subrings of R. Suppose that L has the intersection property. Then*

$$\bigcap_{\alpha \in \Omega} A_\alpha^* = \left(\bigcap_{\alpha \in \Omega} A_\alpha\right)^*.$$

**THEOREM 1.4.28:** *Let $\{A_\alpha \mid \alpha \in \Omega\}$ be a collection of fuzzy subrings of R. Suppose that L has the finite intersection property. Then for all $\beta \in \Omega$*

$$A_\beta^* \cap \sum_{\alpha \in \Omega_\beta} A_\alpha^* = \{0\}$$

*if and only if for all $x \in R$, $x \neq 0$*

$$\left(A_\beta \cap \sum_{\alpha \in \Omega_\beta} A_\alpha\right) = 0.$$

*Proof:* Follows from the definitions.

**THEOREM 1.4.29**: *Let $\{A_\alpha \mid \alpha \in \Omega\} \cup \{A\}$ be a collection of fuzzy subrings of R such that $A = \sum_{\alpha \in \Omega} A_\alpha$. Suppose that L has the finite intersection property, then $A^* = \bigoplus_{\alpha \in \Omega} A_\alpha^*$ if and only if $A = \bigoplus_{\alpha \in \Omega} A_\alpha$.*

*Proof***:** Follows from the very definitions.

**THEOREM 1.4.30:** *Let $\{A_\alpha \mid \alpha \in \Omega\}$ be a collection of fuzzy subrings of R. Then*

$$\left(\sum_{\alpha \in \Omega} A_{\alpha^*}\right) \subseteq \left(\sum_{\alpha \in \Omega} A_\alpha\right)_*$$

*Proof:* Left as an exercise for the reader to prove.

**THEOREM 1.4.31:** *Let $\{ A_\alpha \mid \alpha \in \Omega \}$ be a collection of fuzzy subrings of R. If there exists $t \in L$, $t \neq 1$, such that $t \geq \sup\{A_\alpha(x) \mid x \notin A_{\alpha^*} \text{ for all } \alpha \in \Omega\}$ then*

$$\sum_{\alpha \in \Omega} A_{\alpha^*} = \left(\sum_{\alpha \in \Omega} A_\alpha\right)_*.$$

*Proof:* Follows from the very definition.

**THEOREM 1.4.32:** *Let $\{ A_\alpha \mid \alpha \in \Omega \}$ be a collection of fuzzy subrings of R. Suppose that L has finite intersection property, if*

$$\sum_{\alpha \in \Omega} A_\alpha = \bigoplus_{\alpha \in \Omega} A_\alpha$$

*then*

$$\sum_{\alpha \in \Omega} A_{\alpha^*} = \bigoplus_{\alpha \in \Omega} A_{\alpha^*}.$$

*Proof:* Left for the reader as an exercise.

**THEOREM 1.4.33:** *Let $\{A_\alpha \mid \alpha \in \Omega\}$ be a collection of fuzzy subrings of R. Then*

$$A_\beta \sum_{\alpha \in \Omega} A_\alpha = \sum_{\alpha \in \Omega} A_\beta A_\alpha$$

*and*

$$\left(\sum_{\alpha \in \Omega} A_\alpha\right) A_\beta = \sum_{\alpha \in \Omega} A_\alpha A_\beta.$$

*Proof:* Simple and straightforward for the reader to prove.

**THEOREM 1.4.34:** *If A and B are fuzzy subrings (fuzzy ideals of R). Then*

$$AB(x + y) \geq \inf(AB(x), AB(y))$$
$$AB(x) = AB(-x)$$

*for all $x, y \in R$.*

*Proof:* Left for the reader as an exercise.

**THEOREM 1.4.35:** *Let A and B be fuzzy subrings of R such that*

$$\sup \{\{\sup A(x) \mid x \notin A_\alpha\}, \sup \{B(x) \mid x \notin B_\alpha\}\} < 1.$$

*Then $(AB)_* = A_* B_*$.*

*Proof:* $x \in (AB)_*$ if and only if $(AB)(x) = 1$ if and only if $\sup \{\inf \{\inf (A(y_i) B(z_i))$

$$\Big| i = 1, 2, \ldots, n\} \text{ such that } x = \sum_{i=1}^{n} y_i z_i \text{ with } n \in N\} = 1$$

if and only if

$$x = \sum_{i=1}^{n} y_i z_i$$

for some $y_i \in A_*$ and $z_i \in B_*$, ($i = 1, 2, \ldots, n$) if and only if $x \in A_* B_*$.

The following theorem is left for the reader to prove.

**THEOREM 1.4.36:** *Suppose that R has an identity. Let $\{ A_\alpha \mid \alpha \in \Omega \}$ be a collection of fuzzy subrings of R. Suppose that*

$$\delta_R = \sum_{\alpha \in \Omega} A_\alpha$$

*and for all $x \in R$, $x \neq 0$ $(A_\alpha A_\beta)(x) = 0$ for all $\alpha, \beta \in \Omega$, $\alpha \neq \beta$. Then for all $\alpha \in \Omega$, $A_\alpha$ is a fuzzy ideal of R. If L has the finite intersection property then*

$$\delta_R = \bigoplus_{\alpha \in \Omega} A_\alpha.$$

Now we proceed on to define fuzzy left coset for a fuzzy ideal A of R.

**DEFINITION 1.4.3:** *Let A be a fuzzy ideal of R. For all $r \in R$, define fuzzy left coset $r + A$ by for all $x \in R$ $(r + A)(x) = A(x - r)$. Let $R/A = \{r + A \mid r \in A\}$.*

The following theorem is Straightforward and hence left for the reader to prove.

**THEOREM 1.4.37:** *Let A be a fuzzy ideal of R. Define '+' and '•' on R/A by for all r, s $\in R$ $(r + A) + (s + A) = r + s + A$ and $(r + A) \bullet (s + A) = r \bullet s + A$. Then R/A is a commutative ring and $R/A = R/A_*$.*

**Notation:** Let $R_\alpha$ be a subring of R and let A be a fuzzy ideal of R. Let $A'_\alpha$ denote the restriction of $\delta_{R_\alpha} A$ of $R_\alpha$. Then $A'_\alpha$ is a fuzzy ideal of $R_\alpha$.

**THEOREM 1.4.38:** *Suppose that R has an identity. Let $\{R_\alpha \mid \alpha \in \Omega\}$ be a collection of ideals of R such that $R = \bigoplus_{\alpha \in \Omega} R_\alpha$. Let A be a fuzzy ideal of R. Suppose that L has the finite intersection property. Then there exists a collection of fuzzy ideals of R, $\{A_\alpha \mid \alpha \in \Omega\}$ such that $A = \bigoplus_{\alpha \in \Omega} A_\alpha$ and in fact $A_\alpha = \delta_{R_\alpha} A$ for all $\alpha \in \Omega$. If $1 > \sup \{A(x) \mid x \notin A_\alpha\}$ then $A'_\alpha$ is a fuzzy ideal of $R_\alpha$ and $R/A = \bigoplus_{\alpha \in \Omega} R_\alpha / A'_\alpha$.*

*Proof*: Refer [84]. The following result can be proved as a matter of routine, hence left as an exercise for the reader to prove.

**THEOREM 1.4.39:** *Let $\{A_\alpha \mid \alpha \in \Omega\} \cup \{A\}$ be a collection of fuzzy subrings of R. Suppose that L has the finite intersection property. If $A = \bigoplus_{\alpha \in \Omega} A_\alpha$ then for all $x \in A^*$,*

$$x = \sum_{\alpha \in \Omega} x_\alpha$$

*for unique*

$$x_\alpha \in A_\alpha^* \text{ and } A(x) = \inf\left\{ A_\alpha(x_\alpha) \,\bigg|\, x = \sum_{\alpha \in \Omega} x_\alpha, x_\alpha \in A_\alpha^* \right\}.$$

**THEOREM 1.4.40:** *Let A and B be fuzzy subrings of R. Suppose that L has the finite intersection property. Then $(AB)^* = A^* B^*$.*

*Proof:* Straightforward, hence left for the reader to prove.

**THEOREM 1.4.41:** *Suppose that R has an identity. Let $\{R_\alpha \mid \alpha \in \Omega\}$ be a collection of ideals of R such that $R = \bigoplus_{\alpha \in \Omega} R_\alpha$. Let A be a fuzzy ideal of R such that $1 > \sup \{A(x) \mid x \notin A_\alpha\}$. Let $A_\alpha = \delta_{R_\alpha} A$, for all $\alpha \in \Omega$. Suppose that L has the finite intersection property. Then A is a fuzzy maximal or fuzzy prime or fuzzy primary (L = [0, 1] here) if and only if for all but one of the $A_\alpha$'s is 1 on $R_\alpha$ and the remaining $A'_\alpha$ is fuzzy maximal or fuzzy prime or fuzzy primary in $R_\alpha$.*

*Proof:* Left for the reader as an exercise.

Now we proceed on to recall the notion of pairwise co-maximal.

**DEFINITION 1.4.4:** *Suppose that R has an identity. Let $\{A_\alpha \mid \alpha \in \Omega\}$ be fuzzy ideals of R. The $A_\alpha$ is said to be pair wise co-maximal if and only if $A_\alpha \neq \delta_R$ and $A_\alpha + A_\beta = \delta_R$ for all $\alpha, \beta \in \Omega, \alpha \neq \beta$.*

The following theorem is straightforward and left for the reader to prove.

**THEOREM 1.4.42:** *Suppose that R has an identity. Let $\{A_\alpha \mid \alpha \in \Omega\}$ be a collection of finite valued fuzzy ideals of R. Then the $A_\alpha$ are pairwise co-maximal if and only if $A_\alpha*$ are pairwise co maximal.*

Some results on finite collection of finite valued fuzzy ideals of R is recalled.

**THEOREM 1.4.43:** *Suppose that R has an identity. Let $\{A_\alpha \mid \alpha \in \Omega\} \cup \{B\}$ be a finite collection of finite valued fuzzy ideals of R.*

i. The $A_\alpha$ are pairwise co maximal if and only if the $\sqrt{A_\alpha}$ are pairwise co maximal where $L = [0, 1]$.

ii. If $B$ is co maximal with each $A_\alpha$, then $B$ is co maximal with $\bigcap_{\alpha \in \Omega} A_\alpha$ and $\prod_{\alpha \in \Omega} A_\alpha$ where $L$ has the finite intersection property.

iii. If the $A_\alpha$ are pairwise co maximal then $\bigcap_{\alpha \in \Omega} A_\alpha = \prod_{\alpha \in \Omega} A_\alpha$ where $L$ has the finite intersection property.

Now we proceed on to recall the definition of fuzzy external direct sum.

**DEFINITION 1.4.5:** *Let $\{X_\alpha \mid \alpha \in \Omega\}$ be a collection of non-empty sets and let $A_\alpha$ be a fuzzy subset of $X_\alpha$ for all $\alpha \in \Omega$. Define the Cartesian cross product of the $A_\alpha$ by for all $x_\alpha \in X_\alpha$*

$$\left( \underset{\alpha \in \Omega}{\times} A_\alpha \right)(x) = \inf \{A_\alpha(x_\alpha) \mid \alpha \in \Omega\}$$

*where $x = \langle x_\alpha \rangle$ and $\langle x_\alpha \rangle$ denotes an element of the Cartesian cross product $\underset{\alpha \in \Omega}{\times} X_\alpha$.*

The following theorem can be proved as a matter of routine.

**THEOREM 1.4.44:** *Let $\{R_\alpha \mid \alpha \in \Omega\}$ be a collection of commutative rings and let $A_\alpha$ be a fuzzy subring (fuzzy ideal) of $R_\alpha$ for all $\alpha \in \Omega$. Then $\underset{\alpha \in \Omega}{\times} A_\alpha$ is a fuzzy subring (fuzzy ideal) of $\underset{\alpha \in \Omega}{\times} R_\alpha$ where by $\underset{\alpha \in \Omega}{\times} R_\alpha$ we mean external direct sum of the $R_\alpha$.*

The following theorems are very direct by the use of definitions.

**THEOREM 1.4.45:** *Let $\{R_\alpha \mid \alpha \in \Omega\}$ be a collection of commutative rings and let $A_\alpha$ be a fuzzy subring (fuzzy ideal) of $R_\alpha$ for all $\alpha \in \Omega$. Then for all $\beta \in \Omega$, $\underset{\alpha \in \Omega}{\times} A_\alpha^{(\beta)}$ is a fuzzy subring (fuzzy ideal) of $\underset{\alpha \in \Omega}{\times} R_\alpha$.*

**THEOREM 1.4.46:** *Let $\{A_\alpha \mid \alpha \in \Omega\}$ be a collection of fuzzy subrings of $R$. Suppose that $L$ has the finite intersection property. If*

$$\sum_{\alpha \in \Omega} A_\alpha = \bigoplus_{\alpha \in \Omega} A_\alpha \text{ then } \sum_{\alpha \in \Omega} A_\alpha = \underset{\alpha \in \Omega}{\times} A_\alpha \text{ on } \left( \sum_{\alpha \in \Omega} A_\alpha \right)^*.$$

**THEOREM 1.4.47:** *Let $\{R_\alpha \mid \alpha \in \Omega\}$ be a collection of commutative rings. Let $A_\alpha$ be a fuzzy subring (fuzzy ideal) of $R_\alpha$ for all $\alpha \in \Omega$. Suppose that $L$ has the finite intersection property. Then in $\underset{\alpha \in \Omega}{\times} R_\alpha$, $\bigoplus_{\beta \in \Omega} \underset{\alpha \in \Omega}{\times} A_\alpha^{(\beta)} = \underset{\alpha \in \Omega}{\times} A_\alpha$.*

Now we proceed on to define the concept of complete fuzzy direct sum.

**DEFINITION 1.4.6:** *Let $\{R_\alpha \mid \alpha \in \Omega\}$ be a collection of commutative rings and let $A_\alpha$ be a fuzzy subring of $R_\alpha$, for all $\alpha \in \Omega$. Then the Cartesian cross product $\underset{\alpha \in \Omega}{\times} A_\alpha$ is called the complete fuzzy direct sum of the $A_\alpha$.*

We proceed on to recall the weak fuzzy direct sum.

**DEFINITION 1.4.7:** *Let $\{R_\alpha \mid \alpha \in \Omega\}$ be a collection of commutative rings and let $A_\alpha$ be a fuzzy subring of $R_\alpha$ for all $\alpha \in \Omega$.*

*Let*
$$\sum_{\alpha \in \Omega}^{\oplus} R_\alpha$$

*denote the weak direct sum of the $R_\alpha$. Define the fuzzy subset*
$$\sum_{\alpha \in \Omega}^{\oplus} A_\alpha \text{ of } \sum_{\alpha \in \Omega}^{\oplus} R_\alpha$$

*by for all*
$$x = \langle x \rangle \in \sum_{\alpha \in \Omega}^{\oplus} R_\alpha,$$

$$\left( \sum_{\alpha \in \Omega}^{\oplus} A_\alpha \right)(x) = \underset{\alpha \in \Omega}{\times} A_\alpha(x_\alpha).$$

*Then*
$$\sum_{\alpha \in \Omega}^{\oplus} A_\alpha$$

*is called the weak fuzzy direct sum of the $A_\alpha$.* Now we proceed on to recall some of the extensions of fuzzy subrings and fuzzy ideals as given by [84], that has given a necessary and sufficient conditions for a fuzzy subring or a fuzzy ideal $A$ of a commutative ring $R$ to be extended to one $A^e$ of a commutative ring $S$ containing $R$ as a subring such that $A$ and $A^e$ have the same image.

One of the applications of these results gives a criterion for a fuzzy subring of an integral domain R to be extendable to a fuzzy subfield of the quotient field. To this end we just recall the definition of fuzzy quasi-local subring and some basic properties about collection of ideals.

**THEOREM 1.4.48:** *Let R be a ring and let $\{I_t \mid t \in II\}$ be a collection of ideals of R such that $R = \cup I_t$ for all $s, t \in II$, $s > t$ if and only if $I_s \subset I_t$. Define the fuzzy subset A of R by for all $x \in R$, $A(x) = \sup \{t \mid x \in I_t\}$. Then A is a fuzzy ideal of R.*

*Proof:* Left for the reader as an exercise.

**Notation.** II is a non-empty subset of [0, 1] i.e. we use II as the index set.

On similar lines the following theorems can also be proved.

**THEOREM 1.4.49:** *Let R be a ring and let $\{R_t \mid t \in II\}$ be a collection of subrings of R such that $R = \cup R_t$ for all s, $t \in II$, $s > t$ if and only if $R_s \subset R_t$. Define the fuzzy subset A of R by, for all $x \in R$, $A(x) = \sup \{t \mid x \in R_t\}$. Then A is a fuzzy subring of R.*

**THEOREM 1.4.50:** *Let R be a ring with unity and let A be a fuzzy subring of R. Then for all $y \in R$, y a unit and for all $x \in R$. $A(xy^{-1}) \geq \min \{A(x), A(y)\}$ if and only if $A(y) = A(y^{-1})$ for every unit y in R. In either case $A(1) \geq A(y)$ where y is a unit.*

**THEOREM 1.4.51:** *Let R be a quasi-local ring. If A is fuzzy quasi local subring of R, then for all t such that $0 \leq t \leq A(t)$, $A_t$ is a quasi local ring and $M \cap A_t$ is the unique maximal ideal of $A_t$.*

*(We say that A fuzzy subring of R to be a fuzzy quasi local subring of R if and only if for all $x \in R$ and for all $y \in R$ such that y is a unit, $A(xy^{-1}) \geq \min \{A(x), A(y)\}$ or equivalently $A(y) = A(y^{-1})$).*

**THEOREM 1.4.52:** *Let R be a quasi-local ring. Let A be a fuzzy subset of R. If for all $t \in Im(A)$, $A_t$ is a quasi local ring with unique maximal ideal $M \cap A_t$ then A is fuzzy quasi local subring of R.*

For proof please refer [84] or the reader is advised to take it as a part of research and find the proof as the proof is more a matter of routine using basically the definitions.

**THEOREM 1.4.53:** *Let R be a quasi local ring.*

   i. *If A is a fuzzy quasi local subring of R, then $A^*$ is quasi local and $A^* \cap M$ is the unique maximal ideal of $A^*$.*

   ii. *If R' is a subring of R, then R' is quasi-local with unique maximal ideal $M \cap R'$ if and only if $\delta_R$, is a fuzzy quasi local subring of R.*

Next a necessary and sufficient condition for a fuzzy subset A of a ring R to be fuzzy quasi local is given.

**THEOREM 1.4.54:** *Let R be a quasi local ring. Let $\{R_t \mid t \in II\}$ be a collection of quasi local subrings of R with unique maximal ideals $M_t$ of $R_t$, $t \in II$}, such that $R = \cup R_t$ and for all s, $t \in II$, $s > t$ if and only if $R_s \subset R_t$. Define the fuzzy subset A of R by for all $x \in R$. $A(x) = \sup \{t \mid x \in R_t\}$. Then A is fuzzy quasi local if and only if $M_t = M \cap R_t$ for all $t \in II$.*

Before one defines the notion of extensions of fuzzy subrings and fuzzy ideals it is very essential to see the definition and properties of extension of fuzzy subsets. The following three theorems will give explicitly all properties about extension of fuzzy subsets. The proof of which is once again left for the reader as an exercise.

**THEOREM 1.4.55:** *Let R be a non-empty subset of a set S and let A be a fuzzy subset of R. If B is an extension of A to a fuzzy subset of S, then $A_t \cap B_s = A_s$ for all s, t such that $0 \le t \le s \le 1$.*

**THEOREM 1.4.56:** *Let R be a non-empty subset of a set S and let A be a fuzzy subset of R. If $B = \{B_t \mid t \in Im(A)\}$ is a collection of subsets of S such that*

  i. $\cup B_t = S$.
  ii. *for all $s, t \in Im(A)$, $s > t$ if and only if $B_s \subset B_t$ and*
  iii. *for all $s, t \in Im(A)$, $s \ge t$, $A_t \cap B_s = A_s$, then A has an extension to a fuzzy subset $A^e$ of S such that $(A^e)_t \supseteq B_t$ for all $t \in Im(A)$.*

**THEOREM 1.4.57:** *Let R be a non-empty subset of a set S and let A be a fuzzy subset of R such that A has the sup property. If $\mathcal{B} = \{B_t \mid t \in Im(A)\}$ is a collection of subsets of S which satisfies (i) to (iii) of conditions given in the above theorem then A has a unique extension to a fuzzy subset $A^e$ of S such that $(A^e)_t = B_t$ for all $t \in Im(A)$ and $Im(A^e) = Im(A)$.*

Now we proceed on to define extension of fuzzy subrings and fuzzy ideals. Let R be a subring of S. If I is an ideal of R, we let $I^e$ denote the ideal of S extended by I.

**THEOREM 1.4.58:** *Let R be a subring of S and let A be a fuzzy ideal of R such that A has the sup property. If*

$$\bigcup_{t \in Im(A)} (A_t)^e = S$$

*and for all $s, t \in Im(A)$, $s \ge t$, $A_t \cap (A_s)^e = A_s$, then A has a unique extension to a fuzzy ideal $A^e$ of S such that $(A^e)_t = (A_t)^e$ for all $t \in Im(A)$ and $Im(A^e) = Im(A)$. Let R be a commutative ring with identity. Let M be a multiplicative system in R. Let $N = \{x \in R \mid mx = 0 \text{ for some } m \in M\}$. Then N is an ideal of R. Unless otherwise specified, we assume $N = \langle 0 \rangle$ i.e. M is regular. Let $R_M$ denote the quotient ring of R with respect to M. Since $N = \langle 0 \rangle$, we can assume that $R \subset R_M$. If A is a fuzzy subring of R, we assume $A(1) = A(0)$.*

**THEOREM 1.4.59:** *Let A be a fuzzy subring of R such that A has the sup property. Then A can be extended to a fuzzy subring $A^e$ of $R_M$ such that for all $x, y \in R$, y a unit. $A^e(xy^{-1}) \ge \min\{A^e(x), A^e(y)\}$ if and only if for all $s, t \in Im(A)$, $s \ge t$, $A_t \cap (A_t)_M = A_s$ where $M_s = M \cap A_s$ for all $s \in Im(A)$. If either condition holds, $A^e$ can be chosen so that $(A^e)_t = (A_t)_M$ for all $t \in Im(A)$ and $Im(A^e) = Im(A)$.*

The following two results are very interesting and can be easily verified by the readers.

**THEOREM 1.4.60:** *Let A be a fuzzy subring of R such that A has the sup property. Then A can be extended to a fuzzy quasi local subring $A^e$ of $R_P \Leftrightarrow$ for all s, t $\in$ Im (A), s $\geq$ t, $A_t \cap (A_s)_P = A_s$ where $P_s = P \cap A_s$ for all s $\in$ Im (A).*

**THEOREM 1.4.61:** *Let R be an integral domain and let Q denote the quotient field of R. Let A be a fuzzy subring of R such that A has the sup property. Let $Q_t$ denote the smallest subfield of Q which contains $A_t$, for all t $\in$ Im (A). Then A can be extended to a fuzzy subfield of Q if and only if for all s, t $\in$ Im (A), s $\geq$ t, $A_t \cap Q_s = A_s$.*

Now we recall the definition of the extension of fuzzy prime ideals.

Let R and S be rings and let f be a homomorphism of R into S. Let T denote f (R). If I is an ideal in R, then the ideal $(f(I))^e$ (or simple $I^e$) is defined to be the ideal of S generated by f (I) and is called the extended ideal or extension of I. If J is an ideal of S, the ideal $J^c = f^{-1}(J)$ is called the contracted ideal or the contraction of J.

**DEFINITION [80]:** *Let A and B be fuzzy subsets of R and T respectively. Define the fuzzy subsets f(A) of T and $f^{-1}$(B) of R by f(A) (y) = sup {A(x) | f(x) = y} for all y $\in$ T, $f^{-1}$ (B) (x) = B (f(x)) for all x $\in$ R .*

**THEOREM 1.4.62:** *Suppose A and B are fuzzy ideals of R and T respectively. Then*

  i. *f (A) and $f^{-1}$(B) are fuzzy ideals of T and R respectively.*
  ii. *f (A) (0) = A(0).*
  iii. *$f^{-1}$ (B) (0) = B(0).*

**THEOREM 1.4.63:** *Let A be a fuzzy ideal of R. Then*

  i. *$f(A_*) \subset (f(A))_*$*
  ii. *if A has the sup property then $f(A))_* = f(A_*)$.*

Now we proceed on to recall the definition of f-invariant.

**DEFINITION 1.4.8:** *Let A be a fuzzy ideal of R. A is called f-invariant if and only if for all x, y $\in$ R, f(x) = f(y) implies A(x) = A(y).*

In view of this we have the following nice characterization theorem; the proof of which is left an exercise for the reader.

**THEOREM 1.4.64:** *Let A be a fuzzy ideal of R. Then A is a fuzzy prime ideal of R if and only if A(0) = 1, $|Im(A)| = 2$ and $A_*$ is a prime ideal of R.*

**THEOREM 1.4.65:** *Let A be an f-invariant fuzzy ideal of R such that A has the sup property. If $A_*$ is a prime ideal of R, then $f(A_*)$ is a prime ideal of T.*

**THEOREM 1.4.66:** *Let A be an f-invariant fuzzy ideal of R such that Im(A) is finite. If $A_*$ is a prime ideal of R then $f(A_*) = (f(A))_*$ .*

**THEOREM 1.4.67:** *Let A be an f-invariant fuzzy ideal of R. If A is a fuzzy prime ideal of R then f(A) is a fuzzy prime ideal of T.*

**THEOREM 1.4.68:** *Let B be a fuzzy ideal of T. Then*

   i. $f^{-1}(B_*) = f^{-1}(B)_*$.
   ii. *If $B_*$ is a prime ideal of T, then $f^{-1}(B)_*$ is a prime ideal of R,*
   iii. *If B is a fuzzy prime ideal of T, then $f^{-1}(B)$ is a fuzzy prime ideal of R.*

The following are simple but interesting results on fuzzy ideals.

**THEOREM 1.4.69:** *Let I be a fuzzy ideal of S. Then $(I^c)_* = (I_*)^c$.*

To prove the following theorem we make the assumption that if M is the multiplicative system in R. N = {x∈ R | xm = 0 for some m ∈ M} equals ⟨0⟩.

**THEOREM 1.4.70:** *Let I be a fuzzy ideal of $R_M$ then $(I^{ce})_* = ((I^c)_*)^e = (I_*)^{ce} = I_*$.*

**THEOREM 1.4.71:** *Let A be a fuzzy ideal of R. Then in $R_M$, $(A^{ce})_* = ((A^e)_*)^c = (A_*)^{ec}$ and if M is prime to $A_*$, then $(A_*)^{ec} = A_*$.*

The reader is requested to refer [80].

We replace the interval [0, 1] by a finite lattice L which has 0 to be the least element and 1 to be the largest element. All the while, fuzzy ideals have been defined over [0, 1] when we define it over a lattice L we call them L-fuzzy ideal.

**DEFINITION [61-64]:** *An L-fuzzy ideal is a function $J : R \to L$ (R is a commutative ring with identity L stands for a lattice with 0 and 1) satisfying the following axioms*

   i. $J(x + y) \geq J(x) \wedge J(y)$.
   ii. $J(-x) = J(x)$.
   iii. $J(xy) \geq J(x) \vee J(y)$.

**THEOREM [61-64]:**

   i. *A function $J : R \to L$ is a fuzzy ideal if and only if*

   $$J(x - y) \geq J(x) \wedge J(y) \text{ and}$$
   $$J(xy) \geq J(x) \vee J(y).$$

   ii. *If $J: R \to L$ is a fuzzy ideal then*

   (a) $J(0) \geq J(x) \geq J(1)$ for all $x \in R$.
   (b) $J(x-y) = J(0)$ implies $J(x) = J(y)$ for all $x, y \in R$
   (c) *The level cuts $J_\alpha = \{x \in R \mid J(x) \geq \alpha\}$ are ideals of R. Conversely if each $J_\alpha$ is an ideal then J is a fuzzy ideal.*

*Proof:* Please refer [61-64]. A strict level cut $J_\alpha = J_\alpha = \{x \in R \mid J(x) > \alpha\}$ need not be an ideal unless of course L is totally ordered.

The following result can also be obtained as a matter of routine.

**THEOREM 1.4.72:** *If $f : R \to R'$ is a homomorphism of rings and $J : R \to L$ and $J' : R' \to L$ are fuzzy ideals, then*

i. $f^{-1}(J')$ is a fuzzy ideal which is constant on ker f and
ii. $f^{-1}(J'_{(J')(0')}) = f^{-1}\left(J'_{f^{-1}(J')(0')}\right)$.
iii. $f(J)$ is a fuzzy ideal.
iv. $ff^{-1}(J') = J'$.
v. if J is constant on ker f then $f\left(J_{J(0)}\right) = f(J)_{f(J)(0')}$.
vi. If J is constant on ker f then $f^{-1}(f(J)) = J$.

**THEOREM (CORRESPONDENCE THEOREM):** *If $f: R \to R'$ is an epimorphism of rings, then there is one to one correspondence between the ideals of R' and those of R which are constant on ker f. If J is a fuzzy ideal of R which is constant on ker f, then f (J) is the corresponding fuzzy ideal of R'. If J' is a fuzzy ideal of R', then $f^{-1}(J')$ is the corresponding fuzzy ideal of R.*

*Proof:* Refer [61].

Now we proceed on to recall the definitions of prime fuzzy ideals, primary fuzzy ideals and semiprime fuzzy ideals and also some of it basic properties. For more about these concepts please refer [61-64].

By a prime fuzzy ideal we mean a non-constant fuzzy ideal $P : R \to L$ satisfying the following condition of primeness

$P(xy) = P(x)$ or
$P(xy) = P(y)$ for all $x, y \in R$.

**THEOREM 1.4.73:** *If $P: R \to L$ is a prime fuzzy ideal, then the set P(R) of membership values of P is a totally ordered set with the least element P(1) and the greatest element P(0).*

**THEOREM 1.4.74:** *A fuzzy ideal $P: R \to L$ is prime if and only if every level cut $P_\alpha = \{x \in R \mid P(x) \geq \alpha\}$ is prime for all $\alpha > P(1)$ For $\alpha = P(1)$, $P_\alpha = R$.*

**THEOREM 1.4.75:** *Let Z be a non-empty subset of R. Z is a prime ideal of R if and only if $\chi_z : R \to L$ is a prime fuzzy ideal.*

**THEOREM 1.4.76:** *Let R be a principal ideal domain (PID). If $P : R \to L$ is a prime fuzzy ideal and $P_{P(0)} \neq 0$, then P (R) has two elements. P is properly fuzzy if and only if P(R) has three elements. We see a properly fuzzy prime ideal of a PID R is equivalent to the fuzzy ideal $P : R \to L$ of the following type:*

$$P(0) = 1,$$
$$P(x) = \alpha$$

for all $x \in P_1 \setminus \{0\}$. $P(x) = 0$ for all $x \in R \setminus P_1$ where $P_1$ is a prime ideal of $R$ and $0 < \alpha < 1$.

**DEFINITION 1.4.9:** *A finite strictly increasing sequence of prime ideals of a ring R, $P_0 \subset P_1 \subset P_2 \subset \ldots \subset P_n$ is called a chain of prime ideals of length n. The supremum of the lengths of all chains is called the dimension of R.*

**DEFINITION 1.4.10:** *Let R be a ring. Then $\vee\{ |P(R)| \mid P: R \to [0, 1]$ is a prime fuzzy ideal$\}$ is called the fuzzy dimension of R.*

**THEOREM 1.4.77:**

   i. *The dimension of R is n $(< \infty)$ if and only if its fuzzy dimension is n + 2.*
   ii. *An artinian ring has no properly fuzzy prime ideal.*
   iii. *A Boolean ring has no properly fuzzy prime ideal.*

The following results are relation on homomorphism and epimorphism of rings.

**THEOREM 1.4.78:**

   i. *Let $f : R \to R'$ is an epimorphism of rings. If $P: R \to L$ is a prime fuzzy ideal which is constant on ker f then $f(P) : R' \to L$ is a prime fuzzy ideal.*
   ii. *If $f : R \to R'$ is a homomorphism of rings. If $P' : R' \to L$ is a fuzzy prime ideal then $f^{-1}(P')$ is a prime fuzzy ideal of R.*
   iii. *Let $f : R \to R'$ be an epimorphism of ring.*

     (a) *Let $P : R \to L$ be a fuzzy ideal which is constant on ker f. Then P is prime if and only if $f(P) : R' \to L$ is prime.*
     (b) *Let $P' : R' \to L$ be a fuzzy ideal. Then P' is prime if and only if $f^{-1}(P') : R \to L$ is prime.*

The proof can be got by simple computations and hence left for the reader as exercise. Now we proceed on to give the definitions of primary fuzzy ideals of a ring.

**DEFINITION [149]:** *A fuzzy ideal $Q : R \to L$ is called primary if $Q(xy) = Q(0)$ implies $Q(x) = Q(0)$ or $Q(y^n) = Q(0)$ for some integer $n > 0$.*

**DEFINITION [77]:** *A fuzzy ideal $Q : R \to L$ is called primary, if either Q is the characteristic function of R or*

   i. *Q is non-constant and*
   ii. *$A \circ B \subseteq Q \Rightarrow A \subset Q$ or $B \subseteq \sqrt{Q}$ is the intersection of all prime fuzzy ideals.*

**DEFINITION [146]:** *A fuzzy ideal $Q : R \to L$ is called primary, if*

i. $Q$ is non-constant and
ii. for all $x, y \in R$ and $r, s \in L$ if $x, y \in Q$ then $x_r \in Q$ or $y_s^n \in Q$ for some positive integer $n$.

**DEFINITION 1.4.11:** *A fuzzy ideal $Q : R \to L$ is called primary if $Q$ is non-constant and for all $x, y \in R$, $Q(xy) = Q(x)$ or $Q(y^n)$ for some positive integer $n$.*

**THEOREM 1.4.79:**

i. Let $Q$ be an ideal of $R$. The characteristic function $\chi_Q$ is a primary fuzzy ideal if and only if $Q$ is a primary ideal.
ii. If $Q : R \to L$ is primary then its level cuts $Q_\alpha = \{x \in R \mid Q(x) \geq \alpha\}$, $\alpha \in L$, are primary.
iii. Every prime fuzzy ideal is primary.

*Proof:* Left for the reader to prove.

As in case of prime ideals and its relation with epimorphism and homomorphism, we give here without proof the relation of primary ideals and the epimorphism and homomorphism.

**THEOREM 1.4.80:**

i. Let $f : R \to R'$ be an epimorphism of rings. If $Q : R \to L$ is a primary fuzzy ideal of $R$ which is constant on ker $f$, then $f(Q)$ is a primary fuzzy ideal of $R'$.
ii. Let $f : R \to R'$ be a homomorphism of rings. If $Q' : R' \to L$ is a primary fuzzy ideal of $R'$ then $f^{-1}(Q')$ is a primary fuzzy ideal of $R$.
iii. Let $f : R \to R'$ be an epimorphism of rings and $Q: R \to L$ and $Q' : R' \to L$ be a fuzzy ideals.
   (a) $Q'$ is primary if and only if $f^{-1}(Q')$ is primary.
   (b) If $Q$ is constant on ker $f$, then $Q$ is primary if and only if $f(Q)$ is primary.

Now we proceed on to define the notion of weak primary fuzzy ideals.

**DEFINITION 1.4.12:** *A fuzzy ideal $J: R \to L$ is said to be weak primary or in short w-primary if $J(xy) = J(x)$ or $J(xy) \leq J(y^n)$ for some integer $n > 0$.*

**THEOREM 1.4.81:** *Every primary fuzzy ideal is w-primary. In particular, every prime fuzzy ideal is w-primary.*

*Proof:* Direct by the very definitions, hence left for the reader to prove.

We give the following two nice characterization theorems. The proofs are left for the reader as an exercise.

**THEOREM 1.4.82:** *A fuzzy ideal is w-primary if and only if each of its level cuts is primary.*

**THEOREM 1.4.83:** *Let Q be an ideal of R. The characteristic function $\chi_Q$ is w-primary if and only if Q is primary.*

Finally we give the following result the proof of which is easy by simple computations.

**THEOREM 1.4.84:** *Let $f : R \to R'$ be a homomorphism of rings, and $Q : R \to L$ and $Q' : R' \to L$ be fuzzy ideals.*

   i. *If $\phi$ is w-primary then so is $f^{-1}(Q')$.*

   ii. *Let f be an epimorphism. Then Q is w-primary if and only if f(Q) is w-primary.*

   iii. *Let f be an epimorphism then Q' is w-primary if and only if $f^{-1}(Q)$ is w-primary.*

Let I be an ideal of R, nil-radical defined as $\sqrt{I} = \{x \in R \mid x^n \in I, n > 0\}$.

**DEFINITION [61-64]:** *If $J : R \to L$ is a fuzzy ideal, then the fuzzy set $\sqrt{J} : R \to L$ defined as $\sqrt{J}(x) = \vee \{J(x^n) \mid n > 0\}$ is called the fuzzy nil radical of J.*

This is proved by simple techniques.

**THEOREM 1.4.85:**

   i. *If $J : R \to L$ is a fuzzy ideal, then $\sqrt{J}$ is a fuzzy ideal.*
   ii. *If I is an ideal of R, then $\sqrt{\chi_I} = \chi_{\sqrt{I}}$.*
   iii. *For any $0 \leq \alpha < 1$ and a fuzzy ideal $J : R \to L$, $(\sqrt{J})_\alpha = \sqrt{J_\alpha}$ where L is a totally ordered set, $J_\alpha = (x \in R \mid J(x) > \alpha)$ and $(\sqrt{J})_\alpha = \{x \in R / \sqrt{J(x)} > \alpha\}$.*
   iv. *In case of non-strict level cuts $\sqrt{J_\alpha} \subseteq (\sqrt{J})_\alpha$.*

Using these results and definitions the reader is assigned the task of proving the following theorem.

**THEOREM 1.4.86:**

   i. *If $f : R \to R'$ is an epimorphism of rings and $J : R \to L$ is a fuzzy ideal, then $f(\sqrt{J}) \subseteq \sqrt{f(J)}$. Further if J is constant on ker f then $f(\sqrt{J}) = \sqrt{f(J)}$.*

ii. If $f : R \to R'$ is a homomorphism of rings and $J' : R' \to L$ is a fuzzy ideal then $f^{-1}\left(\sqrt{J'}\right) = \sqrt{f^{-1}(J')}$.

The following theorem is a direct consequence of the definition.

**THEOREM 1.4.87:** *If $J : R \to L$ and $K : R \to L$ are fuzzy ideals, then the following hold:*

i. $\sqrt{\left(\sqrt{J}\right)} = \sqrt{J}$.
ii. If $J \subseteq K$, then $\sqrt{J} \subseteq \sqrt{K}$.
iii. $\sqrt{J \cap K} = \sqrt{J} \cap \sqrt{K}$.
iv. If $J : R \to L$ is a fuzzy ideal with supremum property then $\sqrt{J_\alpha} = \left(\sqrt{J}\right)_\alpha$.
v. If $P : R \to L$ is prime then $\sqrt{P} = P$.
vi. If $Q : R \to L$ is a primary fuzzy ideal with supremum property then $\sqrt{Q}$ is the smallest prime fuzzy ideal containing Q.
vii. If L is a totally ordered set and $Q : R \to L$ is a primary fuzzy ideal, then $\sqrt{Q}$ is the smallest prime fuzzy ideal containing Q.

Now we recall the notion of prime fuzzy radical.

**DEFINITION 1.4.13:** *Let $J : R \to L$ be a fuzzy ideal and $P : R \to L$ denote a prime fuzzy ideal containing J. The fuzzy ideal $r(J) = \cap \{P / J \subseteq P\}$ is called the prime fuzzy radical of J.*

Using the definitions the reader is expected to prove the following two theorems

**THEOREMS 1.4.88:**

i. If $J: R \to L$ is a fuzzy ideal, then $\sqrt{J} \subseteq r(J)$.
ii. If L is a totally ordered set and $J: R \to L$ is a fuzzy ideal then $\sqrt{J} \subseteq r(J)$.

**DEFINITION 1.4.14:** *A fuzzy ideal $S: R \to L$ is called semiprime fuzzy ideal if $S(x^2) = S(x)$ for all $x \in R$.*

The following results are left as an exercise for the reader to prove.

**THEOREM 1.4.89:**

i. Let $S : R \to L$ be a fuzzy ideal, S is semiprime if and only if its level cuts, $S_\alpha = (x \in R | S(x) \geq \alpha)$ are semiprime ideals of R, for all $\alpha \in L$.
ii. Let S be an ideal of R. S is semiprime if and only if its characteristic function $\chi_S$ is a semiprime fuzzy ideal of R.
iii. Let $f : R \to R'$ be a homomorphism. If $S' : R' \to L$ is a semiprime fuzzy ideal of R then $f^{-1}(S')$ is a semiprime fuzzy ideal of R.

iv. Let $f : R \to R'$ be an epimorphism and $S : R \to L$ be a semiprime fuzzy ideal of R which is constant on ker f. Then f(S) is a semiprime fuzzy ideal of R'. Thus by the correspondence theorem between semiprime fuzzy ideals of R' and those of R which are constant on the kernel of f.
v. Every prime fuzzy ideal is semiprime fuzzy ideal.
vi. Intersection of semiprime fuzzy ideal is a semiprime fuzzy ideal. In particular intersection of prime fuzzy ideals is a semiprime fuzzy ideal.
vii. If $S : R \to L$ is a semiprime fuzzy ideal, then the quotient ring R/S is prime.

Now we formulate certain equivalent condition of a fuzzy ideal. The proof of these relations is omitted.

**THEOREM 1.4.90:** *If $S: R \to L$ is a fuzzy ideal then the following are equivalent*

i. *S is semiprime.*
ii. *Each level cut of S is semiprime.*
iii. *$S(x^n) = x$ for all integers $n > 0$ and $x \in R$.*
iv. *$J^2 \subseteq S$ implies $J \subseteq S$ for all fuzzy ideals $J : R \to L$.*
v. *$J^n \subseteq S$ for $n > 0$ implies $J \subseteq S$ for all fuzzy ideals $J : R \to L$.*
vi. *$S = \sqrt{S}$ where $\sqrt{S}$ is the fuzzy nil radical of S when L is totally ordered each of the above statements is equivalent to the following:*

   a. *S coincides with its prime fuzzy radical.*
   b. *$S = \cap \{ P \,/\, P \in C\}$ where C is a class of prime fuzzy ideals.*

## 1.5 Fuzzy birings

In this section we just introduce the notion of birings and fuzzy birings. The very concept of birings is very new [135] so the notion is fuzzy birings is defined only in this book. So we first give some of the basic properties of birings and then give the fuzzy analogue of it. As the main aim of this book is to introduce the notion of Smarandache fuzzy algebraic structures we do not stake in discussing elaborately the concepts other than Smarandache structures.

**DEFINITION 1.5.1:** *A non-empty set $(R, +, \bullet)$ with two binary operations '+' and '$\bullet$' is said to be a biring if $R = R_1 \cup R_2$ where $R_1$ and $R_2$ are proper subsets of R and*

i. *$(R_1, +, \bullet)$ is a ring.*
ii. *$(R_2, +, \bullet)$ is a ring.*

*Example 1.5.1:* Let R = {0, 2, 4, 6, 8, 10, 12, 14, 16, 3, 9, 15} be a non-empty set. (R, +, •) where '+' and '•' are usual addition and multiplication modulo 18. Take $R_1$ = {0, 2, 4, 6, 8, 10, 12, 14, 16, } and $R_2$ = {0, 3, 6, 9, 15, 12}; clearly ($R_1$, +, •) and ($R_2$, +, •) are rings.

A biring R is said to be finite if R contains only finite number of elements. If R has infinite number of elements then we say R is of infinite order. The biring R given in example 1.5.1 is of finite order and has 12 elements. We denote the order of R by $|R|$ or o(R).

**DEFINITION 1.5.2:** *A biring $R = R_1 \cup R_2$ is said to be a commutative biring if both $R_1$ and $R_2$ are commutative rings. Even if one of $R_1$ or $R_2$ is not a commutative ring then we say the biring is a non-commutative biring. We say the biring R has a monounit if a unit exists which is common to both $R_1$ and $R_2$. If $R_1$ and $R_2$ are rings which has separate unit then we say the biring $R = R_1 \cup R_2$ is a biring with unit.*

It is interesting to note that the biring given in example 1.5.1 has no units but is a commutative biring of finite order.

*Example 1.5.2:* Let $R_{2\times 2}$ denote the set of all $2 \times 2$ upper triangular matrices with entries from the ring of integers Z i.e.

$$R_{2\times 2} = \left\{ \begin{pmatrix} a & b \\ 0 & c \end{pmatrix}, \begin{pmatrix} x & 0 \\ y & z \end{pmatrix} \mid x, y, z, a, b, c \in Z \right\}.$$

R is an infinite non-commutative biring with monounit.

For take $R = R_1 \cup R_2$ where

$$\left\{ \begin{pmatrix} a & b \\ 0 & c \end{pmatrix} \mid a, b, c \in Z \right\}$$

is a non-commutative ring with

$$I_{2\times 2} = \begin{pmatrix} 1 & 0 \\ 0 & 1 \end{pmatrix}$$

as the unit

$$R_2 = \left\{ \begin{pmatrix} x & 0 \\ y & z \end{pmatrix} \mid x, y, z \in Z \right\}$$

is a non-commutative ring with

$$I_{2\times 2} = \begin{pmatrix} 1 & 0 \\ 0 & 1 \end{pmatrix}$$

as the unit. Thus $R = R_1 \cup R_2$ is a biring with $I_{2\times 2}$ as a monounit.

**DEFINITION 1.5.3:** *Let $R = R_1 \cup R_2$ be a biring. A non-empty subset S of R is said to be a sub-biring of R if $S = S_1 \cup S_2$ and S itself a biring and $S_1 = R_1 \cap S$ and $S_2 = R_2 \cap S$.*

*Example 1.5.3:* Let $R = R_1 \cup R_2$ be a biring given in example 1.5.2. Take $S = S_1 \cup S_2$ where

$$S_1 = \left\{ \begin{pmatrix} 0 & 0 \\ a & 0 \end{pmatrix} \mid a \in Z \right\}$$

and

$$S_2 = \left\{ \begin{pmatrix} 0 & a \\ 0 & 0 \end{pmatrix} \mid a \in Z \right\}.$$

Clearly S is a sub-biring of R.

The following theorem is very important which is a characterization theorem.

**THEOREM 1.5.1:** *Let R be a biring where $R = R_1 \cup R_2$. A non-empty subset $S = S_1 \cup S_2$ of R is a subring of R if and only if $R_1 \cap S = S_1$ and $R_2 \cap S = S_2$ are sub-biring of $R_1$ and $R_2$ respectively.*

*Proof:* Straightforward by the very definitions.

**DEFINITION 1.5.4:** *Let $R = R_1 \cup R_2$ be a biring. A non-empty subset I of R is said to be a right bi-ideal of R if $I = I_1 \cup I_2$ where $I_1$ is a right ideal of $R_1$ and $I_2$ is a right ideal of $R_2$. I is said to be a left bi-ideal of R if $I = I_1 \cup I_2$ are left ideals of $R_1$ and $R_2$ respectively.*

If $I = I_1 \cup I_2$ is such that both $I_1$ and $I_2$ are ideals of $R_1$ and $R_2$ respectively then we say I is a bi-ideal of R, then we say I is a bi-ideal of R. Now it may happen when $I = I_1 \cup I_2$, $I_1$ may be a right ideal of $R_1$ and $I_2$ may be a left ideal of $R_2$ then how to define ideal structures.

For this case we give the following definition.

**DEFINITION 1.5.5:** *Let $R = R_1 \cup R_2$ be a biring. We say the set $I = I_1 \cup I_2$ is a mixed bi-ideal of R if $I_1$ is a right ideal of $R_1$ and $I_2$ is a left ideal of $R_2$. Thus we see only in case of birings we can have the concept of mixed bi-ideals i.e. an ideal simultaneously having a section to be a left ideal and another section to be right ideal.*

*Example 1.5.4:* Let $R = R_1 \cup R_2$ be a biring given in example 1.5.2 Take $I = I_1 \cup I_2$ where

$$I_1 = \left\{ \begin{pmatrix} a & 0 \\ b & 0 \end{pmatrix} \mid a, b \in Z \right\}$$

$$I_2 = \left\{ \begin{pmatrix} 0 & x \\ 0 & 0 \end{pmatrix} \mid x \in Z \right\}$$

are left ideals of $R_1$ and $R_2$ respectively. Thus I is a left bi-ideal of R.

Now take J = J₁ ∪ J₂ where

$$J_1 = \left\{ \begin{pmatrix} 0 & 0 \\ a & b \end{pmatrix} \mid a, b \in Z \right\}$$

is a right ideal of $R_1$ and

$$J_2 = \left\{ \begin{pmatrix} 0 & b \\ 0 & 0 \end{pmatrix} \mid b \in Z \right\}$$

is right ideal of $R_2$ so J is a right bi-ideal of R.

Take K = K₁ ∪ K₂ where

$$K_1 = \left\{ \begin{pmatrix} a & 0 \\ b & 0 \end{pmatrix} \mid a, b \in Z \right\}$$

and

$$K_2 = \left\{ \begin{pmatrix} 0 & b \\ 0 & 0 \end{pmatrix} \mid b \in Z \right\}$$

are left and right ideals of $R_1$ and $R_2$ respectively; so K is a mixed bi-ideal of R.

**DEFINITION 1.5.6:** *Let $R = R_1 \cup R_2$ be a biring. A bi-ideal $I = I_1 \cup I_2$ is called a maximal bi-ideal of R if $I_1$ is a maximal ideal of $R_1$ and $I_2$ is a maximal ideal of $R_2$. Similarly we can define the concept of minimal bi-ideal as $J = J_1 \cup J_2$ is a minimal ideal if $J_1$ is a minimal ideal of $R_1$ and $J_2$ is a minimal ideal of $R_2$. It may happen in a bi-ideal. I of a ring $R = R_1 \cup R_2$ that one of $I_1$ or $I_2$ may be maximal or minimal then what do we call the structure $I = I_1 \cup I_2$. We call $I = I_1 \cup I_2$ a bi-ideal in which only $I_1$ or $I_2$ is a maximal ideal as quasi maximal bi-ideal. Similarly we can define quasi minimal bi-ideal.*

*We call a bi-ideal $I = I_1 \cup I_2$ to be a prime bi-ideal of $R = R_1 \cup R_2$ the biring if both $I_1$ and $I_2$ are prime ideals of the rings $R_1$ and $R_2$ respectively.*

**DEFINITION 1.5.7:** *Let $R = R_1 \cup R_2$ and $S = S_1 \cup S_2$ be two subrings. We say a map $\phi$ from R to S is a biring homomorphism if $\phi = \phi_1 \cup \phi_2$ where $\phi_1 = \phi \mid R_1$ from $R_1$ to $S_1$ is a ring homomorphism and $\phi_2 : \phi \mid R_2$ is a map from $R_2$ to $S_2$ is a ring homomorphism. We for notational convenience denote by $\phi = \phi_1 \cup \phi_2$ though this union '∪' is not the set theoretic union. We define for the homomorphism $\phi : R \to S$ where $R = R_1 \cup R_2$ and $S = S_1 \cup S_2$ are birings the kernel of the homomorphism $\phi$ as biker $\phi$ = ker $\phi_1$ ∪ ker $\phi_2$ here ker $\phi_1 = \{a_1 \in R_1 \mid \phi_1(a_1) = 0\}$ and ker $\phi_2 = \{a_2 \in R_2 \mid \phi_2(a_2) = 0\}$.*

The following theorem is straightforward hence left for the reader to prove.

**THEOREM 1.5.2:** *Let $R = R_1 \cup R_2$ and $S = S_1 \cup S_2$ be two birings and $\phi$ a biring homomorphism from R to S then biker $\phi$ = ker $\phi_1 \cup$ ker $\phi_2$ is a bi-ideal of the biring R.*

We may have several other interesting results but we advise the reader to refer [135]. Now we proceed on to define fuzzy birings and give some basic and important results about them.

**DEFINITION 1.5.8:** *Let $(R = R_1 \cup R_2, +, \bullet)$ be a biring. The map $\mu : R \to [0, 1]$ is said to be a fuzzy sub-biring of the biring R if there exists two fuzzy subsets $\mu_1$ (of $R_1$) and $\mu_2$ (of $R_2$) such that*

   i. *$(\mu_1, +, \bullet)$ is a fuzzy subring of $(R_1, +, \bullet)$*
   ii. *$(\mu_2, +, \bullet)$ is a fuzzy subring of $(R_2, +, \bullet)$*
   iii. *$\mu = \mu_1 \cup \mu_2$.*

**THEOREM 1.5.3:** *Every t-level subset of a fuzzy sub-biring $\mu$ of a biring $R = R_1 \cup R_2$ need not in general be a sub-biring of the biring R.*

*Proof:* The reader is requested to prove by constructing a counter example.

**DEFINITION 1.5.9:** *Let $(R = R_1 \cup R_2, +, \bullet)$ be a biring and $\mu$ $(= \mu_1 \cup \mu_2)$ be a fuzzy sub-biring of the biring R. The bilevel subset of the fuzzy sub-biring $\mu$ of the biring R is defined as $G_\mu^t = G_{1\mu_1}^t \cup G_{2\mu_2}^t$ for every $t \in \{0, \min \{\mu_1(0), \mu_2(0)\}$.*

The condition $t \in \{0, \min \{\mu_1(0), \mu_2(0)\}$. is essential for the bilevel subset be a sub-biring for if $t \notin \{0, \min \{\mu_1(0), \mu_2(0)\}$ then the bilevel subset-need not in general be a sub-biring of the biring R.

**DEFINITION 1.5.10:** *A fuzzy subset $\mu$ of a ring R is said to be a fuzzy sub-biring of the ring R if there exists two fuzzy subrings $\mu_1$ and $\mu_2$ of $\mu$ ($\mu_1 \neq \mu$ and $\mu_2 \neq \mu$) such that $\mu = \mu_1 \cup \mu_2$.*

**THEOREM 1.5.4:** *Every fuzzy sub-biring of a ring R is a fuzzy subring of the ring R and not conversely.*

*Proof*: Follows by the very definitions.

It is however left for the student to prove that the converse of the above theorem is not true.

**DEFINITION 1.5.11:** *A fuzzy subset $\mu = \mu_1 \cup \mu_2$ of a biring $R = R_1 \cup R_2$ is called a fuzzy bi-ideal of R if and only if $\mu_1$ is a fuzzy ideal of $R_1$ and $\mu_2$ is a fuzzy ideal of $R_2$.*

**DEFINITION 1.5.12:** *A fuzzy subset $\mu = \mu_1 \cup \mu_2$ of a biring $R = R_1 \cup R_2$ is a fuzzy sub-biring (fuzzy bi-ideal) of R if and only if the level subsets $\mu_t$, $t \in Im (\mu) = Im\mu_1 \cup Im\mu_2$ are subrings or ideals of $R_1$ and $R_2$.*

**DEFINITION 1.5.13:** *Let $\mu$ and $\theta$ be any two fuzzy bi-ideals. $R = R_1 \cup R_2$, ($\mu = \mu_1 \cup \mu_2$ and $\theta = \theta_1 \cup \theta_2$). The product $\mu \circ \theta$*

*of $\mu$ and $\theta$ is defined by*

$$(\mu \circ \theta) = \sup_{\substack{x_1 = \sum y_i z_i \\ i < \infty}} (\min_i (\min \mu_1(y_1), \theta_1(z_1))) \cup \sup_{\substack{x_2 = \sum y_j z_j \\ j < \infty}} (\min (\mu_2(y_j), \theta_2(z_j)))$$

*It can be easily verified that $\mu \circ \theta$ is a fuzzy bi-ideal of $R = R_1 \cup R_2$.*

**DEFINITION 1.5.14:** *A fuzzy bi-ideal $\mu = \mu_1 \cup \mu_2$ of a biring $R = R_1 \cup R_2$ is called fuzzy prime if, the bi-ideal $\mu_t = (\mu_1)_t \cup (\mu_2)_t$ where $t = \mu_1(0)$ and $t = \mu_2(0)$ is a prime ideal of $R_1$ and $R_2$ respectively.*

**DEFINITION 1.5.15:** *A non-constant fuzzy bi-ideal $\mu = \mu_1 \cup \mu_2$ of a biring $R = R_1 \cup R_2$ is called fuzzy prime if, for any two bi-ideals $\sigma$ and $\theta$ of $R$ the condition $\sigma \theta \subseteq \mu$ either $\sigma \subset \mu$ or $\theta \subset \mu$ (where $\theta = \theta_1 \cup \theta_2$ and $\sigma = \sigma_1 \cup \sigma_2$ by $\sigma \subset \mu$ we mean $\sigma_1 \subset \mu_1$ and $\sigma_2 \subset \mu_2$ similarly $\sigma \theta \subset \mu$ implies $\sigma_1 \theta_1 \subset \mu_1$ and $\sigma_2 \theta_2 \subset \mu_2$).*

**THEOREM 1.5.5:**

i. *If $\mu = \mu_1 \cup \mu_2$ is a fuzzy prime bi-ideal of a biring $R = R_1 \cup R_2$ then the ideal $\mu_t = (\mu_1 \cup \mu_2)_t = (\mu_1)_t \cup (\mu_2)_t$, $t = \mu_1(0)$ that is $t = \mu_1(0)$ and $t = \mu_2(0)$ is a prime bi-ideal of $R = R_1 \cup R_2$.*

ii. *A bi-ideal $P = P_1 \cup P_2$ of a biring $R = R_1 \cup R_2$, $P \neq R$ is prime if and only if $\chi_P$ ($\chi_P$ is the characteristic function of $P$ i.e. $\chi_P = (\chi_{P_1} \cup \chi_{P_2})$ is a fuzzy prime bi-ideal of $R$.*

**THEOREM 1.5.6:** *A non-constant fuzzy bi-ideal $\mu = \mu_1 \cup \mu_2$ of a biring $R = R_1 \cup R_2$ is fuzzy prime if and only if card Im $\mu_1 = 2$ and card Im$\mu_2 = 2$, $1 \in $ Im $\mu_1$ and $1 \in $ Im $\mu_2$ and the ideal $(\mu_1)_t$ and $(\mu_2)_t$ where $t = \mu_1(0)$ and $t = \mu_2(0)$ is prime.*

*Proof:* Follows as in case of ideals as each of $\mu_i$, $i = 1, 2$ are non-constant fuzzy ideals so is $\mu = \mu_1 \cup \mu_2$.

**DEFINITION 1.5.16:** *Let $R = R_1 \cup R_2$ be a biring. Let $\mu_s$ and $\mu_t$ be two level sub-birings and level bi-ideals (with $s < t$) of a fuzzy sub-biring (fuzzy bi-ideal), $\mu$ of a biring $R$ are equal if and only if there is no $x$ in $R$ such that $s \leq \mu(x) < t$ (that is $s < \mu_1(x_1) < t$, $s < \mu_2(x_2) < t$, $x_1 \in R_1$ and $x_2 \in R_2$).*

The following theorem is left as an exercise for the reader.

**THEOREM 1.5.7:** *Intersection of any family of fuzzy sub-birings (fuzzy bi-ideals) of a ring $R$ is a fuzzy sub-biring (fuzzy bi-ideal) of $R$.*

*Proof:* $\mu_i : R \to [0, 1]$ where $R = R_1 \cup R_2$ and $\mu_i = \mu_{i1} \cup \mu_{i2}$ each $\mu_{i1}$ is a fuzzy subring of $R_1$ $\mu_{i2}$ is a fuzzy sub-biring of $R_2$ then by usual results on fuzzy ideals we see intersection of $\mu_{ij}$ is a fuzzy subring of fuzzy ideal for $j = 1, 2$.

**THEOREM 1.5.8:** *If $\mu$ is a fuzzy bi-ideal of a biring $R = R_1 \cup R_2$. Then $\mu + \mu = \mu$.*

*Proof:* The result follows as in case of fuzzy ideal.

**THEOREM 1.5.9:** *If $\mu$ is any fuzzy sub-biring and $\theta$ is any fuzzy bi-ideal of a biring $R = R_1 \cup R_2$. Then $\mu \cap \theta$ is a fuzzy bi-ideal of the biring $\mu_t = \{ x \in R \mid \mu(x) = \mu(0)\}$. $[\mu(x) = \mu(0)$ will mean for all $x \in R_1$, $\mu_1(x) = \mu_1(0)$ and for all $x \in R_2$, $\mu_2(x) = \mu_2(0)]$.*

**THEOREM 1.5.10:** *Let $\mu$ be any fuzzy subset of a bifield $F = F_1 \cup F_2$. Then $\mu = \mu_1 \cup \mu_2$ is a fuzzy bi-ideal of $F$ if and only if $\mu(x) = \mu(y) \leq \mu(0)$ for all $x, y \in F \setminus \{0\}$. (i.e. $\mu_1(x) = \mu_1(y) \leq \mu_1(0)$ for all $x, y \in F_1 \setminus \{0\}$ and $\mu_2(x) = \mu_2(y) \leq \mu_2(0)$ for all $x, y \in F_2 \setminus \{0\})$.*

*Proof:* Follows as in case of ideals in fields.

**THEOREM 1.5.11:** *Let $R = R_1 \cup R_2$. be any biring. Then $R$ is bifield if and only if $\mu(x) = \mu(y) < \mu(0)$ where $\mu$ is any non-constant fuzzy bi-ideal of $R$ and $x, y \in R \setminus \{0\}$.*

*Proof*: Follows as in case of rings.

**THEOREM 1.5.12**: *If $\mu = \mu_1 \cup \mu_2$ is any fuzzy sub-biring (fuzzy bi-ideal) of a biring $R = R_1 \cup R_2$ and if $\mu(x) < \mu(y)$ for some $x, y \in R$ then $\mu(x - y) = \mu(x) = \mu(y - x)$.*

*Proof:* Left as an exercise for the reader to prove.

**DEFINITION 1.5.17:** *Let $\mu = \mu_1 \cup \mu_2$ be any fuzzy subset of a biring $R = R_1 \cup R_2$. The smallest fuzzy sub-biring (fuzzy bi-ideal) of $R$ containing $\mu$ is called the fuzzy sub-biring (fuzzy bi-ideal) generated by $\mu$ in $R$ and is denoted by $\langle \mu \rangle = \langle \mu_1 \rangle \cup \langle \mu_2 \rangle$.*

**THEOREM 1.5.13:** *Let $\mu = \mu_1 \cup \mu_2$ be any fuzzy subset of a biring. Then the fuzzy subset $\mu^*$ of $R = R_1 \cup R_2$ defined by $\mu^*(x) = \sup \{k / x \in \langle \mu_k \rangle \}$ ($\mu^*(x) = \mu_1^*(x) \cup \mu_2^*(x)$) is a fuzzy sub-biring (fuzzy bi-ideal) generated by $\mu$ in $R$ according as $\langle \mu_k \rangle$ is the sub-biring (bi-ideal) generated by $\mu_k$ in $R$. (In other wards $\mu^*(x) = t$ when ever $x \in \langle \mu_k \rangle$ and $x \notin \langle \mu_s \rangle$ for all $s > t$).*

We just define for any fuzzy bi-ideal $\mu = \mu_1 \cup \mu_2$ of a biring, $R = R_1 \cup R_2$ with $t = \mu(0)$. Then the fuzzy subset $\mu^*$ of $R/\mu_t$ is defined by $\mu^*(x + \mu_t) = \mu(x)$ for all $x \in R$ is a fuzzy bi-ideal of $R/\mu_t$.

$$( \text{Now } R/\mu_t = \frac{R_1 \cup R_2}{(\mu_t)_1 \cup (\mu_t)_2} = \frac{R_1}{(\mu_t)_1} \cup \frac{R_2}{(\mu_t)_2}$$

$$\text{so that } \mu^*(x + \mu_t) = \mu_1^* \cup \mu_2^* = \mu(x)$$

$$= \mu_1(x) \cup \mu_2(x) \text{ where}$$

$$= \mu_1^*(x + (\mu_t)_1) = \mu_1(x) \text{ for } x \in R_1 \text{ and}$$

$$= \mu_2^*(x + (\mu_t)_2) = \mu_2(x) \text{ for } x \in R_2).$$

**DEFINITION 1.5.18:** *Let $\mu = \mu_1 \cup \mu_2$ be any fuzzy bi-ideal of the biring $R = R_1 \cup R_2$. The fuzzy subset $\mu_x^*(=(\mu_x^*)_1 \cup (\mu_x^*)_2)$ of R where $x \in R_1 \cup R_2$ is defined by $\mu_x^*(r) = \mu(r - x)$ for all $r \in R$ i.e. termed as the fuzzy bicoset determined by x and $\mu$.*

*(Here*

$$\mu_x^*(r) = (\mu_x^*)_1(r) \cup (\mu_x^*)_2(r) = \mu_1(r-x) \cup \mu_2(r-x)$$

*where $(\mu_x^*)_1(r) = \mu_1(r - x)$ where $x, r \in R_1$ and $(\mu_x^*)_2(r) = \mu_2(r - x)$ where $x, r \in R_2$ ).*

Several results can be got for fuzzy birings as in fuzzy case fuzzy rings by appropriate and suitable modifications.

**Note:** If $\mu$ is a constant on $R = R_1 \cup R_2$ i.e. $\mu = \mu_1 \cup \mu_2$ where $\mu_1$ is a constant on $R_1$ and $\mu_2$ is a constant on $R_2$ then $R_{\mu_1} = (\mu_1)_0^*$ and $R_{\mu_2} = (\mu_2)_0^*$

with

$$R_\mu = (\mu_0^*), \text{ here}$$

$$\mu_0^* = (\mu_1)_0^* \cup (\mu_2)_0^*$$

Using these definitions the following theorem can be got as a matter of routine.

**THEOREM 1.5.14**: *Let $\mu$ be any fuzzy bi-ideal of a biring $R = R_1 \cup R_2$ then the following holds $\mu(x) = \mu(0)$ if and only $\mu_x^* = \mu_o^*$ where $x \in R$.*

**THEOREM 1.5.15:** *For any fuzzy bi-ideal $\mu$ of a biring R, the following holds good. $R|\mu_t \cong R_\mu$ where $\mu(0) = t$.*

Now we proceed on to define level bi-ideal fuzzy prime bi-ideal, fuzzy maximal bi-ideal and so on and give some interesting results on them.

**DEFINITION 1.5.19:** *Let $\mu = \mu_1 \cup \mu_2$ be any fuzzy bi-ideal of a biring $R = R_1 \cup R_2$ such that each level bi-ideal $\mu_t, t \in \text{Im } \mu$ is prime. If $\mu(x) < \mu(y)$ for some $x, y \in R$ then $\mu(xy) = \mu(y)$ (Here $\mu(x) < \mu(y) \Rightarrow \mu_1(x) < \mu_1(y)$ if $x, y \in R_1$ and $\mu_2(x) < \mu_2(y)$ if $x, y \in R_2$).*

**THEOREM 1.5.16:** *If $\mu = \mu_1 \cup \mu_2$ is any fuzzy prime bi-ideal of a biring $R = R_1 \cup R_2$ then $\mu(xy) = \max(\mu(x), \mu(y))$ for all $x, y \in R$.*

*Proof:* As in case of prime ideal we get the proof.

Now we give a nice characterization theorem, the proof of which is to be supplied by the reader.

**THEOREM 1.5.17:** *If $\mu$ and $\theta$ are any two fuzzy prime bi-ideals of a biring $R$, then $\mu \cap \theta$ is a fuzzy prime bi-ideal of $R$ if and only if $\mu \subseteq \theta$ or $\theta \subseteq \mu$.*

(Hint $\mu = \mu_1 \cup \mu_2$, $\theta = \theta_1 \cup \theta_2$, $\mu \cap \theta = (\mu_1 \cap \theta_1) \cup (\mu_2 \cap \theta_2)$, $R = R_1 \cup R_2$, so $\mu_1 \subseteq \theta_1$ $\mu_2 \subseteq \theta_2$ or $\mu_1 \subseteq \theta_1$ or $\mu_2 \subseteq \theta_2$; $\mu_1 : R_1 \to [0, 1]$, $\mu_2 : R_2 \to [0, 1]$, $\theta_1 : R_1 \to [0, 1]$ and $\theta_2 : R_2 \to [0, 1]$).

**THEOREM 1.5.18:** *Let $\mu = \mu_1 \cup \mu_2$ be any fuzzy bi-ideal of a biring $R$ such that $1 \in \operatorname{Im} \mu$ and let $\theta = \theta_1 \cup \theta_2$ be any fuzzy prime bi-ideal of the biring $R$. Then, $\mu \cap \theta$ is a fuzzy prime bi-ideal of the biring $\mu_t = \{x \in R \mid \mu(x) = 1\} = \{x_1 \in R_1 \mid \mu_1(x_1) = 1\} \cup \{x_2 \in R_2 \mid \mu_2(x_2) = 1\}$, where $\mu = \mu_1 \cup \mu_2$ such that $\mu_1 : R_1 \to [0, 1]$ and $\mu_2 : R_2 \to [0, 1]$ and $\operatorname{Im} \mu = \operatorname{Im} \mu_1 \cup \operatorname{Im} \mu_2$.*

*Proof:* Follows as in case of fuzzy prime ideals.

**THEOREM 1.5.19:** *If $\mu = \mu_1 \cup \mu_2$ is a fuzzy prime bi-ideal of a biring $R = R_1 \cup R_2$ then $R_\mu$ is an integral bidomain.*

*Proof:* Follows from the fact that each of $\mu_1 : R_1 \to [0, 1]$ and $\mu_2 : R_2 \to [0, 1]$ is such that $R_{\mu_1}$ and $R_{\mu_2}$ are integral domains so $R_\mu = R_{\mu_1} \cup R_{\mu_2}$ is an integral bidomain.

**THEOREM 1.5.20:** *Let $\mu$ be a fuzzy bi-ideal of biring $R = R_1 \cup R_2$ such that $\operatorname{Im} \mu = \{1, \sigma\}$ where $\sigma \in [0, 1)$. If $R_\mu$ is an integral bidomain then $\mu$ is a fuzzy prime bi-ideal.*

*Proof:* Follows as in case of rings.

**DEFINITION 1.5.20:** *A fuzzy bi-ideal $\mu$ of a biring $R = R_1 \cup R_2$ is called fuzzy maximal bi-ideal if $\operatorname{Im} \mu = \{1, \sigma\}$ where $\sigma \in [0, 1)$ and the level bi-ideal $\mu_t = \{x \in R \mid \mu(x) = 1\}$ is maximal.*

**THEOREM 1.5.21:** *Let $\mu = \mu_1 \cup \mu_2$ be any fuzzy bi-ideal of a biring $R = R_1 \cup R_2$. Then $\mu$ is a fuzzy maximal bi-ideal if and only if $R_\mu$ is a bifield.*

*Proof:* Straightforward by definitions and using fuzzy ideal analogous of a ring.

**THEOREM 1.5.22:** *Let $R$ be a biring $\alpha \in [0, 1)$ and $X = \{\mu \mid \mu$ is a fuzzy bi-ideal of the biring $R$ with $\operatorname{Im} \mu = \{1, \sigma\}\}$. Then $\mu \in X$ is fuzzy maximal if and only if for each $\beta \in X$ either $\beta \subseteq \mu$ or else $(\mu + \beta)(y) = 1$ for all $y \in R$.*

*Proof:* Follows directly from the definitions and results on rings.

**THEOREM 1.5.23:** *Let $\mu = \mu_1 \cup \mu_2$ and $\sigma = \sigma_1 \cup \sigma_2$ any two distinct fuzzy maximal bi-ideals of a biring R such that $\mu(0) = \sigma(0)$ and Im $\mu$ = Im $\sigma$ then $\mu \cap \sigma = \mu\sigma$.*

*Proof:* Left as an exercise for the reader.

Now we proceed on to define fuzzy semiprime bi-ideal of a biring.

**DEFINITION 1.5.21:** *A fuzzy bi-ideal $\mu = \mu_1 \cup \mu_2$ of a biring $R = R_1 \cup R_2$ is called fuzzy semiprime if for any fuzzy bi-ideal $\theta$ of R the condition $\theta^n \subseteq \mu$ implies $\theta \subseteq \mu$ where $n \in Z_+$.*

**THEOREM 1.5.24:** *Intersection of fuzzy semiprime bi-ideals of a biring R is always a fuzzy semiprime bi-ideal of the biring R.*

*Proof:* Left for the reader as an exercise.

The following theorem is left as an exercise for the reader.

**THEOREM 1.5.25:** *If $\mu$ is any fuzzy semiprime bi-ideal of a bring R then $R_\mu$ is free from non-zero nilpotent elements.*

*Proof:* Follows by results as in case of rings.

**THEOREM 1.5.26:** *Let $\mu$ be any fuzzy bi-ideal of a biring R with Im $\mu = \{t, j\}$ such that $t > j$. If the biring $R_\mu$ has no non-zero nilpotent elements then the fuzzy bi-ideal $\mu$ is fuzzy semiprime.*

*Proof:* Follows by definitions and proofs as in case of rings.

**DEFINITION 1.5.22:** *Let $R = R_1 \cup R_2$ be a biring R is said to be regular if and only if both $R_1$ and $R_2$ are regular rings.*

**THEOREM 1.5.27:** *Let $R = R_1 \cup R_2$ be a biring; R is regular if and only if every fuzzy bi-ideal of R is idempotent i.e. both $R_1$ and $R_2$ are idempotent.*

*Proof:* Follows directly by definitions and results on rings.

**THEOREM 1.5.28:** *A biring $R = R_1 \cup R_2$ is regular if and only if every fuzzy bi-ideal of R is fuzzy semiprime.*

*Proof:* Direct by definitions.

Next we proceed on to define fuzzy primary bi-ideals and semiprimary bi-ideals of birings.

**DEFINITION 1.5.23:** *A fuzzy bi-ideal $\mu$ of a biring R is called fuzzy primary if for any two fuzzy bi-ideals $\sigma$ and $\theta$ of the biring $R = R_1 \cup R_2$ the condition $\sigma\theta \subseteq \sqrt{\mu}$ and* $\sigma \not\subseteq \mu$ together imply $\theta \subseteq \sqrt{\mu}$. (Here for any fuzzy bi-ideal $\mu = \mu_1 \cup \mu_2$ of a biring R. The fuzzy nil radical of $\mu$ symbolized by $\sqrt{\mu} = \sqrt{\mu_1} \cup \sqrt{\mu_2}$ is defined by $(\sqrt{\mu})(x) = t$ whenever $x \in \sqrt{\mu_t}$, $x \notin \sqrt{\mu_s}$ for all $s > t$.

Here we assume $\sqrt{\mu_t} = \sqrt{(\mu_1)_t} \cup \sqrt{(\mu_2)_t}$ ).

The following theorem is straightforward and hence left for the reader to prove.

**THEOREM 1.5.29:** *Every fuzzy prime bi-ideal of a biring $R = R_1 \cup R_2$ is a fuzzy primary bi-ideal of R.*

Thus we are guaranteed of the existence of fuzzy primary bi-ideal in a biring.

**THEOREM 1.5.30:** *If $\mu = \mu_1 \cup \mu_2$ is any fuzzy primary bi-ideal of a biring $R = R_1 \cup R_2$ then $\mu_t = (\mu_1)_t \cup (\mu_2)_t$, $t \in \text{Im } \mu$ is a primary bi-ideal of R.*

*Proof*: Left for the reader as an exercise.

**THEOREM 1.5.31:** *If $A = A_1 \cup A_2$ is a primary bi-ideal of the biring $R = R_1 \cup R_2$, $A \neq R$ then the fuzzy subset $\mu = \mu_1 \cup \mu_2$ of R defined by*

$$\mu_1(x) = \begin{cases} 1 & \text{if } x \in A_1 \\ \alpha & \text{if } x \in R_1 \setminus A_1 \text{ where } \alpha \in [0,1) \end{cases}$$

*and*

$$\mu_2(x) = \begin{cases} 1 & \text{if } x \in A_2 \\ \alpha & \text{if } x \in R_2 \setminus A_2 \text{ where } \alpha \in [0,1) \end{cases}$$

$\mu = \mu_1 \cup \mu_2$ *is a fuzzy primary bi-ideal of R.*

*Proof:* Follows as in case of rings with simple modifications.

Now following the results in the above theorems we get the following characterization theorem; the proof of which is left for the reader to prove.

**THEOREM 1.5.32:** *A necessary and sufficient condition for a bi-ideal $A = A_1 \cup A_2$ of a biring $R = R_1 \cup R_2$ to be fuzzy primary is that $\chi_A$ is a fuzzy primary bi-ideal of R.*

**THEOREM 1.5.33:** *If $\mu = \mu_1 \cup \mu_2$ is any fuzzy primary bi-ideal of a biring R then $\sqrt{\mu}$ is fuzzy prime.*

*Proof*: Follows as in case of rings with appropriate modification. The fact that a biring $R_1 \cup R_2 = R$ is called primary if the zero ideal of both $R_1$ and $R_2$ are primary.

**THEOREM 1.5.34:** *If $\mu$ is any fuzzy primary bi-ideal of a biring $R = R_1 \cup R_2$, then the biring $R_\mu = (R_1)\mu \cup (R_2)\mu$ is primary.*

*Proof*: Straightforward hence left for the reader to prove.

**THEOREM 1.5.35:** *Let $\mu = \mu_1 \cup \mu_2$ be any fuzzy bi-ideal of a biring $R = R_1 \cup R_2$ such that Im $\mu = \{1, \alpha\}$ with $\alpha < 1$. If every zero divisor in $R\mu$ is nilpotent then $\mu$ is fuzzy primary.*

*Proof:* Follows as in case of rings.

We now present a sufficient condition for a fuzzy bi-ideal to be fuzzy primary.

**THEOREM 1.5.36:** *Let $\mu$ be any fuzzy bi-ideal of a biring $R = R_1 \cup R_2$ such that Im $\mu = \{1, \alpha\}$ with $\alpha < 1$. If $\sqrt{\mu}$ is a fuzzy maximal then $\mu$ is fuzzy primary.*

*Proof:* Follows as in case of rings by defining $\mu_t = (\mu_1)_t \cup (\mu_2)_t$ depending on $\mu = \mu_1 \cup \mu_2$ as follows:

$$\mu_1(x) = \begin{cases} 1 & \text{if } x \in (\mu_1)_t \\ \alpha & \text{if } x \in R_1 \setminus (\mu_1)_t \end{cases}$$

$$\mu_2(x) = \begin{cases} 1 & \text{if } x \in (\mu_2)_t \\ \alpha & \text{if } x \in R_2 \setminus (\mu_2)_t \end{cases}$$

so that

$$\left(\sqrt{\mu_1}\right)(x) = \begin{cases} 1 & \text{if } x \in \sqrt{(\mu_1)_t} \\ \alpha & \text{if } x \in R_1 \setminus \sqrt{(\mu_1)_t} \end{cases}$$

$$\left(\sqrt{\mu_2}\right)(x) = \begin{cases} 1 & \text{if } x \in \sqrt{(\mu_2)_t} \\ \alpha & \text{if } x \in R_2 \setminus \sqrt{(\mu_2)_t} \end{cases}$$

$\sqrt{\mu} = \sqrt{\mu_1} \cup \sqrt{\mu_2}$. Using the fact $\sqrt{\mu}$ is fuzzy maximal bi-ideal the rest of the result follows as in case of rings.

**DEFINITION 1.5.24:** *A fuzzy bi-ideal $\mu = \mu_1 \cup \mu_2$ of the biring $R = R_1 \cup R_2$ is called fuzzy semiprimary if $\sqrt{\mu} = \sqrt{\mu_1} \cup \sqrt{\mu_2}$ is a fuzzy prime bi-ideal of the biring R.*

Now we proceed on to give some interesting results on fuzzy semiprimary bi-ideals of a biring.

**THEOREM 1.5.37:** *If $A = A_1 \cup A_2$ is any semiprimary bi-ideal of a biring $R = R_1 \cup R_2$ ($R_1 \neq A_1$, $A_2 \neq R_2$) then the fuzzy subset $\mu = \mu_1 \cup \mu_2$ of R defined by*

$$\mu_1(x) = \begin{cases} 1 & \text{if } x \in A_1 \\ \alpha & \text{if } x \in R_1 \setminus A_1 \text{ where } \alpha < 1 \end{cases}$$

*and*

$$\mu_2(x) = \begin{cases} 1 & \text{if } x \in A_2 \\ \alpha & \text{if } x \in R_2 \setminus A_2 \text{ where } \alpha < 1 \end{cases}$$

*is a fuzzy semiprimary bi-ideal of R.*

*Proof*: Follows as in case of rings with suitable changes.

**THEOREM 1.5.38:** *If $\mu$ is any fuzzy semiprimary bi-ideal of a biring $R = R_1 \cup R_2$ then $\mu_t = (\mu_1)_t \cup (\mu_2)_t$ where $t \in \text{Im } \mu$ is a semiprimary bi-ideal of R.*

*Proof:* Left as an exercise for the reader to prove.

**THEOREM 1.5.39:** *A bi-ideal $A = A_1 \cup A_2$ of a biring $R = R_1 \cup R_2$ is semiprimary if and only if $\chi_A$ is a fuzzy semiprimary bi-ideal of R.*

*Proof*: Follows from definitions and also using the method of rings we can get the proof.

**THEOREM 1.5.40**: *If $\mu = \mu_1 \cup \mu_2$ is any fuzzy semiprimary bi-ideal of a ring R then the biring $R_\mu$ is semiprimary.*

*Proof*: The result follows by proving the zero bi-ideal of $R_\mu = R_{\mu_1} \cup R_{\mu_2}$ is semiprimary.

**THEOREM 1.5.41:** *Let $\mu = \mu_1 \cup \mu_2$ by any fuzzy bi-ideal of a biring $R = R_1 \cup R_2$ such that $\text{Im } \mu = \{1, \alpha\}$ where $\alpha < 1$. If the ring $R_\mu$ is semiprimary then $\mu$ is a fuzzy semiprimary bi-ideal of R.*

*Proof*: Follows by the very definitions.

**THEOREM 1.5.42:** *If $\mu$ is any non-constant fuzzy semiprimary bi-ideal of a regular biring then $\text{Im}\mu = \{1, \alpha\}$, $\alpha \in [0, 1)$.*

*Proof*: Follows directly by the definitions and by routine methods, hence left for the reader as an exercise.

Now we proceed on to define the notions of fuzzy irreducible bi-ideals of a biring R.

**DEFINITION 1.5.25**: *A fuzzy bi-ideal $\mu = \mu_1 \cup \mu_2$ of a biring $R = R_1 \cup R_2$ is called fuzzy irreducible if it is not an intersection of two fuzzy bi-ideals of R properly containing $\mu$; otherwise $\mu$ is termed as fuzzy reducible.*

**THEOREM 1.5.43**: *If $\mu = \mu_1 \cup \mu_2$ is any fuzzy prime bi-ideal of a biring R then $\mu$ is fuzzy irreducible.*

*Proof*: Use the fact $\mu_1: R_1 \to [0, 1]$ and $\mu_2: R_2 \to [0, 1]$ are fuzzy prime ideals of $R_1$ and $R_2$ respectively and they are fuzzy irreducible so $\mu = \mu_1 \cup \mu_2$ is fuzzy irreducible.

**THEOREM 1.5.44:** *Let $\mu$ be any non-constant fuzzy irreducible bi-ideal of a biring R. Then there exists $\alpha \in [0, 1)$ such that*

  i.   *Im $\mu = \{1, \alpha\}$ and*
  ii.  *the level bi-ideal $\{x \in R \,/\, \mu(x) = 1\}$ is irreducible.*

*Proof*: The proof is a matter of routine as in case of rings. Hence it is left for the reader as an exercise.

**THEOREM 1.5.45:** *If $A = A_1 \cup A_2$ is any irreducible bi-ideal of a biring R, $A_1 \neq R_1$, $A_2 \neq R_2$ where $R = R_1 \cup R_2$, then the fuzzy ideal $\mu = \mu_1 \cup \mu_2$ of R defined by*

$$\mu_1(x) = \begin{cases} 1 & \text{if } x \in A_1 \\ \alpha & \text{if } x \in R_1 \setminus A_1 \end{cases}$$

*and*

$$\mu_2(x) = \begin{cases} 1 & \text{if } x \in A_2 \\ \alpha & \text{if } x \in R_2 \setminus A_2 \end{cases}$$

*(where $\mu = \mu_1 \cup \mu_2$) where $\alpha \in [0, 1)$ is a fuzzy irreducible bi-ideal of the biring R.*

*Proof*: Straightforward as in case of rings.

In view of the earlier theorems we have the following nice characterization theorem.

**THEOREM 1.5.46:** *A necessary and sufficient condition for a bi-ideal $A = A_1 \cup A_2$ of a biring $R = R_1 \cup R_2$ to be irreducible is that $\chi_A$ is a fuzzy irreducible bi-ideal of R.*

Now we give condition for a fuzzy irreducible bi-ideal to be fuzzy prime.

**THEOREM 1.5.47:** *If $\mu = \mu_1 \cup \mu_2$ is any fuzzy bi-ideal of a biring R which is both fuzzy semiprime and fuzzy irreducible then $\mu$ is fuzzy prime.*

*Proof*: Straightforward, hence left for the reader to prove.

**THEOREM 1.5.48:** *In a regular biring every fuzzy irreducible bi-ideal is fuzzy prime.*

*Proof*: It is left as an exercise for the reader to prove using earlier results and the definition of fuzzy irreducible bi-ideal.

We say a biring $R = R_1 \cup R_2$ to be Noetherian if both $R_1$ and $R_2$ are Noetherian rings. Every bi-ideal of a Noetherian biring $R = R_1 \cup R_2$ can be represented as a finite intersection of fuzzy primary bi-ideal of R i.e. every ideal of $R_1$ and $R_2$ can be represented as a finite intersection of fuzzy primary ideals of R. Now we give a condition for a biring $R = R_1 \cup R_2$ to be Noetherian.

**THEOREM 1.5.49:** *If the cardinality of the image set of every bi-ideal of a biring $R = R_1 \cup R_2$ is finite then the biring R is Noetherian.*

*Proof:* Using the fact that the biring $R = R_1 \cup R_2$ where $R_1$ and $R_2$ are rings, it is enough if we can prove that each of $R_1$ and $R_2$ satisfies the Noetherian ring property.

**THEOREM 1.5.50:** *A biring R is artinian if and only if every fuzzy ideal is finite valued.*

*Proof:* As in case of rings the result follows.

**THEOREM 1.5.51:** *If $\mu$ is any fuzzy irreducible bi-ideal of a Noetherian biring R then $\mu$ is fuzzy primary.*

*Proof:* Easily follows as in case of rings.

**THEOREM 1.5.52:** *Let $\mu$ be any fuzzy bi-ideal of a Noetherian biring $R = R_1 \cup R_2$ such that $\operatorname{Im} \mu = \{1, \alpha\}$, $\alpha < 1$. Then $\mu = \mu_1 \cup \mu_2$ can be written as a finite intersection of fuzzy irreducible bi-ideals of R.*

*Proof:* Let $\mu_t = (\mu_1)_t \cup (\mu_2)_t = \{x \in R_1 \mid \mu_1(x) = 1\} \cup \{x \in R_2 \mid \mu_2(x) = 1\}$, then by hypothesis we have

$$\mu_1(x) = \begin{cases} 1 & \text{if } x \in (\mu_1)_t \\ \alpha & \text{if } x \in R_1 \setminus (\mu_1)_t \end{cases}$$

$$\mu_2(x) = \begin{cases} 1 & \text{if } x \in (\mu_2)_t \\ \alpha & \text{if } x \in R_2 \setminus (\mu_2)_t \end{cases}$$

We know that there exists irreducible bi-ideals $A_1, A_2, \ldots, A_n$ ($n < \infty$) of $R = R_1 \cup R_2$ satisfying $\mu_t = A_1 \cap A_2 \cap \ldots \cap A_n$ where $A_i = A_{i1} \cup A_{i2}$ for $i = 1, 2, \ldots, n$ so that $\mu_t = (\mu_1)_t \cup (\mu_2)_t$ where $(\mu_1)_t = A_{11} \cap A_{21} \cap \ldots \cap A_{n1}$ and $(\mu_2)_t = A_{12} \cap A_{22} \cap \ldots \cap A_{n2}$. As in case of rings we can define fuzzy bi-ideals $\mu_i = (\mu_1)_i \cup (\mu_2)_i$ by

$$(\mu_1)_i(x) = \begin{cases} 1 & \text{if } x \in A_{i1} \\ \alpha & \text{if } x \in R_1 \setminus A_{i1} \end{cases}$$

$$(\mu_2)_i(x) = \begin{cases} 1 & \text{if } x \in A_{i2} \\ \alpha & \text{if } x \in R_2 \setminus A_{i2} \end{cases}$$

$1 \leq i \leq n$. It easily follows each $\mu_i$ is fuzzy irreducible. Thus $\mu = \mu_1 \cap \ldots \cap \mu_n$ that is $\mu = (\mu_{11} \cap \mu_{21} \cap \ldots \cap \mu_{n1}) \cup (\mu_{12} \cap \mu_{22} \cap \ldots \cap \mu_{n2})$.

Now the following theorem can be easily established by any innovative reader.

**THEOREM 1.5.53:** *If $\mu = \mu_1 \cup \mu_2$ is a fuzzy bi-ideal of a Noetherian bi-ring $R = R_1 \cup R_2$ such that Im $\mu = \{1, \alpha\}$, $\alpha < 1$, then $\mu$ can be expressed as a finite intersection of fuzzy primary bi-ideals of the biring R.*

Several more results on fuzzy subrings and fuzzy sub-birings can be obtained; but as the main aim of the book is only study on the Smarandache fuzzy algebra, we have restrained ourselves from giving several results. Here also only proofs are hinted and it is left for the readers to prove.

## 1.6 Fuzzy Fields and their properties

The concept of fuzzy subfield is recalled in this section. For more about fuzzy subfield please refer [75, 79].

**DEFINITION 1.6.1:** *Let F be a field, A fuzzy subfield of F is a function A from F into the closed interval [0, 1] such that for all $x, y \in F$*

$$(A(x - y) \geq \min \{A(x), A(y)\} \text{ and}$$
$$A(xy^{-1}) \geq \min \{A(x), A(y)\}; y \neq 0.$$

*Let A be a fuzzy subset of F and let $A_{\neq} = \{x \in F / A(x) \geq A(1)\}$ where 1 denotes the multiplicative identity of F. Let K be a subfield of F and let S(F/K) denote the set of all fuzzy subfields, A of F such that $K \subseteq A_{\neq}$. Here we just recall certain properties of field extensions F/K in terms of fuzzy subfields and conversely.*

Let A be a fuzzy subset of the field F. For $0 \leq t \leq 1$, let $A_t = \{x \in F \mid A(x) \geq t\}$. Then $A_{\neq} = A_t$ when $t = A(1)$.

The following results are given without proof for the reader. However the interested reader can refer [75].

**THEOREM 1.6.1:**

  i. If A is a fuzzy subset of F and $s, t \in$ Im (A), the image of A, then $s \leq t$ if and only if $A_s \supseteq A_t$ and $s = t$ if and only if $A_s = A_t$.
  ii. If A is a fuzzy subfield of F, then for all $x \in F$, $x \neq 0$, $A(0) \geq A(1) \geq A(x) = A(-x) = A(x^{-1})$.

**THEOREM 1.6.2:** *Let A be a fuzzy subset of F. If $A_t$ is a subfield of F for all $t \in$ Im (A), then A is a fuzzy subfield of F. Conversely, if A is a fuzzy subfield of F, then for all t such that $0 \leq t \leq A(1)$, $A_t$ is a subfield of F.*

**THEOREM 1.6.3:** *Let S be a subset of F such that $|S|$ (Cardinality of S) $\geq 2$. Then S is a subfield of F if and only if $\chi_S$, the characteristic function of S, is a fuzzy subfield of F.*

Recall if K be a subfield of F i.e. F is an extension field of K then the field extension is denoted by F/K. S(F/K) denotes the set of all fuzzy subfields A of F such that $A_{\neq} \supseteq K$ and $A_{\neq}$ is a subfield of F.

**THEOREM 1.6.4:** *Let $F_1 \subset F_2 \subset \ldots \subset F_i \ldots$ be a strictly ascending chain of subfields of F such that $\cup F_i = F$. Define the fuzzy subset A of F by $A(x) = t_i$, if $x \in F_i \setminus F_{i-1}$ where $t_i > t_{i+1}$ for $i = 1, 2, \ldots$ and $F_o = \phi$. Then A is a fuzzy subfield of F.*

*Proof:* Let $x, y \in F$. Then $x - y \in F_i \setminus F_{i-1}$ for some i. Hence either $x \notin F_{i-1}$ or $y \notin F_{i-1}$. Thus $A(x - y) = t_i \geq \min \{A(x), A(y)\}$, similarly $A(xy^{-1}) \geq \min \{A(x), A(y)\}$ for $y \neq 0$.

**THEOREM 1.6.5:** *Let $F = F_0 \supset F_1 \supset \ldots \supset F_i \supset \ldots$ be a strictly descending chain of subfields of F. Define the fuzzy subset A of F by $A(x) = t_{i-1}$ if $x \in F_{i-1} \setminus F_i$ where $t_{i-1} < t_i < 1$ for $i = 1, 2, \ldots$ and $A(x) = 1$ if $x \in \cap F_i$. Then A is a fuzzy subfield of F.*

*Proof:* Direct, hence left for the reader to prove.

**THEOREM 1.6.6:** *Let F / K be a field extension and let B be a fuzzy subfield of K. Let $r = \inf \{B(x) \mid x \in K\}$. Define the fuzzy subset A of F by $A(x) = B(x)$ for all $x \in K$ and $A(x) = m$ for all $x \in F \setminus K$ where $0 \leq m \leq r$. Then A is a fuzzy subfield of F.*

*Proof:* Left as an exercise for the reader to prove.

The following theorem can be easily proved by any reader. The reader is also advised to refer [75] for more information.

**THEOREM 1.6.7:** *If F is a finite field, then every fuzzy subfield of F is finite valued.*

**THEOREM 1.6.8:** *Let F/K be a field extension. Then $[F : K] < \infty$ if and only if every $A \in S, (F/K)$ is finite valued.*

**THEOREM 1.6.9:**

   i. Suppose that F has characteristic $p > 0$. Then F is finite if and only if every fuzzy subfield A of F is finite-valued.

   ii. Suppose that F has characteristic 0.

Then $[F : Q] < \infty$ if and only if every fuzzy subfield A of F is finite valued.

**THEOREM 1.6.10:** *Suppose that F/K is finitely generated. Then F/K is algebraic if and only if every $A \in S(F/K)$ is finite valued.*

*Proof:* Left as an exercise for the reader to prove.

**THEOREM 1.6.11:** *F/K has no proper intermediate fields if and only if every $A \in S(F/K)$ is three valued or less.*

*Proof:* Please refer [75].

The following theorem which gives equivalent conditions is left as an exercise for the reader to prove.

**THEOREM 1.6.12:** *The following conditions are equivalent.*

  i. *The intermediate fields of F/K are chained.*
 ii. *There exists $C \in S(F/K)$ such that for all $A \in S(F/K)$. $L_A \subseteq L_C$.*
iii. *For all $A, B \in S(F/K)$ and for all $A_t \in L_A$ and $B_s \in L_B$ either $A_t \subseteq B_s$ or $B_s \subseteq A_t$.*

We give a necessary and sufficient condition for F/K to be simple.

**THEOREM 1.6.13:** *F/K is simple if and only if there exists $c \in F$ such that for all $A \in S(F/K)$ and for all $x \in F$, $A(c) \leq A(x)$.*

*Proof:* Direct by the regular computations.

**THEOREM 1.6.14:** *Suppose that $[F : K] < \infty$. Then the following conditions are equivalent.*

  i. *F/K has a finite number of intermediate fields.*
 ii. *There exists $C_1, C_2, \ldots, C_n \in S(F/K)$ such that for all $A \in S(F/K)$. $L_A \subset L_{C_1} \cup \cdots \cup L_{C_n}$.*
iii. *There exists $c \in F$ such that for all $A \in S(F/K)$ and for all $x \in F$, $A(c) \leq A(x)$.*

*Proof:* By the earlier result a
nd from the theorem in [49]. The result can be easily obtained.

**THEOREM 1.6.15:** *Let F/K be a field extension where K has characteristic $p > 0$ and let $c \in F$. Then*

  i. *K(c) / K is separable algebraic if and only if for all $A \in S(F/K)$, $A(c) = A(c^p)$.*
 ii. *K(c) / K is pinely inseparable if and only if there exists a non-negative integer e such that for all $A \in S(F/K)$, $A(c^{p^e}) = A(1)$.*
iii. *K(c) / K is inseparable if and only if there exists $A \in S(F/K)$ such that $A(c) < A(c^p)$ and there exists a positive integer e such that for all $A \in S(F/K)$ $A(c^{p^e}) = A(c^{p^{e-1}})$.*

*Proof:* Left for the reader to find the proof; as the proof does not involve any deeper knowledge of algebra or field theory.

Now we proceed on to define the concept of fuzzy bifields. The concept of the algebraic structure bifields is itself very new [135]. So the notion of fuzzy bifield is defined for the first time in this book. We just for the sake of completeness recall the definition of bifields.

**DEFINITION 1.6.2:** *A biring (R, +, •) where $R = R_1 \cup R_2$ is said to be a bifield if ($R_1$, +, •) and ($R_2$, +, •) are fields. If the characteristic of both $R_1$ and $R_2$ are finite then we say $R = R_1 \cup R_2$ is a bifield of finite characteristic.*

*If in $R = R_1 \cup R_2$ one of $R_1$ or $R_2$ is a field of characteristic 0 and one of $R_1$ or $R_2$ is of finite characteristic we do not associate any characteristic with it. If both $R_1$ or $R_2$ in $R = R_1 \cup R_2$ is zero characteristic then we say R is a field of characteristic zero.*

*Thus unlike in fields we see in case of bifields we can have characteristic prime or characteristic zero or no characteristic associated with it.*

*Example 1.6.1:* Let $R = Z_{11} \cup Q$ the field of rationals Q and $Z_{11}$, the prime field of characteristic 11. R is a bifield with no characteristic associated with it.

*Example 1.6.2:* Let $R = Q(\sqrt{2}) \cup Q$. Clearly R is not a bifield as $Q \subseteq Q(\sqrt{2})$.

*Example 1.6.3:* Let $R = Q(\sqrt{2}) \cup Q(\sqrt{5})$. Clearly R is a bifield of characteristic 0.

*Example 1.6.4:* Let $R = Z_{11} \cup Z_{19}$; R is a bifield of finite characteristic.

**DEFINITION 1.6.3:** *Let $F = F_1 \cup F_2$ be a bifield, we say a proper subset S of F to be a sub-bifield if $S = S_1 \cup S_2$ and $S_1$ is a subfield of $F_1$ and $S_2$ is a subfield of $F_2$. If the bifield has no proper sub-bifield then we call F a prime bifield.*

*Example 1.6.5:* Let $R = Q \cup Z_7$, clearly R is a prime bifield.

*Example 1.6.6:* Let $R = Q(\sqrt{7}, \sqrt{3}) \cup Q(\sqrt{2}, \sqrt{19})$. Clearly R is a non-prime bifield for the subset $S_1 = Q(\sqrt{7}) \cup Q(\sqrt{2})$ and $S_2 = Q(\sqrt{3}) \cup Q(\sqrt{19})$ are sub-bifields of R. Thus R is not a prime bifield.

**DEFINITION 1.6.4:** *Let (D, +, •) with $D = D_1 \cup D_2$ where ($D_1$, +, •) and ($D_2$, +, •) are integral domains then we call D a bidomain.*

**THEOREM 1.6.16:** *Let (R, +, •) be a bifield. $R = R_1 \cup R_2$. Now (R[x], +, •) where $R[x] = R_1[x] \cup R_2[x]$ is a bidomain.*

*Proof:* Straightforward hence left for the reader to prove.

Now having seen the definition of bifield we now define fuzzy bifield.

**DEFINITION 1.6.5:** *Let $(F = F_1 \cup F_2, +, \bullet)$ be a bifield. $\mu : F \to [0, 1]$ is said to be a fuzzy sub-bifield of the bifield F if there exists two fuzzy subsets $\mu_1$ (of $F_1$) and $\mu_2$ (of $F_2$) such that*

    i.    *$(\mu_1, +, \bullet)$ is a fuzzy sub-field of $(F_1, +, \bullet)$.*
    ii.   *$(\mu_2, +, \bullet)$ is a fuzzy sub-field of $(F_2, +, \bullet)$.*
    iii.  *$\mu = \mu_1 \cup \mu_2$.*

**DEFINITION 1.6.6:** *Let $(F = F_1 \cup F_2, +, \bullet)$ be a bifield and $(S = S_1 \cup S_2, +, \bullet)$ be a sub-bifield of F. Let B (F/S) denote the set of all fuzzy sub-bifields; $A = A_1 \cup A_2$ of F, such that $A_{1\neq} \supseteq S_1$ and $A_{2\neq} \supseteq S_2$. We say $A \in B$ (F/S) is finite bivalued if the image of $A_1$ is finite and the image of $A_2$ is finite where $A = A_1 \cup A_2$ is a fuzzy sub-bifield.*

**THEOREM 1.6.17:** *Let $(F = F_1 \cup F_2, +, \bullet)$ be a finite bifield; then every fuzzy sub-bifield of F is finite bivalued.*

*Proof:* Follows from the definition.

In view of this we have the following results.

**THEOREM 1.6.18:** *Suppose that $(F = F_1 \cup F_2, +, \bullet)$ has finite characteristic bifield (i.e. characteristic of the field $F_1$ is a prime $p_1 > 0$ and the characteristic of the field $F_2$ is a prime $p_2 > 0$). Then F is finite if and only if every fuzzy sub-bifield A of F is finite valued.*

*Proof:* Left for the reader to prove, following the steps of finite fuzzy subfield.

**THEOREM 1.6.19:**

    i.   *Let A be a fuzzy subset of the bifield F and $s, t \in$ Im (A) the image of A where $A = A_1 \cup A_2$, then $s \leq t$ if and only if $A_s \supseteq A_t$ and $s = t$ if and only if $A_s = A_t$. (that is if $(A_1)_s \supset (A_1)_t$ and $(A_2)_s \supseteq (A_2)_t$ then only we say $A_s \supseteq A_t$, similarly for $A_s = A_t$).*

    ii.  *If A is a fuzzy sub-bifield of F ($A = A_1 \cup A_2$ and $F = F_1 \cup F_2$) then for all $x \in F = F_1 \cup F_2$, $x \neq 0$ we have $A_i(x) \geq A_i(1) \geq A_i(x) = A_i(-x) = A_i(x^{-1})$ for $i = 1, 2$.*

*Proof:* Direct as in case of fuzzy fields.

**THEOREM 1.6.20:** *Let $A = A_1 \cup A_2$ be a fuzzy subset of $F = F_1 \cup F_2$. If $A_t = (A_1)_t \cup (A_2)_t$ is a sub-bifield of F for all $t \in$ Im$(A)_1$ then A is a fuzzy sub-bifield of the bifield F.*

*Conversely if $A = A_1 \cup A_2$ is a fuzzy sub-bifield of $F = F_1 \cup F_2$ then for all t such that $0 \leq t \leq A_i (1)$, $i = 1, 2$; $A_t$ is a sub-bifield of F. If A is a fuzzy sub-bifield of F, then $A_t$ is called a level sub-bifield of F where $0 < t < A(1)$.*

**THEOREM 1.6.21:** *Let $S = S_1 \cup S_2$ be a subset of the bifield $F = F_1 \cup F_2$ i.e. ($S_1$ is a subset of $F_1$ and $S_2$ is a subset of $F_2$) such that $|S_1| \geq 2$ and $|S_2| \geq 2$. Then S is a sub-bifield of F if and only if $\chi_S = \chi_{S_1} \cup \chi_{S_2}$ '$\cup$' is just only a default of notation, the characteristic bifunction of S is a fuzzy sub-bifield of F.*

*Proof:* Follows as in case of fuzzy fields.

**DEFINITION 1.6.7:** *Let $F = F_1 \cup F_2$ be a bifield. Suppose $S_1 \subset S_2 \subset \ldots \subset S_n \subset \ldots$ be a strictly ascending chain of sub-bifield of F such that $\cup S_i = F$ where $S_i = S_{i1} \cup S_{i2}$ for i = 1, 2, …(i.e. $S_1 = S_{11} \cup S_{12}$ is a sub-bifield) i.e. each $S_{11} \subset S_{21} \subset S_{31} \subset \ldots \subset S_{n1} \subset \ldots$ is a strictly ascending chain of subfields of $F_1$ and $S_{12} \subset S_{22} \subset S_{32} \subset \ldots \subset S_{n2} \subset \ldots$ is a strictly ascending chain of subfields of $F_2$. Define the fuzzy subset $A = A_1 \cup A_2$ of $F = F_1 \cup F_2$ by $A(x) = t_i$ if $x \in S_{ij} \setminus S_{i-1j}$; j = 1, 2; where $t_i > t_{i+1}$ for i = 1, 2, ... and $S_{01} = \phi$ and $S_{02} = \phi$. Then A is a fuzzy sub-bifield of F.*

The proof of the following theorem is left as an exercise for the reader.

**THEOREM 1.6.22:** *Let $F = F_1 \cup F_2$ be a bifield $F = S_0 \supset S_1 \supset \ldots$ be a strictly descending chain of sub-bifields of F. Define the fuzzy subset $A = A_1 \cup A_2$ of F by $A(x) = t_{i-1}$ if $x \in S_{i-1\,j} \setminus S_{ij}$ where $t_{i-1} < t_i < 1$ for i = 0, 1, 2, … and $A(x) = 1$ if $x \in \cap S_{ij}$. j = 1, 2 and $0 \leq i \leq n$. Then A is a fuzzy sub-bifield of F.*

**THEOREM 1.6.23:** *Let F/K be a bifield extension and let B be a fuzzy sub-bifield of K. Let $r = \inf \{B(x) \mid x \in K\}$. Define the fuzzy subset A of F by $A(x) = B(x)$ for all $x \in K$ and $A(x) = m$ for all $x \in F \setminus K$; where $0 \leq m \leq r$. Then A is a fuzzy sub-bifield of F.*

*Proof:* Follows as in case of subfield.

Thus all results regarding fuzzy sub-bifields can be defined and obtained in an analogous way as in case of fuzzy subfields. These concepts will be once again used in case of fuzzy bivector spaces.

## 1.7 Fuzzy Semirings and their generalizations

In this section we introduce the notion of fuzzy semirings and fuzzy sub-birings. The study of fuzzy k-ideals in semirings started in 1985 [31] followed by several authors. Even [19] has studied about it. But the notion of fuzzy bisemirings in literature is totally absent. This concept is defined only in this book. Just for the sake of completeness we start to give the definition of fuzzy semirings and proceed on.

An algebra (S, +, •) is said to be a semiring if (S, +) is a semigroup with identity 0 and (S, •) is a semigroup, satisfying the following conditions

   i. a • (b + c) = a • b + a • c
   ii. (b + c) • a = b • a + b • a

for all a, b, c ∈ S.

A semiring S may have the identity 1 defined by $1 \cdot a = a \cdot 1 = a$ and a zero 0 defined by $0 + a = a + 0 = a$ for all $a \in S$ and $0 \cdot a = a \cdot 0 = 0$. We say the semiring $(S, +, \cdot)$ is abelian if $a \cdot b = b \cdot a$ for all $a, b \in S$. We say the semiring S is strict if $a + b = 0$ forces $a = 0$ and $b = 0$. Let $(S, +, \cdot)$ be a semiring if $a \cdot b = 0$ in S implies $a = 0$ or $b = 0$ for $a, b \in S \setminus \{0\}$ then we say the semiring $(S, +, \cdot)$ is a semidivision ring.

If $(S, +, \cdot)$ is commutative semiring with no zero divisions which is strict then we say $(S, +, \cdot)$ is a semifield.

*Example 1.7.1:* Let $Z^+ \cup \{0\} = S$. $(S, +, \cdot)$ is a semifield ($Z^+$, set of positive integers).

*Example 1.7.2:* Let

$$M_{2 \times 2} = \left\{ \begin{pmatrix} a & b \\ c & d \end{pmatrix} \middle/ a, b, c, d \in Z^+ \cup \{0\} \right\}$$

is a non-commutative strict semiring which is not a semidivision ring.

*Example 1.7.3:* Let

$$M'_{2 \times 2} = \left\{ \begin{pmatrix} a & b \\ c & d \end{pmatrix} \middle/ a, b, c, d \in Z^+ \right\} \cup \left\{ \begin{pmatrix} 0 & 0 \\ 0 & 0 \end{pmatrix} \right\}$$

be a semiring. $M'_{2 \times 2}$ is a semidivision ring.

*Example 1.7.4:* Let $C_{11}$ be the chain lattice with 11 elements, $C_{11}$ is a finite semiring in fact a finite semifield.

We denote the order of the semiring S by $|S|$ or $o(S)$, if $|S| = \infty$ we say the semiring is of infinite order otherwise S is of finite order if $|S| < \infty$.

Now we proceed on to define left (right) ideal in a semiring. Finally we define left k-ideal of the semiring S.

**DEFINITION 1.7.1:** *Let R be a semiring. A non-empty subset I of R is said to be left (resp, right) ideal if $x, y \in I$ and $r \in I$ imply $x + y \in I$ and $r x \in I$ (resp. $x r \in I$). If I is both left and right ideal if R, we say I is a two sided ideal or simply ideal of R. A left ideal I of a semiring R is said to be a left k-ideal if $a \in I$ and $x \in R$ and if $a + x \in I$ or $x + a \in I$ then $x \in I$. Right k-ideal is defined dually and two-sided k-ideal or simply a k-ideal is both a left and a right k-ideal. A mapping $f : R \to S$ is said to be a homomorphism if $f(x + y) = f(x) + f(y)$ and $f(xy) = f(x) f(y)$ for all $x, y \in R$. We note that if $f : R \to S$ is an onto homomorphism and I is left (resp. right) ideal of R, then $f(I)$ is a left (resp. right) ideal of S.*

Now we proceed on to recall the definition of fuzzy ideal of a semiring.

**DEFINITION 1.7.2:** *A fuzzy subset $\mu$ of a semiring R is said to be a fuzzy left (resp. right) ideal of R if $\mu(x + y) \geq \min \{\mu(x), \mu(y)\}$ and $\mu(xy) \geq \mu(y)$ (resp. $\mu(xy) \geq \mu(x)$) for all $x, y \in R$. $\mu$ is a fuzzy ideal of R if it is both a fuzzy left and a fuzzy right ideal of R.*

**DEFINITION 1.7.3:** *A fuzzy ideal $\mu$ of a semiring R is said to be a fuzzy k-ideal of R if $\mu(x) \geq \min \{\max \{\mu(x + y), \mu(y + x)\}, \mu(y)\}\}$ for all $x, y \in R$.*

It is left for the reader to prove the following theorem:

**THEOREM 1.7.1:** *Every fuzzy ideal of a semiring is a fuzzy k-ideal.*

**THEOREM 1.7.2:** *Let I be a non-empty subset of a semiring R and $\lambda_I$ the characteristic function of I. Then I is a k-ideal of R if and only if $\lambda_I$ is a fuzzy k-ideal of R.*

*Proof:* The reader is advised to give the proof as an exercise.

**THEOREM 1.7.3:** *A fuzzy subset $\mu$ of R is a fuzzy left (resp. right) k-ideal of R if and only if for any $t \in [0, 1]$ such that $\mu_t \neq \phi$, $\mu_t$ is a left (resp. right) k-ideal of R, where $\mu_t = \{x \in R \mid \mu(x) \geq t\}$, which is called a level subset of $\mu$.*

*Proof:* For proof please refer [52].

Note that if $\mu$ is a fuzzy left (resp. right) k-ideal of R then the set $R_\mu = \{x \in R \mid \mu(x) \geq \mu(0)\}$ is a left (resp. right) k-ideal of R.

**THEOREM 1.7.4:** *Let I be any left (resp. right) k-ideal of R. Then there exists a fuzzy left (resp. right) k-ideal $\mu$ of R such that $\mu_t = I$ for some $t \in [0, 1]$.*

*Proof:* If we define a fuzzy subset of R by

$$\mu(x) = \begin{cases} t & \text{if } x \in I \\ 0 & \text{otherwise} \end{cases}$$

for some $t \in [0, 1]$, then if follows that $\mu_t = I$ for a given $s \in [0, 1]$. We have

$$\mu_s = \begin{cases} \mu_o (= R) & \text{if } s = 0 \\ \mu_t (= I) & \text{if } s \leq t \\ \phi & \text{if } t < s \leq 1 \end{cases}$$

Since I and R itself are left (resp. right) k-ideals of R, it follows that every non-empty level subset $\mu_s$ of $\mu$ is left (resp. right) k-ideal of R. Thus by the earlier theorem $\mu$ is a fuzzy left (resp. right) k-ideal of R proving the theorem.

Let $\mu$ and $\delta$ be fuzzy subsets of the semiring R. We denote that $\mu \subseteq \delta$ if and only if $\mu(x) \leq \delta(x)$ for all $x \in R$ and $\mu \subset \delta$ if and only if $\mu \subseteq \delta$ and $\mu \neq \delta$.

**THEOREM 1.7.5:** *Let $\mu$ be a fuzzy left (resp. right) k-ideal of R. Then two level left (resp. right) k-ideals $\mu_s$, $\mu_t$ (with $s < t$ in [0, 1] ) of $\mu$ are equal if and only if there is no $x \in R$ such that $s \leq \mu(x) < t$.*

*Proof:* Suppose $s < t$ in [0, 1] and $\mu_s = \mu_t$. If there exists an $x \in R$ such that $s \leq \mu(x) < t$, then $\mu_t$ is a proper subset of $\mu_s$, a contradiction. Conversely, suppose that there is no $x \in R$ such that $s \leq \mu(x) < t$. Note that $s < t$ implies $\mu_t \subseteq \mu_s$. If $x \in \mu_s$, then $\mu(x) \geq s$ and so $\mu(x) \geq t$ because $\mu(x) \not< t$. Hence $x \in \mu_t$ and $\mu_s = \mu_t$. This completes the proof.

For more about these concepts please refer [52]. Now for a given fuzzy k-ideal $\mu$ of a semiring R we denote by Im ($\mu$) the image set of $\mu$.

**THEOREM 1.7.6:** *Let $\mu$ be a fuzzy left (resp. right) k-ideal of R. If $Im(\mu) = \{t_1, t_2, \ldots, t_n\}$ where $t_1 < t_2 < \ldots < t_n$, then the family of left (resp. right) k-ideals $\mu_{t_i}$ ( $i = 1, 2, \cdots, n$ ) constitutes the collection of all left (resp. right) ideals of $\mu$.*

*Proof***:** Left as an exercise for the reader and requested to refer [52].

Now we proceed on to study some more notions.

Given any two sets R and S, let $\mu$ be a fuzzy subset of R and let $f : R \to S$ be any function. We define a fuzzy subset $\delta$ on S by

$$\delta(y) = \begin{cases} \sup_{x \in f^{-1}(y)} \mu(x) & \text{if } f^{-1}(y) \neq \phi, y \in S \\ 0 & \text{otherwise} \end{cases}$$

and we call $\delta$ the image of $\mu$ under f, written $f(\mu)$. For any fuzzy subset $\delta$ on $f(R)$ we define a fuzzy subset $\mu$ on R by $\mu(x) = \delta(f(x))$ for all $x \in R$, and we call $\mu$ the pre-image of $\delta$ under f which is denoted by $f^{-1}(\delta)$.

**THEOREM 1.7.7:** *An onto homomorphic pre image of a fuzzy left (resp. right) k-ideal is a fuzzy left (resp. right) k-ideal.*

*Proof:* Let $f : R \to S$ be an onto homomorphism. Let $\delta$ be a fuzzy left (resp. right) k-ideal on S and let $\mu$ be the pre image of $\delta$ under f. Then it was proved that $\mu$ is a fuzzy left (resp. right) ideal of R.

For any $x, y \in \delta$ we have

$$\begin{aligned} \mu(x) &= \delta(f(x)) \geq \min \{\max \{\delta(f(x) + f(y)), \delta(f(y) + f(x))\}, \delta f(y))\} \\ &= \min\{\max \{ \delta(f(x + y)), \delta(f(x + y))\}, \delta(f(y))\} \\ &= \min\{\max \{ \mu(x + y), \mu(y + x)\}, \mu(y)\}\end{aligned}$$

proving that $\mu$ is a fuzzy left (resp. right) k-ideal of R.

**THEOREM 1.7.8:** *Let f be a mapping from a set X to a set Y, and let $\mu$ be a fuzzy subset of X. Then for every $t \in (0, 1]$ $(f(\mu))_t = \bigcap_{0<s<t} f(\mu_{t-s})$.*

*Proof*: Straightforward, hence left for the reader to prove.

**THEOREM 1.7.9:** *Let $f : R \to S$ be an onto homomorphism and let $\mu$ be a fuzzy left (resp. right) k-ideal of R. Then the homomorphic image $f(\mu)$ of $\mu$ under f is a fuzzy left (resp. right) k-ideal of S.*

*Proof:* The proof can be got as a matter of routine; hence left for the reader to prove.

**DEFINITION 1.7.4:** *A left (resp.right) k-ideal I of R is said to be characteristic, if $f(I) = I$ for all $f \in Aut(R)$, where Aut(R) is the set of all automorphisms of R. A fuzzy left (resp. right) k-ideal $\mu$ of R is said to be fuzzy characteristic if $\mu(f(x)) = \mu(x)$ for all $x \in R$ and $f \in Aut(R)$.*

**THEOREM 1.7.10:** *Let $\mu$ be a fuzzy left (resp. right) k-ideal of R and let $f : R \to R$ be an onto homomorphism. Then the mapping $\mu^f: R \to [0, 1]$ defined by $\mu^f(x) = \mu(f(x))$ for all $x \in R$ is a fuzzy left (resp. right) k-ideal of R.*

*Proof:* It is a matter of routine once we write $\mu^f(x) = \mu(f(x)) \geq \min \{\max \{\mu(f(x) + f(y)), \mu(f(y) + f(x))\}, \mu(f(y))\}$.

The simplifications are direct and simple and hence left for the reader to prove.

The following theorems are also straightforward hence stated without proof; so that the interested reader can do them.

**THEOREM 1.7.11:** *If $\mu$ is fuzzy characteristic left (resp. right) k-ideal of R, then each level left (resp.right) k-ideal of $\mu$ is characteristic.*

**THEOREM 1.7.12:** *Let $\mu$ be a fuzzy left (resp. right) k-ideal of R and let $x \in R$. Then $\mu(x) = t$ if and only if $x \in \mu_t$ and $x \notin \mu_s$ for all $s > t$.*

We recall the proof of the following theorem.

**THEOREM 1.7.13:** *Let $\mu$ be a fuzzy left (resp. right) k-ideal of R. If each level left (resp right) k-ideal of $\mu$ is a characteristic then $\mu$ is fuzzy characteristic.*

*Proof*: Let $x \in R$ and $f \in Aut(R)$. If $\mu(x) = t \in [0, 1]$, then by the just above theorem we have $x \in \mu_t$ and $x \notin \mu_s$ for all $s > t$.

Since each level left (resp right) k-ideal of $\mu$ is characteristic, $f(x) \in f(\mu_t) = \mu_t$. Assume $\mu(f(x)) = s > t$. Then $f(x) \in \mu_s = f(\mu_s)$. Since f is one to one it follows that $x \in \mu_s$, a contradiction. Hence $\mu(f(x)) = t = \mu(x)$ showing that $\mu$ is fuzzy characteristic.

For more about fuzzy semiring properties refer [50, 51, 52, 54].

Now we proceed on to define fuzzy notions on bisemirings. As bisemirings themselves are very new much so is the concept of fuzzy bisemiring. Just we recall the definitions of them.

**DEFINITION 1.7.5:** *Let $(S, +, \bullet)$ be a non-empty set with two binary operations '+' and '$\bullet$'. We call $(S, +, \bullet)$ a bisemiring if*

  i.   $S = S_1 \cup S_2$ where both $S_1$ and $S_2$ are distinct subsets of $S$, $S_1 \not\subset S_2$ and $S_2 \not\subset S_1$.
  ii.  $(S_1, +, \bullet)$ is a semiring.
  iii. $(S_2, +, \bullet)$ is a semiring.

*Example 1.7.5:* $S = Z^0 \cup C_2$ where $C_2$ is a chain lattice of order 2 and $Z^0$ is a semiring, thus $S$ is a bisemiring.

**DEFINITION 1.7.6:** *Let $(S, +, \bullet)$ be a bisemiring. If the number of elements in $S$ is finite we call $S$ a finite bisemiring. If the number of elements in $S$ is infinite we say $S$ is an infinite bisemiring.*

**DEFINITION 1.7.7:** *Let $(S, +, \bullet)$ be a bisemiring. We say $S$ is a commutative bisemiring if both the semirings $S_1$ and $S_2$ are commutative otherwise we say the bisemiring $S$ is a non-commutative bisemiring. Let $(S, +, \bullet)$ be a bisemiring. We say $S$ is a strict bisemiring if both $S_1$ and $S_2$ are strict bisemirings where $S = S_1 \cup S_2$. Let $(S, +, \bullet)$ be a bisemiring if $0 \neq x \in S$ be a zero divisor if there exists a $y \neq 0 \in S$ such that $x \bullet y = 0$. We say a bisemiring $S$ has a unit 1 in $S$ if $a \bullet 1 = 1 \bullet a = a$ for all $a \in S$. Let $(S, +, \bullet)$ is a bisemiring we say $S$ is a bisemifield if $S_1$ is a semifield and $S_2$ is a semifield where $S = S_1 \cup S_2$. If both $S_1$ and $S_2$ are non-commutative semirings with no zero divisors then we call $S = S_1 \cup S_2$ to be a bisemidivision ring. It is to be noted that even if one of $S_1$ or $S_2$ is non-commutative with no zero divisors but other is a semifield still $S$ is a bisemidivision ring.*

*Example 1.7.6:* $(S, +, \bullet)$ is a bisemiring. $S = S_1 \cup S_2$ where $S_1 = C_7$ and $S_2$ is the lattice given by then $S$ is only a bisemiring which is not a semifield.

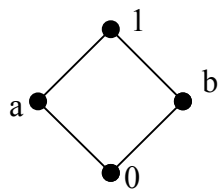

**Figure 1.7.1**

**DEFINITION 1.7.8:** *Let $(S = S_1 \cup S_2, +, \bullet)$ be a bisemiring which is commutative and has unit element. The polynomial bisemiring is denoted by $S[x] = S_1[x] \cup S_2[x]$ where $S_1[x]$ and $S_2[x]$ are polynomial semirings.*

**DEFINITION 1.7.9:** *Let (S, +, •) be a bisemiring. (P, +, •) be a proper sub-bisemiring of S. We call P a bi-ideal if for all $s \in S$ and $p \in P$, sp and ps $\in P$ i.e. $P = P_1 \cup P_2$ with $s_i \in S_i$, $p_i s_i$ and $s_i p_i \in P_i$; i = 1, 2.*

**DEFINITION 1.7.10:** *Let (S, +, •) and (S', +, •) be two bisemirings, where $S = S_1 \cup S_2$ and $S' = S'_1 \cup S'_2$. We say a map $\phi : S \to S'$ is a bisemiring homomorphism if $\phi_1 : S_1 \to S'_1$ and $\phi_2 : S_2 \to S'_2$ is a semiring homomorphism or $\phi_1 : S_1 \to S'_2$ and $\phi_2 : S_2 \to S'_1$ is a semiring homomorphism where we denote $\phi$ just by default of notation as $\phi = \phi_1 \cup \phi_2$.*

As our motivation is to define fuzzy bisemirings as fuzzy bisemirings is not defined till date we give importance only to the study of fuzzy bisemirings and Smarandache fuzzy bisemirings.

**DEFINITION 1.7.11:** *($S = S_1 \cup S_2$, +, •) is a bisemiring. We say $\mu : S \to [0, 1]$ is a fuzzy bisemiring, $\mu = \mu_1 \cup \mu_2$ where both $\mu_1 : S_1 \to [0, 1]$ and $\mu_2 : S_2 \to [0, 1]$ are fuzzy semirings.*

Now we proceed on to define fuzzy bi-ideals in a bisemiring.

**DEFINITION 1.7.12:** *A fuzzy subset $\mu$ of a bisemiring ($S = S_1 \cup S_2$, +, •) is said to be a fuzzy left (resp. right) bi-ideal of S if $\mu : S \to [0, 1]$ such that $\mu(x + y) \geq \min \{\mu(x) + \mu(y)\}$, and $\mu(xy) > \mu(y)$ (resp. $\mu(xy) > \mu(x)$) for all $x, y \in S$. $\mu$ is a fuzzy bi-ideal of S if it is both a fuzzy left and fuzzy right bi-ideal of $S = S_1 \cup S_2$.*

$\mu_1 : S_1 \to [0, 1]$ and
$\mu_2 : S_2 \to [0, 1]$ where $\mu_1$ and $\mu_2$ are fuzzy ideals of the semirings $S_1$ and $S_2$.

**DEFINITION 1.7.13:** *A fuzzy ideal $\mu$ of a bisemiring (S, +, •) is said to be a fuzzy k-ideal of S if $\mu(x) \geq \min \{\{\max \mu(x + y), \mu(y + x)\}, \mu(y)\}$ for all $x, y \in S$. i.e. $\mu(x) \geq \min \{\max \mu_1(x + y), \mu_1(y + x)\}, \mu_1(y)\} \cup \min\{\max \mu_2(x + y), \mu_2(y + x)\}, \mu_2(y)\}$ for all $x, y \in S_1$ ($x, y \in S_2$) where $S = S_1 \cup S_2$. If S is additively commutative then the condition reduces to $\mu(x) \geq \min\{\mu(x + y), \mu(y)\}$ for all $x, y \in S$; where $\mu(x) = \mu_1(x) \cup \mu_2(x) \geq \min \{\mu_1(y + x)\}, \mu_1(y)\} \cup \min \{\mu_2(x + y), \mu_2(y)\}$ for all $x, y \in S_1$ ($x, y \in S_2$).*

**DEFINITION 1.7.14:** *Let $I = I_1 \cup I_2$ be a non-empty subset of a bisemiring $S = S_1 \cup S_2$ and $\lambda_I = \lambda_{I_1} \cup \lambda_{I_2}$ the characteristic function of $I = I_1 \cup I_2$.*

(Here $\lambda_{I_1}$ the characteristic function of $I_1$ and $\lambda_{I_2}$ the characteristic function of $I_2$). Then $I = I_1 \cup I_2$ is a k-bi-ideal of $S = S_1 \cup S_2$ if and only if $\lambda_I$ is a fuzzy k-bi-ideal of $S = S_1 \cup S_2$.

Now we proceed on to give some interesting results on bi-ideals.

**THEOREM 1.7.14:** *A fuzzy subset $\mu = \mu_1 \cup \mu_2$ of $S = S_1 \cup S_2$ is a fuzzy left (resp. right) k-bi-ideal of S if and only if for any $t \in [0, 1]$ such that $\mu_t \neq \phi$ ($\mu_t = (\mu_1)_t \cup$* $(\mu_2)_t$) $\mu_t$ is left (resp. right) k-bi-ideal of $S = S_1 \cup S_2$ where $\mu_t = \{x \in S \mid \mu(x) \geq t\} = \{x \in S_1 \mid \mu_1(x) \geq t\} \cup \{x \in S_2 \mid \mu_2(x) \geq t\}$ which is called a level subset of $\mu$.

*Proof:* Follows as in case of semirings, by taking each $\mu_i$; $i = 1, 2$. Then as $\mu = \mu_1 \cup \mu_2$ we get the desired result.

**THEOREM 1.7.15:** *Let I be any left (resp. right) k-bi-ideal of $S = S_1 \cup S_2$. Then there exists a fuzzy left (resp. right) k-bi-ideal $\mu = \mu_1 \cup \mu_2$ of $S = S_1 \cup S_2$ such that $\mu_t = (\mu_1)_t \cup (\mu_2)_t$; $I = I_1 \cup I_2$ for some $t \in [0, 1]$.*

*Proof:* If we define a fuzzy subset of $S = S_1 \cup S_2$ by $\mu = \mu_1 \cup \mu_2$

$$\mu_1(x) = \begin{cases} t & \text{if } x \in I_1 \\ 0 & \text{otherwise} \end{cases}$$

$$\mu_2(x) = \begin{cases} t & \text{if } x \in I_2 \\ 0 & \text{otherwise} \end{cases}$$

for some $t \in [0, 1]$, then it follows that $\mu_t = I = I_1 \cup I_2$. For a given $s \in [0, 1]$ we have

$$(\mu_1)_s = \begin{cases} (\mu_1)_0 \ (= S_1) & \text{if } s = 0 \\ (\mu_1)_t \ (= I_1) & \text{if } s \leq t. \\ \phi & \text{if } t < s \leq 1. \end{cases}$$

The result follows as in case of fuzzy semiring.

Let $\mu = \mu_1 \cup \mu_2$ and $\delta = \delta_1 \cup \delta_2$ be fuzzy subsets of $S = S_1 \cup S_2$. We say $\mu \subseteq \delta$ if and only if ($\mu_1 \subset \delta_1$, $\mu_2 \subset \delta_2$) if and only if $\mu(x) \leq \delta(x)$, $\mu_1(x) \leq \delta_1(x)$, $\mu_2(x) \leq \delta_2(x)$) for all $x \in S$ and $\mu \subset \delta$ if and only if $\mu \subseteq \delta$ and $\mu \neq \delta$.

**THEOREM 1.7.16:** *Let $\mu = \mu_1 \cup \mu_2$ be a fuzzy left (resp. right) k-bi-ideal of $S = S_1 \cup S_2$. Then two level left (resp. right) k-bi-ideals $\mu_s$, $\mu_t$ (with $s < t$ in $[0, 1]$) of $\mu = \mu_1 \cup \mu_2$ are equal if and only if there is no $x \in S$ such that $s \leq \mu(x) < t$ (i.e. $s \leq \mu_1(x) < t$ and $s \leq \mu_2(x) < t$).*

*Proof:* As in case of semirings.

**THEOREM 1.7.17:** *Let $\mu = \mu_1 \cup \mu_2$ be a fuzzy left (resp. right) k-bi-ideal of $S = S_1 \cup S_2$. If Im $\mu = $ Im $\mu_1 \cup$ Im $\mu_2 = \{t_1, t_2, \ldots, t_n\}$ where $t_1 < t_2 < \ldots < t_n$, then the family of left (resp right) k-bi-ideals $\mu_{t_i}$ ($i = 1, 2, \cdots, n$) constitutes the collection of all level left (resp. right) ideals of $\mu = \mu_1 \cup \mu_2$.*

*Proof:* Direct and the proof is got as a matter of routine as in case of semirings.

Given any two sets $R = R_1 \cup R_2$ and $S = S_1 \cup S_2$, let $\mu$ be a fuzzy subset of $R$ and let $f: R \to S$ be any function $f = f_1 \cup f_2$ where $f_1 : R_1 \to S_1$ and $f_2 : R_2 \to S_2$.

We define a fuzzy subset $\delta = \delta_1 \cup \delta_2$ of $S = S_1 \cup S_2$ by

$$\delta_1(y) = \begin{cases} \sup_{x \in f_1^{-1}(x)} \mu_1(x) & \text{if } f_1^{-1}(y) \neq \phi, y \in S_1 \\ 0 & \text{otherwise} \end{cases}$$

and

$$\delta_2(y) = \begin{cases} \sup_{x \in f_2^{-1}(x)} \mu_2(x) & \text{if } f_2^{-1}(y) \neq \phi, y \in S_2 \\ 0 & \text{otherwise} \end{cases}$$

and we call $\delta = \delta_1 \cup \delta_2$ the image of $\mu = \mu_1 \cup \mu_2$ under $f = f_1 \cup f_2$ written as $f(\mu) = f_1(\mu_1) \cup f_2(\mu_2)$. for any fuzzy subset $\delta$ on $f(R) = f_1(R_1) \cup f_2(R_2)$ we define fuzzy subset $\mu = \mu_1 \cup \mu_2$ on $R = R_1 \cup R_2$ by $\mu(x) = \delta f(x)$ i.e. $\mu(x) = \mu_1(x) \cup \mu_2(x)$ where $\mu_1(x) = \delta_1(f_1(x))$ for $x \in R_1$ and $\mu_2(x) = \delta_2(f_2(x))$ for $x \in R_2$ and we call $\mu = \mu_1 \cup \mu_2$ the pre image of $\delta = \delta_1 \cup \delta_2$ under $f = f_1 \cup f_2$ which is denoted by $f^{-1}(\delta) = f_1^{-1}(\delta_1) \cup f_2^{-1}(\delta_2)$.

Now we state the following theorem on pre image.

**THEOREM 1.7.18:** *An onto homomorphism pre image of a fuzzy left (resp. right) k-bi-ideal is a fuzzy left (resp. right) k-bi-ideal.*

*Proof:* $f : R \to S$ be an onto homomorphism of the bisemirings $R$ and $S$ where $R = R_1 \cup R_2$ and $S = S_1 \cup S_2$, $f = f_1 \cup f_2$. $f_1 : R_1 \to S_1$ and $f_2 : R_2 \to S_2$ be onto semiring homomorphism. Let $\delta$ be a fuzzy left (resp. right) k-bi-ideal on $S = S_1 \cup S_2$ and let $\mu = \mu_1 \cup \mu_2$ be the pre image of $\delta$ under $f$. Then it was proved that $\mu = \mu_1 \cup \mu_2$ is a fuzzy left (resp. right) bi-ideal of $R = R_1 \cup R_2$. For any $x, y \in S = S_1 \cup S_2$ we have $\mu(x) = \delta(f(x))$ proving $\mu = \mu_1 \cup \mu_2$ is a fuzzy left (resp. right) k-bi-ideal.

**THEOREM 1.7.19:** *Let $f = f_1 \cup f_2$ be a mapping from a set $X = X_1 \cup X_2$ to a set $Y = Y_1 \cup Y_2$ and let $\mu = \mu_1 \cup \mu_2$ be a fuzzy subset $X = X_1 \cup X_2$. Then for every $t \in (0, 1]$.*

$$(f(\mu))_t = (f_1(\mu_1))_t \cup (f_2(\mu_2))_t = \left[\bigcap_{0<s<t} f_1((\mu_1))_{t-s}\right] \cup \left[\bigcap_{0<s<t} f_2((\mu_2))_{t-s}\right]$$

*Proof:* Follows as in case of semirings.

**THEOREM 1.7.20:** *Let $f : R \to S$ ($f = f_1 \cup f_2$, $R = R_1 \cup R_2$, $S = S_1 \cup S_2$, $f_1: R_1 \to S_1$; $f_2 : R_2 \to S_2$) be an onto homomorphism, let $\mu = \mu_1 \cup \mu_2$ be a fuzzy left (respectively right) k-bi-ideal of $R = R_1 \cup R_2$. Then the homomorphic image $f(\mu) = f_1(\mu_1) \cup f_2(\mu_2)$ of $\mu = \mu_1 \cup \mu_2$ under $f = f_1 \cup f_2$ is a fuzzy left (resp. right) k-bi-ideal of $S = S_1 \cup S_2$.*

*Proof:* Follows as a matter of routine as in case of semirings.

**DEFINITION 1.7.15:** *A left (resp. right) k-bi-ideal $I = I_1 \cup I_2$ of the bisemiring $R = R_1 \cup R_2$ is said to be characteristic if $f(I) = I$ i.e. $f_1(I_1) = I_1$ and $f_2(I_2) = I_2$ for all $f = f_1 \cup f_2 \in Aut(R)$ where $Aut(R)$ is the set of all automorphism of $R = R_1 \cup R_2$. A fuzzy left (resp. right) k-bi-ideal $\mu = \mu_1 \cup \mu_2$ of $R = R_1 \cup R_2$ is said to be a fuzzy characteristic if $\mu(f(x)) = \mu(x)$ (i.e. $\mu_1(f_1(x)) = \mu_1(x)$; $\mu_2(f_2(x)) = \mu_2(x)$) for all $x \in R_1$ ($x \in R_2$ and $f_1 \in Aut(R_1)$ and $f_2 \in Aut(R_2), f = f_1 \cup f_2$.*

**THEOREM 1.7.21:** *Let $\mu = \mu_1 \cup \mu_2$ be a fuzzy left (resp. right) k-bi-ideal of $R = R_1 \cup R_2$ and let $f : R \to R$ be an onto homomorphism here $f = f_1 \cup f_2$; $f_1 : R_1 \to R_1$ and $f_2 : R_2 \to R_2$. Then the mapping $\mu^f : R \to [0, 1]$;*

$$\left(\mu_1^{f_1} : R_1 \to [0, 1] \text{ and } \mu_2^{f_2} : R_2 \to [0, 1]\right)$$

*defined by $\mu^f(x) = \mu(f(x))$ (where $\mu_1^{f_1}(x) = \mu_1(f_1(x))$ and $\mu_2^{f_2}(x) = \mu_2(f_2(x))$) for all $x \in R_1$ (or $x \in R_2$) is a fuzzy left (resp. right) k-bi-ideal of $R = R_1 \cup R_2$.*

*Proof:* Follows as in case of semirings.

**THEOREM 1.7.22:** *If $\mu = \mu_1 \cup \mu_2$ is a fuzzy characteristic left (resp. right) k-bi-ideal of $R = R_1 \cup R_2$ then each level left (resp. right) k-bi-ideal of $\mu$ is characteristic.*

*Proof:* Left for the reader to work as in case of semiring.

The following theorem is as easy consequence of the above theorem.

**THEOREM 1.7.23:** *Let $\mu = \mu_1 \cup \mu_2$ be a fuzzy left (resp. right) k-bi-ideal of the bisemiring of $R = R_1 \cup R_2$ and let $x \in R$. Then $\mu(x) = t$ if and only if $x \in \mu_t$ and $x \notin \mu_s$ for all $s > t$. ($x \in R_1$ then $\mu_1(x) = t$, $x \in (\mu_1)_t$ and $x \notin (\mu_1)_s$ for all $s > t$), if $x \in R_2$ then $\mu_2(x) = t$, $x \in (\mu_2)_t$ and $x \notin (\mu_2)_s$ for all $s > t$).*

*Proof:* Follows as in case of semirings.

**THEOREM 1.7.24:** *Let $\mu = \mu_1 \cup \mu_2$ be a fuzzy left (resp. right) k-bi-ideal of the bisemiring of $R = R_1 \cup R_2$. If each level left (resp. right) k-bi-ideal of $\mu$ is characteristic then $\mu$ is fuzzy characteristic.*

*Proof:* Follows as in case of semirings.

All results can be got in case of fuzzy bisemirings but as the main motivation is only study of Smarandache fuzzy semirings we do not exhaust all results regarding the fuzzy semirings or fuzzy bisemirings.

Any matrix A in $M_n(F)$ will be called a fuzzy matrix we may call $M_n(F)$ a fuzzy matrix semiring. Several results can be had in this direction.

Further the fuzzy semirings can be utilized in constraint handling rules. For more about this please refer [51-54]. All applications of fuzzy semirings can be had from other research papers.

## 1.8 Fuzzy near-rings and their properties

In this section we recall the definition of fuzzy near-ring, give some new types of fuzzy near-rings and introduce the notion of fuzzy bi near-rings. Also notions like Fuzzy near-ring module and fuzzy congruence of a near-ring module are recalled in this section. For more about fuzzy near-ring literature please refer [28, 38, 56, 57, 71, 122, 130].

**DEFINITION [71]**: *Let R be a near-ring and N a fuzzy set in R. Then N is called a fuzzy near-ring in R if*

  i. $N(x + y) \geq \min\{N(x), N(y)\}$.
  ii. $N(-x) \geq N(x)$.
  iii. $N(xy) \geq \min\{N(x), N(y)\}$ *for all x, y in R.*

**DEFINITION [71]**: *Let R be a near-ring and N a fuzzy near-ring in R. Let Y be a near-ring module over R and M a fuzzy set in Y. Then M is called a fuzzy near-ring module in Y if*

  i. $M(x + y) \geq \min\{M(x), M(y)\}$.
  ii. $M(\lambda x) \geq \min\{N(\lambda), M(x)\}$ *for all x, y $\in$ Y and $\lambda \in$ R.*
  iii. $M(0) = 1$.

*If N is an ordinary near-ring then condition (ii) in the above definition is replaced by*

  ii(a). $M(\lambda x) > M(x)$ *for all $\lambda \in$ N and for all x $\in$ Y.*

**THEOREM [71]**: *Let Y be a near-ring module over a fuzzy near-ring N in R. Then M is a fuzzy near-ring module in Y if and only if $M(\lambda(x) + \mu(x)) \geq \min\{\min\{N(\lambda), M(x)\}, \min\{N(\mu), M(y)\}\}$ for all $\lambda, \mu \in$ N and for all x, y $\in$ Y.*

*If N is an ordinary near-ring then the above condition is replaced by $M(\lambda x + \mu y) \geq \min\{\{M(x), N(y)\}$ for all x, y $\in$ Y.*

*Proof*: Left for the reader as an exercise as it can be got directly by the definitions.

**THEOREM [71]**: *Let Y be a near-ring module over a near-ring R with identity. If M is a fuzzy near-ring module in Y and if $\lambda \in$ R is invertible then $M(\lambda x) = M(x)$ for all x $\in$ Y.*

*Proof*: If $\lambda \in$ M is invertible then we have for all x $\in$ Y. $M(x) = M(\lambda^{-1} \lambda x) \geq M(\lambda x) \geq M(x)$ and so $M(\lambda x) = M(x)$.

**THEOREM [71]**: *Let $\{M_i \mid i \in I\}$ be a family of fuzzy near-ring modules in Y. Then $\bigcap_{i \in I} M_i$ is a fuzzy near-ring module in Y.*

*Proof*: Let $M = \bigcap_{i \in I} M_i$; then we have for all $\lambda \in R$ and for all $x, y \in Y$.

$$\begin{aligned} M(x+y) &= \inf_{i \in I} M_i(x+y) \\ &\geq \inf_{i \in I} \{\min\{M_i(x), M_i(y)\}\} \\ &= \min\left\{\inf_{i \in I} M_i(x), \inf M_i(y)\right\} \\ &= \min\{M(x), M(y)\} \end{aligned}$$

and

$$\begin{aligned} M(\lambda x) &= \inf_{i \in I} M_i(\lambda x) \\ &\geq \inf_{i \in I} \{\min\{N(\lambda), M_i(x)\}\} \\ &= \min\left\{N(\lambda), \inf_{i \in I} M_i(x)\right\} \\ &= \min\{N(\lambda), M(x)\}. \end{aligned}$$

**THEOREM [71]**: *Let Y and W be near-ring modules over a fuzzy near-ring N in a near-ring R and $\theta$ a homomorphism of Y into W. Let M be a fuzzy near-ring module in W. Then the inverse image $\theta^{-1}(M)$ of M is a fuzzy near-ring module in Y.*

*Proof*: For all $x, y \in Y$ and for all $\lambda, \mu \in R$, we have

$$\begin{aligned} \theta^{-1}(M)(\lambda x + \mu y) &= M(\theta(\lambda x + \mu y)) \\ &= M(\lambda \theta(x) + \mu \theta(y)) \\ &\geq \min\{\min\{N(\lambda), M(\theta(x))\}, \min\{N(\mu), M\theta(y))\}\} \\ &= \min\{\min\{N(\lambda), \theta^{-1}(M)(x)\}, \min\{N(\mu), \theta^{-1}(M)(y)\}\} \end{aligned}$$

By earlier theorems $\theta^{-1}(M)$ is a fuzzy near-ring module in Y. Hence the theorem.

We say that a fuzzy set A in M has the sup property if for any subset T of M there exists $t_o \in T$ such that

$$A(t_o) = \sup_{t \in T} A(t).$$

**THEOREM [71]**: *Let Y and W be near-ring modules over a fuzzy near-ring N in a near-ring R and $\theta$ a homomorphism of Y into W. Let W be a fuzzy near-ring module in Y that has the sup property. Then the image $\theta(M)$ of M is a fuzzy near-ring module in W.*

*Proof*: Let $\mu, \nu \in W$. It either $\theta^{-1}(\mu)$ or $\theta^{-1}(\nu)$ is empty, then the result holds good. Suppose that neither $\theta^{-1}(\mu)$ nor $\theta^{-1}(\nu)$ is empty, then we have $\theta(M)(\lambda(u) + \mu\nu) =$

$$\sup_{\omega \in \theta^{-1}(\lambda u + \mu\nu)} M(\omega) \geq \min\{\min\{N(\lambda), \theta(m)(u)\}, \min\{N(\mu), \theta(M)(\nu)\}\}.$$

Hence the result.

Throughout the discussion from now on by a near-ring R we mean a system with two binary operations addition and multiplication such that

i. The elements of R form a group under '+'.
ii. The elements of R form a semigroup under multiplication '•'.
iii. $x \bullet (y + z) = x \bullet y + x \bullet z$ for all $x, y, z \in R$.

An R-module that is near-ring module is a system consisting of an additive group M, a near-ring R, and a mapping $(m, r) \to mr$ of $M \times R$ into M such that

i. $m(x + y) = mx + my$ for all $m \in M$ and for all $x, y \in R$.
ii. $m(xy) = (mx)y$ for all $m \in M$ and for all $x, y \in R$.

**DEFINITION 1.8.1**: *An R-homomorphism f of an R-module M into an R-module M' is a mapping from M to M' such that $(m_1 + m_2)f = m_1 f + m_2 f$ and $(mf)r = (mr)f$ for all $m, m_1, m_2 \in M$ and for all $r \in R$.*

The submodules of an R-module M are defined to be kernels of R-homomorphism.

**THEOREM [22]**: *An additive normal subgroup B of an R-module M is a submodule if and only if, $(m + b)r - mr \in B$ for all $m \in M, b \in B$ and $r \in R$.*

**DEFINITION 1.8.2**: *A relation $\rho$ on an R-module M is called a congruence on M if it is an equivalence relation on M such that $(a, b) \in \rho$ and $(c, d) \in \rho$ imply that $(a + c, b + d) \in \rho$ and $(ar, br) \in \rho$ for all a, b, c, d in M and for all r in R.*

**DEFINITION [98]**: *A non empty fuzzy subset $\mu$ of an additive group G is called a fuzzy normal subgroup of G if*

i. *$\mu(x + y) \geq \min\{\mu(x), \mu(y)\}$ for all x, y in G.*
ii. *$\mu(-x) = \mu(x)$ for all $x \in G$.*
iii. *$\mu(y + x - y) = \mu(x)$ for all x, y in G.*

*Let µ be a fuzzy normal subgroup of an additive group G and x ∈ G. Then the fuzzy subset x + µ of G defined by (x + µ) (x) = µ (y – x) for all y in G is called the fuzzy coset of µ .*

**DEFINITION 1.8.3**: *Let µ be a non-empty fuzzy subset of an R-module M. Then µ is said to be a fuzzy submodule of M if*

    i.  *µ is a fuzzy normal subgroup of M and*
    ii.  *µ {(x + y) r – x r} ≥ µ (y) for all x, y in M and for all r in R .*

The proof of the following two theorems are omitted as they can be easily obtained by the reader using direct methods.

**THEOREM 1.8.1**: *Let B be a non-empty subset of an R-module M. Then the characteristic function $\chi_g$ is a fuzzy submodule of M if and only if B is a submodule of M.*

**THEOREM 1.8.2**: *Let µ be a fuzzy submodule of an R-module M. Then the level subset $µ_t$ = {x ∈ M | µ (x) ≥ t}, t ∈ Im µ is a submodule of M.*

**DEFINITION 1.8.4**: *Let µ be a fuzzy submodule of an R-module M. Then the submodules $µ_I$'s are called level submodules of M.*

**THEOREM 1.8.3**: *For a non-empty fuzzy subset µ of an R-module M, the following assertions are equivalent*

    i.   *µ is a fuzzy submodule of M.*
    ii.  *The level subsets $µ_t$, t ∈ Im µ are submodules of M.*

*Proof*: It is a simple matter of routine.

**THEOREM 1.8.4**: *If µ is a fuzzy normal subgroup of an additive group G. Then x + µ = y + µ if and only if µ(x – y) = µ(0) for all x, y in G.*

*Proof*: Left for the reader as an exercise.

**THEOREM 1.8.5**: *Let µ be a fuzzy submodule of an R-module M. Then the set M / µ of all fuzzy cosets of µ is an R-module with respect to the operations defined by (x + µ) + (y + µ) = (x + y) + µ and (x + µ) r = xr + µr for all x, y in M and for all r in R. If f is a mapping from M to M/ µ, defined by xf = x + µ for all x ∈ M then f is an R-epimorphism with ker f = { x ∈ M /µ(x) = µ (0)}.*

*Proof*: Please refer [28, 39].

**DEFINITION [28]**: *The R-module M / µ is called the quotient R-module of M over its fuzzy submodule µ.*

**DEFINITION [28]:** *Let M be an R-module. A non empty fuzzy relation $\alpha$ on M [i.e. a mapping $\alpha: M \times M \to [0,1]$] is called a fuzzy equivalence relation if*

i. $\alpha(x,x) = \sup_{y,z \in M} \alpha(y,z)$ *for all $x, y, z$ in M (fuzzy reflexive).*
ii. $\alpha(x, y) = \alpha(y, x)$ *for all $x, y$ in M (fuzzy symmetric).*
iii. $\alpha(x,y) \geq \sup_{z \in M} [\min \alpha(x,z), \alpha(z,y)]$ *for all $x, y$ in M (fuzzy transitive).*

**DEFINITION 1.8.5**: *A fuzzy equivalence relation $\alpha$ on an R-module M is called a fuzzy congruence if $\alpha(a + c, b + d) \geq \text{Min}\,[\alpha(a, b), \alpha(c, d)]$ and $\alpha(ar, br) \geq \alpha(a, b)$ for all a, b, c, d in M and all r in R.*

The following theorem is left for the reader as an exercise.

**THEOREM 1.8.6**: *Let $\rho$ be a relation on an R-module M and $\lambda_\rho$ be its characteristic function. Then $\rho$ is a congruence relation on M if and only if $\lambda_\rho$ is a fuzzy congruence on M.*

**DEFINITION [28]**: *Let $\alpha$ be a fuzzy relation on an R-module M. For each $t \in [0, 1]$ the set $\alpha_t = \{(a,b) \in M \times M; \alpha(a,b) \geq t\}$ is called a level relation on $\alpha$.*

**THEOREM [28]**: *Let $\alpha$ be a fuzzy relation on an R-module M. Then $\alpha$ is a fuzzy congruence on M if and only if $\alpha_t$ is a congruence on M for each $t \in \text{Im}\,\alpha$.*

*Proof*: Refer [28].

**THEOREM [28]**: *Let $\alpha$ be a fuzzy congruence on an R-module M and $\mu_\alpha$ be a fuzzy subset of M, defined by $\mu_\alpha(a) = \alpha(a, 0)$, $a \in M$. Then $\mu_\alpha$ is a fuzzy submodule of M.*

*Proof*: Since $\mu_\alpha(0) = \alpha(0, 0) = \sup_{x,y \in M} \alpha(x, y) \neq 0$ (as $\alpha$ is non empty) it follows that $\mu_t$ is non-empty. For a, b in M.

$$\mu_\alpha(a+b) = \alpha(a+b, 0) \geq \min\{\alpha(a,0), \alpha(b,0)\} = \min[\mu_\alpha(a), \mu_\alpha(b)]$$
$$\mu_\alpha(-a) = \alpha(-a, 0) = \alpha(-a+0, -a+a) \geq \text{Min}[\alpha(-a, -a), \alpha(0, a)]$$
$$= \alpha(0, a) = \alpha(a, 0) = \mu_\alpha(a).$$

Similarly $\mu_\alpha(a) \geq \mu_\alpha(-a)$. Thus $\mu_\alpha(-a) = \mu_\alpha(a)$. Again $\mu_\alpha(a + b - a)$

$$= \alpha(a+b-a, 0) = \alpha(a+b-a, a+0-a)$$
$$\geq \alpha(b, 0) = \mu_\alpha(b).$$

So $\mu_\alpha$ is a fuzzy normal subgroup of the R-module M. Now for a, b in M and r in R.

$$\mu_\alpha\{(a+b)r - ar\} = \alpha\{(a+b)r - ar, 0\} = \alpha\{(a+b)r - ar, ar - ar\}$$

$$\geq \text{Min} \left[\alpha\{(a+b)r, ar\}, \alpha(-ar, -ar)\right]$$
$$= \alpha\{(a+b)r, ar\} \geq \alpha(a+b, a) \geq \alpha(b, 0) = \mu_\alpha(b).$$

Thus $\mu_\alpha$ is a fuzzy submodule of M.

**THEOREM 1.8.7**: *Let $\mu$ be a fuzzy submodule of an R-module M. Let $\alpha_\mu$ be the fuzzy relation on M, defined by $\alpha_\mu(x, y) = \mu(x - y)$ for x, y in M. Then $\alpha_\mu$ is a fuzzy congruence on M.*

*Proof*: Since $\mu$ is non-empty it follows that $\alpha_\mu$ is non empty. Now $\alpha_\mu(x, x) = \mu(0) \geq \mu(y - z)$ for all y, z in M $= \alpha_\mu(y, z)$. So

$$\alpha_\mu(x, x) = \sup_{y, z \in M} \alpha_\mu(y, z).$$

Thus $\alpha_\mu$ is fuzzy reflexive. It is clear that $\alpha_\mu$ is fuzzy symmetric. Again $\alpha_\mu(x, y) = \mu(x - y) = \mu(x - z + z - y) \geq \text{Min}[\mu(x - z), \mu(z - y)]$ for all z in M. So

$$\alpha_\mu(x, y) \geq \sup_{z \in M} \left[\text{Min}\left[\alpha_\mu(x, z), \alpha_\mu(z, y)\right]\right].$$

Thus $\alpha_\mu$ is a fuzzy equivalence relation on M.

Now

$$\begin{aligned}
\alpha_\mu(x + u, y + v) &= \mu(x + u - v - y) \\
&= \mu(-y + x + u - v) \\
&\geq \text{Min}[\mu(-y + x), \mu(u - v)] \\
&= \text{Min}[\mu(x - y), \mu(u - v)] \\
&= \text{Min}[\alpha_\mu(x, y), \alpha_\mu(u, v)].
\end{aligned}$$

Again $x_\mu(xr, yr) = \mu(xr - yr) = \mu\{(y - y + x)r - yr\} \geq \mu(-y + x) = \mu(x - y) = \alpha_\mu(x, y)$. Hence $\alpha_\mu$ is a fuzzy congruence on M.

**Note**: $\alpha_\mu$ is called the fuzzy congruence induced by $\mu$ and $\mu_\alpha$ is called the fuzzy submodule induced by $\alpha$.

Here in the following theorem FS (M) set of all fuzzy submodules of M and FC (M) the set of all fuzzy congruences on M. The reader is expected to refer [28] for proof.

**THEOREM 1.8.8**: *Let M be an R-module. Then there exists an inclusion preserving bijection from the set FS(M) of all fuzzy submodules of M to the set FC(M) of all fuzzy congruences on M.*

**THEOREM 1.8.9:** *Let $\alpha$ be a fuzzy congruences on an R-module M and $\mu_\alpha$ be the fuzzy submodule induced by $\alpha$. Let $t \in \text{Im } \alpha$. Then $(\mu_\alpha)_t = \{x \in M \,/\, x \equiv 0 \,(\alpha_t)\}$ is the submodule induced by the congruence $\alpha_t$.*

*Proof*: Let a ∈ M. Now a ∈ $(\mu_\alpha)_t$ if and only if $(\mu_\alpha)$ (a) ≥ t if and only if α(a, 0) ≥ if and only if (a, 0) ∈ $\alpha_t$ if and only if α ∈ $0(\alpha_t)$ if and only if α ∈ {x ∈ M | x ≡ $0(\alpha_t)$}. Hence the theorem.

**THEOREM 1.8.10**: *Let μ be a fuzzy submodule of an R-module M and $\alpha_\mu$ be the fuzzy congruence induced by μ. Let t ∈ Im μ. Then $(\alpha_\mu)_t$ is the congruence on M induced by $\mu_t$.*

*Proof*: Let β be the congruence on the R-module M induced by $\mu_t$. (x, y) ∈ β if and only if x – y ∈ $\mu_t$. Let (x, y) ∈ $(\alpha_\mu)_t$. Therefore $(\alpha_\mu)$ (x, y) ≥ t ⇒ μ (x – y) ≥ t ⇒ (x – y) ∈ $\mu_t$ ⇒ (x, y) ∈ β. Thus $(\alpha_\mu)_t$ ⊆ β. By reversing the above argument we get β ⊆ $(\alpha_\mu)$. Hence $(\alpha_\mu)_t$ = β.

**DEFINITION [28]**: *Let M be an R-module and α be a fuzzy congruence on M. A fuzzy congruence β on M is said to be α-invariant if α (x, y) = α (u, v) implies that β (x, y) = β (u, v) for all (x, y), (u, v) ∈ M × M.*

**THEOREM [28]**: *Let M be a R-module and μ be a fuzzy submodule of M. Let α be the fuzzy congruence on M induced by μ. Then the fuzzy relation α / α on M / μ defined by α / α (x + μ, y + μ) = α (x, y) is a fuzzy congruence on M / μ.*

*Proof:* Assume that x + μ = u + μ and y + μ = u + μ. Then μ (x – u) = μ(0) and μ(y – v) = μ (0).

Thus

$$\alpha(x, u) = \sup_{p,q \in M} \alpha(p,q) \text{ and}$$

$$\alpha(y, v) = \sup_{p,q \in M} \alpha(p,q).$$

Now

| α(x, y) | ≥ | Min [α(x, u), (u, y)] |
|---|---|---|
| | = | α(u, y) |
| | ≥ | Min [α(u, v), α(u, y)] |
| | = | α(u, v). |

Similarly α(u, v) ≥ α(x, y). Thus α(x, y) = α(u, v). Thus, α / α is meaningful. Rest of the proof is a matter of routine verification.

**THEOREM [28]**: *Let M be an R-module and μ be a fuzzy submodule of M. Let α be a the fuzzy congruence on M induced by μ. Then there exists a one to one correspondence between the set $FC_\alpha$ (M) of α-invariant fuzzy congruences on M and the set $FC_{\alpha/\alpha}$ (M / μ) of α /α invariant fuzzy congruences on M / μ.*

*Proof:* Refer [28].

**THEOREM 1.8.11:** *Let M be an R-module and $\mu$ be a fuzzy submodule of M. Let $\alpha$ be the fuzzy congruence on M induced by $\mu$. Let $t = \sup \operatorname{Im} \alpha$ then $M/\mu \cong M/\alpha_t$.*

*Proof*: We define a mapping $\theta: M/\mu \to M/\alpha_t$ by $(x + \mu)\theta = x\alpha_t$ where $x\alpha_t$ denotes the congruence class containing x of the congruence $\alpha_t$, $x + \mu = y + \mu \Rightarrow \mu(x - y) = \mu(0) = \alpha(x, y) = \sup \operatorname{Im} \alpha = t \Rightarrow (x, y) \in \alpha_t \Rightarrow x\alpha_t = y\alpha_t \Rightarrow (x + \mu)\theta = (y + \mu)\theta$. So $\theta$ is well defined.

Now $(x + \mu + y + \mu)\theta = (x + y + \mu)\theta = (x + y)\alpha_t = x\alpha_t + y\alpha_t = (x + \mu)\theta + (y + \mu)\theta$ and $((x + \mu)r)\theta = (xr + \mu)\theta = (xr)\alpha_t = (x\alpha_t)r = ((x + \mu)\theta)r$.

Therefore $\theta$ is an R-homomorphism. Again $(x + \mu)\theta = (y + \mu)\theta \Rightarrow x\alpha_t = y\alpha_t \Rightarrow (x, y) \in \alpha_t \Rightarrow \alpha(x, y) = t \Rightarrow \mu(x - y) = t \Rightarrow \mu(x - y) = \alpha(0, 0) = \mu(0)$ which implies that $x + \mu = y + \mu$. So $\theta$ is injective. Obviously $\theta$ is surjective. Hence the result.

**DEFINITION [56]:** *Let R be a near-ring and let $\mu$ be a fuzzy set in R. We say that $\mu$ is a fuzzy subnear-ring of R if for all $x, y \in R$.*

  i. $\mu(x - y) \geq \min\{\mu(x), \mu(y)\}$.
  ii. $\mu(xy) \geq \min\{\mu(x), \mu(y)\}$.

*If a fuzzy set $\mu$ in a near-ring R satisfies the property $\mu(x - y) \geq \min\{\mu(x), \mu(y)\}$ then letting $x = y$; $\mu(0) > \mu(x)$ for all $x \in R$.*

**DEFINITION [1]:** *Let $(R, +, \bullet)$ be a near-ring. A fuzzy set $\mu$ in R is called a fuzzy right (resp. left) R-subgroup of R if*

  i. $\mu$ *is a fuzzy subgroup of $(R, +)$.*
  ii. $\mu(xr) \geq \mu(x)$ *(resp. $\mu(rx) \geq \mu(x)$) for all $r, x \in R$.*

***Example [56]***: Let $R = \{a, b, c, d\}$ be a set with two binary operations as follows.

| + | a | b | c | d |
|---|---|---|---|---|
| a | a | b | c | d |
| b | b | a | d | c |
| c | c | d | b | a |
| d | d | c | a | b |

| $\bullet$ | a | b | c | d |
|---|---|---|---|---|
| a | a | a | a | a |
| b | a | a | a | a |
| c | a | a | a | a |
| d | a | a | b | b |

Then $(R, +, \bullet)$ is a near-ring. We define a fuzzy set $\mu: R \to [0,1]$ by $\mu(c) = \mu(d) < \mu(b) < \mu(a)$. Then $\mu$ is a fuzzy subgroup of $(R, +)$, and we have that $\mu(xr) \geq \mu(x)$ for all $r, x \in R$. Hence $\mu$ is a fuzzy right R-subgroup of R.

***Example [56]:*** Let $R = \{a, b, c, d\}$ be a set with two binary operations as follows:

| + | a | b | c | d |
|---|---|---|---|---|
| a | a | b | c | d |
| b | b | a | d | c |
| c | c | d | b | a |
| d | d | c | a | b |

| • | a | b | c | d |
|---|---|---|---|---|
| a | a | a | a | a |
| b | a | a | a | a |
| c | a | a | a | a |
| d | a | b | c | d |

Then (R, +, •) is a near-ring. We define a fuzzy set $\mu: R \to [0,1]$ by $\mu(c) = \mu(d) < \mu(b) < \mu(a)$. Thus $\mu$ is also a fuzzy right R subgroup of R.

**THEOREM [56]:** *If $\mu$ is a fuzzy right (resp. left) R-subgroup of a near-ring R, then the set $R_\mu = \{x \in R \mid \mu(x) = \mu(\theta)\}$ is a right (resp. left) R-subgroup of R.*

*Proof*: Let x, y ∈ $R_\mu$. Then $\mu(x) = \mu(y) = \mu(0)$. Since $\mu$ is a fuzzy right (resp. left) R-subgroup, it follows that

$$\mu(x - y) \geq \min\{\mu(x), \mu(y)\}$$
$$= \min\{\mu(0), \mu(0)\}$$
$$= \mu(0).$$

On the other hand $\mu(x - y) \leq \mu(0)$. Hence we have $\mu(x - y) = \mu(0)$ so $x - y \in R_\mu$ Also for any $x \in R_\mu$ and $r \in R$, we get $\mu(xr) \geq \mu(x) = \mu(0)$ (resp. $\mu(rx) \geq \mu(x) = \mu(0)$).

On the other hand $\mu(xr) \leq \mu(0)$ (resp. $\mu(rx) \leq \mu(0)$). Hence we obtain $\mu(xr) = \mu(0)$ (resp $\mu(rx) = \mu(0)$), which shows that $x r \in R_\mu$ (resp. $r x \in R_\mu$) Consequently the set $R_\mu$ is a right (resp. left) R-subgroup of R.

**DEFINITION 1.8.6:** *A fuzzy right (resp. left) R-subgroup $\mu$ of a near-ring R is said to be normal if there exists $x \in R$ such that $\mu(x) = 1$.*

Note that if a fuzzy right (resp. left) R-subgroup $\mu$ of a near-ring R is normal then $\mu(0) = 1$; hence $\mu$ is a normal fuzzy right (resp. left) R-subgroup of a near-ring R if and only if $\mu(0) = 1$.

**THEOREM 1.8.12:** *Let $\mu$ be a fuzzy right (resp. left) R-subgroup of a near-ring R and let $\mu^+$ be a fuzzy set in R defined by $\mu^+(x) = \mu(x) + 1 - \mu(0)$ for all $x \in R$. Then $\mu^+$ is a normal fuzzy right (resp. left) R-subgroup of R containing $\mu$.*

*Proof*: Let x, y ∈ R. Then

$$\min\{\mu^+(x), \mu^+(y)\} = \min\{\mu(x) + 1 - \mu(0), \mu(y) + 1 - \mu(0)\}$$
$$= \min\{\mu(x), \mu(y)\} + 1 - \mu(0)$$
$$\leq \mu(x - y) + 1 - \mu(0)$$
$$= \mu^+(x - y)$$

and for all x, y ∈ R we have

$$\mu^+(xr) = \mu(xr) + 1 - \mu(0)$$
$$\geq \mu(x) + 1 - \mu(0)$$
$$= \mu^+(x).$$

Similarly $\mu^+(rx) = \mu(rx) + 1 - \mu(0) \geq \mu(x) + 1 - \mu(0) = \mu^+(x)$. Hence $\mu^+$ is a fuzzy right (resp. left) R-subgroup of R. Clearly $\mu^+(0) = 1$ and $\mu \subset \mu^+$. This gives the proof.

Using the fact $\mu \subset \mu^+$ we have the following theorem:

**THEOREM 1.8.13**: *If $\mu$ is a fuzzy right (resp. left) R-subgroup of R satisfying $\mu^+(x) = 0$ for some $x \in R$, then $\mu(x) = 0$ also.*

**THEOREM [1]:** *Let $(R, +, \cdot)$ be a near-ring and $\chi_H$ be the characteristic function of a subset $H \subset R$. Then H is a right (resp. left) R-subgroup of R if and only if $\chi_H$ is a fuzzy right (resp. left) R-subgroup of R.*

*Proof*: Refer [1].

**THEOREM 1.8.14:** *For any right (resp. left) R-subgroup H of a near-ring R, the characteristic function $\chi_H$ of H is a normal fuzzy right (resp. left) R-subgroup of R and $R_{\chi_H} = H$.*

*Proof*: From the above theorem the proof of this theorem can be easily obtained.

**THEOREM 1.8.15:** *Let $\mu$ and $\nu$ be fuzzy right (resp. left) R-subgroups of a near-ring R. If $\mu \subset \nu$ and $\mu(0) = \nu(0)$ then $R_\mu \subset R_\nu$.*

*Proof*: Assume that $\mu \subset \nu$ and $\mu(0) = \nu(0)$. If $x \in R_\mu$ then $\nu(x) \geq \mu(x) = \mu(0) = \nu(0)$. Noting that $\nu(x) \leq \nu(0)$ for all $x \in R$ we have $\nu(x) = \nu(0)$ that is $x \in R$. Hence the proof.

The following result is a direct consequence of the above theorem; hence left for the reader.

**THEOREM 1.8.16:** *If $\mu$ and $\nu$ are normal fuzzy right (resp. left) R-subgroups of a near-ring R satisfying $\mu \subset \nu$ then $R_\mu \subset R_\nu$.*

**THEOREM 1.8.17:** *A fuzzy right (resp. left) R-subgroup $\mu$ of a near-ring R is normal if and only if $\mu^+ = \mu$.*

*Proof*: Sufficiency is direct. To prove the necessary condition assume $\mu$ is a normal fuzzy right (resp. left) R-subgroup of R and let $x \in R$. Then $\mu^+(x) = \mu(x) + 1 - \mu(0) = \mu(x)$ and hence $\mu^+ = \mu$.

**THEOREM 1.8.18:** *If $\mu$ is a fuzzy right (resp. left) R-subgroup of a near-ring then $(\mu^+)^+ = \mu^+$.*

*Proof*: For any $x \in R$, we have $(\mu^+)^+(x) = \mu^+(x) + 1 - \mu^+(0) = \mu^+(x)$ completing the proof.

The following result is direct by the above definition.

**THEOREM 1.8.19:** *If $\mu$ is a normal fuzzy right (resp. left) R-subgroup of a near-ring R, then $(\mu^+)^+ = \mu$.*

**THEOREM 1.8.20:** *Let $\mu$ be a fuzzy right (resp. left) R-subgroup of near-ring R. If there exists a fuzzy right (resp. left) R-subgroup $\nu$ of R satisfying $\nu^+ \subset \mu$ then $\mu$ is normal.*

*Proof*: Suppose there exists a fuzzy right (resp. left) R-subgroup $\nu$ of R such that $\nu^+ \subset \mu$. Then $1 = \nu^+(0) \leq \mu(0)$ whence $\mu(0) = 1$. Hence the proof.

Using the above theorem the following result is straightforward.

**THEOREM 1.8.21:** *Let $\mu$ be a fuzzy right (resp. left) R-subgroup of a near-ring R. If there exists a fuzzy right (resp. left) R-subgroup $\nu$ of R satisfying $\nu^+ \subset \mu$ then $\mu^+ = \mu$.*

**THEOREM 1.8.22:** *Let $\mu$ be a fuzzy right (resp. left) R-subgroup of a near-ring R and let $f: [0, \mu(0)] \to [0,1]$ be an increasing function. Define a fuzzy set $\mu_f: R \to [0,1]$ by $\mu_f(x) = f(\mu(x))$ for all $x \in R$. Then $\mu_f$ is a fuzzy right (resp. left) R-subgroup of R. In particular if $f(\mu(0)) = 1$ then $\mu_f$ is normal and if $f(t) \geq t$ for all $t \in [0, \mu(0)]$ then $\mu \subseteq \mu_f$.*

*Proof:* Let $x, y \in R$. Then $\mu_f(x - y) = f(\mu(x - y)) \geq f(\min \mu(x), \mu(y)) = \min \{f(\mu(x)), f(\mu(y))\} = \min\{\mu_f(x), \mu_f(y)\}$ and for all $x, r \in R$, we have $\mu_f(xr) = f(\mu(xr)) \geq f(\mu(x)) = \mu_f(x)$. Similarly $\mu_f(rx) = f(\mu(rx)) \geq f(\mu(x)) = \mu_f(x)$. Hence $\mu_f$ is a fuzzy right (resp. left) R-subgroup of R. If $f(\mu(0)) = 1$ then clearly $\mu_f$ is normal. Assume that $f(t) \geq t$ for all $t \in [0, \mu(0)]$. Then $\mu_f(x) = f(\mu(x)) \geq \mu(x)$ for all $x \in R$, which proves that $\mu \subseteq \mu_f$.

**THEOREM 1.8.23:** *Let $\mu$ be a non constant normal fuzzy right (resp. left) R-subgroup of R, which is maximal in the poset of normal fuzzy right (resp. left) R-subgroups under set inclusion. Then $\mu$ takes only the values 0 and 1.*

*Proof*: Please refer [56].

Another normalization of fuzzy right (resp. left) R-subgroup of a near-ring R as given by [56] is given below.

**THEOREM 1.8.24:** *Let $\mu$ be a fuzzy right (resp. left) R-subgroup of a near-ring R and let $\mu^0$ be a fuzzy set in R defined by $\mu^0(x) = \mu(x) / \mu(0)$ for all $x \in R$. Then $\mu^0$ is a normal fuzzy right (resp.left) R-subgroup of R containing $\mu$.*

*Proof*: For any $x, y \in R$ we have

$$\min \{\mu^0(x), \mu^0(y)\} = \min \left\{\frac{\mu(x)}{\mu(0)}, \frac{\mu(y)}{\mu(0)}\right\}$$

$$= \frac{1}{\mu(0)} \min\{\mu(x), \mu(y)\}$$

$$\leq \frac{1}{\mu(0)} \mu(x-y)$$

$$= \mu^0(x-y)$$

and for all x, y ∈ R we get

$$\mu^0(xr) = \frac{1}{\mu(0)} \mu(xr)$$

$$\geq \frac{1}{\mu(0)} \mu(x)$$

$$= \mu^0(x).$$

Similarly $\mu^0(rx) \geq \mu^0(x)$. Hence $\mu^0$ is a fuzzy right (resp. left) R subgroup of R. Clearly $\mu^0(0) = 1$ and $\mu \subseteq \mu^0$. Hence the claim. Using the fact $\mu \subset \mu^0$ the following results are straightforward.

**THEOREM 1.8.25:** *If $\mu$ is a fuzzy right (resp. left) R-subgroup of a near-ring R satisfying $\mu^0(0) = 0$ for some x ∈ R then $\mu(x) = 0$ also.*

**THEOREM 1.8.26:** *Let H be a right (resp. left) R-subgroup of a near-ring R and let $\mu_H$ be a fuzzy set in R defined by*

$$\mu_H(x) = \begin{cases} 1 & \text{if } x \in H \\ 0 & \text{otherwise} \end{cases}$$

*Then $\mu_H$ is a normal fuzzy right (resp. left) R subgroup of R and $R_{\mu H} = H$.*

*Proof*: Left for the reader as the proof is direct.

**THEOREM 1.8.27:** *Let $\mu$ be a non constant fuzzy right (resp. left) R-subgroup of a near-ring R such that $\mu^+$ is maximal in the poset of normal fuzzy right (resp. left) R-subgroups under set inclusion. Then*

 i. *$\mu$ is normal.*
 ii. *$\mu$ takes only the values 0 and 1.*
 iii. *$\mu_{R_\mu} = \mu$.*
 iv. *$R_\mu$ is a maximal right (resp. left) R-subgroup of R.*

*Proof*: Please refer [56].

We have given several of the proofs verbatim from [56] mainly to make it easy when we do the proofs in case of Smarandache fuzzy near-rings and their ideals and subgroups. For more about fuzzy near-rings please refer [56, 57].

Now we proceed on to define the notion of fuzzy near-rings.

**DEFINITION 1.8.7**: *Let $\mu$ be a non-empty fuzzy subset of a near-ring N (that is $\mu(x) \neq 0$ for some $x \in N$) then $\mu$ is said to be a fuzzy ideal of N if it satisfies the following conditions:*

    i. $\mu(x+y) \geq \min\{\mu(x), \mu(y)\}$.
    ii. $\mu(-x) = \mu(x)$.
    iii. $\mu(x) = \mu(y + x - y)$.
    iv. $\mu(xy) \geq \mu(x)$ and
    v. $\mu\{x(y+i) - xy\} \geq \mu(i)$ for all $x, y, i \in N$.

1. *If $\mu$ is a fuzzy ideal of N then $\mu(x+y) = \mu(y+x)$.*
2. *If $\mu$ is a fuzzy ideal of N then $\mu(0) \geq \mu(x)$ for all $x \in N$. The above two statements can be easily verified for if we put $z = x + y$, then $\mu(x+y) = \mu(z) = \mu(-x+z+x)$ since $\mu$ is a fuzzy ideal $\mu(-x+x+y+x) = \mu(y+x)$ (since $z = x + y$).*

*Likewise for the second statement $\mu(0) = \mu(x-x) \geq \min\{\mu(x), \mu(-x)\} = \mu(x)$ since $\mu$ is a fuzzy ideal (since $\mu(-x) = \mu(x)$ by the very definition of fuzzy ideal).*

**DEFINITION 1.8.8**: *Let I be an ideal of N. we define $\lambda_I : N \to [0, 1]$ as*

$$\lambda_I(x) = \begin{cases} 1 & \text{if } x = 1 \\ 0 & \text{otherwise} \end{cases}$$

*$\lambda_I(x)$ is called the characterize function on 1.*

**THEOREM [130]**: *Let N be a near-ring and $\lambda_I$, the characteristic function on a subset I of N. Then $\lambda_I$ is a fuzzy ideal of N if and only if I is an ideal of N.*

**DEFINITION [130]**: *Let $\mu$ be a fuzzy subset of X. Then the set $\mu_t$ of all $t \in [0, 1]$ is defined by $\mu_t = \{x \in N \,/\, \mu(x) \geq t\}$ is called the level subset of t for the near-ring N.*

**DEFINITION [130]**: *Let N be a near-ring and $\mu$ be a fuzzy ideal of N. Then the level subset $\mu_t$ of N for all $t \in [0, t]$, $t \leq \mu(0)$ is an ideal of N if and only if $\mu$ is a fuzzy ideal of N.*

**DEFINITION [130]**: *A fuzzy ideal $\mu$ of N is called fuzzy prime if for any two fuzzy ideals $\sigma$ and $\theta$ of N $\sigma \circ \theta \subseteq \mu$ implies $\sigma \subseteq \mu$ or $\theta \subseteq \mu$.*

Now we use the concept of fuzziness in $\Gamma$-near-rings as given by [1, 130].

**THEOREM [130]**: *If $\mu$ is a fuzzy ideal of a near-ring N and $a \in N$ then $\mu(x) \geq \mu(a)$ for all $x \in \langle a \rangle$.*

*Proof*: Please refer [1, 130].

**DEFINITION [130]:** *Let $\mu$ and $\sigma$ be two fuzzy subsets of M. Then the product of fuzzy subset ($\sigma o \tau$)(x) = sup {min ($\sigma$ (y), $\tau$ (z))} if x is expressible as a product x = yz where y, z $\in$ M and ($\sigma o \tau$) (x) = 0 otherwise.*

**DEFINITION [130]:** *A fuzzy ideal $\mu$ of N is said to have fuzzy IFP if $\mu$ (a n b) $\geq \mu$(ab) for all a, b, n $\in$ N.*

**DEFINITION [130]:** *Let $\mu$ be a fuzzy ideal of N, $\mu$ has fuzzy IFP if and only if $\mu_k$ is a IFP ideal of N for all $0 \leq k \leq 1$.*

**DEFINITION [130]:** *N has strong IFP if and only if every fuzzy ideal of N has fuzzy IFP.*

**THEOREM [130]:** *If $\mu$ is a fuzzy IFP-ideal of N then $N_\mu = \{x \in N / \mu (x) = \mu(0)\}$ is an IFP ideal of N.*

*Proof*: We have $\mu (0) \geq \mu (x)$ for all $x \in N$; write $t = \mu (0)$. Now $N_\mu = \mu_t$ and by definitions $\mu_t$ has IFP. Therefore $N_\mu$ is an IFP ideal of N.

**Notation:** Let $\mu$ be a fuzzy ideal of N. For any $s \in [0, 1]$ define $\beta_s : N \to [0, 1]$ by

$$_\mu\beta_s(x) = \begin{cases} s \text{ if } \mu(x) \geq s \\ \mu(x) \text{ if } \mu(x) < s. \end{cases}$$

Since $\beta_s$ depends on $\mu$, we also denote $\beta_s$ by $_\mu\beta_s$.

*Results:*

1. $\beta_s (x) \leq s$ for all $x \in N$.
2. $\beta_s$ is a fuzzy ideal of N.
3. If $\mu (0) = t$, then $s \geq t$ if and only if $\mu_t = (\beta_s)_t$.

The proof of the above 3 statements are left as an exercise for the reader to prove.

We can still equivalently define as a definition or prove it as a theorem.

**DEFINITION [130]:** *$\mu$ is a fuzzy IFP ideal of N if and only if $\beta_s$ is a fuzzy IFP ideal for all $s \in [0, 1]$.*

*Result*: A fuzzy ideal $\mu$ has IFP if and only if $N_{(\beta s)}$ has IFP for all $s \in [0, 1]$.

The above result is assigned as an exercise for the reader to prove.

**DEFINITION [130]:** *Let $\mu : M \to [0, 1]$. $\mu$ is said to be a fuzzy ideal of M if it satisfies the following conditions:*

    i. $\mu(x+y) \geq \min\{\mu(x), \mu(y)\}$.
    ii. $\mu(-x) = \mu(x)$.
    iii. $\mu(x) = \mu(y + x - y)$.
    iv. $\mu(x\,\alpha\,y) \geq \mu(x)$ and
    v. $\mu\{(x\,\alpha\,(y+z) - x\,\alpha\,y\} \geq \mu(z)$

for all $x, y, z \in M$ and $\alpha \in \Gamma$.

The following result is left as an exercise for the reader to prove.

**THEOREM 1.8.28**: *Let $\mu$ be a fuzzy subset of M. Then the level subsets $\mu_t = \{x \in M\,/\,\mu(x) \geq t\}$, $t \in \text{Im } \mu$, are ideals of M if and only if $\mu$ is a fuzzy ideal of M.*

All results true in case of fuzzy ideals of N are true in case of fuzzy ideals of M with some minor modifications.

The following result can be proved by routine application of definitions.

**THEOREM 1.8.29**: *Let M and M' be two $\Gamma$-near-rings, $h : M \rightarrow M'$ be an $\Gamma$-epimorphism and $\mu$, $\sigma$ be fuzzy ideals of M and M' respectively then*

    i. $h(h^{-1}(\sigma)) = \sigma$.
    ii. $h^{-1}(h(\mu)) \supseteq \mu$ and
    iii. $h^{-1}(h(\mu)) = \mu$ if $\mu$ is constant on ker h.

**DEFINITION [130]**: *A fuzzy ideal $\mu$ of M is said to be a fuzzy prime ideal of M if $\mu$ is not a constant function; and for any two fuzzy ideals $\sigma$ and $\Gamma$ of M, $\sigma \circ \Gamma \subseteq \mu$, implies either $\sigma \subset \mu$ or $\Gamma \subset \mu$.*

Using these definitions it can be proved.

**THEOREM [130]**: *If $\mu$ is a fuzzy prime ideal of M then $M_\mu = \{x \in M\,/\,\mu(x) = \mu(0)\}$ is a prime ideal of M.*

**PROPOSITION [130]**: *Let I be an ideal of M and $s \in [0, 1)$. Let $\mu$ be a fuzzy subset of M, defined by*

$$\mu(x) = \begin{cases} 1 & \text{if } x \in I \\ s & \text{otherwise.} \end{cases}$$

*Then $\mu$ is a fuzzy prime ideal of M if I is a prime ideal of M.*

*Proof*: Using the fact $\mu$ is a non-constant fuzzy ideal of M we can prove the result as a matter of routine using the basic definitions.

**DEFINITION [130]**: *Let I be an ideal of M. Then $\lambda_I$ is a fuzzy prime ideal of M if and only if I is a prime ideal of M.*

***Result 1***: If μ is a fuzzy prime ideal of M then μ(0) = 1.

*Proof*: It is left for the reader to prove.

***Result 2***: If μ is a fuzzy prime ideal of M then |Im μ| = 2.

We further define fuzzy near matrix ring.

**DEFINITION 1.8.9**: *Let $P_{n \times n}$ denote the set of all $n \times n$ matrices with entries form [0, 1] i.e. $P_{n \times n} = \{(a_{ij}) / a_{ij} \in [0, 1]\}$ for any two matrices $A, B \in P_{n \times n}$ define $\oplus$ as follows:*

$$A = \begin{pmatrix} a_{11} & a_{12} & \ldots & a_{1n} \\ a_{21} & a_{22} & \ldots & a_{2n} \\ \vdots & \vdots & & \vdots \\ a_{n1} & a_{n2} & \ldots & a_{nn} \end{pmatrix}$$

*and*

$$B = \begin{pmatrix} b_{11} & b_{12} & \ldots & b_{1n} \\ b_{21} & b_{22} & \ldots & b_{2n} \\ \vdots & \vdots & & \vdots \\ b_{n1} & b_{n2} & \ldots & b_{nn} \end{pmatrix}$$

$$A \oplus B = \begin{pmatrix} a_{11} + b_{11} & \ldots & a_{1n} + b_{1n} \\ a_{21} + b_{21} & \ldots & a_{21} + b_{2n} \\ \vdots & & \vdots \\ a_{n1} + b_{n1} & \ldots & a_{nn} + b_{nn} \end{pmatrix}$$

*where*

$$a_{ij} + b_{ij} = \begin{cases} a_{ij} + b_{ij} & \text{if } a_{ij} + b_{ij} < 1 \\ 0 & \text{if } a_{ij} + b_{ij} = 1 \\ a_{ij} + b_{ij} - 1 & \text{if } a_{ij} + b_{ij} > 1. \end{cases}$$

*Clearly $(P_{n \times n}, \oplus)$ is an abelian group and*

$$0 = \begin{pmatrix} 0 & 0 & \ldots & 0 \\ 0 & 0 & \ldots & 0 \\ \vdots & \vdots & & \vdots \\ 0 & 0 & \ldots & 0 \end{pmatrix}$$

*is zero matrix which acts as the additive identity with respect to $\oplus$.*

*Define $\odot$ on $P_{n \times n}$ as follows. For $A, B \in P_{n \times n}$*

$$A \odot B = \begin{pmatrix} a_{11} & \ldots & a_{1n} \\ a_{21} & \ldots & a_{2n} \\ \vdots & & \vdots \\ a_{n1} & \ldots & a_{nn} \end{pmatrix} \odot \begin{pmatrix} b_{11} & \ldots & b_{1n} \\ b_{21} & \ldots & b_{2n} \\ \vdots & & \vdots \\ b_{2n} & \ldots & b_{nn} \end{pmatrix} = \begin{pmatrix} a_{11}+\ldots+a_{1n} & \ldots & a_{11}+\ldots+a_{1n} \\ a_{21}+\ldots+a_{2n} & \ldots & a_{21}+\ldots+a_{2n} \\ \vdots & & \vdots \\ a_{n1}+\ldots+a_{nn} & \ldots & a_{n1}+\ldots+a_{nn} \end{pmatrix}$$

*where $a_{ij}.b_{ij} = a_{ij}$ for all $a_{ij} \in A$ and $b_{ij} \in B$. Clearly $(P_{n \times n}, \odot)$ is a semigroup. Thus $(A \oplus B) \odot C = A \odot C \oplus B \odot C$. Hence $(P_{n \times n}, \oplus, \odot)$ is a near-ring, which we call as the fuzzy near matrix ring or fuzzy matrix near-ring.*

**THEOREM 1.8.30**: *The fuzzy near- matrix ring is a commutative near-ring.*

*Proof*: Straightforward.

**THEOREM 1.8.31**: *The fuzzy matrix near-ring is not an abelian near-ring.*

*Proof*: For $A, B \in P_{n \times n}$ we have $A \odot B \neq B \odot A$ in general.

**THEOREM 1.8.32:** *In $\{P_{n \times n}, \oplus, \odot\}$ we have $I_{n \times n} \neq A$ where $I_{n \times n}$ is the matrix with diagonal elements 1 and rest 0.*

*Proof*: Left for the reader to prove.

**DEFINITION 1.8.10**: *Let $\{P_{n \times n}, \oplus, \odot\}$ be a fuzzy near matrix ring we say a subset I of $P_{nxn}$ is a fuzzy left ideal of $P_{nxn}$ if*

    i.    *$(I, +)$ is a normal subgroup of $P_{nxn}$.*
    ii.    *$n(n' + i) + n_r n' \in I$ for each $i \in I$ and $n_r, n, n' \in N$ where $n_r$ denotes the unique right inverse of n.*

*All properties enjoyed by near-rings can be defined and will be true with appropriate modifications.*

Next we proceed on to define the concept of fuzzy complex near-rings.

**DEFINITION 1.8.11:** *Let $V = \{a + ib \,/\, a, b \in [0,1]\}$ define on V the operation called addition denoted by $\oplus$ as follows:*

*For $a + ib, a_1 + ib_1 \in V$, $a + ib \oplus a_1 + ib_1 = a + a_1 + i(b + b_1)$ where $a \oplus a_1 = a + a_1$ if $a + a_1 < 1$ and $a + a_1 = a + a_1 - 1$ if $a + a_1 \geq 1$ where '+' is the usual addition of numbers. Clearly $(V, \oplus)$ is a group. Define $\odot$ on V by $(a + ib) \odot (a_1 + ib_1) = a + ib$ for all $a + ib, a_1 + ib_1 \in V$. $(V, \odot)$ is a semigroup. It is easily verified. $(V, \oplus, \odot)$ is a near-ring, which we call as the fuzzy complex near-ring.*

*Further P = {a / a ∈ [0,1]} and C = {ib / b ∈ [0 1]} are fuzzy complex subnear-rings of (V, ⊕, ⊙).*

**PROPOSITION 1.8.1**: *V has non-trivial idempotent.*

*Proof*: Left for the reader to prove.

**THEOREM 1.8.33**: *Let {V, ⊕, ⊙} be a fuzzy complex near-ring. Every nontrivial fuzzy subgroup of N is a fuzzy right ideal of V.*

*Proof*: Obvious by the fact that if N is a fuzzy subgroup of V then NV ⊆ N.

It is an open question. Does V have nontrivial fuzzy left ideals and ideals? The reader is requested to develop new and analogous notions and definitions about these concepts.

Now a natural question would be can we have the concept of fuzzy non-associative complex near-ring; to this end we define a fuzzy non-associative complex near-ring.

**DEFINITION 1.8.12:** *Let W = {a + ib / a, b ∈ [0, 1] called the set of fuzzy complex numbers}. Define on W two binary operations ⊕ and ⊙ as follows:*

*(W, ⊕) is a commutative loop where for a + ib, c + id ∈ W define a + ib ⊕ c + id = a ~ c + i ( b ~ d) where '~' is the difference between a and b. Clearly (W, ⊕) is a commutative loop.*

*Define ⊙ on W by (a + ib) ⊙ (c + id) = a + ib for all a + ib, c + id ∈ W. (W, ⊕, ⊙) is called the fuzzy complex non-associative near-ring. ([0, 1], ⊕, ⊙) ⊆ (W, ⊕, ⊙) is a fuzzy non-associative subnear-ring.*

Obtain interesting properties about these non-associative fuzzy complex near-rings.

Now we proceed on to define fuzzy polynomial near-rings.

**DEFINITION 1.8.13**: *Let R be the set of reals. The fuzzy polynomial near-ring $R[x^{[0, 1]}]$ consists of elements of the form $p_0 + p_1 x^{\gamma_1} + p_2 x^{\gamma_2} + ... + p_n x^{\gamma_n}$ where $p_0, p_1, ... p_m \in R$ and $\gamma_1, \gamma_2, ..., \gamma_n \in [0, 1]$ with $\gamma_1 < \gamma_2 < ... < \gamma_n$. Two elements $p(x) = q(x) \Leftrightarrow p_i = q_i$ and $\gamma_i = s_i$ where $p(x) = p_0 + p_1 x^{\gamma_1} + ... + p_n x^{\gamma_n}$ and $q(x) = q_0 + q_1 x^{s_1} + ... + q_n x^{s_n}$ Addition is performed as in the case of usual polynomials.*

*Define ⊙ on $R[x^{[0 1]}]$ by $p(x) \odot q(x) = p(x)$ for $p(x), q(x) \in R[x^{[0, 1]}]$.*

*Clearly $\{R[x^{[0, 1]}], +, \odot\}$ is called the fuzzy right polynomial near-ring. $x^0 = 1$ by definition.*

**DEFINITION 1.8.14**: Let $\{R[x^{[0,1]}], +, \odot\}$ be a fuzzy polynomial near-ring. For any polynomial $p(x) \in R[x^{[0,1]}]$ define the derivative of $p(x)$ as follows.

If $p(x) = p_0 + p_1 x^{s_1} + \ldots + p_n x^{s_n}$

$$\frac{dp(x)}{dx} = 0 + s_1 p_1 x^{s_1 \sim 1} + \ldots + s_n p_n x^{s_n \sim 1} = (s_1 p_1) x^{s_1 \sim 1} + \ldots + (s_n p_n) x^{s_n \sim 1}$$

where '~' denotes the difference between $s_i$ and 1.

Clearly if $p(x) \in R[x^{[0,1]}]$ then $\frac{dp(x)}{dx} \in R[x^{[0,1]}]$ Likewise successive derivatives are also defined i.e. product of $s_i p_i \in R$ as $s_i \in [0, 1]$ and $p_i \in R$ i.e. the usual multiplication of the reals.

***Example 1.8.1***: Let R be reals $R[x^{[0,1]}]$ be a polynomial near-ring.

$$p(x) = 5 - 6x^{1/5} + 2x^{3/8} - 15x^{7/9}$$
$$\frac{dp(x)}{dx} = 0 - \frac{1}{5} 6 x^{4/5} + \frac{2 \times 3}{8} x^{5/8} - \frac{15 \times 7}{9} x^{2/9} = -\frac{6}{5} x^{4/5} + \frac{3}{4} x^{5/8} - \frac{35}{3} x^{2/9}.$$

The observation to be made is that no polynomial other than the polynomial x vanishes after differentiation.

**DEFINITION 1.8.15**: Let $p(x) \in \{Rx^{[0,1]}\}$ the fuzzy degree of $p(x)$ is $s_n$ where $p(x) = p_0 + p_1 x^{s_1} + \ldots + p_n x^{s_n}$; $s_1 < s_2 \ldots < s_n$ ($p_n \neq 0$) deg $p(x) = s_n$. The maximal degree of any polynomial $p(x)$ can take is 1. Now it is important to note that as in the case polynomial rings we cannot say deg $[p(x).q(x)]$ = deg $p(x)$ + deg $q(x)$.

But we have always in fuzzy polynomial near-ring.

$$\text{deg } (p(x) q(x)) = \text{deg } p(x) \text{ for}$$
$$p(x), q(x) \in R[x^{[0,1]}]$$

as this degree for fuzzy polynomial near-rings is a fuzzy degree we shall denote them by f(deg (p(x)).

**DEFINITION 1.8.16**: Let $p(x) \in [R[x^{[0,1]}]$, $p(x)$ is said to have a root $\alpha$ if $p(\alpha) = 0$.

***Example 1.8.2***: Let $p(x) = \sqrt{2} - x^{1/2}$ be a fuzzy polynomial in $R\{x^{[0,1]}\}$. The root of $p(x)$ is 2 for $p(2) = \sqrt{2} - 2^{1/2} = 0$.

But as in case of root of polynomial of degree n has n and only n roots which is the fundamental theorem on algebra; we in case of fuzzy polynomial near-rings cannot say the number of roots in a nice mathematical terminology that is itself fuzzy.

A study of these fuzzy polynomial near-rings is left open for any interested researcher. We proceed on to define fuzzy polynomial near-rings when the number of variables is more than one x and y.

**DEFINITION 1.8.17**: *Let R be the reals x, y be two variables we first assume xy = yx. Define the fuzzy polynomial near-ring.*

$$R\,[x^{[0,\,1]},\,y^{[0,\,1]}]\;(by)\;=\;\left\{\sum r_i x^{p_i} y^{q_i}\,/\,r_i \in R;\;p_i \in [0,1]\;\;q_i \in [0,1]\right\}$$

*Define '+' as in the case of polynomial and '⊙' by $p(xy)\odot q(x, y) = p(x, y)$. Clearly $R[x^{[0,1]}, y^{[0,1]}]$ is called as a fuzzy polynomial right near-ring.*

**DEFINITION 1.8.18**: *Let $\{R\,[x^{[0,\,1]},\,y^{[0,\,1]}],\,\oplus,\,\bullet\}$ be a fuzzy polynomial near-ring in the variable x and y.*

*A fuzzy polynomial $p(x, y)$ is said to be homogenous of fuzzy degree t, $t \in [0, 1]$ if $p(x, y) = a_n x^{s_1} y^{t_1} +\ldots+ b_n x^{s_p t_p}$ then $t_i \neq 0$, $s_i \neq 0$ for all $i = 1, 2, \ldots, p$ and $s_i + t_i = t$ for $i = 1, 2, \ldots, p$.*

**DEFINITION 1.8.19**: *Let $R\,\{x^{[0,\,1]},\,y^{[0,\,1]},\,\oplus,\,{}^{\backprime}\bullet{}^{\prime}\}$ be the fuzzy polynomial near-ring in the variables x and y. A symmetric fuzzy polynomial is a homogenous polynomial of fuzzy degree t, $t \in [0, 1]$ such that $p(x, y) = p_1 x^{t_1} y^{s_1} +\ldots+ p_n x^{t_n} y^{s_n}$ where $t_1 < t_2 < \ldots < t_n$, $s_n < s_{n-1} < \ldots < s_1$ with $t_1 = s_n$, $t_2 = s_{n-1}, \ldots, t_n = s_1$ further $p_1 = p_n$, $p_2 = p_{n-1}, \ldots$*

For example $p(x, y) = 3x^{1/2} y^{2/3} + 3x^{2/3} y^{1/2}$. $p(x, y) = x^r + y^r$, $r \in [0, 1]$.

$p(x, y) = x^r + y^r + x^s y^t + y^t s^s$ where $s + t = 1$. $s, t, r \in [0, 1]$. We have like other polynomials we can extend the fuzzy polynomials to any number of variables say $X_1, X_2, \ldots, X_n$. under the assumption $X_i X_j = X_j X_i$ and denote it by $R\,[X_1^{[0,\,1]},\,X_2^{[0,\,1]},\,\ldots X_n^{[0,\,1]}]$ called the fuzzy polynomial near-ring in n variables. The reader is advised to develop new results on these fuzzy polynomial near-rings.

We have introduced the concept of complex near-ring and the non-associative complex near-ring now we just define yet another new notion called fuzzy non-associative near-ring.

**DEFINITION 1.8.20**: *Let $\{W, \oplus, \odot\}$ be the fuzzy non-associative complex near-ring. Let x be an indeterminate. We define the fuzzy non-associative polynomial near-ring as follows:*

*$W\,[x] = \{\Sigma p_i\, x^i\,/\,p_i \in W\}$; we say $p(x)$, $q(x) \in W[x]$ are equal if and only if every coefficient of same power of x is equal i.e. if $p(x) = p_0 + p_1 x +\ldots+ p_n x^n$ and $q(x) = q_0 + q_1 x +\ldots+ q_n x^n$. $p(x) = q(x)$ if and only if $p_i = q_i$. for $i = 1, 2, \ldots, n$. Addition is performed as follows $p(x) \oplus q(x) = p_0 \oplus q_0 +\ldots+ (p_n \oplus q_n)\,x^n$ where $\oplus$ is the operation on W. For $p(x)$, $q(x)$ in $W[x]$ define $p(x) \odot q(x) = p(x)$. Clearly $\{W(x), \oplus, \odot\}$ is a fuzzy non-associative complex polynomial near-ring.*

*Now take $Z^0 = Z^+ \cup \{0\}$. Let $p: Z^0 \to W$ be defined by $p(0) = 0$, $p(x) = \dfrac{1}{x}$ for $0 \neq x \in Z$.*

*Clearly $p(x)$ is a fuzzy non-associative subnear-ring of W. Thus p is a fuzzy non-associative subnear-ring. Let $G = Z^0 \times Z^0$. Define a map $p_0 : G \to W[x]$ by*

$$p(0, 0) = 0.$$
$$p(x, y) = \frac{1}{x} + \frac{1}{y}, \ x \neq 0, y \neq 0.$$
$$p(x, 0) = \frac{1}{x};$$
$$p(0, y) = \frac{1}{y}.$$

*Then the map p is a fuzzy non-associative complex subnear-ring of G.*

Several interesting research in this direction is thrown open for the reader.

Now we proceed on to define a special class of fuzzy near-ring.

**DEFINITION 1.8.21**: *Let $P = [0, 1]$ the interval from 0 to 1. Define $\oplus$ and $\odot$ on P as follows. For $a, b \in P$ define $a \oplus b = a + b$ if $a + b < 1$, $a \oplus b = 0$ if $a + b = 1$ and $a \oplus b = a + b - 1$ if $a + b > 1$. Thus $\oplus$ acts as modulo 1. Define $\odot$ on $a, b \in P = [0, 1]$ by $a \odot b = a$; clearly $(a \oplus b) \odot c = a \odot c + b \odot c = a \oplus b$. Clearly $(P, \oplus)$ is a group and $(P, \odot)$ is a semigroup. Hence $(P, \oplus, \odot)$ is a right near-ring. We call $\{P, \oplus, \odot\}$ the special fuzzy right near-ring.*

**DEFINITION 1.8.22**: *Let $(P, \oplus, \odot)$ be a fuzzy near-ring. $P_0 = \{p \in P \,/\, p.0 = 0\}$ is called the fuzzy zero symmetric part and $P_c = \{n \in P \,/\, n.0 = n\}$ is called the fuzzy constant part.*

**THEOREM 1.8.34**: *The special fuzzy right near-ring $\{P, \oplus, \odot\}$ has no fuzzy invertible elements.*

*Proof*: Left for the reader to prove.

Let $S = \{r/p, 0 \,/\, 1 < r < p\}$ is a fuzzy subnear-ring or to be more specific if $S = \{0, ¼, ½, ¾\}$; S is a fuzzy subnear-ring.

**DEFINITION 1.8.23**: *A fuzzy subnear-ring N of P is called fuzzy invariant if $NP \subseteq N$ and $PN \subseteq N$ we call a fuzzy subnear-ring N of P to be a fuzzy right invariant if $NP \subset N$.*

The following theorem is left as an exercise to the reader.

**THEOREM 1.8.35**: *Every fuzzy subnear-ring N of P is fuzzy right invariant.*

**DEFINITION 1.8.24**: *The set P = {0, 1} with two binary operations ⊕ and ⊙ is called fuzzy right seminear-ring if {P, ⊕} and {P, ⊙} are right seminear-ring.*

All results can be easily extended in case of fuzzy seminear-ring.

***Example 1.8.3***: Let {P, ⊕, ⊙} be the fuzzy seminear-ring. Define ⊕ as p ⊕ q if p + q < 1 and p ⊕ q = 0 if p + q ≥ 1. Then (P, ⊕) is a semigroup. Define ⊙ as p ⊙ q = p for all p, q ∈ P. Clearly {P, ⊕, ⊙} is a special fuzzy seminear-ring.

Now we proceed on to define binear-rings and the concept of fuzzy binear-rings.

**DEFINITION 1.8.25**: *Let (N, +, •) be a non empty set. We call N a binear-ring if N = $N_1 \cup N_2$ where $N_1$ and $N_2$ are proper subsets of N i.e. $N_1 \not\subset N_2$ or $N_2 \not\subset N_1$ satisfying the following conditions:*

*Atleast one of ($N_i$, +, •) is a right near-ring (i = 1, 2) i.e. for preciseness we say*

   i. *($N_1$, +, •) is a near-ring.*
   ii. *($N_2$, +, •) is a ring.*

*We say that even if both ($N_i$, +, •) and ($N_2$, +, •) are right near-rings still we call (N, +, •) to be a binear-ring. By default of notation we mean by a binear-ring only a right binear-ring unless explicitly stated.*

**DEFINITION 1.8.26**: *Let (N, +, •) be a binear-ring. We call (N, +, •) as abelian if (N, +) is abelian i.e. if N = $N_1 \cup N_2$ then ($N_1$, +) and ($N_2$, +) are both abelian. If both ($N_1$, •) and ($N_2$, •) are commutative then we call N a commutative binear-ring. If N = $N_d$ i.e. $N_1 = (N_1)_d$ and $N_2 = (N_2)_d$ then we say N is a distributive binear-ring. If all non-zero elements of N are left (right) cancelable we say that N fulfills the left (right) cancellation law. N is a bi-integral domain if both $N_1$ and $N_2$ has no zero divisors. If N \ {0} = $N_1$ \ {0} and $N_2$ \ {0} are both groups then we call N a binear field.*

**DEFINITION 1.8.27**: *Let (P, +) be a bigroup i.e. (P = $P_1 \cup P_2$) with 0 and let N be a binear-ring μ : N × P → P is called the N-bigroup if for all $p_i \in P_i$ and for all n, $n_i \in N_i$ we have (n + $n_i$) p = np + $n_i$p and (n$n_i$) p = n($n_i$p) for i = 1, 2. $N^P = N_1^P \cup N_2^P$ stands for N-bigroups.*

**DEFINITION 1.8.28**: *A sub-bigroup M of a binear-ring N with M.M ⊂ M is called a sub-binear-ring of N. A bi subgroup S of $N^P$ with NS ⊂ S is a N-sub-bigroup of P.*

**DEFINITION 1.8.29**: *Let N = $N_1 \cup N_2$ be a binear-ring; P a N-bigroup. A binormal subgroup (or equivalently we can call it as normal bisubgroup) I of (I = $I_1 \cup I_2$) (N, +) is called a bi-ideal of N if*

   i. *$I_1 N_1 \subset I_1$ and $I_2 N_2 \subset I_2$.*
   ii. *n, $n_1 \in I_1$, $i_1 \in I_1$, n($n_1 + i_1$) – n $n_1 \in I_1$ and n, $n_2 \in N_2$; $i_2 \in I_2$, n ($n_2 + i_2$) – n$n_2 \in I_2$.*

*Normal sub-bigroup T of (N,+) with (i) is called right bi ideal of N while normal sub-bigroup L of (N,+) with (ii) are called left biideals.*

*A normal sub-bigroup S of P is called bi ideal of $N^P$ if for $s_i \in P_i$ (i = 1, 2) (P = $P_1 \cup P_2$) and $s \in S_i$ (i = 1, 2 and $S = S_1 \cup S_2$) for all $n_i \in N_i$ (i = 1, 2, $N = N_1 \cup N_2$). $n_i (s + s_i) - ns \in S_i$, i = 1,2. Factor binear-ring N / I and factor N-bigroup P/S are defined as in case of birings.*

**DEFINITION 1.8.30:** *A sub-binear-ring M of the binear-ring N is called bi-invariant if $MN_1 \subset M_1$ and $M_2 N_2 \subset M_2$ (where $M = M_1 \cup M_2$ and $N = N_1 \cup N_2$) and $N_1 M_1 \subset M_1$ and $N_2 M_2 \subset M_2$.*

A minimal bi-ideal, minimal right bi-ideal and minimal left bi-ideal and dually the concept of maximal right bi-ideal, left bi-ideal and maximal bi-ideal are defined as in case of bi-rings.

**DEFINITION 1.8.31:** *Let N be a binear-ring ($N = N_1 \cup N_2$, +, •) and S a sub-bisemigroup of ($N = N_1 \cup N_2$, +) (where $S = S_1 \cup S_2$). A binear-ring $N_S$ is called a binear-ring of left (right) quotients with respect to S if*

  i. $N_S = (N_1)_{S_1} \cup (N_2)_{S_2}$ *has identity.*
  ii. *N is embeddable in $N_S$ by a binear-ring homomorphism h.*
  iii. *For all $s_i \in S_i$ (i = 1, 2; $S = S_1 \cup S_2$). $h(s_i)$ is invertible in $((N_i)_{S_i}, •)$, i = 1, 2.*
  iv. *For all $q_i \in N_{S_i}$ there exists $s_i \in S_i$ and there exists $n_i \in N_i$ such that $q_i = h(n_i) h(s_i)^{-1}$; ($q_i = h(s_i)^{-1} h(n_i)$), i = 1, 2.*

**DEFINITION 1.8.32:** *The binear-ring N ($N = N_1 \cup N_2$) is said to fulfill the left (right) ore condition with respect to a given sub-bisemigroup $S_i$ of ($N_i$, •) if for all (s, n) $\in S_i \times N_i$ there exists $n • s_1 = s • n_1$ ($s_1 • n = n_1 • s$); i = 1, 2.*

Let V denote the collection of all binear-rings and X be any non-empty subset.

**DEFINITION 1.8.33:** *A binear-ring $F_X \in V$ is called a free binear-ring in V over a binear-ring N if there exists f: $X \to F_X$ (where X is any non-empty set) for all $N \in V$ and for all g: $X \to N$ there exists a homomorphism $h \in Hom(F_X, N)$; h o f = g*

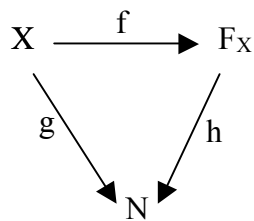

*Let V = {set of all binear-rings}, we simply speak about free binear-ring on X. A binear-ring is called free if it is free over some set X.*

**DEFINITION 1.8.34:** *A finite sequence $N = N_0 \supset N_1 \supset N_2 \supset ... \supset N_t = \{0\}$ of sub-binear-rings $N_i$ of the binear-ring N is called the binormal sequence of N if and only if for all $i \in \{1, 2, ..., t\}$, $N_i$ is a bi-ideal of $N_{i-1}$.*

*In the special case when all the $N_i$'s are bi-ideals of the binear-ring N, we call the binormal sequence an bi-invariant sequence.*

**DEFINITION 1.8.35:** *Let P be a bi-ideal of a binear-ring N. P is called the prime bi-ideal if for all bi-ideals I and J of N; $IJ \subset N$ implies $I \subset P$ or $J \subset P$. The binear-ring N is called a prime binear-ring if $\{0\}$ is a prime ideal.*

**DEFINITION 1.8.36:** *Let S be a bi-ideal of a binear-ring N. S is semiprime if and only if for all bi-ideals I on N, $I^2 \subset S$ implies $I \subset S$. N is called a semiprime binear-ring if $\{0\}$ is a semiprime bi-ideal.*

**DEFINITION 1.8.37:** *Let L be a left bi-ideal in a binear-ring N. L is called modular if and only if there exists $e \in N_1 \cup N_2$ and for all $n \in N_1 \cup N_2$ and for all $n \in N_1 \cup N_2$, $n - ne \in L$. In this case we also say that L is modular by e and that e is a right identity modulo L, since for all $n \in N$, $ne = n \pmod{L}$.*

*Notation*: For $z \in N$, denote the bi-ideal generated by the set $\{n - nz \mid n \in N\}$ by $L_Z$. $L_Z = N$ if $z = 0$. $z \in N$ is called quasi regular if $z \in L_Z$, $S \subset N$ is called quasi-regular if and only if for all $s \in S$, s is quasi regular.

Let N be a binear-ring. An idempotent $e \in N = N_1 \cup N_2$ is called central if it is in the center of $(N_1, \bullet)$ or $(N_2, \bullet)$ i.e. for all $n \in N_1$ (or $n \in N_2$) we have $ne = en$.

**DEFINITION 1.8.38:** *A binear-ring is said to be biregular if there exists some set E of central idempotents with*

    i.    *For all $e \in N_1$ ($e \in N_2$), $N_1 e$ ($N_2 e$) is an ideal of $N_1$ (or $N_2$). i.e. $N_1 e \cup \{0\}$ or $\{0\} \cup N_2 e$ or if $e \in N_1 \cap N_2$ then $N_1 e \cup N_2 e$ is a bi-ideal of N.*

**DEFINITION 1.8.39:** *Let N be a binear-ring; $N = N_1 \cup N_2$ where $N_1$ and $N_2$ are near-rings. The binear-ring N is said to fulfill the insertion factors property (IFP) provided for all $a, b, n \in N_1$ (or $a, b, n \in N_2$) we have $ab = 0$ implies $anb = 0$. The biring has strong IFP property if every homomorphic image of N has the IFP. The binear-ring N has the strong IFP if and only if for all I in N and for all $a, b, n \in N_1$ ($a, b, n \in N_2$), $a, b \in I_1$ implies $anb \in I_1$ ($ab \in I_2$ implies $anb \in I_2$) where $I = I_1 \cup I_2$.*

**DEFINITION 1.8.40:** *Let p be a prime. A binear-ring N is called a p-binear-ring provided that for all $x \in N$, $x^p = x$ and $px = 0$.*

**DEFINITION 1.8.41:** *A biright ideal I of a binear-ring N is called right quasi bireflexive if whenever A and B are bi-ideals of the binear-ring N with $AB \subset I$ then $b(b' + a) - bb' \in I$ for all $a \in A$ and $b, b' \in B$. A binear-ring N is strongly sub-bicommutative if it is right quasi bireflexive.*

**DEFINITION 1.8.42:** *Let N be a binear-ring S a subnormal sub-bigroup of (N, +). S is called a quasi bi-ideal of N if $SN \subset NS \subset S$ whereby NS we mean elements of the form $\{n(n' + s) - nn' /$ for all $s \in S$ and for $n, n' \in N\} = NS$.*

**DEFINITION 1.8.43:** *An infra binear-ring (INR) is a triple (N, +, •), (N = $N_1 \cup N_2$) where*

  i.   *(N, +) is a bigroup.*
  ii.  *(N, •) is a bisemigroup.*
  iii. *(x + y) • z = x • z – 0 • z + y • z for all x, y, z $\in$ N.*

**DEFINITION 1.8.44:** *Let (N, +, •) be a non-empty set, where $N = N_1 \cup N_2 \cup N_3 \cup N_4$ where each $N_i$ is a near-ring or a ring for each i = 1, 2, 3, 4. (N, +, •) is defined as the quad near-ring or bi-binear-ring.*

*Example 1.8.4:* Let (N, +, •) be a non-empty set, where $N = Z_2 \cup Z_{12} \cup Z_7 \cup Z$ where $Z_2$ and Z are near-rings and $Z_{12}$ is the ring of integers modulo 12 and $Z_7$ is the prime field of characteristic 7. N is a quad near-ring.

**DEFINITION 1.8.45:** *Let (N, +, •) be a binear-ring where $N = N_1 \cup N_2$ and $G = G_1 \cup G_2$ be a bigroup. The bigroup binear-ring $NG = N_1G_1 \cup N_1G_2 \cup N_2G_1 \cup N_2G_2$ is a quad near-ring where each $N_iG_j$ is a group near-ring or group ring, i = 1, 2 and j = 1, 2. Similarly we can define bisemigroup binear-ring $NS = N_1S_1 \cup N_2S_1 \cup N_1S_2 \cup N_2S_2$ where S is a bisemigroup ($S = S_1 \cup S_2$) and $N = N_1 \cup N_2$ is a binear-ring. It is easily verified that the bisemigroup binear-ring is also a quad near-ring where $N_iS_j$ are semigroup near-rings or semigroup rings, i = 1, 2 and j = 1, 2. All results connected with binear-rings can be extended in case of quad near-rings.*

**DEFINITION 1.8.46:** *Let (N, +, •) is a non-empty set where $N = N_1 \cup N_2$. We say N is a biseminear-ring if one of ($N_1$, +, •) or ($N_2$, +, •) are seminear-rings. Thus for N to be biseminear-ring we need one of ($N_1$, +, •) or ($N_2$, +, •) or both to be a seminear-ring.*

We can have several such biseminear-rings. All notions studied in case of binear-rings can also be easily studied and generalized for biseminear-rings. The only means to generate several new classes of biseminear-rings are as follows:

**DEFINITION 1.8.47:** *Let (S, +, •) be a biseminear-ring with $S = S_1 \cup S_2$. G any group the group biseminear-ring $SG = S_1G \cup S_2G$ where $S_iG$ are group seminear-rings. Clearly every group biseminear-ring is a biseminear-ring.*

By using this definition we can have infinite class of both commutative or non-commutative and infinite or finite classes of biseminear-rings. If in the definition of group biseminear-rings if we replace the group by a semigroup we get yet another new class of biseminear-rings.

**DEFINITION 1.8.48:** *Let (N, +, •) be a non-empty set. We say (N, +, •) is a quad seminear-ring or bi-biseminear-ring if the following conditions hold good.*

i. $N = N_1 \cup N_2 \cup N_3 \cup N_4$ where $N_i$ are proper subsets of N such that $N_i \not\subset N_j$ for any $i \neq j$, $i, j = 1, 2, 3, 4$.

ii. $(N_i, +, \bullet)$ is at least a seminear-ring for $i = 1, 2, 3, 4$ under the operations of N. Then we call $(N, +, \bullet)$ a quad seminear-ring or bi-biseminear-ring.

**DEFINITION 1.8.49:** Let $(N, +, \bullet)$ be a binear-ring i.e. $N = N_1 \cup N_2$ where $N_1 \not\subset N_1$ or $N_2 \not\subset N_1$ such that $(N_1, +, \bullet)$ is a near-ring and $(N_2, +, \bullet)$ is a ring. The map $\mu : N \rightarrow [0,1]$ is said to be a fuzzy bi near-ring if $\mu = \mu_1 \cup \mu_2$ where $\mu_1 : N_1 \rightarrow [0, 1]$ is a fuzzy near-ring and $\mu_2 : N_2 \rightarrow [0,1]$ is a fuzzy ring; the '$\cup$' is just only a notational convenience.

Equivalently we can define $\mu$ restricted to $N_1$ is a fuzzy near-ring and $\mu$ restricted to $N_2$ is a fuzzy ring where $\mu = \mu_{N_1} \cup \mu_{N_2} = \mu_1 \cup \mu_2$. Here $\mu_{N_1}$ denotes the restriction of $\mu$ to $N_1$ and $\mu_{N_2}$ denotes the restriction of $\mu$ to $N_2$.

All notions as in case of fuzzy near-rings can be extended to fuzzy binear-rings in an analogous way using appropriately the fuzzy near-ring concept and the fuzzy ring concept. However we will be defining in the later chapters the concept of Smarandache fuzzy binear-rings and their properties.

## 1.9 Fuzzy vector spaces and fuzzy bivector spaces

In this section we introduce the concept of fuzzy bivector spaces and recall the definition and properties of fuzzy vector spaces. The study of fuzzy vector spaces started as early as 1977; but till date the study of fuzzy bivector spaces is absent.

For more about fuzzy vector spaces refer [75, 79, 91, 92, 96, 101].

Throughout this section V denotes a vector space over a field F. A fuzzy subset of a non-empty set S is a function from S into [0, 1]. Let A denote a fuzzy subspace of V over a fuzzy subfield K of F and let X denote a fuzzy subset of V such that $X \subseteq A$. Let $\langle X \rangle$ denote the intersection of all fuzzy subspaces of V over K that contain X and are contained in A.

**DEFINITION 1.9.1:**

i. A fuzzy subset K of F is a fuzzy subfield of F, if $K(1) = 1$ and for all $c, d \in F$, $K(c - d) \geq \min \{K(c), K(d)\}$ and $K(cd^{-1}) \geq \min \{K(c), K(d)\}$ where $d \neq 0$.

ii. A fuzzy subset A of V is a fuzzy subspace over a fuzzy subfield K of F, if $A(0) > 0$ and for all $x, y \in V$ and for all $c \in F$, $A(x - y) \geq \min \{A(x), A(y)\}$ and $A(cx) \geq \min \{K(c), A(x)\}$. If K is a fuzzy subfield of F and if $x \in F$, $x \neq 0$, then $K(0) = K(1) \geq K(x) = K(-x) = K(+ x^{-1})$.

In the following we let L denote the set of all fuzzy subfields of F and let $A_L$ denote the set of all fuzzy subspaces of V over $K \in L$.

If A and B are fuzzy subsets of V then A ⊂ B means A(x) ≤ B(x) for all x ∈ A. For 0 ≤ t ≤ 1, let $A_t = \{x \in V \mid A(x) \geq t\}$.

**THEOREM 1.9.1:** *Let A be a fuzzy subset of V and let s, t ∈ Im(A).*

*Then*

    i.    $s \leq t$ *if and only if* $A_s \supseteq A_t$,
    ii.   $s = t$ *if and only if* $A_s = A_t$.

*Proof*: It is a matter of routine hence left for the reader as an exercise.

If A ∈ $A_L$ for every fuzzy subfield K' of K. For S a subset of F, we let $\delta_S$ denote the characteristic function of S.

**THEOREM 1.9.2:** *Let* $A \in A_{\delta_S}$. *Then for all t such that* $0 \leq t \leq A(0)$. $A_t$ *is a subspace of V.*

*Proof:* Straightforward hence left for the reader as an exercise.

**THEOREM 1.9.3:** *Let A be a fuzzy subset of V. If* $A_t$ *is a subspace of V for all t ∈ Im(A), then* $A \in A_{\delta_F}$.

*Proof:* Proof is a matter of routine and the reader is expected to prove. If A is a fuzzy subspace of V, then $A_t$ is called a level subspace of V where $0 \leq t \leq A(0)$.

**THEOREM 1.9.4:** *Let* $V_1 \subset V_2 \subset \ldots \subset V_i \subset \ldots$ *be a strictly ascending chain of subspaces of V. Define the fuzzy subset A of V by* $A(x) = t_i$, *if* $x \in V_i \setminus V_{i-1}$ *where* $t_i > t_{i+1}$ *for i = 1, 2, … and* $V_0 = \phi$ *and A(x) = 0 if* $x \in V \setminus \bigcup_{i=1}^{\infty} V_i$. *Then* $A \in A_{\delta_F}$.

*Proof:* Let c ∈ F. If $x \in V \setminus \bigcup_{i=1}^{\infty} V_i$ then A(cx) ≥ 0 = A(x). Suppose that x ∈ $V_i \setminus V_{i-1}$. Then cx ∈ $V_i$, thus A(cx) ≥ $t_i$ = A(x) ≥ min {$\delta_F$(c), A(x)}.

**THEOREM 1.9.5:** *Let* $V = V_o \supset V_1 \supset \ldots \supset V_i \supset \ldots$ *be a strictly descending chain of subspaces of V. Define the fuzzy subset A of V by* $A(x) = t_{i-1}$ *if* $x \in V_{i-1} \setminus V_i$ *where* $t_{i-1} < t_i < 1$ *for i = 1, 2, …, and A(x) = 1 if* $x \in \bigcap_{i=1}^{\infty} V_i$. *Then* $A \in A_{\delta_F}$.

*Proof:* Let c ∈ F. If $x \in \bigcap_{i=1}^{\infty} V_i$ then $cx \in \bigcap_{i=1}^{\infty} V_i$ and so A(cx) = 1 ≥ min {$\delta_F$(c), A(x)}. Suppose that x ∈ $V_{i-1} \setminus V_i$. Then cx ∈ $V_{i-1}$ thus A(cx) ≥ $t_{i-1}$ = A(x) ≥ min {$\delta_F$(c), A(x)}. If A ∈ $A_K \cap A_{K'}$ for K, K' ∈ L then Im (A) is fixed no matter whether we consider A in $A_K$ or A in $A_{K'}$.

Now using all the above results the following theorem can be proved which is left as an exercise for the reader.

**THEOREM 1.9.6:** *V is finite dimensional over F if and only if every $A \in A_{\delta_F}$ is finite valued.*

The following condition which is given as the theorem is also assigned for the reader to prove.

**THEOREM 1.9.7:** *Let S be a non empty subset of V. Then S is a subspace of V if and only if $\delta_S \in A_{\delta_F}$.*

Now we proceed onto recall the concept of fuzzy spanning.

**DEFINITION 1.9.2:** *Let $A_1, A_2, \ldots, A_n$ be fuzzy subsets of V and let K be a fuzzy subset of F.*

  i. *Define the fuzzy subset $A_1 + \ldots + A_n$ on V by the following: for all $x \in V$, $(A_1 + \ldots + A_n)(x) = \sup \{\min \{A_1(x_1), \ldots, A_n(x_n)\} / x = x_1 + \ldots + x_n, x_i \in V\}$.*
  ii. *Define the fuzzy subset K o A of V by for all $x \in V$, $(K \circ A)(x) = \sup\{\min\{K(c), A(y)\} \mid c \in F, y \in V, x = c\,y\}$.*

**DEFINITION 1.9.3:** *Let S be a set $x \in S$ and $0 \leq \lambda \leq 1$. Define the fuzzy subset $x_\lambda$ of S by $x_\lambda(y) = \lambda$ if $y = x$ and $x_\lambda(y) = 0$ if $y \neq x$. $x_\lambda$ is called a fuzzy singleton.*

**THEOREM 1.9.8:** *Let A be a fuzzy subset of V and let K be a fuzzy subset of F. Let $d \in F$ and $x \in V$. Suppose that $0 \leq \mu, \lambda \leq 1$.*

*Then for all $z \in V$.*

  i. $(d_\mu \circ A)(z) = \min \left\{\mu, A\left(\dfrac{1}{d}z\right)\right\}$ *if $d \neq 0$.*
  ii. $(0_\mu \circ A)(z) = \sup \{\min \{\mu, A(y)\}\, y \in V\}$ *if $z = 0$.*
  iii. $(K \circ x_\lambda)(z) = \sup \{\min \{K(c), \lambda\} \mid c \in F, z = cx\}$ *if $x \neq 0$ and $z \in S_p(x)$; 0 if $x \neq 0$ and $z \notin S_p(x)$.*
  iv. $(K \circ 0_\lambda)(z) = \sup \{\min \{K(c), \lambda\} / c \in F\}$ *if $z = 0$, $(K \circ 0_\lambda)(z) = 0$ if $z \neq 0$.*

*Proof:* Follows by simple and routine work.

**THEOREM 1.9.9:** *Let $c, d \in F$, $x, y \in V$ and $0 \leq k, \lambda, \mu, \nu \leq 1$.*

*Then,*

$d_\mu \circ x_\lambda = (dx)_{\min(\mu, \lambda)}$, $x_\lambda + y_\nu = (x + y)_{\min(\lambda, \nu)}$, $d_\mu \circ x_\lambda + c_K \circ y_\nu = (dx + cy)_{\min\{\mu, \lambda, k, \nu\}}$.

*Proof:* Left for the reader to refer [75].

The following theorem is also expected to be proved by the reader or refer [75].

**THEOREM 1.9.10:** *Let $A \in A_K$ and let B, C be fuzzy subsets of V. Let $b, c \in F$. If $B \subseteq A$ and $C \subseteq A$ then $b_\mu \circ B + C_\nu \circ C \subseteq A$ where $0 \leq \mu \leq K(b)$ and $0 \leq \nu \leq K(c)$.*

**DEFINITION 1.9.4:** *Let $\{A_i \mid i \in I\}$ be a non empty collection of fuzzy subsets of V. Then the fuzzy subset $\bigcap_{i \in I} A_i$ of V is defined by the following; for all $x \in V$*

$$\left(\bigcap_{i \in I} A_i\right)(x) = \inf\{A_i(x) \mid i \in I\}.$$

**THEOREM 1.9.11:** *If $\{A_i \mid A_i \in A_K, i \in I\}$ is non empty then $\bigcap_{i \in I} A_i \in A_K$.*

*Proof:* Let $c \in F$ and $x \in V$. Then

$$\left(\bigcap_{i \in I} A_i\right)(cx) = \inf\{A_i(cx) \mid i \in I\} \geq \inf\{\min\{K(c), A_i(x)\} \mid i \in I\}$$

= either $K(c)$ or $\inf\{A_i(x) \mid i \in I\}$.

Hence

$$\left(\bigcap_{i \in I} A_i\right)(cx) \geq \min\left\{K(c), \bigcap_{i \in I} A(x)\right\}.$$

**DEFINITION 1.9.5:** *Let $A \in A_K$ and let X be a fuzzy subset of V such that $X \subseteq A$. Let $\langle X \rangle$ denote the intersection of all fuzzy subspaces of V (over K) that contain X and are contained in A. Then $\langle X \rangle$ is called the fuzzy subspace of A, fuzzily spanned (or generated) by X.*

**THEOREM 1.9.12:** *Let $A \in A_K$ and let X be a fuzzy subset of V such that $X \subseteq A$. Define the fuzzy subset S of V by the following; for all $x \in V$,*

$$S(x) = \sup\left\{\left(\sum_{i=1}^{n} c_{i\mu_i} \, o \, x_{i\lambda_i}\right)(x) \mid c_i \in F, x_i \in V, K(c_i) = \mu_i,\right.$$

$\left. X(x_i) = \lambda_i, i = 1, 2, \cdots, n, n \geq 1\right\}$. *Then $\langle X \rangle = S$ and $S \in A_K$.*

*Proof:* Follows by routine calculations. For proof refer [75].

Now we proceed on to define the notion of fuzzy freeness. Let $\zeta$ denote a set of fuzzy singletons of V such that $x_\lambda, x_K \in \zeta$, then $\lambda = k > 0$. Define the fuzzy subset $X(\zeta)$ of V by the following; for all $x \in V$, $X(\zeta)(x) = \lambda$ if $x_\lambda \in \zeta$, and $X(\zeta)(x) = 0$, otherwise. Define $\langle \zeta \rangle = \langle X(\zeta) \rangle$. Let X be a fuzzy subset of V. Define $\zeta(X) = \{x_\lambda \mid x \in V, \lambda = X(x) > 0\}$. Then $X(\zeta(x)) = X$ and $\zeta X(\zeta)) = \zeta$. If there are only a finite number of $x_\lambda \in \zeta$, with $\lambda > 0$ we call $\zeta$ finite. If $X(x) > 0$ for only a finite number of $x \in X$, we call X finite. Clearly $\zeta$ is finite if and only if $X(\zeta)$ is finite and X is finite if and only if $\zeta(X)$ is finite. For $x \in V$ let $X \setminus x$ denote the fuzzy subset of V defined by the following; for all $y \in V$ let $(X \setminus x)(y) = X(y)$ if $y \neq x$ and $(X \setminus x)(y) = 0$ if $y = x$.

**DEFINITION 1.9.6:** *Let $A \in A_K$ and let X be a fuzzy subset of V such that $X \subset A$. Then X is called a fuzzy system of generators of A over K if $\langle X \rangle = A$. X is said to be fuzzy free over K if for all $x_\lambda \in X$ where $\lambda = X(x)$, $x_\lambda \not\subset \langle X \setminus x \rangle$. X is said to be a fuzzy basis for A if X is a fuzzy system of generators of A and X is fuzzy free. Let $\zeta$ denote a set of fuzzy singletons of V such that if $x_\lambda, x_K \in \zeta$ then $\lambda = k$ and $x_\lambda \subseteq A$. Then $\zeta$ is called a fuzzy singleton system of generators of A over K. if $\langle \zeta \rangle = A$. $\zeta$ is said to be fuzzy free over K, if for all $x_\lambda \in \zeta$, $x_\lambda \not\subset \langle \zeta \setminus \{x_\lambda\} \rangle$. $\zeta$ is said to be a fuzzy basis of singletons for A if $\zeta$ is a fuzzy singleton system of generators of A and $\zeta$ is fuzzy free.*

**THEOREM 1.9.13:** *Suppose that $A \in A_K$. Then*

  i. *$K^*$ is a subfield of F.*
  ii. *$A^*$ is a subspace of V over $K^*$.*

*Proof*: The work is assigned to the reader as exercise.

**THEOREM 1.9.14:** *Let $A \in A_K$ and let $\varsigma \subseteq \{x_\lambda \mid x \in A^*, 0 < \lambda \leq A(x)\}$ be such that of $x_\lambda, x_k \in \zeta$ then $\lambda = k$ and let $\chi = \{x \mid x_\lambda \in \zeta\}$. Suppose inf $\{k(c) \mid c \in F\} \geq \sup \{A(x) \mid x \in V \setminus \{(0)\}\}$. Then $\zeta$ is fuzzy free over k if and only if x is linearly independent over F.*

*Proof:* It is a matter of routine [75].

**THEOREM 1.9.15:** *Let $A \in A_k$. Let $\varsigma = \{x_\lambda \mid x \in A^*, 0 < \lambda \leq A(x)\}$ be such that if $x_\lambda, x_k \in \zeta$ then $\lambda = k$ and $\chi = \{x \mid x_\lambda \in \zeta\}$. Suppose inf $\{k(c) \mid c \in F\} \geq \sup \{A(x) \mid x \in V \setminus \{0\}\}$. Then $\zeta$ is maximally fuzzy free in A over k if and only if $\chi$ is a basis for $A^*$ over $k^*$.*

*Proof:* Using the fact if $A(x) = 0$ for all $x \in V \setminus \{0\}$ then the result holds with $\zeta$ and $\chi$ empty.

Suppose $\{0\} \subset A^*$. Then $k^* = f$. Suppose $\zeta$ is maximal fuzzy free, then $\chi$ is linearly independent over F. Using earlier theorems we get with simple calculations that $\zeta$ is maximal.

The following is a simple consequence of the above theorem.

**THEOREM 1.9.16:** *Let $A \in A_k$. Suppose that inf $\{k(c) \mid c \in F\} \geq \sup \{A(x) \mid x \in V \setminus \{0\}\}$. Then A has maximally fuzzy free sets over k of fuzzy singletons of V and every such set has the same cardinality.*

The following theorem gives condition for the existence of fuzzy basis.

**THEOREM 1.9.17:** *Let $A \in A_k$. Suppose inf $\{k(c) \mid c \in F\} \geq \sup \{A(x) \mid x \in V \setminus \{0\}\}$. If A is finite valued, then A has a fuzzy basis over k.*

*Proof:* Matter of routine, hence left for the reader as an exercise.

The following results are immediate consequence of the above Theorem hence it is left for the reader to prove.

**THEOREM 1.9.18:** *Let $A \in A_k$. Suppose that inf $\{k(c) \mid c \in F\} \geq \sup A(x)$ such that $x \in V \setminus \{0\}\}$. If V is finite dimensional, then A has a fuzzy basis over k.*

**THEOREM 1.9.19:** *Let $A \in A_k$. Suppose inf $\{k(c) \mid c \in F\} \geq \sup \{A(x) \mid x \in V \setminus \{0\}\}$. If A is finitely fuzzily generated over k, then A has a fuzzy basis over k.*

Now we define fuzzy linearly independent set.

**DEFINITION 1.9.7:** *Let $A \in A_k$ and let $\zeta \subseteq \{x_\lambda \mid x \in A^*, \lambda \leq A(x)\}$ be such that if $x_\lambda, x_k \in \zeta$, then $\lambda = k$. Then $\zeta$ is said to fuzzy linearly independent over k if and only if for every finite subset $\{x_{1_{\lambda_1}}, \ldots, x_{n_{\lambda_n}}\}$ of $\zeta$, whenever $\{(\sum_{i=1}^{n} c_{i\mu_i} \circ x_{i\lambda_i})(x) = 0$ for all $x \in V \setminus \{0\}$ where $c_i \in F$, $0 < \mu_i \leq K(c_i)$ for $i = 1, 2, \ldots, n\}$ then $c_1 = c_2 = \ldots = c_n = 0$.*

Finally we recall the following theorem with proof.

**THEOREM 1.9.20:** *Let $A \in A_K$ and let $\zeta \subseteq \{x_\lambda \mid x \in A^*, 0 < \lambda \leq A(x)\}$ be such that if $x_\lambda, x_K \in \zeta$, then $\lambda = K$ and let $\chi = (x \mid x_\lambda \in \zeta)$. Then $\zeta$ is fuzzy linearly independent over K if and only if $\chi$ is linearly independent over $K^*$.*

*Proof:* Suppose $\zeta$ is fuzzy linearly independent over K. Suppose

$$0 = \sum_{i=1}^{n} c_i x_i \text{ where } c_i \in K^* \text{ and } x_i \in \chi, i = 1, 2, \ldots, n.$$

Let $\lambda = \min \{\mu_1, \ldots, \mu_n, \lambda_1, \ldots, \lambda_n\}$ where $0 < \mu_I \leq K(c_i)$ and $0 < \lambda_i \leq A(x_i)$ for $i = 1, 2, \ldots, n$. Then for all $x \in V \setminus \{0\}$,

$$0 = \left(\sum_{i=1}^{n} c_i x_i\right)_\lambda (x) = \sum_{i=1}^{n} \left(C_{i\mu_i} \circ x_{i\lambda_i}\right)(x).$$

Hence $c_1 = c_2 = \ldots = c_n = 0$. Conversely suppose $\chi$ is linearly independent over $K^*$. Let $x_{1_{\lambda_1}}, \ldots, x_{n_{\lambda_n}} \in \zeta$. Suppose that for all $x \in V \setminus \{0\}$, $0 = \sum (c_{i\mu_i} \circ x_{i\lambda_i})(x)$ then $0 = \left(\sum_{i=1}^{n} c_i x_i\right)_\lambda (x)$. Since $\lambda > 0$, $\sum_{i=1}^{n} c_i x_i = 0$. Hence $c_1 = c_2 = \ldots = c_n = 0$. Hence the theorem.

Now we proceed on to define the notion of fuzzy algebraically independent. Let X be a fuzzy subset of F, we let $\phi(X) = \{x_t \mid t = X(x) > 0\}$. If $\phi$ is a set of fuzzy singletons such that $x_t, x_s \in \phi$ implies $t = s$, then we let $X(\phi)$ denote the fuzzy subset of F defined by $(X(\phi))(x) = t$ if $x_t \in \phi$ and $(X(\phi))(x) = 0$ if $x_t \notin \phi$. Clearly $X(\phi(X)) = X$ and $\phi(X(\phi)) = \phi$. We let $i = (i_1, \ldots, i_n)$. If F is a field F(F) denotes the set of all fuzzy subfields of F. F(A) will denote the set of all fuzzy subsets X of F such that $X \subseteq A$ and F(A/B) the set of all fuzzy subfields C of F such that $B \subseteq C \subseteq A$.

We in the forthcoming definitions and results assume A, B $\in$ F(F) and $B \subset A$.

Now we proceed on to recall the definitions as given by [91]. The following notation and definition is recalled for they are used in the following results.

**DEFINITION 1.9.8:** *Let A, B $\in$ F(R), $B \subset A$, and let X be a fuzzy subset of R such that $X \subseteq A$. Define B[X] to be the intersection of all $C \in F(R)$ such that $B \cup X \subseteq C \subseteq A$. For A, B $\in$ F(F), $B \subseteq A$ and let X be a fuzzy subset of F such that $X \subset A$. Define B(X) to be the intersection of all $C \in F(F)$ such that $B \cup X \subseteq C \subseteq A$. B [X] is a fuzzy subring of R and B [X] is a fuzzy subfield of F.*

**DEFINITION 1.9.9:** *Let $X \in F(A)$. Then X (or $\phi(X)$) is said to be fuzzy algebraically independent over B if and only if for all $(x_1)_{t_1}, \cdots, (x_n)_{t_n} \in \phi(X)$, for all $b_1, b_2, \ldots, b_n \in F$ for all $s \in (0, 1]$, $\sum (b_i)_{u_i} (x^i)_t = .0$, where $B(b_i) \geq u_i$ and $X(x_j) \geq t_j$, ($j = 1, 2, \ldots, n$) implies $b_i = 0$ for all i.*

*If $C_t \subseteq A$ is fuzzy algebraically independent over B, then $C_t$ is also said to be fuzzy transcendental over B.*

**THEOREM 1.9.21:** *Let $X \in F(A)$. Then X is fuzzy algebraically independent over B if and only if for all $s \in (0, 1]$, $X_S$ is algebraically independent over $B_S$.*

*Proof:* Suppose that X is fuzzy algebraically independent over B. Let $s \in (0, 1]$, Suppose that $0 = \sum b_i (x_1)^{i_1} \cdots (x_n)^{i_n}$ where $b_i \in B_s$ and $x_j \in X_s$ ($j = 1, 2, \ldots, n$). Then $0_s = \sum (b_i)_s (x_1)_s^{i_1} \cdots (x_n)_s^{i_n}$ and since $B(b_i) \geq s$ and $X(x_j) \geq s$, $b_j = 0$ for all i. Thus $X_s$ is algebraically independent over $B_s$. Conversely suppose that $X_s$ is algebraically independent over $B_s$ for all $s \in (0, 1]$. Suppose that $0_s = \sum (b_i)_{u_i} (x^i)_t$ where $B(b_i) \geq u_i$ and $X(x_j) \geq t_j$. Then min $\{\min_i\{u_i\}, \min \{t_j / j = 1, 2, \ldots, n\}\} = s$ and $\sum b_i (x_1)^{i_1} \ldots (x_n)^{i_n} = 0$. Thus $B(b_i) \geq s$ and $X(x_j) \geq s$ and so $b_i \in B_s$ and $x_j \in X_s$. Hence $b_i = 0$ for all i. Thus X is fuzzy algebraically independent over B.

The following theorem can also be a proved as a matter of routine.

**THEOREM 1.9.22:** *Let $X \in F(A)$. Then X is fuzzy algebraically independent over B if and only if $X^*$ is algebraically independent over $B^*$.*

**THEOREM 1.9.23:** *Let $X \in F(A)$.*

*i.   For all $t \in (0, 1]$, $B_t (X_t) \subseteq B(X)_t$.*
*ii.  If B(X) has the sup property then for all $t \in (0, 1]$,*

$$B_t (X_t) = B(X)_t.$$

*Proof:*

i. Let $z \in B_t(X_t)$, then $(B(X))(z) \geq t$ follows easily from properties of fuzzy sets.
ii. Let $z \in B(X)_t$. Then $(B(X))(z) \geq t$ and since $B(X)$ has the sup property, $z$ has a representation as an element of $B_t(X_t)$ by results on fuzzy sets. Thus $z \in B_t(X_t)$.

**THEOREM 1.9.24:** *Let $X \in F(A)$. Then $B(X)^* = B^*(X^*)$ and $B[X^*] = B^*[X^*]$.*

*Proof:* $x \in B(X)^*$ if and only if $(B(X))(x) > 0$ if and only if $x \in B^*(X^*)$. Similarly we can prove $B[X]^* = B^*[X^*]$. Let $X \in F(A)$. We say that $X$ is maximally fuzzy algebraically independent over $B$ if and only if $X$ is fuzzy algebraically independent over $B$ and there does not exist $Y \in F(A)$ such that $Y$ is fuzzy algebraically independent over $B$ and $X \subset Y$.

**THEOREM 1.9.25:** *Let $X \in F(A)$. Suppose for all $x \in X^*$ we have $X(x) = A(x)$. Then $X$ is maximally fuzzy algebraically independent over $B$ if and only if $X^*$ is transcendence basis of $A^*/B^*$.*

*Proof:* By earlier results and routine calculations the result can be arrived.

The following theorem is a direct consequence of the above result. Hence the reader is expected to prove the result.

**THEOREM 1.9.26:** *A/B has maximal fuzzy algebraically independent fuzzy subsets of F and the cardinality of each is unique.*

First we give some notational conventions using which we define the concept of fuzzy transcendental and neutral.

Let $c_t \subseteq A$, $t > 0$. Suppose that $c_t$ is not fuzzy algebraically independent over $B$. Then there exists $n \in N$, there exist $b_i \in B^*$, s, $u_i \in (0, 1]$ such that $B(b_i) = u_i$ for $i = 0, 1, \ldots, n$ and such that $0_s = \sum_{i=1}^{n}(b_i)_{u_i}(c_t)^i$ with not all $b_i = 0$. If the only such s that exists for which such an equation holds are strictly less than t, then $C_t$ is not fuzzy algebraic over $B$.

**DEFINITION 1.9.10:** *Let $c_t \subseteq A$, with $t > 0$ Then $c_t$ is called fuzzy algebraical over $B$ if and only if $c_t$ is not fuzzy transcendental over $B$. If every such $c_t$ is fuzzy algebraical over $B$, then $A/B$ is called fuzzy algebraical; otherwise $A/B$ is called fuzzy transcendental.*

**DEFINITION 1.9.11:** *Let $c_t \subseteq A$, with $t > 0$ If there exists $s \in (0, t]$ such that $c_s \subseteq B$, then $c_t$ is called neutral over $B$. If every such $c_t$ is neutral over $B$, then $A/B$ is called neutral.*

Clearly $A/B$ is neutral if and only if $A^* = B^*$. Suppose that for $c \in F$, $t = A(c) > B(c) = s > 0$. Then $c_t + (-c_s) = 0_s$, so that $c_t$ is fuzzy algebraical over $B$. In fact we can think of $c_t$ as being a root of a first degree polynomial in $x$, $x + (-c_s)$ with $c_s \subseteq B$.

**THEOREM 1.9.27:** *A/B is fuzzy algebraical if and only if $A^*/B^*$ is algebraic.*

*Proof:* Suppose that A/B is fuzzy algebraical. Let $c \in A^*$. Then there exist $n \in N$, $K \in F$ not all 0, there exist s, $\upsilon_I \in [0, 1]$ such that

$$(K_n)_{\upsilon_n}(c_t)^n + \cdots + (K_1)_{\upsilon_1}(c_t) + (K_0)_{\upsilon_0} = 0_s \text{ where } s \le t = A(c) \text{ and } (K_i)_{\upsilon_i} \subseteq B$$

for i = 0, 1, 2, …, n. Thus $K_n c^n + \ldots + K_1 c + K_o = 0$ and $K_i \in B^*$ for i = 0, 1, 2, …, n. That is c is algebraic over $B^*$. Conversely, suppose that $A^*/B^*$ is algebraic. Let $c_t \subseteq A$ with t > 0. Then $c \in A^*$ hence there exist $n \in N$, and $K_i \in B^*$, $K_i$ not all 0, i = 0, 1, …, n such that $K_n c^n + \ldots + K_1 c + K_o = 0$. Thus

$$(K_n)_{\upsilon_n}(c_t)^n + \cdots + (K_1)_{\upsilon_1}(c_t) + (K_0)_{\upsilon_0} = 0_s$$

where $B(K_i) = \upsilon_i$ for i = 0, 1, 2, …, n. and s = min $\{t, \upsilon_o, \upsilon_i, \ldots, \upsilon_n\}$. Hence $c_t$ is fuzzy algebraical over B.

The following theorem is an easy consequence of the earlier theorem.

**THEOREM 1.9.28:** *Let $X \in F(A)$. Then $B(X) = B[X]$ if and only if B(X)/B is fuzzy algebraical.*

**THEOREM 1.9.29:** *Let $C \in F(A/B)$. Then A/B is fuzzy algebraical if and only if A/C and C/B are fuzzy algebraical.*

*Proof:* $A^*/B^*$ is algebraic if and only if $A^*/C^*$ and $C^*/B^*$ are algebraic.

**DEFINITION 1.9.12:** *Let $c_t \subseteq A$ with t > 0. Suppose that c is a root of a polynomial $p(x) = K_n x^n + \ldots + K_1 x + K_0$ over $B^*$. We say that $c_t$ is fuzzy algebraical with respect to p(x) over $B^*$ if and only if $0_s = \sum_{i=1}^{n}(K_i)_{\upsilon_i}(c_i)^i$ for some $\upsilon_i \in (0, 1]$ where $B(K_i) \ge \upsilon_i$ for i = 0, 1, 2, …, n and $s \le t$. For s = t we say that $c_t$ is fuzzy algebraic with respect to p(x) over $B^*$.*

**THEOREM 1.9.30:** *Let $K, b \in B^*$, $K \ne 0$. If $B(K) \ne B(b)$ then $B(Kb)$ = min $\{B(K), B(b)\}$.*

*Proof*: Follows as an easy consequence hence left for the reader to prove.

**THEOREM 1.9.31:** *Let $c_t \subseteq A$ with t > 0. Suppose that c is algebraic over $B^*$. Let p(x) be the minimal polynomial of c over $B^*$. If $c_t$ fuzzy algebraical (not fuzzy algebraic) with respect to p(x) then $c_t$ is fuzzy algebraical (not fuzzy algebraic) with respect to every irreducible polynomial over $B^*$ which has c as a root.*

*Proof*: Let $p(x) = x^n + \ldots + K_1 x + K_o$. Then there exists $K_i$ such that $0 < B(K_i) < t$. Let q(x) be any irreducible polynomial over $B^*$ having c as a root. Then q(x) = Kp(x) for some $K \in B^*$. Suppose that $B(K) \ne B(K_i)$. Then $0 < B(KK_i) = B(K_i) < t$ suppose that $B(K) = B(K_i)$, then $0 < B(K) < t$ and the desired result follows in this case since K is the leading coefficient of $q(x)$.

The following result is a direct consequence of the definitions and a matter of routine.

**THEOREM 1.9.32:** *Define the fuzzy subset $B^{(n)}$ of F by $B^{(n)}(x) = A(x)$ if $x \in B^*$ and $B^n(x) = 0$ if $x \notin B^*$. Then $B^n \in F(A/B)$. Suppose that $\inf \{A(x) \mid x \in B^*\} \geq \sup \{A(x) \mid x \notin B^*\}$. If A/B is fuzzy algebraical then $A/B^{(n)}$ is fuzzy algebraic.*

Now we proceed on to define fuzzy transcendence basis.

**DEFINITION 1.9.13:** *Let $X \in F(A)$. Then X is called a fuzzy transcendence basis of A/B if and only if X is fuzzy algebraically independent over B and A/B(X) is fuzzy algebraical.*

**THEOREM 1.9.33:** *A/B has a fuzzy transcendence basis and the cardinality of a fuzzy transcendence basis is unique. In fact X is a fuzzy transcendence basis of A/B if and only if $X^*$ is a transcendence of $A^*/B^*$.*

*Proof*: Follows easily by the very definition.

Now we proceed on to define fuzzy separable and fuzzy pure inseparable.

**DEFINITION 1.9.14:** *Suppose that $c_t \subset A$ with $t > 0$. Then $c_t$ is said to be fuzzy pure inseparable over B if and only if there exists $e \in N \cup \{0\}$, such that $b_o, b_1 \in B^*$, where s, $\upsilon_0, \upsilon_1 \in (0, 1]$. $B(b_i) = \upsilon_i$ for $i = 0, 1$ such that $(K_1)_{\upsilon_1}(c_t)^{p^e} + (K_0)_{\upsilon_0} = 0_s$ A/B is called fuzzy pure inseparable if and only if every $c_t \subseteq A$ with $t > 0$ is fuzzy pure inseparable over B. $c_t$ is said to be fuzzy separable algebraical over B if and only if there exists $n \in N$, such that $b_i \in B^*$, there exist s, $\upsilon_i \in (0, 1]$, $B(b_i) = \upsilon_i$ for $i = 0, 1, 2, \ldots, n$ such that $0_s = \sum (K_i)_{\upsilon_i}(c_i)^i$ and the polynomial $\sum_{i=0}^{n} K_i x^i$ (in x) is separable over $B^*$. A/B is called fuzzy separable algebraical if and only if every $c_t \subset A$ with $t > 0$ is fuzzy separable algebraical over B.*

*Now $c_t \subseteq A$ ($t > 0$) in neutral over B if and only if $c_t$ is fuzzy pure inseparable and fuzzy separable algebraical over B, yet in either event it is not necessarily the case that $c_t \subseteq B$. If $c_t$ is fuzzy pure inseparable and fuzzy separable algebraic over B, then $c_t \subseteq B$.*

The following theorems are a matter of routine and the proofs can be easily supplied by an innovative reader.

**THEOREM 1.9.34:**

  i. *A/B is fuzzy pure inseparable if and only if $A^*/B^*$ is purely inseparable.*
  ii. *A/B is fuzzy separable algebraical if and only if $A^*/B^*$ is separable algebraic.*

**THEOREM 1.9.35:** *Let $C \in F(A/B)$,*

i. *A/B is fuzzy pure inseparable if and only if A/C and C/B are fuzzy pure inseparable.*
ii. *A/B is fuzzy separable algebraical if and only if A/C and C/B are fuzzy separable algebraical.*

**THEOREM 1.9.36:** *Let $c_t \subseteq A$ with $t > 0$. If $c_t$ is fuzzy algebraical (pure inseparable or separable algebraical) over B, then $B(c_t)/B$ is fuzzy algebraical (pure inseparable or separable algebraical).*

*Proof*: Define the fuzzy subset c of F by $c(x) = A(x)$ if $x \in B(c_t)^*$ and $c(x) = 0$ otherwise. Then $c^* = B(c_t)^* = B^*((c_t)^*) = B^*(c)$ and $B^*(c) / B^*$ i.e. $c^* / B^*$ is either algebraic or purely inseparable or separable algebraic according as c has these properties over $B^*$.

Using earlier result the proof of the theorem is a mater or routine.

Now we proceed on to define fuzzy separating transcendence basis of A / B.

**DEFINITION 1.9.15:** *Let $X \in F(A)$. Then X is called a fuzzy separating transcendence basis of A / B if and only if X is fuzzy algebraically independent over B and A / B(X) is fuzzy separable algebraical.*

The following two theorems are also a matter of routine hence left for the reader as an exercise.

**THEOREM 1.9.37:** *A/B has a fuzzy separating transcendence basis if and only if $A^*/B^*$ has a separating transcendence basis. In fact X is a fuzzy separating transcendence basis of A/B if and only if $X^*$ is a separating transcendence basis of $A^*/B^*$.*

**THEOREM 1.9.38:**

i. *B has a fuzzy algebraical (pure inseparable or separable algebraical) closure in A. In fact $B^{(c)}(B^{(i)}, B^{(s)})$ is the fuzzy algebraical pure inseparable, separable algebraical closure of B in A if and only if $B^{(c)^*}(B^{(i)^*}, B^{(s)^*})$ is algebraic (purely inseparable separable algebraic) closure of $B^*$ in $A^*$ i.e. $B^{(c)^*} = B^{*(c)}(B^{(i)^*} = B^{*(i)}, B^{(s)^*} = B^{*(s)})$*

ii. *$B^{(c)} \supseteq B^{(i)}, B^{(s)}$ and $B^{(c)} / B^{(s)}$ is fuzzy pure inseparable.*

The results are as a matter of routine.

**DEFINITION 1.9.16:** *A fuzzy vector space $(V, \eta)$ or $\eta_V$ is an ordinary vector space V with a map $\eta : V \to [0, 1]$ satisfying the following conditions.*

i. $\eta(a + b) \geq \min \{\eta(a), \eta(b)\}$.
ii. $\eta(-a) = \eta(a)$.

    iii.    $\eta(0) = 1$.
    iv.    $\eta(ra) \geq \eta(a)$ for all $a, b \in V$ and $r \in F$ where $F$ is a field.

**DEFINITION 1.9.17:** *For an arbitrary fuzzy vector space $\eta_V$ and its vector subspace $\eta_W$, the fuzzy vector space $(V/W, \hat{\eta})$ or $\eta_{VW}$ determined by*

$$\hat{\eta}(v + W) = \begin{cases} 1 & \text{if } v \in W \\ \sup_{\omega \in W} \eta(v + \omega) & \text{otherwise} \end{cases}$$

*is called the fuzzy quotient vector space, $\eta_V$ by $\eta_W$.*

**DEFINITION 1.9.18:** *For an arbitrary fuzzy vector space $\eta_V$ and its fuzzy vector subspace $\eta_W$, the fuzzy quotient space of $\eta_V$ by $\eta_W$ is determined by*

$$\overline{\eta}(v + W) = \begin{cases} 1 & v \in W \\ \inf_{\omega \in W} \eta(v + \omega) & v \notin W \end{cases}$$

*It is denoted by $\overline{\eta}_{V/W}$.*

Now similar to vector spaces are the notions of modules. We recall the definition and some of the properties of fuzzy modules.

**DEFINITION 1.9.19:** *A fuzzy R-module $\eta_M$ is an ordinary module $M$ with a map $\eta : M \to [0, 1]$ satisfying the following conditions:*

    i.    $\eta(a + b) \geq \min \{\eta(a), \eta(b)\}$.
    ii.    $\eta(-a) = \eta(a)$.
    iii.    $\eta(0) = 1$.
    iv.    $\eta(ra) \geq \eta(a)$ for all $a, b \in V$ and $r \in R$ where $R$ is a commutative ring with 1.

**DEFINITION 1.9.20:** *For an arbitrary fuzzy module $\eta_M$ and its fuzzy module $(M/N, \hat{\eta})$ or $\eta_{M/N}$ determined by*

$$\hat{\eta}(m + N) = \begin{cases} 1 & \text{if } m \in N \\ \sup_{n \in N} \eta(m + n), & \text{otherwise} \end{cases}$$

*is called the fuzzy quotient module of $\eta_M$ by $\eta_N$.*

**THEOREM 1.9.39:** *Let $D$ be an $R$-module and $N$ be its submodule. Then the map $\hat{\eta} : D/N \to [0, 1]$ given by*

$$\hat{\eta}(d + N) = \sup_{n \in N} \eta(d + n)$$

*determines a fuzzy R-module (D/N, $\hat{\eta}$ ).*

*Proof*: Please refer [103].

**DEFINITION 1.9.21:** *For an arbitrary fuzzy module $\eta_D$ and its fuzzy submodule $\eta_N$ the fuzzy module $\hat{\eta}_{D/N}$ given by the above theorem is called the fuzzy quotient module of $\eta_D$ by $\eta_N$.*

**DEFINITION 1.9.22:** *For a fuzzy R-module $\eta_D$ (note R-module D is divisible if for any $d \in D$ and $r \in R$, $r \neq 0$ there exists some $d' \in D$ with $rd' = d$) and its arbitrary fuzzy submodule $\eta_N$ we can define fuzzy quotient module $\hat{\eta}_{D/N}$.*

Its structure is in fact very complicated refer [102].

**DEFINITION 1.9.23:** *For an arbitrary fuzzy divisible module $\eta_D$ and its fuzzy prime submodule $\eta_N$ , the fuzzy quotient module of $\eta_D$ by $\eta_N$ is given by*

$$\overline{\eta}\,(d + N) = \begin{cases} 1 & d \in N \\ \inf_{n \in N} \eta(d+n), & d \notin N \end{cases}$$

*It is denoted $\overline{\eta}_{D/N}$.*

*Let a vector space V be over R, where R is a field or a division ring. For any $v \in V$ and $r \in R$, $r \neq 0$ we choose $v' = r^{-1}v$ and then $rv' = v$. By definition V is a divisible module. Further let W be a vector subspace of V. For any $v \notin W$ and $r \in R$; $r \neq 0$, $rv \in W$ implies $r^{-1}(rv) \in rW = W$ i.e. $v \in W$. This is contrary to the hypothesis. Hence $rv \notin W$ and so W is prime submodule.*

*For any fuzzy module $\eta_M$ we denote by*

$$M_p = \{a \in M \mid \eta(a) \geq p, p \in [0, 1]\} \text{ and}$$
$$M^p = \{a \in M \mid \eta(a) > p, p \in [0, 1)\}$$

*It is well known that $M_p$ is called p-level set of M. Moreover if $M_p$ is a submodule of M and then we obtain a fuzzy submodule $\eta_{M_p}$ for every $p \in [0, 1]$. It is easy to verify that $M^p$ is a submodule of M as well as we call $M^p$ strong p-level set of M.*

The following three results are given without proof. The interested reader is requested to refer [103].

**THEOREM 1.9.40:** *Let $\eta_M$ be a fuzzy module and $N = M_p$ (a non zero p-level set of M). For every $C_t = \{a \in M \mid \eta(a) = t, t \in [0, p)\}$ there exists some index set I (finite or infinite) such that*

$$C_t = \bigcup_{i \in I}(a_i + N) \text{ with } (a_i + N) \cap (a_j + N) = \phi\,(i \neq j).$$

**THEOREM 1.9.41:** *Let $\eta_M$ be an arbitrary fuzzy module over any ring R and $N = M_P$ be its p-level submodule. Then the map $\eta^*; M \to [0, 1]$ given by*

$$\eta^*(a) = \begin{cases} 1 & a \in N \\ \eta(a) & a \notin N \end{cases}$$

*defines a fuzzy module $(M, \eta^*)$ or $\eta^*_M$.*

**THEOREM 1.9.42:** *Let M be an arbitrary R-module over a ring R and $N = M_p$ be its p-level submodule. Then for any fuzzy module $\eta_M$ the map $\hat{\eta}; M/N \to [0, 1]$ given by*

$$\hat{\eta}(a + N) = \sup_{n \in N} \eta^*(a + n)$$

*defines a fuzzy module.*

*Proof*: Please refer [103].

**DEFINITION 1.9.24:** *For an arbitrary fuzzy module $\eta_M$ and its p-level fuzzy submodule $\eta_N$ $(N = M_p)$ the fuzzy module $\hat{\eta}_{M/N}$ given $\hat{\eta}(a + N) = \sup_{n \in N} \eta^*(a + n)$ is called the fuzzy quotient module of $\eta_N$ and $\eta_M$. Especially we denote $N_1 = \{a \in M / \eta(a) = 1\}$. It is the fuzzy singular submodule of $\eta_M$.*

**THEOREM 1.9.43:** *Let M be an arbitrary R-module over the ring R and $Q = M^p$ be its strong p-level submodule then for any fuzzy module $\eta_M$, the map $\hat{\eta}: M/Q \to [0, 1]$ given by $\hat{\eta}(a + Q) = \sup_{q \in Q} \eta^*(a+q)$ defines a fuzzy module.*

*Proof:* Left for the reader to prove.

**DEFINITION 1.9.25:** *For an arbitrary fuzzy module $\eta_M$ and its strong p-level fuzzy submodule $\eta_Q$ $(Q = M^p)$ the fuzzy module $\hat{\eta}_{M/Q}$ given by $\hat{\eta}: M/Q \to [0, 1]$ such that $\hat{\eta}(a + Q) = \sup_{q \in Q} \eta^*(a + q)$ is called the fuzzy quotient module of $\eta_M$ by $\eta_Q$.*

For more about these please refer the 3 papers of [103].

## 1.10 Fuzzy semigroup and their properties

The study of fuzzy semigroups started in 1979. Several results have been obtained in this direction. However to make our study most recent we recall some facts about fuzzy semigroups from the result of [142, 143]. Researchers have studied many classes of semigroups using fuzzy ideals. Here we recall most notions from [143].

**DEFINITION [54]:** *Let S be a non-empty set completely ordered by transitive, irreflexive relation '<'. For all x, y ∈ S, exactly one relations x < y, x = y, y < x holds.*

*For a, b ∈ S define*

$$(a, b) = \{x \in S \,/\, a < x < b\}$$
$$[a, b) = \{x \in S \,/\, a \leq x < b\}$$
$$(a, b] = \{x \in S \,/\, a < x \leq b\} \text{ and}$$
$$[a, b] = \{x \in S \,/\, a \leq x \leq b\}.$$

This notation we have used already we may use them as an algebraic structure viz. semigroups whenever the situation demands.

Now we define a fuzzy semigroup.

**DEFINITION 1.10.1:** *Let (S, •) be a semigroup. A map $\mu: S \to [0, 1]$ is called a fuzzy semigroup if $\mu(x • y) = \min \{\mu(x), \mu(y)\}$ for all x, y ∈ S.*

From this definition we evolve the following. All algebraic structures viz. fuzzy groups, fuzzy rings, fuzzy near-rings, fuzzy semirings are obviously fuzzy semigroups as at least under one operation we have $\mu : X \to [0, 1]$, X an algebraic structure we have $\mu(x * y) = \min\{\mu(x), \mu(y)\}$. $*$ can be some closed operation on X.

*Let S denote a semigroup. A non-empty subset I of S is called a left-ideal of S if $SI \subset I$. I is called a right ideal if $IS \subset I$. I is said to be an ideal of S, if I is both a left and a right ideal of S.*

*We call a non-empty subset L of S which is a left ideal of S to be prime if for any two ideal A and B of S such that $AB \subseteq L$, it implies that $A \subseteq L$ or $B \subseteq L$. L is called quasi prime if for any two left ideals $L_1$ and $L_2$ of S, $L_1L_2 \subseteq L$ then $L_1 \subset L$ or $L_2 \subset L$ and L is called weakly quasi prime if any one of the two left ideals of $L_1$ and $L_2$ of S such that $L \subseteq L_1$, $L_2$ and $L_1L_2 \subset L$ we have $L_1 = L$ or $L_2 = L$.*

*Let S be a semigroup. A function f from S to the unit interval [0, 1] is a fuzzy subset of S. A semigroup S itself is a fuzzy subset of S such that S(x) = 1 for all x ∈ S denoted by S. Let $\mu$ and $\delta$ be any two fuzzy subsets of S. Then the inclusion relation $\mu \subseteq \delta$ is defined by $\mu(x) \leq \delta(x)$ for all x ∈ S. $(1 - \mu)$ is a fuzzy subset of S defined for all x ∈ S.*

*$(1 - \mu)(x) = 1 - \mu(x)$, $\mu \cap \delta$ and $\mu \cup \delta$ are fuzzy subsets of S defined by $(\mu \cap \delta)(x) = \min\{\mu(x), \delta(x)\}$, $(\mu \cup \delta)(x) = \max\{\mu(x), \delta(x)\}$ for all x ∈ S. The product $\mu • \delta$ is defined as follows:*

$$(\mu • \delta)(x) = \begin{cases} \sup_{x=yz}\{min[\{\mu(y), \delta(z)\}]\} \\ 0 \text{ if x is not expressible as } x = yz \end{cases}$$

*'•' is an associative operation. We denote a fuzzy point of S by $a_\lambda$ where*

$$a_\lambda(x) = \begin{cases} \lambda & x = a \\ 0 & otherwise \end{cases}$$

For any fuzzy subset $f$ of $S$ it is obvious that $f = \bigcup_{a_\lambda \in f} a_\lambda$. Let $\lambda f_A$ be a fuzzy subset of $S$ defined as follows:

$$\lambda f_x(x) = \begin{cases} \lambda & x \in A \\ 0 & otherwise \end{cases}$$

Then we have a fuzzy point $a_\lambda$ of $S$ denoted by $\lambda f_{(a)}$.

A fuzzy subset $f$ of $S$ is called a fuzzy left ideal of $S$ if $S \circ f \subseteq f$ and $f$ is called a fuzzy right ideal of $S$ if $f \circ S \subseteq f$. If $f$ is both a fuzzy left and fuzzy right ideal then we call $f$ a fuzzy ideal of $S$. Equivalently if $f(xy) > \max\{f(x), f(y)\}$ for all $x, y \in S$.

**THEOREM 1.10.1:** *Let $f$ and $g$ be fuzzy subsets of $S$. Then the following statements are true:*

i. $f \circ (g \cup h) = (f \circ g) \cup (f \circ h)$.
ii. $f \circ (g \cap h) = (f \circ g) \cap (f \circ h)$.
iii. *If $f_1$, $f_2$ are fuzzy subsets of $S$ such that $f_1 \subseteq f$, $f_2 \subseteq g$. Then $f_1 \circ f_2 \subseteq f \circ g$.*

*Proof:* Straightforward.

**THEOREM 1.10.2:** *Let $a_\lambda$ be a fuzzy point of $S$. Then*

i. *The fuzzy left ideal generated by $a_\lambda$ denoted by $L(a_\lambda)$, is for all $x \in S$,*

$$L(a_\lambda)(x) = \begin{cases} \lambda & x \in L(a) \\ 0 & otherwise \end{cases}$$

*where $L(a)$ is a left ideal of $S$ generated by $a$.*

ii. *The fuzzy ideal generated by $a_\lambda$, denoted by $(a_\lambda)$, is for all $x \in S$,*

$$(a_\lambda)(x) = \begin{cases} \lambda & x \in (a) \\ 0 & otherwise \end{cases}$$

*where $(a)$ is an ideal of $S$ generated by $a$.*

*Proof*: Left as an exercise for the reader.

**THEOREM 1.10.3:** *Let $a_\lambda$ be a fuzzy point of $S$. The following are true:*

i. for all $x \in S$ $(S \text{ o } a_\lambda \text{ o } S)(x) = \begin{cases} \lambda & x \in SaS \\ 0 & \text{otherwise.} \end{cases}$

ii. $(a_\lambda \text{ o } b_\mu) = (ab)_{\lambda \cap \mu}$ for all fuzzy points $a_\lambda$ and $b_\mu$ of S.

iii. $(a_\lambda) = (a_\lambda \cup a_\lambda \text{ o } S \cup S \text{ o } a_\lambda \cup S \text{ o } a_\lambda \text{ o } S)$, $L(a_r) = a_r \cup S \text{ o } a_r$.

iv. $(a_\lambda)^3 \subseteq S \text{ o } a_\lambda \text{ o } S$.

*Proof:* The above theorem is easily verified by using simple calculations.

**THEOREM 1.10.4:** *Let A be a subset of S. Then for any $\lambda \in (0, 1]$ the following statements are true*

i. $\lambda f_A \text{ o } \lambda f_B = \lambda f_{AB}$
ii. $\lambda f_A \cap \lambda f_B = \lambda f_{A \cap B}$
iii. $\lambda f_A = \bigcup_{a \in A} a_\lambda$
iv. $S \text{ o } \lambda f_A = \lambda f_{SA}$
v. *If A is an ideal (right, left ideal) of S then $\lambda f_A$ is a fuzzy ideal (fuzzy right, fuzzy left ideal) of S.*

*Proof:* Easy consequence.

**DEFINITION 1.10.2:** *A fuzzy left ideal f is called prime if for any two fuzzy ideals $f_1$ and $f_2$, $f_1 \text{ o } f_2 \subseteq f$ implies that $f_1 \subseteq f$ or $f_2 \subseteq f$.*

The following theorem can be proved as a matter of routine.

**THEOREM 1.10.5:** *A fuzzy left ideal f of S is prime if and only if for any two fuzzy points $x_r, y_t \in S$ $(rt > 0)$, $x_r \text{ o } S \text{ o } y_t \text{ o } S \subseteq f$ implies that $x_r \in f$ or $y_t \in f$.*

**THEOREM 1.10.6:** *A left ideal L of S is prime if and only if $f_L$ is a prime fuzzy left ideal of S.*

*Proof:* Follows as a matter of routine. The reader who is not able to prove kindly refer [142, 143].

Now we proceed on to recall the definition of quasi prime ideals of semigroup S.

**DEFINITION 1.10.3:** *A fuzzy left ideal f is called quasi prime if for any two fuzzy left ideals $f_1$ and $f_2$; $f_1 \text{ o } f_2 \subseteq f$ implies that $f_1 \subseteq f$ and $f_2 \subseteq f$; f is called quasi semiprime if for any fuzzy left ideal g of S, $g^2 \subseteq f$ implies that $g \subseteq f$.*

Similar to the two theorems just stated we have the following two theorems.

**THEOREM 1.10.7:** *A fuzzy left ideal f of S is quasi prime if and only if for any fuzzy points $x_r, y_t \in S$ $(rt > 0)$, $x_r \text{ o } S \text{ o } y_t \subseteq f$ implies that $x_r \in f$ and $y_t \in f$.*

**THEOREM 1.10.8:** *A left ideal L of S is quasi prime if and only if $f_L$ is a quasi prime fuzzy left ideal of S.*

Now we proceed on to recall the definition of fuzzy m-systems.

**DEFINITION 1.10.4:** *A fuzzy subset f of S is called fuzzy m-system if for any t, s $\in$ [0, 1) and a, b $\in$ S, f(a) > t, f(b) > s imply that there exists an element x $\in$ S such that f(axb) > t $\vee$ s.*

**THEOREM 1.10.9:** *Let M be a subset of S. Then M is a m-system of S if and only if $f_M$ is a fuzzy m-system.*

*Proof:* For any t, s $\in$ [0, 1) and a, b $\in$ S if $f_M(a) > t$, $f_M(b) > s$, then a, b $\in$ M. By hypothesis there exists an element x $\in$ S such that axb $\in$ M that is $f_M(axb) = 1$. Thus $f_M(axb) > t \vee s$.

Conversely let a, b $\in$ M. Then $f_M(a) = f_M(b) = 1$. Thus for any t, s $\in$ [0, 1), $f_M(a) > t$, $f_M(b) > s$ which implies that there exists an element x $\in$ S such that $f_M(axb) > s \vee b$. Therefore axb $\in$ M.

**THEOREM 1.10.10:** *Let f be a fuzzy left ideal of S. Then f is quasi prime if and only if 1 – f is a fuzzy m-system.*

*Proof:* It is easily obtainable by routine calculations, hence left as an exercise for the reader.

**THEOREM 1.10.11:** *A fuzzy left ideal f is called weakly quasi prime if for any two fuzzy left ideals $f_1$ and $f_2$ such that $f \subseteq f_1$ and $f \subseteq f_2$ and $f_1 \circ f_2 \subseteq f$ then $f_1 \subseteq f$ or $f_2 \subseteq f$.*

**THEOREM 1.10.12:** *A left ideal L is weakly quasi prime if and only if f is weakly quasi prime.*

*Proof:* Please refer [143].

**THEOREM 1.10.13:** *Let S be a commutative, f a fuzzy left ideal of S. Then the following statements are equivalent.*

  i.   *f is prime.*
  ii.  *f is quasi prime.*
  iii. *f is weakly quasi prime.*

*Proof*: (i) $\Rightarrow$ (ii) $\Rightarrow$ (iii) is obvious. (iii) $\Rightarrow$ (i). Let $f_1 \circ f_2 \subset f$ for any two fuzzy ideals $f_1$ and $f_2$ of S. Since S is commutative, so we have f is a fuzzy ideal of S and

$$(f_1 \cup f) \circ (f_2 \cup f) \subseteq f_1 \circ f_2 \cup f_1 \circ f \cup f \circ f_2 \cup f^2$$
$$\subseteq f$$

**THEOREM 1.10.14:** *Let f be a fuzzy left ideal of S. Then the following statements are equivalent.*

i. *f is weakly quasi prime fuzzy left ideal of S*

ii. *For any fuzzy ideals $f_1, f_2 \subseteq S$ if $(f_1 \cup f) \circ (f_2 \cup f) \subseteq f$ then $f_1 \subseteq f$ or $f_2 \subseteq f$.*
iii. *For any two fuzzy left ideals $f_1, f_2 \subseteq S$ if $f_1 \subseteq f$ and $f_1 \circ f_2 \subseteq f$ then $f_1 = f$ or $f_2 \subseteq f$.*
iv. *For any two fuzzy ideals $f_1, f_2 \subseteq S$, if $(f_1 \cup f) \circ f_2 \subseteq f$ then $f_1 \subseteq f$ or $f_2 \subseteq f$.*
v. *For any two fuzzy points $a_r, b_t \in S$ $(rt > 0)$ if $(a_r \cup f) \circ S \circ (b_t \cup f) \subseteq f$, then $a_r \in f$ or $b_t \in S$.*

*Proof:* Left for the reader to prove.

**THEOREM 1.10.15:** *Let $f$ be fuzzy left ideal of $S$ and $\mu$ be a fuzzy subset of $S$ satisfying that*

   i. $f \cap \mu = 0$.
   ii. *For any $a_t, b_r \in \mu$ $(a_t \cup f) \circ S(b_r \cup f) \cap \mu \neq 0$.*

*If $g$ is a maximal fuzzy left ideal of $S$ with respect to containing $f$ and $g \cap \mu = 0$. Then $g$ is a weakly quasiprime fuzzy left ideal of $S$.*

*Proof:* Let $f_1$ and $f_2$ be a fuzzy left ideal of $S$ such that $g \subseteq f_1, f_2$ and $f_1 f_2 \subseteq g$. If $g \subset f_1$ and $g_2 \subset f_2$ then there exists $a_t \in f_1 \setminus g$, $b_r \in f_2 \setminus g$, $(rt > 0)$.

Thus $g \subset g \cup L(a_t) \subseteq f_1$, $g \subset g \cup L(b_r) \subseteq f_2$.

By hypothesis, we have

$(g \subset L(a_r)) \cap \mu \neq 0$, $(g \cup L(b_r)) \cap \mu \neq 0$. Let $c_k \in (g \cup L(a_t)) \cap \mu$ $(k > 0)$, $d_l \in (g \cup L(b_r)) \cap \mu$ $(l > 0)$.

Then $(c_k \cup f) \circ S \circ (d_l \cup f) \subseteq (f_1 \cup f) \circ S \circ (f_2 \cup f) \subseteq f_1 \circ S \circ f_2 \cup f_1 \cup S \circ f \cup f \circ S \circ f_2 \cup f \circ S \circ f \subseteq f_1 \circ f_2 \subseteq g$ which contradicts with the fact that $(c_k \cup f) \circ S \circ (d_l \cup f) \cap \mu \neq 0$. Hence the claim.

Fuzzy ideals in i(f) and I(f)

Let $f$ be a fuzzy left ideal of $S$, we define two fuzzy subsets of $S$, denoted by i(f) and I(f) respectively as follows ($\forall x \in S$),

$$i(f)(x) = \vee \{t_\alpha / x_{t_\alpha} \in f, x_{t_\alpha} \circ S \subseteq f, t_\alpha \in [0, 1]\}.$$
$$I(f)(x) = \vee \{t_\alpha / f \circ x_{t_\alpha} \subset f, t_\alpha \in [0, 1]\}.$$

**THEOREM 1.10.16:** *Let $f$ be a left-ideal of $S$. Then i(f) is the largest fuzzy ideal of $S$ contained in $f$.*

*Proof:* Left for the reader to prove.

**THEOREM 1.10.17:** *Let $S$ be a semigroup with an identity $e$ and $f$ a prime fuzzy left ideal of $S$. If $i(f) \neq 0$ then $i(f)$ is a quasi prime fuzzy ideal of $S$.*

*Proof:* Let $a_t$ and $b_r$ (tr > 0) be fuzzy points of S such that $a_t$ o S o $b_r \subseteq i(f)$. Then (S o $a_t$ o S) o (S o $b_r$ o S) $\subseteq i(f) \subseteq f$.

Since f is prime we have S o $a_t$ o S $\subseteq$ f or S o $b_r$ o S $\subseteq$ f. Say S o $a_t$ o S $\subseteq$ f. Using results we have S o $a_t$ o S $\subseteq i(f)$. Since S has an identity e so that $a_t = (eae)_t$, $e_t$ o $a_t$ o $e_t$ $\subseteq$ S o $a_t$ o S $\subseteq i(f)$.

The following theorem is interesting.

**THEOREM 1.10.18:** *Let f be a fuzzy ideal of S. Then I(f) is the largest fuzzy subsemigroup of S such that f is a fuzzy ideal I(f).*

*Proof:* Left as an exercise for the reader. Please refer [143].

**THEOREM 1.10.19:** *Let S be a semigroup with an identity e, f a fuzzy left ideal of S but not a fuzzy ideal of S. Then the following statements are equivalent.*

 i. *f is a weakly quasi prime fuzzy left ideal of S.*
 ii. *If $f_1$ is a fuzzy left ideal of S such that $f o f_1 \subseteq f$ then $f_1 \subseteq f$.*
 iii. *For any fuzzy point $a_t \in S$ if $f o (S o a_t) \subseteq f$, then $a_t \in f$.*
 iv. *f is the largest fuzzy left ideal of S contained in I(f).*

*Proof*: Refer [143].

Now we proceed on to define the concept of strongly semisimple semigroups.

**DEFINITION 1.10.5:** *A semigroup S is called strongly semisimple if every left ideals of S is idempotent. A left ideal f of a semigroup S is called idempotent if f = f o f that is $f^2 = f$.*

The following result gives equivalent formulation of strongly semisimple semigroup.

**THEOREM 1.10.20:** *For a semigroup S the following conditions are equivalent.*

 i. *S is strongly semisimple.*
 ii. *For any $a \in S$, $a \in SaSa$.*

*Proof:* For any $a \in S$, since S is strongly semisimple we have $L(a) = L(a)^2$. Thus $L(a) = L(a)^4$. Since $L(a)^2 = (a \cup Sa)(a \cup Sa) \subseteq Sa \cup SaSa$ then $a \in L(a) = L(a)^4 \subseteq (Sa \cup SaSa)$, $(Sa \cup SaSa) \subseteq SaSa$.

Conversely let L be any left ideal of S. For any $a \in L$ by (ii) we have $a \in SaSa \subseteq L(a)L(a) \subseteq L^2$. Thus $L \subseteq L^2$. Clearly $L^2 \subset L$. Therefore $L^2 = L$.

The following theorem characterizes the strongly semisimple semigroups by means of fuzzy quasi prime left ideals of S.

**THEOREM 1.10.21:** *Let S be a semigroup. Then the following statements are equivalent.*

  i. *Every fuzzy left fuzzy ideal of S is idempotent.*
  ii. *For any two fuzzy left ideals $f_1$ and $f_2$ of S, $f_1 \cap f_2 \subseteq f_1 \circ f_2$.*
  iii. *For any fuzzy point $a_r \in S$, $L(a_r) \subseteq L(a_r)^2$.*
  iv. *For any fuzzy point $a_r \in S$, $a_r \in S \circ a_r \circ S \circ a_r$.*
  v. *Every fuzzy left ideal of S is a quasi semiprime fuzzy left ideal of S.*
  vi. *Every fuzzy left ideal of S is the intersection of all quasi prime fuzzy left ideals of S containing it.*

*Proof:* Since the proof is long, it is left for the reader to prove. The reader is advised to refer [143].

**THEOREM 1.10.22:** *A semigroup S is strongly semisimple if and only if every fuzzy left ideal of S is idempotent.*

*Proof:* Follows easily by the definitions and the theorem.

**THEOREM 1.10.23:** *Let S be commutative. Then the fuzzy left ideals of S are quasi prime if and only if they form a chain and S is strongly semisimple.*

*Proof*: Let g and h be fuzzy left ideals of S. Since g o h is a fuzzy left ideal of S by g o h $\subset$ S o h we have g $\subseteq$ g o h $\subset$ S o h $\subseteq$ h or h $\subseteq$ g o h $\subseteq$ g o S $\subseteq$ g. Thus the fuzzy left ideals of S form a chain. Moreover for any fuzzy left ideal f of S obviously $f^2 \subseteq f$. Since f o f $\subseteq f^2$ we have f $\subseteq f^2$ so that $f^2 = f$.

Conversely let f, g be two fuzzy left ideals of S, and f o g $\subset$ h. Since the fuzzy left ideals of S form a chain i.e. f $\subseteq$ g or g $\subseteq$ f we have $f^2 \subset$ f o g $\subset$ h or $g^2 \subseteq$ f o g $\subseteq$ h. By hypothesis f $\subseteq$ h or g $\subseteq$ h holds.

Let X be an alphabet with $1 \leq |X| \leq \infty$ and $X^+(X^*)$ is the semigroup (free monoid) generated from X with the operation of adjoin. F stands for 'fuzzy', and F(x) denotes the set of all fuzzy subsets of X. We denote the cardinality of $A_{\bar{o}}$ where $A \in F(X^*)$ by card A. 1 is the identity of $X^*$ and we let M be a monoid. For x, y $\in$ M, and a given subset Y of a moniod M, we define

$$x^{-1}y = \{z \in M\ /\ xz = y\}$$
$$xy^{-1} = \{z \in M\ /\ x = zy\}$$

$F(Y) = M^{-1}YM^{-1} = \{m \in M/$ such that there exists u, v $\in$ M, umv $\in$ Y} and we use the notation $\overline{F}(Y)$ to denote the complement of F(Y) in M i.e. $\overline{F}(Y) = M.F(Y)$.

We recall some definitions from [72].

**DEFINITION 1.10.6:** *A word $w \in X^+$ is called unbordered if no proper non-empty left factor of w is a right factor of w. In other words, w is unbordered if and only if $w \in uX^+ \cap X^+u$ implies $u = 1$.*

**PROPOSITION 1.10.1:** *Let X be an alphabet with at least two letters for each word $u \in X^+$, there exists $v \in X^+$ such that uv is unbordered.*

*Proof*: Obvious hence left for the reader to prove.

**DEFINITION 1.10.7:** *Let $A, B \in F(X^*)$. For any $x \in X^*$,*

$$(A - B)(x) = \begin{cases} A(x) & \text{if } B(x) = 0 \\ 0 & \text{if } B(x) > 0 \end{cases}$$

$$(AB)(x) = \sup_{\substack{y,z \in X^* \\ yz \in X}} \min(A(y), B(z))$$

**THEOREM [27]:** *An F-subset $A \in F(X^*)$ of the free monoid $X^*$ is an F-code over X if and only if for all $n, m \geq 1$ and $x_1, ..., x_n, x'_1, ..., x'_m \in X^*$ the condition $x_1 x_2 ... x_n = x'_1 x'_2 ... x'_m = x$ implies $\min(A(x_1), ..., A(x_n), A(x'_1), ..., A(x'_m))_\alpha = \min([n = m], [x'_1, ..., x'_n] ... [x_1, ..., x_n]) \geq A^+(x)$.*

**THEOREM [72]:** *A is a F-prefix if and only if for all $x, x', u$ in $X^*$, $x' = xu$ implies $\min(A(x), A(x')) \leq [x = x']$.*

**DEFINITION [72]:** *A function $\Pi_f : F(X^*) \to R_+ \cup \{\infty\}$ is called a Bernoulli F-distribution on $X^*$ if*

$$\Pi_f(A) = \sum_{\substack{A(x) > 0 \\ x \in X^*}} \Pi(X)$$

*for any $\phi \neq A \in F(X^*)$ and $\Pi_f(\phi) = 0$ where $\Pi : X^* \to R_+$ is a Bernoulli distribution, $\Pi_f$ is positive if $\Pi$ is positive.*

Now we recall definitions from [72].

**DEFINITION 1.10.8:** *Let M be a monoid and let P be a F subset of F(M). An element $m \in M$ is called $\lambda$-completable for P if there exists $u, v$ in M such that $P(umv) > \lambda$ where $\lambda \in [0, 1]$.*

*An 0-completable element for P is called fuzzy completable for P.*

**THEOREM 1.10.24:** *If $m \in M$ is $\lambda$-completable $\Leftrightarrow MmM \cap P_{\bar{\lambda}} \neq \phi$ if and only if $m \in F(P_{\bar{\lambda}}) = M^{-1} P M_{\bar{\lambda}}^{-1}$ where $\lambda \in [0, 1]$.*

*Proof:* Obvious by the very definition.

**DEFINITION 1.10.9:** *A word which is not λ-completable for P is of course $F(P_{\bar{\lambda}})$ the set $\bar{F}(P_{\bar{\lambda}}) = M - \bar{F}(P_{\bar{\lambda}})$ of λ-completable words is a two-sided ideal of M which is disjoint from $P_\lambda$.*

**DEFINITION 1.10.10**: *For a F subset P of F(M), is called λ-dense in M if all elements of M are λ-completable for P, λ ∈ [0, 1].*

**THEOREM 1.10.25**: *P is λ-dense in M if and only if $F(P_{\bar{\lambda}}) = M$ where λ ∈ [0, 1].*

*Proof*: Left as an exercise for the reader to prove.

**DEFINITION 1.10.11**: *An F subset P of F(M) is called λ-complete in M if the F-submonoid generated by P is λ-dense where λ ∈ [0, 1]. P is called F-complete if it is 0-complete.*

**THEOREM 1.10.26**:

  i. *Every λ-dense set is also λ-complete and*
  ii. *An F-subset A of $F(X^*)$ is λ-complete if and only if $F(A^*_{\bar{\lambda}}) = X^*$.*

*Proof*: Suppose that the F-subset P of F(M) is arbitrarily λ-dense. Then for any m ∈ M, there exists u, v ∈ M such that P(umv) > λ. So P(umv) ≥ P(umv) > λ. Thus $P^*$ is λ-complete. (ii) implies, suppose that A of $F(X^*)$ is λ-complete. Then $A^*$ is λ-dense and $F(A^*_\lambda) = X^*$ by earlier results. Suppose that $F(A^*_{\bar{\lambda}}) = X^*$. Then for any m ∈ $X^*$ we have m ∈ $F(A_{\bar{\lambda}}) = X^* A_{\bar{\lambda}} X^*$, that there exists u, v ∈ $X^*$ such that $A^*(umv) > \lambda$. Thus $A^*$ is λ-dense and hence A is λ-complete.

The following theorem is interesting but the proof is left for the reader.

**THEOREM 1.10.27**: *Let A ∈ $F(X^+)$ be an F-code. Let y ∈ $X^*$ be an unordered word such that $X^* y X^* \cap X^* = \phi$. Let $U = X^* - A^* - X^* y X^*$. Then the F-set $B = A \cup y(Uy)^*$ is an F-complete code.*

**THEOREM 1.10.28**: *Any maximal F-code is F-complete.*

*Proof*: Let A ∈ $F(X^+)$ be an F-code, which is not F-complete. If card (X) = 1, then A = ϕ and A is not maximal. If card (X) ≥ 2 consider a word u ∈ $F(A^*_0)$. We know a word v ∈ $X^*$ such that y = uv is unbordered. We still have y ∈ $F(A^*_0)$ (If y ∈ $F(A^*_0)$ then there exists x, z ∈ $X^*$ such that $0 < A^*(xyz) = A^*(xu(vz))$ which contradicts the fact that u ∈ $F(A^*_0)$. It follows that A ∪ {(y, 1)} is an F-code and A is not a maximal F-code.

Now we recall the concept of λ-thin in a monoid.

**DEFINITION 1.10.12:** *An F-subset $p \in F(M)$ of a monoid M that is not $\lambda$-dense ($\lambda \in [0, 1]$) is called $\lambda$-thin. P is called an F-thin set if P is 0-thin.*

**THEOREM 1.10.29:** *If P is $\lambda$-thin then there is at least one element m in M that is $\lambda$-incompletable for P. $P \Leftrightarrow F(P_{\bar{\lambda}}) \neq M$ if and only if $MmM \cap P_{\bar{\lambda}} = \phi$.*

*Proof*: Left for the reader.

**THEOREM 1.10.30:** *Let M be a monoid and $PQ \in F(M)$. Then the $P \cup Q$ is $\lambda$-thin if and only if P and Q are $\lambda$-thin. If R is $\lambda$-dense and P is $\lambda$-thin, then $R - P$ is $\lambda$-dense.*

*Proof*: It is a matter of routine, left for the reader as an exercise.

**THEOREM 1.10.31:** *For $\lambda \in [0, 1]$, any finite F-subset of $X^*$ is clearly $\lambda$-thin. Furthermore if A and B are $\lambda$-thin subsets in $F(X^*)$ then AB is a $\lambda$-thin set.*

*Proof*: Refer [72].

**THEOREM 1.10.32:** *Let $A \in F(X^+)$ be a $\lambda$-thin and $\lambda$-complete. $\lambda \in [0, 1]$. Let w be a $\lambda$-incompletable word for A. Then*

$$X^* = \bigcup_{\substack{d \in D \\ g \in G}} d^{-1} A_{\bar{\lambda}}^* g^{-1} = D^{-1} A_{\bar{\lambda}}^* G^{-1}$$

*where G and D are the sets of right (respectively left) factors of w. Note that the set $D \times G$ is finite.*

*Proof*: Let $z \in X^*$. Since $A^*$ is $\lambda$-dense, the word wzw is $\lambda$-completable for $A^*$, thus for some $u, v \in X^*$, $A^*(uwzwt) > \lambda$. Since w is $\lambda$-incompletable for A, for any $m, n \in X^*$, $A(mwn) \leq \lambda$. Thus w is not a factor of a word in $A_{\bar{\lambda}}$.

Here there exists two factorizations $w = g_1 d = g d_1$ such that $u g_1, dzg, d_1 v \in A_{\bar{\lambda}}^*$. This shows that $z = d^{-1} A_{\bar{\lambda}}^* g^{-1}$.

**THEOREM 1.10.33:** *Let A be an F-thin and F-complete subset in $F(X^*)$. For any positive Bernoulli F-distribution $\Pi_f$ on $X^*$, we have $\Pi_f(A) \geq 1$.*

*Proof*: We have $\Pi_f(X^*) = \infty$. Note that $D \times G$ is finite, then there exists a pair $(d, g) \in D \times G$ such that $\Pi_f(d^{-1} A_0^* g^{-1}) = \Pi_f(d^{-1} A^* g^{-1}) = \infty$.

Now $d(d^{-1} A^* g^{-1})g \subset A^*$. This implies $\Pi_f(d)\Pi_f = (d^{-1} A^* g^{-1})\Pi_f(g) \leq \Pi_f(A^*)$. The positivity of $\Pi_f$ shows that $\Pi_f(dg) \neq 0$. Thus $\Pi_f(A^*) = \infty$. Now $\Pi_f(A^*) \leq \sum_{n>0} \Pi_f(A^n)$ $\leq \sum_{n>0} \Pi_f(A)^n$. Assuming $\Pi_f(A) < 1$, we get $\Pi(A^*) < \infty$. Thus $\Pi_f(A) \geq 1$.

**THEOREM 1.10.34:** *Let $A \in F(X^*)$ be an F-thin set. For any positive Bernoulli F-distribution on $X^*$ we have $\Pi_f(A) < \infty$.*

*Proof:* Let w be a word that is not a factor of a word in $A_{\bar{0}}$, $w \notin F(A_{\bar{0}})$. Set $n = |w|$, we have $n \geq 1$. For $0 \leq i \leq n-1$, consider

$$A_i(x) = \begin{cases} A(x), & \text{if } |x| \equiv i \bmod n \\ 0, & \text{otherwise} \end{cases}$$

It suffices to show that $\Pi_f(A_i)$ is finite for $i = 0, 1, 2, \ldots, n-1$. Now $A_i \subset X^i(X^n - w)^*$.

Since $X^n - w$ is an F-code, we have

$$\Pi_f[(X^n - w)^*] = \sum_{k \geq 0} \left(\Pi_j(X^n - w)\right)^k = \sum_{k \geq 0} (1 - \Pi_f(w))^k.$$

The positivity of $\Pi_f$ implies $\Pi_f(w) > 0$ and consequently $\Pi_f[(X^n - w)^*] = 1/[\Pi_f(w)]$. Thus

$$\Pi_f(A_i) \leq \Pi_f(X^i) \Pi_f[(X^n - w)^*] = \frac{1}{\Pi_f(w)}.$$

$$\Pi_f(A) = \Pi_f\left(\bigcup_{i=0}^{n-1} A_i\right) \leq \sum_{i=0}^{n-1} \Pi_f(A_i) = \frac{1}{\Pi_f(w)} < \infty.$$

**THEOREM 1.10.35**: *Let A be an F-code over X. Then A is F-complete if and only if A is F-dense or maximal.*

*Proof*: The reader is expected to prove.

**THEOREM 1.10.36:** *Let $A \in F(X^*)$ be a finite maximal F-code. For any non-empty subset Y of X, the F-code $A \cap Y^*$ is a maximal F-code over Y. In particular, for each letter $a \in X$, there is an integer n such that $A(a^n) > 0$.*

*Proof*: Refer [72].

**THEOREM 1.10.37:** *Let A be an F-thin code. The following conditions are equivalent.*

 i. *A is a maximal F-code.*
 ii. *there exists a positive Bernoulli F-distribution such that $\Pi_f(A) = 1$.*
 iii. *for any positive Bernoulli F-distribution $\Pi_f$ we have $\Pi_f(A) = 1$.*
 iv. *A is F-complete.*

*Proof:* Easily follows from the earlier results.

**THEOREM 1.10.38:** *Let A be an F-thin subset of $X^*$ and let $\Pi_f$ be a positive Bernoulli F-distribution. Any two among the three following conditions imply the third.*

    i.    *A is an F-code.*
    ii.    *$\Pi_f(A) = 1$.*
    iii.    *A is F-complete.*

*Proof*: Refer [72].

The study of fuzzy topological semigroups are not introduced or studied in this book. As the main aim of this book is fuzzy algebra and Smarandache fuzzy algebra we have restrained ourselves from this study.

## 1.11 Fuzzy Sub-half groupoids and its generalizations

In this section we recall the definition of non-closed fuzzy algebraic structures, namely fuzzy sub-halfgroupoids of a halfgroupoid, which is a generalization of a fuzzy group. Here we give the necessary and sufficient condition for the two fuzzy sub-halfgroupoids to have open fuzzy extension. Also we recall yet another new notion of the concept called anti fuzzy sub halfgroupoid and anti fuzzy extension and obtain some characterization theorems about them. We also recall the notion of normalized fuzzy extension and maximal fuzzy extension and obtain a relation between these two concepts.

A new notion called fuzzy extension chains for half groupoids is studied characterization theorems for normalized fuzzy extension and maximal fuzzy extension are obtained. Also a relationship between fuzzy translation extension and fuzzy multiplication extension is given in this section.

We prove if $\mu$ is a fuzzy sub-halfgroupoid of a halfgroupoid P and $\alpha \in [0, 1 - \sup\{\mu(x) \mid x \in P\}]$ then $\{\mu_{\alpha'}^T\}_{\alpha' \in S}$ (where $S = \{\gamma : 0 \leq \gamma \leq \alpha$ and $\gamma$ rational $\}$) is an open fuzzy extension chain for P. We prove that $\{\mu_{\alpha'}^T\}_{\alpha' \in S}$ generates the halfgroupoid P if and only if $\mu$ is a constant map.

We denote by (P, $*$) the halfgroupoid where '$*$' need not in general be closed or associative binary operation.

**DEFINITION 1.11.1:** *Let P be a halfgroupoid. A map $\mu : P \to [0, 1]$ is said to be a fuzzy sub-halfgroupoid of the halfgroupoid P if $ab = c$ in P implies $\mu(c) \geq \min\{\mu(a), \mu(b)\}$ (Here 'a b= c' in P means that a, b, c $\in$ P and a b = c).*

We illustrate this by an example.

*Example 1.11.1:*

$$P = \left\{\frac{1}{2^0}, \frac{1}{2^1}, \frac{1}{2^2}, \ldots, \frac{1}{2^n}\right\},$$

under usual multiplication as a binary operation P is a half groupoid. Define $\mu: P \to [0, 1]$ by $\mu(x) = 1 - x$ for every $x \in P$. Now, it is easy to verify that $\mu$ is a fuzzy sub-halfgroupoid of the halfgroupoid P.

**THEOREM 1.11.1:** *Let $\mu$ be a fuzzy sub-halfgroupoid of a half groupoid P and a, b, c $\in$ P. Then $\mu(c) \geq \min\{\mu(a), \mu(b)\}$ does not in general imply a b = c in P.*

*Proof:* In the above example 1.11.1 we see $\mu(c) \geq \min\{\mu(a), \mu(b)\}$ but a b $\neq$ c for

$$a = \frac{1}{2^0}, b = \frac{1}{2} \text{ and } c = \frac{1}{2^2}.$$

**DEFINITION 1.11.2:** *Let $\mu_1$ and $\mu_2$ be any two fuzzy sub-halfgroupoids of a halfgroupoid P such that*

 i. $\mu_1(x) \leq \mu_2(x)$ for every $x \in P$.
 ii. if a b = c in P and

*$\mu_1(c) = \min\{\mu_1(a), \mu_1(b)\}$ then $\mu_2(c) = \min\{\mu_2(a), \mu_2(b)\}$. Then we say that $\mu_2$ is a fuzzy extension of $\mu_1$.*

*Example 1.11.2:* Let

$$P = \left\{\frac{1}{2^0}, \frac{1}{2^1}, \frac{1}{2^2}, \ldots, \frac{1}{2^n}\right\}$$

be a halfgroupoid with respect to usual multiplication.

Define $\mu_1, \mu_2 : P \to [0, 1]$ by $\mu_1(x) = 1 - x$ for every $x \in P$ and $\mu_2(x) = 1 - \frac{x}{2}$ for every $x \in P$ and it is easy to verify that $\mu_2$ is a fuzzy extension of $\mu_1$.

**THEOREM 1.11.2:** *Let $\mu_1$ and $\mu_2$ be any two fuzzy sub-halfgroupoid of a halfgroupoid P. If $\mu_2$ is a fuzzy extension of $\mu_1$ and a b = c in P with $\mu_2(c) = \min\{\mu_2(a), \mu_2(b)\}$ then it need not in general imply $\mu_1(c) = \min\{\mu_1(a), \mu_1(b)\}$.*

*Proof:* The proof is by an example. Take P and $\mu_1$ as in example 1.11.2. Define $\mu_2: P \to [0, 1]$ by $\mu_2(x) = 1$ for every $x \in P$. Now if we take

$$a = \frac{1}{2^0} \text{ and } b = \frac{1}{2} \text{ then } ab = \frac{1}{2^0}\frac{1}{2} = \frac{1}{2}$$

It is easy to verify that $\mu_2(c) = 1$, $\min\{\mu_2(a), \mu_2(b)\} = 1$. $\mu_1(c) = \frac{1}{2}$ and $\min\{\mu_1(a), \mu_1(b)\} = 0$. This proves that $\mu_2(c) = \min\{\mu_2(a), \mu_2(b)\}$ and $\mu_1(c) \neq \min\{\mu_1(a), \mu_1(b)\}$. This completes the proof.

**DEFINITION 1.11.3:** *Let $\mu_1$ and $\mu_2$ be any two fuzzy sub-halfgroupoids of a half groupoid P. Then the fuzzy extension $\mu_2$ (of $\mu_1$) is said to be an open fuzzy extension of $\mu$ if a b = c in P and $\mu_1$ (c) > min {$\mu_1$ (a), $\mu_1$ (b)} implies $\mu_2$ (c) > min {$\mu_2$ (a), $\mu_2$ (b)}. $\mu_2$ is said to be closed fuzzy extension of $\mu_1$ if it is not an open fuzzy extension of $\mu_1$.*

The following example illustrates it.

*Example 1.11.3:* Consider the halfgroupoid P = {– n, – (n – 1), …, –1, 0, 1, …, n} under the usual multiplication. Define $\mu_1, \mu_2 : P \to [0, 1]$ by

$$\mu_1(x) = \begin{cases} \dfrac{1}{2} & \text{if } x > 0 \\ \dfrac{1}{4} & \text{if } x \leq 0 \end{cases}$$

and

$$\mu_2(x) = \begin{cases} 1 & \text{if } x > 0 \\ \dfrac{3}{4} & \text{if } x \leq 0 \end{cases}$$

It is easy to check that $\mu_2$ is an open fuzzy extension of $\mu_1$

**THEOREM 1.11.3:** *Every open fuzzy extension is a fuzzy extension but not conversely.*

*Proof:* Every open fuzzy extension is a fuzzy extension, directly follows from the definitions of open fuzzy extension and fuzzy extension.

Converse is not true. This is explicit by the following example.

*Example 1.11.4:* Let P be any halfgroupoid and $\mu_1$ be any fuzzy sub-halfgroupoid of the halfgroupoid P with o (Im ($\mu_1$)) ≥ 2. Then it is easy to verify that $\mu_2 : P \to [0, 1]$ given by $\mu_2$ (x) = 1 for every x ∈ P is a fuzzy extension of $\mu_1$. It is easily seen $\mu_2$ is not an open fuzzy extension of $\mu_1$.

**DEFINITION 1.11.4:** *Let P be a half groupoid. A countable collection of fuzzy sub halfgroupoids of the halfgroupoids P, denoted by $\{\mu_i \mid i = 0, 1, 2, 3, \ldots\}$ or $\{\mu_i\}_{i=0}^{\infty}$ is called fuzzy extension chain for P if $\mu_{i+1}$ is a fuzzy extension of $\mu_i$ for i = 0, 1, 2, … .*

*A fuzzy extension chain $\{\mu_i\}_{i=0}^{\infty}$ is said to generate P if $\bigcup_{i=0}^{\infty} \mu_i = 1_p$ (Here $1_P$ denotes the map $\mu : P \to [0, 1]$ such that $\mu$ (x) = 1 for every x ∈ P).*

*A fuzzy extension chain $\{\mu_i\}_{i=0}^{\infty}$ is said to be a fuzzy extension chain for P if $\mu_{i+1}$ is an open fuzzy extension of $\mu_i$ for i = 0, 1, 2, … .*

*Example 1.11.5:* Let P be any halfgroupoid and $\mu_1: P \to [0, 1]$ be defined by $\mu_i(x) = \dfrac{i}{i+1}$ for every $x \in P$ and $i = 0, 1, 2, \ldots$ .

Then we get a countable collection of fuzzy sub-halfgroupoids $\{\mu_i\}_{i=0}^{\infty}$ of the halfgroupoid P. If $ab = c$ in P and $\mu_{i+1}(c), = \min\{\mu_{i+1}(a), \mu_{i+1}(b)\}$ then we have $\mu_i(c) = \min\{\mu_i(a), \mu_i(b)\}$. Hence $\mu_{i+1}$ is a fuzzy extension of $\mu_i$ for all $i = 0, 1, 2, \ldots$.

This establishes that $\{\mu_i\}_{i=0}^{\infty}$ is a fuzzy extension chain for P.

Now we will prove $\{\mu_i\}_{i=0}^{\infty}$ generates the halfgroupoid P. Let x be an arbitrary element of the halfgroupoid P then we have

$$\left(\bigcup_{i=0}^{\infty}\mu_i\right)(x) = \sup_{i=0,1,2,\ldots}\{\mu_i(x)\}$$

$$= \sup_{i=0,1,2,\ldots}\left\{\dfrac{i}{i+1}\right\}$$

$$= 1 = 1_{P(x)}$$

That is $\left(\bigcup_{i=0}^{\infty}\mu_i\right)(x) = 1_{P(x)}$ for every $x \in P$. Thus we have proved that $\bigcup_{i=0}^{\infty}\mu_i = 1_p$. This implies that $\{\mu_i\}_{i=0}^{\infty}$ generates the half groupoid P.

**THEOREM 1.11.4:** *Let $\mu_1$ and $\mu_2$ be any two fuzzy sub half groupoids of a half groupoid P. If the fuzzy extension $\mu_2$ (of $\mu_1$) is an fuzzy extension of $\mu_1$ then $\mu_1(c) = \min\{\mu_1(a), \mu_1(b)\}$ if and only if $\mu_2(c) = \min\{\mu_2(a), \mu_2(b)\}$.*

We prove the converse of the contra positive method. Let $ab = c$ in P and $\mu_1(c) \neq \min\{\mu_1(a), \mu_2(b)\}$ to prove $\mu_2(c) \neq \min\{\mu_2(a), \mu_1(b)\}$.

If $ab = c$ in P and $\mu_1(c) \neq \min\{\mu_1(a), \mu_1(b)\}$, then we have $\mu_1(c) > \min\{\mu_1(a), \mu_1(b)\}$ (since $\mu_1$ is a fuzzy subhalf-groupoid of the half groupoid P).

This implies that $\mu_2(c) > \min\{\mu_2(a), \mu_2(b)\}$ (since $\mu_2$ is an open fuzzy extension of $\mu_1$). That is $\mu_2(c) \neq \min\{\mu_2(a), \mu_2(b)\}$. Hence we have proved that if $ab = c$ in P and $\mu_1(c) \neq \min\{\mu_1(a), \mu_1(b)\}$ then $\mu_2(c) \neq \min\{\mu_2(a), \mu_2(b)\}$. This proves that if the fuzzy extension $\mu_2$ (of $\mu_1$) is an open fuzzy extension of $\mu_1$ then $\mu_1(c) = \min\{\mu_1(a), \mu_1(b)\}$ if and only if $\mu_2(c) = \min\{\mu_2(a), \mu_2(b)\}$ for every $ab = c$ in P.

We illustrate the above theorem by an example.

*Example 1.11.6:* Take $\mu_1$ and $\mu_2$ given by the half groupoid $P = \{-n_1, -(n-1), \ldots, -1, 0, 1, \ldots, n\}$ then it is easy to verify that if $\mu_2$ is an open fuzzy extension of $\mu_1$ then $\mu_1(c) = \min\{\mu_1(a), \mu_1(b)\}$ if and only if $\mu_2(c) = \min\{\mu_2(a), \mu_2(b)\}$ for every $ab = c$ in P.

Since there are not many books on fuzzy algebra and especially no paper even on fuzzy half groupoids we are forced to give all definitions and results about them.

Now we proceed on to give the definition of anti fuzzy extension.

**DEFINITION 1.11.5:** *A fuzzy subset $\mu$ of a half groupoid P is said to be an anti-fuzzy subhalfgroupoid of the half groupoid P if ab = c in P implies $\mu(c) \leq min\{\mu(a), \mu(b)\}$.*

**DEFINITION 1.11.6:** *Let $\mu_1$ and $\mu_2$ be any two anti fuzzy subhalf-groupoids of a half groupoid P. Then $\mu_2$ is said to be an anti fuzzy extension of $\mu_1$ if the following two conditions hold good.*

 i. *$\mu_1(x) \geq \mu_2(x)$ for every $x \in P$.*
 ii. *if ab = c in P and $\mu_1(c) = max\{\mu_1(a), \mu_1(b)\}$ then $\mu_2(c) = max\{\mu_2(a), \mu_2(b)\}$.*

As in case of fuzzy extension we can prove that if $\mu_2$ is an anti fuzzy extension of $\mu_1$ and ab = c in P with $\mu_2(c) = max\{\mu_2(a), \mu_2(b)\}$ then it need not in general imply that $\mu_1(c) = max\{\mu_1(a), \mu_1(b)\}$.

**DEFINITION 1.11.7:** *Let $\mu_1$ and $\mu_2$ be any two anti fuzzy subhalf-groupoids of a half groupoid P. Then the anti fuzzy extension $\mu_2$ (of $\mu_1$) is said to be an open anti fuzzy extension of $\mu_1$ if a b = c in P and $\mu_1(c) < max\{\mu_1(a), \mu_1(b)\}$ implies $\mu_2(c) < max\{\mu_2(a), \mu_2(b)\}$.*

If $\mu_2$ is not an open anti fuzzy extension of $\mu_1$ then we say that $\mu_2$ is a closed anti fuzzy extension of $\mu_1$.

**DEFINITION 1.11.8:** *Let P be a half groupoid. A countable collection of anti fuzzy sub half groupoids of the half groupoid P, denoted by $\{\mu_i \mid i = 0, 1, 2, \ldots\}$ or $\{\mu_i\}_{i=0}^{\infty}$ is called an anti fuzzy extension chain for P if $\mu_{i+1}$ is an anti fuzzy extension of $\mu_i$ for i = 0, 1, 2, ….*

We give an example of an anti fuzzy sub-half groupoid and anti fuzzy extension.

*Example 1.11.7:* Let

$$P = \left\{\frac{1}{2^0}, \frac{1}{2^1}, \cdots, \frac{1}{2^n}\right\}$$

take usual multiplication as the binary composition in P. Define $\mu: P \to [0, 1]$ by $\mu(x) = x$ for every $x \in P_1$. It is clear that $\mu$ is an anti fuzzy sub half groupoid of the half groupoid P.

Define $\mu_1, \mu_2 : P \to [0, 1]$ by

$\mu_1(x) = x$ for every $x \in P$ and

$\mu_2(x) = \dfrac{x}{2}$ for every $x \in P$

Then it is easy to verify that $\mu_2$ is an anti fuzzy extension of $\mu_1$.

**THEOREM 1.11.5:** *Let P be a half groupoid. Then $\mu$ is a fuzzy subhalf-groupoid of the half groupoid P if and only if $\mu^c$ is an anti fuzzy subhalf-groupoid of the half groupoid P.*

*Proof:* Let $\mu$ be a fuzzy subhalf-groupoid of a half groupoid P and ab = c in P. Since $\mu^c(x) = 1 - \mu(x)$ for every x ∈ P, we have $\mu(c) \geq \min\{\mu(a), \mu(b)\}$ if and only if $\mu^c(c) \leq \max\{\mu^c(a), \mu^c(b)\}$. This proves that $\mu$ is a fuzzy subhalf-groupoid of the half groupoid P if and only if $\mu^c$ is an anti fuzzy sub-half groupoid of the half groupoid P.

**THEOREM 1.11.6:** *Let $\mu_1$ and $\mu_2$ be any two fuzzy subhalf-groupoids of a half groupoid P. Then*

   i. *$\mu_2$ is a fuzzy extension of $\mu_1$ if and only if $\mu_2^c$ is an anti fuzzy extension of $\mu_1^c$ and*
   ii. *$\mu_2$ an open fuzzy extension of $\mu_1$ if and only if $\mu_2^c$ is an open anti fuzzy extension of $\mu_1^c$.*

*Proof:* Let $\mu_1$ and $\mu_2$ be any two fuzzy sub-half groupoids of a half groupoid P.

*Proof of (i):* Let $\mu_2$ be a fuzzy extension of $\mu_1$ to prove $\mu_2^c$ is an anti fuzzy extension of $\mu_1^c$. Since $\mu_2$ is a fuzzy extension of $\mu_1$, we have $\mu_1(x) \leq \mu_2(x)$ for every x ∈ P, that is $\mu_1^c(x) \geq \mu_2^c$ for every x ∈ P.

Let ab = c in P and $\mu_1^c(c) = \max\{\mu_1^c(a), \mu_1^c(b)\}$. That is $1 - \mu_1(c) = \max\{1 - \mu_1(a), 1 - \mu_1(b)\}$ (as $\mu_1^c(x) = 1 - \mu_1(x)$ for every x ∈ P). By the properties of min and max functions we have $1 - \mu_1(c) = 1 - \min\{\mu_1(a), \mu_1(b)\}$. Hence $\mu_1(c) = \min\{\mu_1(a), \mu_1(b)\}$. Since $\mu_2$ is a fuzzy extension of $\mu_1$ and ab = c in P we have $\mu_2(c) = \min\{\mu_2(a), \mu_2(b)\}$. Using the properties of min and max function we have $1 - \mu_2(c) = \max\{1 - \mu_2(a), 1 - \mu_2(b)\}$. That is $\mu_2^c(c) = \max\{\mu_2^c(a), \mu_2^c(b)\}$. Thus $\mu_2^c$ is an anti fuzzy extension $\mu_1^c$. Similarly we can prove the converse using the fact that $1 - \max\{a, b\} = \min\{1 - a, 1 - b\}$ for all a, b ∈ [0, 1].

*Proof of (ii):* Let $\mu_2$ be an open fuzzy extension of $\mu_1$. To prove $\mu_2^c$ is an open anti fuzzy extension of $\mu_1^c$. Since $\mu_2$ is an open fuzzy extension of $\mu_1$ we have $\mu_1(x) \leq \mu_2(x)$ for every x ∈ P that is $\mu_1^c(x) \geq \mu_2^c(x)$ for every x ∈ P.

Let ab = c in P and $\mu_1^c(c) < \max\{\mu_1^c(a), \mu_1^c(b)\}$ then $1 - \mu_1(c) < \max\{1 - \mu_1(a), 1 - \mu_1(b)\}$. (Since $\mu_1^c(x) = 1 - \mu_1(x)$ for every x ∈ P). This implies $1 - \mu_1(c) < 1 - \min\{\mu_1(a), \mu_1(b)\}$. That is $\mu_1(c) > \min\{\mu_1(a), \mu_1(b)\}$.

As $\mu_2$ is an open fuzzy extension of $\mu_1$ and ab = c in P. We have $\mu_2(c) > \min\{\mu_2(a), \mu_2(b)\}$. By the properties of min and max functions we have $1 - \mu_2(c) < \max\{1 - $ $\mu_2(a), 1 - \mu_2(b)\}$. That is $\mu_2^c(c) < \max\{\mu_2^c(a), \mu_2^c(b)\}$. Thus $\mu_2^c$ is an open anti fuzzy extension of $\mu_1^c$. Conversely let $\mu_2^c$ be an open anti fuzzy extension of $\mu_1^c$ then we prove that $\mu_2$ is an open fuzzy extensions of $\mu_1$ using the fact $1 - \max\{a, b\} = \min\{1 - a, 1 - b\}$ for all a, b ∈ [0, 1].

***Example 1.11.8:*** Choose $\mu_1$ and $\mu_2$ as in the earlier example where

$$P = \left\{ \frac{1}{2^0}, \frac{1}{2^1}, \frac{1}{2^2}, \ldots, \frac{1}{2^n} \right\}$$

$\mu_1(x) = 1 - x$ for every $x \in P$ and $\mu_2(x) = 1 - \frac{x}{2}$ for every $x \in P$. It is easy to verify that $\mu_2$ is a fuzzy extension of $\mu_1$. Now we calculate $\mu_1^c$ and $\mu_2^c$ as $\mu_1^c(x) = x$ for every $x \in P$ and $\mu_2^c(x) = \frac{x}{2}$ for every $x \in P$.

Clearly $\mu_2^c$ is an anti fuzzy extension of $\mu_1^c$. This shows that $\mu_2$ is a fuzzy extension of $\mu_1$ if and only if $\mu_2^c$ is an anti fuzzy extension of $\mu_1^c$.

Now we proceed on to define the concept of maximal fuzzy extension.

**DEFINITION 1.11.9:** *Let $\mu_1$ and $\mu_2$ be any two fuzzy sub-half groupoids of a half groupoid P. A fuzzy extension $\mu_2$ (of $\mu_1$) is said to be the maximal fuzzy extension of $\mu_1$ if there does not exist any fuzzy sub half groupoid $\mu$ of the half groupoid P such that $\mu_2 \subset \mu$.*

***Example 1.11.9:*** P = {1, 2, 3, …, n} be a half groupoid with respect to usual multiplication. Define $\mu_1(x) = ½$ for every $x \in P$ and $\mu_2(x) = 1$ for every $x \in P$. It is easy to verify that both $\mu_1$ and $\mu_2$ are fuzzy subhalf groupoids of the groupoid P. Now we cannot find any fuzzy sub half groupoid $\mu$ of the half groupoid P such that $\mu_2 \subset \mu$. Hence $\mu_2$ is the maximal fuzzy extension of $\mu_1$.

**DEFINITION 1.11.10:** *Let $\mu_1$ and $\mu_2$ be any two fuzzy sub-half groupoids of a half groupoid P. A fuzzy extension $\mu_2$ (of $\mu_1$) is said to be a normalized fuzzy extension of $\mu_1$ if $\mu_2(x) = 1$ for some $x \in P$.*

In the above example the fuzzy subhalf-groupoid $\mu_2$ is a normalized fuzzy extension of $\mu_1$.

**THEOREM 1.11.7:** *Let $\mu_1$ and $\mu_2$ be any two fuzzy sub-half groupoids of a half groupoid P. If $\mu_2$ is the maximal fuzzy extension of $\mu_1$ then $\mu_2$ is a normalized fuzzy extension of $\mu_1$.*
*The converse is not true in general.*

*Proof:* The first part of the theorem follows from the definitions of the maximal and normalized fuzzy extensions. The converse is not true. That is every normalized fuzzy extension need not in general be a maximal fuzzy extension which can be seen from the example in which P = {– n, – (n – 1), …, –1, 0, 1, 2, …, n},

$$\mu_1(x) = \begin{cases} \dfrac{1}{2} & \text{if } x > 0 \\ \dfrac{1}{4} & \text{if } x \leq 0 \end{cases}$$

and

$$\mu_2(x) = \begin{cases} 1 & \text{if } x > 0 \\ \dfrac{3}{4} & \text{if } x \leq 0. \end{cases}$$

Now we proceed on to define fuzzy translations and fuzzy multiplications and its relation with anti fuzzy extension is obtained.

**DEFINITION 1.11.11:** *Let $\mu$ be a fuzzy subset of a set X and $\alpha \in [0, 1 - \sup \{\mu(x) : x \in X\}]$. A map $\mu_\alpha^T : X \to [0,1]$ is called a fuzzy translation of $\mu$ if $\mu_\alpha^T(x) = \mu(x) + \alpha$ for every $x \in X$.*

We write this fuzzy translation of $\mu$ by $\mu_\alpha^T$, we illustrate this by an example.

***Example 1.11.10:*** Let $X = \{2, 4, \ldots, 2n\}$ be a set. Define $\mu : X \to [0, 1]$ by $\mu(x) = \dfrac{1}{x}$ for every $x \in X$. Then for $\alpha = \dfrac{1}{4}$ we have $\mu_\alpha^T(x) = \dfrac{1}{x} + \dfrac{1}{4}$ for every $x \in X$. It is easy to verify that $\mu_\alpha^T$ is a fuzzy translation of $\mu$.

**THEOREM 1.11.8:** *Let $\mu$ be a fuzzy sub-half groupoid of a half groupoid P and $\alpha \in [0, 1 - \sup\{\mu(x) ; x \in P\}]$ then every fuzzy translation $\mu_\alpha^T$ of $\mu$ is a fuzzy sub-half groupoid of the half groupoid P.*

*Proof:* Let $\mu$ be a fuzzy sub-half groupoid of the half groupoid P and $\alpha \in [0, 1 - \sup \{\mu(x) ; x \in P\}]$. For ab = c in P we have $\mu(c) \geq \max \{\mu(a), \mu(b)\}$ that is $\alpha + \mu(c) \geq \alpha + \min\{\mu(a), \mu(b)\}$ that is $\alpha + \mu(c) \geq \min\{\alpha + \mu(a), \alpha + \mu(b)\}$ (by min property) that is $\{\mu_\alpha^T(c) \geq \min\{\mu_\alpha^T(a), \mu_\alpha^T(b)\}$ (by the definition of fuzzy translation).

That is $\mu_\alpha^T$ is a fuzzy subhalf-groupoid of the half-groupoid P.

**THEOREM 1.11.9:** *Let $\mu$ be a fuzzy sub-half groupoid of the half groupoid P and $\alpha \in [0, 1 - \sup \{\mu(x); x \in P\}]$. Then every fuzzy translation $\mu_\alpha^T$ (of $\mu$) is a fuzzy extension of $\mu$.*

*Proof:* Let µ be a fuzzy sub-half groupoid of the half groupoid P and $\alpha \in [0, 1 - \sup\{\mu(x) ; x \in P\}]$. We know by the earlier result $\mu_\alpha^T$ is a fuzzy sub-half groupoid of the half groupoid P.

Clearly $\mu_\alpha^T(x) \geq \mu(x)$ for every $x \in P$. Further if $ab = c$ in P and $\mu(c) = \min\{\mu(a), \mu(b)\}$ then $\mu_\alpha^T(c) = \min\{\mu_\alpha^T(a), \mu_\alpha^T(b)\}$. Hence $\mu_\alpha^T$ is a fuzzy extension of µ.

*Example 1.11.11:* Let

$$P = \left\{\frac{1}{3^0}, \frac{1}{3^1}, \cdots, \frac{1}{3^n}\right\}$$

be a half groupoid with respect to usual multiplication. Define $\mu: P \to [0, 1]$ by $\mu(x) = 1 - x$ for every $x \in P$. Now it can be easily checked that µ is a fuzzy sub-half groupoid of the half groupoid P. For $\alpha = \frac{1}{3}$ we have $\mu_\alpha^T(x) = 1 - x + \frac{1}{3}$ for every $x \in P$. It is easy to verify that $\mu_\alpha^T$ is a fuzzy extension of µ.

**THEOREM 1.11.10:** *Let µ be a fuzzy sub half groupoid of a half groupoid P and $\alpha \in [0, 1 - \sup\{\mu(x)/ x \in P\}]$. If $\mu_\alpha^T$ is a fuzzy translation of µ then*

   i. *$\mu_\alpha^T$ is a maximal fuzzy extension of µ if and only if µ is a constant map.*
   ii. *if µ has sup property then $\mu_\alpha^T$ is a normalized fuzzy extension of µ but not conversely.*

*Proof:* Let µ be a fuzzy sub-half groupoid P and $\alpha \in [0, 1 - \sup\{\mu(x) \mid x \in P\}]$. Then by earlier result $\mu_\alpha^T$ is a fuzzy extension of µ.

*Proof of (i):* Let $\mu_\alpha^T$ be a maximal fuzzy extension of µ then $\mu_\alpha^T(x) = 1$ for every $x \in P$, that is $\mu(x) = \sup\{\mu(y) \mid y \in P\}$ for every $x \in P$. Hence µ is a constant map. Conversely let µ be a constant map, that is $\mu(x) = \gamma$ for every $x \in P$ (where γ is a fixed element of [0, 1]. Now consider the fuzzy translation $\mu_\alpha^T$ of µ $\mu_\alpha^T(x) = \mu(x) + \alpha = \gamma + 1 - \gamma = 1$ for every $x \in P$. Hence $\mu_\alpha^T$ is a maximal fuzzy extension of µ.

*Proof of (ii):* Let µ have sup property; that is for every subset A of P there exists $x_0 \in A$ such that $\mu(x_0) = \sup\{\mu(y) \mid y \in A\}$. If possible let as assume that $\mu_\alpha^T$ is not a normalized fuzzy extension of µ then $\mu_\alpha^T(x) < 1$ for every $x \in P$. That is $\mu(x) < \sup\{\mu(y) \mid y \in P\}$ for $x \in P$. Since µ has sup property for every subset A of P there exists $x_0 \in A$ such that $\mu(x_0) = \sup\{\mu(y) \mid y \in A\}$. So we get $\mu(x) < \mu(x_0)$ for every $x \in P$. In particular for $x = x_0$ we have $\mu(x) < \mu(x_0)$ which is a contradiction. Hence $\mu_\alpha^T$ is a normalized fuzzy extension of µ.

If $\mu_\alpha^T$ is a normalized fuzzy extension of $\mu$ then it does not in general imply that $\mu$ has sup property. This can be seen from the example given below.

Consider the half groupoid $P = [0, 1]$ under usual multiplication. Define $\mu : P \to [0, 1]$ by $\mu(x) = 1 - x$ for every $x \in P$. It is easy to see that $\mu$ is a fuzzy sub-half groupoid of the half groupoid P. Let $\alpha \in [0, 1 - \sup \{\mu(x) \mid x \in P\}]$ then $\mu_\alpha^T$ is a normalized fuzzy extension of $\mu$ for there exists $0 \in P$ such that $\mu_\alpha^T (0) = 1$. Take $A = [0, 1]$ then $A \subseteq [0, 1]$ and $\sup \{\mu(x) \mid x \in A\} = 1$.

Now it is easy to verify that there is no $x_0 \in A$ such that $\mu(x_0) = 1$. Hence $\mu$ does not have sup property.

**DEFINITION 1.11.12:** *A fuzzy subset $\mu$ of a set X has the weak sup property if there exist $x_0 \in X$ such that $\mu(x_0) = \sup \{\mu(x) \mid x \in X\}$.*

**THEOREM 1.11.11:** *Let $\mu$ be a fuzzy subset of a set X. If $\mu$ has sup property then $\mu$ has weak sup property but not conversely.*

*Proof:* Let $\mu$ be a fuzzy subset of a set X. If $\mu$ has sup property then by the definition of sup property for every subset A of the set X there exist $x_0 \in A$ such that $\mu(x_0) = \sup \{\mu(y) \mid y \in A\}$. The above condition is true for every subset A of the set X. If we take A as X then there exists $x_0 \in A$ such that $\mu(x_0) = \sup \{\mu(y) \mid y \in X\}$. Hence $\mu$ has weak sup property.

However the converse is not true this can be seen by the following example.

Consider the set $X = [\gamma, \delta]$ where $\gamma$ and $\delta$ are any two arbitrary fixed numbers in the interval $[0, 1]$ with $\gamma < \delta$. Define $\mu: [\gamma, \delta] \to [0, 1]$ by $\mu(x) = x$ for every $x \in X$. Then $\mu(\delta) = \sup \{\mu(x) \mid x \in X\}$. but for $A = [\gamma, \delta)$ we have $\sup \{\mu(x) \mid x \in A\} = \delta$ and there is no $x_0 \in A$ such that $\mu(x_0) = \delta$. That is $\mu$ has weak sup property but $\mu$ does not have sup property.

Hence we give the characterization for normalized fuzzy extension.

**THEOREM 1.11.12:** *Let $\mu$ be a fuzzy sub-half groupoid of a half groupoid P of $\alpha \in [0, 1 - \sup \{\mu(x) \mid x \in P\}]$. Then $\mu_\alpha^T$ is a normalized fuzzy extension of $\mu$ if and only if $\mu$ has weak sup property.*

*Proof:* Let $\mu$ be a fuzzy sub-half groupoid of a half groupoid P and $\alpha \in [0, 1 - \sup \{\mu(x) \mid x \in P\}]$. Let $\mu_\alpha^T$ be a normalized fuzzy extension of $\mu$ then by the definition of normalized fuzzy extension we have $\mu_\alpha^T(x_0) = 1$ for some $x_0 \in P$. That is $\mu(x_0) = \sup \{\mu(x) \mid x \in P\}$ for some $x_0 \in P$. Hence $\mu$ has weak sup property.

Conversely let $\mu$ have weak sup property. Then there exists $x_0 \in P$ such that $\mu(x_0) = \{\mu(x) \mid x \in P\}$.

Now consider the fuzzy translation $\mu_\alpha^T$ of $\mu$. $\mu_\alpha^T(x) = \mu(x) + 1 - \sup\{\mu(x) \mid x \in P\}$. Take $x = x_0$, then we have $\mu_\alpha^T(x_0) = \mu(x_0) + 1 - \mu(x_0) = 1$. Hence $\mu_\alpha^T$ is a normalized fuzzy extension of $\mu$. This completes the proof of the theorem.

**DEFINITION 1.11.13:** *Let $\mu$ be a fuzzy subset of set $X$ and $\beta \in [0, 1]$. A map $\mu_\beta^m : X \to [0, 1]$ is called a fuzzy multiplication of $\mu$ if $\mu_\beta^m(x) = \beta \bullet \mu(x)$ for every $x \in X$. We denote this fuzzy multiplication of $\mu$ by $\mu_\beta^m$.*

We illustrate this by an example.

***Example 1.11.12:*** Let $X = [\gamma, \delta]$ be a set where $\gamma$ and $\delta$ are any two arbitrary fixed numbers in the interval $[0, 1]$ with $\gamma < \delta$. Define $\mu : [\gamma, \delta] \to [0, 1]$ by $\mu(x) = x$ for every $x \in X$. Then for $\beta = \dfrac{1}{2}$ we calculate the fuzzy multiplication of $\mu$ and is given by $\mu_\beta^m(x) = \dfrac{x}{2}$ for every $x \in X$.

**THEOREM 1.11.13:** *Let $\mu$ be an anti fuzzy sub half groupoid of a half groupoid $P$. Then every fuzzy multiplication $\mu_\beta^m$ (of $\mu$) is an anti fuzzy extension of $\mu$.*

*Proof*: Can be proved by similar ways using the techniques used in the earlier theorems.

**THEOREM 1.11.14:** *Let $\mu$ be a fuzzy sub-halfgroupoid of a half group $P$. If $\alpha \in [0, 1 - \sup\{\mu(x) \mid x \in P\}]$ and $\beta \in (0, 1]$ then every fuzzy translation $\mu_\alpha^T$ (of $\mu$) is a fuzzy extension of fuzzy multiplication $\mu_\beta^m$ (of $\mu$).*

*Proof*: Let $\mu$ be a fuzzy sub-halfgroupoid of a half groupoid $P$. If $\alpha \in [0, 1 - \sup\{\mu(x) \mid x \in P\}]$ and $\beta \in (0, 1]$ then fuzzy translation $\mu_\alpha^T$ and the fuzzy multiplication $\mu_\beta^m$ are fuzzy sub-halfgroupoids of the half groupoid $P$.

We see that for every $x \in P$, $\mu_\beta^m(x) = \beta \mu(x) \leq \mu(x) \leq \alpha + \mu(x) = \mu_\alpha^T(x)$. Also if $ab = c$ in $P$ and $\mu_\beta^m(c) = \min\{\mu_\beta^m(a), \mu_\beta^m(b)\}$ then this implies $\beta \bullet \mu(c) = \beta \bullet \min\{\mu(a), \mu(b)\}$ that is $\mu(c) = \min\{\mu(a), \mu(b)\}$ (since $\beta \neq 0$) that is $\alpha + \mu(c) = \min\{\alpha + \mu(a), \alpha + \mu(b)\}$ (by min property) that is $\mu_\alpha^T(c) = \min\{\mu_\alpha^T(a), \mu_\alpha^T(b)\}$ (by the definition of fuzzy translation). Hence $\mu_\alpha^T$ is a fuzzy extension of $\mu_\beta^m$ for every $\alpha \in [0, 1 - \sup\{\mu(x) \mid x \in P\}]$ and $\beta \in (0, 1]$. Hence the theorem. It is to be noted that for $\beta = 0$ the result is not true.

***Example 1.11.13:*** Consider the half groupoid $P = \{-n, -(n-1), \ldots, -1, 0, 1, \ldots, n\}$ under the usual multiplication. Define $\mu : P \to [0, 1]$ by

$$\mu(x) = \begin{cases} \dfrac{1}{2} & \text{if } x > 0 \\ 0 & \text{if } x \leq 0. \end{cases}$$ If $ab = c$

in P with $a < 0$ and $b < 0$ then $c > 0$. Taking $\beta = 0$ we see that $\mu_0^m(c) = 0 = \min\{\mu_0^m(a), \mu_0^m(b)\}$ and $\mu_\alpha^T(c) = 1 > \min\{\mu_\alpha^T(a), \mu_\alpha^T(b)\} = \dfrac{1}{2}$. That is $\mu_\alpha^T$ is not a fuzzy extension of $\mu_\beta^m$ for $\beta = 0$.

**THEOREM 1.11.15:** *Let $\mu$ be a fuzzy sub half groupoid of a half groupoid P and $\alpha \in [0, 1 - \sup\{\mu(x) \mid x \in P\}]$. Then $\{\mu_{\alpha'}^T\}_{\alpha' \in S}$ where (S = $\{\gamma, 0 \leq \gamma \leq \alpha,$ and $\gamma$ rational$\}$) is an open fuzzy extension chain for P. Further $\{\mu_{\alpha'}^T\}_{\alpha' \in S}$ generates the half groupoid P if and only if $\mu$ is a constant map.*

*Proof:* It is easy to see from the proof of the earlier theorems that for every $\alpha' \in S$, $\mu_\alpha^T$, is a fuzzy sub-half groupoid of the half groupoid P.

Choose $\alpha_1, \alpha_2 \in S$ such that $0 \leq \alpha_1 \leq \alpha_2 \leq \alpha$ then $\mu_{\alpha_1}^T(x) \leq \mu_{\alpha_2}^T(x)$ for every $x \in P$. Thus for $\alpha_i, \alpha_j \in S$ ($i \leq j$), $\mu_{\alpha_i}^T(x) \leq \mu_{\alpha_j}^T(x)$ for every $x \in P$ without loss of generality we have $0 \leq \alpha_0 \leq \alpha_1 \leq \alpha_2 \leq \ldots \leq \alpha_i \leq \alpha_{i+1} \leq \ldots \leq \alpha$. Then by the above construction $\mu_{\alpha_0}^T(x) \leq \mu_{\alpha_j}^T(x) \leq \cdots \leq \mu_{\alpha_i}^T(x) \leq \mu_{\alpha_{i+1}}^T \leq \cdots \leq \mu_\alpha^T(x)$ for every $x \in P$. Now if $ab = c$ in P and $\min\{\mu_{\alpha_i}^T(a), \mu_{\alpha_j}^T(b)\}$, then it is easy to see that

$$\mu_{\alpha i+1}^T(c) = \min\{\mu_{\alpha i+1}^T(a), \mu_{\alpha i+1}^T(b)\}$$

for $i = 0, 1, 2, \ldots$. Further if $ab = c$ in P and $\mu_{\alpha_i}^T(c) > \min\{\mu_\alpha^T(a), \mu_\alpha^T(b)\}$ then we can easily prove that

$$\mu_{\alpha+1}^T(c) > \min\{\mu_{\alpha+1}^T(a), \mu_{\alpha+1}^T(b)\}$$

for $i = 0, 1, 2, \ldots$ Hence $\{\mu_{\alpha'}^T\}_{\alpha' \in S}$ is an open fuzzy extension chain for P.

Now we will prove that $\{\mu_{\alpha'}^T\}_{\alpha' \in S}$ generates the half groupoid P if and only if a is a constant map.

Let $\{\mu_{\alpha'}^T\}_{\alpha' \in S}$ generate the half groupoid P. Then we have $\bigcup_{\alpha' \in S} \mu_{\alpha'}^T = 1_P$ that is $\sup(\mu_{\alpha'}^T(x) \mid \alpha' \in S) = 1$ for every $x \in P$, that is $\sup\{\mu(x) + \alpha' \mid \alpha' \in S\} = 1$ for every $x \in P$, that is $\mu(x) + 1 - \sup\{\mu(y) \mid y \in P\} = 1$ for every $x \in P$. That is $\mu(x) = \sup\{\mu(y) \mid y \in P\}$ for every $x \in P$. This proves that $\mu$ is a constant map.

Conversely let μ be a constant map that is $\mu(x) = \gamma$ for every $x \in P$ (where $\gamma$ is fixed element in $[0, 1]$). Consider $\left(\bigcup_{\alpha' \in S} \mu_{\alpha'}^T\right)(x)$ for an arbitrary $x \in P$.

$$
\begin{aligned}
\left(\bigcup_{\alpha' \in S} \mu_{\alpha'}^T\right)(x) &= \sup\left(\mu_{\alpha'}^T(x) \mid \alpha' \in S\right) \\
&= \sup\{\mu(x) + \alpha' \mid \alpha' \in S\} \\
&= \mu(x) + 1 - \sup\{\mu(y) \mid y \in P\} \\
&= \gamma + 1 - \gamma = 1 \\
&= 1_P(x)
\end{aligned}
$$

since x is an arbitrary element of the half groupoid P we have $\left(\bigcup_{\alpha' \in S} \mu_{\alpha'}^T\right)(x) = 1_{P(x)}$ for every $x \in P$. Hence $\{\mu_{\alpha'}^T\}_{\alpha' \in S}$ generates the half groupoid P.

Now we introduce the notions of fuzzy equivalence relation as given by [144] and proceed on to recall the definition of fuzzy relation on the groupoid given by [70].

**DEFINITION 1.11.14:** *Let X and Y be two non empty sets, we call a map $f: X \times X \to Y \times Y$ a semi balanced map if*

  i. *Given $a \in X$ there exists $\upsilon \in Y$ such that $f(a, a) = (\upsilon, \upsilon)$.*
  ii. *$f(a, a) = (\upsilon, \upsilon), f(b, b) = (\upsilon, \upsilon)$ imply that $f(a, b) = (\upsilon, \upsilon)$.*

**DEFINITION [144]:** *Let $\lambda$ and $\mu$ be two fuzzy relations on the set X. The sup-min product $\lambda \circ \mu$ is a fuzzy relation on X defined by*

$$\lambda \circ \mu (a, b) = \sup_{x \in X} \min\{\lambda(a,x), \mu(x,b)\};$$

$a, b \in X$; $\lambda$ is reflexive on X if $\lambda(x, x) = 1$ for all $x \in X$; $\lambda$ is symmetric on X if $\lambda(a, b) = \lambda(b, a)$ for all $a, b \in X$; and $\lambda$ is transitive on X if $\lambda \circ \lambda \subseteq \lambda$.

*A reflexive, symmetric and transitive fuzzy relation on a set is a fuzzy equivalence relation.*

**DEFINITION [112]:** *Let f be a mapping from the set X into the set Y. If $\lambda$ is a fuzzy subset of X, the image $f(\lambda)$ of $\lambda$ is a fuzzy subset of Y defined by*

$$f(\lambda)(\mu) = \begin{cases} \sup_{x \in f^{-1}(x)} \lambda(x), & \text{if } f^{-1}(y) \neq 0 \\ 0 & \text{otherwise} \end{cases}.$$

*If $\mu$ is a fuzzy subset of Y, the preimage or inverse image $f^{-1}(\mu)$ of $\mu$ is the fuzzy subset of X defined by $f^{-1}(\mu)(x) = \mu(f(x))$, $x \in X$.*

**DEFINITION [70]**: *Let $\lambda$ be a fuzzy relation on the groupoid D; $\lambda$ is compatible on D if $\lambda(ac, bd) \geq \min\{\lambda(a, b), \lambda(c, d)\}$ for all $a, b, c, d \in D$.*

*A compatible fuzzy equivalence relation on a groupoid is a congruence.*

**DEFINITION [45]**: *A fuzzy relation $\lambda$ on the set X is G-reflexive if for all $a \neq b$ in X.*

   i. $\lambda(a, a) \geq 0$.
   ii. $\lambda(a, b) < \delta(\lambda)$, where $\delta(\lambda) = \inf_{x \in X} \lambda(x, x)$.

A G-reflexive and transitive fuzzy relation X is a G-pre order on X. A symmetric G-pre order on X is a G-equivalence on X. A compatible G-equivalence fuzzy relation on a groupoid is a G-congruence.

**THEOREM [45]**: *If $\lambda$ is a G-pre order on the set then $\lambda \circ \lambda = \lambda$.*

**DEFINITION [44, 45]**: *Let $\lambda$ be a fuzzy relation on the set X and let $0 < \alpha \leq 1$; $\lambda$ is $\alpha$-reflexive on X if $\lambda(a, a) = \alpha$ and $\lambda(a, a) \leq \alpha$ for all $a, b \in X$. An $\alpha$-reflexive symmetric and transitive fuzzy relation on X is a fuzzy $\alpha$-equivalence relation on X. We call a compatible fuzzy $\alpha$-equivalence relation on a groupoid as $\alpha$-congruence.*

*Clearly a fuzzy $\alpha$-equivalence relation is a fuzzy equivalence, if $\alpha = 1$, and every fuzzy $\alpha$-equivalence relation is $\alpha$ G-equivalence.*

Now we proceed on to recall definition of fuzzy equivalences and congruences on a groupoid.

**THEOREM [44, 45]**: *If $\mu$ is compatible fuzzy relation on the groupoid S and f is a groupoid homomorphism from $D \times D$ into $S \times S$ then $f^{-1}(\mu)$ is a compatible fuzzy relation on D.*

*Proof*: Let $a, b, c, d, \in D$. Then we have $f^{-1}(\mu)(ac, bd) = \mu(f(ac, bd)) = \mu(f(a, b) \bullet f(c, d)) \geq \min\{\mu(f(a, b), \mu(f(c, d))\} = \min\{f^{-1}(\mu)(a, b), f^{-1}(\mu)(c, d)\}$.

Hence $f^{-1}(\mu)$ is a compatible fuzzy relation on D.

**THEOREM [44, 45]**: *If $\lambda$ is a compatible fuzzy relation on the groupoid D and f is a groupoid homomorphism from $D \times D$ in to $S \times S$ then $f(\lambda)$ is a compatible fuzzy relation on S.*

*Proof*: Let $\upsilon, \nu, \omega, r, \in S$. In the case when either $f^{-1}(\upsilon, \nu)$ or $f^{-1}(\omega, r)$ is empty we have $f(\lambda)(\upsilon\omega, \nu r) \geq 0 = \min\{f(\lambda)(\upsilon, \nu), f(\lambda)(\omega, r)\}$.

Now, consider the case when $f^{-1}(\upsilon, \nu)$ and $f^{-1}(\omega, r)$ both are nonempty. Then $f^{-1}(\upsilon\omega, \nu r)$ is non-empty too. We then have

$$
\begin{aligned}
f(\lambda)(\upsilon\omega, \nu r) &= \sup_{(x,x^1)\in f^{-1}(\upsilon\omega,\nu r)} \lambda(x, x^1) \\
&\geq \sup_{(ac,bd)\in f^{-1}(\upsilon\omega,\nu r)} \lambda(ac, bd) \\
&= \sup_{f(a,b);\, f(c,d)=(\mu,\nu)(\omega,r)} \lambda(ac, bd) \\
&\geq \sup_{f(a,b)=(\upsilon,\nu)} \sup_{f(c,d)=(\omega,r)} \min(\lambda(a,b), \lambda(c,d)) \\
&= \min\left\{ \sup_{f(a,b)=(\upsilon,\nu)} \lambda(a,b),\ \sup_{f(c,d)=(\omega,r)} \lambda(c,d) \right\} \\
&= \min\{f(\lambda)(\upsilon, \nu), f(\lambda)(\omega, r)\}.
\end{aligned}
$$

Hence $f(\lambda)$ is a compatible fuzzy relation on S.

**DEFINITION 1.11.15:** *Let X and Y be two non empty sets. A mapping $f : X \times X \to Y \times Y$ is called a semi balanced mapping, if*

  i. *Given $a \in X$ there exists $u \in Y$ such that $f(a, a) = (\upsilon, \upsilon)$.*
  ii. *$f(a, a) = (\upsilon, \upsilon)$ and $f(b, b) = (\nu, \nu)$ where $a, b \in X$, $\upsilon, \nu \in Y$ implies that $f(a, b) = (\upsilon, \nu)$.*

*One can easily verify the following result for the semi balanced mapping f. If $f(a, b) = (\upsilon, \nu)$ where $a, b \in X$, $\upsilon, \nu \in Y$ then*

  i. *$f(a, a) = (\upsilon, \upsilon)$ and $f(b, b) = (\nu, \nu)$.*
  ii. *$f(b, a) = (\nu, \upsilon)$ and*
  iii. *given $x \in X$ there exists a unique $t_x \in Y$ such that $f(a, x) = (\upsilon, t_x)$, $f(x, b) = (t_x, \nu)$ and $f(x, x) = (t_x, t_x)$.*

We give an example of a semibalanced map.

*Example 1.11.14:* Let $X = \{a, b, c\}$ and $Y = \{\upsilon, \nu, \omega, r\}$. Consider the map $f : X \times X \to Y \times Y$ defined by $f(a, b) = f(c, b) = (\upsilon, \nu)$, $f(b, a) = f(b, c) = (\nu, \upsilon)$, $f(a, c) = f(c, a) = f(c, c) = f(a, a) = (\upsilon, \upsilon)$ and $f(b, b) = (\upsilon, \upsilon)$. Then, f is a semi-balanced map.

**THEOREM 1.11.16:** *If f is a semi-balanced map from $X \times X$ into $Y \times Y$ and $\mu$ is an $\alpha$-equivalence fuzzy relation on Y, then $f^{-1}(\mu)$ is an $\alpha$-equivalence fuzzy relation on X.*

*Proof:* Let $a, b \in X$. Then $f^{-1}(\mu)(a, a) = \mu(f(a, a)) = \alpha$, because $f(a, a) = (\upsilon, \upsilon)$ for some $\upsilon \in Y$, and $f^{-1}(\mu)(a, b) = \mu(f(a, b)) \leq \alpha$. Thus $f^{-1}(\mu)$ is an $\alpha$- reflexive fuzzy relation on X.

Now $f^{-1}(\mu)(a, b) = \mu(\upsilon, \nu)$ (where $f(a, b) = (\upsilon, \nu)$) $= \mu(\nu, \upsilon)$ (by symmetry of $\mu$) $= \mu(f(b, a)) = f^{-1}(\mu)(b, a)$.

Further

$$
\begin{aligned}
f^{-1}(\mu) \circ f^{-1}(\mu)(a, b) &= \sup \min \{f^{-1}(\mu)(a, x), f^{-1}(\mu)(x, b)\} \\
&= \sup \min \{\mu f(a, x), \mu f(x, b)\} \\
&= \sup_{x \in X} \min \{\mu(\upsilon, t_x), \mu(t_x, \nu)\} \\
&\leq \sup_{\omega \in Y} \min \{\mu(\upsilon, \omega), \mu(\omega, \nu)\} \\
&= (\mu \circ \mu)(\upsilon, \nu) \leq \mu(\upsilon, \nu) \\
&= \mu(f(a, b)) \\
&= f^{-1}(\mu)(a, b).
\end{aligned}
$$

Hence $f^{-1}(\mu)$ is an $\alpha$-equivalence fuzzy relation on X.

**THEOREM 1.11.17:** *If $\mu$ is an $\alpha$-congruence fuzzy relation on the groupoid S and f is a groupoid homomorphism from $D \times D$ into $S \times S$ which is a semi balanced map, then $f^{-1}(\mu)$ is an $\alpha$-congruence on D.*

*Proof*: Follows by the earlier theorems.

**DEFINITION 1.11.16:** *Let f be a map from $X \times X$ into $Y \times Y$. A fuzzy relation $\lambda$ on X is f invariant if $f(a, b) = f(a_1, b_1)$ implies that $\lambda(a, b) = \lambda(a_1, b_1)$. A fuzzy relation $\lambda$ on X is weakly f-invariant if $f(a, b) = f(a_1, b)$ implies that $\lambda(a, b) = \lambda(a_1, b)$.*

*Clearly if $\lambda$ is f-invariant then $\lambda$ is weakly f-invariant but not conversely.*

**THEOREM 1.11.18:** *Let f be a semi-balanced map from $X \times X$ into $Y \times Y$. If $\lambda$ is a weakly f-invariant symmetric fuzzy relation on X with $\lambda \circ \lambda = \lambda$, then $\lambda$ is f-invariant.*

*Proof*: Let $f(a, b) = f(a_1, b_1) = (\upsilon, \nu)$ where $a, a_1, b, b_1 \in X$, $\upsilon, \nu \in Y$. Given $x \in X$ by definition there exists a unique $t_x \in Y$ such that $f(x, x) = (t_x, t_x)$, $f(a, x) = (\upsilon, t_x) = f(a_1, x)$ and $f(x, b) = (t_x, \nu) = f(x, b_1)$.

But $\lambda$ is weakly f-invariant so we get $\lambda(a, x) = \lambda(a_1, x)$ and $\lambda(x, b) = \lambda(x, b_1)$.

Lastly

$$
\begin{aligned}
\lambda(a, b) &= (\lambda \circ \lambda)(a, b) \\
&= \sup_{x \in X} \{\min \{\lambda(a, x), \lambda(x, b)\}\} \\
&= \sup_{x \in X} \min \{\lambda(a_1, x), \lambda(x, b_1)\} \\
&= (\lambda \circ \lambda)(a_1, b_1) = \lambda(a_1, b_1).
\end{aligned}
$$

Hence $\lambda$ is f-invariant.

The following theorem can be deduced from the earlier theorems.

**THEOREM 1.11.19:** *Let f be a semi balanced map from $X \times X$ into $Y \times Y$. If $\lambda$ is an $\alpha$-equivalence (or G-equivalence) fuzzy relation on X which is weakly f-invariant then $\lambda$ is f-invariant.*

**DEFINITION [44, 45]:** *A mapping $f : X \times X \to Y \times Y$ is a balanced mapping if*

  i.   $f(a, b) = (u, u) \Rightarrow a = b$.
  ii.  $f(a, b) = (u, v) \Rightarrow f(b, a) = (v, u)$.
  iii. $f(a, a) = (u, u)$ and $f(b, b) = (v, v) \Rightarrow f(a, b) = (u, v)$

*for all $a, b \in X$ and $u, v \in Y$.*

A mapping $f : X \times X \to Y \times Y$ is a balanced mapping if and only if it is a one to one semi-balanced mapping.

Now we give a condition under which the image of an $\alpha$-equivalence fuzzy relation is an $\alpha$-equivalence.

**THEOREM 1.11.20:** *Let f be a semi-balanced map from $X \times X$ onto $Y \times Y$. If $\lambda$ is an $\alpha$-equivalence fuzzy relation on X, which is weakly f-invariant, than $f(\lambda)$ is an $\alpha$-equivalence on Y.*

*Proof*: Let $u \in Y$. But, f being an onto semi-balanced map there exist $a, a' \in X$ such that $f(a, a') = (u, u) = f(a, a)$. By the earlier theorem $\lambda$ is f-invariant. Then we have here.

$f(\lambda)(u, u) = \sup_{(x,x') f^{-1}(u,u)} \lambda(x, x') = \lambda(a, a) = \alpha$. If $u, v \in Y$ then there exists $a, b \in F$ such that $f(a, b) = (u, v)$ and $f(b, a) = (v, u)$. We then have $f(\lambda)(u, v) = \sup_{(x_1,x_2) \in f^{-1}(u,v)} \lambda(x_1, x_2) = \lambda(a, b) \leq \infty$ and $f(\lambda)(u, v) = \lambda(a, b) = \lambda(b, a)$.

($\lambda$ is a symmetric fuzzy relation) = $f(\lambda)(v, u)$.

Thus, $f(\lambda)$ is a $\alpha$-reflexive and a symmetric fuzzy relation on Y. Now, given $x \in X$ there exists a unique $t_x \in Y$ such that $f(a, x) = (u, t_x)$, $f(x, b) = (t_x, v)$ and $f(x, x) = (t_x, t_x)$.

Lastly

$$\begin{aligned} f(\lambda)(u, v) &= \lambda(a, b) \geq (\lambda o \lambda)(a, b) \\ &= \sup_{x \in X} \min\{\lambda(a, x), \lambda(x, b)\} \\ &= \sup_{x \in X} \min\{f(\lambda)(u, t_x), f(\lambda)(t_x, v)\}. \\ &= \sup_{\omega \in Y} \min\{f(\lambda)(u, \omega), f(\lambda)(\omega, v)\}. \\ &= (f(\lambda) \circ f(\lambda))(u, v). \end{aligned}$$

Consequently, $f(\lambda)$ is an $\alpha$-equivalence on Y.

The following three results are a direct consequence of the above theorems got by combining a few of them.

The proofs of these are left as an exercise for the reader.

**THEOREM 1.11.21:** *Let f be a balanced map from $X \times X$ onto $Y \times Y$. If $\lambda$ is an $\alpha$-equivalence fuzzy relation on X, then $f(\lambda)$ is an $\alpha$-equivalence on Y.*

**THEOREM 1.11.22:** *Let f be a semi-balanced map and a groupoid homomorphism from $D \times D$ onto $S \times S$. If $\lambda$ is an $\alpha$-congruence fuzzy relation on D, which is weakly f-invariant then $f(\lambda)$ is an $\alpha$-congruence on S.*

**THEOREM 1.11.23:** *Let f be a balanced map and a groupoid homomorphism from $D \times D$ onto $S \times S$. If $\lambda$ is an $\alpha$-congruence fuzzy relation on D, then $f(\lambda)$ is an $\alpha$-congruence on S.*

Images and preimages of fuzzy G-equivalences and G-congruences on a groupoid are recalled from [44] in the following.

**THEOREM 1.11.24:** *Let f be a semi-balanced map from $X \times X$ onto $Y \times Y$. It $\lambda$ is a G-equivalence fuzzy relation on X, which is weakly f-invariant then $f(\lambda)$ is a G-equivalence on Y with $\delta(f(\lambda)) = \delta(\lambda)$.*

*Proof*: Let $\upsilon \neq \nu \in Y$. Then there exists a', $a \neq b \in X$ such that $f(a, a') = (\upsilon, \upsilon) = f(a, a)$ and $f(a, b) = (\upsilon, \nu)$ we then have

$$f(\lambda)(\upsilon, \upsilon) = \sup_{(x,x') \in f^{-1}(\upsilon,\upsilon)} \lambda(x, x') = \lambda(a, a) > 0$$

As $\lambda$ is f-invariant by earlier results and

$$\begin{aligned} f(\lambda)(\upsilon,\nu) &= \lambda(a,b) \leq \delta(\lambda) \\ &= \inf_{x \in X} \lambda(x,x) \\ &= \inf_{x \in X} f(\lambda)(t_x, t_x) \\ &= \inf_{\omega \in Y} f(\lambda)(\omega, \omega) \\ &= \delta(f(\lambda)). \end{aligned}$$

Thus $f(\lambda)$ is a G-reflexive fuzzy relation on Y with $\delta(f(\lambda)) = \delta(\lambda)$. Symmetry and transitivity of $f(\lambda)$ can be proved, as is evident from earlier theorems.

The following three theorems are easily deducible from the earlier theorem.

**THEOREM 1.11.25:** *Let f be a balanced map from $X \times X$ onto $Y \times Y$. If $\lambda$ is a G-equivalence fuzzy relation on X, then $f(\lambda)$ is a G-equivalence on Y with $\delta(f(\lambda)) = \delta(\lambda)$.*

**THEOREM 1.11.26:** *If f is a groupoid homomorphism and a semi-balanced map from $D \times D$ onto $S \times S$ and $\lambda$ is a G-congruence fuzzy relation on D which is weakly f-invariant then $f(\lambda)$ is a G-congruence on S with $\delta(f(\lambda)) = \delta(\lambda)$.*

**THEOREM 1.11.27:** *If f is a groupoid homomorphism and a balanced map from $D \times D$ onto $S \times S$ and $\lambda$ is a G-congruence fuzzy relation on D, then $f(\lambda)$ is a G-congruence on S with $\delta(f(\lambda)) = \delta(\lambda)$.*

**DEFINITION 1.11.17:** *Let $\mu$ be a fuzzy relation on Y and let f be a map from $X \times X$ into $Y \times Y$. We say $\mu$ is f-stable, if $f(a, b) = (\upsilon, \upsilon)$, where $a \neq b \in X$ and $\upsilon \in Y$, implies that $\mu(f(a, b)) \leq \mu(f(x, x))$ for all $x \in X$.*

The f-stable fuzzy relation is illustrated by the following example.

***Example 1.11.15:*** Consider the map f: $X \times X \to Y \times Y$. Define the fuzzy relations $\mu$ and $\sigma$ on Y as follows.

$$\mu(u, u) = \frac{1}{3},$$

$$\mu(\nu, \nu) = \mu(\omega, \omega) = \mu(r, r) = \frac{1}{2}$$

$$\mu(s, t) = \frac{1}{4} \quad \text{for all } s \neq t \in Y; \text{ and}$$

$$\sigma(\upsilon, \nu) = \sigma(\upsilon, \omega)$$
$$= \sigma(\nu, \nu)$$
$$= \sigma(\omega, \omega)$$
$$= \frac{1}{3}$$

$$\sigma(\upsilon, \upsilon) = \sigma(r, r)$$
$$= \frac{1}{4}$$

$$\sigma(\nu, \upsilon) = \sigma(\omega, \upsilon) = \sigma(\nu, \omega)$$
$$= \sigma(\omega, \nu) = \sigma(r, \upsilon)$$
$$= \sigma(\upsilon, r) = \sigma(\nu, r)$$
$$= \sigma(r, \nu) = \sigma(\omega, r)$$
$$= \sigma(r, \omega) = \frac{1}{5}$$

Now,

$$\mu(f(a,c)) = \mu(\upsilon, \upsilon) = \frac{1}{3} \leq \mu(f(x,x)) \text{ and } \sigma(f(a,c)) \leq \sigma(f(x,x))$$

for all $x \in X$. Thus $\mu$ and $\sigma$ both are f-stable. It may be noted that $\mu$ is a G-equivalence whereas $\sigma$ is not a G-equivalence.

**THEOREM 1.11.28:** *Let f be a semi-balanced map from $X \times Y$ into $Y \times Y$ and let $\mu$ be a G-equivalence fuzzy relation on Y which is f-stable. Then $f^{-1}(\mu)$ is a G-equivalence on X with $\delta(f^{-1}(\mu)) \geq \delta(\mu)$. Further if f is onto then $\delta(f^{-1}(\mu)) \geq \delta(\mu)$).*

*Proof*: Let $\alpha \in X$. Then $f^{-1}(\mu)(a, a) = \mu(f(a, a)) = \mu(\upsilon, \upsilon) > 0$ where $\upsilon \in Y$. If $a \neq b \in X$ then $f(a, b) = (\upsilon, v)$ for some $\upsilon, v \in Y$. In case $\upsilon = v$ the f-stability of $\mu$ implies that $f^{-1}(\mu)(a, b) = \mu(f(a,b)) \leq \mu(f(x,x)) = f^{-1}(\mu)(x, x)$ for all $x \in X$.

Now consider the case when $\upsilon \neq v$. Then by the G-reflexivity of $\mu$ we get

$$\begin{aligned} f^{-1}(\mu)(a, b) &= \mu(\upsilon, v) \leq \delta(\mu) \\ &= \inf_{\omega \in Y} \mu(\omega, \omega) \\ &\leq \inf_{x \in X} \mu(f(x, x)) \\ &= \inf_{x \in X} f^{-1}(\mu)((x, x)) = \delta(f^{-1}(\mu)). \end{aligned}$$

Thus $f^{-1}(\mu)$ is a G-reflexive relation on X with $\delta(f^{-1}(\mu)) \geq \delta(\mu)$. Further, if f is onto then obviously $\delta(f^{-1}(\mu)) = \delta(\mu)$. As in the proof of the earlier theorem we can show that $f^{-1}(\mu)$ is a symmetric and transitive fuzzy relation.

The following results also are left for the reader as an exercise.

**THEOREM 1.11.29:** *Let f be a balanced map from $X \times X$ into $Y \times Y$ and let $\mu$ be a G-equivalence fuzzy relation on Y. Then $f^{-1}(\mu)$ is a G-equivalence on X with $\delta(f^{-1}(\mu)) \geq \delta(\mu)$. Further if f is onto then $\delta(f^{-1}(\mu)) = \delta(\mu)$.*

*Proof*: Left for the reader. Please refer [45].

**THEOREM 1.11.30**: *Let f be a groupoid homomorphism and a semi balanced map from $D \times D$ into $S \times S$, and let $\mu$ be a G-congruence on Y which is f-stable. Then $f^{-1}(\mu)$ is a G-congruence on X with $\delta(f^{-1}(\mu)) \geq \delta(\mu)$. Further if f is onto then $\delta(f^{-1}(\mu)) = \delta(\mu)$.*

*Proof*: Follows directly by the earlier theorems.

**THEOREM 1.11.31:** *Let f be a groupoid homomorphism and a balanced map from $D \times D$ into $S \times S$ and let $\mu$ be a G-congruence on Y. Then $f^{-1}(\mu)$ is a G-congruence on Y. Then $f^{-1}(\mu)$ is a G-congruence on X with $\delta(f^{-1}(\mu)) \geq \delta(\mu)$.*

*Further if f is onto then $\delta(f^{-1}(\mu)) = \delta(\mu)$.*

For proof refer [45].

We proceed on to recall the definition of G-groupoid and P-fuzzy correspondence.

**DEFINITION 1.11.18:** *Let $S_1$, $S_2$ and $S_3$ be non empty sets and $A : S_1 \times S_2 \to S_3$ a mapping. Then $(S_1, S_2, S_3, A)$ is a G-groupoid. If $A_1$ and $A_2$ are non empty sets and $(P, \leq)$ is a partially ordered set then a function $\overline{A} : A_1 \times A_2 \to P$ is a binary P-fuzzy correspondence. By a P-fuzzy correspondence sometimes we consider a quadruple $(A_1, A_2, P, \overline{A})$.*

*If $A_1 \ldots, A_n$ are sets and $(P, \leq)$ is a partially ordered set than a mapping $\overline{A} : A_1 \times \ldots \times A_n$ to P is an (n-ary) P-fuzzy correspondence on sets $A_1, \ldots, A_n$. Here we recall a special type of binary P-fuzzy correspondences such that $(A_1, \leq), (A_2, \leq), \ldots, (A_n, \leq)$ are also partially ordered sets. We call such a correspondence PG- fuzzy correspondence. The notion of a binary PG-fuzzy correspondence is a generalization of the notion of a G-groupoid, since we obtain a binary PG-fuzzy correspondence when we equip all sets in a G-groupoid with orderings.*

*If $\overline{A} : A_1 \times \ldots \times A_n \to P$ is a P-(or PG-) fuzzy correspondence for $p \in P$, a p-level correspondence is a function $\overline{A} : A_1 \times A_2 \times \ldots \times A_n \to \{0, 1\}$ defined by $\overline{A}_p (x_1, \ldots, x_n) = 1$ if and only if $\overline{A}(x_1, x_2, \ldots, x_n) \geq p$.*

*By $A_P$ we denote the corresponding p-level subset of $A_1 \times A_2 \times \ldots \times A_n$. Further more we denote the set of all level subsets of $A_P$ by $\mathcal{A}_p = \{A_p \mid p \in P\}$.*

Since a P-fuzzy correspondence is a P-fuzzy set on a product of sets, the following results are recalled from [90].

**THEOREM 1.11.32:** *If $\overline{A} : A_1 \times A_2 \times \ldots \times A_n \to P$ is a P fuzzy correspondence where $(P, \leq)$ is a partially ordered set, then for every $(x_1, x_2, \ldots, x_n) \in A_1 \times A_2 \times \ldots \times A_n$ the following supremum exists in $P \vee (p \in P \mid \overline{A}_p(x_1, x_2, \ldots, x_n) = 1$ and $\overline{A}(x_1, x_2, \ldots, x_n) = \vee\{p \in P$ such that $\overline{A}_p(x_1, x_2, \ldots, x_n) = 1\}$.*

*Proof:* If $\overline{A}(x_1, x_2, \ldots, x_n) = q$ then $\overline{A}_q(x_1, x_2, \ldots, x_n) = 1$. If $\overline{A}_p(x_1, x_2, \ldots, x_n) = 1$ then $\overline{A}(x_1, \ldots, x_n) = q \geq p$. This means that q greater than all p for which $\overline{A}_p(x_1, x_2, \ldots, x_n) = 1$, that is $\overline{A}(x_1, x_2, \ldots, x_n) = q = \vee (p \in P \mid \overline{A}_p(x_1, x_2, \ldots, x_n) = 1)$.

**THEOREM 1.11.33:** *If $\overline{A} = A_1 \times A_2 \times \ldots \times A_n \to P$ is a P-fuzzy correspondence then*

  i. *If $p \leq q$ for $p, q \in P$, then for every $(x_1, x_2, \ldots, x_n) \in A_1 \times A_2 \times \ldots A_n$, $\overline{A}_q(x_1, x_2, \ldots, x_n) \leq \overline{A}_p(x_1, \ldots, x_n)$.*
  ii. *If the supremum of a subset Q of P exists then $\cap (A_p \mid p \in Q) = A_{\vee\{p \mid p \in Q\}}$.*

  *iii.*   $\cup (A_p \mid p \in Q) = A_1 \times \ldots \times A_n$.
  *iv.*   for every $(x_1, x_2, \ldots, x_n) \in A_1 \times A_2 \times \ldots \times A_n \cap (A_p \in \mathcal{A}_p \mid (x_1, \ldots, x_n) \in P)$ belongs to the family $A_P$ of level correspondences of $\overline{A}$.

*Proof*: The results can be easily proved using the definitions.

**THEOREM 1.11.34:** *Let $F = \{F_i \mid i \in I\}$ be a family of subsets of $A_1 \times A_2 \times \ldots \times A_n$ (crisp correspondences) such that*

  *i.*   $\bigcup_{i \in I} F_i = A_1 \times \ldots \times A_n$.
  *ii.*   for every $(x_1, x_2, \ldots, x_n) \in A_1 \times A_2 \times \ldots \times A_n$; $\cap (F_i \in F \mid (x_1, \ldots, x_n) \in F_i) \in F$.

*Let $(F, \leq)$ be the partially ordered set dual to $(F, \subseteq)$. Then a P-fuzzy correspondence $\overline{A}: A_1 \times A_2 \times \ldots \times A_n \to F$ defined with $\overline{A}(x_1, \ldots, x_n) = \cap (F_i \in F \mid (x_1, \ldots, x_n) \in F_i)$ has F as its family of level correspondences. Moreover for every $F_i \in F$, $F_i = A_{F_i}$.*

*Proof*: By the condition (ii), $\overline{A}$ is well defined. The only thing left to prove is that $F_i = A_{F_i}$ for all $F_i \in F$.

Let $F_i \in F$ then $(x_1, x_2, \ldots, x_n) \in A_{F_i}$ if and only if $\overline{A}(x_1, \ldots, x_n) \geq F_i$ if and only if $\overline{A}(x_1, \ldots, x_n) \subseteq F_j$; if and only if $\cap (F_j \in F \mid (x_1, \ldots, x_n) \in F_j) \subseteq F_i$.

Since $(x_1, \ldots, x_n) \in \cap (F_j \in F \mid (x_1, \ldots, x_n)$ is in $F_j)$, then $(x_1, \ldots, x_n) \in F_j$. On the otherhand if $(x_1, x_2, \ldots, x_n) \in F_i$, then $(\cap (F_j \in F \mid (x_1, \ldots, x_n) \in F_j) \subseteq F_i$ and as proven above $(x_1, x_2, \ldots, x_n) \in A_{F_i}$.

Notations as given by [90].

Let $A_1, A_2, A_3, A_4$ and $A_5$ be partially ordered sets and $\overline{B}: A_2 \times A_3 \to A_5$ and $\overline{D}: A_1 \times A_2 \to A_4$, PG fuzzy correspondence. $(\overline{A}, \overline{C})$ where $\overline{A}: A_1 \times A_5 \to P$ and $\overline{C}: A_4 \times A_3 \to P$ are PG-fuzzy correspondences and P is a partially ordered set, of the functional equation $\overline{A}(x, \overline{B}(y, z)) = \overline{C}(\overline{D}(x, y), z)$ for given PC-fuzzy correspondences $\overline{B}$ and $\overline{D}$.

We see the equation $\overline{A}(x, \overline{B}(y, z)) = \overline{C}(\overline{D}(x, y), z)$ always has a nontrivial solution for arbitrary $\overline{B}: A_2 \times A_3 \to A_5$ and $\overline{D}: A_1 \times A_2 \to A_4$. Indeed for $\overline{A}$ and $\overline{C}$ we take arbitrary PG-fuzzy correspondences

$\overline{A}: A_1 \times A_5 \to L$ and $\overline{C}: A_4 \times A_3 \to L$ where L is a distributive lattice (i.e. partially ordered set).

**THEOREM 1.11.35:** *Let $A_1, A_2, A_3, A_4$ and $A_5$ be partially ordered set. Further more let $\overline{B}: A_2 \times A_3 \to A_5$ and $\overline{D}: A_1 \times A_2 \to A_4$ be PG-fuzzy correspondences which are surjections.*

*Let F be a family of ternary correspondences on $A_1 \times A_2 \times A_3$ satisfying the following conditions:*

   i. *For every $y, t \in A_2$ and $z, u \in A_3$, $x \in A_1$ and $E \in F$ if $\overline{B}(y, z) = \overline{B}(t, u)$ then $E(x, y, z) = E(x, t, u)$.*
   ii. *For every $x, z \in A_1$, $y, t \in A_2$, $u \in A_3$ and $E \in F$ if $\overline{D}(x, y) = \overline{D}(z, t)$ then $E(x, y, u) = E(z, t, u)$.*
   iii. *For all $(x, y, z) \in A_1 \times A_2 \times A_3$; $\cap (E \in F \mid E(x, y, z) = 1) \in F$.*
   iv. *$\cup (E / E \in F) = A_1 \times A_2 \times A_3$.*
   v. *For all $E \in F$, let the correspondence $A_E$ on $A_1 \times A_5$ and $C_E$ on $A_4 \times A_3$ be defined in the following way:*
   $A_E(x, y) = E(x, z, t)$ for $\overline{B}(z, t) = y$.
   vi. *$C_E(x, y) = E(\upsilon, v, y)$ for $\overline{D}(\upsilon, v) = x$, then partially ordered sets F, $F_A = \{A_E / E \in F\}$ and $F_C = \{C_E / E \in F\}$ are isomorphic under the mappings $f : F_A \to F$ defined $f(A_E) = E$ and $g : F_C \to F$ defined by $g(C_E) = E$. Solutions of the functional equation $\overline{A}(x, \overline{B}(y, z)) = \overline{C}(\overline{D}(x, y), z)$ are PG-fuzzy correspondences.*

$\overline{A} : A_1 \times A_5 \to F$ and $\overline{C} : A_4 \times A_3 \to F$ given by

$$\widetilde{A}(x,y) = f(\cap (A_E \in F_A \mid A_E(x,y) = 1))$$
$$\widetilde{C}(x,y) = g(\cap (C_E \in F_C \mid C_E(x,y) = 1))$$

where $(F, \leq)$ is the dual of $(F, \subseteq)$.

*Proof:* (Please refer [90]).

Consequent of this we have the following theorem and example given by [90].

**THEOREM 1.11.36:** *All subsets of the functional equation of generalized associativity. $\overline{A}(x, \overline{B}(y, z) = \overline{C}(\overline{D}(x, y), z)$ for PG-fuzzy correspondences where $\overline{B}$ and $\overline{D}$ are subjections given by $\overline{B} : A_2 \times A_3 \to A_5$ and $\overline{D} : A_2 \times A_3 \to A_4$ which are arbitrary, PG-fuzzy correspondences and subjections and $\overline{A}$, $\overline{C}$ are given by (i) to (vi) of the theorem 1.11.35.*

***Example [90]:*** Let $A_1$, $A_2$, $A_3$, $A_4$ $A_5$ be partially ordered sets given by the following figures.

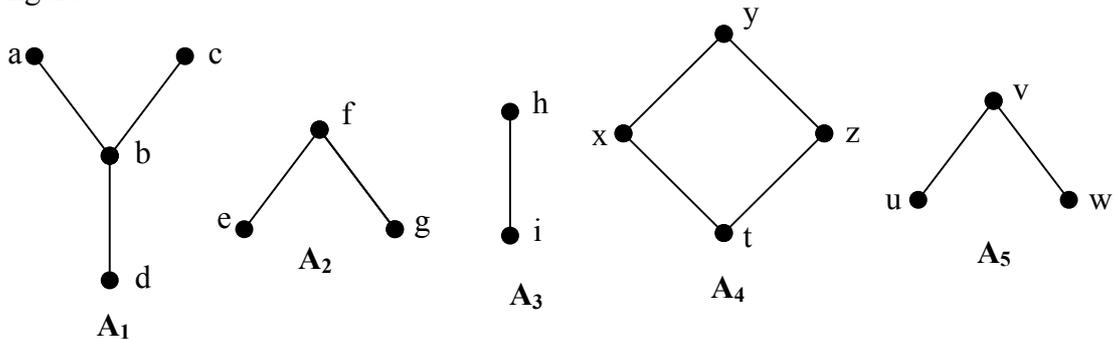

**Figures 1.11.1**

Let P-fuzzy correspondences be given as follows:

| $\overline{B}$ | h | i |
|---|---|---|
| e | υ | ν |
| f | υ | ν |
| g | ν | ω |

| $\overline{D}$ | e | f | g |
|---|---|---|---|
| a | x | x | z |
| b | x | x | z |
| c | y | y | t |
| d | y | y | t |

To find a family of correspondence as on $A_1 \times A_2 \times A_3$ satisfying conditions for all $E \in F$ and all $k \in A_1$ and $j \in A_3$.

$$E(k, e, g) = E(k, f, h)$$
$$E(k, e, i) = E(k, f, i) = E(k, g, h)$$
$$E(a, e, j) = E(a, f, j) = E(b, e, j) = E(b, f, j)$$
$$E(c, e, j) = E(c, f, j) = E(d, e, j) = E(d, f, j)$$
$$E(a, g, j) = E(b, g, j)$$
$$E(c, g, j) = E(d, g, j)$$

such that the family F itself satisfies conditions.

Such a family of correspondences F is given as follows:

| P | eh | ei | fh | fi | gh | gi |
|---|---|---|---|---|---|---|
| a | 0 | 1 | 0 | 1 | 1 | 0 |
| b | 0 | 1 | 0 | 1 | 1 | 0 |
| c | 0 | 0 | 0 | 0 | 0 | 0 |
| d | 0 | 0 | 0 | 0 | 0 | 0 |

| Q | eh | ei | fh | fi | gh | gi |
|---|---|---|---|---|---|---|
| a | 1 | 1 | 1 | 1 | 1 | 1 |
| b | 1 | 1 | 1 | 1 | 1 | 1 |
| c | 1 | 1 | 1 | 1 | 1 | 0 |
| d | 1 | 1 | 1 | 1 | 1 | 0 |

| R | eh | ei | fh | fi | gh | gi |
|---|---|---|---|---|---|---|
| a | 0 | 0 | 0 | 0 | 0 | 1 |
| b | 0 | 0 | 0 | 0 | 0 | 1 |
| c | 0 | 1 | 0 | 1 | 1 | 0 |
| d | 0 | 1 | 0 | 1 | 1 | 0 |

| S | eh | ei | fh | fi | gh | gi |
|---|---|---|---|---|---|---|
| a | 0 | 1 | 0 | 1 | 1 | 1 |
| b | 0 | 1 | 0 | 1 | 1 | 1 |
| c | 0 | 1 | 0 | 1 | 1 | 1 |
| d | 0 | 1 | 0 | 1 | 1 | 1 |

The partially ordered set $(F, \leq)$ is presented in the figure.

The following corresponding families $\{A_E \mid E \in F\}$ and $\{C_E \mid E \in F\}$ are obtained.

| $A_P$ | υ | ν | ω |
|---|---|---|---|
| a | 0 | 1 | 0 |
| b | 0 | 1 | 0 |
| c | 0 | 0 | 0 |
| d | 0 | 0 | 0 |

| $A_Q$ | υ | ν | ω |
|---|---|---|---|
| a | 1 | 1 | 1 |
| b | 1 | 1 | 1 |
| c | 1 | 1 | 0 |
| d | 1 | 1 | 0 |

| $A_R$ | $\upsilon$ | $\nu$ | $\omega$ |
|---|---|---|---|
| a | 0 | 0 | 1 |
| b | 0 | 0 | 1 |
| c | 0 | 1 | 0 |
| d | 0 | 1 | 0 |

| $A_S$ | $\upsilon$ | $\nu$ | $\omega$ |
|---|---|---|---|
| a | 0 | 1 | 1 |
| b | 0 | 1 | 1 |
| c | 0 | 1 | 1 |
| d | 0 | 1 | 1 |

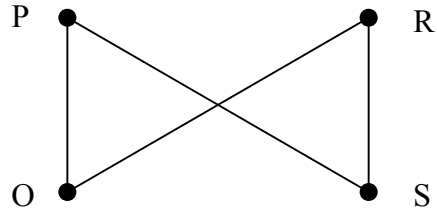

**Figure 1.11.2**

| $C_P$ | h | i |
|---|---|---|
| x | 0 | 1 |
| y | 0 | 0 |
| z | 1 | 0 |
| t | 0 | 0 |

| $C_Q$ | h | i |
|---|---|---|
| x | 1 | 1 |
| y | 1 | 1 |
| z | 1 | 1 |
| t | 1 | 0 |

| $C_R$ | h | i |
|---|---|---|
| x | 0 | 0 |
| y | 0 | 1 |
| z | 0 | 1 |
| t | 1 | 0 |

| $C_S$ | h | i |
|---|---|---|
| x | 0 | 1 |
| y | 0 | 1 |
| z | 1 | 1 |
| t | 1 | 1 |

By the synthesis of these families of correspondences we obtain the following fuzzy correspondences.

| $\overline{A}$ | $\upsilon$ | $\nu$ | $\omega$ |
|---|---|---|---|
| a | Q | P | R |
| b | Q | P | R |
| c | Q | R | S |
| d | Q | R | S |

| $\overline{C}$ | h | i |
|---|---|---|
| x | Q | P |
| y | Q | R |
| z | P | R |
| t | R | S |

$(\overline{A}, \overline{B}, \overline{C}, \overline{D})$ is a solution of the functional equation.

It is important and interesting to note that PG-fuzzy correspondences such that all related partially ordered sets are equal (S, ≤), $\overline{A} : S \times S \to S$, $\overline{B}: S \times S \to S$, $\overline{C}: S \times S \to S$ and $\overline{D}: S \times S \to S$ are semilattices i.e. groupoids satisfying commutative, idempotent and associative laws. We prove that if such $\overline{A}$, $\overline{B}$, $\overline{C}$ and $\overline{D}$ are solutions of the functional equation then they are equal that is $\overline{A} = \overline{B} = \overline{C} = \overline{D}$.

**THEOREM 1.11.37:** *If (S, $\overline{A}$), (S, $\overline{B}$), (S, $\overline{C}$) and (S, $\overline{D}$) are semi-lattices such that for all x, y, z ∈ S, $\overline{A}$ (x, $\overline{B}$ (y, z)) = $\overline{C}$ ($\overline{D}$ (x, y), z) then $\overline{A} = \overline{B} = \overline{C} = \overline{D}$.*

*Proof*: Let x, y ∈ S. Then

$\overline{A}$ (x, y) = $\overline{A}$ (y, x) = $\overline{A}$ (y, $\overline{B}$ (x, x)) = $\overline{C}$ ($\overline{D}$ (y, x), x) = $\overline{C}$ ($\overline{D}$ (x, y), x).
$\overline{A}$ (x, $\overline{B}$ (y, x)) = $\overline{A}$ (x, $\overline{B}$ (x, y)) = $\overline{C}$ ($\overline{D}$ (x, x), y) = $\overline{C}$ (x , y).
$\overline{A}$ (x, y) = $\overline{A}$ (x, $\overline{B}$ (x, y))

$$
\begin{aligned}
&= \overline{A}(x, \overline{B}(\overline{B}(x, y)), \overline{B}(x, y)) \\
&= \overline{C}(\overline{D}(x, \overline{B}(x, y), \overline{B}(x, y)) \\
&= \overline{C}(\overline{D}(\overline{B}(x, y), x), \overline{B}(x, y)) \\
&= \overline{A}(\overline{B}(x, y), \overline{B}(x, \overline{B}(x, y))), \\
&= \overline{A}(\overline{B}(x, y), \overline{B}(x, y)). \\
&= \overline{B}(x, y).
\end{aligned}
$$

Similarly we prove $\overline{C} = \overline{D}$.

This will find several applications in social groups, pattern recognitions, traffic, public relations and correspondences.

Now using [110], we recall some results about lattice of all idempotent fuzzy subsets of a groupoids.

We know from the definition of a fuzzy subset of a set S originally defined by [144] in his classical paper is a mapping from S into the real interval [0, 1]. If λ and μ are fuzzy subsets of S then the equality λ = μ and the inclusion λ ⊆ μ are defined pointwise, endowed with this partial order the set F(S) of all fuzzy subsets of S form a completely distributive lattice. The meet λ ∩ μ and the join λ ∪ μ the arbitrary supremum ∪ $\lambda_i$ and the arbitrary infimum ∩ $\lambda_i$ can be verified to have the following formulas.

    i.    (λ ∩ μ) (x) = min {λ (x) , μ (x)}.
    ii.   (λ ∪ μ) (x) = min {λ(x) , μ(x)}.
    iii.  (∪ $\lambda_i$ ) (x) = $\sup_i$ λ(x).
    iv.  (∩ $\lambda_i$ )(x) = $\inf_i$ λ(x), x ∈ S.

The sup-min product of two fuzzy subsets λ and μ of a groupoid D is the fuzzy subset λ. μ of D defined by

$$(\lambda \bullet \mu)(x) = \begin{cases} \sup_{x=ab} \min\{\lambda(a), \mu(b)\} & \text{if } x \text{ is factorizable in } D \\ 0 & \text{otherwise.} \end{cases}$$

This product is monotone. It is associative and commutative if the composition in D is associative and commutative respectively. The following two statements are equivalent

i. $\lambda \bullet \lambda \subseteq \lambda$
ii. $\min\{\lambda(a), \lambda(b)\} \leq \lambda(ab)$

for all a, b ∈ D.

If D is a semigroup and n is a positive integer, $\lambda^n$ denotes $\lambda \bullet \lambda \bullet \ldots \bullet \lambda$ (n times).

Here we give some of the important subclasses of the class F(G) of all fuzzy subsets of a group G with identity e as follows:

i. $F_{-1}(G)$ : It is the class of all $\lambda \in F(G)$ such that $\lambda \subseteq \lambda \bullet \lambda$.
ii. $F_0(G)$ : It is the class of all $\lambda \in F(G)$ such that $\lambda(x) \leq \lambda(e)$ for all x ∈ G.
iii. $F_1(G)$ : It is the class of all $\lambda \in F(G)$ such that $\lambda \bullet \lambda \subseteq \lambda$. It means $F_1(G)$ is the set of all fuzzy subgroupoids of G.
iv. $F_2(G)$ : It is the class of all $\lambda \in F(G)$ such that $\lambda \bullet \lambda = \lambda$. In other words $F_2(G) : F_1(G) \cap F_{-1}(G)$.
v. $F_3(G) = $ It is defined by $F_3(G) = F_1(G) \cap F_0(G)$.
vi. $F_4(G)$ : It is the class of all fuzzy subgroups of G.

It is easy to see that $F_4(G) \subseteq F_3(G) \subseteq F_2(G) \subseteq F_1(G)$ and that $F_4(G) \subseteq F_3(G) \subseteq F_0(G) \subseteq F_{-1}(G)$. If D is a groupoid the classes $F_{-1}(D)$, $F_1(D)$ and $F_2(D)$ are similarly defined. Besides, if D possesses an identity element the classes $F_0(D)$ and $F_3(D)$ are meaningful the figure shows the inclusion relation among the various $F_i(G)$.

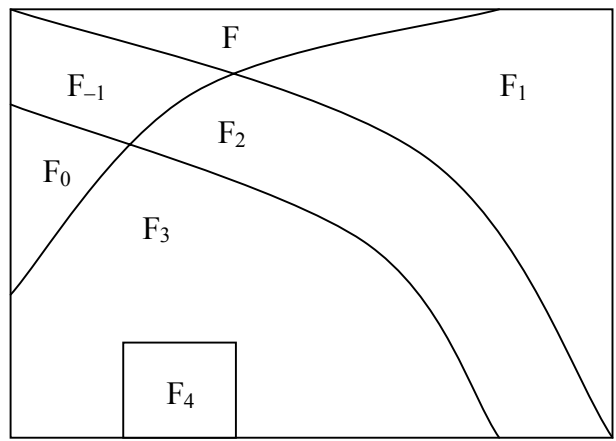

**Figure 1.11.3**

Throughout this note, D denotes a groupoid and G denotes a group with identity e, unless otherwise specifically stated.

**THEOREM 1.11.38:** *Let $\lambda_i, \mu_j, \lambda, \mu \subseteq F(D)$ where $i \in I, j \in J$. Then*

i. $\lambda \bullet \left(\bigcup_j \mu_j\right) = \bigcup_j (\lambda \bullet \mu_j)$

ii. $\left(\bigcup_i \lambda_i\right) \bullet \mu = \bigcup_i (\lambda_i \bullet \mu)$

iii. $\left(\bigcup_i \lambda_i\right) \bullet \left(\bigcup_j \mu_j\right) = \bigcup_{i,j} (\lambda_i, \mu_j)$

*Proof (i):*

$$\left(\lambda \bullet \left(\bigcup_j \mu_j\right)\right)(x) = \sup_{x=ab} \min\left\{\mu(a), \sup_j \mu_j(b)\right\}$$

$$= \sup_{x=ab}\left(\sup_j \min\{\lambda(a), \mu_j(b)\}\right)$$

$$= \sup_j\left(\sup_{x=ab} \min\{\lambda(a), \mu_j(b)\}\right)$$

$$= \sup_j (\lambda, \mu_j)(x)$$

for all $x \in D$. If $x \in D$ and $0 < \alpha \leq 1$, then the fuzzy subset of D that maps x to $\alpha$ and the other elements of D to 0 is the fuzzy point $x_\alpha$ in D.

For the sake of notational convenience we denote it by $(x, \alpha)$. Upon a moments reflection one realizes that the class M (D) of all fuzzy points in D is in fact the set of all join irreducible elements of the completely distributive lattice F (D).

**DEFINITION 1.11.19:** *Given a non empty element $\lambda \in F(S)$, where S is an arbitrary set, there exists a class $\{(x_i, \alpha_i)\}_i$ in M(S) not necessarily unique such that $\lambda = \bigcup_i (x_i, \alpha_i)$.*

*We call such an expression for $\lambda$ a fuzzy-point expression or simply a reduction of $\lambda$ in M(S).*

We illustrate this example from [110].

***Example [110]:*** Let S = {a, b, c}. If $\lambda$ is the fuzzy subset of S that maps a, b, c to 1, $\frac{1}{3}$, 0 respectively then $\lambda = (a,1) \cup \left(b, \frac{1}{3}\right)$ is a reduction of $\lambda$ in M (S). There exist infinitely many reductions of $\lambda$ as follows.

$\lambda = \left(a, \frac{1}{4}\right) \cup (a,1) \cup \left(b, \frac{1}{3}\right)$, $\lambda = \left(b, \frac{1}{3}\right) \cup \left(a, \frac{1}{2}\right) \cup \left(a, \frac{2}{3}\right) \cup \left(a, \frac{3}{4}\right) \cup \cdots$

**THEOREM 1.11.39:** *Let D be a groupoid in which the cancellation laws are valid. If $(x_i, \alpha_i)$, $(y_j, B_j)$, $(y, \beta) \in F(D)$, $i \in I, j \in J$ then*

i. $\left(\bigcap_i (x_i, \alpha_i)\right)((y, \beta)) = \bigcap_i ((x_i, \alpha) \bullet (y, \beta))$

ii. $(y, \beta) \bullet \left(\bigcap_i (x_i, \alpha_i)\right) = \bigcap_i ((y, \beta) \bullet (x_i, \alpha_i))$

iii. $\left(\bigcap_i (x_i, \alpha_i)\right) \bullet \bigcap_j (y_j, B_j) = \bigcap_{i \bullet j} ((x_i, \alpha_i), (y_j, B_j))$

*Proof*: Follows from [4].

**THEOREM 1.11.40:** *Let D be a groupoid in which the cancellation laws are valid. If $\lambda_i, \mu_j, \lambda, \mu \in F(D)$ where $i \in I$ and $j \in J$ then*

i. $\left(\bigcap_i \lambda_i\right) \bullet \mu = \bigcap_i (\lambda_i, \mu)$.

ii. $\lambda \bullet \left(\bigcap_j \mu_j\right) = \bigcap_j (\lambda \bullet \mu_j)$.

iii. $\left(\bigcap_i \lambda_i\right) \bullet \left(\bigcap_j \lambda_j\right) = \bigcap_{i,j} (\lambda_j, \mu_j)$.

*Proof:* Please refer [110]. $F_{-1}(D)$ is not a sublattice of $F(D)$ as is shown in by the example of [110].

*Example [110]*: We show in this example $F_{-1}(D)$ is not closed under finite intersection. Let $D = \{1, -1, i, -i\}$ with the usual multiplication of complex numbers. Let $A = \{1, i\}$, $B = \{-1, i, -i\}$. Let $\chi_A$ denote the characteristic function of A with domain D. We then have $\chi_A \subset \chi_{A^2} = \chi_A \cdot \chi_A$ where C denotes proper inclusion $\chi_B \subset \chi_{B^2} = \chi_B \cdot \chi_B$, $\chi_A \cap \chi_B = \chi_{A \cap B}$ and $(\chi_A \cap \chi_B) \bullet (\chi_A \cap \chi_B) = \chi_{(A \cap B)^2}$. Consequently $\chi_A$, $\chi_B \in F_{-1}(D)$ but $\chi_A \cap \chi_B$ and $(\chi_A \cap \chi_B) \bullet (\chi_A \cap \chi_B)$ are not even comparable. Nevertheless $F_{-1}(D)$ does not form a join – sublattice of $F(D)$. Infact, its members form a complete lattice.

**THEOREM 1.11.41:** *Let D be an arbitrary groupoid. Then*

i. $F_{-1}(D)$ *is closed under the formation of arbitrary unions.*
ii. $F_{-1}(D)$ *is a complete lattice under the meet $\lambda \wedge \mu$ and join $\lambda \vee \mu$ given by $\lambda \wedge \mu = \{\lambda \cap \mu\}_{-1}$ and $\lambda \vee \mu = (\lambda \cup \mu)$.*

*Proof:*

i. Let $\{\lambda_i\}_i \subseteq F_{-1}(D)$. Then $\lambda_i \subseteq \lambda_i \bullet \lambda_i \subseteq \left(\bigcup_i \lambda_i\right) \bullet \left(\bigcup_j \lambda_j\right)$ for all i, which implies that $\cup \lambda_i \subseteq \left(\bigcup_i \lambda_i\right) \bullet \left(\bigcup_i \lambda_i\right)$.

In other words $\cup_i \lambda_i \in F_{-1}$.

ii. By (i) every nonempty subset of $F_{-1}$ has an l.u.b.

in $F_{-1}$. Besides $\phi \in F_{-1}$ Hence by well known results on lattice theory [27]. $F_{-1}$ is a complete lattice under meet $\lambda \cap \mu$ and join $\lambda \cup \mu$ given by $\{\lambda \cap \mu\}_{-1}$ and $\lambda \cup \mu$ respectively.

**THEOREM 1.11.42:** *Let D be a semigroup in which the cancellation laws are valid. If* $\lambda \in F(D)$ *then* $|\lambda|_{-1} = \bigcap_{n=1}^{\infty} \lambda^n$

Proof: Let $\mu = \bigcap_{n=1}^{\infty} \lambda^n \subseteq \lambda$. By earlier theorems we have $\mu \bullet \mu = \lambda^2 \cap \lambda^3 \cap \lambda^4 \cap \ldots \supset \lambda \cap \lambda^2 \cap \lambda^3 \cap \ldots = \mu$. In other words $\mu \in F_{-1}(D)$. Let $\delta \in F_{-1}$ be such that $\delta \subseteq \lambda$. Because $\delta \subseteq \delta^2 \subseteq \lambda^2$, by induction we obtain $\delta \subseteq \lambda^n$ for all positive integers n. Hence $\delta \subseteq \mu$. Consequently $\mu = |\lambda|_{-1}$.

The class $F_1(D)$ is not a sublattice of $F(D)$. This is explained by the example of [110].

*Example [110]:* Here it is proved that $F_1(D)$ is not closed under the formation of finite unions. Let $D = \{a, b, c, d\}$ be the set with composition table:

| • | a | b | c | d |
|---|---|---|---|---|
| a | a | b | a | d |
| b | a | b | d | c |
| c | a | d | a | d |
| d | c | d | c | c |

Let $\lambda$ be the fuzzy subset of D that maps a, b, c, d to $1, \frac{1}{2}, \frac{1}{3}, \frac{1}{3}$ respectively and let $\mu$ be the fuzzy subset of D that maps a, b, c, d to $1, \frac{1}{3}, \frac{1}{2}, \frac{1}{3}$ respectively. Then $\lambda, \mu, \in F_1(D)$, because the level subsets of them are all subgroupoids of D. However $\lambda \cup \mu$ maps a, b, c, d to $1, \frac{1}{2}, \frac{1}{2}, \frac{1}{3}$ respectively. Clearly, $\lambda \cup \mu \notin F_1(D)$ because the level subset $\{a, b, c\}$ is not a subgroupoid of D. Neverthless $F_1(D)$ does not form a meet sublattice of $F(D)$.

**THEOREM [110]:** *Let D be an arbitrary groupoid then*

  i. *$F_1(D)$ is closed under the formation of arbitrary intersections.*

  ii. *$F_1(D)$ is a complete lattice under join $\lambda \vee \mu$ and meet $\lambda \wedge \mu$ given by $\lambda \wedge \mu = \lambda \cap \mu$ and $\lambda \vee \mu = \{\lambda \cup \mu\}$.*

*Proof:*

i. Let $\{\lambda_i\}_i \subseteq F_1(D)$. Then $\left(\bigcap_i \lambda_i\right) \bullet \left(\bigcap_i \lambda_i\right) \subseteq \lambda_i \bullet \lambda_i \subseteq \lambda_i$ for all i, and we get $\left(\bigcap_i \lambda_i\right) \bullet \left(\bigcap_i \lambda_i\right) \subseteq \bigcap_i \lambda_i$. In other words $\bigcap_i \lambda_i \in F_1(D)$.

ii. From (i) every non-empty subset of $F_1$ has a g. l. b in $F_1$. Besides, $l \in F_1$. Consequently, $F_1$ is complete lattice under meet $\lambda \vee \mu$ and join $\lambda \wedge \mu$ given by $\lambda \wedge \mu = \lambda \cap \mu$ and $\lambda \vee \mu = \{\lambda \cup \mu\}$.

**THEOREM [110]:** *Let D be a semigroup. If $\lambda \in F(D)$ then $[\lambda]_i = \bigcup_{n=1}^{\infty} \lambda^n$*

*Proof:* Let $\mu = \bigcup_{n=1}^{\infty} \lambda^n \supseteq \lambda$. By earlier results we have

$$\mu \bullet \mu = \lambda^2 \cup \lambda^3 \cup \lambda^4 \cup \cdots \subseteq \lambda \cup \lambda^2 \cup \lambda^3 \cup \cdots = \lambda.$$

In other words $\mu \in F_1(D)$. Let $\delta \in F_1$ be such that $\lambda \subseteq \delta$. Because $\lambda^2 \subseteq \delta^2 \subseteq \delta$ by induction we obtain $\lambda^n \subseteq \delta$ for all positive integers n. Hence $\mu \subseteq \delta$ consequently $\mu = [\lambda]_l$.

**THEOREM 1.11.43:** *Let D be a semigroup.*

i. *If $\lambda \in F_{-1}(D)$ then $\lambda^n \in F_{-1}(D)$ for all positive integers n. In other words $F_{-1}(D)$ is the union of all the monogenic subsemigroups generated by its own elements.*

ii. *If $\lambda \in F_1(D)$ then $\lambda^n \in F_1(D)$ for all positive integers n. In other words $F_1(D)$ is the union of all the monogenic subsemigroups generated by its own elements.*

*Proof:* Follows easily with simple calculations.

**THEOREM 1.11.44:** *If $\lambda \in F(D)$ then the two statements (i) and (ii) that follow are equivalent.*

i. $\lambda \in F_2(D)$.
ii. *(a) $\lambda \in F_1(D)$ and (b) given $d \in D$ there exists sequences $\{a_n\}, \{b_n\}$ in D such that $d = a_n b_n$ and $\lim_{n \to a} \min\{\lambda(a_n), \lambda(b_n)\} = \lambda(d)$.*

*Proof:* (i) $\Rightarrow$ (ii).

Because

$$\sup_{d=ab} \min\{\lambda(a), \lambda(b)\} = \lambda \bullet \lambda(d) = \lambda(d),$$

given by positive integer n there exists $a_n, b_n, \in D$ such that

$$d = a_n b_n \text{ and } \lambda(d) - \frac{1}{n} < \min \{\lambda(a_n), \lambda(b_n)\} \leq \lambda(d).$$

(ii) $\Rightarrow$ (i). If possible let there be $d \in D$ such that $\lambda \bullet \lambda (d) < \lambda (d)$. There exist sequences $\{a_n\}$, $\{b_n\}$ in D such that $d = a_n b_n$ and that $\lim_{n \to a} \min \{\lambda(a_n), \lambda(b_n)\} = \lambda(d)$.

On the other hand, we have $\min \{\lambda (a_n), \lambda (b_n) \} \leq \lambda \bullet \lambda (d)$ for all n. Therefore $\lambda (d) \leq \lambda \bullet \lambda (d)$ which is a contradiction so we must have $\lambda \bullet \lambda = \lambda$. Hence the claim.

**THEOREM 1.11.45:** *Let $\lambda \in F_1 (D)$ where D has an identity element. If $\alpha$ is the upper bound of $\lambda$ in [0, 1] then the following two statements are equivalent.*

  i.  *Given $d \in D$, there exists sequences $\{a_n\}$, $\{b_n\}$ in D such that $d = a_n b_n$ and $\lim_{n \to a} \min \{\lambda(a_n), \lambda(b_n)\} = \lambda(d)$.*
  ii. $\left(\lambda^{-1}[\tau,\alpha]\right)^2 = \lambda^{-1}[\tau,\alpha]$ *whenever $0 \leq \tau < \alpha$.*

*Proof:* Refer [110].

**THEOREM 1.11.46:** *Let $\lambda \in F(D)$ where D is a semigroup in which the cancellation laws are valid. Then $\lambda \in F_2(D)$ if and only if*

$$\bigcap_{n=1}^{\infty} \lambda^n = \bigcup_{n=1}^{\infty} \lambda^n.$$

*Proof:* Firstly if $\lambda = \lambda^2$, then $\lambda = \lambda^2 = \lambda^3 = \ldots$ and hence $\cap \lambda^n = \lambda = \cup \lambda^n$. Conversely if $\cap \lambda^n \subseteq \cup \lambda^n$ then $\lambda \subseteq \cup \lambda^n = \cap \lambda^n \subseteq \lambda^2 \subseteq \cup \lambda^n = \cap \lambda^n = \lambda$, therefore $\lambda = \lambda^2$.

**THEOREM 1.11.47:** *Let $\{\lambda_i\} \subseteq F_2 (D)$ Then $\bigcup_i \lambda_i \in F_2(D)$ if and only if*

$$\bigcup_{i \neq j}(\lambda_i, \lambda_j) \subseteq \bigcup_i \lambda_i$$

*Proof:* $\cup_i \lambda_i \in F_2(D)$ if and only if

$$\begin{aligned}
\bigcup_i \lambda_i &= \left(\bigcup_i \lambda_i\right) \bullet \left(\bigcup_i \lambda_i\right) \\
&= \bigcup_{i,j}(\lambda_i, \lambda_j) \\
&= \left(\bigcup_i \lambda_i^2\right) \cup \left(\bigcup_{i \neq j}(\lambda_i \bullet \lambda_j)\right) \\
&= \left(\bigcup_i \lambda_i^2\right) \cup \left(\bigcup_{i \neq j}(\lambda_i \bullet \lambda_j)\right)
\end{aligned}$$

if and only if $\bigcup_{i \neq j}(\lambda_i . \lambda_j) \subseteq \bigcup_i \lambda_i$.

**THEOREM 1.11.48:** *Let $\{\lambda_i\}_i \subseteq F_2(D)$ where D is a groupoid in which the cancellation laws are valid. Then $\bigcap_i \lambda_i \in F_2(D)$ if and only if $\left(\bigcap \lambda_i\right) \subseteq \bigcap_{i \neq j}(\lambda_i \bullet \lambda_j)$.*

*Proof:* $\bigcap_i \lambda_i \in F_2(D)$ if and only if

$$\left(\bigcap_i \lambda_i\right) = \left(\bigcap_i \lambda_i\right)\left(\bigcap_i \lambda_i\right) = \bigcap_{i,j}(\lambda_i \bullet \lambda_j)$$
$$= \left(\bigcap_i (\lambda_i^2)\right) \cap \left(\bigcap (\lambda_i \bullet \lambda_j)\right)$$
$$= \left(\bigcap_i \lambda_i\right) \cap \left(\bigcap_{i \neq j}(\lambda_i \bullet \lambda_j)\right)$$

if and only if $\bigcap_i \lambda_i \subseteq \bigcap_{i \neq j}(\lambda_i . \lambda_j)$. $F_2(D)$ is neither a meet sublattice nor a join sublattice of $F(D)$ as given by an example in [110].

We have the following theorem:

**THEOREM 1.11.49:** *Let $\{\lambda_i\}_i \subseteq F_2(D)$ where D is a semigroup in which the cancellation laws are valid. Then*

$$\left(\bigcap_i \lambda_i\right)_2 = \bigcap_{A=1}^{\infty} (\bigcap_i \lambda_i)^n = \left(\bigcap_i \lambda_i\right)_{-1}.$$

*In other words $F_2(D)$ is a complete meet sublattice of $F_{-1}(D)$, if D is a semigroup in which the cancellation laws are valid.*

*Proof*: Left for the reader to prove.

**THEOREM 1.11.50:** *Let $\{\lambda_i\}_i \subseteq F_2(D)$ where D is a semigroup. Then*

$$\left(\bigcup_i \lambda_i\right)_2 = \bigcap_{n=1}^{\infty} (\bigcup_i \lambda_i)^n = \left(\bigcup_i \lambda_i\right)_1.$$

*In other words $F_2(D)$ is a complete join sublattice of $F_1(D)$, if D is a semigroup.*

*Proof*: Refer [110].

Using the above two theorems we get the following theorem which can be easily proved.

**THEOREM 1.11.51:** *If D is a semigroup in which the cancellation laws are valid, then $F_2(D)$ is a completely lattice under meet $\lambda \wedge \mu$ and join $\lambda \vee \mu$ given by*

$$\lambda \wedge \mu = |\lambda \cap \mu|_{-1} = \bigcap_{n=1}^{\infty} (\lambda \cap \mu)^n$$

and

$$\lambda \cup \mu = |\lambda \cup \mu|_{1} = \bigcup_{n=1}^{\infty} (\lambda \cup \mu)^n$$

Now we proceed on to recall the fuzzy subgroupoid as given by [47]. In this study he uses the concept of t-norm as given by [150].

**DEFINITION [150]**: *A t-normal is a function T : [0,1] × [0,1] → [0,1] satisfying for each x, y and z in [0,1]*

    i.   *T (x,1) = x.*
    ii.  *T (x, y) ≤ T (z x) if x ≤ z.*
    iii. *T (x, y) = T (y, x).*
    iv. *T (x, T(y,z)) = T(T(x, y), z).*

*A few t-norms which are frequently encountered are $T_m$, Prod and Min defined by $T_m$ (x, y) = $M_{ax}$ (x + y –1, 0), Prod (x, y) = xy and*

$$\text{Min } (x, y) = \begin{cases} x & \text{if } x \leq y \\ y & \text{if } y < x \end{cases}.$$

**DEFINITION 1.11. 20:** *A t-norm, $T_1$ is stronger than a t-norm $T_2$ if and only if $T_1$ (x, y) ≥ $T_2$ (x, y) for all x, y ∈ [0,1].*

**THEOREM 1.11.52:** *Min is the strongest of all t-norms.*

*Proof*: Directly proved.

Let X = (X, •) be a groupoid we will denote x • y by xy.

**DEFINITION 1.11.21:** *Let X = (X, •) be a groupoid. A function µ: X → [0, 1] is a fuzzy subgroupoid of X with respect to a t-norm T, if and only if for every x, y ∈ X we have µ (xy) ≥ T (µ(x), µ(y)).*

**DEFINITION 1.11.22:** *Let $X_1$ and $X_2$ be groupoid and let $µ_1$ and $µ_2$ be fuzzy subgroupoids of $X_1$ and $X_2$ respectively, with respect to a t-norm T. The fuzzy subgroupoids $µ_1$ and $µ_2$ are homomorphic (isomorphic) if and only if there exists a groupoid homomorphism (isomorphism) ϕ: $X_1$ → $X_2$ such that $µ_1$ = $µ_2$ o ϕ.*

In this situation we say that $µ_1$ is given by pull back of $µ_2$ along ϕ.

**THEOREM 1.11.53:** *Let µ be a fuzzy subgroupoid of X with respect to Min then µ is a fuzzy subgroupoid of X with respect to any t-norm T.*

*Proof*: Left as an exercise for the reader.

**DEFINITION 1.11.23:** *Let (X, •) be a groupoid. Let G denote the family of subgroupoids of X. G is called a generating family if for every element x ∈ X there exists a S ∈ G such that x ∈ S. {X} is a trivial generating family.*

**DEFINITION 1.11.24:** *Let (X, •) be a groupoid. An element e ∈ X is called an identity of X, if ex = xe = x for all x ∈ X.*

**THEOREM 1.11.54:** *Let (X, •) be a groupoid, then X has at most one identity.*

*Proof*: Obvious from the fact if e and e' are two identities e = ee' = e'.

**THEOREM 1.11.55:** *Every groupoid is a subgroupoid of some groupoid with an identity.*

*Proof*: Let (X, •) be a groupoid. Let e be an element such that e ∉ X. Define an operation 'o' on the set X' = X ∪ {e} as follows x o y = x • y if x, y ∈ X; x o e = e o x = x for all x ∈ X, e o e = e. Then obviously (X, o) is a groupoid with identity e and (X, •) is a subgroupoid of (X', o).

**DEFINITION 1.11.25:** *We call (X', o, e) the identity extension of (X, •) by e. Let (X, •) be a groupoid. By the virtue of the above theorem we may assume that X has an identity whenever necessary, without loss of generality.*

**THEOREM 1.11.56:** *Let (X, o) be a groupoid then there exists a unique identity extension of X up to isomorphism.*

*Proof*: Let (X', o, e) be constructed as in the proof of the earlier theorem. Suppose that (X'', ⊙, υ) is also an identity extension of X such that X'' = X ∪ {u} with identity element u. It is easily shown that φ: X' → X'' defined by

$$\phi(x) = \begin{cases} x & \text{if } x \in X \\ u & \text{if } x = e \end{cases}$$

is an isomorphism of groupoid X' onto X''.

**THEOREM 1.11.57**: *Let (X, •) be a groupoid. Let μ : X → [0, 1] be a fuzzy subgroupoid of X with respect to a t-norm T. Then μ can be extended to a fuzzy subgroupoid μ' of X' with respect to the same t-norm T, where X' is the identity extension of X.*

*Proof*: Let (X', o, e) be constructed as in the earlier theorem Let μ' : X' → [0, 1] be defined as

$$\mu'(x) = \begin{cases} \mu(x) & \text{if } x \in X \\ 1 & \text{if } x = e \end{cases}.$$

Now for any x, y ∈ X' if x, y ∈ X then x o y = x • y ∈ X; hence μ'(x o y) = μ(x • y) ≥ T (μ(x), μ(y)) = T (μ'(x), μ'(y)) since μ is a fuzzy subgroupoid of X with respect to the t-norm T. Further more for any x ∈ X'; μ' (xoe) = μ'(x) = T (μ'(x), 1) = T (μ'(x), μ'(e)).

Therefore μ' is a fuzzy subgroupoid of X' with respect to the same t-norm μ' $|_X = μ$. It $(Ω, \mathcal{A}, P)$ is a probability space then P is infact a fuzzy subset of $\mathcal{A}$. Let Ω = [0,1]. Let L be the set of all Lebesgue measurable subsets in [0,1] and P be the Lebesque measure then (Ω, L, P) is a probability space.

**THEOREM 1.11.58:** *Let $(Ω, \mathcal{A}, P)$ be a probability space then $(\mathcal{A}, \cap)$ becomes a groupoid and P is a fuzzy subgroupoid of $\mathcal{A}$ with respect to $T_m$.*

*Proof:* For any A, B ∈ $\mathcal{A}$, we have P(A ∩ B) = P(A) + P(B) − P(A∪B) ≥ P(A) + P(B) − 1. Since P(A ∩ B) ≥ 0, we obtain P (A ∩ B) ≥ Max (P(A) + P(A) − 1,0) = $T_m$(P(A), P(B)).

**THEOREM 1.11.59:** *Let (X, •) be a groupoid and let $(Ω, \mathcal{A}, P)$ be a probability space. Suppose that the mapping, ψ: X → A satisfies ψ (x y) ⊇ ψ(x) ∩ ψ(y). Set μ(x) = P(ψ(x)) for all x ∈ X, then μ is a fuzzy subgroupoid of X with respect to $T_m$. One says that μ is induced by the mapping ψ and the probability space $(Ω, \mathcal{A}, P)$.*

*Proof*: For any x, y ∈ X we have ψ (x y) ⊇ ψ(x) ∩ ψ (y) therefore by the above theorem it follows that μ(xy) = P (ψ (xy)) ≥ P (ψ (x) ∩ ψ (y)) ≥ $T_m$ (P(ψ(x), P (ψ(y))) = $T_m$(μ(x), μ(y)).

The following result is direct hence left for the reader to prove.

**THEOREM 1.11.60:** *Let (X, •) be a groupoid and let $(Ω, \mathcal{A}, P)$ be a probability space. If μ: X → $\mathcal{A}$ is a homomorphism from (X, •) to $(\mathcal{A}, \cap)$. Set μ(x) = P (ψ(x)), then μ is a fuzzy subgroupoid of X with respect to $T_m$ induced by ψ and $(Ω, \mathcal{A}, P)$.*

**THEOREM 1.11.61:** *Let (X, •) be a groupoid and G be a generating family of X. For x ∈ X let $S_x$ = {S ∈ G $/$ x ∈ S} and let W = {$S_x$ $/$ x ∈ X}. Let A be any σ -algebra on G which contains the σ- algebra generated by W and let m be a probability measure on (G, A). Then μ: X → [0,1] defined by μ(x) = m($S_x$) for each x ∈ X is a fuzzy subgroupoid of X with respect to $T_m$. A fuzzy subgroupoid obtained in this manner is called "subgroupoid generated".*

*Proof*: Let x, y ∈ X. Suppose that S ∈ $S_x$ ∩ $S_y$; Then S is a subgroupoid of X containing both x and y. Consequently x y ∈ S and hence S ∈ $S_x$. Therefore $S_{xy}$ ⊇ $S_x$ ∩ $S_y$. Define ψ : X → A by x →$S_x$. Then ψ (xy) = $S_{xy}$ ⊃ $S_x$ ∩ $S_y$ = ψ(x) ∩ ψ(y). Set μ(x) = m (ψ(x)) = m($S_x$) then by earlier theorem μ is a fuzzy subgroupoid of X with respect to $T_m$.

**THEOREM 1.11.62:** *Let (X, •) and all other symbols be the same as described in the theorem 1.11.61. If there exists an element e of X such that e is contained in every subgroupoid of the generated family G, then μ(e) = 1.*

*Proof*: Obviously $S_e = G$ and hence $\mu(e) = m(S_e) = m(G) = 1$.

The following result is straightforward.

**THEOREM 1.11.63:** *Let all symbols be the same as described in the above theorem. If $(X, \bullet)$ is a monoid with the identity e then $\mu(e) = 1$.*

Further we have a still stronger result.

**THEOREM 1.11.64:** *If $(X, \bullet)$ is a group. Let e be its identity, then $\mu(e) = 1$ and $\mu(x^{-1}) = \mu(x)$ for all $x \in X$.*

*Proof*: Note that $S_{x^{-1}} = S_x$ and hence $\mu(x^{-1}) = m(S_{x^{-1}}) = m(S_x) = \mu(x)$. It is pertinent to note that if $(X, \bullet)$ is a group then $\mu$ obtained in the Theorem 1.11.64 is just the "subgroup generated"; fuzzy subgroup of X.

**THEOREM 1.11.65:** *Let $(\Omega, \mathcal{A}, P)$ be a probability space and $\tau$ be a non empty set. Let $\phi : \Omega \to \tau$ be a mapping. Then $\phi$ induces a new probability space $(\tau, A, m)$ where $A = \{B \subseteq \tau \,/\, \phi^{-1}(B) \in \mathcal{A}\}$ and for each $B \in A$, $m(B) = P(\phi^{-1}(B))$.*

*We say that $(\tau, A, m)$ is induced by $\phi$ and $(\Omega, \mathcal{A}, P)$.*

**THEOREM 1.11.66:** *Let $(\Omega, \mathcal{A}, P)$ be a probability space and $\tau$ be a non empty set. $P(\tau)$ denotes the power set of $\tau$. Let $\phi : \Omega \to \tau$ be a mapping. Suppose that $(X, \bullet)$ is a groupoid and there exists a mapping $\psi : X \to P(\tau)$ such that*

   i. $\phi^{-1}(\psi(x)) \in \mathcal{A}$ for all $x \in X$
   ii. $\psi(x y) \supseteq \psi(x) \cap \psi(y)$ for all $x, y \in X$.

*Set $\mu(x) = P(\phi^{-1}(\psi(x)))$ for all $x \in X$, then $\mu$ is a fuzzy subgroupoid of X with respect to $T_m$.*

*Proof*: By the above result $\phi$ and $(\Omega, \mathcal{A}, P)$ induce a new probability space $(\tau, A, m)$. Result follows by the earlier theorems.

The following result is a direct consequence of the theorems proved; hence left for the reader as an exercise.

**THEOREM 1.11.67:** *Let $\Omega = [0,1]$ and L be the set of all Lebesgue measurable sets on $[0,1]$. Let P by the Lebesgue measure. Suppose that we have a non empty set $\tau$, a groupoid $(X, \bullet)$ and two mappings as follows:*
*$\phi : \Omega \to \tau$, $\psi : X \to P(\tau)$ satisfying*

   i. $\phi^{-1}(\psi(x)) \in L$ for all $x \in X$.
   ii. $\psi(xy) \supseteq \psi(x) \cap \psi(y)$ for all $x \in X$.

*Set $\mu(x) = P(\phi^{-1}(\psi(x)))$ for all $x \in X$, then $\mu$ is a fuzzy subgroupoid of X with respect to $T_m$.*

Now we proceed on to recall properties of a fuzzy subgroupoid of a subgroupoid of a direct product.

**THEOREM 1.11.68:** *Let $(X, \bullet)$ be a groupoid and let Y be a fixed subgroupoid of X. Let $(\Omega, \mathcal{A}, P)$ be a probability space and let $(F, \odot)$ be a groupoid of functions mappings $\Omega$ into X and $\odot$ defined by pointwise multiplication in the range space. A further restriction is placed on F by assuming that for each $f \in F$, $X_f = \{w \in \Omega \mid f(w) \in Y\}$ is an element of $\mathcal{A}$. Then $v : F \to [0,1]$ defined by $v(f) = P(X_f)$ for each f in F is a fuzzy subgroupoid of F with respect to $T_m$. A fuzzy subgroupoid obtained in this manner is called "function generated".*

*Proof*: Let f, g $\in$ F. Suppose $\omega \in X_f \cap X_g$. Then $f(\omega) \in$ Y. Since Y is a subgroupoid of X. (f o g)($\omega$) = f($\omega$) o g($\omega$) $\in$ Y and consequently $\omega \in X_{fog}$. Therefore $X_{fog} \supseteq X_f \cap X_g$. Now consider the mapping $\psi : F \to \mathcal{A}$. $f \to X_f$. We have $\psi$ (f o g) = $X_{fog} \supseteq X_f \cap X_g = \psi$ (f) $\cap \psi$(g). Set $v(f) = P(\psi(f)) = P(X_f)$ for all $f \in$ F, then by earlier results $v$ is a fuzzy subgroupoid of F with respect to $T_m$.

It is pertinent to mention that for each $\omega \in \Omega$ by defining $X_\omega = X$, we can construct a direct product $\prod_{\omega \in \Omega} X_\omega$ then (F, $\odot$) is in fact a subgroupoid of $\prod_{\omega \in \Omega} X_\omega$.

If $f \in$ F and $f(\omega) \in$ Y, for all $\omega \in \Omega$ then $v$ (f) = 1.

Now we recall the fuzzy subgroupoid representations.

**THEOREM 1.11.69:** *Every function generated fuzzy subgroupoid is subgroupoid generated.*

*Proof*: The proof is lengthy the reader is requested to refer [47].

The proof of the following theorem is omitted and the reader is advised to refer [47].

**THEOREM 1.11.70:** *Every subgroupoid generated fuzzy subgroupoid is isomorphic to a function generated fuzzy subgroupoid.*

**THEOREM 1.11.71:** *Let $(X, \bullet)$ be a groupoid with the identity e and let $v$ be a fuzzy subgroupoid of X with respect to Min such that $v(e) = 1$. For each t in [0,1] let $X_t = \{x \in X, v(x) \geq t\}$. Let a and b be in [0,1] such that $a \leq b$; then*

   i. $X_t$ *is a subgroupoid of X for every t in [0,1].*
   ii. $X_b \subseteq X_a$.
   iii. *if $x \in X_a - X_b$ for all b in (a, 1]; then $v(x) = a$.*

*Proof*:

   i. For each t in [0,1] since $v(e) = 1 \geq t$ we have $e \in X_t$. Hence $X_t$ is not empty. For any x, y $\in X_t$ it follows that $v(x) \geq t$, $v(y) \geq t$. Therefore $v(xy) \geq$ Min$\{v(x), v(y)\} \geq t$; hence x y $\in X_t$. This shows that $X_t$ is a subgroupoid of X.

ii. For any x ∈ $X_b$ we have ν(x) ≥ b. Since b ≥ a it follows that ν(x) ≥ a and so x ∈ $X_a$. Hence $X_b \subseteq X_a$.

iii. If x ∈ $X_a - X_b$ for all b in (a, 1] then we have a ≤ ν(x) < b for all a < b ≤ 1. Hence ν(x) = a.

Finally we state the following two theorems from [47] the proof is left for the reader to refer.

**THEOREM 1.11.72:** *Every fuzzy subgroupoid with respect to Min is subgroupoid generated.*

**THEOREM 1.11.73:** *Let μ be a subgroupoid generated fuzzy subgroupoid with (X, •) and (G, A, m) and set $S_x$ for x ∈ X, as described in the theorem 1.11.72. If there exists $G^* ∈ A$ which is linearly ordered by set inclusion such that $m(G^*) = 1$ then μ is a fuzzy subgroupoid with respect to Min.*

*The converse of the above theorem holds if (X, •) has an identity e such that μ(e) = 1. In this case we set $G^* = \{\chi_t \mid t \in [0,1]\}$ then from the earlier result it is obvious that $G^*$ is linearly ordered with measure one.*

Some more algebraic properties are discussed in the following:

**DEFINITION 1.11.26:** *Give a groupoid. X = (X, •) a t-norm T and set I, for each i ∈ I, let $μ_i$ be a fuzzy subgroupoid of X with respect to T, we define*

$$\left(\bigcap_{i \in I} \mu_i\right)(x) = \inf_{i \in I} \{\mu_i(x)\}.$$

**THEOREM 1.11.74:** $\bigcap_{i \in I} \mu_i$ *(x) is a fuzzy subgroupoid of X with respect to T.*

*Proof*: For any x, y ∈ X. We have

$$\mu_i(x) \geq \inf_{i \in I} \{\mu_i(x)\}$$

$$\mu_i(y) \geq \inf_{i \in I} \{\mu_i(y)\}$$

hence

$$T(\mu_i(x), \mu_i(y)) \geq T(\inf \mu_i(x), \inf \mu_i(y))$$

for all i ∈ I.

It follows that

$$\inf_{i \in I}[T(\mu_i(x), \mu_i(y))] \geq \left[\inf_{i \in I} \mu_i(x), \inf_{i \in I} \mu_i(y)\right]$$

Therefore we have

$$\left(\bigcap_{i \in I} \mu_i\right)(xy) = \inf_{i \in I} \{\mu_i(xy)\}$$

$$\geq \inf_{i \in I} [T(\mu_i(x), \mu_i(y))]$$

$$= T\left[\left(\bigcap_{i \in I} \mu_i\right)(x), \left(\bigcap_{i \in I} \mu_i\right)(y)\right]$$

Hence the result.

It is evident that fuzzy subgroupoids of X with respect to a t-norm T form a complete lattice. In this lattice, the inf of a set of fuzzy subgroupoids $\mu_i$ with respect to T is just $\cap \mu_i$, while their sup is the least $\mu$ (i.e. the $\cap$ of all $\mu$'s) which $\supseteq \cup \mu_i$, where $\cup \mu_i (x) = \sup \mu_i (x)$ for all $x \in X$.

**DEFINITION 1.11.27:** *Let $f: X \to X'$ be a homomorphism from groupoid X into groupoid X'. Suppose that $\nu$ is a fuzzy subgroupoid of X' with respect to a t-norm T. Then the fuzzy set $\mu = \nu \circ f$ (defined by $\mu(x) = \nu(f(x))$ for all $x \in X$) is called the pre image of $\nu$ under f.*

**THEOREM 1.11.75:** *Let $\mu = \nu \circ f$ be the preimage of $\nu$ under f as described in the above definition, then $\mu$ is a fuzzy subgroupoid of X with respect to T.*

*Proof*: For $x, y \in X$ we have $\mu(xy) = (\nu \circ f)(x y) = \nu(f(x y)) = \nu(f(x), f(y)) \geq T(\nu(f(x)), \nu(f(y))) = T(\mu(x), \mu(y))$. So by definition $\mu$ is a fuzzy subgroupoid of X with respect to T.

**DEFINITION 1.11.28:** *Let $f: X \to X'$ be a homomorphism from groupoid X onto groupoid X'. Suppose that $\mu$ is a fuzzy subgroupoid of X with respect to a t-norm T, then the fuzzy set $\nu$ in $X' = f(X)$ defined by*

$$\nu(y) = \sup_{x \in f^{-1}(y)} \mu(x)$$

*for all $y \in X'$ is called the image of $\mu$ under f.*

**DEFINITION [112]:** *We say that a fuzzy set $\mu$ in X has the sup property if for any subset $A \subseteq X$ there exists $a_o \in A$ such that $\mu(a_0) = \sup_{a \in A} \mu(a)$*

**THEOREM 1.11.76:** *Let $f: X \to X'$ be a homomorphism from groupoid X on to groupoid X. Suppose that $\mu$ is a fuzzy subgroupoid of X with respect to a t-norm T and that $\mu$ has the sup property. Let $\nu$ be a homomorphic image of $\mu$ under f as described in the earlier definition then $\nu$ is a fuzzy subgroupoid of X' with respect to T.*

*Proof:* Given $f(x), f(y)$ in $f(X) = X'$, let $x_o \in f^{-1}(f(x))$, $y_o \in f^{-1}(f(x))$ be such that

$$\mu(x_o) = \sup_{t \in f^{-1}(f(x))} \mu(t)$$

$$\mu(y_o) = \sup_{t \in f^{-1}(f(y))} \mu(t)$$

respectively. Then we have

$$\nu(f(x)) = \sup_{t \in f^{-1}(f(x))} \mu(t) = \mu(x_o)$$

$$\nu(f(y)) = \sup_{t \in f^{-1}(f(y))} \mu(t) = \mu(y_o).$$

since

$$f(x_o, y_o) = f(x_o)\, f(y_o)$$
$$= f(x)\, f(y)$$

we have $x_o, y_o \in f^{-1}(f(x), f(y))$ hence

$$\sup_{z \in f^{-1}(f(x),f(y))} \mu(z) \geq \mu(x_o, y_o).$$

It follows that

$$\nu(f(x), f(y)) = \sup_{z \in f^{-1}(f(x),f(y))} \mu(z) \geq \mu(x_o, y_o)$$
$$\geq T(\mu(x_o), \mu(y_o))$$
$$\geq T(\nu f(x), \nu f(y)).$$

Therefore $\nu$ is a fuzzy subgroupoid of X' with respect to T.

Now we proceed on to define fuzzy bigroupoid. The concept of bigroupoid is itself very new so the notion is fuzzy bigroupoid is totally absent in literature. Here we define fuzzy bigroupoid and recall just the definition of bigroupoid.

**DEFINITION 1.11.29:** *Let $(G, +, \bullet)$ be a non empty set we call G a bigroupoid if $G = G_1 \cup G_2$ and satisfies the following:*

    i.    *$(G_1, +)$ is a groupoid.*
    ii.   *$(G_2, \bullet)$ is a semigroup.*

*Example 1.11.16:* Let $(G, +, \bullet)$ be a groupoid where $G = G_1 \cup G_2$ with $G_1$ a groupoid given by the following table:

| + | $X_1$ | $X_2$ | $X_3$ |
|---|---|---|---|
| $X_1$ | $X_1$ | $X_3$ | $X_2$ |
| $X_2$ | $X_2$ | $X_1$ | $X_3$ |
| $X_3$ | $X_3$ | $X_2$ | $X_1$ |

And $G_2 = S(3)$ the symmetric semigroup of mappings of (123) to itself.

Clearly (G, +, •) is a bigroupoid.

Now we define fuzzy bigroupoid as follows.

**DEFINITION 1.11.30:** *Let (G, +, •) be a bigroupoid; $G = G_1 \cup G_2$ proper subsets of G; with ($G_1$, +) a groupoid and $G_2$ a semigroup. $\mu : G \rightarrow [0, 1]$ is said to be a fuzzy bisubgroupoid (or by default of notation fuzzy bigroupoid) if and only if $\mu = \mu_1 \cup \mu_2$ where $\mu_1$ from $G_1 \rightarrow [0, 1]$ is a fuzzy subgroupoid and $\mu_2: G_2 \rightarrow [0, 1]$ is a fuzzy subsemigroup.*

Thus almost all properties true in case of fuzzy groupoids and fuzzy semigroups can be easily extended to the case of fuzzy bisubgroupoids.

Now we proceed on to define fuzzy loops and fuzzy biloops.

**DEFINITION [118]:** *Let G be a group. A fuzzy subset $\mu : G \rightarrow [0, 1]$ is called a fuzzy subloop if for at least a triple x, y, z $\in$ G we have $\mu((xy)z) \neq \mu(x(yz))$.*

Study of this type was carried out by [118].

Since all groups are loops we can define fuzzy subloop in a more generalized way as

**DEFINITION [118]:** *Let L be a loop A: $L \rightarrow [0, 1]$ is a fuzzy subloop of L if atleast for a triple x, y z $\in$ L we have $A((xy)z) \neq A(x(yz))$.*

**DEFINITION [118]:** *Let L be a group or a loop. The fuzzy subloop V: $L \rightarrow [0, 1]$ is called a fuzzy Bruck subloop if $P(x(yx)z) = P(x(y(xz)))$ and $P(xy)^{-1} = P(x^{-1} y^{-1})$ for all x, y $\in$ L.*

**DEFINITION [118]:** *Let G be a group or a loop. A fuzzy subloop P of G (P: $G \rightarrow [0, 1]$ is a fuzzy Bol subloop of G if $P(((xy)z)y) = P(x((yz))y)$ for all x, y, z $\in$ G.*

**DEFINITION [118]:** *Let L be a loop or a group. A fuzzy subloop P: $L \rightarrow [0, 1]$ is called a fuzzy Moufang subloop of L if $P((xy(zx)) = P((x(yz))x)$ for all x, y, z $\in$ L.*

**THEOREM 1.11.77:** *Every fuzzy subloop of a loop L need not in general be a fuzzy Moufang subloop of L.*

*Proof*: Consider the loop L given by the following table:

| * | e | 1 | 2 | 3 | 4 | 5 |
|---|---|---|---|---|---|---|
| e | e | 1 | 2 | 3 | 4 | 5 |
| 1 | 1 | e | 5 | 4 | 3 | 2 |
| 2 | 2 | 3 | e | 1 | 5 | 4 |
| 3 | 3 | 5 | 4 | e | 2 | 1 |
| 4 | 4 | 2 | 1 | 5 | e | 3 |
| 5 | 5 | 4 | 3 | 2 | 1 | e |

Define $V : L \to [0\ 1]$ by $V(0) = 0$, $V(i) = 0.i$ for all $i \in V$. $V((14)(21)) = 0.2 = V(2)$ $V((1(42)1)) = 0 = V(0)$. So V is not fuzzy Moufang subloop of L.

**THEOREM 1.11.78:** *Let L be a loop or a group. The fuzzy subloop V of $L \to [0\ 1]$ is a right alternative fuzzy subloop of L if $V((xy)y) = V(x(yy))$ for all $x, y \in L$.*

*Proof*: Straightforward by the very definition.

**THEOREM 1.11.79:** *Let L be a loop or a group. A fuzzy subloop V of L is a weak inverse property loop if $V((xy)z) = 0$ imply $V(x(yz)) = 0$ for all $x, y, z \in L$.*

*Proof*: Given $V: L \to [0, 1]$ is a fuzzy subloop of L, clearly if $V((xy)z) = 0$ imply $V(x(yz)) = 0$ then it is the weak inverse property as the result is true for all $x, y, z \in L$ we have V to be a fuzzy subloop which satisfies weak inverse property.

As we are not able to find any other means to define fuzzy subloop as to the best of the authors knowledge we do not have any other definition for fuzzy subloops we take this as the basic definition. Also we wish to mention the definition of fuzzy subloops given by us is distinctly different from the classical definitions of other algebraic structure.

Also it is pertinent to mention here that all fuzzy subloops are not subloops for they do not satisfy the axiom of a subloop.

Now we proceed on to define the notion of fuzzy biloop; for which first we give the definition of biloop.

**DEFINITION 1.11.31:** *Let $(L, +, \bullet)$ be a non empty set with two binary operations. L is said to be a biloop if L has 2 nonempty finite proper subsets $L_1$ and $L_2$ of L such that*

    i.    $L = L_1 \cup L_2$.
    ii.    $(L_1, +)$ *is a loop.*
    iii.    $(L_2, \bullet)$ *is a loop or a group.*

We define fuzzy biloop and request the reader to refer [135] for more properties about biloops.

**DEFINITION 1.11.32:** *Let $(L, +, \bullet)$ be a biloop. A map $\mu : L \to [0, 1]$ is called fuzzy biloop if we have $(L = L_1 \cup L_2)$ with $(L_1, +)$ a loop and $(L_2, \bullet)$ group or a loop) $\mu = \mu_1 \cup \mu_2$ where $\mu_1$ denotes the fuzzy subloop of $L_1$ and $\mu_2$ is the fuzzy subloop of $L_2$ (i.e.*

*µ restricted to $L_1$ is denoted by $µ_1$ and µ restricted to $L_2$ denoted by $µ_2$ the symbol '∪' in $µ = µ_1 \cup µ_2$ is only by default of notation).*

All fuzzy properties studied for loops can be easily extended to fuzzy biloops in an analogous and routine ways.

## 1.12 Miscellaneous Properties in Fuzzy algebra

This section is devoted to the miscellaneous properties about fuzzy algebra which has not been covered in the earlier eleven sections like fuzzy continuous map on groups, fuzzy polynomial ring, fuzzy polynomial semiring, and finally the concept of fuzzy opposite sets and systems. Now we proceed on to recall the concept of fuzzy continuous map on groups. For more about these refer [139].

**DEFINITION 1.12.1:** *A fuzzy topology τ on a group G is called a g fuzzy topology. The pair (G, τ) is called a g-fuzzy topological space.*

*Example 1.12.1:* Let G = {1, –1} be the group with respect to usual multiplication and $\tau = \{\phi_G, 1_G, \lambda, \mu\}$ where λ, µ : G → [0, 1] are given by

$$\lambda(x) = \begin{cases} 1 & \text{if } x = 1 \\ 0 & \text{if } x = -1 \end{cases}$$

and

$$\mu(x) = \begin{cases} 0 & \text{if } x = 1 \\ 1 & \text{if } x = -1. \end{cases}$$

Now the empty fuzzy set $\phi_G$ and the whole fuzzy set $1_G$ are in τ. Further, it is easily verified that the intersection of any two members of τ is a member of τ and arbitrary union of members of τ is a member of τ. Hence τ is a g-fuzzy topology on G.

**DEFINITION 1.12.2:** *Let $\tau_1$ and $\tau_2$ be g-fuzzy topologies on the group $G_1$ and $G_2$ respectively. A function f : ($G_1$ $\tau_1$) → ($G_2$ $\tau_2$) is said to be a g-fuzzy continuous map from $G_1$ to $G_2$ if it satisfies the following conditions:*

  i. *For every $\mu \in \tau_2, f^{-1}(\mu) \in \tau_1$ and*
  ii. *For every fuzzy subgroup µ (of $G_2$) in $\tau_2$, $f^{-1}$ (µ) is a fuzzy subgroup (of $G_1$) in $\tau_1$.*

This definition is illustrated by the following example.

*Example 1.12.2:* Let $G_1$ = {1, –1, i, –i} be a group with respect to usual multiplication and $G_2$ = {e, a, b, ab} be the Klein four group, where $a^2 = e = b^2$, ab = ba and e is the identity element of the group G. The corresponding g-fuzzy topologies are given by $\tau_1 = \{\phi_{G_1}, 1_{G_1}, \lambda_1, \mu_1\}$ and $\tau_2 = \{\phi_{G_2}, 1_{G_2}, \lambda_2, \mu_2\}$ where $\lambda_1, \mu_1 : G_1 \to [0, 1]$ and $\lambda_2, \mu_2 : G_2 \to [0, 1]$ are defined as follows ;

$$\lambda_1(x) = \begin{cases} 1 & \text{if } x = 1, -1 \\ 0 & \text{if } x = i, -i \end{cases}$$

$$\mu_1(x) = \begin{cases} 1 & \text{if } x = i, -i \\ 0 & \text{if } x = 1, -1 \end{cases}$$

$$\lambda_2(x) = \begin{cases} 1 & \text{if } x = e \\ 0 & \text{if } x = a, b, ab \end{cases}$$

$$\mu_2(x) = \begin{cases} 0 & \text{if } x = e \\ 1 & \text{if } x = a, b, ab. \end{cases}$$

Define f: $(G_1 \tau_1) \to (G_2 \tau_2)$ by

$$f(x) = \begin{cases} e & \text{if } x = 1, -1 \\ b & \text{if } x = i \\ ab & \text{if } x = -i. \end{cases}$$

Then for every $x \in G_1$, we have calculated the following:

$$(f^{-1}(\phi_{G_2}))(x) = \phi_{G_1}(x)$$

$$(f^{-1}(1_{G_2}))(x) = 1_{G_1}(x)$$

$$(f^{-1}(\lambda_2))(x) = \lambda_1(x)$$

$$(f^{-1}(\mu_2))(x) = \mu_1(x).$$

Hence we have

$$(f^{-1}(\phi_{G_2})) = \phi_{G_1}, (f^{-1}(1_{G_2})) = 1_{G_1}, (f^{-1}(\lambda_2)) = \lambda_1 \text{ and } (f^{-1}(\mu_2)) = \mu_1.$$

This proves $f^{-1}(\mu) \in \tau_1$ for every $\mu \in \tau_2$. Further, it is easy to verify that $f^{-1}(\mu)$ is a fuzzy subgroup of the group $G_1$ whenever $\mu$ is a fuzzy subgroup of the group $G_2$. Hence f is a g-fuzzy continuous map from $G_1$ to $G_2$.

Both the conditions given in definition 1.12.2 are essential for in the above example if we take $\tau_1 = \{\phi_{G_1}, 1_{G_1}, \lambda_1\}$ instead of $\tau_1 = \{\phi_{G_1}, 1_{G_1}, \lambda_1, \mu_1\}$ then clearly condition (ii) of definition 1.12.2 holds good, but not condition (i) as $f^{-1}(\mu_2) \notin \tau_1$. Hence both the conditions given in the definition are essential.

**THEOREM 1.12.1:** *Let $G_1$ and $G_2$ be any two groups. If $\tau_1$ is a g-fuzzy topology on the group $G_1$ and $\tau_2$ is an indiscrete g-fuzzy topology on the group $G_2$ then every function $f: (G_1, \tau_1) \to (G_2, \tau_2)$ is a g-fuzzy continuous map.*

*Proof*: A fuzzy topology $\tau$ is said to be an indiscrete fuzzy topology if its only elements are the empty fuzzy set and the whole fuzzy set. Let $\tau_1$ be a g-fuzzy topology on the group $G_1$ and $\tau_2$ be an indiscrete g-fuzzy topology on the group $G_2$. Since $\tau_2$ is an indiscrete g-fuzzy topology we have $\tau_2 = \{\phi_{G_2}, 1_{G_2}\}$.

Let $f : (G_1, \tau_1) \to (G_2, \tau_2)$ be any function. We see that every member of $\tau_2$ is a fuzzy subgroup of the group $G_2$. So it is enough to prove that for every $\mu \in \tau_2$, $f^{-1}(\mu) \in \tau_1$.

**Case i:** Let $\phi_{G_2} \in \tau_2$. Then for any $x \in G_1$,

$$f^{-1}(\phi_{G_2})(x) = \phi_{G_2}(f(x))$$
$$= 0 \ (\text{as } f(x) \in G_2)$$
$$= \phi_{G_1}(x)$$

(by the definition of empty fuzzy set) That is $(f^{-1}(\phi_{G_2}))(x) = \phi_{G_1}(x)$ for every $x \in G_1$. Thus we have $(f^{-1}(\phi_{G_2})) = \phi_{G_1} \in \tau_1$.

**Case ii:** Let $1_{G_2} \in \tau_2$ and $x \in G_1$, then we have

$$(f^{-1}(1_{G_2}))(x) = 1_{G_2}(f(x))$$
$$= 1 \ (\text{as } f(x) \in G_2)$$
$$= 1_{G_1}(x)$$

(by the definition of whole fuzzy set). Hence $(f^{-1}(\phi_{G_2}))(x) = 1_{G_1}(x)$ for every $x \in G_1$. This proves $(f^{-1}(1_{G_2})) = 1_{G_1} \in \tau_1$. Hence f is a g-fuzzy continuous map from $G_1$ to $G_2$.

The following example will illustrate the above theorem.

*Example 1.12.3:* Let $G_1 = \{1, -1, i, -i\}$ and $G_2 = \{1, -1\}$ be two groups with respect to usual multiplication. The corresponding g-fuzzy topologies are given by $\tau_1 = \{\lambda_1, \mu_1\}$ and $\tau_2 = \{\lambda_2, \mu_2\}$ where $\lambda_1(x) = 0$ for every $x \in G_1$, $\mu_1(x) = 1$ for every $x \in G_1$, $\lambda_2(x) = 0$ for every $x \in G_2$ and $\mu_2(x) = 1$ for every $x \in G_2$.

Let as consider the function $f: G_1 \to G_2$ defined by $f(x) = x^2$ for every $x \in G_1$. Thus we have for every $x \in G_1$

$$
\begin{aligned}
\left(f^{-1}(\lambda_2)\right)(x) &= \lambda_2(f(x)) \\
&= \lambda_2(x^2) \\
&= 0 \\
&= \lambda_1(x).
\end{aligned}
$$

Thus $\left(f^{-1}(\lambda_2)\right)(x) = \lambda_1(x)$ for every $x \in G_1$. Hence we have $f^{-1}(\lambda_2) = \lambda_1$. Further we have $f^{-1}(\mu_2) = \mu_1$. It is easy to check that for every function f, $f^{-1}(\lambda_2) = \lambda_1$ and $f^{-1}(\mu_2) = \mu_1$.

Clearly $\lambda_1$ and $\mu_1$ are fuzzy subgroups of the group $G_1$. Thus f is a g-fuzzy continuous map from $G_1$ to $G_2$. Clearly $\lambda_1$ and $\lambda_2$ are fuzzy subgroups of the group $G_1$. Thus f is a g-fuzzy continuous map from $G_1$ to $G_2$.

**THEOREM 1.12.2:** *Let $\tau_1$ and $\tau_2$ be any two discrete g-fuzzy topologies on the groups $(G_1, \bullet)$ and $(G_2, *)$ respectively. Then every group homomorphism f: $(G_1, \tau_1) \to (G_2, \tau_2)$ is a g-fuzzy continuous map but not conversely.*

*Proof*: We say a fuzzy topology $\tau$ on a set X to be a discrete fuzzy topology if it contains all fuzzy subsets of X. Let $\tau_1$ and $\tau_2$ be any two discrete g-fuzzy topologies on the groups $(G_1, \bullet)$ and $(G_2, *)$ respectively and f be a group homomorphism from $G_1$ to $G_2$. Since $\tau_1$ and $\tau_2$ are discrete g-fuzzy topologies by the definition of g-fuzzy topology we have for every $\mu \in \tau_2$; $f^{-1}(\mu) \in \tau_1$, here it is important to note that $f^{-1}$ is not the usual inverse homomorphism from $G_2$ to $G_1$.

Let $\mu$ be a fuzzy subgroup (of $G_2$) in $\tau_2$ then for x, y $\in$ G, we have

$$
\begin{aligned}
\left(f^{-1}(\mu)\right)(x \bullet y) &= \mu(f(x \bullet y)) \\
&= \mu(f(x) * f(y)) \text{ (since f is a group homomorphism)} \\
&\geq \min\{\mu(f(x)), \mu(f(y))\} \text{ (for } \mu \text{ is a fuzzy subgroup of } G_2) \\
&= \min\{\left(f^{-1}(\mu)\right)(x), \left(f^{-1}(\mu)\right)(y)\}.
\end{aligned}
$$

Hence we have $\left(f^{-1}(\mu)\right)(x \bullet y) \geq \min\{\left(f^{-1}(\mu)\right)(x), \left(f^{-1}(\mu)\right)(y)\}$ for every x, y $\in G_1$.

Further

$$
\begin{aligned}
\left(f^{-1}(\mu)\right)(x^{-1}) &= \mu(f(x^{-1})) \\
&= \mu(f(x)^{-1}) \text{ (since f is a group homomorphism)} \\
&= \mu(f(x)) \text{ (since } \mu \text{ is a fuzzy subgroup of the group } G_2) \\
&= \left(f^{-1}(\mu)\right)(x).
\end{aligned}
$$

Thus $\left(f^{-1}(\mu)\right)(x^{-1}) = \left(f^{-1}(\mu)\right)(x)$ for every $x \in G_1$. This proves that $f^{-1}(\mu)$ is a fuzzy subgroup (of $G_1$) in $\tau_1$ and hence f is a g-fuzzy continuous map from $G_1$ to $G_2$.

The following example shows that a g-fuzzy continuous map need not in general be a group homomorphism.

*Example 1.12.4:* Take $G_1 = \{1, -1\}$ and $G_2 = \{1, \omega, \omega^2\}$ to be groups under the usual multiplication with discrete g-fuzzy topologies $\tau_1$ and $\tau_2$ on $G_1$ and $G_2$ respectively.

Now define $f: (G_1, \tau_1) \to (G_2, \tau_2)$ by $f(1) = 1$ and $f(-1) = \omega$.

Since $\tau_1$ and $\tau_2$ are discrete g-fuzzy topologies, clearly for every $\mu \in \tau_2$, $f^{-1}(\mu) \in \tau_1$. Now we will prove that if $\mu$ is a fuzzy subgroup (of $G_2$) in $\tau_2$ then $f^{-1}(\mu)$ is a fuzzy subgroup (of $G_1$) in $\tau_1$. Let $\mu$ be a fuzzy subgroup (of $G_2$) in $\tau_2$ then we have the following two cases.

**Case (i):** If $\mu = \phi_{G_2}$ or $\mu = 1_{G_2}$ then clearly $f^{-1}(\mu)$ is a fuzzy subgroup (of $G_1$) in $\tau_1$.

**Case (ii):** If $\mu = \phi_{G_2}$ and $\mu \neq 1_{G_2}$ then $\mu$ has the following form:

$$\mu(x) = \begin{cases} t_1 & \text{if } x = 1 \\ t_2 & \text{if } x = \omega, \omega^2 \end{cases}$$

where $1 \geq t_1 > t_2 \geq 0$.

Now for every $x \in G_1$,

$$\left(f^{-1}(\mu)\right)(x) = \mu(f(x)) = \begin{cases} t_1 & \text{if } x = 1 \\ t_2 & \text{if } x = -1, \end{cases}$$

where $1 \geq t_1 > t_2 \geq 0$.

Hence $f^{-1}(\mu)$ is a fuzzy subgroup (of $G_1$) in $\tau_1$. Thus f is a g-fuzzy continuous map from $G_1$ to $G_2$. It is easy to verify that $f(xy) \neq f(x)(f(y))$ for $x = y = -1 \in G_1$. This proves that f is not a group homomorphism.

**THEOREM 1.12.3:** *Let $\tau_1$ and $\tau_2$ be any two g-fuzzy topologies on groups $(G_1, \bullet)$ and $(G_2, *)$ respectively. Then every group homomorphism $f: (G_1, \tau_1) \to (G_2, \tau_2)$ need not in general be a g-fuzzy continuous map.*

*Proof*: To prove this theorem it is sufficient if we prove the result to be false for a particular $\tau_1$ and $\tau_2$ defined on any group G as in our definition of g-fuzzy continuous map we have not assumed $G_1$ and $G_2$ to be distinct.

Let G be any group. Define two g-fuzzy topologies $\tau_1$ and $\tau_2$ on the group G as follows:

$\tau_1 = \{\phi_G, 1_G, \lambda\}$ and $\tau_2 = \{\phi_G, 1_G, \mu\}$ where $\lambda, \mu : G \to [0, 1]$ is as given below:

$$\lambda(x) = \begin{cases} 1 & \text{if } x = e \\ 0 & \text{if } x \neq e \end{cases}$$

and

$$\mu(x) = \begin{cases} 1 & \text{if } x \neq e \\ 0 & \text{if } x = e \end{cases}$$

Hence e is the identity element of the group G.

Define $f : (G, \tau_1) \to (G, \tau_2)$ by $f(x) = x^{-1}$ for every $x \in G$. It can be easily verified that f is a group homomorphism. For $x \in G$ and $\mu \in \tau_2$ we have

$$\begin{aligned}(f^{-1}(\mu))(x) &= \mu(f(x)) \\ &= \mu(x^{-1}) \\ &= \mu(x) \text{ (for } x = e, \text{if and only if } x^{-1} = e).\end{aligned}$$

This gives $(f^{-1}(\mu))(x) = \mu(x)$ for every $x \in G$. That is $f^{-1}(\mu) = \mu$. Thus $f^{-1}(\mu) \notin \tau_1$, as $\mu \notin \tau_1$. Hence f is not a g-fuzzy continuous map on G.

Now we proceed on to define g-fuzzy homeomorphism.

**DEFINITION 1.12.3:** *Let $(G_1, \tau_1)$ and $(G_2, \tau_2)$ be any two g-fuzzy topological spaces. A function $f: (G_1, \tau_1) \to (G_2, \tau_2)$ is said to be a g-fuzzy homeomorphism if it satisfies the following three conditions.*

  i. *f is one to one and onto.*
  ii. *f is a g-fuzzy continuous map from $G_1$ to $G_2$ and*
  iii. *$f^{-1}$ is a g-fuzzy continuous map from $G_2$ to $G_1$.*

We now give an example of a g-fuzzy homeomorphism.

The following example will illustrate the g-fuzzy homeomorphism.

***Example 1.12.5:*** Take $G_1$ and $G_2$ to be the set of all integers and even integers respectively (they are groups with respect to usual addition) with corresponding topologies $\tau_1 = \{\phi_{G_1}, 1_{G_1}, \lambda\}$ and $\tau_2 = \{\phi_{G_2}, 1_{G_2}, \mu\}$ where $\lambda: G_1 \to [0, 1]$ is defined as

$$\lambda(x) = \begin{cases} 1 & \text{if } x = 0 \\ 0 & \text{if } x = \pm 1, \pm 2, \pm 3, \cdots \end{cases}$$

and $\mu: G_2 \to [0, 1]$ is defined as

$$\mu(x) = \begin{cases} 1 & \text{if } x = 0 \\ 0 & \text{if } x \pm 2, \pm 4, \pm 6, \cdots \end{cases}$$

The function $f : (G_1, \tau_1) \to (G_2, \tau_2)$ given by $f(x) = 2x$ is clearly one to one and onto. It is easy to verify that $f^{-1}(\phi_{G_2}) = \phi_{G_1}$, $f^{-1}(1_{G_2}) = 1_{G_1}$ and $f^{-1}(\mu) = \lambda$.

Hence f is a g-fuzzy continuous map from $G_1$ and $G_2$. Further we have checked that $f(\phi_{G_1}) = \phi_{G_2}$, $f(1_{G_1}) = 1_{G_2}$ and $f(\lambda) = \mu$. This proves $f^{-1}$ is a g-fuzzy continuous map from $G_2$ to $G_1$. Thus f is a g-fuzzy homeomorphism from $G_1$ to $G_2$.

We know that a fuzzy subset of a set X is said to be a fuzzy point if and only if it takes the value 0 for all y ∈ X. expect on one element, say, x ∈ X.

If its value at x is t ($0 < t \leq 1$) then we denote this fuzzy point by $x_1$.

**DEFINITION 1.12.4:** *A fuzzy subset $\mu$ in a g-fuzzy topological space, $(G, \tau)$ is called a $Q_g$-neighborhood of the fuzzy point $x_t$ (for $x \in G$) if and only if there exists fuzzy subgroup $\lambda$ (of G) in $\tau$ such that $\lambda \subset \mu$ and $x_t \, q \, \lambda$ (where $x_t \, q \, \lambda$ means $t + \lambda(x) > 1$ and $x_t$ is quasi coincident with $\lambda$).*

This definition is illustrated by the following example:

***Example 1.12.6:*** Let $G = \{1, \omega, \omega^2\}$ be the group with respect to the usual multiplication, where $\omega$ denotes the cube root of unity.

Let $\tau = \{\phi_G, 1_G, \lambda, \mu\}$ where $\mu, \lambda: G \to [0, 1]$ are given by

$$\lambda(x) = \begin{cases} 0.9 & \text{if } x = 1 \\ 0.7 & \text{if } x = \omega, \omega^2 \end{cases}$$

and

$$\mu(x) = \begin{cases} 0.9 & \text{if } x = 1 \\ 0.8 & \text{if } x = \omega \\ 0.7 & \text{if } x = \omega^2 \end{cases}$$

It is easy to verify that $\tau$ is a g-fuzzy topology on the group G. We observe that $\mu$ is a $Q_g$-neighborhood of the fuzzy point $\omega_{0.5}$ as there exists a fuzzy subgroup $\lambda \in \tau$ such that $\lambda \subseteq \mu$ and the fuzzy point $\omega_{0.5}$ is quasi coincident with $\lambda$.

**DEFINITION 1.12.5:** *A g-fuzzy topological space $(G, \tau)$ is said to be a g-fuzzy Housdorff space if and only if for any two fuzzy points $x_t$ and $y_s$ ($x, y \in G$ and $x \neq y$) there exist $Q_g$-neighborhoods $\lambda$ and $\mu$ of $x_t$ and $y_s$ respectively such that $\lambda \cap \mu = \phi_G$.*

The following is a nice characterization theorem.

**THEOREM 1.12.4:** *Let $f : (G_1, \tau_1) \to (G_2, \tau_2)$ be a g-fuzzy homeomorphism. Then $(G_1, \tau_1)$ is a g-fuzzy Housdorff space if and only if $(G_2, \tau_2)$ is a g-fuzzy Housdorff space.*

*Proof*: Let $f : (G_1, \tau_1) \to (G_2, \tau_2)$ be a g-fuzzy homeomorphism. Suppose we assume that $(G_1, \tau_1)$ is a g-fuzzy Housdorff space we prove $(G_2, \tau_2)$ is a g-fuzzy Housdorff space. Let $x_t$ and $y_t$ be any two fuzzy points in $\tau_2$ with $x \neq y$ ($x, y \in G_2$) then $f^{-1}(x) \neq f^{-1}(y)$ as f is one to one.

Now consider for $z \in G_1$

$$\left(f^{-1}(x_t)\right)(z) = x_t(f(z))$$

$$= \begin{cases} t & \text{if } f(z) = x \\ 0 & \text{if } f(z) \neq x \end{cases}$$

$$= \begin{cases} t & \text{if } z = f^{-1}(x) \\ 0 & \text{if } z \neq f^{-1}(x) \end{cases}$$

$$= \left(f^{-1}(x)\right)_t (z).$$

That is $\left(f^{-1}(x_t)\right)(z) = \left(f^{-1}(x)\right)_t (z)$ for every $z \in G_1$.

From the above equality we have $f^{-1}(x_t) = \left(f^{-1}(x)\right)_t$.

Similarly we can prove that $f^{-1}(y_s) = \left(f^{-1}(y)\right)_s$, just by replacing the fuzzy point $x_t$ by $y_s$. Since $\left(f^{-1}(x)\right)_t$ and $\left(f^{-1}(y)\right)_s$, are fuzzy points in $\tau_1$ we have $f^{-1}(x_t)$ and $f^{-1}(y_s)$ are also fuzzy points in $\tau_1$ with $f^{-1}(x) \neq f^{-1}(y)$.

By the definition of a g-fuzzy Hausdorff space there exists $Q_g$-neighborhoods $\mu_x$ and $\mu_y$ of $f^{-1}(x_t)$ and $f^{-1}(y_s)$ respectively such that $\mu_x \cap \mu_y = \phi_{G_1}$. That is there exist fuzzy subgroups $\lambda_x, \lambda_y \in \tau_1$ such that

i. $\lambda_x \subset \mu_x$ and $f^{-1}(x_t)$ q $\lambda_x$.
ii. $\lambda_y \subset \mu_y$ and $f^{-1}(y_s)$ q $\lambda$ y and
iii. $\mu_x \cap \mu_y = \phi_{G_1}$.

Since f is a g-fuzzy continuous map from $G_1$ to $G_2$ and $f^{-1}$ is a g-fuzzy continuous map from $G_2$ to $G_1$ there exists $Q_g$ – neighborhoods $f(\mu_x)$ and $f(\mu_y)$ of $x_t$ and $y_s$ respectively such that $f(\mu_x) \cap f(\mu_y) = \phi_{G_2}$. Hence by the definition of g-fuzzy Hausdorff space $(G_2, \tau_2)$ is a g-fuzzy Hausdorff space.

Conversely let $(G_2, \tau_2)$ be a g-fuzzy Hausdorft space. By a similar argument and by also using the fact that both f and $f^{-1}$ are g-fuzzy continuous maps we can prove that $(G_1, \tau_1)$ is a g-fuzzy Hausdorff space.

Now we proceed on to define fuzzy polynomial rings and fuzzy polynomial semirings.

We define it in a very different way and in fact it is also quite new.

**DEFINITION 1.12.6:** *Let [0, 1] be a closed interval of the real line. Let R be a commutative ring with 1 or the field of reals. The fuzzy polynomial ring in the variable x with coefficients from R denoted by $R[x^{[0, 1]}]$ consists of elements of the form*

$$\left\{ a_o + a_1 x^{\gamma_1} + \cdots + a_n x^{\gamma_n} \,\middle|\, \begin{array}{l} a_o, a_1, \cdots, a_n \in R \text{ and} \\ \gamma_1, \gamma_2, \cdots, \gamma_n \in [0, 1] \text{ with } \gamma_1 < \gamma_2 < \cdots < \gamma_n \end{array} \right\}.$$

*In order to make a ring out of $R[x^{[0, 1]}]$ we must be able to recognize when two elements in it are equal, we must be able to add and multiply elements in $R[x^{[0, 1]}]$ so that axioms defining a ring hold true for $R[x^{[0, 1]}]$. If $p(x) = a_o + a_1 x^{\gamma_1} + \cdots + a_n x^{\gamma_n}$ and $q(x) = b_o + b_1 x^{s_1} + \cdots + b_m x^{s_m}$ are in $R[x^{[0, 1]}]$; $p(x) = q(x)$ if and only if (1) m = n, $\gamma_i = s_i$ and $a_i = b_i$, for each i.*

*If $p(x) = a_o + a_1 x^{\gamma_1} + \cdots + a_n x^{\gamma_n}$ and $q(x) = b_o + b_1 x^{s_1} + \cdots + b_m x^{s_m}$ then $p(x) + q(x) = c_o + c_1 x^{t_1} + \cdots + c_k x^{t_k}$ where $\gamma_1 < \gamma_2 < \ldots < \gamma_n$, $s_1 < s_2 < \ldots < s_m$, and $t_1 < t_2 < \ldots < t_K$ with $\gamma_i, s_j, t_p \in [0, 1]$, $1 \leq i \leq n$, $1 \leq j \leq m$ and $1 \leq p \leq k$. i.e. we add two polynomials as in case of usual polynomials by adding the like powers of x. Now it is easily verified $R[x^{[0, 1]}]$ is an abelians group under '+'; with $0 = 0 + 0 x^{\gamma_0} + \cdots + 0 x^{\gamma_n}$ as the additive identity which for short will be denoted by 0.*

*Now we have to define multiplication 'o' of two fuzzy polynomials in $R[x^{[0, 1]}]$. For $x^{s_1}$ and $x^{t_1} \in R[x^{[0, 1]}]$*

$$\begin{array}{rcl} x^{s_1} \,o\, x^{t_1} & = & x^{s_1+t_1} \text{ if } s_1 + t_1 \leq 1 \\ & = & x^{s_1+t_1-1} \text{ if } s_1 + t_1 > 1. \end{array}$$

*We extend this way of multiplication for any two polynomial $p(x), q(x) \in R[x^{[0, 1]}]$.*

*Clearly $R[x^{[0, 1]}, o]$ is a semigroup under 'o' and $\{R[x^{[0, 1]}], +, o\}$ is defined to be a fuzzy polynomial ring.*

Now unlike in a polynomial ring we see in case of fuzzy polynomial rings the degree of fuzzy polynomial is also fuzzy.

**DEFINITION 1.12.7:** *Let $p(x) = p_o + p_1 x^{\gamma_1} + \cdots + p_n x^{\gamma_n} \in R[x^{[0, 1]}]$ we call p(x); a fuzzy polynomial, the degree of p(x) is $\gamma_n$ provided $p_n \neq 0$. $\gamma_n$ always lie between 0 and 1.*

It is pertinent to note that if $p(x), q(x) \in R[x^{[0, 1]}]$ of degrees $\gamma_n$ and $s_n$ respectively then the degree of $p(x)$ o $q(x)$ can be $\gamma_n + s_n$ or $\gamma_n + s_n - 1$ depending on the fact whether $\gamma_n + s_n \leq 1$ or $\gamma_n + s_n > 1$ respectively. Thus deg $(p(x) o q(x))$ in general are not equal to deg $p(x)$ + deg $q(x)$ in fuzzy polynomials.

We have nice results analogous to usual polynomial rings.

**THEOREM 1.12.5:** *Let $R[x^{[0, 1]}]$ be the polynomial ring over the field of reals R. Then $R[x^{[0, 1]}]$ is an infinite dimensional vector space over R.*

*Proof*: Follows from the fact $R[x^{[0, 1]}]$ is an additive abelian group for any fuzzy polynomial ring and we have for all $r \in R$ and $p(x) \in R[x^{[0, 1]}]$, $rp(x) \in R[x^{[0, 1]}]$. $R[x^{[0, 1]}]$ is a vector space over R since the interval [0, 1] is of infinite cardinality we say $R[x^{[0, 1]}]$ is a vector space with infinite basis.

Can we have fuzzy polynomial rings which are finite dimensional or has a finite basis? The answer is yes.

*Example 1.12.7:* Let R be the ring of integers $R[x^{[0, 1]}]$ = {all polynomial formed by $x^0$ = 1 $x^{0.2}$, $x^{0.4}$, $x^{0.6}$, $x^{0.8}$ and x with coefficients from R} = {$\alpha_0 + \alpha_1 x^{0.2} + \alpha_2 x^{0.4} + \alpha_3 x^{0.6} + \alpha_4 x^{0.8} + \alpha_5 x \mid \alpha_1, \alpha_2, \alpha_3, \alpha_4, \alpha_5 \in R$}.

Clearly $R[x^{[0, 1]}]$ is a finite dimensional vector space over R.

Dimension of $R[x^{[0, 1]}]$ is 6. One of the basis is [1, $x^{0.2}$, $x^{0.4}$, $x^{0.6}$, $x^{0.8}$ and x]

**DEFINITION 1.12.8:** *Let $R[x^{[0, 1]}]$ be a fuzzy polynomial ring. A non empty subset P of $R[x^{[0, 1]}]$ is said to be a fuzzy polynomial subring of $R[x^{[0, 1]}]$ if P itself is a fuzzy polynomial ring under the operations inherited from $R[x^{[0, 1]}]$.*

*One more major difference between a polynomial ring and fuzzy polynomial ring is that $R[x^{[0, 1]}]$ is not an integral domain even if R is a field.*

**THEOREM 1.12.6:** *Let R be a field. $R[x^{[0, 1]}]$ be a fuzzy polynomial ring $R[x^{[0, 1]}]$ is not a field or an integral domain.*

*Proof*: To show $R[x^{[0, 1]}]$ is not a field or an integral domain we have to show $R[x^{[0, 1]}]$ contains two polynomials $0 \neq p(x)$, $0 \neq q(x) \in R[x^{[0, 1]}]$ with $p(x) \circ q(x) = 0$. Take $p(x) = x + x^{0.5}$ and $q(x) = 1 - x^{0.5}$.

$$
\begin{aligned}
p(x) \circ q(x) &= (x + x^{0.5})(1 - x^{0.5}) \\
&= x + x^{0.5} - x^{0.5} - x^1 \\
&= x + x^{0.5} - x^{0.5} - x \\
&= 0.
\end{aligned}
$$

Hence the claim.

**DEFINITION 1.12.9:** *Let $R[x^{[0, 1]}]$ be a fuzzy polynomial ring. A proper subset I of $R[x^{[0, 1]}]$ is said to be a fuzzy polynomial ideal of the fuzzy polynomial ring if*

    i. *I is a fuzzy polynomial subring of $R[x^{[0, 1]}]$*
    ii. *For every $p(x) \in R[x^{[0, 1]}]$ and $q(x) \in I$; $p(x) \circ q(x) \in I$.*

***Example 1.12.8:*** Let $Z_2 = [0, 1]$ be the prime field of characteristic two and $Z_2[x^{[0, 1]}]$ be the fuzzy polynomial ring. Can $Z_2[x^{[0, 1]}]$ have fuzzy polynomial ideal?

Several interesting results can be obtained in case of fuzzy polynomial rings.

Now we proceed on to define fuzzy polynomial semirings. In this section S will denote a commutative semiring with unit or a semifield.

**DEFINITION 1.12.10:** *Let S be a semiring. Let x be an indeterminate. We call $S[x^{[0, 1]}]$ to be a fuzzy polynomial semiring if*

$$S[x^{[0, 1]}] = \{s_o + s_1 x^{\gamma_1} + \cdots + s_n x^{\gamma_n} \mid s_o, s_1, \cdots, s_n \in S$$

*and $\gamma_1, \gamma_2, \cdots, \gamma_n \in [0,1]$ with $\gamma_1 < \gamma_2 < \cdots < \gamma_n \}$.*

*Let $p(x), q(x) \in S[x^{[0, 1]}]$, where*

$$p(x) = p_o + p_1 x^{\gamma_1} + \cdots + p_n x^{\gamma_n}$$
$$q(x) = q_o + q_1 x^{s_1} + \cdots + q_m x^{s_m}$$

*$p(x) = q(x)$ if and only if $p_i = q_i$, $\gamma_i = s_i$, $i = 1, ..., n$ with $n = m$.*

*Let $p(x), q(x) \in S[x^{[0, 1]}]$. $p(x) + q(x) = c_o + c_1 x^{t_1} + \cdots + c_k x^{t_k}$ where for each $c_i = p_i + q_i$ provided $\gamma_i = s_i$ i.e. we add two polynomials by adding terms of x which have equal powers.*

*Define multiplication of two polynomial $p(x), q(x) \in S[x^{[0, 1]}]$,*

$$p(x) \circ q(x) = \{p_o + p_1 x^{\gamma_1} + \cdots + p_n x^{\gamma_n}\} \circ \{q_o + q_1 x^{s_1} + \cdots + q_m x^{s_m}\}$$
$$= p_o q_o + p_1 q_0 x^{\gamma_1} + p_2 q_0 x^{\gamma_2} + \cdots + p_n q_0 x^{\gamma_n} +$$
$$p_o q_1 x^{s_1} + p_1 q_1 x^{\gamma_1 + s_1} + \cdots + p_n q_1 x^{\gamma_1 + s_m} + \cdots + p_n q_m x^{\gamma_n + s_m};$$

$$\gamma_i + s_j = \begin{cases} p & \text{if } \gamma_i + s_j \leq 1 \\ p-1 & \text{if } \gamma_i + s_j > 1. \end{cases}$$

***Example 1.12.9:*** Let $p(x), q(x)$ be fuzzy polynomials in a fuzzy polynomial semiring $Z[x^{[0, 1]}]$. $p(x) = x^{0.11} + 18x^{0.2} + x^{0.7} + 3x^{0.9}$, $q(x) = 3 + x^{0.21} + 8x^{0.71} + 11x$. We say degree $p(x)$ is 0.9 and degree of $q(x)$ is 1, degree of $p(x) \circ q(x) \neq \deg p(x) + \deg q(x)$.

We define fuzzy semifields in a different way not as fuzzy subsets.

**DEFINITION 1.12.11:** *Let [0, 1] be a unit interval. Define operation '+' and '•' on [0, 1] by for all $a, b \in [0, 1]$,*

$$a + b \;=\; \begin{cases} a+b & \text{if } a+b \le 1 \\ (a+b)-1 & \text{if } a+b > 1. \end{cases}$$

*Then ([0, 1], +) is a commutative semigroup with 0 as its identity. Define '•' on [0, 1] by a • b as the usual multiplication. Clearly ([0, 1], +, •) is a semifield. We call this as a fuzzy semifield; as in our opinion when entries are from the unit interval we can call them as fuzzy semifield as fuzziness allows us to do so.*

We can also define fuzzy semirings in a different and classical way as follows.

**DEFINITION 1.12.12:** *Let R be any ring and [0, 1] be the unit interval a map $p : R \to [0, 1]$ is said to be a fuzzy semiring if*

  i. $p(x + y) = \begin{cases} p(x) + p(y) & \text{if } p(x) + p(y) \le 1 \\ p(x) + p(y) - 1 & \text{if } p(x) + p(y) > 1. \end{cases}$

  ii. $p(x\,y) = p(x) \bullet p(y)$.
  iii. $p(0) = 0$.
  iv. $p(1) = 1$ if R is a ring with 1.

*We say p is a fuzzy semifield if in $p(x\,y) = p(x) \bullet p(y) = 0$ i.e. $p(x) \bullet p(y) = 0$ if and only if $p(x) = 0$ or $p(y) = 0$.*

*It is important to note that if R is taken as a semiring instead of a ring then any semiring homomorphism p from R to [0, 1] is a fuzzy subsemiring. If p is not a homomorphism then p is not a semiring.*

*We have for $p : R \to [0, 1]$ which is a fuzzy subsemiring then p is said to be a strict fuzzy semiring if $x + y = 0$ implies $x = 0$ and $y = 0$ then $p(x) + p(y) = 0$ implies $p(x) = 0$ and $p(y) = 0$.*

Once again we mention the fuzzy semiring in the usual sense has been carried out in the section on fuzzy semirings. Here we devote to see fuzzy semiring in a different perspective.

**DEFINITION 1.12.13:** *Let $V = [0, 1]$ be a fuzzy semifield. An additive abelian semigroup P with 0 is said to be a fuzzy semivector space over [0, 1] if for all $x, y \in P$ and $c \in [0, 1]$; $c\,x$ and $x\,c \in P$ i.e. $c[x + y] = c\,x + c\,y \in P$. In short $[0, 1]\,P \subseteq P$ and $P\,[0, 1] \subseteq P$.*

We define for the fuzzy semivector space defined in this manner the fuzzy semivector transformation.

**DEFINITION 1.12.14:** *Let V be a semivector space over a semifield. Let F and P be fuzzy semivector spaces over [0, 1]. A map $p : V \to P$ is called a fuzzy semivector transformation if for all $v \in V$, $p(v) \in V$. For every $c \in F$, $p(c) \in [0, 1]$ such that $p(cv + u) = p(c)p(v) + p(u)$ where $p(c) \in [0, 1]$; $p(u), p(v) \in P$.*

$$\text{Further } p(c+d) = \begin{cases} p(c) + p(d) & \text{if } p(c) + p(d) \leq 1 \\ p(c) + p(d) - 1 & \text{if } p(c) + p(d) > 1. \end{cases}$$

$$p(cd) = p(c) \bullet p(d)$$
$$p(0) = 0$$
$$p(1) = 1$$

for all $c, d \in F$.

**DEFINITION 1.12.15:** *Let P be a fuzzy semivector space over [0, 1]. The fuzzy dimension of P over [0, 1] is the minimum number of elements in P required to generate P.*

As in case of semivector spaces [134] several results in this direction can be developed and defined. But as in case of classical fuzzy semivector space we do not view fuzzy semivector spaces as a fuzzy subset. As once again we wish to state that our main motivation is the study of Smarandache fuzzy algebra we leave the fuzzy algebra development to the reader.

Now we proceed on to define a system called fuzzy opposite sets and systems. This system happens to be opposite of the system [0, 1] denoted [0, 1] opposite. Here 1 happens to be the least element and 0 happens to be largest element i.e. for example 0.99 is an element closer to largest element 1 in the classical fuzzy sets and systems but in fuzzy opposite sets and systems 0.99 is an element closer to the least element.

**DEFINITION 1.12.16:** *Let S be a set. A fuzzy subset A of S is a function $A : S \rightarrow [0, 1]$. The fuzzy opposite subset or opposite fuzzy subset for the same function. A denoted by $A^{opp}$ is a map from S to $[0, 1]^{opp}$.*

*In most of the algebraic developments the opposite fuzzy subset which is a fuzzy subgroup or a fuzzy subring or a fuzzy ideal or a fuzzy module or a fuzzy vector space happens to give way to fuzzy opposite structures i.e. if A is a fuzzy subgroup of S with a modification $A^{opp}$ happens to be opposite fuzzy subgroup.*

**DEFINITION 1.12.17:** *Let G be a group. A fuzzy subset A of G is a fuzzy subgroup of G. If $A : G \rightarrow [0, 1]$ such*

    i.   $A(xy) \geq \min\{A(x), A(y)\}$ *for all $x, y \in G$.*
    ii.  $A(x) = A(x^{-1})$ *for all $x \in G$.*

*Now $A^{opp}$ is a fuzzy opposite subgroup if $A^{opp} : G \rightarrow [0, 1]$ such that*

    i.   $A^{opp}(xy) \leq \max\{A(x), A(y)\}$
    ii.  $A^{opp}(x) = A^{opp}(x^{-1})$ *for all $x, y \in G$.*

Thus the notions of fuzzy normality and other notions can be developed using $A^{opp}$ given A to be fuzzy algebraic structure.

**DEFINITION 1.12.18:** *Let A be a fuzzy subset of S i.e. $A : S \to [0, 1]$. For $t \in [0, 1]$ the set $A_t = \{x \in S / A(x) \geq t\}$ is called the level subset of the fuzzy subset A. Now $A^{opp}$ is a fuzzy subset of the set S. For $t \in [0, 1]^{opp}$, the set $A_t^{opp} = \{x \in S / A(x) \leq t\}$ is called the level subset of the fuzzy opposite subset A.*

Most results which hold good for fuzzy level subset can be carried out for level subset of fuzzy opposite subset A.

# PART TWO

# Smarandache Fuzzy Algebra

# Chapter Two

# SMARANDACHE FUZZY SEMIGROUPS AND ITS PROPERTIES

This chapter has five sections. In section one we give the definition of Smarandache fuzzy semigroup, describe several Smarandache definitions and give results. In section two we study the Smarandache fuzzy semigroups. Here we have given 28 Smarandache fuzzy definitions of a S-semigroup and we have developed these properties in about 67 theorems. Section three analyzes about elementwise S-fuzzy properties in S-semigroups. About 30 new definitions in this direction are given and the properties enunciated in 59 theorems. The fourth section is devoted to the study of Smarandache fuzzy bisemigroups and bigroups. Thirty-two new definitions about S-fuzzy semigroups are introduced and its properties are analyzed in 23 theorems. The chief attraction is the fifth section which proposes 54 problems for the reader to solve.

## 2.1. Definition of Smarandache fuzzy semigroups with examples:

In this section we for the first time introduce the concept of Smarandache fuzzy semigroups. Smarandache semigroups are thoroughly studied in the year [123]. Fuzzy semigroups were introduced in the late 1970s. To get the concept of Smarandache fuzzy semigroups we go for the concept of fuzzy groups; as Smarandache groups do not exist in literature we get almost all properties of groups to be present in Smarandache semigroups. Just for the sake of reader we recall the definition of Smarandache semigroup then proceed on to define Smarandache fuzzy semigroup and give examples.

**DEFINITION 2.1.1:** *Let S be a semigroup. S is said to be a Smarandache semigroup (S-semigroup) if S has a proper subset P such that P is a group under the operations of S.*

*Example 2.1.1:* Let $S(3)$ be a permutation semigroup. $S(3)$ is a S-semigroup as $S_3 \subset S(3)$ and $S_3$ is the permutation group of degree 3.

*Example 2.1.2:* Let $M_{n \times n} = \{(a_{ij}) \mid a_{ij} \in Q\}$ be the set of all $n \times n$ matrices with entries from Q. $M_{n \times n}$ is a semigroup under matrix multiplication. $P_{n \times n} = \{(a_{ij}) = A \mid |A| \neq 0\}$ is a proper subset of $M_{n \times n}$ which is a group. So $M_{n \times n}$ is a S-semigroup.

The semigroup given in example 2.1.1 is a non-commutative S-semigroup of finite order where as in example 2.1.2 the semigroup is an infinite non-commutative S-semigroup.

*Example 2.1.3:* Let $Z_{20} = \{0, 1, 2, \ldots, 19\}$ be a semigroup under multiplication. Clearly, $Z_{20}$ is a S-semigroup for $P = \{1, 19\}$ is a subgroup in $Z_{20}$. This is a commutative S-semigroup of finite order.

*Example 2.1.4:* Let Q[x] be the ring of polynomials; Q[x] is a semigroup under multiplication. P = Q \ {0} ⊂ Q[x] is a group under multiplication; so Q[x] is an infinite commutative S-semigroup.

As our main motivation in this paper is the study of Smarandache fuzzy semigroup we stop at this stage and request the reader to refer [123] for complete information regarding S-semigroup.

Fuzzy semigroups are dealt in chapter one section ten. Now we proceed on to define Smarandache fuzzy semigroup.

**DEFINITION 2.1.2:** *Let S be a S-semigroup. A fuzzy subset A of S is said to be a Smarandache fuzzy semigroup (S-fuzzy semigroup) if $A : S \to [0, 1]$ is such that A restricted to atleast one subset P of S which is a subgroup is a fuzzy group. That is for all $x, y \in P \subset S$. $A(x, y) \geq \min \{A(x), A(y)\}$ and $A(x) = A(x^{-1})$ for all x in P.*

*This S-fuzzy semigroup is denoted by $A_P$ i.e. $A_P : P \to [0, 1]$ is a fuzzy group.*

*Thus in case of S-fuzzy semigroup we face several situations which are enlisted below. Here S denotes a S-semigroup. Suppose $P_1, P_2, \ldots, P_n$ are proper subsets of S which are groups under the operations of S.*

*A fuzzy subset $A : S \to [0, 1]$ may be a S-fuzzy semigroup or sometimes may not be S-fuzzy semigroup. So we call all fuzzy subsets from S to [0, 1] which are not S-fuzzy semigroups where S is a S-semigroup as Smarandache non-fuzzy semigroup (S-non-fuzzy semigroup).*

*If $A : S \to [0, 1]$ such that there exist atleast one subset $P_i \subset S$ which is a group and $A_{P_i}$ is a fuzzy subgroup then we call the fuzzy subset P a S-fuzzy semigroup.*

*If the fuzzy subset $A : S \to [0, 1]$ in such that $A_{P_i} : P_i \to [0, 1]$, for $i = 1, 2, \ldots, n$ are fuzzy subgroups then we call the fuzzy subset $A : S \to [0, 1]$ a Smarandache strong fuzzy semigroup (S-strong fuzzy semigroup).*

*Thus we see in general a fuzzy subset $A : S \to [0, 1]$ where S is a S-semigroup may be a S-fuzzy semigroup or a S-strong fuzzy semigroup or S-non-fuzzy semigroup.*

*Example 2.1.5:* Let S(3) be a S-semigroup. The map $\mu : S(3) \to [0, 1]$ defined by

$$\mu(x) = \begin{cases} 0.5 & \text{if } x = \begin{pmatrix} 1 & 2 & 3 \\ 1 & 2 & 3 \end{pmatrix} \\ 0.4 & \text{if } x = \begin{pmatrix} 1 & 2 & 3 \\ 2 & 3 & 1 \end{pmatrix} \text{ and } \begin{pmatrix} 1 & 2 & 3 \\ 3 & 1 & 2 \end{pmatrix} \\ 0.7 & \text{otherwise.} \end{cases}$$

Clearly μ(x) is a S-fuzzy semigroup. It is easily verified that μ restricted to the subset

$$P = \left\{ \begin{pmatrix} 1 & 2 & 3 \\ 1 & 2 & 3 \end{pmatrix}, \begin{pmatrix} 1 & 2 & 3 \\ 2 & 3 & 1 \end{pmatrix}, \begin{pmatrix} 1 & 2 & 3 \\ 3 & 1 & 2 \end{pmatrix} \right\}$$

which is a subgroup in S(3) is a fuzzy subgroup. i.e. $\mu_P : P \to [0, 1]$ is a fuzzy subgroup.

Thus $\mu : S(3) \to [0, 1]$ is S-fuzzy semigroup. Further, it is worthwhile to note $\mu : S(3) \to [0, 1]$ is not a fuzzy subgroup. For if

$$x = \begin{pmatrix} 1 & 2 & 3 \\ 1 & 3 & 2 \end{pmatrix} \text{ and } y = \begin{pmatrix} 1 & 2 & 3 \\ 2 & 1 & 3 \end{pmatrix}$$

$$\mu(xy) = \mu \begin{pmatrix} 1 & 2 & 3 \\ 2 & 3 & 1 \end{pmatrix} \not\geq \min \{\mu(x), \mu(x)\}$$

$$= \min \{0.7, 0.7\} = 0.7.$$

Thus μ is only a S-fuzzy semigroup and not a S-strong fuzzy semigroup.

*Example 2.1.6:* Let S(3) be a S-semigroup. $\lambda : S(3) \to [0, 1]$ be defined by

$$\lambda(x) = \begin{cases} 0.7 & \text{if } x = \begin{pmatrix} 1 & 2 & 3 \\ 1 & 2 & 3 \end{pmatrix} \\ 0.5 & \text{if } x = \begin{pmatrix} 1 & 2 & 3 \\ 1 & 3 & 2 \end{pmatrix} \\ 0.8 & \text{otherwise} \end{cases}$$

λ is a S-fuzzy semigroup and is not a S-strong fuzzy semigroup.

We call these S-fuzzy semigroups as level I S-fuzzy semigroups. Now we proceed on to define Smarandache fuzzy semigroup of level II.

**DEFINITION 2.1.3:** *Let S be a S-semigroup A fuzzy semigroup $A : S \to [0, 1]$ is called the Smarandache fuzzy semigroup of level II (S-fuzzy semigroup of level II) if A restricted to one of the subsets P of S where P is a group is a fuzzy group that is $A_P : P \to [0, 1]$ is a fuzzy group.*

**THEOREM 2.1.1:** *All S-fuzzy semigroup of level II are S-fuzzy semigroup of level I.*

*Proof*: Direct by the very definition.

It is important to note that the converse is not true.

**THEOREM 2.1.2:** *Let S be a S-semigroup, $\mu : S \to [0, 1]$ be a S-fuzzy semigroup I. Then $\mu$ in general need not be a S-fuzzy semigroup II.*

*Proof*: Follows by a counter example. Take S(3) which is a S-semigroup. Define $\mu : S(3) \to [0, 1]$ by

$$\mu(x) = \begin{cases} 0.5 & \text{if } x = \begin{pmatrix} 1 & 2 & 3 \\ 1 & 2 & 3 \end{pmatrix} \\ 0.3 & \text{if } x = \begin{pmatrix} 1 & 2 & 3 \\ 2 & 1 & 3 \end{pmatrix} \\ 0.8 & \text{otherwise} \end{cases}$$

Clearly $\mu$ is a S-fuzzy semigroup I and is not a S-fuzzy semigroup II.

Now we reformulate the definition of S-fuzzy semigroup II in Smarandache language.

Define S-strong fuzzy semigroup II.

**DEFINITION 2.1.4:** *Let S be a S-semigroup. A fuzzy semigroup $A : S \to [0, 1]$ is called the Smarandache strong fuzzy semigroup of level II (S-strong fuzzy semigroup II) if A restricted to every proper subsets $P_i$ which are subgroups in S i.e. $A_{P_i} : S \to [0, 1]$ is a fuzzy group for every proper subgroup $P_i$ in S.*

**THEOREM 2.1.3:** *Let S be a S-semigroup. Every S-strong fuzzy semigroup II is a S-fuzzy semigroup and not conversely.*

*Proof:* Straightforward.

**DEFINITION 2.1.5:** *Let S be a S-semigroup $\mu : S \to [0, 1]$ be a fuzzy semigroup. $\mu$ is said to be a Smarandache fuzzy semigroup III (S-fuzzy semigroup III) if $\mu$ has a proper subset $\sigma$ such that $\sigma \leq \mu$ and $\sigma$ is a fuzzy subgroup.*

The study of inter relation between these three types of S-fuzzy semigroup is left as an exercise for the reader.

*Example 2.1.7:* Let

$$M_{2\times 2} = \left\{ \begin{pmatrix} a & b \\ c & d \end{pmatrix} \mid a, b, c, d \in Z_2 = \{0,1\} \right\}$$

is a S-semigroup. Define $\mu : M_{2 \times 2} \to [0, 1]$ by

$$\mu(A) = \begin{cases} .5 & \text{if } A = \begin{pmatrix} 1 & 0 \\ 0 & 1 \end{pmatrix} \\ .4 & \text{if } |A| \neq 0 \\ .9 & \text{if } |A| = 0 \end{cases}$$

Clearly μ is a S-fuzzy semigroup.

## 2.2. Substructures of S-fuzzy semigroups and their properties

In this section we introduce the notion of S-fuzzy ideals in a S-fuzzy semigroup and also define and introduce several properties enjoyed by fuzzy groups as all S-fuzzy semigroups contain S-fuzzy subgroups.

**DEFINITION 2.2.1:** *Let G be a S-semigroup. $\mu : G \to [0, 1]$ is said to be a Smarandache fuzzy ideal (S-fuzzy ideal) of the S-fuzzy semigroup μ if $\mu(x) = \mu(yxy^{-1})$ for all $x, y \in A$ where $A \subset G$ and A is a subgroup of G.*

**DEFINITION 2.2.2:** *Let G be a S-semigroup and $\mu : G \to [0, 1]$ be a S-fuzzy semigroup of G. For $t \leq \mu(e)$, $t \in [0, 1]$ the set $\mu_t = \{x \in A \,/\, \mu(x) \geq t\}$ is called the Smarandache level semigroup (S-level semigroup) of the S-fuzzy semigroup μ where $A \subset G$ and A is the subgroup of the S-semigroup G.*

The following results are important.

**THEOREM 2.2.1:** *Let G be a S-semigroup with $A \subset G$ is a proper subset of G which is a subgroup of G $\mu : G \to [0, 1]$ is a S-fuzzy semigroup if and only if $\mu(xy^{-1}) \geq \min(\mu(x), \mu(y))$ for all $x, y \in A$.*

*Proof*: The result is straightforward by the very definitions of fuzzy subgroup, S-semigroup and S-fuzzy semigroup.

**THEOREM 2.2.2:** *Let G be a S-semigroup. For $t \leq \mu(e)$, $t \in [0, 1]$ for a given t and given G and fixed μ we can have several S-level semigroups depending on the number of subgroups A in G.*

**DEFINITION 2.2.3:** *Let G be a S-semigroup. If A be a proper subset of G which is a subsemigroup of G and A contains the largest subgroup of G; then we call the fuzzy subset $\mu : G \to [0, 1]$ to be a Smarandache fuzzy hyper subsemigroup (S-fuzzy hypersub semigroup) if μ restricted to P, i.e. $\mu : P \to [0, 1]$, $P \subset A$ and P is the largest subgroup of A is a fuzzy subgroup of G.*

*Example 2.2.1:* Let S (3) be a S-semigroup. Clearly

$$A = \left\{ \begin{pmatrix} 1 & 2 & 3 \\ 1 & 1 & 1 \end{pmatrix}, \begin{pmatrix} 1 & 2 & 3 \\ 2 & 2 & 2 \end{pmatrix}, \begin{pmatrix} 1 & 2 & 3 \\ 3 & 3 & 3 \end{pmatrix} \right\} \cup S_3$$

is a semigroup. μ : S (3) → [0, 1]. The restricted map μ : $S_3$ → [0, 1] which is a fuzzy group. Then μ is a S-fuzzy hyper subsemigroup.

**THEOREM 2.2.3:** *Every S-fuzzy subsemigroup of a S-semigroup G in general need not be a S-fuzzy hyper subsemigroup of G.*

*Proof*: Follows from the very definition.

**COROLLARY 2.2.1:** *Every S-fuzzy hyper subsemigroup of a S-semigroup G is a S-fuzzy subsemigroup of G.*

*Proof:* Direct hence left as an exercise for the reader.

**DEFINITION 2.2.4:** *Let G be a S-semigroup. We say G is a Smarandache fuzzy simple semigroup (S-fuzzy simple semigroup) if G has no S-fuzzy hyper subsemigroup.*

Now we proceed on to define Smarandache fuzzy symmetric semigroup.

**DEFINITION 2.2.5:** *Let S(n) denote the symmetric semigroup i.e. the semigroup of all mapping of a set of n elements {1, 2, …, n} to itself S(n) under composition of mappings is a semigroup.*

*Let SF(S(n)) denote the set of all Smarandache fuzzy subsemigroups of S(n). If μ ∈ SF(S(n)) then Im μ = { f(x) /x ∈ A ⊂ S(n)} where A is a proper subset of S(n) which is a subgroup of S(n) under the operations of S(n) Let μ, σ ∈ SF (S(n)) we say if |Im(μ)| < |Im(σ)| then we write μ < σ. By this rule we can define max SF (S(n)). Let f be a S-fuzzy subsemigroup of S(n). If f = max SF(S(n)) then we say that f is a Smarandache fuzzy symmetric semigroup of S(n) (S-fuzzy symmetric semigroup of S(n)).*

Now we introduce a new concept called Smarandache co-fuzzy symmetric semigroup. To do this we just recall the definition of co-fuzzy symmetric group.

**DEFINITION [139]:** *Let G ($S_n$) = {g /g is a fuzzy subgroup of $S_n$ and g(C(Π)) is a constant for every Π ∈ $S_n$} where C (Π) is the conjugacy class of $S_n$ containing Π, which denotes the set of all y ∈ $S_n$ such that y = x Π $x^{-1}$ for x ∈ $S_n$. If g = max G($S_n$) then we call g as co-fuzzy symmetric subgroup of $S_n$.*

**DEFINITION 2.2.6:** *Let SG (S(n)) = {g /g is a S-fuzzy subsemigroup of S(n) and g(C(Π)) is constant for every Π ∈ $S_n$ } where C (Π) is the conjugacy class of $S_n$ ⊂ S(n) containing Π which denotes the set of all y ∈ $S_n$ such that y = x Π $x^{-1}$ (for x ∈ $S_n$). If g = max SG(S(n)) then we call g as Smarandache co-fuzzy symmetric subsemigroup of S(n) (S-co-fuzzy symmetric subsemigroup of S(n)).*

**THEOREM 2.2.4:** *Every S-co-fuzzy symmetric subsemigroup of a symmetric semigroup S(n) is a S-fuzzy symmetric subsemigroup of the symmetric semigroup S(n).*

*Proof*: Follows from the definitions of S-fuzzy symmetric subsemigroup and S-co fuzzy symmetric subsemigroup.

Now we proceed on to define yet another new concept called Smarandache fuzzy normal subsemigroup of a S-semigroup.

**DEFINITION 2.2.7:** *A S-fuzzy subsemigroup $\mu$ of a S-semigroup G is said to be a Smarandache fuzzy normal subsemigroup (S-fuzzy normal subsemigroup) of the S-semigroup G if $\mu(xy) = \mu(yx)$ for all $x, y \in A$, A any proper subset of G which is a subgroup and $\mu$ restricted to A i.e. $\mu : A \to [0, 1]$ is a fuzzy subgroup of A.*

Now we define Smarandache fuzzy cosets of a S-semigroup.

**DEFINITION 2.2.8:** *Let $\mu$ be a S-fuzzy subsemigroup of S-semigroup G. For any $a \in A \subset G$ (where A is the subgroup associated with this $\mu$), a $\mu$ defined by $(a\mu)(x) = \mu(a^{-1}x)$ for every $x \in A$ is called the Smarandache fuzzy coset (S-fuzzy coset) of the S-semigroup (S-fuzzy coset of the S-semigroup) G determined by a and $\mu$.*

**DEFINITION 2.2.9:** *Let $\mu$ be a S-fuzzy subsemigroup of a S-semigroup G and H = $\{x \in A \subset G \mid \mu(x) = \mu(0)\}$ (where $\mu$ is the associated S-fuzzy subsemigroup with A) then $o(\mu)$, order of $\mu$ is defined as $o(\mu) = o(H)$.*

**DEFINITION 2.2.10:** *Let $\lambda$ and $\mu$ be two S-fuzzy subsemigroup of a S-semigroup G. Then $\lambda$ and $\mu$ are said to be Smarandache conjugate fuzzy subsemigroups (S-conjugate fuzzy subsemigroups) of G (relative to the same A) if for some $g \in A \subset G$ (A proper subset of G which is a subgroup) we have $\lambda(x) = \mu(g^{-1}xg)$ for every $x \in A$.*

The main result on S-conjugate fuzzy subsemigroups of a S-semigroup G is as follows.

**THEOREM 2.2.5:** *If $\lambda$ and $\mu$ are S-conjugate fuzzy subsemigroups relative to a subgroup $A \subset G$; where G is a S-semigroup, then $o(\lambda) = o(\mu)$.*

*Proof*: Let $\lambda$ and $\mu$ be S-conjugate fuzzy subsemigroups of the S-semigroup G. By the very definition of S-conjugate fuzzy subsemigroups of the S-semigroup there exists $g \in A \subset G$ such that $\lambda(x) = \mu(g^{-1}xg)$ for every $x \in A \subset G$ (A a proper subgroup of G). Now let us define $H = \{x \in A \mid \lambda(x) = \lambda(e)\}$ and $K = \{x \in A \mid \mu(x) = \mu(e)\}$ where e is the identity element of $A \subset G$.

Clearly H is a subgroup of the subgroup A, for H is a t-level subset of the subgroup A where $t = \lambda(e)$. Similarly K is also a subgroup of the subgroup $A \subset G$.

To prove that $o(\lambda) = o(\mu)$ it is sufficient to prove $o(H) = o(K)$, by the definition of the order of the fuzzy subgroup of the subgroup A.

To prove that $o(H) = o(K)$ as the first step we prove $H \subset gKg^{-1}$ for some $g \in A \subset G$. In the second step we prove that for the same $g \in A$; $K \subset g^{-1}Hg$.

Let x be any element in H. Since λ and μ are S-conjugate fuzzy subsemigroups of the subgroup A ⊂ G (equivalently of the S-semigroup G) we have for some g ∈ G.

$$\begin{aligned}
\mu(g^{-1} x g) &= \lambda(x) \\
&= \lambda(e) \text{ (since } x \in H\text{)} \\
&= \mu(g^{-1} e g) \text{ (since } \lambda \text{ and } \mu \text{ are S-conjugate fuzzy subsemigroups of G)} \\
&= \mu(e).
\end{aligned}$$

Hence there exists g ∈ A such that $\mu(g^{-1} x g) = \mu(e)$. Now by the definition of K we have $g^{-1} x g \in K$ implies $x \in gKg^{-1}$. Therefore $H \subseteq gKg^{-1}$.

To prove the other inclusion: Take x an arbitrary element in K. Since λ and μ are S-conjugate fuzzy subsemigroups of the S-semigroup G for the same g ∈ A ⊆ G (A a proper subset which is a subgroup of the S-semigroup and λ and μ are conjugate fuzzy subgroup relative to this A ⊂ G; i.e. S-conjugate fuzzy subsemigroups of G) used in the earlier result to prove $H \subseteq g K g^{-1}$ we have

$$\begin{aligned}
\lambda(gxg^{-1}) &= \mu(x) \\
&= \mu(e) \text{ (since } x \in K\text{)} \\
&= \lambda(geg^{-1}) \text{ (since } \lambda \text{ and } \mu \text{ are conjugate fuzzy subgroups of A} \subset G\text{)} \\
&= \lambda(e).
\end{aligned}$$

Hence $gxg^{-1} \in H$ for the same g ∈ A ⊆ G, that is $x \in g^{-1} H g$. Hence $K \subseteq g^{-1} H g$.

From the first and the second steps we have $H \subseteq gKg^{-1}$ and $K \subseteq g^{-1}Hg$, i.e. $Hg \subseteq gK$ and $gK \subset Hg$ Thus $Hg = gK$ so $H = gKg^{-1}$. Since K is a subgroup of the group A contained in the S-semigroup G we have $o(xKx^{-1}) = o(K)$ for every x ∈ A. Now choose x = g. Then we have $o(gKg^{-1}) = o(K)$. Therefore $o(H) = o(K)$. Hence $o(\lambda) = o(\mu)$. Hence the theorem.

Recall G a S-semigroup. A and B are S-fuzzy subsemigroups of G related to the subgroup P in G (P ⊂ G, P a proper subset of G which is a group) such that B ⊆ A. Let $x_t \subseteq A$. Then the fuzzy subset $x_t$ o B (B o $x_t$) is called the Smarandache fuzzy left (right) coset (S-fuzzy left (right) coset) of B in A with representative $x_t$ where by the operation 'o' we mean the following for all x ∈ P ($x_t$ o B) (x) = sup {inf {$x_t$ (y), B(z)} such that x = yz}.

We give the following results which can be easily extended to S-fuzzy left (right) cosets in an analogous way.

**THEOREM 2.2.6:** *Let G be a S-semigroup. A be a S-fuzzy semigroup of G relative to a subgroup P ⊂ G. B be a S-fuzzy subsemigroup relative to P, such that B ⊂ A. Then for all z ∈ P, ($x_t$ o B) (x) = inf {t, B ($x^{-1}$z)} and (B o $x_t$) (z) = inf {t, B (z $x^{-1}$)} for $x_t$ ∈ A.*

*Proof*: As in case of fuzzy subgroups the proof is straightforward.

**THEOREM 2.2.7:** *Let A, B, and G as in the Theorem 2.2.6 Let $x_t$, $y_s$ ∈ A. Then*

i.   $(x_t \circ B) = (y_s \circ B)$ if and only if inf $\{t, (B(e)\} = $ inf $\{s, B\ (y^{-1}x)\}$ and inf $\{s, B(e)\} = $ inf $\{t, B\ (x^{-1}y)\}$.
ii.  $B \circ x_t = B \circ y_s$ if and only if inf $\{t, B(e)\ \} = $ inf $\{s, B(xy^{-1})\}$ and inf $\{s, B(e)\} = $ inf $\{t, B\ (yx^{-1})\}$.

*Proof*: Straightforward; hence left for the reader as an exercise.

**THEOREM 2.2.8:** *Let A, B, G be as in Theorem 2.2.6 Let $x_t, y_t \subseteq A$. If $B(y^{-1}x) = B(e)$ then $x_t \circ B = y_t \circ B$.*

*Proof*: Since $B(x^{-1}y) = B(y^{-1}x) = B(e)$; inf$\{t, B(e)\} = $ inf$\{t, B(x^{-1}y)\} = $ inf$\{t, B(y^{-1}x)\}$. Hence by earlier results $x_t \circ B = y_t \circ B$.

**THEOREM 2.2.9:** *Let A, B, G be as in Theorem 2.2.6 Let $x_t, y_t \subseteq A$. Then the following conditions are equivalent.*

i.   $x_t \circ B = y_t \circ B$.
ii.  $(y^{-1}x)_t \circ B = e_t \circ B$.
iii. $(x^{-1}y)_t \circ B = e_t \circ B$.

*Proof*: By earlier results we have $x_t \circ B = y_t \circ B$ iff inf$\{t, B(e)\} = $ inf$\{t, B(y^{-1}x)\}$ and inf$\{t, B(e)\} = $ inf $\{t, B(x^{-1}y)\}$. The latter conditions are equivalent to (ii) and (iii).

**THEOREM 2.2.10:** *Let A, B, G and P be as in Theorem 2.2.6. Let $x, y \in P$ and $s, t \in [0, A(e)]$. Suppose that $B(e) = A(e)$; then*

i.   $x_t \circ B = y_s \circ B$ if and only if $t = $ inf $\{s, B\ (y^{-1}x)\}$, $s = $ inf $\{\ t, B\ (x^{-1}y)\}$.
ii.  $x_t \circ B = y_t \circ B$ if and only if $(y^{-1}x)_t \subseteq B$.
iii. $x_t \circ B = y_s \circ B$ if and only if $t = s \leq B\ (x^{-1}y)$
iv.  $x_t \circ B = x_s \circ B$ if and only if $t = s$.

*Proof*:

i.   By earlier results $x_t \circ B = y_s \circ B$ if and only if $t = $ inf $\{s, B(y^{-1}x)\}$ and $s = $ inf$\{t, B(x^{-1}y)\}$.
ii.  By (i), $x_t \circ B = y_t \circ B$ if and only if and only if $t = $ inf $\{t, B\ (y^{-1}x)\}$, $t = $ inf$\{t, B(x^{-1}y)\}$ if and only if $B(y^{-1}x) \geq t$, $B(x^{-1}y) \geq t$.
iii. By (i) and the fact that $B((y^{-1}x) = B(x^{-1}y)$; $x_t \circ B = y_s \circ B$ if and only if $t = s < B(x^{-1}y)$.
iv.  Follows immediately from (iii).

The proof of the following theorem is true hence only stated and left for the reader to prove.

**THEOREM 2.2.11:** *Let $s, t \in [0, A(e)]$. Suppose that $B(e) = A(e)$. If $t \neq s$, then $\{x_t \circ B\ /\ x_t \subseteq A\} \cap \{y_s \circ B\ /y_s \subseteq A\} = \phi$.*

We define Smarandache fuzzy normal subsemigroup of a S-semigroup G.

**DEFINITION 2.2.11:** *Let G be a S-semigroup. P a proper subset of G which is a group. A and B S-fuzzy subsemigroup of G relative to the subgroup $P \subset G$, $B \subset A$; B is said to be Smarandache fuzzy normal (S-fuzzy normal) in A if and only if for all $x_t \subseteq A$, $x_t \circ B = B \circ x_t$.*

**THEOREM 2.2.12:** *Let A, B and G be as in Definition 2.2.11 Let $x_t, y_s \subseteq A$. If B is fuzzy normal in A then $(x_t \circ B) \circ (y_s \circ B) = (xy)_r \circ B$ where $r = \inf \{t, s\}$.*

*Proof*: 'o' operation is associative and $B \circ B = B$. The proof follows from [73, 74].

**THEOREM 2.2.13:** *Let A, B, P and G be as in Theorem 2.2.12. Let $A / B = \{x_1 \circ B \mid x_1 \subset A$ and $x \in P\}$. Suppose that B is S-fuzzy normal in A. Then (A / B, o) is a semigroup with identity. If $B(e) = A(e)$ then A / B is completely regular, that is A / B is union of (disjoint) groups.*

*Proof*: If $x_t \circ B$, $y_s \circ B \in A / B$, then clearly $(xy)_r \circ B \in A / B$ where $r = \inf \{t, s\}$. Clearly, $e_{A(e)}$ is the identity of A / B. 'o' is an associative operation. For fixed $t \in [0, A(e)]$, let $(A / B)^{(t)} = \{x_t \circ B \mid x_t \subset A, x \in P\}$. Then $(A / B)^{(t)}$ is closed under o, $e_t \circ B$ is the identity of $(A / B)^{(t)}$ and $(x^{-1})_t \circ B$ is the inverse of $x_t \circ B$. Hence $(A / B)^{(t)}$ is a group.

Clearly $A / B = \bigcup_{t \in [0, A(e)]} (A / B)^{(t)}$.

*Example 2.2.2:* Let $G = \{e, a, b, c, d, f, g\}$ be a semigroup given by the table

| o | e | a | b | c | d | f | g |
|---|---|---|---|---|---|---|---|
| e | e | a | b | c | d | f | g |
| a | a | e | c | b | a | a | a |
| b | b | c | e | a | b | b | b |
| c | c | b | a | e | c | c | c |
| d | d | a | a | a | d | f | g |
| f | f | b | b | b | d | f | g |
| g | g | c | c | c | d | f | g |

G is S-semigroup. Define fuzzy subsets A and B of $P \subset G$ where $P = \{a, b, e, c\}$ the klein four group by $A(e) = A(a) = 1$, $A(b) = A(c) = ¾$ and $B(e) = B(a) = 1$, $B(b) = B(c) = ½$ with $A(x) = B(x) = 0$ for all other $x \in G \setminus P$. Then A and B are S-fuzzy subsemigroups of G such that $B \subset A$ and B is S-fuzzy normal in A. Now $e_1 \circ B$ is the identity of A/B but $e_{3/4} \circ B$ does not have an inverse. Hence A/B is not a group.

**THEOREM 2.2.14:** *Let G be a S-semigroup. A and B be S-fuzzy semigroup of G related to the subgroup P in G. (i.e $P \subset G$ is a proper subset of G which is a group under the operations of G) such that $B \subset A$. B is S-fuzzy normal in A if and only if for all $t \in [0, B(e)]$ $B_t$ is normal in $A_t$.*

*Proof*: As in case of fuzzy subgroups.

**THEOREM 2.2.15:** *Let A, B, P and G be as in Theorem 2.2.14. Suppose $0 \leq t \leq B(e)$, $x_s \subseteq A$ and $t \leq s$. Then $(x_s \circ B)_t = xB_t$ and $(B \circ x_s)_t = B_t x$.*

*Proof*: $y \in (x_s \circ B)_t$ if and only if $(x_s \circ B)(y) > t$ if and only if $\inf \{s, B(x^{-1}y)\} \geq t$ if and only if $B(x^{-1}y) \geq t$ if and only if $x^{-1} y \in B_t$ if and only if $y \in xB_t$.

**THEOREM 2.2.16:** *Let $t \in [0, B(e)]$. Suppose that B is fuzzy normal in A. (A, B, G and P as in Theorem 2.2.14); then $A_t / B_t \cong (A / B)^{(t)}$.*

*Proof*: Refer [86] and the definition to prove this theorem.

**DEFINITION 2.2.12:** *Let G be a S-commutative semigroup; i.e. every proper subset P of G which is a group is commutative. A, a S-fuzzy subsemigroup of G related to P; P a subgroup of G. B is S-fuzzy normal in A. We say A is Smarandache bounded (S-bounded) over B if there exists $n \in N$ such that for all $x_t \subset A$, $(x_t)^n \subseteq B$. A is bounded over B if and only if $A_t / B_t$ is uniformly bounded for all $t \in [0, A(e)]$. Hence if A is bounded over B then $A_t / B_t$ is a direct product of cyclic groups for all $t \in [0, A(e)]$.*

**DEFINITION 2.2.13:** *Let G be a S-commutative semigroup. A be a S-fuzzy semigroup of G related to a subgroup $P \subset G$, such that $C \subset A$ and $A = B \otimes C$ the Smarandache fuzzy direct product (S-fuzzy direct product) of B and C i.e. $A = B \circ C$ for all $x \in P$, $(B \cap C)(x) = 0$. Then $A_t = B_t \otimes C_t$ for all $t \in [0, A(e)]$ Thus $A / B \cong \bigcup_{t \in (0, A(e)]} C_t \cup \{e_0 \circ B\}$ we give some conditions for B to be a Smarandache fuzzy direct factor (S-fuzzy direct factor) of A.*

*Let SF (A) denote the set of all S-fuzzy subsemigroups C of G (G a S-semigroup) such that $C \subseteq A$ and $C(e) = A(e)$. Let $C^* = \{x \in P / C(x) > 0\}$. Then $C^*$ is a subgroup of P.*

*We say that B is Smarandache compatible (S-compatible) in A if and only if $A(e) = B(e)$ and for all $s, t \in (0, A(e)]$, $s \leq t$, $A_s = A_t B_s$ and $A_t \cap B_s = B_t$. It is shown that if B is Smarandache divisible (S-divisible) i.e. for all $x_t \subseteq B$ with $t > 0$ and for all $n \in N$ there exists $y_t \subset B$ such that $(y_t)^n = x_t$, then it need not be the case that B is a S-fuzzy direct factor of A. If B is S-compatible in A and B is divisible, then B is a S-fuzzy direct factor of A.*

**THEOREM 2.2.17:** *The following conditions are equivalent.*

  i. $A(e) = B(e)$ *and there exist subgroup H of G (G a S-semigroup) such that for all $t \in (0, A(e)]$, $A_t = B_t \otimes H$.*
  ii. *There exists $C \in SF(A)$ such that $A = B \otimes C$ and $C^* = C_*$ where $C_* = C_{C(0)}$.*
  iii. *B is S-compatible in A and there exists $C \in SF(A)$ such that $A_* = B_* \otimes C_*$.*

*Proof*: Left for the reader to prove.

**THEOREM 2.2.18:** *Let A and B be S-fuzzy subsemigroups of the S-semigroup G such that $B \subseteq A$, we say that B is S-pure in A if and only if for all $x_t \subset B$ with $t > o$, for all $n \in N$ and for all $y_t \subset A$, $(y_t)^n = x_t$ implies that there exists $b_t \subseteq B$ such that $(b_t)^n = x_t$.*

**THEOREM 2.2.19:** *Suppose that B is S-compatible and S-pure in A.*

   i. *If $A_n/B_n$ is a direct product of S-cyclic groups, then B is a S-fuzzy direct factor of A.*

   ii.   *If B is S-bounded then B is a S-fuzzy direct factor of A.*

*Proof:*

   i. Since B is S-pure in A, $B_*$ is S-pure in $A_*$. Hence there exists a subgroup H if $A_*$ such that $A_* = B_* \otimes H$. Since B is S-compatible in A, $A_t = B_t \otimes H$ for all $t \in (0, A(e)]$, the desired results follows from earlier results.

   ii. $B_*$ is S-pure in $A_*$ and $B_*$ is S-bounded. Hence there exists subgroup H of $A_*$ such that $A_* = B_* \otimes H$. The result of the proof is direct.

**THEOREM 2.2.20:** *Suppose that B is S-compatible in A. If B is divisible, then B is a S-direct factor of A.*

*Proof***:** Left as an exercise for the reader to prove.

**THEOREM 2.2.21:** *Let G be a S-semigroup. H and K be subgroups of G such that $H \not\subset K$ and $K \not\subset H$, then there exists a S-fuzzy subsemigroup $\mu$ of $P \times P$ such that $\sigma_\mu$ is not a S-fuzzy subsemigroup of G. ($P \subset G$ is a proper subset which is a group).*

*Proof***:** Let $H_o = \{(e, e)\}$, $H_1 = H \times G$ and $H_2 = P \times P$. Let $1 \geq t_0 \geq t \geq t_2 \geq 0$. Let $\mu$ be a fuzzy subset of $P \times P$ such that

$$\mu(H_o) = t_o$$
$$\mu(H_1 \setminus H_0) = t_1 \text{ and}$$
$$\mu(H_2 \setminus H_1) = t_2.$$

Then $\mu$ is a S-fuzzy subsemigroup of $P \times P$, since $\mu_t$ is a subgroup of $P \times P$ for all $t \in \text{Im}(\mu)$. Now there exists $h \in H$, $k \in K$ such that $k \notin H$ and $h \notin K$ then $hk \notin H$ and $hk \notin K$. Hence $(hk, x), (x, hk) \notin H_1$ for all $x \in G$.

Then $\sigma_\mu(hk) = t_2$. Now since $(h, k) \in H_1 \setminus H_0$, $\sigma_\mu(h) = t_1 = \sigma_\mu(k)$. Thus $\sigma_\mu(hk) = t_2 < t_1 = \min \{\sigma_\mu(h), \sigma_\mu(k)\}$. Hence $\sigma_\mu$ is not a S-fuzzy subsemigroup of G.

**THEOREM 2.2.22:** *Let $\sigma_\mu$ be a S-fuzzy semigroup of the S-semigroup G for all S-fuzzy subgroups $\mu$ of $P \times P \subset G \times G$ (P a subgroup of G). Then either G is S a S- cyclic semigroup of order $p^m$ for some $m \geq 1$ or $P = Z(p^\infty)$.*

*Proof***:** Follows as in case of fuzzy subgroup.

Interested reader may refer [86].

**THEOREM 2.2.23:** *Let $\mu$ be a S-fuzzy semigroup of $P \times P \subset G \times G$. Then $\sigma_\mu (e) = \mu(e,e) \geq \sigma_\mu (x)$ for all $x \in P$ and $\sigma_\mu(x) = \sigma_\mu (x^{-1})$ for all $x \in P$.*

*Proof*: Left as an exercise for the reader.

**THEOREM 2.2.24:** *Let G be a S-cyclic semigroup having cyclic subgroups of order $p^m$ for some $m \geq 1$, where p is a prime. If $\mu$ is a S-fuzzy semigroup of $P \times P \subset G \times G$ then $\sigma_\mu$ is a fuzzy subgroup of P.*

*Proof*: Follows from [86].

**THEOREM 2.2.25:** *Let G be a S-semigroup having a subgroup $Z(p^\infty) = P$. If $\mu$ is a S-fuzzy semigroup of $P \times P \subset G \times G$, then $\sigma_\mu$ is a fuzzy subgroup of P.*

*Proof*: Refer [86].

**DEFINITION 2.2.14:** *Let G be a S-semigroup under addition. A be a S-fuzzy semigroup of G relative to a subgroup of G relative to a subgroup P, $P \subset G$. We assume $A(0) > 0$. If $x \in P$ then $A(0) \geq A(x)$. For $t \in [0, A(0)]$, let $A_t = \{x \in P \mid A(x) \geq t\}$, then $A_t$ is a subgroup of P. For $x \in P$ and $t \in [0, 1]$ we let $x_t$ denote the fuzzy subset of P defined by $x_t(y) = t$ if $y = x$ and $x_t(y) = 0$ if $y \neq x$.*

*Then $x_t$ is called a Smarandache fuzzy singleton (S-fuzzy singleton). If $x_t$ and $y_s$ are S-fuzzy singletons, then $x_t + y_s$ is defined to be the fuzzy subset $(x + y)_r$ where $r = \min\{t,s\}$.*

**DEFINITION 2.2.15:** *Let G be a S-semigroup. A a S-fuzzy semigroup of G relative to the subgroup P, $P \subset G$. A is called a Smarandache torsion fuzzy subsemigroup (S-torsion fuzzy subsemigroup) of G if and only if for all S-fuzzy singletons $x_t \subset A$ with $t > 0$ there exists $n \in N$ (N – set of natural integers) such that $n(x_t) = 0$.*

**DEFINITION 2.2.16:** *Let G be a S-semigroup. A be a S-fuzzy subsemigroup of G. A is called a Smarandache p-primary fuzzy subsemigroup (S-p-primary fuzzy subsemigroup) of G if and only if there exists a prime p such that for all S-fuzzy singletons $x_t \subset A$ with $t > 0$ there exists $n \in N$ (N-natural integers) such that $p^n(x_t) = 0_t$.*

**THEOREM 2.2.26:** *A is S-p-primary if and only if $A^*$ is S-p-primary.*

*Proof*: Left for the reader to prove.

**THEOREM 2.2.27:** *A is a S-p-primary if and only if $A_t$ is S-p-primary for all $t \in (0, A(0)]$ If A is S-p-primary then $A_*$ is S-p-primary.*

*Proof*: Easily follows by the very definitions and a matter of routine.

**THEOREM 2.2.28:** *Let p be a prime Define the fuzzy subset $A^{(p)}$ of G, G a S-semigroup by for all $x \in P \subset G$ (P a subgroup of G), $A^{(p)}(x) = A(x)$ if $x \in (A^*)_p$ and $A^{(p)}(x) = 0$ otherwise. Then $A^{(p)}$ is a S-p-primary fuzzy subsemigroup of $P \subset G$ and $(A^{(p)})^* = (A^*)_p$.*

*Furthermore $(A^{(p)})_t = (A_t)_p$ for all $t \in (0, A(0)]$.*

*Proof:* As in case of groups refer [83].

**THEOREM 2.2.29:** *Let p be a prime. Then $A^{(p)}$ is the unique S-maximal p-primary fuzzy subsemigroup of P ($P \subset G$) contained in A a S-fuzzy subsemigroup related to P.*

*Proof***:** Refer [83].

For p a prime $A^{(p)}$ is the Smarandache p-primary component of A. Now we proceed on to define Smarandache divisible fuzzy semigroup of the S-semigroup.

**DEFINITION 2.2.17:** *Let G be a S-semigroup. A be a S-fuzzy subsemigroup of G. A is called Smarandache divisible fuzzy subsemigroup (S-divisible fuzzy subsemigroup) of G if and only if for all S-fuzzy singletons $x_t \subset A$ with $t > 0$ and for all $n \in N$ there exists a S-fuzzy singleton $y_t \subset A$ such that $n(y_t) = x_t$.*

**THEOREM 2.2.30:** *A is S-divisible if and only if $A_t$ is S-divisible for all $t \in (0, A(0)]$.*

*Proof***:** Follows by the very definition.

**THEOREM 2.2.31:** *If A is S-divisible then $A^*$ and $A_*$ are S-divisible.*

*Proof***:** Direct by the very definition.

**THEOREM 2.2.32:** *If $A^*$ is S-divisible and A is constant on $A^* \setminus \{0\}$ then A is S-divisible.*

*Proof***:** Let $x_t \subseteq A$. with $t > 0$ and $n \in N$. Then $x \in A^*$ and so there exists $y \in A^*$ such that $x = ny$. If $y = 0$ then $x = 0$ the result is true. Let $y \neq 0$. Since A is constant on $A^* \setminus \{0\}$, $A(y) = A(x) \geq t$. Then $y_t \subset A$ and $x_t = n(y_t)$ Hence A is S-divisible.

**THEOREM 2.2.33:** *Let G be a S-semigroup. P a subgroup of G. For all $x, y \in P \subset G$ and $n \in N$, $ny = x$ implies that $A(x) = A(y)$ for all S-divisible fuzzy subsemigroups of $P \subset G$ if and only if P is torsion free.*

*Proof:* Follows as in case of groups by taking the S-semigroup G to contain a proper subset P such that P is a torsion free group and $A : G \to [0, 1]$ is such that A restricted to P is a divisible fuzzy subgroup A of P.

Recall for G a S-semigroup A a S-fuzzy subsemigroup of G related to P.   SF (A) = {B | B is a S-fuzzy subsemigroup P such that $B \subseteq A$} = {B | B is a fuzzy subgroup of P such that $B \subset A$}. Both are equal as S-fuzzy subsemigroup of G say A acts as a fuzzy subgroup on P, $P \subset G$.

**DEFINITION 2.2.18:** *Let $B \in SF(A)$ Then B is said to Smarandache pure (S-pure) in A if and only if for all S-fuzzy singletons $x_t \subseteq B$ with $t > 0$, for all $n \in N$, for all $y_t \subseteq A$, $n(y_t) = x_t$ implies that there exists $b_t \subset B$ such that $n(b_t) = x_t$.*

Throughout this portion G is an additive abelian group and A fuzzy subgroup of G.

**THEOREM 2.2.34:** *Let $B \in SF(A)$. Then B is S-pure in A if and only if $B_t$ is S-pure in $A_t$ for all $t \in (0, B(0)]$.*

*Proof*: Follows from the very definition.

**THEOREM 2.2.35:** *Let $B \in SF(A)$. If B is S-pure in A, then $B_*$ is S-pure in $A_*$.*

*Proof*: Direct by the very definition.

**DEFINITION 2.2.19:** *Let G be a S-semigroup. P a proper subset of G which is a subgroup of G. Let X be a fuzzy subset of P and let $n \in N$. Define the Smarandache fuzzy subset (S-fuzzy subset) nX of $P \subset G$ by for all $x \in P$ $(nX)(x) = 0$ if $x \notin nP$ and $(nX)(x) = \sup\{X(y) \mid y \in P, x = ny\}$ if $x \in nP$, We say X has Smarandache sup property (S-sup property) if and only if Im (X) has a maximal element.*

**THEOREM 2.2.36:** *Let $n \in N$ then*

  i. $nA(0) = A(0)$.
  ii. $nA \subseteq A$.
  iii. $nA$ is a S-fuzzy subsemigroup of P; $P \subset G$.
  iv. If A has the sup property, then Im $(n, A) \subseteq$ Im $(A)$ where A is a S-fuzzy subsemigroup of G relative to the group P; $P \subset G$.

*Proof*: Follows as in case of fuzzy subgroups.

**THEOREM 2.2.37:** *Let B and C be S-fuzzy subsemigroups of the S-semigroup G related to the subgroup P. Then*

  i. $(B \cap C)^* = B^* \cap C^*$.
  ii. $(B \cap C)_t = B_t \cap C_t$ for all $t \in (0, \min\{B(0), C(0)\}]$.

*Proof*: Now B and C be S-fuzzy subsemigroup of the S-semigroup G related to the subgroup $P \subset G$. $x \in (B \cap C)^*$ if and only if $(B \cap C)(x) > 0$ if and only if inf $\{B(x), C(x)\} > 0$ if and only if $B(x) > 0$ and $C(x) > 0$ if and only if $x \in B^* \cap C^*$.

ii. $x \in (B \cap C)_t$ if and only if $(B \cap C)(x) \geq t$ if and only if inf $\{B(x), C(x)\} \geq t$ if and only if $B(x) \geq t$ and $C(x) \geq t$ if and only if $x \in B_t \cap C_t$.

**DEFINITION 2.2.20:** *Let A be a S-fuzzy subsemigroup of G, G a S-semigroup. If $A^*$ is Smarandache torsion free (S-torsion free) and for $x \in nA^*$ ($n \in N$) we have $(n A)(x) = A(\omega)$ for some unique $\omega \in A^*$ such that $x = n\omega$*

**THEOREM 2.2.38:** *Let $n \in N$. Then*

i. $(nA)^* = nA^*$.
ii. $nA_t \subseteq (nA)_t$, *for all $t \in (0, A(0)]$ (A is a S-fuzzy subsemigroup of the S-semigroup G related to the subgroup P, $P \subset G$).*
iii. *Let $B \in SF(A)$ and $B(0) = A(0)$.*

*If A has the sup property and B is pure in A, then $nB_t = (nB)_t$ for all $t \in (0, A(0)]$.*

*Proof*: Exactly as in case of fuzzy group as A is a S-fuzzy semigroup related to $P \subseteq G$.

**THEOREM 2.2.39:** *Let $B, C \in SF(A)$ such that $B(0) = C(0) = A(0)$ Then we have the following results to be true.*

i. *Suppose that A has the sup property. If B is S-pure in A, then for all $n \in N$, $nB = B \cap nA$.*
ii. *Suppose that A and B have sup property. If for all $n \in N$, $nB = B \cap nA$, then B is S-pure in A.*
iii. *Suppose that $A^*$ is S-torsion free, then B is S-pure in A if and only if for all $n \in N$, $nB = B \cap nA$.*
iv. *For $C \subseteq B$. If C is S-pure in B and B is S-pure in A, then C is pure in A.*
v. *If B is S-divisible then B is S-pure in A.*
vi. *Suppose that A is S-divisible. Then B is S-pure in A if and only if B is S-divisible.*
vii. *Let $B \in SF(A)$ be such that $B(0) = A(0)$. If for all $x_t \subseteq A$ such that $x_t \not\subset B$, there does not exist $n \in N$ such that $n(x_t) \subseteq B$, then B is S-pure in A.*

*Proof*: The proof is for the reader to prove using the definitions.

**THEOREM 2.2.40:** *Suppose that $A^*$ is S-torsion free. Let $\{B_\alpha / \alpha \in \Omega\}$ be a collection of S-fuzzy subsemigroups of the S-semigroup G such that $B_\alpha \subset A$, $B_\alpha$ is S-pure in A and $B_\alpha(0) = A(0)$ for all $\alpha \in \Omega$. Then $\bigcap_{\alpha \in \Omega} B_\alpha$ is S-pure in A.*

*Proof*: Matter of routine using the fact $B_\alpha$ is a S-fuzzy subsemigroup of S-semigroup G and $B_\alpha$ is S-pure in A.

**THEOREM 2.2.41:** *Suppose that $A^*$ is S-torsion free. Let $B \in SF(A)$ be such that $B(0) = A(0)$, then B is contained in a unique smallest S-pure fuzzy subsemigroup in A.*

*Proof*: Straightforward.

**THEOREM 2.2.42:** *Let $\{B_\alpha / \alpha \in \Omega\}$ and $\{C_\alpha / \alpha \in \Omega\}$ be chains of S-fuzzy subsemigroups of S-semigroup. Let $n \in N$ then*

$$\text{i.} \quad n\left(\bigcup_{\alpha \in \Omega} B_\alpha\right) = \bigcup_{\alpha \in \Omega} nB_\alpha$$

ii. $\bigcup_{\alpha \in \Omega}(B_\alpha \cap C_\alpha) = \left(\bigcup_{\alpha \in \Omega} B_\alpha\right) \cap \left(\bigcup_{\alpha \in \Omega} C_\alpha\right)$

*Proof*: Using the fact each $B_\alpha$ and $C_\alpha$ are S-fuzzy subsemigroups hence also fuzzy subgroups of a S-semigroup the proof is analogous to fuzzy subgroup of a group.

**THEOREM 2.2.43:** *Let $\{A_\alpha / \alpha \in \Omega\}$ be a chain of S-fuzzy subsemigroup of the S-semigroup G such that $A_\alpha \subseteq A$, $A_\alpha$ is S-pure in A and $A_\alpha(0) = A(0)$ for all $\alpha \in \Omega$. Suppose that either A and $\bigcup_{\alpha \in \Omega} A_\alpha$ have the sup property of $A^*$ is S-torsion free, then $\bigcup_{\alpha \in \Omega} A_\alpha$ is a S-pure fuzzy subsemigroup of G in A.*

*Proof*: For $n \in N$; we have

$$n\left(\bigcup_{\alpha \in \Omega} A_\alpha\right) = \bigcup_{\alpha \in \Omega} nA_\alpha = \bigcup_{\alpha \in \Omega}(A_\alpha \cap nA) = \left(\bigcup_{\alpha \in \Omega} A_\alpha\right) \cap nA$$

by earlier results we get the theorem.

Next we define a new notion called Smarandache fuzzy weak direct sum of a S-semigroup G.

**DEFINITION 2.2.21:** *Let $\{A_\alpha / \alpha \in \Omega\}$ be a collection of S-fuzzy subsemigroups of a S-semigroup G. Then A is said to be Smarandache fuzzy weak direct sum (S-fuzzy weak direct sum) of the $A_\alpha$ if and only if $A = \Sigma A_\alpha$ and $\forall x \in P \subset G$ (P a subgroup of G)*

$$x \neq 0, \left(A_\beta \cap \sum_{\alpha \in \Omega \setminus \{\beta\}} A_\alpha\right)(x) = 0.$$

*If A is the S-fuzzy weak direct sum of the $A_\alpha$ then we write $A = \bigoplus_{\alpha \in \Omega} A_\alpha$.*

**THEOREM 2.2.44:** *Let $\{A_\alpha / \alpha \in \Omega\}$ be a collection of S-fuzzy subsemigroups of the S-semigroup G. Such that $A_\alpha \subset A$ and $A_\alpha(0) = A(0)$ for all $\alpha \in \Omega$. Then $A = \bigoplus_{\alpha \in \Omega} A_\alpha$ if and only if $A^* = \bigoplus_{\alpha \in \Omega} A_\alpha^*$ and $A = \sum_{\alpha \in \Omega} A_\alpha$.*

*Proof*: As in case of fuzzy subgroups.

It is pertinent to mention here that if $\{A_\alpha | \alpha \in \Omega\}$ is a collection of S-fuzzy subsemigroups of G. A necessary condition for $A = \bigoplus_{\alpha \in \Omega} A_\alpha$ is that $A = A_\alpha$ or $A_\alpha^*$ for all $\alpha \in \Omega$.

**THEOREM 2.2.45:** *Let $\{A_\alpha \,/\, \alpha \in \Omega\}$ be a collection of S-fuzzy subsemigroups of a S-semigroup G such that $A = A_\alpha$ on $A_\alpha^*$ for all $\alpha \in \Omega$ and $A^* = \sum A_\alpha^*$. If for all $\alpha, \beta \in \Omega$ such that $A_\alpha \neq A_\beta$; $Im(A_\alpha) \cap Im(A_\beta) \subseteq \{0,1\}$ then $A = \sum_{\alpha \in \Omega} A_\alpha$.*

*Proof*: Directly as in case of fuzzy subgroups of a group G.

**THEOREM 2.2.46:** *Let $\{A_\alpha \,/\, \alpha \in \Omega\}$ be a collection of S-fuzzy subsemigroups of a S-semigroup G such that $A_\alpha \subseteq A$ for all $\alpha \in \Omega$. Suppose that either*

i. $\bigcup_{\alpha \in \Omega} Im(A_\alpha)$ is finite or

ii. $A^* = \bigoplus_{\alpha \in \Omega} A_\alpha^*$

*Then $A = \sum_{\alpha \in \Omega} A_\alpha$ if and only if $A_\alpha(0) = A(0)$ and $A_t = \sum_{\alpha \in \Omega}(A_\alpha)_t$ for all $t \in [0, A(0)]$.*

*Proof*: As in case of fuzzy subgroups.

The following theorems for S-fuzzy subsemigroups of a S-semigroups G can be easily proved on similar lines using the techniques of fuzzy subgroups.

**THEOREM 2.2.47:** *Let $\{A_\alpha \,/\, \alpha \in \Omega\}$ be a collection of S-fuzzy subsemigroups of the S-semigroup G such that $A_\alpha \subseteq A$ for all $\alpha \in \Omega$. Then $A = \bigoplus_{\alpha \in \Omega} A_\alpha$ if and only if $A_\alpha(0) = A(0)$ for all $\alpha \in \Omega$ and $A_t = \bigoplus_{\alpha \in \Omega}(A_\alpha)_t$ for all $t \in (0, A(0)]$.*

**THEOREM 2.2.48:** *Let $\{A_\alpha \,/\, \alpha \in \Omega\}$ be a collection of S-fuzzy subsemigroups of the S-semigroup G such that $A_\alpha \subseteq A$ for all $\alpha \in \Omega$. Suppose that either $\bigcup_{\alpha \in \Omega} Im(A_\alpha)$ is finite or $A^* = \bigoplus_{\alpha \in \Omega} A_\alpha^*$ holds. If $A = \sum_{\alpha \in \Omega} A_\alpha$ then $Im(A) \subseteq \bigcup_{\alpha \in \Omega} Im(A_\alpha)$.*

**THEOREM 2.2.49:** *Let $\{A_\alpha \,/\, \alpha \in \Omega\}$ be a collection of S-fuzzy subsemigroup of the semigroup G such that $A_\alpha \subseteq A$ for all $\alpha \in \Omega$. Suppose that $|\Omega| \geq 2$. If $A = \bigoplus_{\alpha \in \Omega} A_\alpha$, then $Im(A) \cup \{0\} = \bigcup_{\alpha \in \Omega} Im(A_\alpha)$.*

*Proof*: Using the above theorem we have $Im(A) \subseteq \bigcup_{\alpha \in \Omega} Im(A_\alpha)$ since $A = \bigoplus_{\alpha \in \Omega} A_\alpha$, $A = A_\alpha$ on $A_\alpha^*$.

Thus $\bigcup_{\alpha \in \Omega} Im(A_\alpha) \subseteq Im(A) \cup \{0\}$.

**THEOREM 2.2.50:** *Let $B, C, D \in SF(A)$ be such that $D \subseteq B$ and $D(0) = A(0)$. Suppose that $A = B \oplus C$. Then*

i.         *for all n ∈ N, n A = n B ⊕ nC.*
ii.         *D + C = D ⊕ C and B ∩ (D ⊕ C) = D.*

*Proof*: Matter of routine once we take B, C, D to be S-fuzzy subsemigroups we get the proof as in case of fuzzy subgroups.

**THEOREM 2.2.51**: *Let B and C ∈ SF(A). If A = B ⊕ C then B and C are S-pure in A.*

*Proof*: Suppose that $ny_t = x_t$ where $n \in N$, $y_t \subseteq A$ and $x_t \in B$. Since $A = B \oplus C$, $A_t = B_t \oplus C_t$. Thus $B_t \cap n A_t = nB_t$. Now $ny_t = (ny)_t$. Hence $ny = x$, $y \in A_t$ and $x \in B_t$. Thus there exists $b \in B_t$ such that $nb = x$. Hence $nb_t = x_t$. Thus B is a S-fuzzy subsemigroup in A is S-pure in A.

Thus from these results we chiefly observe that those fuzzy results true in case of just groups can be very easily studied in case of S-semigroups by using S-fuzzy subsemigroup. Thus it is important to note that except for Smarandache structure such study would be certainly impossible. Now we give some interesting results on S-fuzzy semigroups.

**THEOREM 2.2.52**: *Let G be a finite S-semigroup. Suppose that there exists a S-fuzzy subsemigroup μ of G (for some proper subset P of G with P a subgroup of G) satisfying the following conditions for x, y ∈ P ⊆ G.*

    i.     *μ(x) = μ(y) implies ⟨x⟩ = ⟨y⟩.*
    ii.     *μ(x) > μ(y) implies ⟨x⟩ ⊂ ⟨y⟩.*

*Then G is a S-cyclic semigroup.*

*Proof*: For $x \in P \subset G$, P a subgroup of a S-semigroup G. ⟨x⟩ denotes the subgroup generated by x.

$\mu : G \to [0, 1]$ such that μ restricted to P denoted by $\mu_P$ is a fuzzy subgroup of P, μ is a S-fuzzy subsemigroup of G i.e. $\mu_P : P \to [0, 1]$.

If $\mu_P$ is constant on P then $\mu_P(x) = \mu_P(y)$ for all $x, y \in P$ and so ⟨x⟩ = ⟨y⟩ by (i) consequently P = ⟨x⟩. Let us assume that $\mu_P$ is not constant on P. Let Im $\mu_P$ = $\{t_0, \ldots, t_n\}$ with $t_0 > t_1 > \ldots > t_n$. Then the chain of level subgroups of $\mu_P$ in given by

$$(\mu_P)_{t_o} \subset (\mu_P)_{t_1} \subset \cdots \subset (\mu_P)_{t_n} = P$$

Let $x \in P \setminus (\mu_P)_{t_{n-1}}$. We assert that P = ⟨x⟩. Consider any element g in $P \sim (\mu_P)_{t_{n-1}}$; Then $\mu_P(g) = \mu_P(x) = t_{n-1}$ so that ⟨g⟩ = ⟨x⟩. Hence $P \setminus (\mu_P)_{t_{n-1}} \subseteq$ ⟨x⟩. Now at the next stage let $g \in (\mu_P)_{t_{n-1}}$. Then $\mu_P(g) \geq t_{n-1} > t_n = \mu_P(x)$. By condition (ii) we have ⟨g⟩ ⊆ ⟨x⟩. Thus $(\mu_P)_{t_{n-1}} \subset$ ⟨x⟩. This yields P = ⟨x⟩.

**THEOREM 2.2.53:** *Let G be a S-cyclic semigroup having cyclic group P ($P \subset G$) of order $p^\alpha$ where p is a prime. Then there exists a S-fuzzy subsemigroup $\mu_P$ of G satisfying $\mu_P(x) = \mu_P(y)$ implies $\langle x \rangle = \langle y \rangle$ and $\mu_P(x) > \mu_P(y) \Rightarrow \langle x \rangle \subset \langle y \rangle$ where $\mu_P : P \to [0, 1]$ and $x, y \in P \subset G$.*

*Proof*: Consider the chain of subgroups of $P \subset G$; $(e) = P_o \subset P_1 \subset \ldots \subset P_n = P$ where $P_i$ is the subgroup of P generated by an element of order $p^i$, $i = 0, 1, 2, \ldots, n$ and e is the identity element of P. Define a fuzzy subset $\mu_P$ of P by $\mu_P(e) = t_o$ and $\mu_P(x) = t_i$ if $x \in P_i \setminus P_{i-1}$ for all $i = 1, 2, \ldots, n$; where $t_o, t_1, \ldots, t_n$ are numbers lying in the interval [0, 1] such that $t_o > t_1 > \ldots > t_n$. It is a routine matter to confirm that $\mu_P$ is a S-fuzzy subsemigroup of G satisfying the desired conditions.

The result which follows is a consequence of the above two theorems.

**THEOREM 2.2.54:** *Let G be a S-semigroup having a subgroup P of order $p^\alpha$. Then G is S-cyclic if and only if there exists a S-fuzzy subsemigroup $\mu_P$ of G such that for all $x, y \in P$, $\mu(x) = \mu(y) \Rightarrow \langle x \rangle = \langle y \rangle$ and $\mu(x) > \mu(x) \Rightarrow \langle x \rangle \subset \langle y \rangle$.*

From now onwards till we explicitly specify L denotes a completely distributive lattice with minimum 0 and maximum 1. X a set, the elements of the direct power $L^X$ are called L-subsets of fuzzy subsets of X. $f_0$ and $f_1$ denotes the minimum and maximum respectively of $L^X$.

**DEFINITION 2.2.22:** *Let G be a S-semigroup. $P \subset G$ be a subgroup in G. The Smarandache normal fuzzy subgroups (S-normal fuzzy subgroups) of the S-semigroup G i.e. the L-subsets of P whose cuts are normal subgroups of P.*

**THEOREM 2.2.55:** *An S-L subsemigroup f of $P \subset G$ is normal if and only if $f(x^{-1} z x) \geq f(z)$ for all $x, z \in P$. Several interesting results in this direction can be developed for L-subsets.*

*Recall for $\mu$ a fuzzy set in a set S. Then the level subset $\mu_t$ and strong level subset $\mu_t^>$ of $\mu$ are defined by*

  i. $\mu_t = \{x \in S / \mu(x) \geq t\}$ for $t \in [0, 1]$ ([0, 1] just the unit interval).
  ii. $\mu_t^> = \{x \in S \mid \mu(x) > t\}$ for $t \in [0, 1]$.

*We know for a group G $\mu$ be a fuzzy set of G. Then $\mu$ is a fuzzy subgroup of G if the following conditions hold*

  i. $\mu(xy) \geq \min \{\mu(x), \mu(y)\}$ for all $x, y \in G$.
  ii. $\mu(x^{-1}) = \mu(x)$ for $x \in G$.

*If $\mu$ is a fuzzy subgroup then it attains its supremum at 'e' the identity of G that is*

$$\sup_{x \in G} \mu(x) = \mu(e).$$

*We agree to call μ(e) to be the tip of the fuzzy subgroup μ. On the other hand a fuzzy subgroup may or may not attain its infimum. We shall write in short inf μ for $\inf_{x \in G} \mu(x)$ and refer to it as the tale of the fuzzy subgroup μ. Two fuzzy subgroup μ and η are said to be similarly bounded if they have the same tip and same tale i.e.*

$$\mu(e) = \eta(e) \text{ and}$$
$$\inf \mu = \inf \eta.$$

*The range set of fuzzy subgroup we shall denote by Im μ. It is well known that for a fuzzy subgroup μ the level subset $\mu_t$ for each t ∈ [0, μ(e)] is a subgroup of the given group and is called a level subgroup of μ. The set of all level subgroups of a fuzzy subgroup forms a chain. For μ a fuzzy subgroup of G, the level subset $\mu_t$, for t ∈ Im μ is a subgroup of G and is called the level subgroup of μ.*

Now we proceed on to give a Smarandache analogue for these using S-semigroups as we do not have S-groups.

**DEFINITION 2.2.23:** *Let μ be a fuzzy set in a S-semigroup G. If μ is a S-fuzzy subsemigroup of G associated with a proper subset P ⊂ G, P a subgroup of G. '$e_P$' the identity element of P.*

*That is*

$$\sup_{x \in P} \mu(x) = \mu(e_P)$$

*we call μ($e_P$) to be the Smarandache tip (S-tip) of the S-fuzzy subsemigroup μ of G relative to the subgroup P ⊂ G.*

*Unlike in the case of fuzzy subgroup of a group G we have in case of S-fuzzy subsemigroups of the S-semigroup G several S-tips associated even with a single μ : G → [0, 1] i.e. if $P_1$, ..., $P_K$ are K subgroups in G and if μ is such that its restriction on each $P_i$ happens to be a S-fuzzy subsemigroup of G and if*

$$\sup_{x \in P_i} \mu(x) = \mu(e_{P_i})$$

*where $e_{P_i}$ is the identity element of $P_i$ for i = 1,2,..., K and if they are distinct i.e. $e_{P_i} \neq e_{P_j}$, if j ≠ i. Then we see for a given μ : G → [0 1] we have several S-tips for a given μ. We call such maps μ and such S-semigroups G to be Smarandache multi-tiped fuzzy semigroups (S-multi-tiped fuzzy semigroups).*

**Example 2.2.3:** Let $Z_{12}$ = {0, 1, 2,…, 11}; $Z_{12}$ is a semigroup under multiplication modulo 12. In fact $Z_{12}$ is a S-semigroup. The subgroups of $Z_{12}$ are $P_1$ = {1, 11}, $P_2$ = {3, 9} $P_3$ = {1, 5} , $P_4$ = {1, 7} and $P_5$ = {4, 8}.

Clearly one can define a μ so that $Z_{12}$ has several S-tips for 1, 3 and 4 are the units of the subgroup.

**DEFINITION 2.2.24:** *Let G be a S-semigroup. μ : G → [0, 1] restricted to some subset P of G be a S-fuzzy subsemigroup of G. We see a S-fuzzy subsemigroup as in case of a fuzzy subgroup may or may not attain its infimum. We shall write inf $μ_P$ for $\inf_{x \in P} μ_P(x)$ and refer to it as Smarandache tale (S-tale) of the S-fuzzy subsemigroup μ. Let μ and η be two S-fuzzy subsemigroups of the S-semigroup G related to the same subset P ⊂ G, P a subgroup of G. We say μ and η are Smarandache similarly bounded (S-similarly bounded) if they have the same S-tip and the same S-tale that is $μ_P(e) = η_P(e)$ and inf $μ_P$ = inf $η_P$.*

Now we define Smarandache equivalent of S-fuzzy subsemigroups.

**DEFINITION 2.2.25:** *Let G be a S-semigroup. Two S-fuzzy subsemigroup η and μ of the S-semigroup G related to the same subgroup P, P ⊂ G are said to be Smarandache equivalent (S-equivalent) denoted by $μ \underset{S}{\approx} η$ if μ and η have the same chain of S-level subgroups in P. Thus, we have contrary of groups in case of S-semigroups several S-fuzzy subsemigroups and S-equivalences depending on the subgroups P ⊂ G. The S-relation $\underset{S}{\approx}$ can be proved to be a Smarandache equivalence relation (S-equivalence relation) on the set of S-fuzzy subsemigroups only related to the same subgroup P ⊂ G. Thus on G a S-semigroup we can have several S-relations which are S-equivalence relations depending on P, the proper subset of G which is a subgroup of G. Thus we can say the maximum number of S-relations will correspond to the number of proper subsets in G which are subgroups under the operations of G.*

*All properties pertaining to fuzzy subgroups of a group can be extended to S-fuzzy subsemigroups of a S-semigroup G.*

We proceed on to define Smarandache penultimate subsemigroup of a S-semigroup G.

**DEFINITION 2.2.26:** *Let μ be a fuzzy set in G then the penultimate subset P (μ) of μ in G, defined by P(μ) = {x ∈ G / μ(x) > inf μ}. In case μ is a fuzzy subgroup of a group G then P (μ) is a subgroup of G provided μ is a non-constant and P (μ) is called the penultimate subgroup of μ in G.*

For more about penultimate subgroups please refer [11]. Now on similar lines we define Smarandache penultimate subsemigroups of a S-fuzzy subsemigroup μ of a S-semigroup G.

**DEFINITION 2.2.27:** *Let G be a S-semigroup. Let μ be a S-fuzzy subsemigroup of the S-semigroup G related to a proper subset E of G, E a subgroup of G. $P(μ_E) = \{x ∈ E / μ_E(x) > \text{Im } μ_E\}$.*

*In case $\mu_E$ is a S-fuzzy subsemigroup of the S-semigroup G ($E \subset G$) then $P(\mu_E)$ is a S-fuzzy subsemigroup of E provided $\mu_E$ is non-constant and $P(\mu_E)$ is called the Smarandache penultimate subsemigroup (S- penultimate subsemigroup) of $\mu_E$ in E $\subset G$.*

*It is to be noted by S-penultimate subsemigroup of a S-fuzzy subsemigroup may or may not be a proper subgroup of the group E ($E \subset G$ a subgroup of the S-semigroup G).*

**THEOREM 2.2.56:** *Let $\mu_E$ be a non-constant S-fuzzy subsemigroup of the S-semigroup G. Then the S-penultimate subsemigroup $P(\mu_E)$ is a proper subgroup of $E \subset G$ if and only if $\mu_E$ attains its infimum.*

*Proof*: $\mu : G \to [0, 1]$ is a fuzzy set of G; E a proper subset of G which is a subgroup under the operation of G. $\mu_E$ denotes the restriction map i.e. $\mu_E : E \to [0\ 1]$ such that $\mu_E$ is fuzzy subgroup of G or $\mu$ is a S-fuzzy subsemigroup of G related to E. Now the rest of the proof is as in case of fuzzy subgroups.

**THEOREM 2.2.57:** *Let $\eta_1$, $\eta_2$ and $\mu$ be fuzzy sets in a S-semigroup G such that $\eta_1 \leq \eta_2$. Then $\mu \circ \eta_1 \leq \mu \circ \eta_2$ and $\mu*\eta_1 \leq \mu*\eta_2$ and $\mu \bullet \eta_1 \leq \mu \bullet \eta_2$ on some fixed subgroup E of G.*

*Proof*: First G is a S-semigroup. E a proper subset of G, E a subgroup under the operations of G. $\eta_1$, $\eta_2$ and $\mu$ are fuzzy sets on G satisfying the conidition $\eta_1 \leq \eta_2$ restricted to the subgroup $E \subset G$.

Now using the proof of [7, 10] the result follows:

**THEOREM 2.2.58:** *Let $\eta$ and $\mu$ be S-fuzzy subsemigroups of a S-semigroup of G. Then the set product $\eta \circ \mu$ contains $\eta$ and $\mu$ if and only if $\eta$ and $\mu$ have the same S-tip that is $\eta(e_P) = \mu(e_P)$ for a subgroup P in G.*

*Proof*: As in case of groups. (Hint: Using $P \subset G$ and P is a subgroup $\mu$ and $\eta$ restricted to P serves the purpose).

**THEOREM 2.2.59:** *Let $\eta$ and $\mu$ be S-fuzzy subsemigroups with the same S-tip of a S-semigroup G relative P ($P \subset G$ a subgroup of G). Then the set product $\eta \circ \mu$ is a S-fuzzy subsemigroup generated by the union of $\eta$ and $\mu$ if $\eta \circ \mu$ is a S-fuzzy subsemigroup of G.*

*Proof*: Direct hence left as an exercise for the reader.

**THEOREM 2.2.60:** *Let $\eta$ and $\mu$ be S-fuzzy subsemigroups of a S-semigroup G. Then the S-penultimate product $\eta \circ \mu$ contains $\eta$ and $\mu$ if $\eta$ and $\mu$ are similarly bound with respect to the same subgroup $P \subset G$.*

*Proof*: Please refer [11] using the subgroup P of the S-semigroup G instead of using the whole of G.

**THEOREM 2.2.61:** *Let $\eta$ and $\mu$ be S-fuzzy subsemigroups in a S-semigroup G. Then $\eta \circ \mu = \mu \circ \eta$ if either $\eta$ or $\mu$ is a S-fuzzy ideal of the S-fuzzy subsemigroup.*

*Proof*: If the S-fuzzy subsemigroup $\mu$ related to P, P $\subset$ G a subgroup of G then we see $\mu$ can also be fuzzy normal so that $\mu$ is a S-fuzzy ideal of the S-fuzzy subsemigroup.

**THEOREM 2.2.62:** *Let $\eta$ and $\mu$ be S-fuzzy subsemigroups of a S-semigroup G. Then $\eta \bullet \mu = \mu \bullet \eta$ if either $\eta$ or $\mu$ is S-fuzzy ideal of the S-fuzzy subsemigroup.*

*Proof*: Analogous to the proof in case of fuzzy subgroups as S-fuzzy subsemigroup $\mu$ of G is also a fuzzy subgroup of G.

**THEOREM 2.2.63:** *Let $\eta$ and $\mu$ be S-fuzzy subsemigroups of a S-semigroup G. Then the set product $\mu \circ \eta$ is a S-fuzzy subsemigroup if and only if $\eta \circ \mu = \mu \circ \eta$.*

*Proof*: Straightforward, hence left for the reader to prove.

**THEOREM 2.2.64:** *Let $\eta$ and $\mu$ be S-fuzzy subsemigroup of a S-semigroup G. Then the S-penultimate product $\eta \bullet \mu$ is a S-fuzzy subsemigroup if and only if $\eta \bullet \mu = \mu \bullet \eta$.*

*Proof*: As in case of fuzzy subgroups by using the very definition.

*Example 2.2.4:* Let S(4) be the S-symmetric semigroup. Let $\mu$ and $\eta$ be S-fuzzy subsemigroups of S(4) relative to the proper subset $S_4$ of S(4) which is a subgroup.

For

$$\eta(e) = \frac{3}{16},$$
$$\eta\begin{pmatrix} 1 & 2 & 3 & 4 \\ 3 & 4 & 1 & 2 \end{pmatrix} = \frac{9}{16},$$
$$\eta\begin{pmatrix} 1 & 2 & 3 & 4 \\ 2 & 1 & 4 & 3 \end{pmatrix} = \frac{3}{16},$$
$$\eta\begin{pmatrix} 1 & 2 & 3 & 4 \\ 4 & 3 & 2 & 1 \end{pmatrix} = \frac{3}{16},$$
$$\eta(p) = \frac{1}{16}$$

for all odd permutations of (1234) and $\eta(\beta) = 0$ for any other permutation $\beta$ in $S_4$.

$$\mu\left(\begin{pmatrix} 1 & 2 & 3 & 4 \\ 3 & 4 & 1 & 2 \end{pmatrix}\right) = \frac{11}{16},$$

$$\mu\left(\begin{pmatrix} 1 & 2 & 3 & 4 \\ 2 & 1 & 4 & 3 \end{pmatrix}, \begin{pmatrix} 1 & 2 & 3 & 4 \\ 4 & 3 & 2 & 1 \end{pmatrix}\right) = \left\{\frac{7}{16}\right\}$$

$$\mu\left\{\begin{pmatrix} 1 & 2 & 3 & 4 \\ 2 & 3 & 4 & 1 \end{pmatrix}, \begin{pmatrix} 1 & 2 & 3 & 4 \\ 4 & 3 & 2 & 1 \end{pmatrix}, \begin{pmatrix} 1 & 2 & 3 & 4 \\ 3 & 2 & 1 & 4 \end{pmatrix}, \begin{pmatrix} 1 & 2 & 3 & 4 \\ 1 & 4 & 3 & 2 \end{pmatrix}\right\} = \left\{\frac{5}{16}\right\}$$

and $\mu(\gamma) = 0$ for any other permutation $\gamma$ in $S_4$.

We define fuzzy subgroups $\overline{\eta}$ and $\overline{\mu}$ as follows:

$$\overline{\eta}(e) = \frac{13}{16},$$

$$\overline{\eta}\left(\begin{pmatrix} 1 & 2 & 3 & 4 \\ 3 & 4 & 1 & 2 \end{pmatrix}\right) = \frac{9}{16}$$

$$\overline{\eta}\left\{\begin{pmatrix} 1 & 2 & 3 & 4 \\ 2 & 1 & 4 & 3 \end{pmatrix}, \begin{pmatrix} 1 & 2 & 3 & 4 \\ 4 & 3 & 2 & 1 \end{pmatrix}\right\} = \frac{3}{16}$$

$\overline{\eta}(\alpha) = \frac{2}{16}$ for each 3 cycle. $\overline{\eta}(\beta) = \frac{1}{16}$ for any other permutation $\beta$ in $S_4$.

$$\overline{\mu}\left\{\begin{pmatrix} 1 & 2 & 3 & 4 \\ 1 & 2 & 3 & 4 \end{pmatrix}, \begin{pmatrix} 1 & 2 & 3 & 4 \\ 3 & 4 & 1 & 2 \end{pmatrix}\right\} = \frac{11}{16}$$

$$\overline{\mu}\left\{\begin{pmatrix} 1 & 2 & 3 & 4 \\ 2 & 1 & 4 & 3 \end{pmatrix}, \begin{pmatrix} 1 & 2 & 3 & 4 \\ 4 & 3 & 2 & 1 \end{pmatrix}\right\} = \frac{7}{16}$$

$$\overline{\mu}\left\{\begin{pmatrix} 1 & 2 & 3 & 4 \\ 2 & 3 & 4 & 1 \end{pmatrix}, \begin{pmatrix} 1 & 2 & 3 & 4 \\ 4 & 3 & 2 & 1 \end{pmatrix}, \begin{pmatrix} 1 & 2 & 3 & 4 \\ 3 & 4 & 1 & 2 \end{pmatrix}\right\} = \frac{5}{16}$$

and $\overline{\mu}(\gamma) = \frac{1}{16}$ for any other permutation $\gamma$ in $S_4$.

Free product: The free product of $\overline{\eta}$ and $\overline{\mu}$ turns out to be constant fuzzy set

$$(\overline{\eta} * \overline{\mu})(g) = \frac{11}{16}$$

for each $g \in S_4$. Set product can be verified to be

$$(\overline{\eta} o \overline{\mu}) \left\{ \begin{pmatrix} 1 & 2 & 3 & 4 \\ 1 & 2 & 3 & 4 \end{pmatrix}, \begin{pmatrix} 1 & 2 & 3 & 4 \\ 3 & 4 & 1 & 2 \end{pmatrix} \right\} = \frac{11}{16},$$

$$(\overline{\eta} o \overline{\mu}) \left\{ \begin{pmatrix} 1 & 2 & 3 & 4 \\ 2 & 1 & 4 & 3 \end{pmatrix}, \begin{pmatrix} 1 & 2 & 3 & 4 \\ 4 & 3 & 2 & 1 \end{pmatrix} \right\} = \left\{ \frac{7}{16} \right\}$$

$$(\overline{\eta} o \overline{\mu})(\alpha) = \frac{2}{16}$$

for any permutation $\alpha$ in $S_4$.

It can be verified that neither $\overline{\eta}$ nor $\overline{\mu}$ is a S-fuzzy ideal of the S-fuzzy subsemigroup $\overline{\mu}$ and $\overline{\eta}$ of $S_4$ but $\overline{\eta} o \overline{\mu} = \overline{\mu} o \overline{\eta}$.

The S-permulate subsemigroup $P(\overline{\eta})$ and $P(\overline{\mu})$ of $\overline{\eta}$ and $\overline{\mu}$ respectively are $A_4$ and $D_4$.

**THEOREM 2.2.65:** *Let $\eta$ and $\mu$ be the fuzzy set with sup property in a S-semigroup G relative to a subgroup $P \subset G$. Then for each $t \in [0, 1]$ ; $(\hat{\eta} o \hat{\mu})_t = \hat{\eta}_t \hat{\mu}_t$ where $\hat{\eta}$ and $\hat{\mu}$ denote the restriction of $\eta$ and $\mu$ respectively over $P \subset G$ i.e. $\hat{\eta} : P \to [0,1]$ and $\hat{\mu} : P \to [0,1]$ with a sup property in the group P.*

*Proof*: Analogous to proof in case of sup property in a group G.

**THEOREM 2.2.66:** *Let $\eta$ and $\mu$ be fuzzy sets in a S-semigroup G with S-penultimate subsets $P(\eta)$ and $P(\mu)$ respectively. Then*

i. $(\eta o \mu)_t^> = \eta_t^> P(\mu) \cap P(\eta) \mu_t^>$ for $t \in [0, 1]$
ii. $(\eta o \mu)_t = \eta_t P(\mu) \cap P(\eta) \mu_t$ for $t \in [0, 1]$ provided $\eta$ and $\mu$ have sup property.

*Proof*: Easily proved by using all the properties and definition.

**THEOREM 2.2.67:** *Let $\eta$ and $\mu$ be S-fuzzy ideals of the S-fuzzy subsemigroup of the S-semigroup G. Then the free product $\eta * \mu$ is a S-fuzzy ideal.*

*Proof:* As $\eta * \mu$ is a S-fuzzy subsemigroup of the S-semigroup G relative to P we see $\eta * \mu$ is a fuzzy subgroup of P. The rest of the proof follows as in case of groups.

We just recall the notion of Smarandache normal subgroup of a S-semigroup.

**DEFINITION 2.2.28:** *Let S be a S-semigroup. Let A be a proper subset of S which is a group under the operations of S. We say A is a Smarandache normal subgroup (S- normal subgroup) of the S-semigroup S if $xA \subset A$ and $Ax \subseteq A$ or $xA = \{0\}$ and $Ax = \{0\}$ for all $x \in S$ and if 0 is an element in S then we have $xA = \{0\}$ and $Ax = \{0\}$.*

*Let S be a S-semigroup μ : S → [0, 1] be a S-fuzzy semigroup of S relative to P, P ⊂ S; P a subgroup under the operations of S. i.e. μ$_P$ : P → [0, 1] is a fuzzy group i.e. μ$_P$ is nothing but μ restricted to P. If in addition we have μ (gx) = μ (xg) for every g ∈ P and x ∈ S then we call μ a Smarandache fuzzy normal subgroup (S-fuzzy normal subgroup) of the S-semigroup S.*

## 2.3 Element-wise properties of S-fuzzy subsemigroups

In this section we define several of the element-wise properties in S-fuzzy subsemigroups.

**DEFINITION 2.3.1:** *Let G be a S-semigroup. A fuzzy subset λ of G of the form*

$$\lambda(y) = \begin{cases} t(\neq 0) & \text{if } y = x \\ 0 & \text{if } y = x. \end{cases}$$

*For y in P ⊂ G, P a subgroup in G is called the Smarandache fuzzy point (S-fuzzy point) with Smarandache support (S-support) x and value t is denoted by $x_t$.*

*It is very important to note that all elements in the S-semigroup G need not in general have S-fuzzy point. That is to be more precise only those elements in G which fall into one or more subgroups of G will have S-fuzzy point that too depending on the fuzzy subset λ of G. Thus unlike in a group G we will not be in a position to associate S-fuzzy points.*

**DEFINITION 2.3.2:** *A fuzzy subset λ of a S-semigroup G is said to have Smarandache sup property (S-sup property) if for every non empty subset T of G where T is a subgroup and there exists a ∈ T such that λ(a) = sup {λ(t) | t ∈ T}.*

**DEFINITION 2.3.3:** *A fuzzy subset λ of a S-semigroup G is said to be an Smarandache (∈, ∈∨q) fuzzy subsemigroup (S-(∈, ∈∨q) fuzzy subsemigroup) of G if for any x, y ∈ P ⊂ G (P a subgroup of G) and t, r ∈ (0, 1]*

    i.    $x_t, y_t \in \lambda \Rightarrow (xy)_{M(t,r)} \in \vee q\lambda$ *and*
    ii.    $x_t \in \lambda \Rightarrow (x^{-1})_t \in \vee q\lambda.$

**DEFINITION 2.3.4:** *A S-fuzzy subsemigroup λ of G is said to be*

    i.    *Smarandache (∈,∈) fuzzy normal (S-(∈,∈) fuzzy normal) if for all x, y ∈ P ⊂ G and t ∈ (0,1], $x_t \in \lambda \Rightarrow (x^{-1}y x)_t \in \lambda$.*

    ii.    *Smarandache (∈,∈∨q) fuzzy normal (S-(∈,∈∨q) fuzzy normal) (or simply S-fuzzy normal) if for any x, y ∈ P ⊂ G and t ∈ (0, 1], $x_t \in \lambda \Rightarrow (y^{-1}x y)_t \in \vee q\lambda.$*

Several relations and results in this direction can be evolved using Smarandache fuzzy notions.

**THEOREM 2.3.1:**

i. *A fuzzy subset $\lambda$ of a S-semigroup G is an S- $(\in, \in \vee q)$ fuzzy subsemigroup (respectively S- $(\in, \in \vee q)$ fuzzy normal subgroup) related to a subgroup P of G if and only if the $(\in \vee q)$- level subset $\lambda_q$ of $P \subset G$ is a subgroup of P (normal subgroup of P) for all $t \in (0, 1]$ where $\lambda_q = \{ x \in P \mid \lambda(x) \geq t$ or $\lambda(x) + t > 1\}$.*

ii. *Let $f : G \to H$ be a S-semigroup homomorphism. Let $\lambda$ and $\mu$ be S-fuzzy subsemigroups of G and H respectively. Then $f(\lambda)$ and $f^{-1}(\mu)$ are S-fuzzy subsemigroups of $f(P)$ and P respectively ($P \subset G$ is a subgroup of G).*

*Proof:* Refer [25].

**THEOREM 2.3.2:** *For any S-$(\in, \in \vee q)$ fuzzy subsemigroup $\lambda$ of S-semigroup G the S-$(\in \vee q)$-level subgroup $\lambda_q = P$ ($P \subset G$, P a subgroup of G) for all $t \in I$ if and only if $\lambda(x) \geq 0.5$ for all $x \in P \subset G$.*

*Proof:* Straightforward.

**THEOREM 2.3.3:** *Let G be a S-semigroup. Then given any chain of subgroups $P_0 \subset P_1 \subset ... \subset P_r = P$ (P a subgroup of G) there exists a S- $(\in, \in \vee q)$- fuzzy subsemigroup of $P \subset G$ whose S- $(\in \vee q)$-level subgroups are precisely the members of the chain.*

*Proof:* Let $\{t_i \mid t_i \in (0, 0.5); i = 1, 2, ..., r\}$ be such that $t_1 > t_2 > ... > t_r$. Let $\lambda : P \to I$ be defined as follows.

$$\lambda(x) = \begin{cases} t > 0.5 & \text{if } x = e \\ \mu > t & \text{if } x \in P_0 - \{e\} \\ t_1 & \text{if } x \in P_1 - P_0 \\ t_2 & \text{if } x \in P_2 - P_1 \\ \vdots & \\ t_r & \text{if } x \in P_r - P_{r-1} \end{cases}$$

Then $\lambda$ is an S-$(\in, \in \vee q)$-fuzzy subsemigroup of $P \subset G$ which is not a S- $(\in, \in)$ fuzzy subsemigroup. Note that $\lambda_{0.5} = P_0$ and $\lambda_{t_i} = P_i$ for $i = 1, 2, ..., r$ as $x_t \subset \vee q\lambda \Rightarrow x_t \subset \lambda$, if $t \in (0, 0.5)$.

It is important to mention here that for a given S-semigroup G we can have several subgroups P, Q, R,..., X, Y, Z such that we have a chain of subgroups associated with each of the subgroups P, Q, R,..., X, Y, Z leading to many S-$(\in, \in \vee q)$ fuzzy subgroup P, Q, R,..., X, Y, Z of G.

Thus only the Smarandache notion alone can give a nice spectrum of S-($\in,\in \vee q$) fuzzy subsemigroups relative to a S-semigroup G having several distinct proper subsets which are subgroups of G.

*Example 2.3.1:* Let G be the set of integers together with $\left\{0, \pm\frac{1}{2}, \pm\frac{1}{2^2}, \cdots, \pm\frac{1}{2^n} \pm \cdots\right\}$. G is a S-semigroup. Take P = additive group of all integers. P is a subgroup of G. Let nP = {additive group of all integers multiple of n}. Then $16P \subset 8P \subset 4P \subset 2P \subset P$ be a chain of subgroups of G.

Let $\lambda : P \to [0, 1]$ defined by

$$\lambda(x) = \begin{cases} 0.6 & \text{if } x = 0 \\ 0.9 & \text{if } x \neq 0, x \in 16P \\ 0.7 & \text{if } x \in 8P - 16P \\ 0.5 & \text{if } x \in 4P - 8P \\ 0.2 & \text{if } x \in 2P - 4P \\ 0.1 & \text{if } x \in P - 2P. \end{cases}$$

Then $\lambda$ is a S-fuzzy subsemigroup of P. Note that $\lambda_{0.25} = 4P = \lambda_{0.6} = \lambda_{0.7}$, $\lambda_{0.9} = 2P = \lambda_{0.2}$, $\lambda_{0.1} = P$.

**THEOREM 2.3.4:** *Let $\lambda$ be a S-($\in, \in \vee q$)-fuzzy subgroup of G with Im $\mu = \{t, r\}$ where $0 < t < r < 0.5$. If $\lambda = \mu \cup \nu$ where $\mu, \nu$ are S-($\in,\in \vee q$)-fuzzy subsemigroups of G then either $\mu \leq \nu$, or $\nu \leq \mu$.*

*Proof:* Follows as in case of ($\in,\in \vee q$)-fuzzy subgroup. We say S-($\in, \in \vee q$)-fuzzy subgroup of a S-semigroup G is said to be proper if it is not constant in G.

**THEOREM 2.3.5:** *Let G be a S-semigroup that has several proper subgroup. A proper S-($\in,\in$) fuzzy subsemigroup $\lambda$ of G such that card. Im $\lambda \geq 3$ can be expressed as the union of two proper non-equivalence S-($\in,\in$)- fuzzy subsemigroups of G.*

*Proof:* As in case of ($\in,\in$)-fuzzy subgroups.

**THEOREM 2.3.6:** *Let $\lambda$ be a proper S-($\in,\in \vee q$)-fuzzy subgroup of G such that the cardinality of $\{\lambda(x) \mid \lambda(x) < 0.5 \} \geq 2$. Then $\lambda$ can be expressed as the union of two proper non-($\in \vee q$)-equivalent S-($\in,\in \vee q$) fuzzy subsemigroups of G.*

*Proof:* Follows by the same arguments as in case of fuzzy subgroups.

**THEOREM 2.3.7:** *Let $\lambda$ be an S-($\in,\in \vee q$)-fuzzy subsemigroup of the S-semigroup G such that $\lambda(x) \geq 0.5$ for all $x \in G$. Then $\lambda$ can be expressed as the union of two non-($\in \vee q$) equivalent S- ($\in,\in \vee q$)-fuzzy subsemigroups if and only if $P = H \cup K$, where H and K are proper S- ($\in,\in \vee q$)-fuzzy subsemigroups of P (where P is a proper subset of G which is a subgroup of G).*

*Proof:* For proof please refer [23].

**DEFINITION 2.3.5:** *Let $\mu$ be a fuzzy subset of a S-semigroup G. An S- $(\in, \in \vee q)$ fuzzy subsemigroup $\xi$ of G is said to be the S-$(\in, \in \vee q)$-fuzzy subsemigroup generated by $\mu$ in G if $\xi \geq \mu$ and for any other S- $(\in, \in \vee q)$-fuzzy subsemigroup $\eta$ of G with $\eta \geq \mu$ it must be $\eta \geq \xi$.*

**THEOREM 2.3.8:** *Let $\mu$ be a fuzzy subset of the S-semigroup G where card Im $\mu$ is finite. Define S-subsemigroups $G_i$ of G as follows.*

$$\begin{aligned} G_0 &= \langle \{x \in G \,/\, \mu(x) \geq 0.5 \} \rangle \\ G_1 &= \langle \{G_0 \cup \{x \in G \,/\, \mu(x) = \sup(\mu(x) \,;\, x \in G \setminus G_0)\}\} \rangle \\ &\vdots \\ G_i &= \langle \{G_{i-1} \cup \{x \in G \,|\, \mu(x) = \sup\{\mu(z), z \in G - G_{i-1}\}\}\} \rangle \end{aligned}$$

*$i = 1, 2, ..., k$ where $k \leq$ Card Im $\mu$ and $G_K = G$. $P_i \subset G_i$ is a proper subset of $G_i$ which is a subgroup. Then the fuzzy subset $\theta$ of G defined by*

$$\theta(x) = \begin{cases} \mu(x) & \text{if } x \in P_o \subset G_o \\ \sup\{\mu(z); z \in G \setminus G_{i-1}\} & \text{if } x \in P_i - P_{i-1} \subseteq G_i - G_{i-1} \end{cases}$$

*(where $P_i$ is a subgroup in $G_i$, for $i = 0, 1, 2, ..., k$ is the S- $(\in, \in \vee q)$-fuzzy subsemigroup generated by $\mu$ in G.*

*Proof:* Left as an exercise for the reader to prove.

**THEOREM 2.3.9:** *Let $f : G \to K$ be a S-semigroup homomorphism. If $\lambda$ and $\mu$ are S-fuzzy subsemigroups of G and K respectively, then $f(\lambda)$ and $f^{-1}(\mu)$ are S-fuzzy subsemigroups of K and G respectively.*

*Proof:* Please refer [23].

**THEOREM 2.3.10:** *Let $f : G \to H$ be a S-semigroup homomorphism. Let $\lambda$ and $\mu$ be two S-fuzzy subsemigroups of G and H respectively. Then*

  i.   *$f^{-1}(\mu)$ is a S-fuzzy subsemigroup of G.*
  ii.  *$f(\lambda)$ is a S-fuzzy subsemigroup of $f(G)$.*

*Proof:* Straightforward by the earlier result and definitions.

Now we proceed on to define Smarandache $(\in, \in)$-fuzzy left (resp. right) cosets.

**DEFINITION 2.3.6:** *Let $\lambda$ be a S-fuzzy subsemigroup of a S-semigroup G. For $x \in P \subset G$, P-a subgroup of G.*

$\lambda_x^l$ (resp. $\lambda_x^r$) : $G \to [0,1]$ be defined by $\lambda_x^l(g) = \lambda(gx^{-1})$ (resp. $\lambda_x^r(g) = \lambda(x^{-1}g)$ for all $g \in G$ is called Smarandache ($\in, \in$) fuzzy left (resp. right coset) (S-($\in,\in$) fuzzy left (resp. right coset)) of P in G determined by x and $\lambda$.

Let $\lambda$ be a S-fuzzy subsemigroup of the S-semigroup G. Then $\lambda$ is a S-($\in,\in$)-fuzzy normal if and only if $\lambda_x^l = \lambda_x^r$ for all $x \in P \subsetneq G$.

However if $\lambda$ is S-($\in, \in \vee q$) fuzzy normal, then $\lambda_x^l$ may not be equal to $\lambda_x^r$.

This is illustrated by the following example.

*Example 2.3.2:* Let G be a S-semigroup given by the following table.

| * | 0 | e | a | b | c | d | f | g |
|---|---|---|---|---|---|---|---|---|
| 0 | 0 | 0 | 0 | 0 | 0 | 0 | 0 | 0 |
| e | 0 | e | a | b | c | d | f | e |
| a | 0 | a | b | e | f | c | d | d |
| b | 0 | b | e | a | d | f | c | b |
| c | 0 | c | d | f | e | a | b | a |
| d | 0 | d | f | c | b | e | a | g |
| f | 0 | f | c | d | a | b | e | c |
| g | 0 | g | a | b | d | e | c | 0 |

H = {e, a, b, c, d, f} is a subgroup of G. Define $\lambda$ : G $\to$ [0, 1] by $\lambda$, (e a b c d f) = $\lambda$(0.7, 0.75, 0.8. 0.4, 0.4, 0.4). $\lambda$(0) = 0 and $\lambda$(g) = 0. Thus $\lambda$ is an S- ($\in,\in \vee q$)- fuzzy normal subgroup of G. But since $\lambda$(fd) = $\lambda$(b) = 0.8 $\neq$ 0.75 = $\lambda$(df) = $\lambda$ (a). $\lambda$ is not S-($\in,\in$) fuzzy normal.

**DEFINITION [23]:** *Let $\lambda$ be a S-fuzzy subsemigroup of the S-semigroup G. For any $x \in P \subset G$, $\hat{\lambda}_x$ (resp. $\tilde{\lambda}$) : G $\to$ [0, 1] is defined by $\hat{\lambda}_x(g) = M(\lambda(gx^{-1}, 0.5)$ [resp. $\hat{\lambda}_x(g) = M(\lambda(x^{-1}g, 0.5))$] for all $g \in P \subseteq G$ and is called S- ($\in,\in \vee q$)- fuzzy left (resp. right coset) of $P \subset G$ determined by x and $\lambda$.*

**THEOREM [23]:** *Let $\lambda$ be a S-fuzzy subsemigroup of a S-semigroup G. Then $\lambda$ is S-($\in,\in \vee q$)- fuzzy normal if and only if $\hat{\lambda}_x = \tilde{\lambda}_x$ for all $x \in P \subset G$.*

*Proof:* Refer [23].

**THEOREM [23]:** *Let $\lambda$ be a S-fuzzy normal subgroup of the S-semigroup G. Let F be the set of all S-fuzzy cosets of $\lambda$ in G. Then F is a S-semigroup of all fuzzy cosets of G determined by $\lambda$ where the multiplication is defined by $\lambda_x o \lambda_y = \lambda_{xy}$ for all x, y $\in P$, ($P \subset G$).*

Let $\bar{\lambda}$ : $F \to [0,1]$ be defined by $\bar{\lambda}(\lambda_x) = M(\lambda(x^{-1}), 0.5)$ for all $x \in P$. Then $\lambda$ is a S-fuzzy normal subgroup in F.

*Proof:* Refer [23].

**THEOREM 2.3.11:** *Let $\lambda$ be an S-$(\in, \in \vee q)$- fuzzy normal subgroup of the S-semigroup G and $\lambda_t = H$ where $\lambda_t$ is the S-$(\in \vee q)$ level subgroup of G for $t \in (0, 1]$. If a, b $\in P \subset G$ are such that $\lambda_a = \lambda_b$; then $H_a = H_b$.*

*Proof:* Refer [25].

**DEFINITION 2.3.7:** *A S-subsemigroup H of a S-semigroup G is said to be Smarandache quasi normal (S-quasi normal) if for every S-subsemigroup K of H, we have $XY = YX$, $X \subset H$ and $Y \subset K$ are subgroups of X and Y respectively.*

**THEOREM 2.3.12:** *The S-homomorphic pre-image of a S-quasi normal subgroup is S-quasi normal.*

*Proof:* Follows by very definitions.

**THEOREM 2.3.13:** *A maximal S-quasi normal subgroup of a S-semigroup G is normal.*

*Proof:* Left as an exercise for the reader.

**THEOREM 2.3.14:** *Let G be a finite S-semigroup and Q be a S-quasi normal subgroup T of G. If Q is core-free then every Sylow subgroup in T is S-quasi normal in $T \subset Q$.*

*Proof:* Follows by simple computations.

**DEFINITION 2.3.8:** *A S-fuzzy subsemigroup $\xi$ of G (G a S-subsemigroup) relative to a subgroup $X \subset G$ is called Smarandache fuzzy quasi normal I (S-fuzzy quasi normal I) in G if $\xi \circ \eta = \eta \circ \xi$ for every S-fuzzy subsemigroup $\eta$ of X in G.*

**DEFINITION 2.3.9:** *A S-fuzzy subsemigroup $\mu$ of a S-semigroup G is said to be Smarandache fuzzy maximal (S-fuzzy maximal) if $\mu$ is not constant for any S-fuzzy subsemigroup $\eta$ of G whenever $\mu \leq \eta$ either $P_\mu = P_\eta$ ($P \subset G$, P a subgroup of G relative to which $\mu$ and $\eta$ are defined or $\eta = \chi_P$ where $P_\mu = \{x \in P \mid \mu(x) = \mu(e_p)\}$, $e_P$ identity element of the subgroup P in G.*

**DEFINITION 2.3.10:** *Let $\mu$ be a S-fuzzy subsemigroup of a finite S-semigroup G and let $S_p$ be a p-sylow subgroup of P, P a subgroup in G. Define a fuzzy subset $\mu S_p$ in G as follows*

$$\mu S_p(x) = \begin{cases} \mu(x) & \text{if } x \in S_p \\ 0 & \text{if } x \notin S_p. \end{cases}$$

*Clearly $\mu S_p$ is a S-fuzzy subsemigroup called the Smarandache p-fuzzy sylow subgroup (S-p-fuzzy Sylow subgroup) of $\mu$. The following result can be proved as a matter of routine.*

**THEOREM 2.3.15:** *Let $\xi$ be a S-fuzzy subsemigroup of the S-semigroup G with S-sup property. Then each level subset $\mu_t$, for $t \in (0, \xi(e))$ is a S-quasi normal subgroup of $P \subset G$ ($\xi$ is a S-fuzzy subsemigroup related to $P \subset G$) if and only if $\xi$ is S-fuzzy quasi normal (I).*

**DEFINITION 2.3.11:** *A S-fuzzy semigroup of a S-semigroup G is called a S-fuzzy quasi normal II if its S-level subsemigroups are S-quasi normal subsemigroups of G.*

Now a characterization theorem for this $\xi$ to be S-fuzzy quasi normal II is given by the following Theorem; the proof of which is left for the reader.

**THEOREM 2.3.16:** *Let $\xi$ be a S-fuzzy subsemigroup of G, G a S-semigroup. Then $\xi$ is S-fuzzy quasi normal II if and only if $\xi \circ \eta = \eta \circ \xi$ (restricted to some subgroup $P \subset G$ relative to which $\xi$ is defined) for any S-subsemigroup $\eta$ of G.*

*From now onwards G will denote a S-semigroup; P a subgroup of G ($P \subset G$) under the operations of G. Let $e_P$ be the identity element of P and $\lambda_P$ an S-$(\in, \in \vee q)$- fuzzy subsemigroup of G. By a S-fuzzy subsemigroup of G we shall mean a S-$(\in, \in \vee q)$- fuzzy subsemigroup of G. For a fuzzy subgroup $\lambda$ of G there exists x and $y \in G$ such that $\lambda(x) \geq 0.5$ and $\lambda(y) < 0.5$.*

**DEFINITION 2.3.12:**

i. $\lambda$ is said to be a S-$(\in, \in)$-fuzzy quasi normal if for any S-fuzzy subsemigroup $\mu$ of G and for all $z \in P \subset G$, $t \in (0,1]$, $z_t \in (\lambda \circ \mu)$ if and only if $z_t \in (\mu \circ \lambda)$.

ii. $\lambda$ is said to be a S-$(\in, \in \vee q)$ fuzzy quasi normal if for any S-fuzzy subsemigroup $\mu$ of G and for all $z \in P \subset G$, $t \in (0,1]$, $z_t \in (\lambda \circ \mu)$ implies $z_t \in \vee q (\mu \circ \lambda)$ and $z_t \in (\mu \circ \lambda)$ implies $z_t \in \vee q (\lambda \circ \mu)$.

**THEOREM 2.3.17:** *$\lambda$ is S-$(\in, \in)$ fuzzy quasi normal if and only if for any S-fuzzy subsemigroup $\mu$ of G and for all $z \in P \subset G$, $(\lambda \circ \mu) z = (\mu \circ \lambda) z$.*

*Proof*: By way of contradiction if $(\lambda \circ \mu)(z) \neq (\mu \circ \lambda)(z)$ for some $z \in P \subseteq G$ and we arrive at a contradiction.

***Remark***: Let $\lambda$ be an S-$(\in, \in)$ fuzzy subsemigroup of a S-semigroup G. If $\lambda$ is a S-$(\in, \in)$-fuzzy quasi normal then $\lambda$ is a S-fuzzy quasi normal subsemigroup of type I or II.

As an S-($\in,\in\vee q$) fuzzy subsemigroup is different from the S-fuzzy semigroup it follows that an S-($\in,\in$) fuzzy quasi normal for an S-($\in,\in\vee q$) fuzzy subsemigroup is different from a S-fuzzy quasi normal subsemigroup of type I or II.

**THEOREM 2.3.18**: *$\lambda$ is S-($\in,\in\vee q$)-fuzzy quasi normal if and only if for any S-fuzzy subsemigroup $\mu$ of G relative to a subgroup $P \subset G$ and for all $z \in P$, $(\mu \circ \lambda)(z) = (\lambda \circ \mu)(z)$ if $(\lambda \circ \mu)(z) < 0.5$ and $(\mu \circ \lambda)(z) \geq 0.5$ if $(\lambda \circ \mu)(z) \geq 0.5$.*

*Proof*: Let $\lambda$ be S-($\in,\in\vee q$) fuzzy quasi normal. Let $z \in P \subset G$.

Case 1 (i): Let $(\lambda \circ \mu)(z) < 0.5$. If possible let $(\mu \circ \lambda)(z) \neq (\lambda \circ \mu)(z)$. Suppose (i) $(\mu \circ \lambda)(z) = t < (\lambda \circ \mu)(z) = r$. Then $z_r \in (\lambda \circ \mu)$ but $z_r \overline{\vee} q\, (\mu \circ \lambda)$ a contradiction.

Case 1 (ii) : Let $r < t$. Then there exists $\delta > 0$ such that $r + \delta < t$ and $r + \delta + t \leq 1$.

Now $z_{r+\delta} \in (\mu \circ \lambda)$ but $z_{r+\delta} \overline{\in \vee q}\, (\lambda \circ \mu)$ a contradiction. So $(\lambda \circ \mu)(z) = (\mu \circ \lambda)(z)$ for all $z \in P \subset G$.

Case 2: Let $(\lambda \circ \mu)(z) = t > 0.5$. If possible let $(\mu \circ \lambda)(z) < 0.5$. Then $z_{0.5} \in (\lambda \circ \mu)$ but $z_{0.5} \overline{\in \vee q}\, (\mu \circ \lambda)$ a contradiction. So $(\mu \circ \lambda)(z) > 0.5$. The converse part follows similarly.

**THEOREM 2.3.19:** *If $\mu$ is a S-fuzzy quasi normal subsemigroup (type I or II) of $P \subset G$, then $\lambda$ is an S- ($\in,\in\vee q$)-fuzzy quasi normal subsemigroup of G.*

*Proof*: Straightforward as in case of subgroups. The following theorems are stated and the proof is left as an exercise for the reader.

**THEOREM 2.3.20:** *Let H be any non-empty subset of the S-semigroup G. H is a S-quasi normal subgroup of G if and only if $\chi_H$ ( the characteristic function of H) is an S-($\in,\in\vee q$) fuzzy quasi normal subgroup of G.*

**THEOREM 2.3.21:** *Let $\lambda$ be a S-($\in,\in\vee q$)-fuzzy quasi normal subsemigroup of G, G a S-semigroup with the sup property. Then $(\in\vee q)$-level subset $\lambda_t = \{x \in P \mid x_t \in\vee q\lambda\}$ is a S-quasi normal subsemigroup G relative to $P \subset G$ ( P a subgroup of G) for all $t \in (0, 1]$.*

**THEOREM 2.3.22:** *Let $\lambda$ be a fuzzy subset of the S-semigroup G with "sup-property" and for all $t \in (0, 1]$, S-($\in\vee q$)-level subset $\lambda_t$ be a S-quasi normal subsemigroup of G.*

*Proof*: $\lambda$ is a S-fuzzy subsemigroup of G relative to a proper subset $P \subset G$; P a subgroup of G. Let $\mu$ be any S-fuzzy subsemigroup of G and $z \in P$. Then $z_t \in (\lambda \circ \mu)$ implies $z \in (\lambda \circ \mu)_t$ implies $z \in \lambda_t \bullet \mu_t$ implies $z \in \mu_t \bullet \lambda_t$. (Since $\lambda_t$ is S-quasi normal) implies $z \in (\mu \circ \lambda)_t$ implies $z_t \in\vee q\, (\mu \circ \lambda)$. Similarly $z_t \in (\mu \circ \lambda) \Rightarrow z_t \in\vee q\, (\lambda \circ \mu)$. Therefore $\lambda$ is S-($\in,\in\vee q$)-fuzzy quasi normal.

**THEOREM 2.3.23:** *Let G be a S-semigroup such that for all x, y, υ, ν ∈ G, xy • υν = xυ • yν. If λ μ are S-(∈, ∈∨q) fuzzy quasi normal subsemigroups of G then so also is λ o μ..*

*Proof*: Refer [25] and obtain the proof analogous for S-semigroup G with suitable modifications.

*Remark*: Let G be a S-semigroup such that for all x, y, υ,ν ∈ P (P a subgroup of G), xy • υν = xυ • yν and Q(P) is the set of all S-(∈,∈∨q)-fuzzy quasi-normal subsemigroup of P then (Q(P), •) is a commutative semigroup.

Thus for a single S-semigroup G we have several commutative semigroups associated with them in fact depending on the number of nontrivial subsets P of G which are subgroups of G.

**THEOREM 2.3.24:** *Let H be a S-semigroup. G a S-semigroup, suppose f is a S-semigroup homomorphism from G onto H (i.e. f : P ⊂ G → L ⊂ H i.e. P and L are subgroups of G and H respectively and f is a onto group homomorphism i.e. equivalently f is a S-semigroup homomorphism); then*

  i. *f (λ o θ ) = f(λ) o f(θ ) where λ and θ are S-fuzzy subsemigroups of P ⊂ G.*
  ii. *f (f$^{-1}$(θ ) = θ where θ is a S-fuzzy subsemigroup of f(P).*

*Proof*: Direct as in case of fuzzy subgroups.

**THEOREM 2.3.25:** *Let f : G → H (H a S-semigroup), be a S-semigroup homomorphism). If λ is a S-(∈,∈∨ q)-fuzzy quasi normal subgroup of P ⊂ G then f (λ) is an S-(∈,∈∨ q)-fuzzy quasi normal subgroup of f (P). (P ⊂ G is a subgroup of the S-semigroup G).*

*Proof*: As in case of (∈,∈∨q) fuzzy normal subgroups of G.

**THEOREM 2.3.26**: *Let f : G → H be a S-semigroup homomorphism and μ be a S-(∈,∈∨q)-fuzzy quasi normal subsemigroup of G with sup property. Then f $^{-1}$(μ) is an S-(∈,∈∨q)-fuzzy quasi normal subsemigroup of G.*

*Proof*: Refer [25] and obtain an analogous proof by studying the subgroup P in the S-semigroup G.

**THEOREM 2.3.27:** *S-(∈,∈∨q) fuzzy normality implies S (∈,∈∨q) fuzzy quasi normality.*

*Proof*: Direct as in case of groups.

It is important to note that S-(∈,∈∨q) fuzzy quasi normality need not imply S-(∈,∈∨q) fuzzy normality in case of S-semigroups also.

**DEFINITION 2.3.13:** *Let $\lambda$ be a S-fuzzy subsemigroup of a S-semigroup G. For any $g \in P \subset G$ let $S(\lambda^g) : P \subset G \to I$ defined by $S(\lambda^g(x)) = \lambda(g^{-1}xg)$ for all $x \in P \subset G$.*

*This definition is true and can be defined for every proper subset P in G, P a subgroup of G.*

**THEOREM 2.3.28:** *Let $\lambda$ be a S-fuzzy subsemigroup of the S-semigroup G. Then*

  i.  *For all $g \in P \subset G$, $\lambda^g$ is a S-fuzzy subsemigroup of G.*
  ii. *$\lambda$ is a S- $(\in, \in \vee q)$-fuzzy normal if and only if $\lambda^g(x) \geq M(\lambda(x), 0.5)$ for all $x \in P \subset G$.*

*Proof*:

  i.  For any $x, y \in P \subset G$, $\lambda^g(xy) = \lambda(g^{-1}xyg) = \lambda(g^{-1}xygg^{-1}yg) \geq M(\lambda(g^{-1}xg), \lambda^{-1}(g^{-1}yg), 0.5) = M(\lambda^g(x), \lambda^g(y), 0.5)$. Again $\lambda^g(x^{-1}) = \lambda(g^{-1}x^{-1}g) = \lambda(g^{-1}xg)^{-1} \geq M(\lambda(g^{-1}xg), 0.5)) = M(\lambda^g(x), 0.5)$. So $\lambda^g$ is a S-fuzzy subsemigroup relative the subgroup $P \subset G$.

  ii. Follows directly by the very definition.

**DEFINITION 2.3.14:** *Let $\lambda$ be a S-fuzzy subsemigroup of the S-semigroup G. Then the Smarandache core (S-core) of $\lambda$ denoted by $S(\lambda_G^P)$ and defined by*

$$S(\lambda_G^P) = \cap \{\lambda^g \mid g \in P \subset G\}.$$

*Thus for a given S-semigroup G we can have as many S-core as the number of proper subset of G which are subgroups of G.*

**THEOREM 2.3.29:** *Let $\lambda$ be a S-fuzzy subsemigroup of G. Then $(\cap \lambda^g)_t = \cap g\lambda_t g^{-1} = S(\lambda_G^P)_t$ for all $t \in (0,1]$, for each subgroup $P \subset G$.*

*Proof*: Let $t \in (0,1]$.

Then $x \in S(\lambda_G^P)_t = (\cap \lambda^g)_t$ if and only if

$$(\cap \lambda^g)(x) \geq t \text{ or } (\cap \lambda^g)(x) + t > 1$$

if and only if inf $\{\lambda(g^{-1}xg) \mid g \in P\} \geq t$ or inf $\{\lambda(g^{-1}xg) \mid g \in G\} + t > 1$ if and only if $(\lambda(g^{-1}xg) \geq t$ or $(\lambda(g^{-1}xg) + t \geq 1$ for all $g \in P$ if and only if $(g^{-1}xg)_t \in \vee q \, \lambda$ if and only if $g^{-1}xg \in \lambda_t$ if and only if $x \in g \lambda_t g^{-1}$ for $g \in G$ if and only if $x \in \underset{g}{\cap} g \lambda_t g^{-1}$.

So $(\cap \lambda^g)_t = \cap g \lambda_t g^{-1}$ for all $t \in (0,1]$.

**THEOREM 2.3.30:** $\lambda_g$ *is an S-* $(\in, \in \vee q)$ *fuzzy normal subgroup of G.*

*Proof*: Matter of routine, hence left for the reader as an exercise.

**DEFINITION 2.3.15:** *A S-fuzzy subsemigroup $\lambda$ is said to be Smarandache core free (S-core free) if there exists some $\alpha \in (0,1]$ such that $\lambda_P = e_\alpha$, $P \subset G$. P a subgroup of G.*

**DEFINITION 2.3.16:** *Let $\lambda$ be a S-fuzzy subsemigroup of a finite S-semigroup G and let $S_p$ be a p-sylow subgroup of $P \subset G$. Let a fuzzy subset $\xi_p$ of $P \subset G$ such that $\xi_p \leq \lambda$ be defined by*

$$\xi_p(x) = \begin{cases} \lambda(x) & \text{if } x \in S_p \\ 0 & \text{if } x \notin S_p. \end{cases}$$

**THEOREM 2.3.31:** $\xi_p$ *is a S-fuzzy subsemigroup of G.*

*Proof*: Follows by the very definition.

**DEFINITION 2.3.17:** $\xi_P$ *is called S-*$(\in, \in)$ *fuzzy p-sylow subgroup of $\lambda$.*

**DEFINITION 2.3.18:** *Let $\lambda$ be a S-fuzzy subsemigroup of a finite S-semigroup G and $S_p$ be a p-sylow subgroup of $S_p \subset G$. Let a fuzzy subset $\mu_p$ of G be defined as follows:*

$$\forall x \in S_p \; \mu_p(x) = \begin{cases} \lambda(x) & \text{if } \lambda(x) < 0.5 \\ \geq 0.9 & \text{if } \lambda(x) \geq 0.5 \end{cases}$$

$\mu_p(x) = 0$ if $x \notin S_p$.

**THEOREM 2.3.32:** $\mu_p$ *is a S-fuzzy subsemigroup of a S-semigroup of G.*

*Proof*: As in case of groups.

**DEFINITION 2.3.19:** $\mu_p$ *is called a Smarandache $(\in, \in \vee q)$ fuzzy p-sylow subgroup (S-*$(\in, \in \vee q)$ *fuzzy p-sylow subgroup) of $\lambda$.*

**THEOREM 2.3.33:** *If $\mu_p$ is a S-*$(\in, \in \vee q)$ *fuzzy p-sylow subgroup of $\lambda$ then $(\mu_p)_t = \lambda_t \cap S_p$ for all $t \in (0, 1]$.*

*Proof*: Analogous to group. Refer [24].

**THEOREM 2.3.34:** *Let $\lambda$ be an S-* $(\in, \in \vee q)$ *fuzzy quasi normal subsemigroup of a finite S-semigroup G. If $\lambda$ is S-core-free then every S-*$(\in, \in \vee q)$-*fuzzy p-sylow subgroup of $\lambda$ is S-*$(\in, \in \vee q)$ *fuzzy quasi normal.*

*Proof*: Left as an exercise for the reader to prove.

**DEFINITION 2.3.20:** *Let $\lambda$ be a S-fuzzy subsemigroup of a S-semigroup G. $\lambda$ is said to be a S-($\in$, $\in$) fuzzy maximal. [respectively S-($\in,\in \vee q$)-fuzzy maximal] if $\lambda$ is not constant for any other S-fuzzy subsemigroup $\mu$ of $P \subset G$, whenever $\lambda \leq \mu$ either $[\in - \lambda_{0.5}] = [\in - \mu_{0.5}]$ or $\mu = \chi_G$ [resp. either $[(\in \vee q) - \lambda_{0.5}] = [(\in \vee q) - \mu_{0.5}]$ or $\mu = \chi_G$] where $[\in - \lambda_t]$ denotes the level subset or denotes $[(\in \vee q)$-level subset].*

*It is important to mention here following the definition of [32] it is a $[\in - \lambda_t]$-level subset and $[(\in \vee q) - \lambda_t ]$ denotes $(\in \vee q)$-level subset following definition of [24], $(\in \vee q)$-level subset].*

It is left for the reader to verify $[(\in \vee q) - \lambda] = [\in -\lambda_t]$ and $[(\in \vee q) - \lambda_{0.5}] = [\in -\lambda_{0.5}]$ for all $t \in (0, 0.5)$.

**THEOREM 2.3.35:** *Let $\lambda$ be a S-fuzzy subsemigroup of a S-semigroup G. If $\lambda$ is S-($\in$, $\in$)-fuzzy maximal then $\lambda(x) \geq 0.5$ for all $x \in P \subset G$ ( P subgroup of G).*

*Proof*: If possible let there exist $a \in P \subset G$ such that $\lambda(a) < 0.5$. Define $\mu : P \subset G \to I$ by

$$\mu(x) = \begin{cases} \lambda(x) & \lambda(x) \geq 0.5 \\ 0.5 & \lambda(x) < 0.5 \end{cases}$$

Then clearly $\mu$ is a S-fuzzy subsemigroup of G such that $\lambda < \mu$ but $[\in - \lambda_{0.5}] \neq [\in - \mu_{0.5}]$. Clearly $\mu \neq \chi_G$. So $\mu$ is not S- $(\in,\in)$ maximal a contradiction.

Therefore $\lambda(x) \geq 0.5$ for all $x \in P \subseteq G$.

**DEFINITION 2.3.21:** *A S-fuzzy subsemigroup $\lambda$ of a S-semigroup G is said to be q-fuzzy maximal if for any other S-fuzzy subsemigroup $\mu$ of G whenever $\lambda \leq \mu$ either $\mu = \chi_a$ or*

$$\overline{P}_\lambda = \{x \in P \subset G \mid \lambda(x) \geq \lambda(e)\} = \overline{P}_\mu = \{x \in P \subset G \mid \mu(x) \geq \mu(e)\}.$$

Following the definitions of [5].

**THEOREM 2.3.36:** *If $\lambda$ is a S-fuzzy maximal subsemigroup of the S-semigroup G then $\lambda$ is q-fuzzy maximal.*

*Proof*: The result follows from the fact that

$$\{x \in P \subset G \mid \lambda(x) \geq \lambda(e)\} = \{x \in P \subset G \mid \lambda(x) = \lambda(e)\}$$

for any S-fuzzy subsemigroup $\lambda$ of G.

**THEOREM 2.3.37:** *If $\lambda$ is a S-q fuzzy maximal subsemigroup of the S-semigroup G such that $\lambda(x) < 0.5$ for some $x \in P \subset G$, then Im $\lambda = 2$.*

*Proof*: Let $\lambda$ be a S-q-fuzzy maximal subsemigroup of the S-semigroup G; let Im $\lambda >$ 2. Let t, r, $\lambda$(e) $\in$ Im $\lambda$ and t $\neq$ r $\neq \lambda$ (e). If r < 0.5 then define a fuzzy subset $\mu$ of P $\subset$ G as follows:

$$\mu(xs) = \begin{cases} \lambda(e) & \text{if } x = e \\ t(> \lambda(e)) & \text{if } x \neq e. \end{cases}$$

Then $\mu$ is a S-fuzzy subsemigroup of G and $\mu \geq \lambda$. Clearly $\overline{G}_\lambda \neq \overline{G}_\mu$ and $\mu \neq \chi_G$ which contradicts that $\lambda$ is q-fuzzy maximal.

Hence r $\geq$ 0.5. Similarly t > 0.5. Again t, r $\geq$ 0.5 implies that $\lambda$ (x) $\geq$ 0.5 for all x $\in$ P $\subset$ G which is a contradiction. Moreover Im $\lambda$ = 1 only when $\lambda$ is constant. Therefore Im $\lambda$ = 2.

**THEOREM 2.3.38**: *If $\lambda$ is a q-fuzzy maximal or S-($\in, \in \vee$ q) fuzzy maximal subsemigroup of P $\subset$ G. Then $\overline{P}_\lambda$ and $\lambda_{0.5}$ are maximal subgroups of P $\subset$ G.*

*Proof*: Left for the reader as an exercise.

**THEOREM 2.3.39:** *Let $\lambda$ be a S-($\in, \in \vee q$)-fuzzy quasi normal subsemigroup of P $\subset$ G with sup property. If $\lambda$ is q-fuzzy maximal or S-($\in, \in \vee q$)-fuzzy maximal, then $\lambda$ is S-($\in, \in \vee q$)-fuzzy normal.*

*Proof*: Refer [25] for analogous proof.

Now we define the concept of Smarandache primary fuzzy subsemigroup of a S-semigroup G. We for the sake of completeness recall some of the definitions.

**DEFINITION [14]**: *Let A be a fuzzy subset of a group G. Then A is called a fuzzy subgroup of G under a t-norm T (T-fuzzy subgroup) if and only if for all x, y $\in$ G.*

    i.    *A(xy) $\geq$ T (A(x), A(y)).*
    ii.    *A(e) = 1 where e is the identity of G.*
    iii.    *A(x) = A($x^{-1}$).*

*A is called a Min-fuzzy subgroup if A satisfies conditions (i) and (iii) only by replacing T with Min Let A be a fuzzy subset of the group G. Then the subset {x $\in$ G / A(x) $\geq$ t }, t $\in$ [0, 1] of G is called a t-level subset of G under A and denoted by $A_t$.*

**DEFINITION 2.3.22:** *Let A be a fuzzy subset of the S-semigroup G. Then A is called a Smarandache fuzzy subsemigroup of P $\subset$ G under a t-norm T (S-T-fuzzy subsemigroup) if and only if for all x, y $\in$ P $\subset$ G.*

    i.    *A (xy) $\geq$ T (A(x), A(y)).*
    ii.    *A(e) = 1 where e is the identity element of P.*
    iii.    *A(x) = A($x^{-1}$) for all x $\in$ P.*

*A is called a S-Min-fuzzy subsemigroup if A satisfies conditions (i) and (ii) only by replacing T by Min.*

*Let A be a S-fuzzy subsemigroup of $P \subset G$ is called a t-level subset of $P \subset G$ under A and is denoted by $S(A_t)$.*

**THEOREM 2.3.40:** *Let A be a S-min fuzzy subsemigroup of the S-semigroup G. Then every S-t-level subset $S(A_t)$ of $P \subset G$, $t \in [0, A(e)]$ is a S-subsemigroup of G.*

*Proof:* As in case of subgroups.

**THEOREM 2.3.41:** *Let A be a S-fuzzy subset of S-subsemigroup G such that every t-level subset $S(A_t)$ of $P \subset G$, $t \in Im (A)$ and $A(e) = 1$. Then A is a S-Min-fuzzy subsemigroup of $P \subset G$.*

*Proof:* Similar to results in groups.

**THEOREM 2.3.42:** *Let $f : G \to H$ be a S-semigroup homomorphism and A be a S-T-fuzzy subsemigroup of G. Then f(A) is a S-T-fuzzy subsemigroup of H.*

*Proof:* As in case of groups, hence left for the reader to prove.

**DEFINITION 2.3.23:** *Let A be a S-T-fuzzy subsemigroup of a S-semigroup $P \subset G$, xA (Ax) be left (right) fuzzy coset of A in $P \subset G$ such that $xA(g) = A(x^{-1} g) A(xg) = A(gx^{-1}))$, $g \in P \subset G$. Then A is said to be a Smarandache normal fuzzy subsemigroup (S-normal fuzzy subsemigroup) of G if and only if xA = Ax for all $x \in P \subset G$ and hence $P/A = \{x A \mid x \in P\}$ is a group (if A is normal) with the operation $xAyA = xyA$, $x, y \in P \subset G$.*

**Notation**: Let A be a S-T-fuzzy subsemigroup of G, $x \in P \subset G$ $r, n, \in N$, $r > n$, such that $A(x^r) = A(x^n) = 1$. Since $r > n_1$ then there exists $q_1, q_2 \in N$ such that $r = q_1 n_1 + n_2$, $0 \leq n_2 < n_1$. Therefore, we have two cases either $n_2 = 0$ or $n_2 \neq 0$. If $n_2 = 0$ then $n_1 \mid r$ i.e. $(n_1, r) = n_1$. If $n_2 \neq 0$ we get $A\left(x^{n_2}\right) = A\left(x^{r-q_1 \mid n_1}\right) \geq T(A(x^r), A(x^{q/n})) = 1$. Again, since $n_2 < n_1$, then there exists $q_2, n_3 \in N$ such that $n_1 = q_2 n_2 + n_3$, $0 \leq n_3 < n_2$. Also here we have two cases either $n_3 = 0$ or $n_3 \neq 0$. If $n_3 = 0$ we get $n_2 \mid r$ and $(n_2, r) = n_2$. If $n_3 \neq 0$, we get

$$A\left(x^{n_1}\right) = A\left(x^{n_1 - q_2 n_2}\right) \geq T\left(A(x^{n_1}), A\left(x^{q_2 n_2}\right)\right) = 1$$

and so on. By the division algorithm there exists $q_{i+1}, n_{i+2} \in N$ such that $n_i = q_{i+1} n_{i+1} + n_{i+2}$, $n_{i+2} = 0$ and hence $A\left(x^{n_1}\right) = A\left(x^{q_{i+1} n_{i+1}}\right) = 1$ that is $n_{i+1} \mid r$.

The following theorem is evident from the above definition.

**THEOREM 2.3.43**: *Let A be a S-T-fuzzy subsemigroup of the S-semigroup G; $x \in P \subset G$, $r, n_1 \in N$, $r > n$ such that $A(x^r) = A (x^n) = 1$. Then there exists $d \in N$ such that $A(x^d) = 1$ where $d = (r, n)$.*

**DEFINITION 2.3.24:** *Let A be a T-fuzzy subsemigroup of a S-semigroup of G and $x \in P \subset G$. Then the least positive integer n satisfying the condition $A(x^n) = 1$ is called the fuzzy order of x with respect to A and we use the notation $A\ o\ (x) = n$. If n does not exist we say x is of infinite fuzzy order with respect to A and write $A\ o\ (x) = \infty$.*

**THEOREM 2.3.44:** *Let A be a S-T-fuzzy subsemigroup of a S-semigroup of $P \subset G$, $x \in P \subset G$. Then*

  i. *If $A(x^r) = 1$ then $A\ o\ (x)\ /\ r$*
  ii. *If $A\ o\ (x) < \infty$ then $A\ o\ (x)\ /\ O(x)$.*

*(Here $A\ o(x) = O(x)$ iff $A_1 = \{e\}$ such that $x \in P \subset G$ and $O(x)$ is the order x).*

*Proof*: As in case of T-fuzzy subgroups.

**THEOREM 2.3.45:** *Let A be a S-T-fuzzy subgroup of $P \subset G$, G a S-semigroup and $x \in P \subset G$ such that $A\ o\ (x) < \infty$. Then $A\ o\ (x) = A\ o\ (x^{-1})$.*

*Proof*: Let $A\ o\ (x) = n$ and $A\ o\ (x^{-1}) = m$. Then $A((x^{-1})^n) = A((x^n)^{-1}) = A(x) = 1$. From the above theorem we get $m\ /\ n$. On the other hand, $A(x^m) = A((x^m)^{-1}) = A((x^{-1})^m) = 1$ and hence $n\ /\ m$. Therefore $n = m$.

**THEOREM 2.3.46:** *Let A be a S-T-fuzzy subsemigroup an S-commutative semigroup G. Then the set of all elements in $P \subset G$ (P a proper subset of G which is an abelian group of G) whose fuzzy order with respect to A are finite is a subgroup of P.*

*Proof*: Direct as in case of T-fuzzy subgroup of an abelian group.

**DEFINITION 2.3.25:** *Let A be a S-T fuzzy subsemigroup of an S abelian semigroup G. Then the subgroup of all elements of P in G whose fuzzy orders with respect to A are finite is called the torsion part of P with respect to A and is denoted by $A^P$.*

**THEOREM 2.3.47:** *If A is a S-fuzzy subsemigroup of G, G a S-semigroup such that for $x \in P$, $A\ o\ (x) < \infty$, $n \in N$ then $A\ o\ (x^n) = [A\ o\ (x)]\ /\ [(n, A\ o\ (x))]$.*

*Proof*: As in case of crisp case.

**THEOREM 2.3.48:** *Let A be a S-T fuzzy subsemigroup of G and $x \in P \subseteq G$ such that $A\ o\ (x) = n_1 n_2 \ldots n_r$, $n_i \in N$, $i = 1, 2, \ldots, n_r$, $(n_i, n_j) = 1$ if $i \neq j$. Then there exists $x_1, x_2, \ldots, x_r \in P$ such that $A\ o\ (x_i) = n_i$ for every $i = 1, 2, \ldots, r$; $x = x_1 x_2 \ldots x_r$, $x_i x_j = x_j x_i$ for all i, j.*

*Proof*: Let $q_i = \Pi\{n_j; i \neq j\} > 1$ since the greatest common divisor of $q_1, q_2, \ldots, q_r$ is 1 then there exists $a_i \in Z$ such that $a_1 q_1 + a_2 q_2 + \ldots + a_r q_r = 1$. Let $b_i = (x)^{a_i q_i}$

$$Ao(b_i) = \frac{\pi n_j}{(\pi n_j, a_i q_i)} = \frac{n_i}{(n_i, a_i)} = n_i$$

If $p \neq 1$ is a prime number such that $p|n_1$ and $p|a_i$ then $p \mid a_1 q_1 + \ldots + a_r q_r$. Therefore, $p \mid 1$ which is a contradiction. Put $b_i = x_i$, $i = 1, 2, \ldots, r$. Therefore $x_1 x_2 \ldots x_r = (x)^{a_1 q_1} (x)^{a_2 q_2} \ldots (x)^{a_r q_r} = (x)^{a_1 q_1 + \ldots + a_r q_2}$. It is clear $x_i x_j = x_j x_i$. Hence the claim.

**THEOREM 2.3.49:** *Let $f: G \to H$ be a S-homomorphism of S-semigroups and A be a S-T-fuzzy subsemigroup of G, $x \in P \subset G$. Then $f(A)of(x)) \mid A \circ (x)$.*

*Proof:* Let $Ao(x) = n$ and $(f(A)) \circ (f(x)) = m$. Then $(f(A)((f(x))^n)) = f(A)(f(x^n)) = \sup\{A(y) \mid y \in f^{-1}(f(x^n)) = \sup \{A(y) \mid y \in f^{-1}(f(x^n))\} = 1$; that is because one of the elements y is $x^n$ and consequently $A(x^n) = 1$. From earlier results we get m/n. Hence the result.

**DEFINITION 2.3.26:** *Let A be a S-T-fuzzy subsemigroup of $P \subset G$, G a S-semigroup. Then the least common multiple of the fuzzy order of the elements of P in G with respect to A is called the Smarandache order (S-order) of the S-fuzzy subsemigroup A and is denoted by $|SA|_F$. If it does not exists $|SA|_F = \infty$.*

**THEOREM 2.3.50:** *Let A be a S-T-fuzzy subsemigroup of a finite S-semigroup G; then $|SA|_F \mid |P|$ but $|SA|_F \nmid |G|$ ($P \subset G$ relative to which A is defined).*

*Proof:* Straightforward; to prove $|SA|_F \mid |P|$, the reader is requested to construct an example for $|SA|_F \mid |G|$ by considering symmetric semigroups S(n).

Thus in general analogous of Lagranges theorem in case of S-fuzzy subsemigroups is not always possible.

Recall: Let A be a T fuzzy subgroup of G, p a prime. Then A is called primary fuzzy subgroup of G if for every x in G there exists a natural number r such that $A \circ (x) = p^r$. Now we proceed on to define S-primary fuzzy subsemigroups of a S-semigroup G.

**DEFINITION 2.3.27:** *Let A be a S-T-fuzzy subsemigroup of $P \subset G$, G a S-semigroup and p a prime. Then A is called Smarandache primary fuzzy subsemigroup (S-primary fuzzy subsemigroup) of $P \subset G$ if for every $x \in P$ there exists a natural number r such that $Ao(x) = p^r$.*

*Clearly every S-T-fuzzy subsemigroup of $P \subset G$ is S-primary if and only if there exists $n \in N$ such that $|SA|_F = p^n$, p a prime.*

**THEOREM 2.3.51:** *If A is a fuzzy subset of a S-semigroup G and every $A_t$, $t \in Im(A)$ is a S-primary subsemigroup of $P \subset G$ then A is a S-primary Min fuzzy subsemigroup of $P \subset G$.*

*Proof*: Let $x \in P \subset G$ such that $A(x) = t$, $t \in \text{Im}(A)$. Then $x \in A_t$ and hence there exists $m \in Z$ such that $O(x) = p^m$ where p is a prime. Since $A \circ (x) \mid O(x)$ there exists $0 \leq s \leq m$ such that $Ao(x) = p^s$ i.e. $A\left(x^{p^s}\right) = 1$.

Therefore A is a S-primary fuzzy subsemigroup of $P \subset G$.

**THEOREM 2.3.52:** *Every S-T-fuzzy subsemigroup of a S-primary semigroup is also S-primary.*

*Proof*: Direct from the earlier theorem and the fact if the semigroup G is S-primary then every proper subset $P \subset G$ where P is a subgroup and is S-primary.

**THEOREM 2.3.53:** *Let A be a S-primary Min- fuzzy subsemigroup of the S-semigroup G. Then every t-level subset $A_t$, $t \in \text{Im}(A)$, of G is S-primary if $A_\alpha = \{e\}$ where $\alpha = A(e)$, e is the identity element of $P \subset G$ relative to which A is defined.*

*Proof*: Let $x \in P \subset G$, $A(x) = t$, $A \circ (x) = P^r$, where p is prime $r \in N$. Since $A_{A(e)} = \{e\}$, then $O(x) = A \circ (x) = p^r$. Hence $A_t$ is S-primary.

**THEOREM 2.3.54:** *If A is a S-normal T-fuzzy subsemigroup of a S-primary semigroup G then P/ A is a S-primary semigroup (A is defined relative to $P \subset G$).*

*Proof*: As in case of T-fuzzy subgroup of a primary group G.

**THEOREM 2.3.55:** *Let $f : G \to H$ be a S-semigroup homomorphism and A be S-primary T-fuzzy sub semigroup of $P \subset G$ and $A_{A(e)} = \{e\}$. Then f(A) is S-primary T-fuzzy subsemigroup of $T \subset P$.*

*Proof*: Left for the reader as an exercise.

**THEOREM 2.3.56:** *Let A be a S-T-fuzzy subsemigroup of a S-cyclic primary semigroup G such that $A_1 = \{e\}$ for $x, y \in P \subset G$ Then*

    i.    *If $Ao(x) = Ao(y)$ then $A(x) = A(y)$.*
    ii.    *If $Ao(x) > Ao(y)$ then $A(y) > A(x)$.*

*Proof*: Since $A_1 = \{e\}$. Therefore $Ao(x) = O(x)$ and hence the proof follows by earlier results.

**THEOREM 2.3.57:** *Let A be a S-T-fuzzy subgroup of an S-abelian semigroup G. Then the set of all elements whose fuzzy order is a power of a prime p is a subgroup of $P \subset G$.*

*Proof*: Let $x, y \in P \subset G$ such that $A \circ (x) = p^r$ and $A \circ (y) = p^s$, $r, s \in N$, $s \geq r$. Then $A\left((xy^{-1})^{p^s}\right) = A\left(x^{p^s}(y^{-1})^{p^s}\right) \geq T\left(A(x^{p^s}), A((y^{-1})^{p^s})\right) = 1$ that is $A \circ (xy^{-1}) = P^s$ or $A \circ (xy^{-1}) \mid p^s$. If $A \circ (xy^{-1}) \mid p^s$ then there exists $t \in N$, $r \leq t \leq s$ such that $Ao(xy^{-1}) = p^t$. Hence $xy^{-1}$ has a fuzzy order with respect to A as a power of p. Hence the claim.

**DEFINITION 2.3.28:** *Let A be a S-T-fuzzy subsemigroup of G. A is said to be Smarandache abelian (S-abelian) if for all a, b ∈ P ⊂ G, A(a) > 0, A(b) > 0 implies ab = ba.*

**DEFINITION 2.3.29:** *A a S-Min-fuzzy subsemigroup of P ⊂ G is called the Smarandache p-sylow fuzzy subgroup (S-p-Sylow fuzzy subgroup) of G if it has the following properties:*

   i.  *Card (ImA) ≤ 2.*
   ii. *$A_{A(e)}$ is a sylow p-subgroup of P (P ⊂ G).*

**THEOREM 2.3.58:** *If A is a S-p-sylow normal Min fuzzy subsemigroup of G then A is not a S-primary fuzzy subgroup of P ⊂ G with order as a power of p.*

*Proof*: A matter of routine as in case of subgroups.

Now we proceed on to define a Smarandache p-component of a S-semigroup G.

**DEFINITION 2.3.30:** *If A is a S-T-fuzzy subsemigroup of P ⊂ G then the set of all elements in P ⊂ G of a fuzzy order as a power of a prime p with respect to A is called Smarandache p-component (S-p-compotent) of P ⊂ G with respect to A and it is denoted by $P_{P(A)}$.*

**THEOREM 2.3.59:** *Let A be a S-T-fuzzy subsemigroup of P ⊂ G, G a S-semigroup. Then $T_A$ is the sum of $P_{P(A)}$, p ∈ V where V is the set of all prime and it is a direct sum if $A_1 = \{e\}$.*

*Proof*: Refer [17] and prove it as in case of S-semigroups.

In the following section we define Smarandache fuzzy bisemigroup.

## 2.4 Smarandache fuzzy bisemigroups

In this section we introduce the concept of Smarandache fuzzy bisemigroups. The study of Smarandache bisemigroups and Smarandache bigroups is very recent [135]. Several interesting properties of these concepts are studied in this section.

**DEFINITION 2.4.1:** *Let (G, •, ∗) be a non empty set such that $G = G_1 \cup G_2$ where $G_1$ and $G_2$ are proper subsets of G, (G, •, ∗) is called a Smarandache bigroup (S-bigroup) if the following conditions are true:*

   i.  *(G, •) is a group.*
   ii. *(G,∗) is a S-semigroup.*

Now we have got several results about them which we just recall without proof.

**THEOREM [135]:** *Let G be a S-bigroup then G need not be a bigroup.*

*Proof*: Construct an example to prove.

**THEOREM [135]:** *If $G = G_1 \cup G_2$ is a S-bigroup then G contains a bigroup.*

*Proof*: Left for the reader to refer [135].

**DEFINITION [135]:** *Let $G = G_1 \cup G_2$ be a S-bigroup, we say G is a Smarandache commutative bigroup (S-commutative bigroup) if $G_1$ is a commutative group and every proper subset S of $G_2$ which is a group is a commutative group. If both $G_1$ and $G_2$ happens to be commutative trivially G becomes a S-commutative bigroup.*

*We say G is S-weakly commutative if the S-semigroup $G_2$ has atleast one proper subset which is commutative.*

The following theorem is straight forward and hence left for the reader to prove.

**THEOREM [135]:** *Let $G = G_1 \cup G_2$ be a S-commutative bigroup then G is a S-weakly commutative bigroup and not conversely.*

*Proof*: Left for the reader as an exercise.

**DEFINITION 2.4.2:** *Let $G = G_1 \cup G_2$ be a S-bigroup. We say G is a Smarandache cyclic bigroup (S-cyclic bigroup) if $G_1$ is a cyclic group and $G_2$ is a S-cyclic semigroup. We call G a S-weakly cyclic bigroup if every subgroups of $G_1$ is cyclic and $G_2$ is a S-weakly cyclic bigroup.*

The following theorem is straight forward hence left for the reader to prove.

**THEOREM 2.4.1:** *Let $G = G_1 \cup G_2$ be a S-cyclic bigroup then G is a S-weakly cyclic bigroup.*

Now we proceed on to define Smarandache hyper bigroup.

**DEFINITION 2.4.3:** *Let $G = G_1 \cup G_2$ be a S-bigroup. If A be a proper subset of $G_2$ which is the largest group contained in $G_2$ and $G_1$ has no subgroups, then we call $P = G_1 \cup A$ the Smarandache hyper bigroup (S-hyper bigroup). We say G is Smarandache simple bigroup (S-simple bigroup) if G has no S-hyper bigroup.*

We define Smarandache double cosets in bigroups.

**DEFINITION 2.4.4:** *Let $G = G_1 \cup G_2$ be a S-bigroup. Let $A = A_1 \cup A_2$ and $B = B_1 \cup B_2$ be any S-sub-bigroups of G. We define the Smarandache double coset (S-double coset) as $AxB = \{A_1 x_1 B_1 \cup A_2 x_2 B_2 / A_1 x_1 B_1$ provided $x = x_1 \in G_1$ if $x \neq x_1 \notin G_1$ then we just take $A_1 B_1$ if $x = x_2 \in G_2$ then take $A_2 x_2 B_2$ otherwise take $A_2 B_2$ for every $x \in G\} : AxB = \{A_1 B_1 \cup A_2 x B_2$ if $x \in G_2\}$ or $= \{A_1 x B_1 \cup A_2 B_2$ if $x \in G_1\}$ or $= \{A_1 x B_1 \cup A_2 x B_2 / $ if $x \in G_1 \cup G_2\}$.*

We just recall the definition of Smarandache normal bigroup.

**DEFINITION 2.4.5:** *Let $G = G_1 \cup G_2$ be a S-bigroup we call a S-sub-bigroup A of G to be a Smarandache normal sub-bigroup (S-normal sub-bigroup) of G where $A = A_1 \cup A_2$ if $xA_1 \subset A_1$ and $A_1x \subset A_1$ if $x \in G_1$ and $xA_2 \subset A_2$ and $A_2 x \subset A_2$ or $xA_2 = \{0\}$ or $Ax_2 = \{0\}$ if 0 is an element in $G_2$ and $x \in G_2$ for all $x \in G_1 \cup G_2$. We say G is Smarandache pseudo simple bigroup (S-pseudo simple bigroup) if G has no S-normal sub-bigroup.*

Now we proceed on to define Smarandache maximal S-bigroup.

**DEFINITION 2.4.6:** *Let $G = G_1 \cup G_2$ be a S-bigroup. We say a proper subset $A = A_1 \cup A_2$ to be Smarandache maximal S-bigroup (S-maximal S-bigroup) of G if $A_1$ is the largest normal subgroup in $G_1$ and $A_2$ is the largest proper subset of $G_2$ which is the subgroup of $G_2$.*

Now we proceed on to define Smarandache fuzzy bigroup.

**DEFINITION 2.4.7:** *Let $(G = G_1 \cup G_2, +, \bullet)$ be a S-bigroup. Then $\mu : G \to [0, 1]$ is said to be a Smarandache fuzzy sub-bigroup (S-fuzzy sub-bigroup) of G if there exists two fuzzy subsets $\mu_1$ (of $G_1$) and $\mu_2$ (of $G_2$) such that*

  i.   *$(\mu_1, +)$ is a fuzzy subgroup of $(G_1, +)$*
  ii.  *$(\mu_2, \bullet)$ is a S-fuzzy subsemigroup of $(G_2, \bullet)$ and*
  iii. *$\mu = \mu_1 \cup \mu_2$.*

We illustrate this by an example.

**Example 2.4.1:** Consider the S-bigroup $G = \{0, \pm1, \pm2, \cdots, \frac{1}{2}, \frac{1}{2^2}, \frac{1}{2^4}, \cdots\}$, $G_1 = \{0, \pm1, \pm2, \cdots\}$ is a group under usual multiplication and $G_2 = \{1, \frac{1}{2}, \frac{1}{2^2}, , \cdots\}$ is a S-semigroup under product. Define $\mu : G \to [0, 1]$ by

$$\mu(x) = \begin{cases} \frac{1}{3} & \text{if } x \in \left\{\frac{1}{2}, \frac{1}{2^2} \cdots\right\} \\ 1 & \text{if } x \in \{0, \pm 2, \cdots\} \\ \frac{1}{2} & \text{if } x \in \{\pm 1, \pm 3, \cdots\} \end{cases}$$

It is easily verified that $\mu_1 : G_1 \to [0, 1]$ is given by

$$\mu_1(x) = \begin{cases} 1 & \text{if } x \in \{0, \pm2, \pm4, \cdots\} \\ \frac{1}{2} & \text{if } x \in \{\pm1, \pm3, \cdots\} \end{cases}$$

and $\mu_2 : G_2 \to [0, 1]$ given by

$$\mu_2(x) = \begin{cases} \dfrac{1}{2} & \text{if } x = 1 \\ \dfrac{1}{3} & \text{if } x \in \left\{ \dfrac{1}{2}, \dfrac{1}{2^2}, \ldots \right\} \end{cases}$$

Thus $\mu = \mu_1 \cup \mu_2$ is a S-fuzzy sub-bigroup of $G = G_1 \cup G_2$.

**THEOREM 2.4.2:** *Every t-level subset of a S-fuzzy sub-bigroup $\mu$ of a S-bigroup G need not in general be a S-sub bigroup of the S-bigroup G.*

*Proof*: The proof is by an example. Take $G = \{0, 1, -1, 2, 4, 6, \ldots\}$ to be a S-bigroup. $G = G_1 \cup G_2 = \{-1, 1\} \cup \{0, 2, 4, 6, \ldots\}$. $(G_1, \bullet)$ is a group and $(G_2, +)$ is a S-semigroup.

$\mu : G \to [0, 1]$

$$\mu(x) = \begin{cases} \dfrac{1}{2} & \text{if } x = 1 \text{ or } -1 \\ \dfrac{1}{4} & \text{if } x = 0 \\ \dfrac{1}{8} & \text{if } x = 2, 4, \cdots \end{cases}$$

Clearly $(\mu, +, \bullet)$ is a S-fuzzy sub bigroup of $(G, +, \bullet)$. Now consider the level subset $G_\mu^{\frac{1}{2}}$ of the fuzzy sub-bigroup $\mu$

$$G_\mu^{\frac{1}{2}} = \left\{ x \in G \,\bigg|\, \mu(x) \geq \dfrac{1}{2} \right\} = \{-1, 1\}.$$

It is easily verified that $\{-1,1\}$ is not a S-sub-bigroup of the S-bigroup $(G, +, \bullet)$. Hence the t-level subset $G_\mu^t \left[\text{for } t = \dfrac{1}{2}\right]$ of the S-fuzzy sub-bigroup $\mu$ is not a sub-bigroup of the bigroup $(G, +, \bullet)$.

Now we proceed on to define the bilevel subset of the Smarandache fuzzy sub-bigroup of a S-bigroup G.

**DEFINITION 2.4.8:** *Let $(G = G_1 \cup G_2, +, \bullet)$ be a bigroup and $\mu = (\mu_1 \cup \mu_2)$ be a S-fuzzy sub-bigroup of G. The bilevel subset of the S-fuzzy sub-bigroup $\mu$ of the S-bigroup G is defined as $G_\mu^t = G_{1\mu_1}^t \cup G_{2\mu_2}^t$ for every $t \in [0, \min \{\mu_1(e_1), \mu_2(e_2)\}]$ where $e_1$ denotes the identity elements of the group $(G_1, +)$ and $e_2$ denotes the identity element of $P_2 \subset G_2$ where $P_2$ is the proper subset of $G_2$ which is the subgroup of $G_2$ relative to which $\mu_2$ is defined.*

It is left as an exercise for the reader to prove.

**THEOREM 2.4.3:** *Every bilevel subset of a S-fuzzy sub-bigroup μ of a S-bigroup G is a S-sub-bigroup of the S-bigroup G.*

Now we define the concept of Smarandache fuzzy sub-bigroup of a S-semigroup G.

**DEFINITION 2.4.9:** *A fuzzy subset μ of a S-semigroup G is said to be a Smarandache fuzzy sub-bisemigroup (S-fuzzy sub-bisemigroup) of the S-semigroup G if there exists two S-fuzzy subsemigroups $\mu_1$ and $\mu_2$ of μ ($\mu_1 \neq \mu$ and $\mu_2 \neq \mu$) such that $\mu = \mu_1 \cup \mu_2$. Here by the term S-fuzzy subsemigroup of μ we mean that λ is a S-fuzzy subsemigroup of the S-semigroup G and $\lambda \subseteq \mu$ (where μ is also a fuzzy subgroup of G).*

*Example 2.4.2:* Consider the additive S-semigroup G; G = $\{0, \pm 1, \pm 2, \ldots, n/2^t \mid n = 0, 1, 2, \ldots$ and $t = 1, 2, \ldots\}$. Clearly G is a S-semigroup.

Define $\mu : G \to [0, 1]$

$$\mu(x) = \begin{cases} 1 & \text{if } x \in \{0, \pm 2, \pm 4, \ldots\} \\ 0.5 & \text{if } x \in \{\pm 1, \pm 3, \ldots\} \\ 0.25 & \text{if } x \in \{n/2^t \mid \begin{matrix} n=1,2,\ldots \\ t=1,2,\ldots \end{matrix}\} \end{cases}$$

It is easily verified μ is a S-fuzzy subsemigroup of the S-semigroup G for we have two S-fuzzy subsemigroups $\mu_1$ and $\mu_2$ of μ ($\mu_1 \neq \mu$ and $\mu_2 \neq \mu$) such that $\mu = \mu_1 \cup \mu_2$ where $\mu_1$ and $\mu_2$ are as given below:

$$\mu_1(x) = \begin{cases} 1 & \text{if } x \in \{0, \pm 2, \pm 4, \ldots\} \\ 0.5 & \text{if } x \in \{\pm 1, \pm 3, \ldots\} \\ 0.25 & \text{if } x \in \{n/2^t \mid \begin{matrix} n=1,2,\ldots \\ t=1,2,\ldots \end{matrix}\} \end{cases}$$

$$\mu_2(x) = \begin{cases} 0.75 & \text{if } x \in \{0, \pm 2, \pm 4, \ldots\} \\ 0.5 & \text{if } x \in \{\pm 1, \pm 3, \ldots\} \\ 0 & \text{if } x \in \{n/2^t \mid \begin{matrix} n=1,2,\ldots \\ t=1,2,\ldots \end{matrix}\} \end{cases}$$

In view of these definitions and examples we have the following theorem.

**THEOREM 2.4.4:** *Let $\mu = \mu_1 \cup \mu_2$ be a S-fuzzy subsemigroup of a S-semigroup G where $\mu_1$ and $\mu_2$ are S-fuzzy subsemigroups of G. For $t \in [0, \min \{\mu_1 (e_1), \mu_2 (e_2)\}]$ where $e_1$ and $e_2$ are identity elements of $P_1$ and $P_2$ ($P_1 \subset G_1$, $P_2 \subset G_2$ subgroups of $G_1$* *and $G_2$ relative to which $\mu_1$ and $\mu_2$ are defined. The level subset $G_\mu^t$ of $\mu$ can be represented as the union of two subsemigroups of the S-semigroup G, that is $G_\mu^t$.*

*Proof*: Left as an exercise for the reader.

**THEOREM 2.4.5:** *Let $\mu$ be a Smarandache fuzzy subsemigroup of a S-semigroup G with $3 \leq o(Im(\mu)) < \infty$ then there exists two S-fuzzy subsemigroups $\mu_1$ and $\mu_2$ of $\mu$, ($\mu_1 \neq \mu$, $\mu_2 \neq \mu$) such that $\mu = \mu_1 \cup \mu_2$.*

*Proof*: Let $\mu$ be a S-fuzzy subsemigroup of the S-semigroup G. Suppose Im ($\mu$) = $\{a_1, a_2,\ldots, a_n\}$ where $3 \leq n < \infty$ and $a_1 > a_2 > \ldots > a_n$. Choose $b_1, b_2 \in [0, 1]$ be such that $a_1 > b_1 > a_2 > b_2 > a_3 > b_3 > \ldots > a_n$ and define $\mu_1, \mu_2: P \to [0, 1]$, P a subgroup of the S-semigroup G relative to which $\mu$ is defined.

$$\mu_1(x) = \begin{cases} a_1 & \text{if } x \in \mu_{a_1} \\ b_2 & \text{if } x \in \mu_{a_2} \setminus \mu_{a_1} \\ \mu(x) & \text{otherwise} \end{cases}$$

and

$$\mu_2(x) = \begin{cases} b_1 & \text{if } x \in \mu_{a_1} \\ a_2 & \text{if } x \in \mu_{a_2} \setminus \mu_{a_1} \\ \mu(x) & \text{otherwise} \end{cases}$$

Thus it can be easily verified that both $\mu_1$ and $\mu_2$ are S-fuzzy subsemigroups of $\mu$. Further $\mu_1 \neq \mu$, $\mu_2 \neq \mu$ and $\mu = \mu_1 \cup \mu_2$.

It is very evident that the condition $3 \leq o(Im\, \mu) < \infty$ cannot be dropped.

The following theorem is left as an exercise for the reader to prove.

**THEOREM 2.4.6:** *Every S-fuzzy sub bisemigroup of a S-semigroup G is a S-fuzzy subsemigroup of the S-semigroup G and not conversely.*

**THEOREM 2.4.7:** *Let $\mu$ be a fuzzy subset of a S-semigroup G with $3 \leq o(Im\, \mu) \leq \infty$ then $\mu$ is a S-fuzzy sub-bisemigroup of the S-semigroup G if and only if $\mu$ is a S-fuzzy sub-bisemigroup of the S-semigroup G.*

*Proof*: The $\mu$ once defined, it is relative to a particular fixed proper subset P of G where P is a subgroup under the operations of G. The rest of the result is identical with the fuzzy sub-bigroup. Almost all the definitions and properties derived in case of S-fuzzy semigroups can be easily defined and derived for S-fuzzy bigroups and S-fuzzy bisemigroups. We just give very few of the related results S-fuzzy bigroups and S-fuzzy bisubsemigroup notions in this section.

**DEFINITION 2.4.10:** *Let $(G, +, \bullet)$ be a S-bigroup. $G = G_1 \cup G_2$ where $(G_1, +)$ is a group and $(G_2, \bullet)$ is a S-semigroup. $A: G \to [0, 1]$ where $A = A_1 \cup A_2$ with $A_1$ a divisible fuzzy subgroups of $G_1$ and $A_2$ is a divisible S-fuzzy subsemigroup of $G_2$ (i.e. $A_2 : G_2 \to [0, 1]$ is such that $A_2$ restricted to the subgroup $P_2$ of $G_2$ is a fuzzy subgroup $G_2$). Then we call $A = A_1 \cup A_2: G \to [0, 1]$ is a Smarandache divisible fuzzy subgroup (S-divisible fuzzy sub-bigroup) of G.*

Several interesting results in this direction can be obtained.

*Let G be a S-bigroup; $(G = G_1 \cup G_2)$ $A = A_1 \cup A_2: G \to [0, 1]$ be a S-fuzzy sub-bisemigroup of G; i.e., $A_1: G_1 \Rightarrow [0, 1]$ is a fuzzy subgroup of $G_1$ and $A_2 : G_2 \to [0, 1]$ is a S-fuzzy subsemigroup of $G_2$ related to S-sub-bigroup P of G.*

$S(F(A)) = \{B \mid B$ *is a S-fuzzy sub-bisemigroup of G such that* $B \subseteq A$ *i.e.* $B_1 \subseteq A_1$ *and* $B_2 \subseteq A_2\}$

**DEFINITION 2.4.11:** *Let $(G = G_1 \cup G_2, +, \bullet)$ be a S-bigroup. Let $A : G \to [0, 1]$ be a S-fuzzy bisubsemigroup of G. Let $B \in S(F(A))$. Then B is said to be Smarandache pure (S-pure) in A if and only if for all fuzzy singletons $x_t \subset B$ with $t > 0$, for all $n \in N$ for all $y_t \in A$, $n(y_t) = x_t$ implies that there exists $b_t \subset B$ such that $n(b_t) = x_t$.*

Next we define Smarandache $(\in, \in \vee q)$- fuzzy sub-bigroup of a S-bigroup G.

**DEFINITION 2.4.12:** *Let $(G, +, \bullet)$ be a S-bigroup. A fuzzy subset $\lambda = \lambda_1 \cup \lambda_2 : G = G_1 \cup G_2 \to [0, 1]$ is said to be a Smarandache $(\in, \in \vee q)$ fuzzy sub-bigroup (S-$(\in, \in \vee q)$ fuzzy sub-bigroup) of G if*

    i. *For any $x, y \in G_1$ and $t, r \in (0,1]$ $\lambda_1 : G_1 \to [0, 1]$, $(G_1, +)$ a group; $x_t, y_r \in \lambda_1 \Rightarrow (xy)_{M(t,r)} \in \vee q \lambda_1$ and*
    ii. *$x_t \in \lambda_1 \Rightarrow (x^{-1})_t \in \vee q \lambda_1$*
    iii. *$\lambda_2 : G_2 \to [0, 1]$, $\lambda_2$ defined on $P_2 \subset G_2$ ($P_2$ subgroup of the S-semigroup $(G_2, \bullet))$ is such that for any $x, y \in P_2$ and $t, r \in (0,1]$, $x_t, y_r \in \lambda_2 \Rightarrow (xy)_{M(t,r)} \in \vee q \lambda_2$ and*
    iv. *$x_t \in \lambda_2 \Rightarrow (x^{-1})_t \in \vee q \lambda_2$.*

*The conditions in the above definition are respectively equivalent to*

    i.    *$\lambda_1 (xy) \geq M(\lambda_1 (x), \lambda_1 (y), 0.5)$ for $x, y \in G_1$.*
    ii.    *$\lambda_2 (xy) \geq M(\lambda_2 (x), \lambda_2 (y), 0.5)$ for $x, y \in P_2 \subset G_2$ ($P_2$ subgroup of $G_2$).*
    iii.    *$\lambda_1(x^{-1}) \geq M(\lambda_1 (x), 0.5)$ for all $x \in G_1$.*
    iv.    *$\lambda_2(x^{-1}) \geq M(\lambda_2 (x), 0.5)$ for all $x \in P_2 \subset G_2$.*

*Further we have for any S-$(\in, \in \vee q)$ - fuzzy sub-bisemigroup $\lambda = \lambda_1 \cup \lambda_2$ of $G = G_1 \cup G_2$ with $\lambda_1 (x) \geq 0.5$ for some $x \in G_1$, $\lambda_1 (e_1) \geq 0.5$, e is the identity element of $G_1$ and $\lambda_2 (x) \geq 0.5$ for some $x \in P_2 \subset G_2$, $\lambda_2 (e_2) \geq 0.5$; $e_2$ identity element of $P_2$ relative to which $\lambda_2$ is defined.*

*Further if $\lambda$ is a S-fuzzy subsemigroup of G then $\lambda$ is an S-($\in, \in \vee q$)-fuzzy subsemigroup.*

**DEFINITION 2.4.13:** *Let $(G = G_1 \cup G_2, +, \bullet)$ be a S-bigroup. $\lambda$ a S-fuzzy sub-bisemigroup of G. $\lambda = \lambda_1 \cup \lambda_2 : G = G_1 \cup G_2 \to [0, 1]$ such that $\lambda_1 : G_1 \to [0, 1]$ is fuzzy subgroup and $\lambda_2 : G_2 \to [0, 1]$ is a S-fuzzy subsemigroup relative to $P \subset G_2$; $\lambda$ is said to be a*

i.      *S- ($\in, \in$)-fuzzy normal if for all $x, y \in G_1$ and $t \in (0,1]$ $x_t \in \lambda_1 \Rightarrow (y^{-1} xy)_t \in \lambda_1$. For $x, y \in P_2 \subset G_2$ and $t \in (0,1]$, $x_t \in \lambda_2 \Rightarrow (y^{-1} xy)_t \in \lambda_2$.*

ii.      *S-($\in, \in \vee q$)-fuzzy normal (or simply fuzzy normal).*

     a      *if for any $x, y \in G_1$ and $t \in (0,1]$, $x_t \in \lambda_1 \Rightarrow (y^{-1} xy)_t \in \vee q \lambda_1$.*
     b      *if for any $x, y \in G_2$ and $t \in (0,1]$, $x_t \in \lambda_2 \Rightarrow (y^{-1} xy)_t \in \vee q \lambda_2$.*

*Several results derived for groups using S-($\in, \in \vee q$)-fuzzy subgroup and S-($\in \vee q$)-fuzzy level subgroups can be derived in an analogous way with appropriate modifications in case of S-bigroups.*

The reader is requested to refer [24] and [135] for more innovative notions and concepts. Now we just define the concept of Smarandache ($\in, \in$)-fuzzy left (resp. right) coset of a S-bigroup G.

**DEFINITION 2.4.14:** *Let G be a S-bigroup i.e. $G = G_1 \cup G_2$, $(G_1, +)$ is a group and $(G_2, \bullet)$ is a S-semigroup $\lambda: G = G_1 \cup G_2 \to [0, 1]$ be the S-fuzzy subsemibigroup of G; i.e. $\lambda = \lambda_1 \cup \lambda_2$ where $\lambda_1: G_1 \to [0, 1]$ is a S-fuzzy semigroup and $\lambda_2 : G_2 \to [0, 1]$ is a S-fuzzy subsemigroup relative to a subgroup $P_2$ of $G_2$. For $x \in G = G_1 \cup G_2$ we define the following if $x_1 \in G_1$ then $\lambda^l_{x_1}$ ( resp $\lambda^r_{x_1}$ ) : $G_1 \to I$ is defined by*

$$\lambda^l_x(g) = \lambda(gx_1^{-1}) \text{ and } \lambda^r_{x_1}(g) = \lambda(x_1^{-1}g)$$

*for all $g \in G_1$ and for $x_2 \in P_2 \subset G_2$ we have $\lambda^l_{x_2}$ ( resp $\lambda^r_{x_2}$ ) : $P_2 \subset G \to I$ defined by*

$$\lambda^l_{x_2}(g) = \lambda(gx_2^{-1}) \text{ and } \lambda^r_{x_2}(g) = \lambda(x_2^{-1}g)$$

*for all $g \in P_2$ $\lambda^l_x = \lambda^l_{x_1} \cup \lambda^l_{x_2}$ ( resp. $\lambda^r_x = \lambda^r_{x_1} \cup \lambda^r_{x_2}$ ) is called the Smarandache ($\in, \in$)-fuzzy left (resp. right) bicoset (S-($\in, \in$)-fuzzy left (resp. right) bicoset) of G determined by $x = x_1 \cup x_2$ and $\lambda$.*

**DEFINITION 2.4.15:** *Let $\lambda = \lambda_1 \cup \lambda_2$ be a S-fuzzy sub-bisemigroup of a S-bigroup G. For $x \in G = G_1 \cup G_2$, $\hat{\lambda}_{x_1}$ (resp. $\tilde{\lambda}_{x_1}$) : $G_1 \to [0, 1]$ is defined by* $\hat{\lambda}_{x_1}(g) = M(\lambda(gx_1^{-1}, 0.5))$ [resp. $\tilde{\lambda}_{x_1}(g) = M(\lambda(x_1^{-1}g, 0.5))$] for all $x_1 \in G_1$ and $\hat{\lambda}_{x_2}$ (resp. $\tilde{\lambda}_{x_2}$): $P_2 \to [0,1]$ is defined by $\tilde{\lambda}_{x_2}(g) = M(\lambda(gx_2^{-1}, 0.5))$] [resp. $\tilde{\lambda}_{x_2}(g) = M(\lambda(x_2^{-1}g, 0.5))$] for all $g \in P_2 \subset G_2$.

Then $\hat{\lambda}$ (resp $\tilde{\lambda}$) is called the Smarandache ($\in, \in \vee q$)-fuzzy left (resp. right) bicoset (S-($\in, \in \vee q$)-fuzzy left (resp. right) bicoset) of G determined by $x = x_1 \cup x_2$ and $\lambda$.

The following theorem is straightforward by the very definition; hence left for the reader to prove.

**THEOREM 2.4.8:** *Let $\lambda$ be a S-fuzzy sub-bisemigroup of a S-bigroup G. Then $\lambda$ is a S-($\in, \in \vee q$) fuzzy normal if and only if $\tilde{\lambda}_{x_1} = \hat{\lambda}_{x_1}$ for all $x_1 \in G_1$ and $\tilde{\lambda}_{x_2} = \hat{\lambda}_{x_2}$ for all $x_2 \in P_2 \subseteq G_2$ where $\lambda$ is the restriction of $P_2$ in $G_2$.*

**DEFINITION 2.4.16:** *Let G be a S-bigroup. H be a S-sub-bigroup of G. H of G is said to be Smarandache quasi normal (S-quasi normal) if for every S-sub-bigroup K of G; $H \bullet K = K \bullet H$.*

**DEFINITION 2.4.17:** *Let G be a S-bigroup. A S-fuzzy sub-bisemigroup $\xi$ of $G = G_1 \cup G_2$ is called S-fuzzy quasi normal sub-bisemigroup if $\xi o \eta = \eta o \xi$ for every S-fuzzy sub-bisemi group $\eta$ of G.*

*Unless otherwise mentioned G will denote a S-bigroup with $e_1$ identity of $G_1$ and $e_2$ the identity of $P_2 \subset G_2$ and $\lambda$ an S-($\in, \in \vee q$) fuzzy sub-bisemigroup of G. By a S-fuzzy sub-bisemigroup of G we shall mean an S-($\in, \in \vee q$) fuzzy sub-bisemigroup of G.*

*For a S-fuzzy subsemigroup $\lambda$ of $G = G_1 \cup G_2$ for all $x, y \in G_1$ or $x, y \in P_2 \subset G_2$ such that $\lambda_i(x) \geq 0.5$, $i = 1,2$ and $\lambda_i(y) < 0.5$, $i = 1,2$.*

**DEFINITION 2.4.18:**

  i. *$\lambda$ is said to be S-($\in, \in$) fuzzy quasi normal if for any S-fuzzy sub-bisemigroup $\mu$ of $G = G_1 \cup G_2$ for all $x \in G_1$, $t \in (0,1]$ and if $x \in P_2 \subseteq G_2$,*

  $$z_t \in (\lambda o \mu) \Leftrightarrow z_t \in (\mu \ o \ \lambda).$$

  ii. *$\lambda$ is said to be a S-($\in, \in \vee q$) fuzzy quasi normal if for any S-fuzzy sub-bisemigroup $\mu$ of $G = G_1 \cup G_2$ and for all $z \in G_1$, $t \in (0, 1]$, $z \in P_2 \subset G_2$, $t \in (0, 1]$, $z_t \in (\lambda \ o \ \mu)$ implies $z_t \in \vee q \ (\mu \ o \ \lambda)$ and if $z \in P_2 \subset G_2$.*

  *($P_2$ is the subgroup relative to which $G_2$ is defined) and $z_t \in (\mu \ o \ \lambda)$ implies $z_t \in \vee q \ (\lambda \ o \ \mu)$.*

The following theorem is left an exercise for the reader to prove.

**THEOREM 2.4.9:** *$\lambda$ is S-($\in,\in$)-fuzzy quasi normal if and only if for any S-fuzzy sub-bisemigroup $\mu$ of G and for all $z \in G_1$ (or $z \in P_2 \subset G_2$), $(\lambda \circ \mu)(z) = (\mu \circ \lambda)(z)$.*

**THEOREM 2.4.10:** *$\lambda$ is a S-($\in, \in \vee q$)-fuzzy quasi normal if and only if for any S-fuzzy sub-bisemigroup $\mu$ of G and for all $z \in G_1$ and for all $z \in P_2 \subset G_2$ ($P_2$ a subgroup of $G_2$) we have $(\mu \circ \lambda)(z) = (\lambda \circ \mu)(z)$ if $(\lambda \circ \mu)(z) < 0.5$ and $(\mu \circ \lambda)(z) \geq 0.5$ if $(\lambda \circ \mu)(z) \geq 0.5$.*

*Proof*: As in case of fuzzy subgroups.

**THEOREM 2.4.11:** *If $\lambda$ is a S-fuzzy quasi normal subsemigroup of the S-bigroup $G = G_1 \cup G_2$, then $\lambda$ is an S-($\in,\in \vee q$)-fuzzy quasi normal sub-bisemigroup of G.*

*Proof*: Using definitions and the fact that $\lambda = \lambda_1 \cup \lambda_2$ and $\lambda_1$ is fuzzy quasi normal subgroup of $G_1$ and $\lambda_2$ is a fuzzy quasi normal subgroup of $P_2$ contained in $G_2$.

**THEOREM 2.4.12:** *Let H be any non-empty subset of $G = G_1 \cup G_2$ i.e. $H = H_1 \cup H_2$ ($H_1 \subset G_1$ and $H_2 \subset G_2$). Then H is a S-quasi normal sub-bisemigroup of G if and only if $X_H$ (the characteristic function H) is an S-($\in,\in \vee q$)-fuzzy quasi normal sub-bisemigroup of G.*

*Proof*: As in case of fuzzy subgroups.

**THEOREM 2.4.13:** *Let $\lambda$ be a S-($\in, \in \vee q$)-fuzzy quasi normal sub-bisemigroup of G with the sup property. Then S-($\in \vee q$)-level subset $\lambda_t = \{x_1 \in G_1 | (x_1)_t \in \vee q \lambda_1\} \cup \{x_2 \in P_2 \subset G_2 | (x_2)_t \in \vee q \lambda_2\}$ is a S-quasi normal sub-bi semigroup of $G = G_1 \cup G_2$ for all $t \in (0,1]$.*

*Proof*: Left as an exercise for the reader to prove.

**THEOREM 2.4.14:** *Let $\lambda$ be a fuzzy subset of a S-bigroup G with sup property and for all $t \in (0,1]$, S-($\in \vee q$)-level subset $\lambda_t$ be a S-quasi normal sub-bisemigroup of G. Then $\lambda$ is an S-($\in, \in \vee q$)-fuzzy quasi normal sub-bisemigroup of S-bigroup G.*

*Proof*: As in case of fuzzy subgroup G the result can be easily proved for S-bigroup G.

**THEOREM 2.4.15:** *Let $f : G \to H$ be S-bigroup homomorphism from the S-bigroup G to S-bigroup H. Then*

   i. *$f(\lambda \circ \theta) = f(\lambda) \circ f(\theta)$ where $\lambda$ and $\theta$ are S-fuzzy sub-bisemigroups of G.*
   ii. *$f(f^{-1}(\theta)) = \theta$ where $\theta$ is a S-fuzzy sub-bigroup of $f(G)$.*

*Proof*: Follows by the very definitions.

**THEOREM 2.4.16:** *Let $f : G \to H$ be a S-bigroup homomorphism of the S-bigroups G and H If $\lambda$ is an S-($\in, \in \vee q$) fuzzy quasi normal sub-bisemigroups of G then $f(\lambda)$ is an S-($\in,\in \vee q$) fuzzy quasi normal sub-bisemigroup of $f(G)$.*

*Proof*: Follows by the very definitions.

**THEOREM 2.4.17:** *Let $f : G \to H$ be a S-bigroup homomorphism and $\mu$ be an $(\in, \in \vee q)$-fuzzy quasi normal sub-bisemigroup of H with sup-property. Then $f^{-1}(\mu)$ is a S-$(\in, \in \vee q)$-fuzzy quasi normal sub-bisemigroup of G.*

*Proof*: Left for the reader to prove, using the fact $f : G = G_1 \cup G_2 \to H = H_1 \cup H_2$ where $f = f_1 \cup f_2$ where $f_1 : G_1 \to H_1$ is a group homomorphism $f_2 : G_2 \to H_2$ is a S-semigroup homomorphism i.e. $f_2 : P_2 \subset G_2 \to Q_2 \subset H_2$ where $P_2$ and $Q_2$ are subgroups of $G_2$ and $H_2$ respectively.

**DEFINITION 2.4.19:** *Let $\lambda$ be a S-fuzzy sub-bisemigroup of a S-bigroup $G = G_1 \cup G_2$. For any $g_1 \in G_1$ let $\lambda_1^{g_1} : G_1 \to [0,1]$ and for $g_2 \in P_2 \subset G_2$ ($\lambda$ defined related to this $P_2$) $\lambda_2^{g_2} : G_2 \to [0,1]$ defined by $\lambda_1^{g_1}(x) = \lambda_1(g_1^{-1} x g_1)$ for all $x \in G_1$ and $\lambda_2^{g_2}(x) = \lambda_2(g_2^{-1} x g_2)$ for all $x \in P_2 \subset G_2$. So $\lambda^g = \lambda_1^{g_1} \cup \lambda_2^{g_2}$.*

**THEOREM 2.4.18:** *Let $\lambda$ be a S-fuzzy sub-bisemigroup of the S-bigroup $G = G_1 \cup G_2$.*

*Then*

i. For all $g_1 \in G_1$, $\lambda_1^{g_1}$ is a fuzzy subgroup of $G_1$ and $\lambda_2^{g_2}$ is a S-fuzzy subsemigroup of $G_2$.

ii. $\lambda$ is S- $(\in, \in \vee q)$ fuzzy normal if and only if $\lambda_1^{g_1}(x_1) \geq M(\lambda_1(x), 0.5)$ for all $x_1 \in G_1$, $g_2 \in G_2$, $\lambda_2^{g_2}(x_2) \geq M(\lambda_2(x), 0.5)$ for all $x_2 \in P_2$ and $g_2 \in G_2$, where $P_2$ is the subgroup associated with the S-semigroup $G_2$ and $\lambda = \lambda_1 \cup \lambda_2$.

*Proof:* Direct using definitions.

**DEFINITION 2.4.20:** *Let $\lambda$ be a S-fuzzy sub-bisemigroup of the S-bigroup $G = G_1 \cup G_2$. Then the Smarandache bicore (S-bicore) of $\lambda = \lambda_1 \cup \lambda_2$ denoted by $\lambda_{G_1} \cup \lambda_{G_2}$ is defined by*

$$\lambda_{G_1} = \left( \cap \left\{ \lambda_1^{g_1} \mid g_1 \in G_1 \right\} \right)$$

*and*

$$\lambda_{G_2} = \left( \cap \left\{ \lambda_2^{g_2} \mid g_2 \in P_2 \subset G_2 \right\} \right).$$

**DEFINITION 2.4.21:** *Let $\lambda$ be a S-fuzzy sub-bisemigroup of a S-bigroup $G = G_1 \cup G_2$. $\lambda$ is said to be Smarandache bicore free (S-bicore free) if there exists $\alpha \in (0, 1]$ such that $\lambda_{G_1} = e_\alpha$, e identity element of $G_1$ and $\lambda_{G_2} = e_\alpha^1$. $e_\alpha^1$ -identity element of $P_2 \subset G_2$ ($P_2$ – a subgroup of $G_2$).*

**DEFINITION 2.4.22:** *Let $\lambda$ be a S-fuzzy sub-bisemigroup of the S-bigroup G. $\lambda$ is said to be a Smarandache $(\in,\in)$-fuzzy maximal (S-$(\in,\in)$-fuzzy maximal) if $\lambda$ is not constant and for any other S-fuzzy sub-bisemigroup $\mu$ of the S-bigroup G whenever $\lambda \le \mu$ either $[\in - \lambda_{0.5}] = [\in - \mu_{0.5}]$ (resp. either $[(\in \vee q) - \lambda_{0.5}] = [(\in \vee q) - \mu_{0.5}]$ where $[\in - \lambda_t]$- denotes the level subset and $[(\in \vee q) - \lambda_t]$ denotes the $(\in \vee q)$-level subset.*

**DEFINITION 2.4.23:** *A S-fuzzy sub-bisemigroup $\lambda$ of a S-bigroup G is said to be a S-q-fuzzy maximal if for any other S-fuzzy sub-bisemigroup $\mu$ of G whenever $\lambda \le \mu$ either $\mu = \chi_G$ of $\overline{G}_\lambda = \{x \in G_1 \cup G_2 \mid \lambda(x_1) \ge \lambda(e_1)$ if $x_1 \in G_1$, if $x_2 \in P_2 \subset G_2$ then $\lambda(x_2) \ge \lambda(e_2)$, $e_2$ identity element of the subgroup $P_2 \subset G_2\}$.*

$$\overline{G}_\mu = \{x \in G_1 \cup G_2 \mid \mu(x_1) \ge \mu(e_1),\ x_1 \in G_1\ \text{and}\ \mu(x_2) > \mu(e_2)\ \text{if}\ x_2 \in P_2 \subset G_2\}$$

Several interesting results in this direction can be obtained, some of the results are proposed as problems in section 2.5.

**DEFINITION 2.4.24:** *Let A be a fuzzy subset of a S-bigroup G. Then A is called a S-fuzzy sub-bisemigroup of G and a t-norm T (T-fuzzy bigroup) if and only if for all $x, y \in G$.*

i. $A_1(xy) \ge T(A_1(x), A_1(x), A_1(y))$ $x, y \in G_1$.
ii. $A_1(e) = 1$ where e is the identity of $G_1$.
iii. $A_1(x) = A_1(x^{-1})$.
iv. $A_2(xy) \ge T(A_2(x), A_2(y))$, $x, y \in P_2 \subset G_2$ ($P_2$ a subgroup of $G_2$ relative to which $A = A_1 \cup A_2$ is defined).
v. $A_2(e) = 1$ where e is the identity of $P_2 \subset G_2$.
vi. $A_1(x) = A_1(x^{-1})$ for $x \in P_2 \subset G_2$, A is called the Smarandache Min fuzzy sub-bisemigroup (S-Min fuzzy sub-bisemigroup) if A satisfies conditions (i) to (vi) only by replacing T with Min.

**DEFINITION 2.4.25:** *Let A be a fuzzy subset of a S-bigroup G. Then the subset $\{x \in G \mid A(x) \ge t\}$ of $G = G_1 \cup G_2$ is called the Smarandache t-level subset (S-t-level subset) of G under A denoted by $A_t$ (i.e. $A = A_1 \cup A_2$ and $A_t = (A_1)_t \cup (A_2)_t$ with $A_1(x_1) \ge t$ if $x_1 \in G_1$ and $A_2(x_2) \ge t$ if $x_2 \in G_2$).*

**THEOREM 2.4.19:** *Let A be a S-Min fuzzy sub-bisemigroup of the S-bigroup G. Then every t-level subset $A_t$ of G, $t \in [0, A(e)] = [0, A_1(e_1) \cup A_2(e_2)]$ where $e_1$ is the identity of the group $G_1$ and $e_2$ is the identity of the subgroup $P_2$ in $G_2$ relative to which A is defined is a S-sub-bigroup of G.*

*Proof*: Follows by the very definitions; hence left for the reader as an exercise to prove.

**THEOREM 2.4.20:** *Let A be a fuzzy subset of a S-bigroup G such that every t-level subset $A_t$ of G, $t \in \text{Im}(A)$ and $A(e_1) = 1 = A(e_2)$ where $e_1$ is the identity element of $G_1$* and $e_2$ is the identity element of $P_2$ in $G_2$ relative to which $A$ is defined. Then $A$ is a S-Min fuzzy sub-bisemigroup of $G$.

*Proof*: As in case of groups. Hence the reader is requested to obtain the proof.

**THEOREM 2.4.21**: *Let $f : G \to H$ be a S-bigroup homomorphism and $A$ be a S-T-fuzzy sub-bisemigroup of $G$. Then $f(A)$ is a S-T-fuzzy sub-bisemigroup of $H$.*

**DEFINITION 2.4.26**: *Let $A$ be a S-T-fuzzy sub-bisemigroup of the S-bigroup $G$. $x_1 A_1$ ($A_1 x_1$) be the left (right) fuzzy coset of $A_1$ in $G_1$ such that $x_1(A_1(g_1)) = A_1(x_1^{-1} g_1)$, $(A_1(x_1(g_1))) = A_1(g_1 x_1^{-1})$, $g_1 \in G_1$, $x_1 \in G_1$ and for $x_2 \subset P_2 \subset G_2$ where $A$ is defined relative to this subgroup $P_2 \subset G_2$ we have $x_2 A_2$ ($A_2 x_2$) be left (right) fuzzy coset of $A_2$ in $G_2$ such that*

$$x_2 (A_2(g_2)) = A_2(x_2^{-1} g_2)$$

$$(A_2(x_2(g_2))) = A_2(g_2 x_2^{-1}), \; g_2 \in G_2,$$

*then $A = A_1 \cup A_2$ is said to be a Smarandache normal fuzzy sub-bisemigroup (S-normal fuzzy sub-bisemigroup) of $G$ if and only if $x_1 A_1 = A_1 x_1$ ($x_2 A_2 = A_2 x_2$) for all $x_1 \in G_1$ and for all $x_2 \in P_2 \subset G$.*

Hence

$$\frac{G}{A} = \left( \frac{G_1}{A_1} \cup \frac{G_2}{A_2} \right) = \left( x_1 A_1 \cup x_2 A_2 \mid x_1 \in G_1 \text{ and } x_2 \in P_2 \subset G_2 \right)$$

is a S-bigroup with operation.

$$xAyA = x_1 y_1 A_1 \cup x_2 y_2 A_2, \; x_1, y_1 \in G_1 \text{ and } x_2, y_2 \in P_2 \subset G_2.$$

**THEOREM 2.4.22**: *Let $A$ be a S-T-fuzzy sub-bisemigroup of the S-bigroup $G = G_1 \cup G_2$. For $x \in G = G_1 \cup G_2$ there exists $r_1, n_1 \in N$, $r_1 > n_1$ such that $A_1(x_1^{r_1}) = A_1(x_1^{n_1}) = 1$, $x_1 \in G_1$ and $r_2, n_2 \in N$, $r_2 > n_2$ such that $A_2(x_2^{r_2}) = A_2(x_2^{n_2}) = 1$ for $x_2 \in P_2 \subset G_2$ relative to which $A$ is defined and there exists $d_1, d_2 \in N$ such that $A_1(x^{d_1}) = 1$ and $A_2(x^{d_2}) = 1$ where $d_1 = (r_1, n_1)$ and $d_2 = (r_2, n_2)$.*

*Proof*: Follows directly using definitions.

**DEFINITION 2.4.27**: *Let $A$ be a S-T-fuzzy sub-bisemigroup of a S-bigroup $G$ and $x \in G = G_1 \cup G_2$. Then the least positive integers $n_1$, $n_2$ which satisfies the condition $A_1(x_1^{n_1}) = 1$ and $A_2(x_2^{n_2}) = 1$, $x_1 \in G_1$ and $x_2 \in P_2 \subset G_2$ ($P_2$ subgroup relative to which $A = A_1 \cup A_2$ is defined) is called the Smarandache fuzzy order (S-fuzzy order) of $x \in G_1 \cup G_2$ relative to $A$ and use the notation $A_1 o (x_1) = n_1$ and $A_2 o (x_2) = n_2$.*

*If $n_1$ and $n_2$ does not exist we say $x_1$ and $x_2$ is of Smarandache infinite fuzzy order (S-infinite fuzzy order) with respect to $A_1$ and $A_2$ and write $A_1$ o $(x_1) = \infty$ and $A_2$ o $(x_2) = \infty$. From the above definition we see $A_1$ o $(x_1) = O(x_1)$ if and only if $A_1 = \{e_1\}$ such that $x_1 \in G_1$ and $O(x_1)$ is the order of $x_1$ and $A_2$ o $(x_2) = O(x_2)$ if and only if $A_2 = \{e_2\}$ such that $x_2 \in P_2$ ($P_2 \subset G_2$ is a subgroup of $G_2$) and $O(x_2)$ is the order of $x_2$.*

*Thus we see the S-fuzzy order is a complete generalization of fuzzy order which is the generalization of the usual order of the S-bigroup. But that is not true for all $x_1 \in G_1$ and $x_2 \in P_2 \subset G_2$ as $G_1 (x_1) = 1$; hence $G_1$ o $(x_1) = 1$ which is clearly impossible. $P_2(x_2) = 1$, hence $P_2$ o $(x_2) = 1$ which is once again impossible.*

Now we proceed on to define Smarandache torsion part of a S-bigroup G.

**DEFINITION 2.4.28:** *Let $G = G_1 \cup G_2$ be a S-bigroup. A be a S-T-fuzzy sub-bisemigroup of the bigroup G. Then the S-sub-bigroup of all elements in A whose S-fuzzy orders with respect to $A = A_1 \cup A_2$ are finite is called the Smarandache torsion part (S-torsion part) of G relative to A and is denoted by $S(T_A)$.*

Several interesting results in this direction can be developed. We propose some problems for the reader in section 2.5.

Now we proceed on to define Smarandache order of the S-fuzzy sub-bisemigroup A.

**DEFINITION 2.4.29:** *Let A be a S-T-fuzzy sub-bisemigroup of the S-bigroup G. Then the least common multiple of the S-fuzzy order of the elements of G with respect to A is called the Smarandache order (S-order) of the S-fuzzy sub-bisemigroup A and it is denoted by $|S(A)|_G = \left(|S(A_2)|_G, |S(A_1)|_G\right)$ where $A = A_1 \cup A_2$. If it does not exist then $|S(A)|_G = \infty$.*

*We call a S-bigroup $G = G_1 \cup G_2$ to be Smarandache primary (S-primary) if $G_1$ is a primary group and every proper subset $P_2$ of $G_2$ which is a subgroup of $G_2$ is a primary group. We say the S-bigroup $G = G_1 \cup G_2$ to be Smarandache weakly primary (S-weakly primary) if $G_1$ is primary and atleast a proper subset which is a proper subgroup in $G_2$ is primary.*

**THEOREM 2.4.23:** *Let G be a S-bigroup. If G is a S-primary bigroup then G is a S-weakly primary bigroup.*

*Proof*: Direct by the very definition.

**DEFINITION 2.4.30:** *Let A be a S-T-fuzzy sub-bisemigroup of the S-bigroup $G = G_1 \cup G_2$ and p a prime. Then A is called Smarandache primary fuzzy sub-bisemigroup (S-primary fuzzy subsemigroup) of G if for every $x \in G_1 \cup G_2$ if there exists natural numbers $r_1, r_2 \in N$ such that $A_1$ o $(x_1) = p^{r_1}$ and $A_2$ o $(x_2) = p^{r_2}$, $x_1 \in G_1$ and $x_2 \in P_2 \subset G_2$ (A defined relative to $P_2$).*

**DEFINITION 2.4.31:** *Let A be a S-T-fuzzy sub-bisemigroup of a S-bigroup $G = G_1 \cup G_2$. A is said to be abelian if for all $a_1, b_1 \in G_1, A_1(a_1) > 0, A_1(b_1) > 0$ implies* $a_1 b_1 = b_1 a_1$ for all $a_2, b_2 \in P_2 \subset G_2$ ($P_2$ is the subgroup relative to which $A = A_1 \cup A_2$ is defined) $A_2(a_2) \geq 0$, $A_2(b_2) \geq 0$ implies $a_2 b_2 = b_2 a_2$.

**DEFINITION 2.4.32:** *If A is a S-T fuzzy sub-bisemigroup of a S-bigroup $G = G_1 \cup G_2$, then the set of all elements of G is a S-fuzzy order as a power of prime numbers $(p_1, p_2)$ with respect to $A = A_1 \cup A_2$ is called the Smarandache $(p_1, p_2)$ components of $G = G_1 \cup G_2$ with respect to $A = A_1 \cup A_2$ denoted by*

$$(G)_{p_1(A_1) \cup p_2(A_2)} = (G_1)_{p_1(A_1)} \cup (P_2)_{p_2(A_2)}$$

*($P_2 \subset G_2$ a subgroup relative to which A is defined).*

## 2.5 Problems

In this section we give fifty-four problems on S-fuzzy semigroups and S-fuzzy bisemigroups for the reader to solve.

The problem section of each chapter happens to be an integral part of the book, for the solutions to these problems will throw a lot of light about these new Smarandache fuzzy structures.

**Problem 2.5.1:** If in S(3) σ is a S-fuzzy symmetric semigroup what is o(Im(σ))?

**Problem 2.5.2:** Give an example of a S-fuzzy symmetric semigroup in S(5).

**Problem 2.5.3:** Obtain some interesting results about S-fuzzy symmetric group S(n).

**Problem 2.5.4:** If g is a S-co-fuzzy symmetric subsemigroup of the symmetric semigroup S(3); find o(Im g).

**Problem 2.5.5:** If g is a S-co-fuzzy symmetric subsemigroup of S(n) then find O(Im(g)).

**Problem 2.5.6:** Give an example to show that every S-fuzzy symmetric subsemigroup of a symmetric subsemigroup S(n) need not in general be a S-co-fuzzy symmetric subsemigroup of S(n).

**Problem 2.5.7:** Obtain some interesting results about S-fuzzy normal subsemigroup of the S-semigroup G.

**Problem 2.5.8:** Find conditions so that (A/B, •) is a S-semigroup; or Is (A/B, •) always a S-semigroup? Justify your answer.

**Problem 2.5.9:** Find conditions so that (A/B, •) is a group? Is this possible? If so illustrate by an example.

**Problem 2.5.10:** Give an example of a S-simple semigroup of finite order.

**Problem 2.5.11:** Obtain some interesting results about SF (A).

**Problem 2.5.12:** Characterize those S-semigroups G which are S-torsion fuzzy subsemigroup.

**Problem 2.5.13:** Give an example of a S-p-primary fuzzy subsemigroup of a S-semigroup G.

**Problem 2.5.14:** Characterize those S-semigroups G in which A is S-p-primary if and only if A* is S-p-primary.

**Problem 2.5.15:** Can all S-semigroups have S-maximal p-primary fuzzy subsemigroup?

**Problem 2.5.16:** Give an example of a S-divisible fuzzy semigroup and characterize those S-semigroups G which are S-divisible fuzzy semigroups.

**Problem 2.5.17:** Characterize those S-semigroups G which are S-pure.

**Problem 2.5.18:** Does there exist any relation between S-pure and S-divisible S-semigroups G?

**Problem 2.5.19:** Can all S-semigroups which has a S-fuzzy subsemigroup which is a S-fuzzy weak direct sum a fuzzy weak direct sum?

**Problem 2.5.20:** Give an example of a S-semigroup G which has S-normal fuzzy subgroups.

**Problem 2.5.21:** Illustrate by an example S-tips in a S-semigroup.

**Problem 2.5.22:** Characterize those S-semigroups which has S-multi-tiped fuzzy semigroups.

**Problem 2.5.23:** Do all S-semigroups have S-penultimate subsemigroups?

**Problem 2.5.24:** Characterize those fuzzy sets and the S-semigroup G in which the concept of S-fuzzy semigroup and S-fuzzy ideal coincide.

**Problem 2.5.25:** Characterize those fuzzy subset $\lambda$ of S-semigroup G which are S-($\in$, $\in \vee q$) fuzzy subsemigroup.

**Problem 2.5.26:** Give an example of a S-semigroup G which has S-($\in$, $\in \vee q$) fuzzy subsemigroup.

**Problem 2.5.27:** Illustrate by an example a S-($\in$, $\in$) fuzzy normal subsemigroup of a S-semigroup G.

**Problem 2.5.28:** Characterize those $\lambda : G \to [0,1]$ which are S-($\in$, $\in \vee q$) fuzzy normal.

**Problem 2.5.29:** Give an example of a S-fuzzy quasi normal subsemigroups I and II.

**Problem 2.5.30:** Characterize those S-core S-fuzzy subsemigroup.

**Problem 2.5.31:** Illustrate by examples S-q-fuzzy maximal subsemigroups of a S-semigroup G.

**Problem 2.5.32:** If $\lambda$ is a S-q-fuzzy maximal subsemigroup of a S-semigroup G. Then prove $\overline{G}_\lambda$ and $\lambda_{0.5}$ are maximal subgroups of G.

**Problem 2.5.33:** Prove S-($\in$, $\in \vee q$)-fuzzy normality implies S-($\in$, $\in \vee q$)-fuzzy quasi normality in a S-bigroup $G = G_1 \cup G_2$.

**Problem 2.5.34:** Prove for $\lambda$ a S-fuzzy sub-bisemigroup of G, that $\left(\cap \lambda_1^{g_1}\right)_t = \cap g_1 (\lambda_1)_t g_1^{-1}$ for all $t \in (0,1]$ and $\left(\cap \lambda_2^{g_2}\right)_t = \cap g_2 (\lambda_2)_t g_2^{-1}$, for all $t \in (0, 1]$, for all $g_2 \in P_2 \subset G_2$. $\lambda_G$ is an S- ($\in$, $\in \vee q$) fuzzy normal sub-bisemigroup of the S-bigroup $G = G_1 \cup G_2$.

**Problem 2.5.35:** G be a S-bigroup. Let $\lambda$ be S-($\in$, $\in$) fuzzy maximal. Can we prove if $\lambda$ is not constant and for any other S-fuzzy sub-bisemigroup $\mu$ of G. $\mu = \chi_G$ ?

**Problem 2.5.36:** Under the conditions of the above problem can we prove $t \in (0,0.5)$ $[(\in \vee q) - \lambda_t] = [\in - \lambda_t]$ and $[(\in \vee q) - \lambda_{0.5}] = [\in - \lambda_{0.5}]$?

**Problem 2.5.37:** Prove if $\lambda$ is a S-fuzzy sub-bisemigroup of a S-bigroup G. If $\lambda$ is S- ($\in$, $\in$)-fuzzy maximal then $\lambda(x) \geq 0.5$ for all $x \in G_1$ ($G_1$ a group and $x \in P_2 \subset G_2$, $P_2$ a subgroup of $G_2$).

**Problem 2.5.38:** Let $\lambda$ be a S-fuzzy maximal fuzzy sub-bisemigroup of the S-bigroup G. Is $\lambda$ a S-q-fuzzy maximal?

**Problem 2.5.39:** If $\lambda$ is a S-q-fuzzy maximal sub-bisemigroup of the S-bigroup G such that $\lambda(x) < 0.5$ for some $x \in G_1 \cup G_2$ (if $x \in G_1$ or $x \in P_2 \subset G_2$ relative to $P_2$, $\lambda$ is defined and $P_2$ is a subgroup of $G_2$). Is Im $\lambda = 2$?

**Problem 2.5.40:** If $\lambda$ is a S-q-fuzzy maximal or S-($\in$, $\in \vee q$) fuzzy maximal sub-bisemigroup of the bigroup G. Are $G_\lambda$ and $\lambda_{0.5}$ S-maximal subgroups of G?

**Problem 2.5.41:** Prove if A is a S-T-fuzzy sub-bisemigroup of the S-bigroup G then if $A_1(x_1^{r_1}) = 1$ then $A_1 \circ (x_1) / r_1 \in G$ and if $A_2(x_2^r) = 1$ then $A_2 \circ (x_2) / r_2$, $x_2 \in P_2 \subset G_2$ (A = $A_1 \cup A_2$ defined relative to $G_2$) and if $A_1 \circ (x_1) < \infty$ then $A_1 \circ (x_1) | O(x_1)$, $x_1 \in G_1$ and if $A_2 \circ (x_2) < \infty$ then $A_2 \circ (x_2) | O(x_2)$, $x_2 \in P_2 \subset G_2$.

**Problem 2.5.42:** Prove if A is S-T-fuzzy sub-bisemigroup of a S-bigroup $G = G_1 \cup G_2$ and $x_1 \in G_1$ ($x_2 \in P_2 \subset G_2$) such that $A_1 \circ (x_1) < \infty$ then $A_1 \circ (x_1) = A_1 \circ (x_1^{-1})$. $A_2 \circ (x_2) < \infty$ then $A_2 \circ (x_2) = A_2 \circ (x_2^{-1})$, $x_2^{-1} \in P_2 \subset G_2$.

**Problem 2.5.43:** Prove if A is a S-T-fuzzy sub-bisemigroup of a S-bigroup G where G is a S-commutative bigroup, then the set of all elements in G whose S-fuzzy orders with respect to A = $A_1 \cup A_2$ are finite is a S-sub-bigroup of G.

**Problem 2.5.44:** Prove if A is a S-T-fuzzy sub-bisemigroup of G then $S(T_A)$ is the torsion subgroup of A if and only if $A_1 = (A_1 \cup A_2)_1 = (A_1)_1 \cup (A_2)_1 = \{e_1\} \cup \{e_2\}$, $e_1$ the identity element of $G_1$ and $e_2$ is the identity element of $P_2 \subset G_2$; $P_2$ the subgroup of $G_2$ relative to which A is defined.

**Problem 2.5.45:** Prove if A is a S-T-fuzzy sub-bisemigroup of the S-bigroup $G = G_1 \cup G_2$, $x \in G$ such that $A_1 \circ (x_1) < \infty$, $n_1 \in N$, ($A_2 \circ (x_2)) < \infty$, $x_1 \in G_1$, $x_2 \in P_2 \subset G_2$ then $A_1 o(x_1^{n_1}) = [A_1 o(x_1)] / [n_1, A_1 o(x_1)]$ and $A_2 o(x_2^{n_2}) = [A_2 o(x_2)] / [n_2, A_2 o(x_2)]$, $n_2 \in N_2$. Let $f : G \to H$ be a S-bigroup homomorphism and A be a S-T-fuzzy sub-bisemigroup of G, for $x \in G_1 \cup G_2$; prove $f(A_1) \circ (f(x_1)) | A_1 \circ (x_1)$ if $x_1 \in G_1$ and $f(A_2) \circ (f(x_2)) | A_2 \circ (x_2)$ if $x_2 \in P_2 \subset G_2$. ($P_2$ is the subgroup in $G_2$ relative to which A = $A_1 \cup A_2$ is defined).

**Problem 2.5.46:** Prove if A is a S-T-fuzzy sub-bisemigroup of a finite S-bigroup G then $|S(A_1)|_F \, | \, |G_1|$ and $|S(A_2)|_F \, | \, |P_2|$ where $P_2 \subset G_2$ is the subgroup relative to which A is defined.

**Problem 2.5.47:** Give an example of a S-primary bigroup.

**Problem 2.5.48:** Give an example of a S-weakly primary bigroup. Characterize those bigroups which are

        i. S-primary
       ii. S-weakly primary

**Problem 2.5.49:** Prove if G is a S-primary bigroup and A is a S-normal Min fuzzy sub-bisemigroup of G and [G : A] is finite; then there exists natural numbers $r_1$, $r_2$ in N such that $[G_1; A_1] = p^{r_1}$ and $[P_2; A_2] = p^{r_2}$ where $G = G_1 \cup G_2$ and $A = A_1 \cup A_2$, $A_1 : G_1 \to [0, 1]$ is a Min fuzzy normal subgroup of $G_1$ and $A_2 : P_2 \to [0, 1]$ is a Min fuzzy normal subgroup of $G_2$.

**Problem 2.5.50:** Prove every S-T-fuzzy sub bisemigroup of a S-bigroup G is S-primary if and only if there exists $n \in N$ such that $|SA|_F = p^n$, p is a prime $\left(\text{i.e.} |(SA_1)|_F = p^{n_1} \text{ and } |SA_2|_F = p^{n_2}\right)$.

**Problem 2.5.51:** Is it true that every S-T-fuzzy primary sub bisemigroup is a S-primary bigroup?

**Problem 2.5.52:** If A is a S-primary Min fuzzy sub bisemigroup of a S-bigroup G, are the t-level sub bisemigroup $A_t$ S-Primary? Justify your claim.

**Problem 2.5.53:** If A is a S-normal T-fuzzy sub bisemigroup of a S-primary bigroup G can G/A be a S-primary bigroup?

**Problem 2.5.54:** Is the homomorphic image of a S-primary T-fuzzy sub bisemigroup a primary S-T-fuzzy sub bisemigroup?

**Chapter Three**

# SMARANDACHE FUZZY GROUPOIDS AND THEIR GENERALIZATION

This chapter mainly introduces the notions of Smarandache fuzzy groupoid and Smarandache fuzzy bigroups. Several interesting results in this direction are given. The chapter has four sections. In section one we give some results on Smarandache fuzzy groupoids and around 20 definitions related to Smarandache fuzzy groupoids are given. Groupoids are the generalizations of loops; so, as all loops are groupoids, the Smarandache fuzzy loop study has become inevitable. In section two we define Smarandache fuzzy loop and study these notions. The third section gives in 26 definitions the concepts about Smarandache fuzzy bigroupoids and Smarandache fuzzy biloops. In order to make this book appealing to researchers, in the final section we give 111 problems, which would make a reader master the Smarandache fuzzy groupoid and its generalizations.

## 3.1 Some results on Smarandache Fuzzy groupoids

In this section we for the first time define the notion of Smarandache fuzzy groupoids and illustrate it with examples. Fuzzy groupoids are studied extensively by several researchers like [11, 45, 47, 89, 90, 110].

**DEFINITION 3.1.1**: *A Smarandache groupoid (S-groupoid) G is a groupoid which has a proper subset S, $S \subset G$ such that S under the operations of G is a semigroup.*

*Example 3.1.1:* G = {0, 1, 2, …, 5} be a groupoid given by the following table:

| o | 0 | 1 | 2 | 3 | 4 | 5 |
|---|---|---|---|---|---|---|
| 0 | 0 | 5 | 4 | 3 | 2 | 1 |
| 1 | 3 | 2 | 1 | 0 | 5 | 4 |
| 2 | 0 | 5 | 4 | 3 | 2 | 1 |
| 3 | 3 | 2 | 1 | 0 | 5 | 4 |
| 4 | 0 | 5 | 4 | 3 | 2 | 1 |
| 5 | 3 | 2 | 1 | 0 | 5 | 4 |

Clearly A = {0, 3} is a semigroup of G. Thus G is a S-groupoid.

For several properties about S-groupoids please refer [128].

**DEFINITION 3.1.2:** *Let G be a S-groupoid. A fuzzy subset µ from G to [0, 1] is said to be a Smarandache fuzzy groupoid (S-fuzzy groupoid) G if µ restricted to at least one of the proper subsets $P \subset G$, P a semigroup under the operations of G we have µ (xy) $\geq$ min {µ(x), µ(y)} for x, y $\in$ P, we denote this by $\mu_P$. If for every proper subset $P_i$ of the S-groupoid G. $P_i$ a semigroup, the fuzzy subset µ : G $\to$ [0, 1] is a such that $\mu_{P_i}$* is a fuzzy semigroup, then we call $\mu : G \rightarrow [0, 1]$ a Smarandache strong fuzzy groupoid (S-strong fuzzy groupoid) of G.

**THEOREM 3.1.1:** *Let G be a S-groupoid. Every S-strong fuzzy groupoid $\mu$ of G is a S-fuzzy groupoid $\mu$ of G.*

*Proof*: Follows by the very definitions.

The definition of S-fuzzy groupoids using S-groupoids will be called as the type I or level I S-fuzzy groupoids, by default of notation we just denote it by S-fuzzy groupoids.

We define level II or type II S-fuzzy groupoids in the following.

**DEFINITION 3.1.3**: *Let $\mu : X \rightarrow [0, 1]$ be a fuzzy subgroupoid of the groupoid X with respect to a t-norm T. We say $\mu$ is a Smarandache fuzzy subgroupoid of type II (S-fuzzy subgroupoid of type II) of X if and only if X is a S-groupoid.*

**THEOREM 3.1.2:** *Every S-fuzzy subgroupoid of type II (S-subgroupoid of type II) of X is a fuzzy subgroup of X.*

*Proof*: Direct by the very definition.

**THEOREM 3.1.3:** *Let X be a groupoid. Every fuzzy subgroupoid X need not be a S-fuzzy subgroupoid of X.*

*Proof*: Follows from the fact if X is not a S-groupoid then certainly even if X is a fuzzy subgroupoid of X, X will not be a S-fuzzy subgroupoid of X.

**DEFINITION 3.1.4:** *Let $X_1$ and $X_2$ be two S-groupoids. Let $\mu_1$ and $\mu_2$ be S-fuzzy subgroupoids of $X_1$ and $X_2$ respectively, with respect to a t-norm T. The S-fuzzy subgroupoids $\mu_1$ and $\mu_2$ are Smarandache homomorphic (isomorphic) (S-homomorphic (isomorphic)) if and only if there exists a S-groupoid homomorphic (isomorphism) $\phi : X_1 \rightarrow X_2$ such that $\mu_1 = \mu_2 \circ \phi$. In this situation we say that $\mu_1$ is given by a Smarandache pull back (S-pull back) of $\mu_2$ along $\phi$.*

**DEFINITION 3.1.5:** *Let $(X, \bullet)$ be a S-groupoid. Let G be the family of S-subgroupoids of X. G is called a Smarandache generating family (S-generating family) if for every element $x \in S \subset X$ (S a semigroup of X) there exists $P \in X$ such that $x \in P$.*

**DEFINITION 3.1.6:** *Let $(X, \bullet)$ be a S-groupoid and let Y be a fixed S-subgroupoid of X. Let $(\Omega, A, P)$ be a probability space and let $(V, \odot)$ be a S-groupoid of functions mapping $\Omega$ into X with $\odot$ defined by point wise multiplication in range space. A further restriction in placed on V by assuming that for each $f \in V$; $X_f : \{\omega \in \Omega \mid f(\omega) \in Y\}$ is an element of A. Then $\upsilon : V \rightarrow [0, 1]$ defined by $\upsilon(f) = P(X_f)$ for each $f \in V$ is a fuzzy subgroupoid of V with respect to $T_M$. If $P(X_f)$ is a S-fuzzy subgroupoid of V for each $f \in V$ then we say the S-fuzzy subgroupoid obtained in this manner is called Smarandache function generated (S-function generated).*

An innovative reader can obtain several interesting results in this direction.

**DEFINITION 3.1.7:** *Let $f : X \to X'$ be a S-homomorphism from a S-groupoid X into a S-groupoid X. Suppose that $\nu$ is a S-fuzzy subgroupoid of X' with respect to a t-norm T. Then the fuzzy let $\mu = \nu \circ f$ (defined by $\mu(x) = \nu(f(x))$ for all $x \in X$) is called the Smarandache pre image (S-pre image) of $\nu$ under f.*

**DEFINITION 3.1.8:** *Let $f : X \to X'$ be a S-homomorphism from the S-groupoid X onto S-groupoid X'. Suppose that $\mu$ is a S-fuzzy subgroupoid of X with respect to a t-norm T, then the fuzzy set $\nu$ in $X' = f(X)$ defined by*

$$\nu(y) = \sup_{x \in f^{-1}(y)} \mu(x)$$

*for all $y \in X'$ is called the Smarandache image (S-image) of $\mu$ under f.*

**DEFINITION 3.1.9:** *A fuzzy relation $\mu$ on a S-groupoid X is said to be fuzzy left (resp. right) Smarandache compatible (S-compatible) if for all $x, y, s \in P \subset X$, P a semigroup in X. $\mu(sx, sy)$ [resp $\mu(xs, ys)] \geq \mu(x, y)$. A S-fuzzy C-equivalence relation $\mu$ on a S-groupoid X is said to be a S-fuzzy C-congruence relation on $P \subset X$ if $M(\mu(x, y), \mu(z, \omega)) \leq \mu(xz, y\omega)$ for all $x, y, z, \omega \in P \subset X$.*

**DEFINITION 3.1.10:** *A S-fuzzy subgroupoid $\lambda$ of a S-groupoid X is called a S-fuzzy regular subgroupoid if for all $x \in P \subset X$ their exist $x_1 \in R_\alpha$ such that $(x, t) \in_\alpha \lambda$ implies $(x_1, t) \in_\alpha \vee q_{a,b} \lambda$ for all $t \in (a, c]$ or equivalently for all $x \in \lambda_{sa}$ there exists $x_1 \in R_\alpha$ such that $\lambda(x_1) \geq M(\lambda(x), k)$.*

For more about notations please refer [34].

**DEFINITION 3.1.11:** *A fuzzy inverse subsemigroup $\lambda$ of a S-groupoid X is called a S-fuzzy normal subsemigroup if*

  i. $\lambda(e) \geq k$ for all idempotent e in $P \subset X$ (P-semigroup in X).
  ii. $\lambda$ is fuzzy closed.
  iii. $(x, t) \in_\alpha \lambda \Rightarrow (yxy^{-1}, t) \in_\alpha \vee q_{ab} \lambda$ for all $x, y \in P$ and all $t \in (a, c]$.

*(iii) is equivalent to $\lambda(yxy^{-1}) \geq M(\lambda(x), k)$ for all $x \in \lambda_{s\alpha}$ for all $y \in P \subset S$.*

**DEFINITION 3.1.12:** *Let X be a S-groupoid. A fuzzy subset f of X is called S-fuzzy left ideal of X if $P \circ f \subseteq f$ and f is called the S-fuzzy right ideal of X if $f \circ P \subseteq f$, P a semigroup of X. If f is both a S-fuzzy left and S-fuzzy right ideal of X then f is called a S-fuzzy ideal i.e. $f(xy) \geq \max \{f(x), f(y)\}$ for all $x, y \in P \subset X$.*

**DEFINITION 3.1.13:** *A fuzzy left ideal f is called Smarandache prime (S-prime) if for any two S-fuzzy ideals $f_1$ and $f_2$, $f_1 \circ f_2 \subseteq f$ implies $f_1 \subseteq f$ or $f_2 \subseteq f$.*

**DEFINITION 3.1.14:** *A S-fuzzy left ideal f is called Smarandache quasi prime (S-quasi prime) if for any two S-fuzzy left ideals $f_1$ and $f_2$, $f_1 \circ f_2 \subset f$ implies $f_1 \subset f$ or $f_2 \subset f$. f is* called *Smarandache quasi semiprime (S-quasi semiprime) if for any S-fuzzy ideal g of X; $g^2 \subset f$ implies $g \subset f$.*

**DEFINITION 3.1.15:** *A fuzzy subset f of a S-groupoid X is called Smarandache fuzzy m-system (S-fuzzy m-system) if for any t, s $\in$ [0, 1) and a, b $\in$ P $\subset$ X, f(a) $\geq$ t, f(b) > s imply that there exists an element x $\in$ P such that f (a x b) > t $\vee$ s.*

**DEFINITION 3.1.16:** *A S-fuzzy left ideal f is called Smarandache weakly quasi-prime (S-weakly quasi-prime) if for any two S-fuzzy left ideals $f_1$ and $f_2$ such that $f \subseteq f_1$, $f \subseteq f_2$ and $f_1$ o $f_2 \subseteq f$ then $f_1 \subseteq f$ or $f_2 \subseteq f$.*

**Notation:** Let f be a S-fuzzy left ideal of a S-groupoid X, we define two fuzzy subsets of X denoted by i(f) and i(f) respectively as follows : For all $x \in P \subset S$

i. $i(f)(x) = \vee (t_\alpha \mid x_{t_\alpha} \in f, x_{t_\alpha} \text{ o } P \subset f, t_\alpha \in [0, 1])$.
ii. $i(f)(x) = \vee \{t_\alpha \mid f \text{ o } x_{t_\alpha} \subseteq f, t_\alpha \in [0, 1]\}$.

**DEFINITION 3.1.17:** *A S-groupoid X is called Smarandache strongly semisimple (S-strongly semisimple) if left ideal of X is idempotent. A S-fuzzy ideal f of a S-groupoid X is called idempotent if f = f o f that is $f = f^2$.*

**DEFINITION 3.1.18:** *Let X be a S-groupoid. We call X a Smarandache fuzzy multiplicative groupoid (S-fuzzy multiplicative groupoid) if for any two S-fuzzy ideals $\mu$ and $\delta$ of X satisfying $\mu \leq \delta$ there exists a S-fuzzy ideal $\lambda$ of X such that $\mu = \delta \lambda$.*

(Recall for X a S-groupoid; $\lambda$ a fuzzy subset of X is called a S-fuzzy ideal of X if $\lambda(xy) \geq \max\{\lambda(x), \lambda(y)\}$ for all x, y $\in P \subset X$ (P a proper subset of X which is a semigroup).

**DEFINITION 3.1.19:** *Let I = [0, 1], $I^X$ will denote the set of all mappings $\lambda: X \to I$, M(x, y) will denote the minimum of x and y. Let $\lambda, \mu \in I^X$. Then $\lambda \mu: X \to I$ is defined by $\lambda\mu(x) = \sup \{M (\lambda(y), \mu(z)) \mid x = yz\}$ for all $x \in P \subset X$. A S-fuzzy ideal $\lambda$ of X is called a proper S-fuzzy ideal if $\lambda = \chi_X$. Also for any two –S-fuzzy ideals $\lambda, \mu$, $\lambda < \mu$ will mean $\lambda \leq \mu$, $\lambda \neq \mu$.*

A S-fuzzy ideal $\lambda$ of the S-groupoid X is called S-prime if for any two S-fuzzy ideals $\mu, \delta$ of X $\mu\delta \leq \lambda \Rightarrow \mu \leq \lambda$ or $\delta \leq \lambda$.

For a S-fuzzy ideal $\lambda$ of a S-groupoid X, the Smarandache fuzzy radical (S-fuzzy radical) of $\lambda$ denoted by S(rad $\lambda$) is defined by

$$S(\text{rad } \lambda(x)) = \sup \{ \lambda (x^n) \mid n \in N\} \text{ for all } x \in P \subset X ;$$

$\lambda$ defined relative to P in X. Clearly S(rad $\lambda(x)$) is a S-fuzzy ideal of X.

**DEFINITION 3.1.20:** *A S-fuzzy ideal $\mu$ of a S-groupoid X is said to be idempotent if $\mu = \mu^2$. A S-fuzzy ideal of form $x_r\chi_P$ where $x_r$ is a fuzzy point of the S-groupoid X is* *called a Smarandache principal fuzzy ideal (S-principal fuzzy ideal) of X. (P ⊂ X; P is a semigroup).*

Several other interesting results in this direction can be obtained. We define Smarandache compatible fuzzy relation on a S-groupoid X.

**DEFINITION 3.1.21:** *Let X be a S-groupoid. Let $P_1, P_2, \ldots, P_r$ be r (r >1) proper subset of X which are semigroups under the operations of X. Let λ be a fuzzy relation on the groupoid X. We say λ is a Smarandache compatible (S-compatible) on X, if λ(ac, bd) ≥ min{λ(a, b), λ(c, d)} for all a, b, c, d in each $P_i \subset X$, i.e. for i = 1, 2, …, r.*

*A Smarandache compatible fuzzy equivalence relation (S-compatible fuzzy equivalence relation) on a S-groupoid X is a Smarandache congruence (S-congruence) on X. If λ is compatible on atleast one of the $P_i$'s then we say λ is Smarandache weakly compatible (S-weakly compatible).*

Any interested reader can obtain several exciting relations in this direction. Let λ be a fuzzy relation on a set X and let 0 ≤ α ≤ 1; λ is α-reflexive on X if λ (a, a) = α and λ(a, b) ≤ α for a, b ∈ X. An α-reflexive, symmetric and transitive fuzzy relation on X is a fuzzy α-equivalence relation on X. We call a compatible fuzzy α-equivalence relation on a groupoid, an α-equivalence. A fuzzy α-equivalence relation is a fuzzy equivalence if α = 1 and every fuzzy α-equivalence relation is a G-equivalence. A G-reflexive and transitive fuzzy relation on X is a G-preorder on X.

**DEFINITION 3.1.22:** *Let X and Y be two nonempty set. A mapping $f : X \times X \to Y \times Y$ is called a semibalanced mapping if*

  i. *Given a ∈ X there exists a u ∈ Y such that f (a, a) = (u, u)*
  ii. *f(a, a) = (u, u) and f (b, b) = (v, v) where a, b ∈ X, u, v ∈ Y implies f (a, b) =(u, v).*

Several important properties on semibalanced mapping can be derived.

If f is a semibalanced map from X × X into Y × Y and μ is an α-equivalence fuzzy relation on Y then $f^{-1}(\mu)$ is an α-equivalence fuzzy relation on X.

**DEFINITION 3.1.23:** *Let f be a map from $X \times X$ into $Y \times Y$. A fuzzy relation λ on X is f-invariant if f(a, b) = f($a_1, b_1$) implies that λ(a, b) = λ($a_1, b_1$). A fuzzy relation λ on X is weakly f-invariant if f(a, b) = f ($a_1, b$) implies λ(a, b) = λ($a_1, b$).*

We have several important properties about f-invariant fuzzy relation λ. A mapping is a balanced mapping if

  i. f (a, b) = (u, u) ⇒ a = b.
  ii. f (a, b) = (u, v) ⇒ f (b, a) = (v, u).
  iii. f (a, a) = (u, u) and f(b, b) = (v, v) ⇒ f(a, b) = (u, v)

for all a, b ∈ X and u, v ∈ Y.

A mapping f : X × X → Y × Y is a balanced mapping if and only if it is one to one semibalanced mapping. Let μ be a fuzzy relation on Y and f be a map from X × X into Y × Y. We say μ is f-stable if f (a, b) = (u, u) where a ≠ b ∈ X and u ∈ Y implies that μ(f(a, b)) ≤ μ f (x, x)) for all x ∈ X.

Several innovative results in this direction are derived by [45].

The application of Smarandache fuzzy groupoids to the study of fuzzy automaton will find its place in the seventh chapter of this book.

## 3.2 Smarandache Fuzzy loops and its Properties

In this section we introduce the concept of Smarandache fuzzy loops. Smarandache loops was introduced in [129] and the concept of fuzzy loops was introduced and studied in [118]. The notion of Smarandache fuzzy loops enjoys more properties than that of the fuzzy loops:

Now we just recall the definition of Smarandache loop for the sake of completeness.

**DEFINITION 3.2.1:** *A Smarandache loop (S-loop) is defined to be a loop L such that a proper subset A of L is a subgroup (with respect to the same induced operation) that is $\phi \neq A \subset L$.*

*Example 3.2.1:* Let L = {e, 1, 2, 3, 4, 5} be a loop given by the following table.

| * | e | 1 | 2 | 3 | 4 | 5 |
|---|---|---|---|---|---|---|
| e | e | 1 | 2 | 3 | 4 | 5 |
| 1 | 1 | e | 3 | 5 | 2 | 4 |
| 2 | 2 | 5 | e | 4 | 1 | 3 |
| 3 | 3 | 4 | 1 | e | 5 | 2 |
| 4 | 4 | 3 | 5 | 2 | e | 1 |
| 5 | 5 | 2 | 4 | 1 | 3 | e |

It is easily verified L is a S-loop.

**DEFINITION 3.2.2:** *Let L be a loop. If L has no subloops but only subgroups then we call L a Smarandache subgroup loop (S-subgroup loop).*

*Example 3.2.2:* Let L = {e, 1, 2, 3, 4, 5} be the loop given in example 3.2.1. All subloops in L are subgroups. So L is a S-subgroup loop.

**DEFINITION 3.2.3:** *Let L be a loop. A proper subset A of L is said to be a Smarandache subloop (S-subloop) of L if A is a subloop of L and A itself is a S-loop. i.e. A contains a proper subset B such that B is a group under the operations of L.*

**DEFINITION 3.2.4:** *Let L be a loop. We say a non-empty subset A of L is a Smarandache normal subloop (S- normal subloop) of L if*

i. *A is itself a normal subloop of L.*
ii. *A contains a proper subset B where B is a subgroup under the operations of L.*

*If L has no S-normal subloop we say the loop L is Smarandache simple (S-simple).*

Several results in this direction can be had from [129].

**DEFINITION 3.2.5:** *Let L be a S-loop. A fuzzy subset $\mu$ of L ($\mu : L \to [0, 1]$) is said to be a Smarandache fuzzy loop (S-fuzzy loop) if $\mu$ restricted to at least one of the proper subsets of P of L where R is a group under the operations of L is such that $\mu(xy) \geq \min \{\mu(x), \mu(y)\}$ for every $x, y \in P$. $\mu(x^{-1}) = \mu(x)$ for every $x \in P$. i.e. we can denote this S-fuzzy loop by $\mu_P$ for $\mu$ is defined relative to P. If $\mu : L \to [0, 1]$ is a fuzzy subset and $\mu$ is such that for every proper subset P of L which is a subgroup of L, $\mu : P \to [0, 1]$ is a fuzzy group; then we call $\mu$ a Smarandache strong fuzzy subloop (S-strong fuzzy subloop) of L.*

**THEOREM 3.2.1:** *Every S-strong fuzzy subloop of L is a S-fuzzy subloop of L.*

*Proof*: Left as an exercise for the reader as the proof is straightforward.

This definition of S-fuzzy subloop and S-strong fuzzy subloop will be known as the level I or type I S-fuzzy subloops or S-strong fuzzy subloops. But by default of notation we do not mention it as type I or level I we just say S-fuzzy subloop.

We give the results related with S-fuzzy subloops.

**DEFINITION 3.2.6:** *Let L be a S-loop. A fuzzy subset $\mu$ of L is said to be a Smarandache fuzzy normal subloop (S-fuzzy normal subloop) if*

i. *$\mu : L \to [0, 1]$ is a S-fuzzy subloop relative to P. (P a subgroup in L relative to which $\mu$ is defined).*
ii. *$\mu(xy) = \mu(yx)$ for all $x, y \in P$.*

*Let $\mu$ be a S-fuzzy normal subloop of the loop L.*

*For $t \in [0, 1]$ the set $\mu_t = \{(x, y) \in P \times P \mid \mu(xy^{-1}) \geq t\}$ is called a Smarandache t-level relation (S-t-level relation) of $\mu$.*

Now we define Smarandache fuzzy coset of a S-loop L.

**DEFINITION 3.2.7:** *Let L be a S-loop $\mu$ be a S-fuzzy subloop of L. For any $a \in P \subset L$, P a subgroup of L, a $\mu$ defined by $(a\mu)(x) = \mu(a^{-1}x)$ for every $x \in P \subset L$ is called the Smarandache fuzzy coset (S-fuzzy coset) of the S-loop L defined by a and $\mu$. Thus if $\mu$ happens to be a S-strong fuzzy subloop of L then we see related to this $\mu$ we will have as many S-fuzzy cosets as the number of proper subsets of G which are subgroups of G.*

**DEFINITION 3.2.8:** *Let $\lambda$ and $\mu$ be two S-fuzzy subloops of a S-loop L. Then $\lambda$ and $\mu$ are said to be Smarandache conjugate fuzzy subloops (S-conjugate fuzzy subloops) of L if for some $p \in P$. (P subgroup in L relative to which $\lambda$ and $\mu$ are defined) we have $\lambda(x) = \mu(p^{-1}xp)$ for every $x \in P \subset L$.*

*If $\lambda$ and $\mu$ are S-conjugate fuzzy subloops of the S-loop L then $o(\lambda) = o(\mu)$.*

**DEFINITION 3.2.9:** *Let $\mu$ be a S-fuzzy subloop of a S-loop L then for any $a, b \in P \subseteq L$ the Smarandache fuzzy middle coset (S-fuzzy middle coset) $a\mu b$ of the S-loop L is defined by $(a\mu b)(x) = \mu(a^{-1}xb^{-1})$ for every $x \in P \subset L$. Further we see almost all results true in case of fuzzy groups can be extended in the case of S-fuzzy subloops of level I.*

Now we define Smarandache fuzzy relations on the S-loop L.

**DEFINITION 3.2.10:** *Let $R_\lambda$ and $R_\mu$ be any two Smarandache fuzzy relations (S-fuzzy relations) on a S-loop L. Then $R_\lambda$ and $R_\mu$ are said to be Smarandache conjugate fuzzy relations (S-conjugate fuzzy relations) on the S-loop L if there exists $(g_1, g_2) \in P \times P \subset L \times L$ (P a subgroup in L relative to which fuzzy relations $R_\mu$ and $R_\lambda$ are defined) such that $R_\lambda(x, y) = R_\mu\left(g_1^{-1}xg_1, g_2^{-1}yg_2\right)$ for every $(x, y) \in P \times P$.*

**DEFINITION 3.2.11:** *Let $\mu$ be a S-fuzzy subloop of a loop L. For $a \in P$ (P $\subset$ L a subgroup of L relative to which $\mu$ is defined). The Smarandache pseudo fuzzy coset (S-pseudo fuzzy coset) $(a\mu)^p$ is defined by $((a\mu)^P)(x) = p(a)\mu(x)$ for every $x \in P \subset L$ and for some prime p.*

**DEFINITION 3.2.12:** *A S-fuzzy subloop $\mu$ of a S-loop L is said to be a Smarandache positive fuzzy subloop (S-positive fuzzy subloop) of the loop L if $\mu$ is a positive fuzzy subset of the S-loop L.*

Several interesting and innovative Smarandache analogue results can be established in this direction.

**DEFINITION 3.2.13:** *A S-fuzzy subloop $\mu$ of a S-loop L is Smarandache normalized (S-normalized) if and only if $\mu(e) = 1$ where e is the identity element of the subgroup P, $P \subset L$, P a subgroup of L relative to which $\mu$ is defined.*

**DEFINITION [138]:** *A fuzzy topology $\tau$ on a group G is called a g-fuzzy topology. The pair $(G, \tau)$ is a g-fuzzy topological spaces.*

As loops are non-associative loops the concept of analogous l-topological spaces cannot be directly built using loops instead of groups. To overcome this problem the concept of S-loops becomes handy. To this end we next define Smarandache fuzzy topology on a S-loop L.

**DEFINITION 3.2.14:** *A fuzzy topology $\tau$ on a S-loop L is called a Smarandache l-fuzzy topology (S-l-fuzzy topology) if $\tau$ is a g-fuzzy topology relative to atleast one proper subset P of L where P is a subgroup. We call $\tau$ a Smarandache strong l-fuzzy topology* (S-strong l-fuzzy topology) if $\tau$ is a S-g fuzzy topology with respect to every proper subset P of L.

**THEOREM 3.2.2:** *Every S-strong l-fuzzy topology is a S-l fuzzy topology.*

*Proof:* Direct by the very definition hence left as an exercise for the reader.

**DEFINITION 3.2.15:** *Let $\tau_1$ and $\tau_2$ be S-l-fuzzy topologies on the S-loops $L_1$ and $L_2$ respectively. A function f: $(L_1, \tau_1) \to (L_2, \tau_2)$ is said to be a Smarandache l-fuzzy continuous map (S-l-fuzzy continuous map) from $L_1$ to $L_2$ if*

  i. *f from two proper subsets $P_1$ and $P_2$ of $L_1$ and $L_2$ respectively is such that $f : (P_1, \tau_1) \to (P_2, \tau_2)$ is a g-fuzzy continuous map from $P_1$ to $P_2$.*
  ii. *for every $\mu \in \tau_2, f^{-1}(\mu) \in \tau_1$ and*
  iii. *for every S-fuzzy subloop $\mu$ of $L_2$ in $\tau_2$, $f^{-1}(\mu)$ is a S-fuzzy subloop of $G_1$ in $\tau_1$.*

**DEFINITION 3.2.16:** *Let $(L_1, \tau_1)$ and $(L_2, \tau_2)$ be any two S-l-fuzzy topological spaces. A function f : $(L_1, \tau_1) \to (L_2, \tau_2)$ is said to be a Smarandache l-fuzzy homomorphism (S-l-fuzzy homomorphism) if it satisfies the following 3 conditions.*

  i. *$f : (P_1, \tau_1) \to (P_2, \tau_2)$ is one to one and onto. ($P_i$ is a subgroup of $L_i$ relative to which $\tau_1$ is defined for i = 1, 2).*
  ii. *f is a S-l-fuzzy continuous map from $P_1$ to $P_2$ and*
  iii. *$f^{-1}$ is S-l-fuzzy continuous map from $P_2$ to $P_1$.*

**DEFINITION 3.2.17:** *A fuzzy subset $\mu$ of a S-l-fuzzy topological space $(L, \tau)$ is called a Smarandache $Q_l$-neighborhood (S-$Q_l$-neighborhood) of the fuzzy point $x_t$ (for $x \in P \subset L$, P a subgroup relative to which $\tau$ is defined) if and only if there exists S-fuzzy subloop $\lambda$ (of L) in $\tau$ such that $\lambda \subseteq \mu$, $x_t q \lambda$ (where $x_t q \lambda$ means $t + \lambda (x) > 1$ and $x_t$ is quasi coincident with $\lambda$).*

**DEFINITION 3.2.18:** *A S-l-fuzzy topological space $(L, \tau)$ is said to be a S-l-fuzzy Hausdorff space if and only if for any two fuzzy points $x_t$ and $y_s$ ($x, y \in P$, P a subset of L which is a subgroup relative to which $\tau$ is defined $x \neq y$) there exists S-$Q_l$-neighborhoods $\lambda$ and $\mu$ of $x_t$ and $y_s$ respectively such that $\lambda \cap \mu = \phi_G$.*

We have earlier defined the notion of fuzzy singletons; now we proceed on to define Smarandache p-primary fuzzy subloop of a S-loop L.

**DEFINITION 3.2.19:** *Let L be a S-loop. $\mu$ is called a Smarandache p-primary subloop (S-p-primary subloop) of a S-loop L $\Leftrightarrow$ there exits a prime p such that for all fuzzy singleton $x_t \subset \mu$ with $t > 0$ there exists $n \in N$ (N-integers) such that $p^n(x_t) = 0_t$.*

**DEFINITION 3.2.20:** *Let L be a S-loop $\mu$ is called a Smarandache divisible fuzzy subloop (S-divisible fuzzy subloop) of L if and only if for all singletons $x_t \subset \mu$ with $t > 0$, there exists $n \in N$ there exists a fuzzy singleton $y_t \subset \mu$ such that $n(y_t) = x_t$.*

**DEFINITION 3.2.21:** *Let $\{\mu_\alpha \mid \alpha \in \Omega\}$ be a collection of S-fuzzy subloops of a S-loop L. Then $\mu$ is said to be a Smarandache fuzzy weak direct sum (S-fuzzy weak direct sum) of the $\mu_\alpha$ if and only if*

$$\mu = \sum_{\alpha \in \Omega} \mu_\alpha$$

*and for all $x \in P_\alpha \subset L$ ($P_\alpha$ is a S-fuzzy subloop relative to which $\mu_\alpha$ is defined)*

$$x \neq o\left(\mu_\beta \cap \sum_{\alpha \in \Omega\{\beta\}} \mu_\alpha \right)(x) = 0.$$

*If $\mu$ is a S-fuzzy weak direct sum of $\mu_\alpha$ then we write*

$$\mu = \bigoplus_{\alpha \in \Omega} \mu_\alpha.$$

**DEFINITION 3.2.22:** *A fuzzy subset $\lambda$ of a S-Loop L is said to be a Smarandache $(\in, \in \vee q)$ fuzzy subloop (S-$(\in, \in \vee q)$ fuzzy subloop) of L relative to a proper subset $P \subset L$, P a subgroup in L if for any $x, y \in P$ and $t, r \in (0, 1]$*

   i.  $x_r, y_s \in \lambda$ implies $(xy)_{M(t, r)} \in \vee q \, \lambda$ and
   ii. $x_t \in \lambda \Rightarrow x^{-1} \in \vee q \, \lambda$

*Clearly it is left for the reader to verify that conditions (i) and (ii) are equivalent to*

   i.  $\lambda(xy) \geq M(\lambda(x), \lambda(y), 0.5)$ for all $x, y \in G$ and
   ii. $\lambda(x^{-1}) \geq M(\lambda(x), 0.5)$ for all $x \in G$.

**DEFINITION 3.2.23:** *A S-fuzzy subloop $\lambda$ of a S-loop L is said to be*

i. *Smarandache $(\in, \in)$ -fuzzy normal (S-$(\in, \in)$-fuzzy normal) if for all $x, y \in P \subset L$ (P a subloop of L) and for $t \in (0, 1]$, $x_t \in \lambda \Rightarrow (y^{-1}xy)_t \in \lambda$.*

ii. *Smarandache $(\in, \in \vee q)$-fuzzy normal subloop (S-$(\in, \in \vee q)$-fuzzy normal subloop) if for any $x, y \in P \subset L$ and $t \in (0, 1]$; $x_t \in \lambda$ implies $(y^{-1} xy)_t \in \vee q \lambda$.*

**DEFINITION 3.2.24:** *Two S-fuzzy subloops of a S-loop L are said to be Smarandache equivalent (S-equivalent) if they have the same family of level subloops otherwise they are said to be Smarandache non equivalent (S-non equivalent).*

**DEFINITION 3.2.25:** *Let $\mu$ be a fuzzy subset of a S-loop L. an S-$(\in, \in \vee q)$ fuzzy subloop $\xi$ of L is said to be the S-$(\in, \in \vee q)$ fuzzy subloop generated by $\mu$ in L if $\xi \geq \mu$ and for any other S-$(\in, \in \vee q)$-fuzzy subloop $\eta$ of L with $\eta \geq \mu$ it must be $\eta \geq \xi$.*

**DEFINITION 3.2.26:** *Let $\lambda$ be a S-fuzzy subloop of a S-loop L. For $x \in P \subset L$ $\lambda_x^l$ (resp. $\lambda_x^l$): $L \to I$ defined by $\lambda_x^l(g) = \lambda(gx^{-1})$ (resp. $\lambda_x^r(g) = \lambda(x^{-1}g))$ for all g, f $\in P \subset L$ by $\lambda_x^l(g) = \lambda(gx^{-1})$ (resp. $\lambda_x^r(g) = \lambda(x^{-1}g)$) for all $g \in P \subset L$ (P a subgroup in L relative to which $\lambda$ is defined) is called a S-($\in, \in$) fuzzy left (resp. right) coset of $P \subset L$ determined by x and $\lambda$.*

**DEFINITION 3.2.27:** *Let $\lambda$ be a S-fuzzy subloop of a S–loop L. For any $x \in P \subset L$, $\overline{\lambda}_x$ (resp.($\widetilde{\lambda}_x$): $P \subset L \to [0, 1]$ is defined by $\hat{\lambda}_x(g) = M(\lambda(gx^{-1}), 0.5)$ [resp. $\widetilde{\lambda}_x(g) = M(\lambda(x^{-1}g), 0.5)$] for all $g \in P$ (P a subgroup in L related to which $\lambda$ is defined) is called Smarandache ($\in, \in \vee q$) fuzzy left (resp. right) coset (S-($\in, \in \vee q$) fuzzy left (resp. right) coset) of L determined by x and $\lambda$.*

**DEFINITION 3.2.28:** *A S-subloop H of a loop L is said to be Smarandache quasi normal (S-quasi normal) if for every S-subloop K of L, H.K = K.H.*

**DEFINITION 3.2.29:** *A S-fuzzy subloop $\xi$ of a S-loop L is called Smarandache fuzzy quasi normal (S-fuzzy quasi normal) in G if $\xi o \eta = \eta o \xi$ for every S-fuzzy subloop $\eta$ of P in L.*

*A S-fuzzy subloop $\mu$ of a S-loop L is said to be Smarandache fuzzy maximal (S-fuzzy maximal) if $\mu$ is not constant and for any S-fuzzy subloop $\eta$ of L whenever $\mu \leq \eta$ either $P_\mu = P_\eta$ or $\eta = \chi_P$ where P is a subgroup in L.*

*$P_\mu = \{x \in P \mid \mu(x) = \mu(e)$; e identity element of the subgroup P in G and P is the subgroup relative to which $\mu$ is defined$\}$.*

**DEFINITION 3.2.30:** *Let $\mu$ be a S-fuzzy subloop of a finite S-loop L and let $S_P$ be a S-p-sylow subloop of L. Define a fuzzy subset $\mu_{S_P}$ in L as follows:*

$$\mu_{S_P}(x) = \begin{cases} \mu(x) & \text{if } x \in S_p \\ 0 & \text{if } x \notin S_p \end{cases}$$

*Clearly $\mu_{S_P}$ is a S-fuzzy subloop called a Smarandache p-fuzzy Sylow subloop (S-p-fuzzy Sylow subloop) of $\mu$. A S-fuzzy subloop of L is called Smarandache fuzzy quasi normal (S-fuzzy quasi normal) if its level subloops are S-quasi normal subloop of L.*

The application of S-fuzzy groupoids and S-fuzzy bigroupoids will be described and discussed in Chapter VII of this book. The major application of them being in the fuzzy automaton theory.

### 3.3 Smarandache fuzzy bigroups and Smarandache fuzzy biloops

This section is fully devoted to the introduction of Smarandache fuzzy biloops and Smarandache fuzzy bigroupoid. In chapter I we have introduced the notion of both fuzzy biloops and fuzzy bigroupoids.

Study of just biloops and bigroupoids is elaborately done is [135]. So here we define several of the properties of Smarandache biloops and Smarandache bigroupoids and proceed on to sketch their applications.

**DEFINITION [128]:** *Let $(G, +, \bullet)$ be a non empty set. We call G a Smarandache bigroupoid (S-bigroupoid) if $G_1 \cup G_2$ where $G_1$ and $G_2$ are proper subsets of G satisfying the following conditions*

      i.    *$(G_1, +)$ is a S-groupoid.*
      ii.   *$(G_2, \bullet)$ is a S-semigroup.*

Several interesting and important properties used about them can be had from [135].

**DEFINITION 3.3.1:** *Let $(G = G_1 \cup G_2, +, \bullet)$ be a S-bigroupoid. A fuzzy subset $\mu : G \to [0, 1]$ is called the Smarandache fuzzy sub-bigroupoid (S-fuzzy sub-bigroupoid) of G if $\mu = \mu_1 \cup \mu_2$ where $\mu_1 : G_1 \to [0, 1]$ is a S-fuzzy subgroupoid of $G_1$ ($\mu$ restricted to $G_1$ denoted by $\mu_1$) $\mu_2 : G_2 \to [0, 1]$ is a S-fuzzy subsemigroup of $G_2$ ($\mu$ restricted to $G_2$ is denoted by $\mu_2$). Thus we say if $\mu : G \to [0, 1]$ that is $\mu = \mu_1 \cup \mu_2 : G_1 \cup G_2 \to [0, 1]$ is a S-fuzzy bigroupoid then $\mu_1 : G_1 \to [0, 1]$ which is a S-fuzzy subgroupoid is defined related to a semigroup $P_1$ of $G_1$ and $\mu_2 : G_2 \to [0, 1]$ which is a S-fuzzy subsemigroup is defined related to a subgroup $P_2$ of $G_2$.*

Thus we can give yet another definition of the Smarandache fuzzy sub-bigroupoids which we choose to call it as Smarandache fuzzy sub-bigroupoids of type II.

**DEFINITION 3.3.2:** *Let G be a bigroupoid. A fuzzy subset $\mu$ of G is said to be a Smarandache fuzzy sub-bigroupoid II (S-fuzzy sub-bigroupoid II) of G if there exists $\eta < \mu$ such that $\eta$ is a fuzzy sub-bisemigroup of G.*

**THEOREM 3.3.1:** *Every S-fuzzy sub-bigroupoid $\mu$ of a bigroupoid $G = G_1 \cup G_2$ is a S-fuzzy sub-bigroupoid of type II.*

*Proof:* Straightforward by the very definitions.

It is left for the reader to verify whether a S-fuzzy sub-bigroupoid II, $\mu$ of G is a S-fuzzy sub-bigroupoid of G or will $\mu$ be a S-fuzzy sub-bigroupoid II of G will imply G is a S-bigroupoid.

**DEFINITION 3.3.3:** *Let $(T, +, \bullet)$ be a S-biquasi group. i.e. $(T, +)$ is a S-semigroup and $(T, \bullet)$ is a semigroup. A fuzzy subset $\mu: T \to [0, 1]$ is said to be a Smarandache fuzzy sub-biquasi group (S-fuzzy sub-biquasi group) if $\mu : (T, +) \to [0, 1]$ is a fuzzy subsemigroup and $\mu : (T, \bullet) \to [0, 1]$ is a fuzzy subsemigroup. (Recall the set $(T, +, \bullet)$ is a Smarandache biquasi group (S-biquasi group) if $(T, +)$ is a S-semigroup and $(T, \bullet)$ is just a semigroup).*

We proceed on to define type II Smarandache sub-biquasi group.

**DEFINITION 3.3.4:** *Let (X, +, •) be a biquasi group. A fuzzy subset µ of X is said to be a Smarandache fuzzy sub-biquasi group of type II (S-fuzzy sub-biquasi group of type II) of X if µ : (X, +) → [0, 1] is a S-fuzzy subsemigroup and µ : (X, •) → [0, 1] is a fuzzy subsemigroup of X.*

Thus we have the following result, which is direct by the very definitions.

**THEOREM 3.3.2:** *Let (T, +, •) be a S-biquasi group. Let µ: T → [0, 1] be a S-fuzzy sub-biquasi group then µ is a S-fuzzy sub-biquasi group of type II.*

But if µ : T → [0, 1] is a S-fuzzy sub-biquasi group of type II, will µ be a S-fuzzy sub-biquasi group?

This problem is left for the reader as an exercise.

The applications of S-fuzzy bigroupoids will be dealt in the chapter VII.

**DEFINITION 3.3.5:** *Let $X_1$ and $X_2$ be two bigroupoids. Let $\mu_1$ and $\mu_2$ be S-fuzzy sub-bigroupoids of $X_1$ and $X_2$ respectively with respect to a t-norm T. The S-fuzzy sub-bigroupoids $\mu_1$ and $\mu_2$ are homomorphic if and only if there exists a S-groupoid homomorphism (isomorphism) $\phi : X_1 \to X_2$ such that $\mu_1 = \mu_2 \circ \phi$.*

**DEFINITION 3.3.6:** *Let $(X = X_1 \cup X_2, +, •)$ be a S-bigroupoid. A fuzzy subset $\lambda = \lambda_1 \cup \lambda_2$ of X is called a Smarandache fuzzy bi-ideal (S-fuzzy bi-ideal) of X if*

*$\lambda_1(x, y) \geq \max \{\lambda_1(x), \lambda_1(y)\}$ for all $x, y \in P_1 \subset X_1$; $P_1$ relative to which $\lambda_1$ is defined and $\lambda_2(x, y) \geq \max \{\lambda_2(x), \lambda_2(y)\}$ for all $x, y \in P_2 \subset X_2$; $P_2$ relative to which $\lambda_2$ is defined.*

*A S-fuzzy bi-ideal $\mu = \mu_1 \cup \mu_2$ of $X = X_1 \cup X_2$ is called prime if for any two S-fuzzy bi-ideals $\lambda$ and $\delta$ of X, $\lambda\delta \leq \mu \Rightarrow \lambda \leq \mu$ or $\delta \leq \mu$. Let µ be a S-fuzzy bi-ideal of S. The Smarandache fuzzy biradical (S-fuzzy biradical) of µ denoted by S(birad µ) = $\{\sup \mu_1(x_1^{n_1}) \mid n_1 \in N\}$ for all $x_1 \in P_1 \subset X_1$ relative to which $\mu_1$ is defined} $\cup$ $\{\sup \mu_2(x_2^{n_1}) \mid n_2 \in N\}$ for all $x_2 \in P_2 \subset X_2$ relative to which $\mu_2$ is defined}. S(birad µ) is a S-fuzzy bi-ideal of X.*

*A S-fuzzy bi-ideal $\mu = \mu_1 \cup \mu_2$ of $X = X_1 \cup X_2$ is called Smarandache primary (S-primary) if for any two S-fuzzy bi-ideals $\delta$, $\lambda$ of X. $\delta\lambda \leq \mu \Rightarrow \delta \leq \mu$ or $\lambda \leq S(birad \mu)$. A S-fuzzy bi-ideal $\mu = \mu_1 \cup \mu_2$ of $X = X_1 \cup X_2$ is called Smarandache maximal (S-maximal) if there does not exist any proper S-fuzzy bi-ideal $\lambda$ of X such that $\mu < \lambda$.*

**DEFINITION 3.3.7:** *Let X be a S-bigroupoid, X is called a Smarandache multiplication bigroupoid (S-multiplication bigroupoid) if for any two S-fuzzy bi-ideals µ and $\delta$ of X satisfying $\mu \leq \delta$ there exists a S-fuzzy bi-ideal $\lambda$ of X such that $\mu = \delta\lambda$.*

**DEFINITION 3.3.8:** *A S-fuzzy bi-ideal of the form $x_r \chi_X$ where $x_r$ is a fuzzy point of X is called the Smarandache fuzzy bi-ideal (S-fuzzy bi-ideal) of X. A S-bigroupoid X is said* to be a S-fuzzy principal bi-ideal bigroupoid if every S-fuzzy bi-ideal of X is a S-principal fuzzy bi-ideal.

We proceed on to define Smarandache function generated in case of S-bigroupoids.

**DEFINITION 3.3.9:** *Let $(X, +, \bullet)$ be a S-bigroupoid and let Y be a fixed S-bigroupoid of X. Let $(\Omega, \mathcal{A}, P)$ be a probability space and let $(F, \odot)$ be a S-bigroupoid of functions mapping $\Omega$ into X with $\odot$ defined by pointwise multiplication in the range space. A further restriction is placed on F by assuming that for each $f \in F$, $X_f = \{\omega \in \Omega \,|\, f(\omega) \in Y\}$ is an element of $\mathcal{A}$.*

*Then $\nu : F \to [0, 1]$ defined by $\nu(f) = P(X_f)$ for each f in F is a S-fuzzy sub-bigroupoid of F with respect to $T_m$. A S-fuzzy sub-bigroupoid obtained in this manner is called the Smarandache function generated (S-function generated).*

**DEFINITION 3.3.10:** *Given a S-bigroupoid $(X = X_1 \cup X_2, +, \bullet)$ a t-norm T and a set I, for each $i \in I$ and let $\mu_I$ be a S-fuzzy sub-bigroupoid of X with respect to T define.*

$$\left(\bigcap_{i \in I} \mu_i\right)(x) = \inf_{i \in I} \left(\mu_i(x)\right).$$

**DEFINITION 3.3.11:** *Let $f : X \to X'$ be a S-homomorphism from a S-bigroupoid X into a S-bigroupoid X'. Suppose that $\nu$ is a S-fuzzy sub-bigroupoid of X' with respect to a t-norm T. Then the fuzzy set $\mu = \nu \circ f$ is called the Smarandache preimage (S-preimage) of $\nu$ under f.*

**DEFINITION 3.3.12:** *Let $f : X \to X'$ be a S-homomorphism from S-bigroupoid X onto S-bigroupoid $X^1$. Suppose that $\mu$ is a S-fuzzy sub-bigroupoid of X with respect to a t-norm T, then the fuzzy set $\nu$ in $X' = f(X)$ defined by*

$$\nu(y) = \sup_{x \in f^{-1}(y)} \mu(x)$$

*for all $y \in P' \subset X'$ (P' the bisemigroup contained the S-bigroupoid X' relative to which $\mu$ is defined) is called the Smarandache image (S-image) of $\mu$ under f.*

**DEFINITION 3.3.13:** *We say that a fuzzy set $\mu$ in X (X a S-bigroupoid) has the Smarandache sup property (S-sup property) if for any subset $A \subseteq X$ there exists $a_o \in A$ such that*

$$\mu(a_o) = \sup_{a \in A} \mu(a).$$

**DEFINITION 3.3.14:** *Let $\eta$ and $\mu$ be fuzzy subset in a S-bigroupoid $(G = G_1 \cup G_2, +, \bullet)$ with penultimate bisubset $BP(\eta)$ and $BP(\mu)$ respectively. Then for $g_1 \in G_1, g_2 \in G_2$ we define*

$$\eta_i^l (g_i, \mu_i) = \{x_i \in P_i (\eta_i), g_i = x_i y_i \text{ for some } y_i \in P_i (\mu_i)\}$$

$$\left(\eta^l = \eta_1^l \cup \eta_2^l,\, BP(\eta) = P_1(\eta_1) \cup P_2(\eta_2),\, \eta^r = \eta_1^r \cup \eta_2^r\right)$$

$(\eta = \eta_1 \cup \eta_2,\, \mu = \mu_1 \cup \mu_2\ \mu_i\ \eta_i\ \text{defined on}\ G_i;\ i = 1, 2)$

$\mu_i^r(g_i,\, \eta_i) = \{y_i \in P_i(\mu_i)\,;\, g_i = x_i y_i\ \text{for some}\ x_i \in P_i(\eta_i)\}.$

$\eta_i^r(g_i,\, \mu_i) = \{x_i \in P_i(\eta_i),\, g_i = x_i y_i\ \text{for some}\ y_i \in P_i(\mu_i)\}$

$\mu_i^l(g_i,\, \eta_i) = \{y_i \in P_i(\mu_i)\, /g_i = x_i y_i\ \text{for some}\ x_i \in P_i(\eta_i)\};\ i = 1, 2.$

*The fuzzy set $\eta_1\ o\ \mu_1$ in G defined by*

$$\eta_i\ o\ \mu_i\,(g_i) = \begin{cases} \min \sup\limits_{x \in \eta^l(g\mu)} \eta(x)\ \sup\limits_{y \in \mu^r(g\eta)} \mu(y) & \text{if}\ g \in BP(\eta)BP(\mu) \\ \min(\eta(g),\, \mu(g)) & \text{if}\ g \notin BP(\eta)\,BP(\mu) \end{cases}$$

*is called Smarandache penultimate biproduct (S-penultimate biproduct) of $\eta$ and $\mu$.*

For notations about fuzzy regular and inverse subsemigroups please refer [34] we will be adopting these notions and notations in S-bigroupoids as S-bigroupoids contain proper subset which is a bisemigroup. Using the concept of fuzzy relation $\rho$ on a semigroup S given by [117]. We develop the concept of Smarandache bifuzzy relation.

**DEFINITION 3.3.15:** *Let $X = X_1 \cup X_2$ be a S-bigroupoid. A Smarandache bifuzzy relation (S-bifuzzy relation) $\rho = \rho_1 \cup \rho_2$ on a S-bigroupoid X is one for which $\rho_1$ is a fuzzy relation on $P_1 \subset X_1$; where $P_1$ is a semigroup and $\rho_2$ is a fuzzy relation on $P_2 \subset X_2$ where $P_2$ is a semigroup. A S-bifuzzy relation $\rho = \rho_1 \cup \rho_2$ on a S-bigroupoid $X = X_1 \cup X_2$ is said to be a Smarandache fuzzy C-bicongruence relation (S-fuzzy C-bicongruence relation) on X if $M(\rho_1(x_1,\, y_1),\, \rho(z_1,\, \omega_1)) \leq \rho_1(x_1 z_1,\, y_1 \omega_1)$ for all $x_1,\, y_1,\, z_1,\, \omega_1 \in P_1 \subset X_1$ ($P_1$ semigroup of $X_1$) and $M(\rho_2(x_2,\, y_2),\, \rho_2(z_2,\, \omega_2)) \leq \rho_2(x_2,\, z_2,\, y_2,\, \omega_2)$ for all $x_2,\, z_2,\, y_2,\, \omega_2 \in P_2 \subset X_2$.*

*A S-bifuzzy relation $\rho = \rho_1 \cup \rho_2$ on a S-bigroupoid $X = X_1 \cup X_2$ is said to be a S-fuzzy left (resp. right) bicompatible if for all $x_1, y_1, s_1 \in X_1$ and $x_2, y_2, s_2 \in X_2$ we have $\rho_1(s_1 x_1,\, s_1 y_1)$ (resp. $\rho_1(x_1 s_1,\, y_1 s_1)) \geq \rho_1(x_1,\, y_1)$ and $\rho_2(s_2 x_2,\, s_2 y_2)$ (resp. $\rho_2(x_2 s_2,\, y_2 s_2)) \geq \rho_2(x_2,\, y_2)$.*

**DEFINITION 3.3.16:** *Let $X = X_1 \cup X_2$ be a S-bigroupoid such that X is a S-regular bigroupoid (i., e. $X_1$ has a proper subset $R_1$ which is a regular semigroup and $X_2$ has a proper subset $R_2$ which is a regular semigroup).*

*A S-fuzzy sub-bigroupoid $\lambda$ of X is called a S-fuzzy regular sub-bigroupoid if for all $x_1 \in R_1 \subset X_1$ and $x_2 \in R_2 \subset X_2$ there exists $x'_1, x'_2 \in T_x$ where $T_x = \{x'_1 \in R_1 \subset X_1\ |\ x_1 x'_1 x_1 = x_1\} \cup \{x'_2 \in R_2 \subset X_2\ |\ x_1 x'_2 x_1 = x_1\}$ (for all $x_1 \in R_1$ and for all $x_2 \in R_2$) such that $(x_1, t) \in_a \lambda_1 \Rightarrow (x'_1, t) \in_a \vee q_{(a,\, b)} \lambda_1$ and $(x_2, t) \in_a \lambda_2 \Rightarrow (x'_2, t) \in_a \vee q_{(a,\, b)} \lambda_2$ for all $t \in (a, c]$.*

Several results in this direction can be had in an analogous way. We suggest some problems in section 3.4. We just define Smarandache fuzzy prime bi-ideal of a S-bigroupoid X; Smarandache fuzzy quasi prime bi-ideal of a S-bigroupoid and Smarandache fuzzy weakly quasi prime bi-ideal of a S-bigroupoid X.

**DEFINITION 3.3.17:** *Let $X = X_1 \cup X_2$ be a S-bigroupoid, a S-fuzzy left bi-ideal f is called Smarandache prime (S-prime) if for any two S-fuzzy bi-ideals $f_1$ and $f_2$, $f_1 \circ f_2 \subset f$ implies $f_1 \subset f$ or $f_2 \subset f$.*

**DEFINITION 3.3.18:** *Let $X = X_1 \cup X_2$ be a S-bigroupoid. A S-fuzzy left bi-ideal f of X is called Smarandache quasi prime (S-quasi prime) if for any two S-fuzzy left bi-ideals $f_1$ and $f_2$, $f_1 \circ f_2 \subseteq f$ implies $f_1 \subseteq f$ or $f_2 \subset f$; f is called Smarandache quasi-semi prime (S-quasi semiprime) if for any S-fuzzy left bi-ideal g of X; $g^2 \subseteq f$ implies $g \subseteq f$.*

**DEFINITION 3.3.19:** *Let X be a S-bigroupoid. A S-fuzzy left bi-ideal f is called Smarandache weakly quasi prime (S-weakly quasi prime) if for any two S-fuzzy left bi-ideals $f_1$ and $f_2$ such that $f \subseteq f_1, f \subseteq f_2$ and $f_1 \circ f_2 \subseteq f$ then $f_1 \subseteq f$ or $f_2 \subseteq f$.*

Several important results in this direction can be obtained by using [133, 128]. We have given some problems for any researcher to tackle; as the book is for researchers to work more on Smarandache notions we restrain ourselves by giving elaborate proofs.

Now we proceed on to define Smarandache biloops. Several properties related to bigroups can be easily formulated and obtained we define certain concepts about these new structures like S-bitale, S-bitip and S-penultimate biloops and give some results about them.

**DEFINITION 3.3.20:** *Let (L, +, •) be a biloop. We call L a Smarandache biloop (S-biloop) if L has a proper subset P which is a bigroup.*

Now we define Smarandache fuzzy sub-biloop of a biloop L.

**DEFINITION 3.3.21:** *Let (L, +, •) be a S-biloop. We call $\mu : L \to [0, 1]$ a Smarandache fuzzy sub-biloop (S-fuzzy sub-biloop) of L if $\mu$ related to the sub-bigroup P in L, is a fuzzy bigroup i.e. for $\mu_P : P \to [0, 1]$ is a fuzzy sub-bigroup where $\mu_P$ denotes the restriction of $\mu$ to P.*

This S-fuzzy sub-biloop of L we call as type I but by default of notation we just call the biloop L as S-fuzzy sub-biloop of L. We now proceed on to give the Smarandache fuzzy sub-biloop of type II.

**DEFINITION 3.3.22:** *Let L be a biloop. We call $\mu$ a Smarandache fuzzy sub-biloop of type II (S-fuzzy sub-biloop of type II) if the fuzzy subset $\mu$ of L is such that exists a fuzzy subset $\eta$ of L with $\eta < \mu$ and $\eta$ is a Smarandache fuzzy sub-bigroup of L. That is $\eta$ on a fuzzy subset of L where $\eta$ is a fuzzy sub-bigroup of L.*

We are requesting the reader to find relations between S-fuzzy sub-biloops of type I and type II.

**DEFINITION 3.3.23:** *Let $(L = L_1 \cup L_2, +, \bullet)$ be a S-biloop. Two S-fuzzy sub-biloops $\eta$ and $\mu$ of the S-biloop L are said to be Smarandache equivalent (S-equivalent) denoted by $\mu \approx \eta$ if both $\mu$ and $\eta$ have the same S-chain of level subgroups. (If $\mu$ is a S-fuzzy sub-biloop of the S-biloop X; $\mu_t$ for each $t \in [0, \mu(e)]$ ($e_1$ – identity element of $P_1 \subset X_1$, $P_1$ a subgroup of $X_1$ relative to $\mu_1$; $P_2$ is a subgroup of $x_2$ where $e_2$ is the identity element of $P_2 \subset X_2$ relative to $\mu_2$ where $\mu = \mu_1 \cup \mu_2$), here $\mu(e) = \mu_1(e_1) \cup \mu_1(e_2)$ is a S-sub-biloop of the given S- biloop L called the Smarandache level sub-biloop (S-level sub-biloop) of $\mu$. The set of all S-level sub-biloop of a S-fuzzy sub-biloop forms a chain.*

**DEFINITION 3.3.24:** *Two S-fuzzy subloops $\eta$ and $\mu$ of a S-biloop L are said to be S-equivalent denoted by $\mu \approx_B \eta$ if $\mu$ and $\eta$ have the same chain of S-level subloops. Let $\mu$ be a fuzzy subset of a S-biloop L. Then the Smarandache penultimate (S-penultimate) subset SP ($\mu$) of $\mu$ in L defined by*

*SP ($\mu$) = $\{x \in B \subset L \mid \mu(x) > \inf \mu\}$ ($B \subset L$ is a sub-bigroup in the biloop L). In case $\mu$ is a S-fuzzy sub-biloop of a S-biloop L then SP ($\mu$) is a S-sub-biloop of L provided $\mu$ is non -constant and SP ($\mu$) is called the S-penultimate sub-biloop of $\mu$ in L.*

Now we define yet another interesting property about S-biloop L.

**DEFINITION 3.3.25:** *Let $\eta$ and $\mu$ be S-fuzzy sub-biloops of a loop L. Then $\eta$ is said to be obtained by a Smarandache shift (S- shift) of the range set Im $\mu$ if there exists a one to one order preserving map $\phi$ from Im $\mu$ onto Im$\eta$ such that $\eta = \phi \circ \mu$ where 'o' is the composition of mappings.*

**DEFINITION 3.3.26:** *Let L be a S-biloop, $\mu$ be a S-fuzzy sub-biloop of L. i.e. $\mu = \mu_1 \cup \mu_2 : L_1 \cup L_2 \to [0, 1]$ is a map, that $\mu_I$ restricted to a subgroup $P_i$ in $L_i$ is such that*

$$\mu_i(x_i y_i) \geq \min \{\mu_I(x_i), \mu_I(y_i)\}, x_i, y_i \in P_i \subset L_i \ (i = 1, 2)$$
$$\mu_I(x_i^{-1}) = \mu_I(x_i) \text{ for } x_i \in P_i \subset L_i \ (i = 1, 2).$$

*If $\mu$ is a S-fuzzy sub-biloop then it attains its suppremum at $e = e_1 \cup e_2$ of L, that is*

$$\sup_{x_i \in P_i} \mu_i(x_i) = \mu(e_i),$$

*where $i = 1, 2$.*

*We call $\mu(e) = \mu_1(e_1) \cup \mu_2(e_2)$ to be the Smarandache bitip (S-bitip) of the S-fuzzy sub-biloop $\mu$ of L.*

*The S-fuzzy sub-biloop may or may not attain its infimum. We shall write in short inf $\mu$ for*

$$\bigcup_i \inf_{x_i \in P_i} \mu(x_i)$$

*and refer to if as Smarandache bitale (S-bitale) of the S-biloop $\mu = \mu_1 \cup \mu_2$ of L if*

$$\inf \mu = \left[\inf_{x_i \in P_1} \mu_1(x_1) \cup \inf_{x_i \in P_2} \mu_2(x_i)\right]$$

*Two S-fuzzy sub-biloops $\mu$ and $\eta$ are bisimilarly bounded if they have the same S-bitip and same S-bitale that is $\mu(e) = \eta(e)$ i.e. $\mu_1(e_1) \cup \mu_2(e_2) = \eta_1(e_1) \cup \eta_2(e_2)$ and $\inf \mu = \inf \eta$.*

Refer bigroup paper[89].

### 3.4 Problems

This section is purely devoted to the problems. We have proposed 111 problems for the reader to solve. The problems are of different levels: some require examples/ counter examples, others are routine like theorems some of them are characterizations. Thus tackling of these problems will make the reader not only strong in Smarandache notions but will induce him to discover more in Smarandache notions.

**Problem 3.4.1:** Let $\mu$ be a S-fuzzy subgroupoid of X with respect to Min. Is $\mu$ a S-fuzzy subgroupoids of X with respect to any t norm T? Justify / illustrate your claim.

**Problem 3.4.2:** Let $(X, \bullet)$ be a S-groupoid. Let $\mu : X \to [0, 1]$ be a S-fuzzy subgroupoid of X with respect to a t norm T. Can $\mu$ be extended to a S-fuzzy subgroupoid $\mu'$ of X' with respect to the same t-norm T, where X' is the identity extension of X? Substantiate your answer.

**Problem 3.4.3:** Let $(\Omega, A, P)$ be a probability space with P a fuzzy subset of A, $(A, \cap)$ is a groupoid and P is a fuzzy subgroupoid of A with respect to $T_M$. Will P be a S-fuzzy subgroupoid of A?

**Problem 3.4.4:** Let $(\Omega, A, P)$ be a probability space and $\tau$ be a non-empty set $\wp(\tau)$ denotes the power set of $\tau$. Let $\phi: \Omega \to \tau$ be a mapping. Suppose $(X, \bullet)$ is a S-groupoid $\mu: X \to [0, 1]$, $\mu(x) = P(\phi^{-1}(\psi(x)))$ for all $x \in X$ where $\psi: X \to \wp(\tau)$, such that $\phi^{-1}(\psi(x)) \in A$ for all $x \in X$, $\psi(xy) \supseteq \psi(x) \cap \psi(y)$ for all $x, y \in X$. Prove $\mu$ is a S-fuzzy subgroupoid of X.

**Problem 3.4.5:** Is every S-function generated S-fuzzy subgroupoid a S-subgroupoid generated?

**Problem 3.4.6:** Can we say every S-subgroupoid, S-generated fuzzy subgroupoid is isomorphic to a S-function generated S-fuzzy groupoid?

**Problem 3.4.7:** Will every S-fuzzy subgroupoid with respect to Min a S-subgroupoid generated? Justify.

**Problem 3.4.8:** Let $\mu = \nu \circ f$ be the S-preimage of $\nu$ under f, can $\mu$ be a S-fuzzy subgroupoid of X with respect to T?

**Problem 3.4.9:** Prove a S-fuzzy C-equivalence relation $\mu$ on a S-groupoid X is a S-fuzzy C-congruence relation if and only if it is both fuzzy left and right compatible.

**Problem 3.4.10:** Prove a fuzzy subset $\lambda$ of X is a S-fuzzy regular subgroupoid of X if and only if $\lambda_t$ is a regular S subgroupoid of X for all $t \in (a, k]$.

**Problem 3.4.11:** Prove a fuzzy subset $\lambda$ of X, X a S-groupoid is a fuzzy normal subgroupoid of X if and only if $\lambda_t$ is normal subgroupoid of S for all $t \in (a, k]$.

**Problem 3.4.12:** Prove a S-fuzzy left ideal f of a S-groupoid X is prime if and only if for any two fuzzy points $x_r, y_t \in P \subset X$ (r t > 0) $x_r \circ P \circ y_t \circ S \subset f$ implies that $x_r \in f$ or $y_t \in f$.

**Problem 3.4.13:** Prove a left ideal L of X, X a S-groupoid is prime if and only if $f_L$ is a S-prime fuzzy left ideal of X.

**Problem 3.4.14:** Prove a S-fuzzy left ideal f of a S-groupoid X is S-quasi prime if and only if for any two fuzzy points $x_r, y_t \in P \subset X$, (r t > 0), $x_r \circ P \circ y_t \subseteq f$ implies $x_r \in f$ or $y_t \in f$. (P is a semigroup contained in X relative to which f is defined).

**Problem 3.4.15:** Prove a left ideal L of a S-groupoid X is S-quasi prime if and only $f_L$ is a S-quasi prime fuzzy left ideal of X.

**Problem 3.4.16:** Let M be a subset of X. Prove M is a S-m-system of X if and only if $f_m$ is a S-fuzzy m-system.

**Problem 3.4.17:** Let f be a S-fuzzy left ideal of a S-groupoid X, prove f is S-quasi prime if and only if $1 - f$ is a S-fuzzy m-system.

**Problem 3.4.18:** Prove a left ideal L is S-weakly quasi prime if and only if $f_L$ is S-weakly quasi prime.

**Problem 3.4.19:** Prove if X is a S-commutative groupoid, f a S-fuzzy left ideal of X, then

    i.    f is S-prime if and only if f is S-quasi prime.
    ii.   f is S-prime if and only if f is S-weakly quasi prime.

**Problem 3.4.20:** Let f be a S-left ideal of a S-groupoid X. Then prove i(f) is the S-largest fuzzy ideal of X contained in f.

**Problem 3.4.21:** Let X be a S-groupoid with identity e and f a S-prime fuzzy left ideal of S. If $i(f) \neq 0$, is i(f) a S-quasi prime fuzzy ideal of S?

**Problem 3.4.22:** Prove a S-groupoid X is strongly semisimple if and only if every S-fuzzy left ideal of X is idempotent.

**Problem 3.4.23:** Let X be a S-commutative groupoid. Is the S-fuzzy left ideals of X S-quasi prime if and only if they form a chain and X is strongly semisimple?

**Problem 3.4.24:** Let X be a S-groupoid which is a S-fuzzy multiplication groupoid. If $f : X \to T$ is an epimorphism then will T be a S-fuzzy multiplication S-groupoid ?

**Problem 3.4.25:** Can we prove a S-fuzzy ideal $\lambda$ of a S-groupoid X is S-prime if and only if for any two S-fuzzy points $x_r$, $y_s$, $x_r y_s \in \lambda \Rightarrow x_r \in \lambda$ or $y_s \in \lambda$?

**Problem 3.4.26:** A S-fuzzy ideal $\lambda$ of a S-groupoid X is S-primary if and only for any two fuzzy points $x_r$, $y_s$, $x_r y_s \in \lambda \Rightarrow x_r \in \lambda$ or $y_s^n \in \lambda$ for some n > 0. (Is this statement a valid one?).

**Problem 3.4.27:** Let $\lambda$ be a S-fuzzy ideal of a S-groupoid X. Is $\lambda^n$ also a S-fuzzy ideal of X for all n > 0 and $\lambda^n \geq \lambda^{n+1}$?

**Problem 3.4.28:** Can we say $\lambda$ is S-prime then $\lambda = S(\text{rad } \lambda)$?

**Problem 3.4.29:** For any S-fuzzy ideal $\lambda$ of a S-groupoid X. Will $S(\text{rad } \lambda) = \cap \{\mu \mid \mu$ is a S-prime fuzzy ideal of X such that $\mu \geq \lambda\}$?

**Problem 3.4.30:** If $\lambda$ is a S-primary fuzzy ideal of a S-groupoid X then prove, $S(\text{rad } \lambda)$ is a S-prime fuzzy ideal of X.

**Problem 3.4.31:** Let $f: X \to T$ be an epimorphism of S-groupoids and $\lambda$, $\mu$ be S-fuzzy ideals of X and T respectively. Prove $f^{-1}(\mu)$ is a S-fuzzy ideal of X, and $f(\lambda)$ is a S-fuzzy ideal of T.

**Problem 3.4.32:** If $S(\text{rad } \lambda)$ is S-prime then $\lambda$ is S-primary – prove !

**Problem 3.4.33:** Let $\lambda$ be S-prime. Then prove for all positive integer m, $\lambda^m$ is S-primary and its S-fuzzy radical is $\lambda$.

**Problem 3.4.34:** Let $\lambda$ be S-prime and $\lambda^m \neq \lambda^{m+1}$ for all m > 0. Is $\lambda^m$ S-prime? Justify.

**Problem 3.4.35:** If $\lambda$ is S-primary, prove that $\lambda = \mu^n$ for some positive integer n where $\mu = S(\text{rad } \lambda)$.

**Problem 3.4.36:** If $\lambda$ is a proper S-prime and $\mu$ is a S-fuzzy ideal of X, X-S-groupoid such that $\mu \leq \lambda^n$, $\mu \leq \lambda^{n+1}$ for some n > 0 then $\lambda^n = \mu$; $y_t \chi_X$ where $y_i \notin \lambda$ - Prove.

**Problem 3.4.37:** If $\gamma$ is a non –idempotent S-prime fuzzy ideal of X then prove that there is at most one S-prime fuzzy ideal $\mu < \gamma$ such that there is no S-prime ideal between. $\gamma$ and $\mu$.

**Problem 3.4.38:** If S-prime fuzzy ideal of a S-groupoid X are linearly ordered with respect to the inclusion relation '≤' then every non-idempotent S-prime fuzzy ideal is S-principal – prove.

**Problem 3.4.39:** Let $\tau$ be the unique S-maximal fuzzy ideal of X, X a S-groupoid having the sup property. Then for every S-fuzzy ideal $\mu$ of X either $\mu = \tau^n$ for some n > 0 or $\mu \leq \tau^n$, prove or disprove?

**Problem 3.4.40:** Let X be a S-commutative groupoid with every proper subset, which is a semigroup, is a monoid. Prove every S-fuzzy principal ideal is S-fuzzy multiplication groupoid.

**Problem 3.4.41:** If $\mu$ is a S-compatible fuzzy relation on a S-groupoid X and f is a S-groupoid homomorphism form $X \times X$ into $Y \times Y$ will $f^1(\mu)$ be a S-compatible relation on X?

**Problem 3.4.42:** Let $\lambda$ be a S-compatible fuzzy relation on the S-groupoid X and f is a S-groupoid X and f is a S-groupoid homomorphism from $X \times X$ into $Y \times Y$ will $f(\lambda)$ be a S-compatible fuzzy relation on Y?

**Problem 3.4.43:** If $\lambda$ is a G-pre order on the set X prove $\lambda \circ \lambda = \lambda$.

**Problem 3.4.44:** If $\mu$ is an $\alpha$-congruence fuzzy relation on the S-groupoid X and f is a S-groupoid homomorphism from $D \times D$ into $T \times T$ which is a semibalanced map then prove $f^{-1}(\mu)$ is an $\alpha$-congruence on D.

**Problem 3.4.45:** Let f be a semibalanced map and a S-groupoid homomorphism from $D \times D$ onto $S \times S$. Prove if $\lambda$ is an $\alpha$-congruence fuzzy relation on D which is weakly f-invariant then $f(\lambda)$ is an $\alpha$-congruence on S.

**Problem 3.4.46:** If f is a S-groupoid homomorphism and a semibalanced map from $D \times D$ onto $T \times T$ and $\lambda$ is a G-congruence fuzzy relation on D which is weakly f-invariant then prove or disprove $f(\lambda)$ is a G-congruence on S with $\delta(f(\lambda)) = \delta(\lambda)$.

**Problem 3.4.47:** If f is a S-groupoid homomorphism and a balanced map from $D \times D$ into $T \times T$ and $\lambda$ is a G-congruence fuzzy relation on D then prove $f(\lambda)$ is a G-congruence on S with $\delta(f(\lambda)) = \delta(\lambda)$.

**Problem 3.4.48:** Obtain some applications of S-fuzzy groupoids to S-fuzzy automatons.

**Problem 3.4.49:** Let $\mu$ be a S-fuzzy normal subloop of a loop L and t ∈ [0, 1]. Prove $\mu_t$ is a congruence relation on the subgroup P in L relative to which $\mu$ is defined.

**Problem 3.4.50:** Let $\mu$ be a S-fuzzy normal subloop of a S-loop L and x ∈ P ⊂ L. Then prove $\mu(xy) = \mu(y)$ for every y ∈ P if and only if $\mu(x) = \mu(e)$, e the identity element of the subgroup P relative to which $\mu$ is defined.

**Problem 3.4.51:** Let L be a S-loop $\mu$ and $\lambda$ be S-conjugate subloop of L relative to the subgroup P of L then prove $o(\lambda) = o(\mu)$.

**Problem 3.4.52:** Let $\mu$ be a S-fuzzy subloop of a S-loop L then for any $a \in P \subset L$ prove the S-fuzzy middle coset $a \mu a^{-1}$ of the S-loop L is also a S-fuzzy subloop of L.

**Problem 3.4.53:** Let $\mu$ be a S-fuzzy subloop a of S-loop L and $a \mu a^{-1}$ be a S-fuzzy middle coset of the S-loop L relative to the subloop P in L then $o(a \mu a^{-1}) = o(\mu)$ for any $a \in P$.

**Problem 3.4.54:** Show if $\mu$ is a S-fuzzy subloop of a finite S-loop then prove in general $o(\mu) \neq o(L)$.

**Problem 3.4.55:** Let $\lambda$ be a S-fuzzy subloop of a S-loop L relative to a proper subset $P \subset L$ (P a subgroup in L). If $\lambda$ and $\mu$ are conjugate fuzzy subsets of the S-loop L then prove $\mu$ is a S-fuzzy subloop of the S-loop L.

**Problem 3.4.56:** Let $R_\lambda$ and $R_\mu$ be any two fuzzy relations on a S-loop L. If $R_\lambda$ and $R_\mu$ are generalized conjugate fuzzy relations on the subgroup P of L; then prove $R_\lambda$ and $R_\mu$ are conjugate fuzzy relations on the S-loop L.

**Problem 3.4.57:** Let $\mu$ be a S-fuzzy normal subloop of the S-loop L. Then prove for any $g \in P$ we have $\mu(gxg^{-1}) = \mu(g^{-1}xg)$ for every $x \in P \subset L$ (P a subgroup L relative to $\mu$).

**Problem 3.4.58:** Let $\lambda$ and $\mu$ be S-conjugate fuzzy subloops of a S-loop L relative to the subgroup P in L. Prove $\lambda \times \mu$ and $\mu \times \lambda$ are S-conjugate fuzzy relations on L.

**Problem 3.4.59:** Prove if $\mu$ is a S-positive fuzzy subloop of a S-loop L then any two S-pseudo fuzzy cosets of $\mu$ are either identical or disjoint.

**Problem 3.4.60:** Let $\mu$ be a S-fuzzy subloop of a S-loop L then prove the S-pseudo fuzzy coset $(a\mu)^p$ is a S-fuzzy subloop of the S-loop L for every $a \in P$, $P \subset L$ is the subgroup relative to which $\mu$ is defined.

**Problem 3.4.61:** Let $\mu$ be a S-fuzzy subloop of a S-loop L and $R_\mu: L \times L \to [0, 1]$ be given by $R_\mu(x, y) = \mu(xy^{-1})$ for every $x, y \in P \subset L$ (P a subgroup of L relative to which $\mu$ is defined). If $\lambda$ is a fuzzy subset of the S-loop L such that $\lambda \subseteq \mu$ then prove $(a \lambda)^p$ is a pre class of $R_\alpha$ for any $a \in P \subset L$.

**Problem 3.4.62:** Let $\mu$ and $\lambda$ be any two S-fuzzy subloops of a S-loop L and $R_{\mu \cap \lambda} : L \times L \to [0, 1]$ given by $R_{\lambda \cap \mu}(x, y) = (\lambda \cap \mu)(xy^{-1})$ for every $x, y \in P \subset L$ (P a subgroup of L relative to which $\mu$ and $\lambda$ are defined). Prove a $R_{\lambda \cap \mu}$ is a similarity relation on P only when both $\mu$ and $\lambda$ are normalized.

**Problem 3.4.63:** Let $L_1$ and $L_2$ be any two S-loops. If $\tau_1$ is a S-l-fuzzy topology on the S-loop $L_1$ and $\tau_2$ is an indiscrete S-l-fuzzy topology on the S-loop $G_2$, then prove every function is a S-l-fuzzy continuous map.

**Problem 3.4.64:** Let $\tau_1$ and $\tau_2$ be two discrete S-l-fuzzy topology on the S-loops, $(L_1, \bullet)$ and $(L_2, *)$ respectively. Then prove every S-loop homomorphism from $(L_1, \tau_1) \to (L_2, \tau_2)$ is a S-l-fuzzy continuous map but not conversely.

**Problem 3.4.65:** Let $\tau_1$ and $\tau_2$ be any two S-l-fuzzy topologies on the S-loops $(L_1, \bullet)$ and $(L_2, *)$ respectively. Then prove every S-loop homomorphism need not in general be a S-l-fuzzy continuous map.

**Problem 3.4.66:** Let $f : (L_1, \tau_1) \to (L_2, \tau_2)$ is a S-l-fuzzy homomorphism. Then prove $(L_1, \tau_1)$ is a S-l-fuzzy Hausdorff space if and only if $(L_2, \tau_2)$ is a S-l-fuzzy Hausdorff space.

**Problem 3.4.67:** Formulate the Existence theorem for S-loops. That is "Let $\mu$ be a S-fuzzy subloop of a S-loop L. The congruence class $[x]_\mu$ of $\mu_t$ containing the element x of the subgroup P in the S-loop L (relative to which $\mu$ is defined) exists only when $\mu$ is a S-fuzzy normal subloop of the S-loop L". Is this the existence theorem for S-loops? Justify.

**Problem 3.4.68:** Prove $\mu$ is S-fuzzy divisible if and only if $\mu_t$ is S-fuzzy divisible.

**Problems 3.4.69:** Prove for all $x, y \in P \subset L$ (P a subgroup of the S-loop L) and $n \in N$ $ny = x$ implies $\mu(x) = \mu(y)$ for all S-divisible fuzzy subloops $\mu$ of L if and only if P is torsion free.

**Problem 3.4.70:** Prove if B and C be S-fuzzy subloops of a S-loop L. Then $(B \cap C)^* = B^* \cap C^*$ and $(B \cap C)_t = B_t \cap C_t$ for all $t \in (0, \min \{B(0), C(0)\})$.

**Problem 3.4.71:** Let $\{\mu_\alpha \mid \alpha \in \Omega\}$ be a collection of all S-fuzzy subloops of the S-loop L such that $\mu_\alpha \subseteq \mu$ and $\mu_\alpha (0) = \mu (0)$ for all $\alpha \in \Omega$. Prove $\mu = \bigoplus_{\alpha \in \Omega} \mu_\alpha$ if and only if $\mu^* = \bigoplus_{\alpha \in \Omega} \mu_\alpha^*$ and $\mu = \sum_{\alpha \in \Omega} \mu_\alpha$.

**Problem 3.4.72:** Let $\mu$ be a S-fuzzy subloop of a S-loop L; $\mu$ defined relative to P, P a subgroup contained in L. If $x_1, \ldots , x_n \in P \subset L$ are such that $\min \{\mu(x_1), \ldots, \mu(x_{n-1})\} > \mu(x_n)$ then $\mu(x_1 + \ldots + x_n) = \mu(x_n)$.

**Problem 3.4.73:** Is an S-($\in, \in$)-fuzzy normal subloop a S-($\in, \in \vee q$)-fuzzy normal subloop? Is the converse possible?

**Problem 3.4.74:** Is a S-($\in, \in \vee q$)-fuzzy subloop of a S-loop L is said to be proper if it is not constant on G? Justify.

**Problem 3.4.75:** Let $f : L \to L_1$ be a S-loop homomorphism of any two S-subloops. If $\lambda$ and $\mu$ are S-fuzzy subloops of L and $L_1$ respectively; Prove $f(\lambda)$ and $f^{-1}(\mu)$ are S-fuzzy subloops of $L_1$ and L respectively.

**Problem 3.4.76:** Let $f : L \to L_1$ be a S-loop homomorphism. Let $\lambda$ and $\mu$ be any two S-fuzzy subloops of L and $L_1$ respectively. Prove

i. $f^{-1}(\mu)$ is a S-fuzzy subloop of L.
ii. $f(\lambda)$ is a S-fuzzy subloop of $f(L)$.

**Problem 3.4.77:** Let $\lambda$ be a S-fuzzy subloop of a S-loop L. Then prove $\lambda$ is S-($\in$, $\in \vee q$) – fuzzy normal if and only if $\hat{\lambda}_x = \tilde{\lambda}_x$ for all $x \in P \subset L$ where $\lambda$ is defined relative to P.

**Problem 3.4.78:** Let $\lambda$ be a S-fuzzy normal subloop of a S-loop S. Let F be the set of all cosets of $\lambda$ in $P \subset L$ (P is the subgroup of L relative to which $\lambda$ is defined). Prove F is a group of all S-fuzzy cosets of P in L determined by $\lambda$ where the multiplication is defined by $\lambda_x \lambda_y = \lambda_{xy}$ for all x, y $\in$ P $\subset$ L. Let $\overline{\lambda} : F \to [0,1]$ be defined by $\overline{\lambda}(\lambda x) = M(\lambda(x^{-1}), 0.5)$ for all $x \in P \subset L$. Prove $\lambda$ is a S-fuzzy normal subgroup of F.

Several interesting problems can be defined in this direction.

**Problem. 3.4.79:** Will a homomorphic preimage of a S-quasi normal subloop be S-quasi normal?

**Problem 3.4.80:** Let $\mu: G \to [0, 1]$ be S-fuzzy sub-bigroupoid of type II. When will G be a S-bigroupoid?

**Problem 3.4.81:** Let $\mu : T \to [0, 1]$ is a S-fuzzy sub-biquasi group of type II. Is (T, +, •) a S-biquasi group?

**Problem 3.4.82:** Let $\mu$ be a S-fuzzy sub-biloop of a S-biloop L. Can $\mu$ be a S-fuzzy sub-biloop of type II? Justify your claim.

**Problem 3.4.83:** Construct an example of a biloop L to distinguish between the S-fuzzy sub-biloops of type I and type II.

**Problem 3.4.84:** Prove A S-fuzzy bi-ideal $\lambda$ of a S-bigroupoid X is S-prime if and only if for any two fuzzy points $x_r$, $y_s$, $x_r y_s \in \lambda$ implies $x_r \in \lambda$ or $y_s \in \lambda$.

**Problem 3.4.85:** A S-fuzzy bi-ideal of a S-bigroupoid X is S-primary if and only if for any two fuzzy points $x_r$, $y_s$, $x_r y_s \in \lambda$ implies $x_r \in \lambda$ or $y^n_s \in \lambda$ for some n > 0. Prove.

**Problem 3.4.86:** Prove for any S-fuzzy bi-ideal $\lambda$ of a S-bigroupoid X, S(birad $\lambda$) = $\cap \{\mu \mid \mu$ is a S-fuzzy prime bi-ideal of X such that $\mu \geq \lambda\}$.

**Problem 3.4.87:** Prove if $\lambda$ is a S-primary fuzzy bi-ideal of X, X a S-bigroupoid then S(birad $\lambda$) is a S-prime fuzzy bi-ideal of X.

**Problem 3.4.88:** Let X be a S-fuzzy multiplication bigroupoid. If $f : X \to T$ is an S-epimorphism of S-bigroupoids then prove T is a S-fuzzy multiplication bigroupoid.

**Problem 3.4.89:** If S(birad $\lambda$) is S-prime then show $\lambda$ is S-primary.

**Problem 3.4.90:** Let $\lambda$ be S-prime. Then prove for all positive integer n, $\lambda^n$ is S-primary and its S(birad) is $\lambda$.

**Problem 3.4.91:** Let $\lambda$ be S-prime then $\lambda^n \neq \lambda^{n-1}$ for all n > 0. Prove $\lambda^n$ is S-prime.

**Problem 3.4.92:** If $\lambda$ is S-primary; prove $\lambda = \mu^n$ for some positive integer n where $\mu$ = S(birad $\lambda$).

**Problem 3.4.93:** If $\lambda$ is a proper S-prime bi-ideal of X and $\mu$ is a S-fuzzy bi-ideal of a S-bigroupoid X such that $\mu \leq \lambda^n$, $\mu \leq \lambda^{n+1}$ for some n > 0 then prove $\lambda^n = \mu : y_r\chi_x$ where $y_r \notin \lambda$.

(Here $\lambda : \mu$ is defined for any two S-fuzzy bi-ideals of a S-bigroupoid X as $\lambda : \mu = \cup\{\delta | \delta$ is a S-fuzzy bi-ideal of X such that $\delta\mu \leq \lambda\}$; $\lambda : \mu$ is a S-fuzzy bi-ideal of X).

**Problem 3.4.94:** Prove every S-function generated fuzzy sub-bigroupoid is sub-bigroupoid generated.

**Problem 3.4.95:** Every sub-bigroupoid generated S-fuzzy sub-bigroupoid is isomorphic to a S-function generated fuzzy sub-bigroupoid. Prove.

**Problem 3.4.96:** Prove every S-sub-bigroupoid with respect to Min is S-sub-bigroupoid generated.

**Problem 3.4.97:** Let $\mu = \nu$ o f be the S-pre image of $\nu$ under f then $\mu$ is a S-fuzzy sub-bigroupoid of X with respect to T.

**Problem 3.4.98:** Let f : X $\to$ X' be a S-homomorphism from S-bigroupoid X onto S-bigroupoid X'. Suppose that $\mu$ is a S-fuzzy sub-bigroupoid of X with respect to a t-norm T and that $\mu$ has the S-sup property. Let $\nu$ be a S-homomorphic image of $\mu$ under f, then prove $\nu$ is a S-fuzzy sub-bigroupoid of X' with respect to T.

**Problem 3.4.99:** Is a S-penultimate sub-biloop of a S-fuzzy sub-biloop a proper S-sub-biloop of the given S-biloop L?

**Problem 3.4.100:** Let $\mu$ be a non constant S-fuzzy sub-biloop of a S-biloop L. Then, prove $\mu$ is S-penultimate sub-biloop, SP ($\mu$) is a proper S-sub-biloop of L if and only if $\mu$ attains its infimum.

**Problem 3.4.101:** Prove if $\eta$ and $\mu$ are S-fuzzy sub-biloops of a S-biloop L then the set product $\eta o\mu$ contains $\eta$ and $\mu$ if and only if $\eta$ and $\mu$ have the same S-bitip that is $\eta(e) = \mu(e)$ i.e. $\eta_1(e_1) \cup \eta_2(e_2) = \mu_1(e_1) \cup \mu_2(e_2)$ where $\eta = \eta_1 \cup \eta_2$ and $\mu = \mu_1 \cup \mu_2$.

**Problem 3.4.102:** Let $\eta$ and $\mu$ be S-fuzzy sub-biloops with the same bitip of a S-biloop L. Then prove the set product $\eta$ o $\mu$ is a S-fuzzy sub-biloop generated by the union of $\eta$ and $\mu$ if $\eta o\mu$ is a S-fuzzy sub-biloop of S-biloop.

**Problem 3.4.103:** Let $\eta$ and $\mu$ be S-fuzzy sub-biloops of a S-biloop L. Prove the S-penultimate product $\eta \circ \mu$ contains $\eta$ and $\mu$ if $\eta$ and $\mu$ are S- bisimilarly bounded.

**Problem 3.4.104:** Let $\eta$ and $\mu$ be S-fuzzy sub-biloops of a S-sub-biloop L. Then prove that the set product $\mu \circ \eta$ is a S-fuzzy sub-biloop if and only if $\eta \circ \mu = \mu \circ \eta$.

**Problem 3.4.105:** Let $\eta$ and $\mu$ be S-fuzzy sub-biloops of a S-biloop L. Prove the S-penultimate biproduct $\eta \odot \mu$ is a S-fuzzy sub-biloop if and only if $\eta \odot \mu = \mu \odot \eta$.

**Problem 3.4.106:** Prove a fuzzy subset $\lambda$ of a S-bigroupoid X is a S-fuzzy sub-bigroupoid of X if and only if $\lambda_t$ is a S-sub-bigroupoid of S for all $t \in (a, k]$

**Problem 3.4.107:** Let f be a S-bigroupoid epimorphism of X onto a S-bigroupoid T. Let $\lambda$ and $\mu$ be S-fuzzy sub-bigroupoids of X and T respectively. Then prove

  i. $f(\lambda)$ is a S-fuzzy sub-bigroupoid of T.
  ii. $f^{-1}(\mu)$ is a S-fuzzy sub-bigroupoid of X.

**Problem 3.4.108:** Prove a fuzzy subset $\lambda$ of a S-bigroupoid $X = X_1 \cup X_2$ is said to be a S-fuzzy regular bigroupoid of X if and only if $\lambda_t$ is a regular S-bigroupoid of X for all $t \in (a, k]$.

**Problem 3.4.109:** Prove a S-fuzzy left bi-ideal f of a S-bigroupoid X is S-prime if and only if for any two fuzzy points $x_r^1, y_t^1 \in P_1 \subset X_1$ and $x_r^2, y_t^2 \in P_2 \subset X_2$ ( $X = X_1 \cup X_2$) (rt > 0) $x_r^1 \circ P_1 \circ y_t^1 \circ P_1 \subset f_1$ (and $x_r^2 \circ P_2 \circ y_t^2 \circ P_2 \subset f_2$ here $f = f_1 \cup f_2$) implies $x_r^1 \in f_1, x_r^2 \in f_2$ or $y_t^1 \in f_2, y_t^2 \in f_2$.

**Problem 3.4.110:** Prove a S-fuzzy left bi-ideal f of a S-bigroupoid X is S-quasi prime if and only if for any two fuzzy points $x_r^1, y_t^1 \in X_1$ and $x_r^2, y_t^2 \in X_2$ (r t > 0)

$$x_r^1 \circ X_1 \circ y_t^1 \subset f_1 \text{ and } x_r^2 \circ X_2 \circ y_t^2 \subseteq f_2$$

implies that $x_r \in f$ or $y_t \in f$.

$$\left(\text{i.e. } y_t^1, x_r^1 \in f_1 \text{or } y_t^2, x_r^2 \in f_2\right).$$

**Problem 3.4.111:** Will the following statements be equivalent where X a S-commutative bigroupoid. F a S-fuzzy left bi-ideal of X.

  i. f is S-prime bi-ideal.
  ii. f is S-quasi prime bi-ideal.
  iii. f is S-weakly quasi prime.

# CHAPTER FOUR

# SMARANDACHE FUZZY RINGS AND NON-ASSOCIATIVE RINGS

In this chapter we introduce the notion of Smarandache fuzzy rings and Smarandache-fuzzy non-associative rings. This chapter has five sections. In the first section we introduce the concept of Smarandache fuzzy rings exhibited in 38 definitions and thirty theorems. In section two the notion of Smarandache fuzzy vector spaces are defined. Section three deals with Smarandache fuzzy non-associative rings. Section four deals with Smarandache fuzzy birings and several of its properties are defined and studied. In section five we have given 125 problems for the researcher. As this book is a research book we propose more problems.

## 4.1 Smarandache fuzzy rings definitions and properties

In this section we introduce the notion of Smarandache fuzzy rings. The study fuzzy rings [36, 37, 40, 109] is itself recent and still recent is the study of Smarandache rings [132]. So in this section for the first time we introduce the notion of Smarandache fuzzy rings. Here we give several examples and give some important properties about them.

**DEFINITION 4.1.1:** *Let $R$ be a S-ring. A fuzzy subset $\mu$ of $R$ is said to be a Smarandache fuzzy ring (S-fuzzy ring) relative to a subset $P$ of $R$ where $P$ is a field if $\mu: P \to [0, 1]$ is such that*

$$\mu(x-y) \geq \min\{\mu(x), \mu(y)\}$$
$$\mu(xy^{-1}) \geq \min\{\mu(x), \mu(y)\}$$
$$y \neq 0 \text{ for all } x, y \in P.$$

*This S-fuzzy ring will be called as type I; S-fuzzy ring. The condition of S-fuzzy ring I can be still tightened.*

**DEFINITION 4.1.2:** *Let $R$ be a S-ring. Let $\mu: R \to [0, 1]$ be a fuzzy subset of $R$. If $\mu$ is such that every proper subset $P_i \subset R$ say $i = 1, 2, ..., n$ where $P_i$ is a subfield of $R$; we have $\mu: P_i \to [0, 1]$ is a fuzzy subfield of $R$ for every $P_i$ ($i = 1, 2, ..., n$) then we call $\mu$ a Smarandache fuzzy strong ring of type I (S- fuzzy strong ring of type I).*

**THEOREM 4.1.1:** *Every S-fuzzy strong ring of type I is a S-fuzzy ring of type I.*

*Proof*: Direct by definitions, the converse is not true.

**THEOREM 4.1.2:** *A S-fuzzy ring of type I need in general be a S-fuzzy strong ring of type I.*

*Proof:* We prove this by an example. Let $Z_6 = \{0, 1, 2, 3, 4, 5\}$ be a S-ring. $P = \{0, 2, 4\}$ is a subset of $Z_6$ which is a field.

Define $\mu : Z_6 \to [0, 1]$. By
$$\mu(x) = \begin{cases} 0 & \text{if } x \in [0, 1, 3, 5] \\ 0.5 & \text{if } x = 2 \\ 1 & \text{if } x = 4. \end{cases}$$

Clearly μ restricted to P is a S-fuzzy ring I which is not S-fuzzy strong ring I as μ on the set $P_1 = \{0, 3\}$ which is a field in $Z_6$ and μ is trivial on $P_1$. Hence the claim.

We define Smarandache fuzzy ring of type II.

**DEFINITION 4.1.3:** *Let R be any ring $\mu: R \to [0, 1]$ be a fuzzy subset of R. We say μ is a Smarandache fuzzy subring of type II (S-fuzzy subring of type II) if there exist a fuzzy subset σ on R such that $\sigma \subset \mu$ and $\sigma: R \to [0, 1]$ is a fuzzy subfield of R. If for a given μ on a ring R no such σ exists such that $\sigma \subset \mu$ we say μ is a Smarandache non-fuzzy subring (S-non-fuzzy subring).*

*Even if $\sigma \subset \mu$ is a fuzzy subring of R still we call μ only a S-non-fuzzy ring.*

**DEFINITION 4.1.4:** *Let R be a S-ring. A fuzzy subset μ of the ring R is called a Smarandache fuzzy ideal (S-fuzzy ideal) of R if the following conditions are true*

   i. *R has a proper subset X such that*
      *$\mu(x - y) \geq \min\{\mu(x), \mu(y)\}$*
      *$\mu(xy) \geq \max\{\mu(x), \mu(y)\}$ for all $x, y \in X$. i.e. μ on x is a fuzzy ideal.*

   ii. *X contains a proper subset P such that P is a field under the operations of R and $\mu: P \to [0, 1]$ is a fuzzy subfield.*

*Example 4.1.1:* Let $Z_{12} = \{0, 1, 2, \ldots, 11\}$ be a S-ring under multiplication and addition modulo 12.

$X = \{0, 2, 4, 6, 8, 10\}$, $\mu : Z_{12} \to [0, 1]$

$$\mu(x) = \begin{cases} 0.1 & \text{if } x \text{ is odd} \\ 1 & \text{if } x = 4, 8, 0 \\ 0.5 & \text{if } x = 2, 6, 10. \end{cases}$$

It is easily verified μ is a S-fuzzy ideal of $Z_{12}$.

**Notation**: If μ and θ are S-fuzzy ideals of a ring R relative to the same ideal and fuzzy subfield; then the product μ o θ of μ and θ defined by

$$(\mu \circ \theta)(x) = \sup_{x = \sum x_i z_i} \left( \min_i \left( \min(\mu((y_i)), \theta(z_i)) \right) \right)$$

where $x, y_i, z_i \in X \subset R$, X a S-subring.

It can be easily verified that μ o θ is a S-fuzzy ideal of R.

**DEFINITION 4.1.5:** *A S-fuzzy ideal μ of a ring R a Smarandache fuzzy prime (S-fuzzy prime) if σ and θ are S-fuzzy ideals of R the condition $\sigma\theta \subseteq \mu$ implies either $\sigma \subset \mu$ or $\theta \subset \mu$.*

Now we proceed on to define Smarandache level subring.

**DEFINITION 4.1.6:** *Let μ be a S-fuzzy subring of a S-ring R; $t \in [0, 1]$ and $t \leq \mu(0)$. The S-fuzzy subring (ideal) $\mu_t$ is called a Smarandache level subring (S-level subring) of μ.*

**DEFINITION 4.1.7:** *Let μ by any fuzzy subset of a ring R. The smallest S-fuzzy ideal of R containing μ is called the S-fuzzy ideal generated by μ and is denoted by $\langle \mu \rangle$.*

**THEOREM 4.1.3:** *Let μ be any fuzzy subset of a S-ring R. Then $\langle \chi_\mu \rangle = \chi_{\langle \mu \rangle}$.*

*Proof:* Left as an exercise for the reader to prove.

**DEFINITION 4.1.8:** *Let μ be a fuzzy subset of a S-ring R. Then the fuzzy subset $\mu^*$ of R defined by $\mu^*(x) = \sup\{k \mid x \in \langle \mu_k \rangle\}$ is a S-fuzzy ideal generated by μ in R according as $\langle \mu_k \rangle$ is a S-fuzzy ideal generated by $\mu_k$ in R. i.e. $\mu^*(x) = t$ whenever $x \in \langle \mu_t \rangle$ and $x \notin \langle \mu_s \rangle$ for all $s > t$.*

**DEFINITION 4.1.9:** *Let μ be any S-fuzzy ideal of a S-ring R. The fuzzy subset $\mu_x^*$ of R where $x \in P \subset R$ (P a subfield of R) is defined by $\mu_x^*(r) = \mu(r - x)$ for all $r \in R$, is termed as the Smarandache fuzzy coset (S-fuzzy coset) determined by x and μ.*

**DEFINITION 4.1.10:** *A fuzzy subset λ of a S-ring R is said to be an Smarandache ($\in$, $\in \vee q$) fuzzy subring (S-($\in$, $\in \vee q$) fuzzy subring) of R if for all $x, y \in P \subset R$ (P is a subfield relative to which λ is defined) and $t, r \in (0, 1]$*

  i. $x_t, y_r \in \lambda \Rightarrow (x + y)_{M(t, r)} \in \vee q \lambda$.
  ii. $x_t \in \lambda \Rightarrow (-x)_t \in \vee q \lambda$.
  iii. $x_t, y_r \in \lambda \Rightarrow (x y)_{M(t, r)} \in \vee q \lambda$.

**THEOREM 4.1.4:** *A fuzzy subset λ of a S-ring R is an S-($\in$, $\in \vee q$) fuzzy subring of R if and only if $\lambda(x - y), \lambda(xy) \geq M(\lambda(x), \lambda(y), 0.5)$ for all $x, y \in X \subset R$, X a proper subset of R. which is a field.*

*Proof:* Direct by the definitions.

**DEFINITION 4.1.11:** *A fuzzy subset λ of a ring R is said to be a Smarandache ($\in$, $\in \vee q$)-fuzzy ideal (S-($\in$, $\in \vee q$)-fuzzy ideal) of R if*

  i. *λ is an S-($\in$, $\in \vee q$) fuzzy subring of R.*
  ii. $x_t \in \lambda$ and $y \in R \Rightarrow (xy)_t, (yx)_t \in \vee q \lambda$.

*(Just we once again recall the notation a fuzzy point, $x_t$ is said to belong to (resp be quasi coincident with) a fuzzy set $\lambda$ written as $x_t \in \vee \lambda$ (resp $x_t q \lambda$) if $\lambda(x) \geq t$ (resp. $\lambda(x) + t > 1$), $x_t \in \lambda$ or $x_t q \lambda$ will denoted by $x_t \in \vee q \lambda$) [25].*

In view of these definition the reader is expected to prove.

**THEOREM 4.1.5:** *A fuzzy subset $\lambda$ of a S-ring R is an $(\in, \in \vee q)$-fuzzy ideal of a S-ring R if and only if*

   i. $\lambda(x - y) \geq M(\lambda(x), \lambda(y), 0.5)$.
   ii. $\lambda(xy), \lambda(yx) \geq M(\lambda(x), 0.5)$ for all $x, y \in X \subset R$.

(X a proper S-subring of R i.e. X contains a subfield as a proper subset).

**DEFINITION 4.1.12:** *An $S$-$(\in, \in \vee q)$ fuzzy ideal of a S-ring R is said to be*

  i. *Smarandache $(\in, \in \vee q)$ fuzzy semiprime (S-$(\in, \in \vee q)$-fuzzy semiprime) if for all $x, y \in X \subset R$ and $t \in (0,1]$, $(x^2)_t \in \lambda \Rightarrow x_t \in \vee q \lambda$.*

  ii. *Smarandache $(\in, \in \vee q)$-fuzzy prime (S-$(\in, \in \vee q)$-fuzzy prime) if for all $x, y \in X \subset R$ and $t \in (0, 1]$ $(xy)_t \in \lambda \Rightarrow x_t \in \vee q \lambda$ or $y_t \in \vee q \lambda$.*

  iii. *Smarandache $(\in, \in \vee q)$-fuzzy semiprimary (S-$(\in, \in \vee q)$-fuzzy semiprimary) if for all $x, y \in X \subset R$ and $t \in (0, 1]$, $(xy)_t \in \lambda$ implies $x_t^n \in \vee q \lambda$ or $y_t^m \in \vee q \lambda$ for some $n, m \in N$.*

  iv. *Smarandache $(\in, \in \vee q)$-fuzzy primary (S-$(\in, \in \vee q)$-fuzzy primary) if for all $x, y \in X \subset R$ and $t \in (0, 1]$, $(xy)_t \in \lambda \Rightarrow x_t \in \vee q \lambda$ or $y_t^n \in \vee q \lambda$ for some $n \in N$.*

The following three theorems are left as an exercise for the reader.

**THEOREM 4.1.6:** *A fuzzy subset $\lambda$ of a S-ring R is a S-$(\in, \in \vee q)$ fuzzy subring (ideal) of R if and only if $\lambda_t$ is a S-subring (S-ideal) of R for all $t \in (0, 0.5]$.*

**THEOREM 4.1.7:** *A S-$(\in, \in \vee q)$ fuzzy ideal of R is S-$(\in, \in \vee q)$ fuzzy prime if and only if Max $\{\lambda(x), \lambda(y)\} \geq M(\lambda(xy), 0.5)$ for all $x, y \in X \subset R$.*

**THEOREM 4.1.8:** *A S-fuzzy ideal $\lambda$ of a S-ring R is a S-$(\in, \in \vee q)$ fuzzy semiprime (or S-prime or S-semiprimary or S-primary) if and only if $\lambda_t$ is S-semiprime (or S-prime or S-semiprimary or S-primary) for all $0 < t \leq 0.5$.*

Now we proceed on to define Smarandache $(\in, \in \vee q)$ fuzzy radical of $\lambda$.

**DEFINITION 4.1.13:** *Let $\lambda$ be an S-$(\in, \in \vee q)$ fuzzy ideal of R. The fuzzy subset Rad $\lambda$ of R is defined by*

$$(S\ Rad\ \lambda)\ (x) = \begin{cases} M(\sup \lambda(x^n); n \in N\}, 0.5\ )\ if\ \lambda(x) < 0.5 \\ \lambda(x)\ if\ \lambda(x) \geq 0.5 \end{cases}$$

*is called the Smarandache $(\in, \in \vee q)$ fuzzy radical (S-$(\in, \in \vee q)$ fuzzy radical) of $\lambda$.*

Using the definitions and results of [23] the reader is expected to prove.

**THEOREM 4.1.9:** *Let $\lambda$ be an S-$(\in, \in \vee q)$ fuzzy ideal of the S-ring R. Then*

  i.      *S Rad $\lambda$ is an S-$(\in, \in \vee q)$ fuzzy ideal of R.*
  ii.     *If $\lambda$ has the sup-property then S Rad $\lambda_t$ = (S Rad $\lambda)_t$ for all $0 < t \leq 0.5$.*

**DEFINITION 4.1.14:** *The subset $\lambda_t = \{x \in X\ |\ \lambda(x) \geq t\ or\ \lambda(x) + t > 1\} = \{x \in X, x_t \in \vee q\lambda\}$ is called $(\in \vee q)$-level subset of X determined by $\lambda$ and t.*

**DEFINITION 4.1.15:** *A Smarandache L-fuzzy ring (S-L fuzzy ring) is a function $\mu : R \to L$, L a distribution lattice where (R, +, •) is a S-ring that satisfies,*

  i.    *$\mu \not\equiv 0$.*
  ii.   *$\mu(x - y) \geq \mu(x) \wedge \mu(y)$ for every $x, y \in P \subset R$*
        *(P a subfield relative to which $\mu$ is defined).*
  iii.  *$\mu (xy) \geq \mu (x) \wedge \mu (y)$ for every $x, y \in P \subset R$.*
  iv.   *If R is unitary then $\mu(1) = \mu(0)$.*

*$S(\mu)$ is a S-subring of R where $\mu|_S$ denotes the restriction of $\mu$ to $S(\mu)$. Let $\mu: R \to L$ and $\mu': R' \to L$ are S-L-fuzzy ring (R and R' S-ring).*

*A Smarandache homomorphism between $\mu$ and $\mu'$ is a S-ring homomorphism $f: R \to R'$ that satisfies $\mu'(f(x)) = \mu(x)$ for all x in $P \subset R$ (P is a subfield relative to which $\mu$ is defined) i.e. $f^{-1}(\mu') = \mu$.*

**DEFINITION 4.1.16:** *A Smarandache L-fuzzy subring (S-L-fuzzy subring) of a S-L-fuzzy ring; $\mu$ is a S-fuzzy ring $\mu': R \to L$ (R a S-ring) satisfying $\mu'(x) \leq \mu(x)$ for all $x \in P \subset R$ (where $\mu$ is defined relative to the subfield P in R).*

*For $\mu'$ a S-fuzzy subring of the S-L-fuzzy ring $\mu$ with L totally ordered. We shall call $\mu'$ the restriction of $\mu'$ to $S(\mu)$. It can be easily proved $\mu'|_S$ is a S-fuzzy subring of $\mu|_S$.*

**DEFINITION 4.1.17:** *Let $\mu: R \to L$ be a S-L-fuzzy ring. A Smarandache L-fuzzy ideal (S-L-fuzzy ideal) of $\mu$ is a map $\delta: R \to L$ such that the following properties hold.*

  i.      *$\delta \not\equiv 0$.*
  ii.     *$\delta(xy) \geq \delta(x) \wedge \delta(y)$ for every $x, y \in P \subset R$.*
  iii.    *$\delta(xy) \geq \mu(x) \wedge \delta(y)$ for every $x, y \in P \subset R$.*
  iv.     *$\delta(x) \leq \mu(x)$ for every $x \in P \subset R$.*

*(P is the subfield relative to which the S-L fuzzy ring $\mu$ is defined).*

*As in case of S-L-fuzzy subring, for $\delta$ a S-L-fuzzy ideal of the S-L-fuzzy ring $\mu$ with L totally ordered, $\delta|_S$ the restriction of $\delta$ to $S(\mu)$. It holds that $\delta|_S$ is a S-L-fuzzy ideal of $\mu|_S$ where $S(\mu)$ is a S-subring of R, $\mu|_S$ denotes the restriction of $\mu$ to $S(\mu)$.*

**DEFINITION 4.1.18:** *Let $\mu$ be a fuzzy subset and $t \in L$. The t-level cut $\mu_t$ will be $\{x \in R \mid \mu(x) \geq t\}$ and the Smarandache t-level cut (S-t-level cut) $S\mu_t$ will be $\{x \in P \subset R \mid \mu(x) \geq t$; P is the subfield of the S-ring R relative to which $\mu$ is defined$\}$ Analogously the t-level strong cut will be $\{x \in R \mid \mu(x) > t\}$ and Smarandache t-level strong cut (S-t-level strong cut) $S\mu_t$ will be $\{x \in P \subset R / \mu(x) > t$ – P defined above$\}$ and finally $\mu^0$ will stand for $\{x \in R / \mu(x) = \mu(0)\}$ and $S\mu^0 = \{x \in P \subset R / \mu(x) = \mu(0)$, '0' the additive identity of P$\}$.*

**DEFINITION 4.1.19:** *Let S be a $\delta$-L-fuzzy ideal of $\mu$. We define Smarandache fuzzy quotient ring ( S-fuzzy quotient ring ) $\mu/s \mid \delta/s$ as $\mu/s \mid \delta/s : S(\mu)|_\alpha \to L$ with $\mu/s \mid \delta/s (x + \alpha) = \mu(x)$ where $\alpha$ is the S-ideal of $S(\mu)$, S-maximal among those contained in $S\delta^0 \cap S\mu^0$.*

Now the notion of Smarandache prime L-fuzzy ideals is defined.

**DEFINITION 4.1.20:** *A S-L-fuzzy ideal $\delta$ of a S-L-fuzzy ring $\mu$ is said to Smarandache prime (or Smarandache L-prime) if $\delta \not\equiv 1$ and it satisfies $\delta(xy) \wedge \mu(x) \wedge \mu(y) < \delta(x) \vee \delta(y)$ for every $x, y \in P \subset R$ (P a subfield of R relative to which $\mu$ is defined).*

*Clearly if $\mu \not\equiv 1$ then we have $\delta(xy) \leq \delta(x) \vee \delta(y)$ for every $x, y \in P \subset R$. Further a S-L-fuzzy ideal $\delta$ of a S-L-fuzzy ring R is S-L-prime if for every $x, y \in P \subset R$, $\delta(xy) = \delta(0)$ implies $\delta(x) = \delta(0)$ or $\delta(y) = \delta(0)$. We see if L is totally ordered then it is evident that $\delta(xy) \leq \delta(x) \vee \delta(y)$ implies $\delta(xy) = \delta(o)$ implies $\delta(x) = \delta(0)$ or $\delta(y) = \delta(0)$.*

*A Smarandache strongly prime L-fuzzy ideal (S-strongly prime L-fuzzy ideal) is defined as a non-constant S-fuzzy ideal $\delta$ of a S-ring $\mu$ that satisfies $\delta(xy) = \delta(x)$ or $\delta(xy) = \delta(y)$ for all $x, y \in P \subset R$.*

*A Smarandache weakly completely fuzzy prime ideal (S-weakly completely fuzzy prime ideal) is defined when $L = [0, 1]$ as a non constant S-fuzzy ideal $\delta$ of R that satisfies $\delta(xy) = \max \{\delta(x), \delta(y)\}$ for every $x, y \in P \subset R$. (P is a subfield relative to which $\delta$ is defined).*

**Note**: This notion can be easily generalized to the case in which we have as arbitrary lattice L, changing maximum by supremum.

**DEFINITION 4.1.21:** *A S-fuzzy ring $\mu$ is said to an Smarandache fuzzy integral domain (S-fuzzy integral domain) if $x, y = 0$ and $\mu(x) \wedge \mu(y) > 0$ implies $x = 0$ or $y = 0$, $x, y \in P \subset R$ (P is the subfield of the S-ring R relative to which $\mu$ is defined).*

**DEFINITION 4.1.22:** *Let $\delta$ be a S-L-fuzzy ideal of a S-L-fuzzy ring $\mu$. The Smarandache L-fuzzy radical (S-L-fuzzy radical) of $\delta$ is the fuzzy subset of R defined* by $Sr(\delta)(x) = \vee_{n \in N} \delta(x^n) \wedge \mu(x)$, when $\mu = 1$, definition transforms $Sr(\delta)(x) = \vee_{n \in N} \delta(x^n)$, $Sr(\delta)(x) = \sup \{t \mid x \in Sr(\delta_t)\}$.

The following theorems are left as exercises for the reader to prove.

**THEOREM 4.1.10:** *Let $\delta$ be a S-L-fuzzy ideal of $\mu$. Then $Sr(\delta)$ is S-L-fuzzy ideal of $\mu$.*

**THEOREM 4.1.11:** *Let $\delta$ be a S-L-fuzzy ideal of a S-L-fuzzy ring $\mu$. Then $Sr(r(\delta)) = Sr(\delta)$.*

**THEOREM 4.1.12:** *Let $\delta$ be a S-prime L-fuzzy ideal of $\mu$ then $Sr(\delta) = \delta$.*

**THEOREM 4.1.13:** *Let $\delta$ be a S-L-fuzzy ideal of $\mu$ then $\delta \subseteq Sr(\delta)$.*

**THEOREM 4.1.14:** *Let $\delta_1, \delta_2$ be S-L-fuzzy ideal of a S-fuzzy ring $\mu$*

  i. *If $\delta_1 \subseteq \delta_2$ then $Sr(\delta_1) \subset Sr(\delta_2)$.*
  ii. *$Sr(\delta_1 \cap \delta_2) = Sr(\delta_1) \cap Sr(\delta_2)$.*
  iii. *If L is totally ordered then*

$$Sr(\delta) = \bigcap_{\substack{\delta' \delta-\text{prime} \\ \delta \subseteq \delta'}} \delta'.$$

Now we proceed on to define Smarandache-primary L-fuzzy ideals.

**DEFINITION 4.1.23:** *A S-L-fuzzy ideal $\delta$ is said to be Smarandache L-primary (S-L-primary) if $\delta(xy) \wedge \mu(x) \wedge \mu(y) \leq \delta(x) \vee (\{\vee_{n \in N} \delta(y^n)\}$ for all $x, y \in X \subset R$ (X a proper subset of R which has a proper subset which is a subfield relative to which $\delta$ is defined).*

In case L is a totally order set i.e. L = [0, 1] then we have S-L-primary will coincide with the notion of S-primary.

**DEFINITION 4.1.24:** *A S-fuzzy ideal $\lambda$ of a S-ring R is called*

  i. *Smarandache semiprime (S-semiprime) if $\lambda(x^2) = \lambda(x)$ for all $x \in P \subset R$ (P is a subfield of R relative to which $\lambda$ is defined).*

  ii. *Smarandache prime (S-prime) if $\lambda(xy) = \lambda(x)$ or $\lambda(xy) = \lambda(y)$ for all $x, y \in P \subset R$.*

  iii. *Smarandache semiprimary (S-semiprimary) if for all $x, y \in P \subset R$, $\lambda(xy) \leq \lambda(x^n)$ for some positive integer n or $\lambda(xy) \leq \lambda(y^n)$ for some positive integer n.*

  iv. *Smarandache primary (S-primary) if for all $x, y \in P \subset R$, $\lambda(xy) = \lambda(x)$ or $\lambda(y^n)$ for some positive integer n.*

v.   Maximal if the level S-ideal $\{ x \in P \subset R \mid \lambda (x) = 1\}$ is maximal.

**DEFINITION 4.1.25:** *Let $\lambda$ be a S-fuzzy ideal of a S-ring R. The fuzzy subset S rad $\lambda$ of R defined by $S(Rad\ \lambda)(x) = \sup \{\lambda (x^n) \mid n > 0\}$ is called Smarandache fuzzy radical (S-fuzzy radical) of $\lambda$.*

Several inter-relations between S-fuzzy prime, S-fuzzy primary, S-fuzzy semiprimary etc can be obtained as a matter of routine as in of fuzzy rings.

**DEFINITION 4.1.26:** *A S-fuzzy ideal $\mu$ of a S-ring R is called Smarandache fuzzy irreducible (S-fuzzy irreducible) if it is not a finite intersection of two S-fuzzy ideals of R properly containing $\mu$. If S-fuzzy ideal $\mu$ of a S-ring R is called Smarandache fuzzy weakly irreducible (S-fuzzy weakly irreducible) if it not a finite intersection of a S-fuzzy ideal and a fuzzy ideal of R properly containing $\mu$ otherwise $\mu$ is termed as Smarandache fuzzy reducible (S-fuzzy reducible).*

**THEOREM 4.1.15:** *If $\mu$ is any S-fuzzy prime ideal of a S-ring R, then $\mu$ is S-fuzzy irreducible.*

*Proof:* Straightforward as in case of fuzzy prime ideals.

It is left for the reader as an exercise to prove the following theorem.

**THEOREM 4.1.16:** *Let $\mu$ be a S-fuzzy ideal of a S-ring R. If for some $x, y \in P \subset R$ (P is the subfield relative to which $\mu$ is defined) $\mu(x) < \mu (y)$ then $\mu (x- y) = \mu(x) = \mu(y - x)$.*

Just for the sake of completeness we recall the definition of coset of a fuzzy ideal and its Smarandache analogue.

**DEFINITION 4.1.27:** *Let $J: R \to [0, 1]$ be a S-fuzzy ideal where R is a S-ring. The fuzzy subset $x + J : R \to L = [0, 1]$ defined by $(x + J)y = J(y - x)$ is called a Smarandache coset of the S-fuzzy ideal (S-coset of the S-fuzzy ideal) J.*

**THEOREM 4.1.17:** *If $J : R \to [0, 1]$ is a S-fuzzy ideal, then $x + J = y + J$ if and only if $J(x - y) = J(0)$. In that case $J(x) = J(y)$.*

*Proof:* As in case of fuzzy ideal. Hence left for the reader to prove.

**THEOREM 4.1.18:** *Every S-coset of a fuzzy ideal $J : R \to [0, 1]$ is constant on every S-coset of a (non-fuzzy) ideal $R_J$ of R. More specifically $x + J (z) = J(y - x)$ for all $z \in y + R_J$. In particular $x + J (z) = J (x)$ for all $x \in R_J$.*

*Proof:* As in case of fuzzy ideals the proof is straightforward.

**DEFINITION 4.1.28:** *The ring R/J of the S-cosets of the S-fuzzy ideal J is called the Smarandache factor ring (S-factor ring) or the Smarandache quotient ring (S-quotient ring) of R by J.*

Let SF denote the set of all Smarandache fuzzy ideals of a S-ring R.

**DEFINITION 4.1.29:** *Let $A \in SF$. A representation of A is a finite intersection $A = Q_1 \cap Q_2 \cap ... \cap Q_n$ of S-fuzzy primary ideals of R is called a Smarandache fuzzy primary representation (S-fuzzy primary representation) of A. If is called irredundant or reduced if no $Q_i$ contains*

$$\bigcap_{\substack{j=1 \\ i \neq j}}^{n} Q_j$$

*and the $Q_i$ have distinct S-radicals, $i = 1, 2, ..., n$.*

Now we proceed on to define Smarandache fuzzy prime ideal divisor of A.

**DEFINITION 4.1.30:** *Let $\mu \in SF$. A S-fuzzy prime ideal $\sigma$ of R is called a Smarandache fuzzy prime ideal divisor (S-fuzzy prime ideal divisor) of $\mu$ if $\sigma \supseteq \mu$ and $\sigma_* \supseteq \mu_*$.*

*A Smarandache fuzzy prime ideal divisor $\sigma$ of $\mu$ is called minimal or isolated if there does not exist a S-fuzzy prime ideal divisor $\sigma'$ of $\mu$ such that $\sigma' \subseteq \sigma$ and $\sigma' \neq \sigma$.*

**DEFINITION 4.1.31:** *Let $\mu$ be a S-fuzzy subring of a S-ring R. We say that $\mu$ is a Smarandache fuzzy quasi-local subring (S-fuzzy quasi-local subring) of R if and only if for all $x \in R$ and for all $y \in R$ such that y is a unit.*

$$\mu(xy^{-1}) \geq \min\{\mu(x), \mu(y)\} \text{ or equivalently } \mu(y) = \mu(y^{-1}).$$

The following theorem is straightforward hence left for the reader to prove.

**THEOREM 4.1.19**: *Let R be a S-quasi local ring. If $\mu$ is a S-fuzzy quasi local subring of R, then for all t such that $0 \leq t \leq \mu(1)$. $\mu_t$ is a S-quasi local ring and $M \cap \mu_t$ is the unique maximal ideal of $\mu_t$. (M denotes a unique S-maximal ideal of R).*

*Proof*: $\mu_t$ is a S-fuzzy subring of R. Using $\mu$ is S-quasi local $\mu(x) = \mu(x^{-1})$ leading to $x^{-1} \in \mu_t$.

**Notation**: For the sake of completeness we recall that $\mu^* = \{x \in P \subset R \mid \mu(x) > 0\}$ is a S-subring of the S-ring R. Let R be a non empty subset of a set S and let A be a fuzzy subset of R. If B is an extension of A to a fuzzy subset of S, then $A_t \cap B_s = A_s$ for all s, t such that $0 \leq t \leq s \leq 1$. A is an extension to a fuzzy subset $A^e$ of S such that $(A^e)_t \supseteq B_t$ for all $t \in \text{Im}(A)$.

$A^e(y) = \sup\{t \mid y \in B_t\}$ where $B_t \in B = \{B_t \mid t \in \text{Im}(A)\}$, $A^e(y) = \max\{m \mid y \in B_m\}$.

Let I be an ideal of a ring R, S be a subring of R. $I^e$ denotes the ideal of S generated by I. If A is a fuzzy ideal of R; S a subring of R, such that A has sup property. If

$$\bigcup_{t \in \text{Im}(A)} (A_t)^e = S$$

and for all s, t ∈ Im (A) , s ≥ t. $A_t \cap (A_s)^e = A_s$ then A has a unique extension to a fuzzy ideal $A^e$ of S such that $(A^e)_t = (A_t)^e$ for all t ∈ Im (A) and Im ($A^e$) = Im (A).

Further $A^e(xy^{-1}) \geq$ min {$A^e(x)$, $A^e(y)$} if and only if for all s, t ∈ Im (A) , s ≥ t. Further A can be extended to a fuzzy quasi local subring $A^e$ of $R_P$ if and only if for all s, t ∈ Im (A), s ≥ t $A_t \cap (A_s)_{P_s} = A_s$ where $P_s = P \cap A_s$ for all s ∈ Im (A). Let R and S be S-rings, f a S-ring homomorphism of R onto S. Let T denote f (R).

**DEFINITION 4.1.32:** *Let A and B be fuzzy subsets of R and T respectively. Define the fuzzy subsets f (A) of T and $f^{-1}(B)$ of R by f(A) (y) = Sup {A(x) | f(x) = y} for all y ∈ T. $f^{-1}(B)$ (x) = Bf(x)) for all x ∈ R.*

**DEFINITION 4.1.33:** *Suppose A and B are S-fuzzy ideals of R and T respectively. Then we have f(A) and $f^{-1}(B)$ are S-fuzzy ideals of T and R respectively, with f(A)(0) = A(0), $f^{-1}(B)(0) = B(0)$; where R and S are S-rings, f a S-ring homomorphism of R into S; T denotes f (R). If I is an S-ideal of R , then the ideal f (I) $^e$ or $I^e$ is defined to be the S-ideal of S degenerated by f (I) and it is called the Smarandache extended ideal (S-extended ideal) or Smarandache extension (S-extension) of I. If J is an S-ideal of S, the S-ideal $J^c = f^{-1}(J)$ is called the Smarandache contracted ideal (S-contracted ideal) or the Smarandache contraction (S-contraction) of J.*

The following result is straightforward from the very definitions.

**THEOREM 4.1.20:** *Let A be a S-fuzzy ideal of the S-ring R. The*

   i.   $f(A_*) \subseteq f(A)_*$
   ii.  *If A has the sup property then $f(A)_* = f(A_*)$.*

We proceed on to define the notion of Smarandache f-invariant ideal.

**DEFINITION 4.1.34:** *Let A be a S-fuzzy ideal of S-ring R. A is called Smarandache f-invariant (S-f-invariant ideal) if and only if for all x ∈ P ⊂ R, f(x) = f(y) implies A(x) = A(y) . (P the subfield relative to which A is defined).*

The following theorems are straightforward; hence left for the reader to prove.

**THEOREM 4.1.21:** *Let A be a S-fuzzy ideal of R. Then A is a S-fuzzy ideal of R if and only if A(0) = 1, |Im(A)| = 2 and $A_*$ is a prime ideal of R.*

**THEOREM 4.1.22:** *Let A be an Smarandache f-invariant fuzzy ideal of the S-ring R , such that A has the sup property. If $A_*$ is a S-prime ideal R, then f ($A_*$) is a S-prime ideal of T.*

**THEOREM 4.1.23:** *Let A be a S-f-invariant fuzzy ideal of the S-ring R such that Im (A) is finite.*

*If $A_*$ is a S-prime ideal of R then $f(A_*) = f(A)_*$*

**THEOREM 4.1.24:** *Let A be an S-f-invariant fuzzy ideal of R. If A is a S-fuzzy prime ideal of R then f(A) is a S-fuzzy prime ideal of T.*

(Here f : R → S; R and S, S-rings and f is a S-ring homomorphism and f(S) = T).

Several other interesting results in this direction can be obtained using the results of [84]. The innovative reader is expected to find several Smarandache analogous of these results and prove the above theorems; as the proof of these theorems can be done as a matter of routine using mainly the definitions.

Throughout this section $\Omega$ denote a non empty set we take R to be always a S-ring and all fuzzy ideals are S-fuzzy ideals of R.

**DEFINITION 4.1.35:** *Let $\{A_\alpha / \alpha \in \Omega\}$ be a collection of fuzzy subsets of a S-ring R. Define the fuzzy subset $\sum_{\alpha \in \Omega} A_\alpha$ of R and for all $x \in R$ by*

$$\left(\sum_{\alpha \in \Omega} A_\alpha\right)(x) = \sup\{\inf\{\{A_\alpha(x_\alpha) / \alpha \in \Omega\} \mid x = \sum_{\alpha \in \Omega} A_\alpha\}.$$

*For $x, x_\alpha \in R$, $x = \sum_{\alpha \in \Omega} x_\alpha$.*

**DEFINITION 4.1.36:** *Let $\{A_\alpha / \alpha \in \Omega\}$ be a collection of all S-fuzzy subrings (S-ideals) of a S-ring R. Then $\sum_{\alpha \in \Omega} A_\alpha$ is a S-fuzzy subring (S-fuzzy ideals) of R and $A_\alpha \subseteq \sum_{\alpha \in \Omega} A_\alpha$ for all $\alpha \in \Omega$. Let $\{A_\alpha / \alpha \in \Omega\} \cup \{A\}$ be a collection of S-fuzzy subrings of R. Then A is said to be the Smarandache fuzzy weak direct sum (S-fuzzy weak direct sum) of the $A_\alpha$ if and only if $A = \sum_{\alpha \in \Omega} A_\alpha$ and for all $x \in X \subset R$ ($x \neq 0$, A defined relative to X in R), $x \neq 0$, $\left(A_\beta \cap \sum_{a \in \Omega\beta} A_\alpha\right)(x) = 0$. If A is the S-fuzzy weak direct sum of $A_\alpha$ then we write $A = \bigoplus_{\alpha \in \Omega} A_\alpha$ We have $A^* = \{x \in R \mid A(x) > 0\}$ and $A_* = \{x \in R \mid A(x) = A(0)\}$.*

**THEOREM 4.1.25:** *Let A be S –fuzzy subring (S-fuzzy ideal) of a S-ring R. Then $A_*$ is a S-subring (S-ideal) of R. If L has the finite intersection property then $A^*$ is a S-subring (S-ideal) of R.*

*Proof*: Follows by the very definitions, hence left for the reader to prove.

**Notation**: Let A and B be fuzzy subsets of a S-ring R. Define the fuzzy subset AB of R for all $x \in R$, by $(AB)(x) = \sup\{\inf(\inf\{A(y_i), B(z_i)\}\})$ such that i = 1, 2, …, n |

$$x = \sum_{i=1}^{n} y_i z_i, n \in N\}.$$

Using this notation prove the following theorem:

**THEOREM 4.1.26:** *Let A and B be S-fuzzy subrings (S-fuzzy ideals) of R. Then AB is a S-fuzzy subring (S-fuzzy ideal) of R.*

**THEOREM 4.1.27:** *Let $\{A_\alpha / \alpha \in \Omega\}$ be a collection of S-fuzzy subrings of R. (i.e. $A_\alpha$: R $\to$ L where L has the finite intersection property $A_\alpha$'s defined relative to a fixed subfield X of R. Then*

$$\sum_{\alpha \in \Omega} A_\alpha^* = \left( \sum_{\alpha \in \Omega} A_\alpha \right)^*.$$

*Proof:* $x \in \left( \sum_{\alpha \in \Omega} A_\alpha \right)^*$ if and only if $\left( \sum_{\alpha \in \Omega} A_\alpha \right)(x) > 0$ if and only if $\sup\{\inf\{A_\alpha(x_\alpha) \mid \alpha \in \Omega\} \mid x = \sum_{a \in \Omega} x_\alpha > 0$ if and only if $x = \sum_{a \in \Omega} x_\alpha$ for sum $x_\alpha \in A_\alpha^*$ if and only if $x \in \sum_{\alpha \in \Omega} A_\alpha^*$.

Prove on similar lines.

**THEOREM 4.1.28:** *Let $\{A_\alpha / \alpha \in \Omega\}$ be a collection of S-fuzzy subrings of R; where $A_\alpha: R \to L$, L has intersection property, then prove.*

$$\bigcap_{\alpha \in \Omega} A_\alpha^* = \left( \bigcap_{\alpha \in \Omega} A_\alpha \right)^*$$

*(all $A_\alpha$'s defined relative to one fixed subfield).*

**THEOREM 4.1.29:** *Let $\{A_\alpha / \alpha \in \Omega\}$ be a collection of S-fuzzy subrings of R relative to a complete distributive lattice L with finite intersection property. Then for all $\beta \in \Omega$.*

$A_\beta^* \cap \sum_{\alpha \in \Omega_\beta} A_\alpha^* = \{0\}$ *if and only if for all* $x \in X \subset R; x \neq 0$, $\left( A_\beta \cap \sum_{\alpha \in \Omega_\beta} A_\alpha \right)(x) = 0.$

*Proof:* Straightforward using the definitions.

Now we define Smarandache pairwise co maximal fuzzy ideals of a S-ring with identity.

**DEFINITION 4.1.37:** *Suppose that R is a S-ring with 1. Let $A_\alpha$, $\alpha \in \Omega$ be S-fuzzy ideals of R. $A_\alpha$ are said to be Smarandache pairwise co-maximal (S-pairwise co-maximal) if and only if $A_\alpha \neq \delta_R$ for all $\alpha \in \Omega$ and $A_\alpha + A_\beta = \delta_R$ for all $\alpha, \beta \in \Omega, \alpha \neq \beta$.*

It is left for the reader to prove using the definitions.

**THEOREM 4.1.30:** *Suppose that R is a S-ring with identity. Let $\{A_\alpha \mid \alpha \in \Omega\}$ be a collection of finite valued S-fuzzy ideals of R. $A_\alpha$ are S-pairwise co-maximal if and only if $A_{\alpha*}$ are S- pairwise co-maximal.*

We define the notion of Smarandache complete fuzzy direct sum.

**DEFINITION 4.1.38:** *Let $\{R_\alpha \mid \alpha \in \Omega\}$ be a collection of S-commutative ring and let $A_\alpha$ be a S-fuzzy subring of $R_\alpha$ for all $\alpha \in \Omega$ then the Cartesian cross-product $\bigtimes_{\alpha \in \Omega} A_\alpha$ is called Smarandache complete fuzzy direct sum (S-complete fuzzy direct sum) of the $A_\alpha$. Define the fuzzy subset $\sum_{\alpha \in \Omega}^{\oplus} A_\alpha$ of $\sum_{\alpha \in \Omega}^{\oplus} R_\alpha$ by for all $x \in_\alpha \langle x_\alpha \rangle \in \sum_{\alpha \in \Omega}^{\oplus} R_\alpha$ . $\sum_{\alpha \in \Omega}^{\alpha} A_\alpha (x)$ = $\bigtimes_{\alpha \in \Omega} A_\alpha (x_\alpha)$. Then $\sum_{\alpha \in \Omega}^{\oplus} A_\alpha$ is called the Smarandache weak fuzzy direct sum (S-weak fuzzy direct sum) of the $A_\alpha$.*

Several interesting results in this direction can be obtained by any studious reader. Here it is pertinent to mention that the definition of several fuzzy concepts are defined in different ways by different researchers hence in this text one may find that we would have at many instances recalled those definitions. Also in some places we would have recalled the definition it is a pertinent repetition, mainly for easy reading.

Finally we mention that our study is not totally exhaustive about all fuzzy algebraic structures. We give here only those fuzzy algebraic structures for which we are in a position to obtain a Smarandache analogue. Every chapter on Smarandache notions ends with a set of problems several of them are routine theorems and some of them are research problems but an innovative researcher can certainly work on them.

### 4.2. Smarandache Fuzzy vector spaces and its properties

The study of Smarandache fuzzy vector spaces is new as the notion of Smarandache vector spaces is meager. A little insight about S-vector spaces can be had from [135]. In this section we define S-fuzzy vector spaces and its properties. Nearly 15 definitions about fuzzy vector spaces are given in this section.

**DEFINITION [135]:** *The Smarandache K-vectorial space (S-K-Vectroial space) is defined to be a S-vectorial space $(A, +, \bullet)$ such that a proper subset of A is a K-algebra (with respect to the same induced operation and another '$\times$' operation internal on A where K is a commutative field.*

**DEFINITION [135]:** *Let A be a K-vectorial space. A proper subset X of A is said to be a Smarandache K-vectorial subspace (S-K-vectorial subspace) of A if X itself is a S-K-vectorial space.*

**DEFINITION [135]:** *Let V be a finite dimensional vector space over a field K. Let B = $\{v_1, ..., v_n\}$ be a basis of V. We say B is a Smarandache basis (S-basis) of V if B has a proper subset say A, $A \subset B$ and $A \neq \phi$, $A \neq B$ such that A generates a subspace which* is linear algebra over K that is W is the subspace generated by A then W must be a K-algebra with the same operation of V.

**THEOREM [135]:** *Let A be a K-vectorial space. If A has a S-K-vectorial subspace then A is a S-K-vectorial space.*

*Proof*: Straightforward by the very definition.

**THEOREM [135]:** *Let V be a vector space over the field K. If B is a S-basis then B is a basis of V.*

*Proof*: Left for the reader as an exercise.

**DEFINITION [135]:** *Let V be a finite dimensional vector space over a field K. Let B = {$v_1$, ..., $v_n$} be a basis of V. If every proper subset of B generates a linear algebra over K then we call B a Smarandache strong basis (S-strong basis) for V.*

**DEFINITION [135]:** *Let V be any vector space over the field K. We say L is Smarandache finite dimensional vector space (S-finite dimensional vector space) of K if every S-basis has only finite number of elements in it. It is interesting to note that if L is a finite dimensional vector space then L is a Smarandache finite dimensional space (S-finite dimensional space) provided L has a finite S-basis.*

**THEOREM [135]:** *Let V be a vector space over the field K. If A = {$v_1$, ..., $v_n$} is a S-strong basis of V then A is a S-basis of V.*

*Proof*: Direct by the very definitions hence left for the reader to prove.

**DEFINITION 4.2.1**: *Let V be a vector space over a field F and let T be a linear operator from V to V. T is said to be a Smarandache linear operator (S-linear operator) on V if V has a S-basis which is mapped by T onto another basis V.*

**DEFINITION 4.2.2:** *Let T be a S-linear operator defined on the space V. A characteristic value C in F associated with T is said to be a Smarandache characteristic value (S-characteristic value) of T if the characteristic vector of T associated with C is in a linear algebra. So the eigen vector associated with the S-characteristic values will be called a Smarandache eigen vectors (S-eigen vectors) or Smarandache characteristic vector (S-characteristic vectors).*

As we do not have any book or paper on Smarandache vector spaces except a section of it covered by [135] we have just recalled what is very basic. All other results can be carried out as a matter of routine.

Now we proceed on to define Smarandache fuzzy vector spaces. From now onwards we take only Smarandache fields.

**DEFINITION [135]:** *A Smarandache field is defined to be a field (A, +, ×) such that a proper subset of A is a k-algebra (with respect to the same induced operations and external operations).*

But we in the paper [135] fields define differently.

**DEFINITION 4.2.3:** *A finite ring (i.e. a ring having finite number of elements) is said to be a Smarandache Galois field (S-Galois field) if S contains a proper subset A, A $\subset$ S such that A is a field under the operations of S.*

This definition is partly justified as all finite fields are Galois fields.

Now for us to develop the notion of Smarandache fuzzy vector spaces we make some more amendments in the definition of Smarandache vector space. We will call the definition of Smarandache vector space given here as the classical definition of Smarandache vector space. We define Smarandache vector space in a Smarandache style so that lucid expansion of the theory is possible. This Smarandache vector space we call it as Smarandache vector space of type II or non classical Smarandache vector space.

**DEFINITION 4.2.4:** *Let R be a S-ring. We call an additive abelian group V to be a Smarandache vector space of type II (S-vector space of type II) over R relative to F $\subset$ R (where F is a proper subset of R which is a field) if V is a vector space over F.*

*Thus at the first sight we may have V to be a S-vector space of type II relative to one subfield say F $\subset$ R; but V may fail to be a S-vector space of type II over some other subfield $F_1$ in R. Secondly when will it continue to be S-vector space of type II over all subfields in R.*

As the definition is very new we are constrained to give some examples of them.

***Example 4.2.1:*** Let $Z_{12}$ = {0,1, 2,..., 11} be a S-ring. Let V = P [x] where P is the field [0, 4, 8]. P [x] is the polynomial ring in the variable x over P = [0, 4, 8]. P[x] is a S-vector space of type II over the $Z_{12}$.

As P[x] is a vector space over P = [0,4, 8] $\subset Z_{12}$. We use this concept of S-vector space II to define S-fuzzy vector space.

**DEFINITION 4.2.5:** *Let R be a S-ring. V be a S-vector space of type II over R relative to P (P $\subset$ R). We call a fuzzy subset $\mu$ of V to be a Smarandache fuzzy vectorspace over the S-fuzzy ring (S-fuzzy vectorspace over the S-fuzzy ring) $\sigma$ of R (i.e. $\sigma: R \to [0, 1]$ is a fuzzy subset of R such that $\sigma: P \to [0, 1]$ is a fuzzy field where P $\subset$ R is a subfield of the S-ring R) if $\mu(0) > 0$ and for all x, y $\in$ V and for all c $\in$ P $\subset$ R, $\mu(x - y) \ge$ min {$\mu(x), \mu(y)$} and $\mu(cx)$ = min {$\sigma(c), \mu(x)$}.*

**DEFINITION 4.2.6:** *Let R be a S-ring having n-subfields in it say $P_1, ..., P_n$ (i.e.) each $P_i \subset$ R and $P_i$ is a subfield under the operations of R). Let V be a S-vector space over R. If V is a S-vector space over R relative to every subfield $P_i$ in R then we call V the Smarandache strong vector space (S-strong vector space) over R.*

Thus we have the following straightforward result.

**THEOREM 4.2.1:** *Let R be a S-ring. If V is a S-strong vector space over R then V is a S-vector space over R.*

**THEOREM 4.2.2:** *A S-vector space V over a S-ring R in general is not a S-strong vector space over R.*

*Proof*: By an example. Consider $Z_6 = \{0, 1, 2, 3, 4, 5\}$ be the S-ring under addition and multiplication modulo 6. Take $P_1 = \{0,3\}$ and $P_2 = \{0, 2, 4\}$; these are the only subfields of $Z_6$. Let $V = P_1[x]$; $P_1[x]$ is a S-vector space over $P_1$ but V is not a S-vector space over $P_2$. Thus V is not a S-strong vector space. Hence the claim.

On similar lines we can define S-strong fuzzy space.

**DEFINITION 4.2.7:** *Let R be a S-ring. V be a S-strong vector space of type II over R. We call a fuzzy subset μ of V to be a Smarandache strong fuzzy vector space (S-strong fuzzy vector space) over the S-fuzzy ring $\sigma_i$ of R; $\sigma_i: P_i \subset R \to [0, 1]$ where $P_i \subset R$ are subfields for i = 1, 2, ..., n if μ(0) > 0 and for all x, y ∈ V and for all $c \in P_i \subset R$ (i = 1, 2,..., n), $\mu(x - y) \geq \min \{ \mu(x), \mu(y)\}$ and $\mu(cx) = \min \{\sigma_i(c), \mu(x)\}$, i = 1, 2,..., n.*

As in case of S-strong vector spaces we have the following theorem.

**THEOREM 4.2.3:** *Let V be a S-strong fuzzy vector space over the S-ring R. Then V is a S-fuzzy vector space over R..*

*Proof*: Direct by the very definitions.

The reader is expected to construct examples and counter examples to study them.

If A and B are fuzzy subset of V then $A \subseteq B$ means $A(x) \leq B(x)$ for all $x \in V$. For $0 \leq t \leq 1$, let $A_t = \{x \in V \mid A(x) \geq t\}$, where V is a S-fuzzy vector space over the S-ring R. Let $\wp$ denote the set of all S-fuzzy subfield of a S-ring R. Let $\mathcal{A}_P$ denote the set of all S-fuzzy subspaces of V over $P \subset R$.

**DEFINITION 4.2.8:** *Let A, $A_1$,..., $A_n$ be fuzzy subsets of a S-vector space V and let K be a fuzzy subset of the S-ring R. Define Smarandache fuzzy subset $A_1 + ... + A_n$ of V by the following, for all $x \in V$. $(A_1 + ... + A_n)(x) = \sup \{\min \{A_1(x_1),..., A_n(x_n)\} / x = x_1 + ... + x_n, x_i \in V\}$. Define the fuzzy subset K o A of V by for all $x \in V$, $(K o A)(x) = \sup\{\min \{K(c), A(y)\} / c \in P \subset R, y \in V, x = cy / P \text{ is a subfield in R relative to which V is defined}\}$.*

**DEFINITION 4.2.9:** *Let $\{A_i / i \in I\}$ be a non empty collection of fuzzy subsets of V, V a S-fuzzy subspace of V. Then the fuzzy subset $\bigcap_{i \in I} A_i$ of V is defined by the following for all $x \in V$, $\left(\bigcap_{i \in I} A_i\right)(x) = \inf \{A_i(x) / i \in I\}$.*

Let $A \in \mathcal{A}_k$, K a fuzzy subset of the S-ring R relative to a subfield $P \subset R$. X be a fuzzy subset of V such that $X \subset A$. Let $\langle X \rangle$ denote the intersection of all fuzzy subspaces of

*V (over K) that contain X and are contained in A then ⟨X⟩ is called the Smarandache fuzzy subspace (S-fuzzy subspace) of A fuzzily spanned or generated by X.*

*We give some more notions and concepts about S-fuzzy vector spaces. Let ξ denote a set of fuzzy singletons of V such that $x_\lambda$, $x_k \in \xi$ then $\lambda = k > 0$. Define the fuzzy subset of X (ξ) of V by the following for all $x \in V$. $X(\xi)(x) = \lambda$ if $x_\lambda \in \xi$ and $X(\xi)(x) = 0$ otherwise. Define ⟨ξ⟩ = ⟨X (ξ)⟩. Let X be a fuzzy subset of V, define $\xi(X) = \{x_\lambda \mid x \in V, \lambda = X(x) > 0\}$.*

*Then $X(\xi(X)) = X$ and $\xi(X(\xi)) = \xi$. If there are only a finite number of $x_\lambda \in \xi$ with $\lambda > 0$ we call ξ finite (or Smarandache finite). If $X(x) > 0$ for only a finite number of $x \in X$, we call X finite. Clearly ξ is finite if and only if X (ξ) is S-finite and X is finite if and only if ξ (X) is S-finite. For $x \in V$ let $X \setminus \{x\}$ denote the fuzzy subset of V defined by the following; for all $y \in V$. $(X \setminus \{x\})(y) = X(y)$ if $y \neq x$ and $(X \setminus \{x\})(y) = 0$ if $y = x$. Let $A \in A_K$ (K a fuzzy subset defined relative to a subfield P in a S-ring R) and let X be a fuzzy subset of V such that $X \subseteq A$. Then X is called the Smarandache fuzzy system of generators (S-fuzzy system of generators) of A over K if ⟨X⟩ = A. X is said to be Smarandache fuzzy free (S-fuzzy free) over K if for all $x_\lambda \subset X$ where $\lambda = X(x)$, $x_\lambda \not\subset \langle X \setminus x \rangle$. X is said to be a Smarandache fuzzy basis {S-fuzzy basis} for A relative to a subfield $P \subset R$, R a S-ring if X is a fuzzy system of generators of A and X is S-fuzzy free. Let ξ denote a set of fuzzy singletons of V such that if $x_\lambda$, $x_k \in \xi$ then $\lambda = k$ and $x_\lambda \subseteq A$. Then ξ is called a Smarandache fuzzy singleton system of generators (S-fuzzy singleton system of generators) of A over K if ⟨ξ⟩ = A. ξ is said to be S-fuzzy free over K if for all $x_\lambda \in \xi$, $x_\lambda \not\subset \langle \xi \setminus \{x_\lambda\} \rangle$. ξ is said to be a S-fuzzy basis of singletons for A if ξ is a fuzzy singleton system of generators of A and ξ is S-fuzzy free.*

Several interesting results in the direction of Smarandache fuzzy free can be obtained analogous to results on fuzzy freeness.

Now we proceed on to define Smarandache fuzzy linearly independent over $P \subset R$.

**DEFINITION 4.2.10**: *Let R be a S-ring, V be a S-vector space over $P \subset R$ (P a field in R). Let $A \in \mathcal{A}_k$, K a fuzzy field of P or K is a S-fuzzy ring R relative to P and let $\xi \subseteq \{x_\lambda \mid x \in A^*\, \lambda \leq A(x)\}$ be such that if $x_\lambda$, $x_k \in \xi$ then $\lambda = k$. Then ξ is said to be Smarandache fuzzy linearly independent (S-fuzzy linearly independent) over K of P if and only if for every finite subset $\left(x_{1_{\lambda_1}}, \cdots, x_{n_{\lambda_n}}\right)$ of ξ whenever $\left(\sum_{i=1}^{n} c_{i\mu} o\, x_{i_{\lambda_j}}\right)(x) = 0$ for all $x \in V \setminus \{0\}$ where $c_i \in P \subset R$, $0 < \mu_i \subseteq K(c_i)$ for $i = 1, 2, \ldots, n$ then $c_1 = c_2 = \ldots = c_n = 0$.*

Now we proceed on to define Smarandache fuzzy modules.

**DEFINITION 4.2.11:** *The Smarandache R-module (S-R-module) is defined to be an R-module $(A, +, \times)$ such that a proper subset of A is a S-algebra (with respect to the same induced operations and another '×' operations internal on A) where R is a commutative unitary S-ring and S its proper subset which is field.*

*Example 4.2.2:* Let R[x] be the polynomical ring in the variable x with coefficients from the real field R. Q[x] is a S-R-module for it is a S-algebra.

We just recall the definition of Smarandache right (left) module.

**DEFINITION 4.2.12:** *Let R be a S-ring I. A non-empty set B which is an additive abelian group is said to be a Smarandache right (left) module I (S-right(left) module I) relative to a S-subring I, A if $D \subset A$ where D is a field then $DB \subset B$ and $BD \subset B$ i.e. bd (and db) are in B with $b(d + c) = bd + dc$ for all $d, c \in D$ and $b \in B$ $((d + c)b = db + cb)$. If B is simultaneously a S-right module I and S-left module I over the same relative S-subring I then we say B is a Smarandache module I (S-module I).*

**DEFINITION 4.2.13:** *Let R be a S-ring II. We say a non-empty set B which is an additive abelian group is said to be a Smarandache right (left) module II (S-right (left) module II) relative to a S-subring II, A if $D \subset A$ where D is a division ring or an integral domain, then $DB \subset B$ and $BD \subset B$; i.e., bd(and db) are in B. with $b(d + c) = bd + bc$ $\forall d, c \in D$ and $b \in B$ $((d + c) b = db + cb)$. If B is simultaneously a S-right module II and S-left module II over the same relative S-subring II then we say B is a Smarandache module II (S-module II).*

**DEFINITION 4.2.14:** *Let $(A, +, \bullet)$ be a S-ring. B be a proper subset of A ($B \subset A$) which is a field. A set M is said to be a Smarandache pseudo right (left) module (S-pseudo right(left) module) of A related to B if*

  i.  *$(M, +)$ is an additive abelian group.*
  ii. *For $b \in B$ and $m \in M$ $m.b \in M$ $(b.m \in M)$.*
  iii. *$(m_1 + m_2)b = m_1b + m_2b$, $(b.(m_1+m_2)=bm_1+bm_2)$ for $m_1, m_2 \in M$ and $b \in B$. If M is simultaneously a S-pseudo right module and S-pseudo left module, we say M is a Smarandache pseudo module (S-pseudo module) related to B.*

Here also we wish to state if $M_1$ is a S-pseudo module related to B, $M_1$ need not be S-pseudo module related to some other subfield $B_1$ of A. Thus we see we can have different S-pseudo modules associated with different subfields in a ring.

*Example 4.2.3:* Let $Z_2S_4$ be the group ring of the symmetric group of degree 4 over the field $Z_2$. $M = \{0, \Sigma g, g \in S_4\}$ ($\Sigma g$ denotes the sum of all elements from $S_4$). M is a S-module II over $Z_2$. M is a S-M-module II over $Z_2A_4$. Clearly M is also a S-ideal II and S-pseudo ideal of $Z_2S_4$.

But to define Smarandache fuzzy module and Smarandache fuzzy algebra we have to strike a mediocre between the S-modules and fuzzy modules and S-algebras and fuzzy algebras. As the concept of fuzziness involves lots of fuzziness and Smarandache property involves adaptability of both richer structure as a substructure or weaker structure as a substructure of a richer structure we define Smarandache fuzzy modules and Smarandache fuzzy algebras in a very different way.

**DEFINITION 4.2.15:** *Let M be a S-left (or S-right) R-module over a S-ring R. ($\eta$, M) is called a Smarandache fuzzy left (Smarandache fuzzy right) R-module if there is a map $\eta : M \to [0, 1]$ satisfying the following conditions:*

  i.   $\eta (a + b) \geq \min \{\eta(a), \eta(b)\}$ $a, b \in M$.
  ii.  $\eta (-a) = \eta(a)$ for all $a \in M$.
  iii. $\eta (0) = 1$ (0 is the zero element of M),
  iv.  $\eta (ra) \geq \eta (a)$ ($a \in M$ and $r \in P \subset R$, P a field) It is denoted by $\eta_M^P$ where the S-fuzzy module is defined relative to $P \subset R$; P a subfield of R.

*$\eta(M)$ expresses the cardinal number of all fuzzy values $\eta(a)$ ($a \in M$). We are interested in the situation $\eta(V) < + \infty$. It is well known fact that Noetherian modules are a large class of modules which play an important role in modular theory.*

Several important and interesting results can be obtained in this direction.

### 4.3. Smarandache fuzzy non-associative rings

In this section we introduce the notion of Smarandache fuzzy non associative rings and recall the basic notion about Smarandache non associative ring which for short we denote as SNA-rings. Lot of information about SNA-rings can be had from [131].

Throughout this book by a non associative ring $(R, +, \bullet)$ we mean a non empty set R endowed with two binary operations '+' and '•' such that $(R, +)$ is an additive abelian group with '0' acting as the additive identity and $(R, \bullet)$ is a non associative semigroup or a groupoid such that the distributive laws $a \bullet (b + c) = a \bullet b + a \bullet c$ and $(a + b) \bullet c = a \bullet c + b \bullet c$ for all $a, b, c \in R$ holds true.

Now we recall the definition of Smarandache non-associative ring (SNA-ring).

**DEFINITION 4.3.1:** *Let $(R, +, \bullet)$ be a non-associative ring R is said to be a Smarandache non-associative ring (SNA-ring) if R contains a proper subset P such that P is an associative ring under the operations of R.*

*Example 4.3.1:* Let R be a field and L be a loop given by the following table. RL be the loop ring of the loop L over R. RL is a SNA-ring. The loop L is given by the following table:

| *     | e     | $a_1$ | $a_2$ | $a_3$ | $a_4$ | $a_5$ | $a_6$ | $a_7$ |
|-------|-------|-------|-------|-------|-------|-------|-------|-------|
| e     | e     | $a_1$ | $a_2$ | $a_3$ | $a_4$ | $a_5$ | $a_6$ | $a_7$ |
| $a_1$ | $a_1$ | e     | $a_5$ | $a_2$ | $a_6$ | $a_3$ | $a_7$ | $a_4$ |
| $a_2$ | $a_2$ | $a_5$ | e     | $a_6$ | $a_3$ | $a_7$ | $a_4$ | $a_1$ |
| $a_3$ | $a_3$ | $a_2$ | $a_6$ | e     | $a_7$ | $a_4$ | $a_1$ | $a_5$ |
| $a_4$ | $a_4$ | $a_6$ | $a_3$ | $a_7$ | e     | $a_1$ | $a_5$ | $a_2$ |
| $a_5$ | $a_5$ | $a_3$ | $a_7$ | $a_4$ | $a_1$ | e     | $a_2$ | $a_6$ |
| $a_6$ | $a_6$ | $a_7$ | $a_4$ | $a_1$ | $a_5$ | $a_2$ | e     | $a_3$ |
| $a_7$ | $a_7$ | $a_4$ | $a_1$ | $a_5$ | $a_2$ | $a_6$ | $a_3$ | e     |

Clearly RL is a SNA-ring. We can get class of SNA-rings using loops over rings and groupoid over rings which we choose to call as loop rings and groupoid rings.

**DEFINITION 4.3.2:** *Let R be a non-associative ring. A non-empty subset S of R is said to be a SNA-subring of R if S contains a proper subset P such that P is an associative ring under the operations of R.*

The following theorem is straightforward.

**THEOREM 4.3.1:** *Let R be a non-associative ring. If R has a proper subset S such that S is a SNA subring then R is a SNA-ring.*

Further it is left for the reader to verify that every subring of a SNA-ring R in general need not be a S-subring of R.

**DEFINITION 4.3.3:** *Let R be any non associative ring. A proper subset I of R is said to be a SNA right / left ideal of R if*

  i. *I is a SNA-subring of R say J ⊂ I, J is a proper subset of I which is an associative subring under the operations of R.*
  ii. *If I is simultaneous by both a SNA right ideal and SNA left ideal then we say I is a SNA-ideal of R.*

**THEOREM 4.3.2:** *Let R be any non associative ring. If R has a SNA-ideal then R is a SNA-ring.*

*Proof*: Obvious from the fact that if R has a SNA-ideal say I then we see R has a subset which is an associative ring. Hence the claim.

We have SNA-rings which satisfies special identities like Bruck, Bol, Moufang, right/ left alternative.

**DEFINITION 4.3.4:** *Let R = (R , +, ∗) be a SNA-ring we say R is a SNA-Moufang ring if R contains a subring S where S is a SNA subring and for all x, y z in S we have*

$$(x * y) * (z * x) = (x * (y * z)) * x$$

*that is the Moufang identity is true.*

*(R, +, ∗) is said to be SNA Bol ring if R satisfies the Bol identity*

$$((x * y) * z) * y = x * ((y * z) * y)$$

*for all x, y, z ∈ S ⊂ R; where S is a proper SNA-subring. We say the SNA-ring R is said to be SNA-right (left) alternative ring if (xy)y = x(yy) [(xx)y = x(xy)]. If R is a simultaneously both a right and left alternative NA ring then we call R a Smarandache alternative (S-alternative ring) ring.*

**DEFINITION 4.3.5:** *Let R be a non-associative ring. R is said to be a SNA-commutative ring if R has a subring S such that a proper subset P of S is a commutative associative ring with respect to the operation of R. R is said to be a SNA strong commutative ring if every subset P of R which is an associative subring is commutative.*

Now we proceed on to define Smarandache fuzzy non-associative rings.

**DEFINITION 4.3.6:** *Let R be a SNA-ring. A fuzzy subset $\mu : R \rightarrow [0, 1]$ is said to be a Smarandache non-associative fuzzy ring (SNA-fuzzy ring) if $\mu : P \rightarrow [0, 1]$ is a fuzzy ring where P is a proper subset of R which is an associative ring under the operations of R and denote it by $\mu_P$. We call the fuzzy subset $\mu : R \rightarrow [0, 1]$ to be a Smarandache non-associative strong fuzzy ring (SNA-strong fuzzy ring) if $\mu : P_i \rightarrow [0, 1]$ is a fuzzy subring of R for every proper subset $P_i \subset R$ which is a subring of R. Thus we denote the SNA strong fuzzy ring by just $\mu$ as $\mu$ is a SNA-fuzzy ring for every subring $P_i \subset R$ where $P_i$ is a subring of R.*

We have the following interesting theorem.

**THEOREM 4.3.3:** *Let R be a SNA-ring. Every SNA-strong fuzzy ring is a SNA-fuzzy ring.*

*Proof*: Straightforward by the very definitions.

**DEFINITION 4.3.7:** *Let R be a SNA-ring. A fuzzy subset $\mu$ of the SNA-ring R is said to be Smarandache non associative fuzzy subring (S-non associative fuzzy subring) of $\mu$, if $\mu : T \rightarrow [0, 1]$ is a SNA-fuzzy ring where T is a SNA-subring of R. i.e. $\mu_L : L \rightarrow [0, 1]$ where L is a proper subset of T and L is an associative subring of $T \subset R$.*

**DEFINITION 4.3.8:** *Let R be a SNA-ring. A fuzzy subset $\mu$ of R is said to be Smarandache non associative fuzzy ideal (S-non associative fuzzy ideal) of R if for all $x, y \in P \subset R$ (where P is a subset of R which is an associative subring under the operations of R). $\mu(x - y) \geq \min(\mu(x), \mu(y))$ and $\mu(xy) \geq \max(\mu(x), \mu(y))$. Recall if $\mu$ and $\theta$ be any two SNA-fuzzy ideals of a SNA-ring R. The product $\mu \circ \theta$ of $\mu$ and $\theta$ is defined by*

$$(\mu \circ \theta)(x) = \sup_{\substack{x = \sum y_i z_i \\ i < \infty}} \{\min(\min(\mu(y_i), \theta(z_i))\}$$

*where $x, y_i, z_i \in P \subset R$ }. It is easily verified that $\mu \circ \theta$ is a SNA-fuzzy ideal of R and $\mu \circ \theta = \langle \mu\theta \rangle$.*

Further $\mu\theta$ in general is not a SNA-fuzzy ideal. It is left as an exercise for the reader to prove the above statement.

**DEFINITION 4.3.9:** *An SNA-ideal P of a SNA-ring R, $P \neq R$ is called SNA-prime if $a, b \in P$, implies either $a \in P$ or $b \in P$.*

**DEFINITION 4.3.10:** *A SNA fuzzy ideal µ of a SNA-ring R is called SNA-fuzzy prime ideal if the ideal $\mu_t$ where $t = \mu(0)$ is a SNA-prime ideal of R.*

*A non constant SNA-fuzzy ideal µ of a SNA-ring R is called SNA-fuzzy prime if for any two SNA-fuzzy ideals $\sigma$ and $\theta$ of R (defined relative to the same $P \subset R$) the condition $\sigma\theta \subset \mu$ implies either $\sigma \subset \mu$ or $\theta \subseteq \mu$.*

*Recall a SNA-ring R is SNA-regular if for each element x in $P \subset R$ (P an associative subring of R) there exists y in P such that xyx = x. We say R is SNA-strongly regular if for each subset $P \subset R$, P an associative ring; P is regular. It can be easily verified.*

**THEOREM 4.3.4:** *If R is a SNA strongly regular ring then R is a SNA-regular ring.*

We have the notion of fuzzy primary ideal in case of SNA-rings.

**DEFINITION 4.3.11:** *Let R be a SNA-ring. A SNA-fuzzy ideal µ of R relative to $P \subset R$ (P an associative subring of R) is SNA-fuzzy primary if for any two SNA fuzzy ideals $\sigma$ and $\theta$ of R relative to the same associative subring P the conditions $\sigma\theta \subseteq \sqrt{\mu}$ and $\sigma \not\subset \mu$ together imply $\theta \subseteq \sqrt{\mu}$.*

Several other properties and results for SNA-rings can be developed analogous to associative rings.

Now we just define SNA-fuzzy semiprimary ideal of a SNA-ring R.

**DEFINITION 4.3.12:** *Let R be a SNA-ring. A SNA fuzzy ideal µ of the ring R relative to $P \subset R$ is called SNA-fuzzy semiprimary if $\sqrt{\mu}$ is a SNA-fuzzy prime ideal of R relative to P.*

Almost all results in case of fuzzy semiprimary ideals of an associative ring can be deduced/ derived in case of SNA-rings.

Now we proceed on to define the concept of fuzzy irreducible ideals in case of SNA-rings.

**DEFINITION 4.3.13:** *Let R be a SNA-ring. A fuzzy ideal µ of the ring R relative to the associative subring $P \subset R$ is called SNA-fuzzy irreducible if it is not an intersection of any two SNA-fuzzy ideals of R (relative to the same P) properly containing µ; otherwise µ is termed SNA-fuzzy reducible.*

**THEOREM 4.3.5:** *Let R be a SNA-ring. If µ is a SNA-fuzzy prime ideal of R then µ is SNA-fuzzy irreducible.*

*Proof*: Using the fact µ is a SNA-fuzzy prime ideal of R relative to the associative subring P in R we prove µ is SNA-fuzzy irreducible relative to the same subring P.

Several results true in case of S-ring can be also derived in case of SNA-rings. This task is left as an exercise for the reader. This study of SNA-rings is more complicated and interesting as μ a fuzzy structure can be defined relative to every associative subring of a SNA-ring. Each subring may possess various properties accordingly the μ will also enjoy distinct properties. So further research in this direction will yield many interesting results.

## 4.4. Smarandache fuzzy birings and its properties

In this section we just introduce the notions of Smarandache fuzzy birings and its properties. We deal with both birings which are associative as well as non associative. We define some concepts so that the reader can develop further properties that are already studied in case of associative and non-associative rings. The concept of birings and fuzzy birings are introduced in chapter 1.

**DEFINITION 4.4.1:** *A Smarandache biring (S-biring) $(R, +, \bullet)$ is a non empty set with two binary operations '+' and '$\bullet$' such that $R = R_1 \cup R_2$ where $R_1$ and $R_2$ are proper subsets of R and*

    i. *$(R_1, +, \bullet)$ is a S-ring.*
    ii. *$(R_2, +, \bullet)$ is a S-ring.*

*If only one of $R_1$ or R is a S-ring then we call $(R, +, \bullet)$ a Smarandache weak biring (S- weak biring).*

We define Smarandache conventional biring as follows:

**DEFINITION 4.4.2:** *Let $(R, +, \bullet)$ be a biring. R is said to be a Smarandache conventional biring (S-conventional biring) if and only if R has a proper subset that is a bifield.*

**DEFINITION 4.4.3:** *Let $(R, +, \bullet)$ be a biring. A proper subset P of R is said to be a Smarandache sub-biring (S-sub-biring) if $(P, +, \bullet)$ is itself a S-biring. Similarly we say for a biring $(R, +, \bullet)$ a proper subset P of R is said to be as Smarandache conventional sub-biring (S-conventional sub-biring) if $(P, +, \bullet)$ is itself a S-conventional biring under the operations of R.*

The following theorem is straightforward.

**THEOREM 4.4.1:** *If a biring R has a S-sub-biring then R is itself a S-biring. If a biring R has S-conventional sub-biring then R itself is a S-conventional biring.*

**DEFINITION 4.4.4:** *Let R be a biring. We say R is a Smarandache commutative biring (S-commutative biring) if every S-sub-biring is commutative. If atleast one of the S-sub-biring of R is commutative then we say R is a Smarandache weakly commutative biring (S-weakly commutative biring).*

The following theorem is direct hence left for the reader to prove.

**THEOREM 4.4.2:** *Let (R, +, •) be a biring. If R is a S-commutative biring then R is a S-weakly commutative biring.*

Now we define type II Smarandache biring.

**DEFINITION 4.4.5:** *Let (R, +, •) be a non empty set. We say (R, +, •) is a Smarandache biring II (S-biring II) if R contains a proper subset P such that P is a bidivision ring. If R has no S-zero divisors then we call R a Smarandache division biring (S-division biring) if R is non commutative and R is commutative then we call R a Smarandache integral domain (S-integral domain).*

We recall the definition of Smarandache bi-ideal and Smarandache pseudo bi-ideal of S-birings.

**DEFINITION 4.4.6:** *Let (R, +, •) be a biring. The Smarandache bi-ideal (S-bi-ideal) P of R is defined as an ideal of the biring R where P is a S-sub-biring.*

**DEFINITION 4.4.7:** *Let (R, +, •) be a S-conventional biring. $B \subset R$ is a bifield of R . A non-empty subset P of R is said to be Smarandache pseudo right bi-ideal (S-pseudo right bi-ideal) of R related to P if*

  i. *(P,+) an abelian bigroup.*
  ii. *For $b \in B$ and $p \in P$ we have $p \bullet b \in P$.*

*On similar lines we define Smarandache pseudo left bi-ideal (S-pseudo left bi-ideal); if P is both a S-pseudo right bi-ideal and S-pseudo left bi-ideal.*

*We call a S-biring (R, +, •) to be a Smarandache simple biring (S-simple biring) if R has no S- bi-ideals. If R has no S-pseudo bi-ideals then we call R Smarandache pseudo simple biring (S-pseudo simple biring).*

**DEFINITION 4.4.8:** *Let (R, +, •) be a biring. A non-empty subset M is said to be a Smarandache bimodule (S-bimodule) if M is an S-abelian semigroup with 0 under the operation + i.e. $M = M_1 \cup M_2$ and Let $R = R_1 \cup R_2$ be a biring, a non-empty set M is said to be a R-bimodule (or a bimodule over the biring R) if M is an abelian bigroup, under addition '+' say with $M = M_1 \cup M_2$ such that for every $r_1 \in R_1$ and $m_1 \in M_1$ there exists an element $r_1 m_1$ in $M_1$ and for every $r_2 \in R_2$ and $m_2 \in M_2$ there exists an element $r_2 m_2$ in $M_2$ subject to*

  i. $r_1 (a_1 + b_1) = r_1 a_1 + r_1 b_1$ ; $r_1 \in R_1$ and $a_1, b_1 \in M_1$.
  ii. $r_2 (a_2 + b_2) = r_2 a_2 + r_2 b_2$ ; $r_2 \in R_2$ and $a_2, b_2 \in M_2$.
  iii. $r_1 (s_1 a_1) = (r_1 s_1) a_1$ for $r_1, s_1 \in R_1$ and $a_1 \in M_1$.
  iv. $r_2 (s_2 a_2) = (r_2 s_2) a_2$ for $r_2, s_2 \in R_2$ and $a_2 \in M_2$.
  v. $(r_1 + s_1) a_1 = r_1 a_1 + s_1 a_1$, $r_1, s_1 \in R_1$ and $a_1 \in M_1$.
  vi. $(r_2 + s_2) a_2 = r_2 a_2 + s_2 a_2$ where $r_2, s_2 \in R_2$ and $a_2 \in M_2$.

*In short if M is a bigroup ($M = M_1 \cup M_2$) under '+' then $M_1$ is a $R_1$-module and $M_2$ is a $R_2$-module then we say M is a R-bimodule over the biring $R = R_1 \cup R_2$.*

*If the biring R has unit 1 and if $1.m_1 = m_1$ and $1.m_2 = m_2$ for every $m_1 \in M_1$ and $m_2 \in M_2$ then we call M a unitial R-bimodule. Thus the concept of R-bimodule forces both the structures to be bistructures i.e. we demand the group should be a bigroup and also the ring must be a biring then only we can speak of a R-bimodule.*

*Let M be a R-bimodule an additive sub-bigroup, A of M i.e. $A = A_1 \cup A_2$ where $A_1$ is a subgroup of $M_1$ and $A_2$ is a subgroup of $A_2$ is called the sub-bimodule of the bimodule M if when ever $r_1 \in R_1$ and $r_2 \in R_2$ and $a_1 \in A_1$ and $a_2 \in A_2$ we have $r_1 a_1 \in A_1$ and $r_2 a_2 \in A_2$.*

*A bimodule M is cyclic if there is an element $m_1 \in M_1$ and $m_2 \in M_2$ such that for every $m \in M_1$ is of the form $m = r_1 m_1$ where $r_1 \in R_1$ and for every $m' \in M_2$ is of the form $m' = r_2 m_2$ where $r_2 \in R_2$. Thus cyclic bimodules is nothing but bicyclic groups, that is the bigroup $G = G_1 \cup G_2$ is a bicyclic group if both $G_1$ and $G_2$ are cyclic groups.*

**DEFINITION 4.4.9**: *Let (R, +, o) be a biring. We say R is a Smarandache semiprime biring (S-semiprime biring) if R contains no non zero S-bi-ideal with square zero.*

Now we proceed on to define S-na biring.

**DEFINITION 4.4.10:** *Let (R, +, •) with $R = R_1 \cup R_2$ be a na-biring, R is said to be a Smarandache na-biring (S-na-biring) if R has a proper subset P such that $P = P_1 \cup P_2$ and P is an associative biring under the operations of R. Let (R, +, •) be a na-biring.*

*A proper subset $P \subset R$ is said to be a Smarandache na sub-biring (S-na-sub-biring) if P is itself a na-sub-biring and P is a S-na-biring.*

**THEOREM 4.4.3:** *If (R, +, •) has a na-sub-biring which is a S-na sub-biring then we say R is a S-na-biring.*

*Proof*: Follows directly from the definitions, hence left for the reader to prove.

**DEFINITION 4.4.11:** *Let (R, +, •) be a na-biring. We say R is a Smarandache commutative na-biring (S-commutative na-biring) if every S- na sub-biring of R is commutative. If at least one S-na sub-biring is commutative then we call R a Smarandache weakly commutative na-biring (S-weakly commutative na-biring).*

**THEOREM 4.4.4:** *Let $R = R_1 \cup R_2$ be a na-biring. If R is a strongly subcommutative biring then R is a S-commutative biring.*

*Proof:* Straightforward, hence left for the reader to prove.

**DEFINITION 4.4.12:** *Let (R, +, •) be a na-biring, we call R a Smarandache Moufang biring (S-Moufang biring) if R has a proper subset P; such that P is a S-sub-na-biring of P and every triple of R satisfies the Moufang identity*

$$(xy)(zx) = (x(yz))x$$

*for all x, y, z ∈ P.*

*Similarly $R = R_1 \cup R_2$ is a Smarandache Bruck biring (S-Bruck biring) if R has proper subset P such that*

    i.   *P is a S-na sub-biring of R .*
    ii.  *(x(yx)) z = x (y(xz)) and $(xy)^{-1} = x^{-1}y^{-1}$*

*is true for all x, y, z ∈ P. We call (R, +, •) a Smarandache Bol biring (S-Bol biring) if ((xy z) y = x((yx) y) for all x, y, z ∈ P, P ⊂ R is a S-na sub-biring of R.*

On similar lines we define Smarandache left (right) alternative birings and Smarandache WIP-birings.

**DEFINITION 4.4.13:** *Let $R = (R_1 \cup R_2, +, •)$ be a na-biring. If every S-na sub-biring P of R satisfies Moufang, Bol, Bruck, alternative (right / left) or WIP then we call R a Smarandache strong Moufang, Bol, Bruck, alternative (right/ left) or WIP biring respectively.*

Several interesting results can be obtained in this direction.

Now we proceed on to define the concepts of Smarandache fuzzy birings.

**DEFINITION 4.4.14:** *Let (R, +, •) be a S-biring where $R = R_1 \cup R_2$. A fuzzy subset $\mu : R \to [0, 1]$ is said to be a Smarandache fuzzy biring (S-fuzzy biring) if the following conditions are true.*

$\mu = \mu_1 \cup \mu_2 : R_1 \cup R_2 \to [0, 1]$ is such that $\mu : R_1 \to [0, 1]$ is a S-fuzzy ring of $R_1$ relative to $P_1 \subset R_1$ and $\mu_2 : R_2 \to [0, 1]$ is a S-fuzzy ring of $R_2$ relative to $P_2 \subset R_2$. (Here $P_1$ and $P_2$ are proper subsets of $R_1$ and $R_2$ such that $P_1$ and $P_2$ are subfields of $R_1$ and R respectively) i.e. $\mu_i : P_i \to [0, 1]$ is a fuzzy subfield; i = 1, 2. If $\mu : R \to [0, 1]$ the fuzzy subset of R is such that $\mu$ is a S-fuzzy biring relative to every subfield in R then we call $\mu$ a Smarandache strong fuzzy biring (S-strong fuzzy biring).

The following theorem is direct hence left for the reader as an exercise.

**THEOREM 4.4.5:** *Every S-strong fuzzy biring R is a S-fuzzy biring.*

Now we proceed on to define Smarandache fuzzy sub-biring and Smarandache fuzzy bi-ideal of a S-biring R.

**DEFINITION 4.4.15:** *Let (R, +, •) be a S-biring. A fuzzy subset $\mu : R \to [0, 1]$ is said to be a Smarandache fuzzy ideal (S-fuzzy ideal) of the biring R if the following conditions are satisfied.*

*$\mu = \mu_1 \cup \mu_2 : R_1 \cup R_2 \to [0, 1]$ is such that $\mu_1 : R_1 \to [0, 1]$ is a S-fuzzy ideal of $R_1$ and $\mu_2 : R_2 \to [0, 1]$ is a S-fuzzy ideal of $R_2$. We call the fuzzy subset $\mu =$* $\mu_1 \cup \mu_2 : R \to [0, 1]$ to be a Smarandache fuzzy weak bi-ideal (S-fuzzy weak bi-ideal) of $R = R_1 \cup R_2$ if atleast one $\mu_1$ or $\mu_2$ is a S-fuzzy ideal and the other is just a fuzzy ideal.

**THEOREM 4.4.6:** *Let $(R = R_1 \cup R_2, +, \bullet)$ be a S-biring. $\mu : R \to [0, 1]$ be a S-fuzzy bi-ideal of R then $\mu$ is a S-fuzzy weak bi-ideal of R.*

*Proof*: Obvious by the very definitions.

Now just we define the notion of Smarandache fuzzy sub-biring.

**DEFINITION 4.4.16:** *Let $(R = R_1 \cup R_2, +, \bullet)$ be a S-biring. The fuzzy subset $\mu = \mu_1 \cup \mu_2 : R \to [0, 1]$ is a said to be a Smarandache fuzzy sub-biring (S-fuzzy sub-biring) of R if $\mu = \mu_1 \cup \mu_2 : R \to [0, 1]$ is such that $\mu_1 : R_1 \to [0, 1]$ a S-fuzzy subring of $R_1$ and $\mu_2 : R_2 \to [0, 1]$ is a S-fuzzy subring of $R_2$.*

*Now if $\mu = \mu_1 \cup \mu_2 : R_1 \cup R_2 \to [0, 1]$ is said to be a S-fuzzy weak sub-biring if only one of $\mu_1$ of $\mu_2$ is a S-fuzzy subring and the other is just fuzzy subring.*

The following theorem is left as an exercise for the reader to prove.

**THEOREM 4.4.7:** *Let $(R = R_1 \cup R_2, +, \bullet)$ be a S-biring. The fuzzy subset $\mu = \mu_1 \cup \mu_2 : R_1 \cup R_2 \to [0, 1]$ is a S-fuzzy sub-biring of R then $\mu$ is a S-fuzzy weak sub-biring of R.*

**DEFINITION 4.4.17:** *Let $(R = R_1 \cup R_2, +, \bullet)$ be a S-biring. A S-fuzzy bi-ideal $\mu$ of the S-biring R is said to be Smarandache fuzzy prime (S-fuzzy prime) if the ideal $\mu_t = (\mu_1)_t \cup (\mu_2)_t$ where $t = \mu(0) = \mu_1(0) \cup \mu_2(0)$ is a S-prime bi-ideal of R.*

Following the definition of [100] we have this analogous definition for S-fuzzy prime bi-ideals.

**DEFINITION 4.4.18:** *Let $(R = R_1 \cup R_2, +, \bullet)$ be a S-biring. A non constant S-fuzzy bi-ideal $\mu$ of the biring R is called Smarandache fuzzy prime bi-ideal (S-fuzzy prime bi-ideal) of R if for any two S-fuzzy bi-ideals $\sigma$ and $\theta$ of R the condition $\sigma\theta \subseteq \mu$ implies either $\sigma \subset \mu$ or $\theta \subseteq \mu$.*

*A S-fuzzy bi-ideal $\mu$ of a biring R not necessarily non constant is called Smarandache fuzzy prime if for any two S-fuzzy bi-ideals $\sigma$ and $\theta$ of R the condition $\sigma\theta \subset \mu$ implies either $\sigma \subset \mu$ or $\theta \subset \mu$.*

**Notation**: Let $\mu$ and $\theta$ be any two S-fuzzy bi-ideals of a S-biring R. $\mu, \theta : R_1 \cup R_2 \to [0, 1]$ the sum $\mu + \theta = (\mu_1 \cup \theta_1) \cup (\mu_2 \cup \theta_2) : R_1 \cup R_2 \to [0, 1]$ is defined by

$(\mu + \theta)(x) = (\mu_1 \cup \theta_1)(x_1) \cup (\mu_2 \cup \theta_2)(x_2) =$

$$\sup_{\substack{x_1 = y_1 + z_1 \\ x_2 = y_2 + z_2}} [\min(\mu_1(y_1), \min \theta_1(z_1)) \cup (\min \mu_2(y_2), \min \theta_2(z_2))],$$

where $x_1, y_1, z_1 \in R_1$ and $x_2, y_2, z_2 \in R_2$.

**DEFINITION 4.4.19:** *Let $\mu$ be a S-fuzzy sub-biring (S-fuzzy bi-ideal) of a S-biring $R = R_1 \cup R_2$, $t \in [0, 1]$ and $t \leq \mu(0) = \mu_1(0) \cup \mu_2(0)$.*

*The sub-biring (bi-ideal) $\mu_t$ is called a S-level sub-biring (S-level bi-ideal) of $\mu$.*

**DEFINITION 4.4.20:** *Let $\mu$ be a fuzzy subset of a S-biring $R$. The Smarandache smallest fuzzy sub-biring (Smarandache smallest fuzzy bi-ideal) of $R = R_1 \cup R_2$ containing $\mu$ is called the S-fuzzy sub-biring (S-fuzzy bi-ideal) generated by $\mu$ in $R$ and is denoted by $\langle \mu \rangle$.*

**DEFINITION 4.4.21:** *Let $F$ be a S-homomorphism from a S-biring $R$ onto a S-biring $R'$ if and only if for $\mu$ and $\sigma$ any S-fuzzy bi-ideals of $R$. We have*

$$\begin{aligned} f(\mu + \sigma) &= f(\mu) + f(\sigma) \\ f(\mu\sigma) &= f(\mu)f(\sigma) \\ f(\mu \cap \sigma) &\subseteq f(\mu) \cap f(\sigma) \end{aligned}$$

*with equality if at least one of $\mu$ or $\sigma$ is f-invariant.*

Now we proceed on to define Smarandache fuzzy bi coset of a S-biring $R$.

**DEFINITION 4.4.22:** *Let $\mu$ be any S-fuzzy bi-ideal of the S-biring $R$. The fuzzy subset $\mu_x^* = \left(\mu_1^* \cup \mu_2^*\right)_x$ of $P = P_1 \cup P_2 \subset R_1 \cup R_2$ defined by for $x_1 \in P_1$ and $x_2 \in P_2$*

$$\left(\mu_1^*\right)_{x_1}(r_1) = \mu_1(r_1 - x_1) \text{ for all } r_1 \in P_1 \text{ and}$$

$$\left(\mu_2^*\right)_{x_2}(r_2) = \mu_2(r_2 - x_2) \text{ for all } r_2 \in P_2$$

*is termed as the Smarandache fuzzy bicoset (S-fuzzy bicoset) determined by $(x_1, x_2)$ and $\mu = \mu_1 \cup \mu_2$.*

We give the S-fuzzy prime ideal in yet another form.

**DEFINITION 4.4.23:** *Let $\mu = \mu_1 \cup \mu_2$ be any S-fuzzy bi-ideal of a S-biring $R$ is said to be a Smarandache fuzzy bilevel ideal (S-fuzzy bilevel ideal) if and each level bi-ideal $\mu_t = (\mu_1 \cup \mu_2)_t$, $t \in Im\ \mu$ is prime $\mu_1(x_1) < \mu_1(y_1)$ and $\mu_2(x_2) < \mu_2(y_2)$ for some $x_1, y_1 \in P_1$ ($\mu_1$ is defined relative to $P_1 \subset R_1$ and $\mu_2$ is defined relative to $P_2 \subset R_2$) then $\mu_1(x_1 y_1) = \mu_1(y_1)$ and $\mu_2(x_2 y_2) = \mu_2(y_2)$.*

**DEFINITION 4.4.24:** *A S-fuzzy bi-ideal of a S-biring $R$ is called a Smarandache fuzzy maximal bi-ideal (S-fuzzy maximal bi-ideal) if $Im\ \mu = \{1, \alpha\}$ where $\alpha \in [0,1)$ and the level bi-ideal $(\mu_1 \cup \mu_2)_t = \{x_1 \in P_1 \subset R_1 \text{ and } x_2 \in P_2 \subset R_2\ /\ \mu_1(x_1) = 1, \mu_2(x_2) = 1\}$ is maximal.*

**DEFINITION 4.4.25:** *A fuzzy bi-ideal $\mu = \mu_1 \cup \mu_2$ of a S-biring, $R = R_1 \cup R_2$ is called Smarandache fuzzy semiprime (S-fuzzy semiprime) if for any S-fuzzy bi-ideal $\theta$ of R the condition $\theta^n \subset \mu$ (i.e. $\theta = \theta_1 \cup \theta_2$; $\theta_1^n \subset \mu_1$ and $\theta_2^n \subset \mu$) implies that $\theta \subset \mu$, where $n \in Z_+$.*

**DEFINITION 4.4.26:** *Let R be a S-biring. R is biregular if for each element $x_1$, $y_1$ in $P_1 \subset R_1$ and $x_2$, $y_2$ in $P_2 \subset R_2$ where $P_1$ and $P_2$ are subfields in $R_1$ and $R_2$ respectively; we have $x_1 y_1 x_1 = x_1$, $x_2 y_2 x_2 = x_2$.*

**DEFINITION 4.4.27:** *A S-fuzzy bi-ideal $\mu = \mu_1 \cup \mu_2$ of a S-biring R is called Smarandache fuzzy primary (S-fuzzy primary) if for any two S-fuzzy bi-ideals $\sigma$ and $\theta$ of R the conditions $\sigma\theta \subset \sqrt{\mu}$ and $\sigma \not\subset \mu$ together imply $\theta \subset \sqrt{\mu}$.*

**DEFINITION 4.4.28:** *A S-fuzzy bi-ideal $\mu$ of a S-biring R is called S-fuzzy semiprimary if $\sqrt{\mu}$ is a S-fuzzy prime bi-ideal of R.*

**DEFINITION 4.4.29:** *A S-fuzzy bi-ideal $\mu$ of a S-biring R is called S-fuzzy irreducible if it is not an intersection of two S-fuzzy bi-ideals of R properly containing $\mu = \mu_1 \cup \mu_2$ otherwise $\mu$ is termed as S-fuzzy reducible.*

Now one can study all the properties which we have defined and described for S-birings can be easily extended to the case of Smarandache non associative birings.

Here we mention a few of the concepts and the rest of the development is left for the reader as research.

From now onwards we denote a Smarandache non-associative biring by S-na-biring. We give only few definitions.

**DEFINITION 4.4.30:** *Let $(R = R_1 \cup R_2, +, \bullet)$ be a S-na-biring that is $P = P_1 \cup P_2$ be a S-biring of R. A fuzzy subset $\mu : R \rightarrow [0, 1]$ is said to be a S-na-fuzzy sub-biring of R if for all x, y $\in P = P_1 \cup P_2$ we have*

    i.   $\mu (x - y) \geq \min (\mu (x), \mu (y))$ *and*
    ii.  $\mu (xy) \geq \min (\mu (x), \mu (y))$

*or equivalently we can say if $\mu = \mu_1 \cup \mu_2$ then $\mu_1$ is a S-fuzzy subring of $P_1$ and $\mu_2$ is a S-fuzzy subring of $P_2$.*

*So the S-na-fuzzy sub-biring is associated with an associative S-biring contained in the S-na-biring.*

*If in the definition of S-na-fuzzy sub-biring $\mu = \mu_1 \cup \mu_2$ only one of $\mu_1$ or $\mu_2$ (say) $\mu_1$ is a S-fuzzy subring and $\mu_2$ is just a fuzzy subring then we call $\mu$ a Smarandache na-fuzzy weakly sub-biring (S-na-fuzzy weakly sub-biring) of R.*

The following theorem is direct and hence left for the reader as an exercise.

**THEOREM 4.4.8:** *Let $(R = R_1 \cup R_2, +, \cdot)$ be a S-na biring. $\mu$ be a S-na-fuzzy sub-biring of R then $\mu$ is a S-na –fuzzy weakly sub-biring of R.*

**DEFINITION 4.4.31:** *Let R be a S-na-biring. A fuzzy subset $\mu$ of a S-biring R is called a Smarandache na fuzzy bi-ideal (S-na fuzzy bi-ideal) of R if $\mu = \mu_1 \cup \mu_2$ is defined relative to $P = P_1 \cup P_2$, P an associative biring in R.*

*If $\mu_1(x_1 - y_1) \geq \min(\mu_1(x_1), \mu_1(y_1))$, $\mu_1(x_1 y_1) \geq \max(\mu_1(x_1), \mu_1(y_1))$ and $\mu_2(x_2 - y_2) \geq \min(\mu_2(x_2), \mu_2(y_2))$, $\mu_2(x_2 y_2) \geq \max(\mu_2(x_2), \mu_2(y_2))$ for all $x_1, y_1 \in P_1$ and $x_2, y_2 \in P_2$ or equivalently we can say $\mu_1$ is a S-fuzzy ideal of $P_1$ and $\mu_2$ is a S-fuzzy ideal of $P_2$ or still we can say $\mu_1$ is a S-fuzzy bi-ideal of $P_1 \cup P_2$.*

*If $\mu : R_1 \cup R_2 \to [0, 1]$ is a fuzzy subset such that $\mu_1 : P_1 \to [0, 1]$ is a S-fuzzy ideal of $P_1$ and $\mu_2 : P_2 \to [0, 1]$ is a fuzzy ideal of $P_2$ then we call $\mu = \mu_1 \cup \mu_2 : P = P_1 \cup P_2 \to [0, 1]$ as a S-na-fuzzy weakly bi-ideal of R.*

*Thus we can say if in the definition of S-na-fuzzy bi-ideal of R only if $\mu_1$ or $\mu_2$ is just a fuzzy ideal and not a S-fuzzy ideal then we call $\mu$ a S-na – fuzzy weakly bi-ideal of R.*

**THEOREM 4.4.9:** *Let $R = R_1 \cup R_2$ be a S-na-biring. $\mu : R \to [0, 1]$ be a S-na-fuzzy bi-ideal of R then $\mu$ is a S-na-fuzzy weakly bi-ideal of R.*

*Proof*: Straightforward by very definitions, hence left for the reader as an exercise to prove.

**DEFINITION 4.4.32:** *Let $(R = R_1 \cup R_2, +, \cdot)$ be a S-na-biring. A fuzzy set $\mu : R \to [0, 1]$ is called a Smarandache na-fuzzy prime bi-ideal (S-na-fuzzy prime bi-ideal) relative to $P \subset R$ if $\mu = \mu_1 \cup \mu_2$ is a such that for any two S-na-fuzzy bi-ideals $\sigma = \sigma_1 \cup \sigma_2$ and $\theta = \theta_1 \cup \theta_2$ of R relative to the same P, the condition $\sigma_1 \theta_1 \subset \mu_1$ and $\sigma_2 \theta_2 \subset \mu_2$ implies that either $\sigma_1 \subset \mu_1$, $\sigma_2 \subset \mu_2$ or $\theta_1 \subset \mu_1$ and $\theta_2 \subset \mu_2$.*

All notions like S-na-fuzzy primary bi-ideal, S-na-fuzzy maximal bi-ideal, S-na-fuzzy semiprimary bi-ideal and S-na-fuzzy irreducible bi-ideal can be defined and all properties studied as in case of S-birings. As the aim of the text is to make a researcher work on these notions in several places we just define and give the possible properties for the reader to prove. The last section of each chapter serves this purpose.

### 4.5. Problems

This section gives problems about S-fuzzy rings and their generalizations. About 125 problems regarding S-fuzzy rings and S-fuzzy birings are given.

**Problem 4.5.1:** Can the S-ring $Z_6$ have S-fuzzy ideals?

**Problem 4.5.2:** Find all S-fuzzy ideals of the S-ring $Z_{30}$.

**Problem 4.5.3:** Let $Z_n = \{0,1,2,\ldots, n-1\}$ be a S-ring, $n = p_1^{\alpha_1} p_2^{\alpha_2} \cdots p_t^{\alpha_t}$, $t > 2$, $\alpha_i > 1$; $i = 1, 2, \ldots, t$.

How many S-fuzzy ideals can $Z_n$ have?

**Problem 4.5.4:** Characterize those S-rings which has no S-fuzzy ideals.

**Problem 4.5.5:** Characterize those S-rings which has always non trivial S-fuzzy ideals.

**Problem 4.5.6:** Obtain some interesting results between S-fuzzy ideals and level subsets of a S-ring R.

**Problem 4.5.7:** Is it true that if $\mu$ is a S-fuzzy ideal of R, then $\mu$ is S-fuzzy prime if the ideal $\mu_t$, where $t = \mu(0)$ is a S-prime ideal of R?

**Problem 4.5.8:** A S-ideal P of a S-ring R, $P \neq R$ is S-prime if and only if $\chi_P$ is a S-fuzzy prime ideal of R.

**Problem 4.5.9:** Prove or disprove a non constant S-fuzzy ideal $\mu$ of a S-ring R is S-fuzzy prime if and only if card Im $\mu = 2$, $1 \in$ Im $\mu$ and the S-ideal $\mu_t$ where $t = \mu(0)$ is S-prime.

**Problem 4.5.10:** Prove two S-level subrings $\mu_s$ and $\mu_t$ (with $s < t$) of a S-fuzzy subring $\mu$ of a S-ring R are equal if and only if there is no x in R such that $s \leq \mu(x) < t$.

**Notation:** $F_\mu = \{\mu_t \mid t \in \text{Im } \mu\}$ If Im $\mu = \{t_0, \ldots, t_n\}$ with $t_0 > t_1 > \ldots > t_n$ then we have the chain $\mu_{t_0} \subset \mu_{t_1} \subset \mu_{t_2} \subset \cdots \subset \mu_{t_n} = R$.

**Problem 4.5.11:** Prove for two S-fuzzy subrings (S-fuzzy ideals) $\mu$ and $\theta$ of a S-ring such that card Im $\mu < \infty$, card Im $\theta < \infty$ are equal if and only if Im$\mu$ = Im$\theta$ and $F_\mu = F_\theta$.

**Problem 4.5.12:** If A is a S-subring (S-ideal) of a S-ring R, $A \neq R$ then the fuzzy subset $\mu$ of R defined by

$$\mu(x) = \begin{cases} s & \text{if } x \in A \\ t & \text{if } x \in R \setminus A \end{cases}$$

where $s, t \in [0, 1]$, prove for $s > t$ is a S-fuzzy subring (S-fuzzy ideal) of the S-ring.

**Problem 4.5.13:** Prove or disprove if a non-empty subset P of a S-ring R is a S-subring (S-ideal) of R if and only if $\chi_P$ is a S-fuzzy subring (S-fuzzy ideal) of R.

**Problem 4.5.14:** Prove the intersection of any family of S-fuzzy subrings (S-fuzzy ideals) of a S-ring R is a S-fuzzy subring, S-fuzzy ideal of R.

**Problem 4.5.15:** If μ is a S-fuzzy ideal of a S-ring R will μ + μ = μ?

**Problem 4.5.16:** Let μ be a S-fuzzy subring and θ is any S-fuzzy ideal of a S-ring R then prove μ ∩ θ is a S-fuzzy ideal of the S-ring.

**Problem 4.5.17:** Let μ be any S-fuzzy ideal of a S-ring R such that Imμ = {t} or {0, s} where t ∈ [0, 1] and s ∈ (0, 1]. If μ = σ ∪ θ, where σ and θ are S-fuzzy ideals of R then prove either σ ⊆ θ or θ ⊆ σ.

**Problem 4.5.18:** Let μ be any S-fuzzy ideal of a S-ring R such that Im μ has three or more elements or Im μ = {t, t'} where t, t' ∈ (0, 1], t > t'. Does there exists S-fuzzy ideals σ and θ of R such that μ = σ ∪ θ, σ ⊄ θ and θ ⊄ σ?

**Problem 4.5.19:** Let f be a S-ring homomorphism from R to R'. If μ and σ are S-fuzzy ideals of R prove the following are true.

    i.    $f(\mu + \sigma) = f(\mu) + f(\sigma)$.
    ii.    $f(\mu\sigma) = f(\mu) f(\sigma)$.
    iii.    $f(\mu \cap \sigma) \subset f(\mu) \cap f(\sigma)$

with equality if atleast one of μ or σ is f-invariant.

**Problem 4.5.20:** Prove if μ is a S-fuzzy ideal of a S-ring R; then $R_\mu$ the set of all S-fuzzy cosets of μ in R is a S-ring under the compositions

$$\mu_x^* + \mu_y^* = \mu_{x+y}^*$$
$$\mu_x^* + \mu_y^* = \mu_{xy}^*$$

for all x ∈ P ⊂ R (P a subfield of R).

**Problem 4.5.21:** Prove or disprove if μ is a constant $R_\mu = (\mu_0^*)$ where μ is a S-fuzzy ideal relative to a fixed subfield P in R.

**Problem 4.5.22:** Prove or disprove if μ is a S-fuzzy ideal of a S-ring R then μ(x) = μ(0) ⇔ $\mu_x^* = \mu_0^*$ where x ∈ P ⊂ R (μ defined relative to the subfield P in R).

**Problem 4.5.23:** For any S-fuzzy ideal μ of a S-ring R. Prove $R\big|\mu_t \cong R_\mu$, where t = μ(0).

**Problem 4.5.24:** Can we say if μ is any S-fuzzy ideal of a S-ring R and R is S-Noetherian or S-Artinian then so is $R_\mu$?

**Problem 4.5.25:** Prove if f is a S-homomorphism from a S-ring R onto a S-ring R' and μ' is any S-fuzzy prime ideal of R' then prove $f^{-1}(\mu')$ is a S-fuzzy prime ideal of R.

**Problem 4.5.26:** Prove

i. $(\lambda \cup \mu)_t = \lambda_t \cup \mu_t$.
ii. $(\lambda \cap \mu)_t = \lambda_t \cap \mu_t$.
iii. $\{\lambda \cup (\mu \cap \nu)\}_t = (\lambda \cup \mu)_t \cap (\lambda \cup \nu)_t$.
iv. $\{\lambda \cap (\mu \cup \nu)\}_t = (\lambda \cap \mu)_t \cup (\lambda \cap \nu)_t$.
v. $(\overline{\lambda}_t) \subset \{(\overline{\lambda})_t \cap (\overline{\lambda})_{1-t}\}$ where $\overline{\lambda}$ denotes the complement of λ.
vi. $(\overline{\lambda}_t) \subset \{(\overline{\lambda})_t \cup (\overline{\lambda})_{1-t}\}$.
vii. $(\overline{\lambda_1 \cup \mu_t}) \subset (\lambda)_{1-t} \cap (\mu)_{1-t}$.

**Problem 4.5.27:** Is μ|s | δ|s a S-L-fuzzy ring? Justify your claim.

**Problem 4.5.28:** Prove or disprove if δ is a S-prime L fuzzy ideal of μ, then $\delta(x_1, \ldots, x_n) \wedge \mu(x_1) \wedge \mu(x_2) \ldots \wedge \mu(x_n) \leq \delta(x_1) \vee \ldots \vee \delta(x_n)$ for all $x_1, \ldots, x_n \in P \subset R$ (P is a subfield of R relative to which μ is defined).

**Problem 4.5.29:** If L is a totally ordered set and δ is a S-prime L-fuzzy ideal of μ, then δ is a S-prime ideal of $\mu_t$ for every t in L such $t \leq \delta(0)$ prove.

**Problem 4.5.30:** Let μ be a S-L-fuzzy ring and δ a S-prime fuzzy ideal of μ. If $\mu_1$ is a S-fuzzy subring of μ and $\delta(0) \wedge \mu_t(0) \neq 0$ then $\delta \cap \mu_t$ is a S-prime fuzzy ideal of $\mu_t$. Justify.

**Problem 4.5.31:** Let μ and μ' be S-fuzzy rings f: R → R' be a S-ring homomorphism and δ' is a S-prime fuzzy ideal of μ'; prove $f^{-1}(\delta^{-1})$ is a S-prime fuzzy ideal of μ.

**Problem 4.5.32:** Is $S_r(\delta)(x) = \bigvee_{n \in N} \delta(x^n)$ equivalent to $Sr(\delta)(x) = \sup\{t \mid x \in Sr(\delta_t)\}$?

**Problem 4.5.33:** Prove if δ is a S-primary L-fuzzy ideal and L if finite and totally ordered, then prove $\delta_t$ is a S-primary ideal of $\mu_t$ for every $t \in L$.

**Problem 4.5.34:** Prove if δ is a S-prime fuzzy ideal then δ is S-primary.

**Problem 4.5.35:** If δ is a S-primary prove Sr (δ) is a S-fuzzy prime ideal.

**Problem 4.5.36:** If μ is any non constant S-fuzzy irreducible ideal of a S-ring R, are the following results true?

   i. $1 \in \text{Im } \mu$.
   ii. There exists $a \in [0,1)$ such that $\mu(x) = \alpha$ for all $x \in P' \subset R$ where $P' = P \setminus \{x \in P \mid \mu(x) = 1\}$ here P is a subfield of the S-ring R relative to which μ is defined.
   iii. The S-ideal $\{x \in R \mid \mu(x) = 1\}$ is irreducible.

**Problem 4.5.37:** Prove a S-fuzzy ideal of a S-commutative ring with unity is S-fuzzy semiprime if and only if each of its S-level ideals is a semiprime ideal of the given ring.

**Problem 4.5.38:** Prove a S-commutative ring with unity is S-regular if and only if each of its S-fuzzy ideals is S-fuzzy semiprime.

**Problem 4.5.39:** Let R be a S-commutative ring with unity. If μ is any S-fuzzy ideal of R which is both S-fuzzy semiprime and S-fuzzy irreducible prove μ is S-fuzzy prime.

**Problem 4.5.40:** Prove in a S-commutative regular ring with unity every S-fuzzy irreducible ideal is S-fuzzy prime.

**Problem 4.5.41:** $J : R \to L = [0, 1]$ be a S-fuzzy, ideal and x, y, u, v be any element in R. If $x + J = u + J$ and $y + J = v + J$, then prove $(x + y) + J = (u + v) + J$ and $xy + J = uv + J$.

**Problem 4.5.42:** Prove every S-fuzzy ideal A of R such that $A(0) = 1$ and A is finite valued has a S-fuzzy primary representations if and only if every S-ideal of R has a S-primary representation.

**Problem 4.5.43:** Prove every S-fuzzy ideal A of R such that $A(0) = 1$ has a S-fuzzy primary representation if and only if R is S-artinian.

**Problem 4.5.44:** Let $A \in S\ F$. Suppose that a has a S-reduced fuzzy primary representation $A = Q_1 \cap \ldots \cap Q_n$. Prove A has a finite set $(P_1, \ldots, P_m)$, $m \leq n$, of isolated S-fuzzy prime ideal divisors. Further more. $\sqrt{A} = P_1 \cap \ldots \cap P_n$.

**Problem 4.5.45:** Let $A \in SF$. Suppose that A has a S-reduced fuzzy primary representation $A = Q_1 \cap \ldots \cap Q_n$. Prove $\text{Im}(A) = \bigcup_{i=1}^{n} \text{Im}(Q_i)$.

**Problem 4.5.46:** Let R be a S-quasi local ring.

i. If μ is a S-fuzzy quasi local subring of R then prove μ* is quasi local and μ* ∩ M is the unique S-maximal ideal of μ*.

ii. If R' is a S-subring of R, prove R' is a S-quasi local with unique S-maximal ideal M ∩ R' if and only if $\delta_{R'}$ is a S-fuzzy quasi local subring of R.

**Problem 4.5.47:** Let $\{A_\alpha \, / \, \alpha \in \Omega\} \cup A$ be the collection of S-fuzzy subrings of a S-ring R such that $A = \sum_{\alpha \in \Omega} A_\alpha$. Suppose L has finite intersection property. Prove $A^* = \bigoplus_{\alpha \in \Omega} A_\alpha^*$ if and only if $A = \bigoplus_{\alpha \in \Omega} A_\alpha$.

**Problem 4.5.48:** Let $\{A_\alpha / \alpha \in \Omega\}$ be a collection of S-fuzzy subrings of a S-ring R. Prove $\sum_{\alpha \in \Omega}(A_\alpha)_* \subseteq \left(\sum_{\alpha \in \Omega} A_\alpha\right)_*$.

**Problem 4.5.49:** Let $\{A_\alpha \, / \, \alpha \in \Omega\}$ be a collection of all S-fuzzy subrings of a S-ring. If there exists $t \in L$, $t \neq 1$, such that $t \geq \sup \{A_\alpha(x)| \, x \notin A_{\alpha*}$ for all $\alpha \in \Omega\}$, then prove

$$\sum_{\alpha \in \Omega}(A_{\alpha*}) = \left(\sum_{\alpha \in \Omega} A_\alpha\right)_*.$$

**Problem 4.5.50:** Let $\{A_\alpha| \, \alpha \in \Omega\}$ be a finite collection of S-fuzzy subring of R. If $A_\alpha$ is finite-valued for all $\alpha \in \Omega$, then prove

$$\sum_{\alpha \in \Omega}(A_{\alpha*}) = \left(\sum_{\alpha \in \Omega} A_\alpha\right)_*.$$

**Problem 4.5.51:** Let $\{A_\alpha| \, \alpha \in \Omega\}$ be a finite collection of S-fuzzy subrings of R. Suppose that L has finite intersection. If

$$\sum_{\alpha \in \Omega} A_\alpha = \bigoplus_{\alpha \in \Omega} A_\alpha$$

Prove

$$\sum_{\alpha \in \Omega} A_{\alpha^*} = \bigoplus_{\alpha \in \Omega} A_{\alpha^*}$$

**Problem 4.5.52:** Let $\{A_\alpha| \, \alpha \in \Omega\}$ be a collection of S-fuzzy subrings of R. Then prove

$$A_\beta \sum_{\alpha \in \Omega} A_\alpha = \sum_{\alpha \in \Omega} A_\beta A_\alpha$$

and

$$\left(\sum_{\alpha \in \Omega} A_\alpha\right) A_\beta = \sum_{\alpha \in \Omega} A_\alpha A_\beta$$

**Problem 4.5.53:** Let A and B be S-fuzzy subrings of R such that sup{sup{A(x) | x ∉ $A_\alpha$}, sup {B(x)| α ∉ $B_\alpha$}} < 1. Prove $(AB)_* = A_* B_*$.

**Problem 4.5.54:** Let A be a fuzzy subset of V and let s, t ∈ Im (A). Prove

  i.  s ≤ t if and only if $A_s \supseteq A_t$.
  ii. s = t if and only if $A_s = A_t$.

**Problem 4.5.55:** Let A ∈ $\mathcal{A}_{\delta\mathcal{F}}$. Then for all t such that 0 ≤ t ≤ A(0), prove $A_t$ is a S-subspace of V.

**Problem 4.5.56:** Let A be fuzzy subset of the S-vector space V. If $A_t$ is a S-subspace of V for all t ∈ Im(A) will A ∈ $\mathcal{A}_{\delta\mathcal{S}}$.

[Hint: For S a subset of F, we let $\delta_S$ denote the characteristic function of S].

**Problem 4.5.57:** Let A be a fuzzy subset of a S-fuzzy vector space V and let K be a fuzzy subset of the S-ring R. For d ∈ P ⊂ R and x ∈ V (V defined relative to the subfield P in R). Suppose that 0 ≤ μ, λ ≤ 1. Then for all z ∈ V, prove

  i.   $(d_\mu \circ A)(z) = \min\{\mu, A(\frac{1}{d}z)\}$ if d ≠ 0.
  ii.  $(O_\mu \circ A)(z) = \sup\{\min\{\mu, A(y)\}/ y \in V\}$ if z = 0; $(O_\mu \circ A)(z) = 0$ if z ≠ 0.
  iii. $(k \circ x_\lambda)(z) = \sup\{\min\{k(c), \lambda\} / c \in P \subset R, z = cx\}$ if x ≠ 0 and z ∈ sp(x); 0 if x ≠ 0 and z ∉ sp(x).
  iv.  $(k \circ O_\lambda)(z) = \sup\{\min\{k(c), \lambda / c \in F\}$ if z = 0; $(k \circ 0_\lambda)(z)$ if z ≠ 0.

**Problem 4.5.58:** Prove the following. Let c, d ∈ P ⊂ R, x, y ∈ V and 0 ≤ κ, λ, μ, ν ≤ 1. Then $d_\mu \circ x_\lambda = (dx)_{\min(\mu, \lambda)}$. $x_\lambda + y_\nu = (x + y)_{\min(\lambda, \nu)}$, $d_\mu \circ x_\lambda + c_\kappa \circ y_\nu = (dx + cy)_{\min\{\mu, \lambda, \kappa, \nu\}}$.

**Problem 4.5.59:** Show if M is a δ-Noetherian module over a S-ring R, prove for any S-fuzzy module $\eta_M^P$, η(M) < +∞.

**Problem 4.5.60:** Is it true if η(M) < ∞ for any S-module M over a S-ring R then M is a S-Noetherian module?

**Problem 4.5.61:** Prove or disprove if M is a S-left (or S-right) R module where R is an S-ring. Then prove or disprove η(M) < ∞ for every S-fuzzy evaluation η if and only if M is a S-Noetherian module.

**Problem 4.5.62:** Let V be a S-vector space II. Then prove $\eta(V) < \infty$ for every S-fuzzy vector space $\eta_V$ if and only if V is a S-finite dimensional vector space.

**Problem 4.5.63:** Let $\eta_M^P$ is a S-fuzzy module over R, ($P \subset R$ is a subfield relative to which $\eta_M$ is defined) and $N_1$ be its singular S-fuzzy submodule i.e. $N_1 = \{a \in M / \eta(a) = 1\}$. Then prove or disprove the following statements:

i. In the S-fuzzy quotient module $\left(\overline{\eta}, M/N_1\right), \overline{\eta}(\overline{a}) = 1$, if and only if $\overline{a} = \overline{0}$.
ii. In the S-fuzzy module $(\eta, M)$, $\eta(a') = \eta(a)$ if a and a' belong to the same coset.

**Problem 4.5.64:** For any S-fuzzy module $\eta_M^P$ and its singular S-fuzzy module $N_1$, prove if $|M/N_1| < \infty$ implies $\eta(M) < \infty$.

**Problem 4.5.65:** Give an example of a SNA-fuzzy ring which is not a SNA-strong fuzzy ring.

**Problem 4.5.66:** Prove if $\mu$ is a SNA-fuzzy ideal of R then

i. $\mu(x) = \mu(-x) \leq \mu(0)$ for all $x \in P \subset R$.
ii. $\mu(x - y) = \mu(0) \Rightarrow \mu(x) = \mu(y)$, $x, y \in P \subset R$.

P relative to which $\mu$ is defined.

**Problem 4.5.67:** Let R be a SNA-ring. A fuzzy subset $\mu$ of R is a SNA-fuzzy ideal of R if and only if the level subsets $\mu_t$, $t \in \text{Im } \mu$ are SNA-fuzzy ideals of R relative to the same subring $P \subset R$ – Prove.

**Problem 4.5.68:** If f is a S-homomorphism of a SNA-ring R onto a SNA-ring R', then for each SNA-fuzzy subring $\mu$ of R, $f(\mu)$ is a SNA-fuzzy subring of R' and for each SNA-fuzzy subring $\mu'$ of R', $f^{-1}(\mu')$ is a SNA-subring of R – Prove.

**Problem 4.5.69:** Verify whether the above statement is true, if SNA-fuzzy subring is replaced by the SNA-fuzzy ideal of R.

**Problem 4.5.70:** R is a SNA-ring, is it true that R is SNA-regular ring if and only if $\sigma\theta = \sigma \cap \theta$ where $\sigma$ and $\theta$ are any two SNA-fuzzy ideal of R relative to the same $P \subset R$.

**Problem 4.5.71:**

i. If μ is any SNA-fuzzy prime ideal of a SNA-ring R then prove the SNA-ideal $\mu_t$, $t = \mu(0)$, is a prime ideal of R.
ii. Prove an SNA-ideal P of a SNA-ring R, R ≠ P is prime if and only if characteristic function $\chi_P$ is a SNA-fuzzy prime ideal of R.
iii. Will Im μ = 2? Justify your claim.

**Problem 4.5.72:** SNA-fuzzy ideal μ of a SNA-ring. Prove μ is a SNA-fuzzy semiprime if for any SNA-fuzzy ideal θ of R the condition $\theta^n \subseteq \mu$ implies that $\theta \subset \mu$ where $n \in Z_+$.

**Problem 4.5.73:** Prove an ideal I of a SNA-ring R is semiprime if and only if $\chi_I$ is a SNA-fuzzy semiprime ideal of R.

**Problem 4.5.74:** If μ is any SNA-fuzzy ideal of a SNA-ring such that Im μ = $\{t_0, t_1, \ldots, t_n\}$ and $t_0 > \ldots > t_n$, then prove μ is SNA-fuzzy semiprime if and only if $\mu_{t_i}$ is a semiprime ideal of R for all i = 0, 1, 2, …, n.

**Problem 4.5.75:** Is the intersection of a SNA-fuzzy semiprime ideal of a SNA-ring relative to a fixed $P \subset R$, a SNA-fuzzy semiprime of R relative to the same fixed P? (P an associative subring or R).

**Problem 4.5.76:** Is it true that a SNA-regular ring every SNA-fuzzy ideal of R is idempotent?

**Problem 4.5.77:** Suppose the SNA-fuzzy ideal of the SNA-ring is idempotent relative to an associative subring P of R. Then P is a regular subring of R.

**Problem 4.5.78:** Prove in a SNA-ring R every SNA-fuzzy prime ideal of R relative to the associative subring P of R is a SNA-fuzzy primary ideal of R relative to the same P.

**Problem 4.5.79:** Let A be any SNA-primary ideal of a SNA-ring R, relative to $P \subset R$, A ≠ R, then the fuzzy subset μ of R defined by

$$\mu(x) = \begin{cases} 1 & \text{if } x \in A \\ \alpha & \text{if } x \in P \setminus A \text{ where } \alpha \in [0,1) \end{cases}$$

is a SNA-fuzzy primary ideal of R relative to P.

**Problem 4.5.80:** If μ is any SNA-fuzzy semiprimary ideal of a SNA-ring relative to P then prove $\mu_s$ where $s \in$ Im μ is a SNA-semiprimary ideal of R relative to P.

**Problem 4.5.81:** Let R be a SNA-ring. If μ is any SNA-fuzzy semiprimary ideal of R prove the ring $R_\mu$ is semiprimary relative to the same P.

**Problem 4.5.82:** Let μ be any SNA-fuzzy ideal of a SNA-ring R μ is SNA-fuzzy semiprime and SNA-fuzzy irreducible relative to the same associative subring P ⊂ R. Prove μ is SNA-fuzzy prime relative to P.

**Problem 4.5.83:** Give an example of a S-weakly commutative biring which is not a S-commutative biring.

**Problem 4.5.84:** Let (R = Q[x] ∪ (R × R), +, •) be a biring. Define μ : R → [0, 1] a S-fuzzy biring of R.

[Hint: R is a biring as R = $R_1$ ∪ $R_2$ where $R_1$ = Q[x] is a S-ring as Q ⊂ Q[x] is a field in Q[x], $R_2$ = $\Re \times \Re$, $\Re \times \{0\}$ ⊂ $\Re \times \Re$ is a field in $R_2$).

**Problem 4.5.85:** Give an example of a S-strong fuzzy biring.

**Problem 4.5.86:** Give an example of a S-fuzzy biring which is not a S-strong fuzzy biring.

**Problem 4.5.87:** Prove a fuzzy subset μ = $μ_1$ o $μ_2$; R = $R_1$ ∪ $R_2$ → [0, 1] of a S-fuzzy sub-biring ( S-fuzzy bi-ideal ) of R if and only if $μ_t$ = $(μ_1)_t$ ∪ $(μ_2)_t$ ∈ Im μ = Im $μ_1$ ∪ Im$μ_2$ are S-fuzzy bisubrings (S-fuzzy bi-ideals) of R.

**Problem 4.5.88:** Prove a non-constant S-fuzzy bi-ideal μ of a biring R is S-fuzzy prime if and only if Card Im μ = 2, 1 ∈ Im μ and bi-ideal $μ_t$, where t = μ(0) is prime.

**Problem 4.5.89:** Prove two levels S-sub-birings (level S-bi-ideals) $μ_s$ and $μ_t$ (with s < t) of a S-fuzzy sub-biring (S-fuzzy bi-ideal) μ of a S-biring R are equal if and only if there is no x in R such that s ≤ μ(x) < t; (s < $μ_1$(x) < t; s < $μ_2$(x) < t).

**Problem 4.5.90:** Will the intersection of any family of S-fuzzy sub-birings (S-fuzzy bi-ideals) of a S-biring a S-fuzzy sub-biring (S-fuzzy bi-ideal) of R.

**Problem 4.5.91:** If μ is any S-fuzzy bi-ideal of a S-biring R then will μ + μ = μ? Justify your claim.

**Problem 4.5.92:** If μ is any S-fuzzy sub-biring and θ is any S-fuzzy bi-ideal of R then will μ ∩ θ a S-fuzzy bi-ideal of R?

**Problem 4.5.93:** Let μ be any fuzzy subset of a S-bifield F. Then prove μ is a S-fuzzy bi-ideal of F if and only if μ(x) = μ(y) < μ(0) (i.e. $μ_1(x_1)$ = $μ_1(y_1)$ < $μ_1(0)$ and $μ_2(x_2)$ = $μ_2(y_2)$ < $μ_2(0)$ for all $x_1$, $x_2$ ∈ $F_1$ \ {0} and $y_1$, $y_2$ ∈ $F_2$\{0}; where F = $F_1$ ∪ $F_2$).

**Problem 4.5.94:** Let μ be any S-fuzzy sub-biring (S-fuzzy bi-ideal) of a S-biring R = $R_1$ ∪ $R_2$ and if $μ_1(x_1)$ < $μ_1(y_1)$ and $μ_2(x_2)$ < $μ_2(y_2)$ for some $x_1$, $y_1$ ∈ $R_1$ and $x_2$, $y_2$ ∈ $R_2$ then prove $\mu_1(x_1 - y_1) = \mu_1(x_1) = \mu_1(y_1 - x_1)$ and $\mu_2(x_2 - y_2) = \mu_2(x_2) = \mu_2(y_2 - x_2)$ where $\mu = \mu_1 \cup \mu_2$; $\mu_1 : R_1 \to [0, 1]$ and $\mu_2 : R_2 \to [0, 1]$.

**Problem 4.5.95:** Prove if $\{\mu^i = \mu_1^i \cup \mu_2^i / i \in Z_+\}$ is a collection of S-fuzzy bi-ideals of S-biring R such that $\mu_1^1 \subset \mu_1^2 \subset \mu_1^3 \subset ... \subset \mu_1^n \subset ...$ and $\mu_2^1 \subset \mu_2^2 \subset \mu_2^3 \subset ... \subset \mu_2^n \subset ...$ then

$$\bigcup_{i \in Z_+} \mu^i = \left(\bigcup_{i \in Z_+} \mu_1^i\right) \cup \left(\bigcup_{i \in Z_+} \mu_2^i\right)$$

is a S-fuzzy bi-ideal of R.

**Problem 4.5.96:** If $\mu$ is any S-fuzzy bi-ideal of a S-biring R then $\mu(x) = \mu(0) \Leftrightarrow \mu_x^* = \mu_0^*$ where $x \in R_1 \cup R_2$.

(Hint: To be more precise $\mu = \mu_1 \cup \mu_2$, $\mu_1 : P_1 \to [0, 1]$ and $\mu_2 : P_2 \to [0, 1]$ with $\mu_1(x_1) = \mu_1(0)$ and $\mu_2(x_2) = \mu_2(0)$ if and only if $\mu_{1x_1}^* = (\mu_1)_0^*$ and $\mu_{2x_2}^* = (\mu_2)_0^*$, $x_1 \in P_1$ and $x_2 \in P_2$).

**Problem 4.5.97:** Prove if $\mu$ is any S-fuzzy prime bi-ideal of a S-biring R then $\mu_1(x_1y_1)$ = max $\{\mu_1(x_1), \mu_1(y_1)\}$, $\mu_2(x_2y_2)$ = max $\{\mu_2(x_2), \mu_2(y_2)\}$ for all $x_1, y_1 \in P_1 \subset R_1$ and $x_2, y_2 \in P_2 \subset R_2$ ($R = R_1 \cup R_2$).

**Problem 4.5.98:** Prove if $\mu$ and $\theta$ are any two S-fuzzy prime bi-ideals of a S-biring R then $\mu \cap \theta$ is a S-fuzzy prime bi-ideal of R if and only if $\mu \subseteq \theta$ or $\theta \subseteq \mu$.

**Problem 4.5.99:** Prove I a S-bi-ideal of a S-biring R is semiprime if and only if $\chi_I$ is a S-fuzzy semiprime bi-ideal of R.

**Problem 4.5.100:** If $\mu$ is any S-fuzzy bi-ideal of a S-biring R such that Im $\mu = \{t_0, t_1, ..., t_n\}$ and $t_0 > t_1 > ... > t_n$ then $\mu$ is S-fuzzy semiprime if and only if $\mu_{t_i}$ is a S-semiprime bi-ideal of R for all $i = 0, 1, 2, ..., n$.

**Problem 4.5.101:** Prove the intersection of S-fuzzy semiprime (S-fuzzy prime, S-fuzzy maximal) bi-ideals of a S-biring R is always a S-fuzzy semiprime bi-ideal of R.

**Problem 4.5.102:** Prove a S-biring R is regular if and only if every S-fuzzy bi-ideal of R is idempotent.

**Problem 4.5.103:** Prove every S-fuzzy prime bi-ideal of a S-biring R is S-fuzzy primary bi-ideal of R.

**Problem 4.5.104:** Prove if $\mu$ is any S-fuzzy primary bi-ideal of a S-biring R then $\mu_t$, $t \in$ Im $\mu$ is a S-primary bi-ideal of R.

**Problem 4.5.105:** Prove if μ is any S-fuzzy primary bi-ideal of a S-biring R then √μ is S-fuzzy prime.

**Problem 4.5.106:** Prove if μ is any S-fuzzy primary bi-ideal of a S-biring R then the S-biring $R_\mu$ is primary.

**Problem 4.5.107:** If μ is any S-fuzzy semiprimary bi-ideal of a S-biring R then prove $\mu_s$ where s ∈ Im μ is a semi-primary bi-ideal of R.

**Problem 4.5.108:** Prove if μ is any S-fuzzy semiprimary bi-ideal of a S-biring R then the biring $R_\mu$ is semiprimary.

**Problem 4.5.109:** Prove if μ is any S-fuzzy prime bi-ideal of a S-biring R then μ is S-fuzzy irreducible.

**Problem 4.5.110:** Prove if R is a S-regular biring then every S-fuzzy irreducible bi-ideal is fuzzy prime.

**Problem 4.5.111:** Obtain some innovative results on S-fuzzy bi-modules on S-birings.

**Problem 4.5.112:** Give an example of a S-na-fuzzy weakly sub-biring, which is not a S-na-fuzzy sub-biring.

**Problem 4.5.113:** Give an example of a S-na-fuzzy bi-ideal.

**Problem 4.5.114:** Give an example of a S-na-fuzzy weakly bi-ideal, which is not a S-na-fuzzy bi-ideal.

**Problem 4.5.115:** Give a necessary and sufficient condition for a S-na-biring to have only S-na-fuzzy weakly bi-ideal.

**Problem 4.5.116:** Define S-na-fuzzy irreducible bi-ideal of a S-na-biring.

**Problem 4.5.117:** Give an example of

  i. a S-na-fuzzy irreducible bi-ideal.
  ii. a S-na-fuzzy reducible bi-ideal.

**Problem 4.5.118:** Define S-na-fuzzy semi-primary bi-ideal of a S-na-biring. Illustrate it by an example.

**Problem 4.5.119:** Prove if μ is any S-na-fuzzy semiprimary bi-ideal of a S-na-biring R then the biring Rμ is semiprimary.

**Problem 4.5.120:** Prove/disprove if μ is any S-na-fuzzy ideal which is both S-na-fuzzy semiprime and S-na-fuzzy irreducible then μ is S-na-fuzzy prime bi-ideal.

**Problem 4.5.121:** If R is a S-na-regular biring will every S-na-fuzzy irreducible bi-ideal of R be S-na-fuzzy prime?

**Problem 4.5.122:** Obtain a necessary and sufficient condition for a S-na-fuzzy bi-ideal to be S-na-fuzzy irreducible bi-ideal of R.

**Problem 4.5.123:** Obtain a necessary and sufficient condition for a S-na-fuzzy bi-ideal of a S-na-biring R to be S-na-fuzzy primary bi-ideal of R.

**Problem 4.5.124:** Let R be a S-na-biring. If μ is any S-na-fuzzy bi-ideal of R, will $\sqrt{\mu}$ a S-na- fuzzy maximal bi-ideal imply μ is a S-na-fuzzy primary bi-ideal of R?

**Problem 4.5.125:** If μ is a S-na-fuzzy primary bi-ideal of a S-na-biring R. Prove $\sqrt{\mu}$ is the smallest S-na-fuzzy prime bi-ideal of R containing μ.

**CHAPTER FIVE**

# SMARANDACHE FUZZY SEMIRINGS AND THEIR GENERALIZATIONS

Study of semirings and fuzzy semirings are carried out by several researchers. But the notion of Smarandache semirings is new [134]. The concept of Smarandache fuzzy semiring is introduced in this chapter. This chapter has five sections. In section one we introduce the concept of Smarandache semirings and its properties. Smarandachaic fuzzy properties and Smarandache properties are given in semirings. We give about 25 definitions recalling the Smarandache-ic fuzzy properties in semirings. In section two Smarandache fuzzy semivector spaces are defined and its properties are derived. Smarandache fuzzy non-associative semirings are introduced and studied in section three. Smarandache fuzzy bisemirings both associative and non-associative are introduced and studied in section four which contains about forty definitions. The final part of the chapter, the fifth section contains fifty-four problems of interest.

## 5.1 Smarandache semirings and its properties

The study of fuzzy semirings is dealt in chapter I. In this section we just recall the definition of Smarandache semirings and give a few of its properties. This section defines the notions of Smarandache fuzzy semirings and gives some of its basic properties. As the main motivation of the text is to introduce the Smarandache fuzzy concepts and make the reader work on these concepts we have restrained from giving several theorems with proofs, instead we propose them as problems in section five. Now we start to recollect some of the basic definitions about Smarandache semirings.

**DEFINITION 5.1.1:** *The Smarandache semiring (S-semiring) S which will be denoted from here on wards as S-semiring is defined to be a semiring S such that a proper subset B of S is a semifield with respect to the same induced operations.*

***Example 5.1.1***: $R^0$ = {set of all positive reals with 0} is a S-semiring. For $Z^0 \subset R^0$ is a semifield. Clearly we can prove all semirings in general are not S-semirings.

**DEFINITION 5.1.2**: *Let S be a semiring. A non-empty proper subset A of S is said to be a Smarandache subsemiring (S-subsemiring) if A is a S-semiring.*

In view of these definitions we have the following theorem which is left as an exercise for the reader to prove.

**THEOREM 5.1.1:** *If S is a semiring which has a non-trivial S-subsemiring then S is a S-semiring.*

**DEFINITION 5.1.3:** *Let S be a S-semiring. We say S is a Smarandache commutative semiring (S-commutative semiring) if S has a S-subsemiring which is commutative. If the S-semiring has no commutative S-subsemiring then we say S is a Smarandache non-commutative semiring (S-non-commutative semiring).*

**DEFINITION 5.1.4:** *Let S be a semiring. A non-empty subset P of S is said to be a Smarandache right (left) ideal (S-right (left) ideal) of S if the following conditions are satisfied*

  i.   *P is a S-subsemiring.*
  ii.  *For every $p \in P$ and $A \subset P$ where A is the semifield of P we have for all $a \in A$ and $p \in P$, ap (pa) is in A.*

*If P is simultaneously both S-right ideal and a S-left ideal then we say P is a Smarandache ideal (S-ideal) of S.*

The following theorem is of interest.

**THEOREM 5.1.2:** *Let S be a S-semiring. Every S-ideal of S is a S-subsemiring of S but every S-subsemiring of S need not in general be a S-ideal of S.*

*Proof*: By an example.

**Example 5.1.2**: Let $M_{2\times 2} = \left\{ \begin{pmatrix} a & b \\ c & d \end{pmatrix} \Big/ a,b,c,d \in C_2 = (0,1) \right\}$ = *set of all 2 × 2 matrices with entries from the chain lattice $C_2$.*

$$M_{2\times 2} = \left\{ \begin{pmatrix} 0 & 0 \\ 0 & 0 \end{pmatrix}, \begin{pmatrix} 0 & 1 \\ 0 & 0 \end{pmatrix}, \begin{pmatrix} 1 & 0 \\ 0 & 0 \end{pmatrix}, \begin{pmatrix} 1 & 1 \\ 0 & 0 \end{pmatrix}, \right.$$
$$\begin{pmatrix} 0 & 1 \\ 0 & 1 \end{pmatrix}, \begin{pmatrix} 0 & 0 \\ 1 & 1 \end{pmatrix}, \begin{pmatrix} 1 & 0 \\ 0 & 1 \end{pmatrix}, \begin{pmatrix} 0 & 1 \\ 1 & 0 \end{pmatrix},$$
$$\begin{pmatrix} 1 & 1 \\ 1 & 0 \end{pmatrix}, \begin{pmatrix} 0 & 1 \\ 1 & 1 \end{pmatrix}, \begin{pmatrix} 1 & 0 \\ 1 & 1 \end{pmatrix}, \begin{pmatrix} 1 & 1 \\ 0 & 0 \end{pmatrix},$$
$$\left. \begin{pmatrix} 0 & 1 \\ 1 & 1 \end{pmatrix}, \begin{pmatrix} 1 & 0 \\ 1 & 1 \end{pmatrix}, \begin{pmatrix} 1 & 1 \\ 0 & 1 \end{pmatrix}, \begin{pmatrix} 1 & 1 \\ 1 & 1 \end{pmatrix} \right\}.$$

*Clearly $M_{2\times 2}$ is a S-semiring with 16 elements in it. For*

$$A = \left\{ \begin{pmatrix} 1 & 0 \\ 0 & 1 \end{pmatrix}, \begin{pmatrix} 0 & 0 \\ 0 & 0 \end{pmatrix} \right\}$$

is a semifield. To find ideals in $M_{2\times 2}$. The set

$$S = \left\{ \begin{pmatrix} 0 & 0 \\ 0 & 0 \end{pmatrix}, \begin{pmatrix} 1 & 0 \\ 0 & 0 \end{pmatrix} \right\}$$

is a subsemiring of $M_{2\times 2}$. Clearly S is not a S-subsemiring so S cannot be a S-ideal of $M_{2\times 2}$.

$$A = \left\{ \begin{pmatrix} 1 & 0 \\ 0 & 0 \end{pmatrix}, \begin{pmatrix} 1 & 0 \\ 0 & 1 \end{pmatrix}, \begin{pmatrix} 0 & 0 \\ 0 & 0 \end{pmatrix} \right\}$$

is a S-subsemiring of $M_{2\times 2}$, which is not a S-ideal of $M_{2\times 2}$.

***DEFINITION 5.1.5:*** Let S be a semiring. A non-empty proper subset A of S is said to be Smarandache pseudo-subsemiring (S-pseudo subsemiring) if the following condition is true.

If there exists a subset P of S such that $A \subset P$; where P is a S-subsemiring i.e. P has a subset B such that B is a semifield under the operations of S or P itself is a semifield under the operations of S.

***THEOREM 5.1.3:*** Let S be a semiring every proper subset of S need not in general be a S-pseudo subsemiring.

*The reader is expected to prove by an example.*

***THEOREM 5.1.4:*** Let S be a semiring if S contains a S-pseudo subsemiring then S is a S-semiring.

Proof: *S is a semiring and S contains a S-pseudo subsemiring A; i.e. A is contained in a semifield P, P contained in S. So S is a S-semiring.*

*The concept of S-pseudo subsemiring leads to the definition of S-pseudo ideals in the semiring.*

***DEFINITION 5.1.6:*** Let S be a semiring. A non-empty subset P of S is said to be a Smarandache pseudo right (left) ideal (S-pseudo right (left) ideal) of the semiring S if the following conditions are true.

  i.   P is a S-pseudo subsemiring i.e. $P \subset A$, A a semifield in S.
  ii.  For every $p \in P$ and every $a \in A$, $ap \in P$ ($pa \in P$).

If P is simultaneously both a S-pseudo right ideal and a S-pseudo left ideal we say P is a Smarandache pseudo ideal (S-pseudo ideal).

*Now we define two new notions about S-semiring viz. Smarandache dual ideal and Smarandache pseudo dual ideal of a semiring S.*

***DEFINITION 5.1.7:*** Let S be a semiring. A non-empty subset P of S is said to be a Smarandache dual ideal (S-dual ideal) of S if the following conditions hold good.

  i.   P is a S-subsemiring.
  ii.  For every $p \in P$ and $a \in A\setminus\{0\}$ $a + p$ is in A, where $A \subset P$.

***DEFINITION 5.1.8:*** Let S be a semiring. A non-empty subset P of S is said to be a Smarandache pseudo dual ideal (S-pseudo dual ideal) of S if the following conditions are true.

i. P is a S-pseudo subsemiring i.e. P ⊂ A, A is a semifield in S or A contains a semifield.

ii. For every p ∈ P and a ∈ A, p + a ∈ P. Clearly P is simultaneously left and right S-pseudo dual ideal of S as S is additively commutative.

*The reader is assigned the work of finding examples. Of course, several examples exist. Here we introduce yet another nice substructure in a semiring called a Smarandache semidivision ring.*

**DEFINITION 5.1.9:** Let S be a S-semiring. S is said to be a Smarandache semidivision ring (S-semidivision ring) if the proper subset A ⊂ S is such that

i. A is a Smarandache subsemiring which is non-commutative.
ii. A contains a subset P such that P is a semidivision ring that is P has no zero divisors and P is a non-commutative semiring.

Now we proceed on to define the notion of Smarandache fuzzy semirings.

**DEFINITION 5.1.10:** *A fuzzy subset µ of a S-semiring S is called a Smarandache fuzzy semiring (S-fuzzy semiring) relative to P ⊂ S where P is a field if for all x, y ∈ P*

$$\mu(x + y) \geq \min(\mu(x), \mu(y)) \text{ and}$$
$$\mu(xy) \geq \min(\mu(x), \mu(y)).$$

*Thus every S-fuzzy semiring µ will be associated with a semifield P contained in S. Further we see µ need not be a S-fuzzy semiring relative to all fuzzy subsets µ on a S-semiring S.*

*We say a fuzzy subset µ of a S-semiring S to be a Smarandache strong fuzzy semiring (S-strong fuzzy semiring) if µ relative to every proper subset P which is a subsemifield of S is a S-fuzzy semiring.*

**THEOREM 5.1.5:** *Let µ be a fuzzy subset of a S-semiring S. If µ is a S-strong fuzzy semiring then µ is a S-fuzzy semiring.*

*Proof***:** Direct by the very definition.

Now we define the notion of Smarandache fuzzy ideal of a semiring.

**DEFINITION 5.1.11:** *Let S be a S-semiring. A fuzzy subset µ of a S-semiring R is called a Smarandache fuzzy ideal (S-fuzzy ideal) of R relative to a proper subset P of R where P is a semifield if µ satisfies the following conditions.*

$$\mu(x + y) \geq \min(\mu(x), \mu(y))$$
$$\mu(xy) \geq \max(\mu(x), \mu(y))$$

*for all x, y ∈ P ⊂ R.*

***Example 5.1.3:*** Consider the S-semiring which is a distributive lattice whose Hasse diagram is as follows:

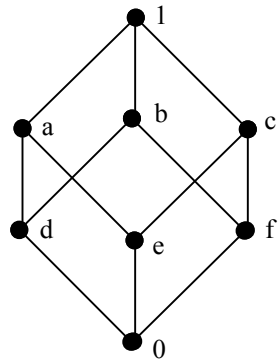

**Figure 5.1.1**

This has several semifields and we can define S-fuzzy ideals accordingly. Product of S-fuzzy ideals of a S-semiring S is defined as in case of S-rings relative to a proper subset which is a semifield.

Now we proceed on to define Smarandache fuzzy prime ideals of a S-semiring S.

**DEFINITION 5.1.12:** *A S-fuzzy ideal μ of a S-semiring S is called a Smarandache fuzzy prime (S-fuzzy prime) if the ideal $\mu_t$ where t = (0) is a prime ideal of P ⊂ S (P a subset of S which is a semifield) or following [116] we define the S-fuzzy prime ideal of a S-semiring in a different way.*

*A non-constant-S-fuzzy ideal μ of a S-semiring S is called S-fuzzy prime if for any two S-fuzzy ideals σ and θ of P ⊂ S (μ defined only relative to the semifield P in S) the condition σ θ ⊂ μ implies σ ⊂ μ or θ ⊂ μ.*

*Recall, μ be any S-fuzzy subsemiring of the S-semiring S, t ∈ [0, 1] and t ≤ μ (0). The S-subsemiring (S-ideal) $\mu_t$ is called a S-level subsemiring (S-level ideal) of μ. We just recall the definition of S-fuzzy quotient ideal of a S-semiring S.*

**DEFINITION 5.1.13:** *Let μ be a S-fuzzy ideal of a S-semiring S relative to a subsemifield P in S. The S-fuzzy ideal μ' of $R_\mu$ defined by $\mu'(\mu_x^*) = \mu(x)$ for all x ∈ P is designated as the Smarandache fuzzy quotient ideal (S-fuzzy quotient ideal) of the S-semiring S relative to the semifield P, P ⊂ S. From now onwards we will denote by L = { L, ≤, ∧, ∨} the completely distributive lattice which has '0' to be the least element and '1' to be the greatest element.*

*Let X be a non-empty set. A L-fuzzy set in X is a map μ : X→ L and F(X) will denote the set of all L-fuzzy sets in X. if μ, ν ∈ F (X) then μ ⊂ ν if and only if μ (x) ⊂ ν (x) for all x ∈ X and μ ⊂ ν if and only if μ ⊆ ν and ν ≠ μ.*

*It is easily seen* $(F(X), \subseteq, \vee, \wedge) = F(X)$ *is again a complete distributive lattice which has the least and the greatest element say* $\bar{0}$ *and* $\bar{1}$. $\bar{0}(x) = 0$ *and* $\bar{1}(x) = 1$ *for all* $x \in X$.

We recall these concepts in order to define the notion of Smarandache normal L-fuzzy ideals in a S-semiring.

**DEFINITION [50]:** *Let* $\mu \in F(S)$, *S a semiring. Then* $\mu$ *is called on L-fuzzy left (resp. right) ideal of S if for all* $x, y \in S$.

  i. $\mu$ *is an L-fuzzy subsemigroup of* $(S,+)$ *that is* $\mu(x-y) \geq \min(\mu(x), \mu(y))$ *and*
  ii. $\mu(x, y) \geq \mu(y)$ *(resp.* $(\mu(x) > \mu(xy))$*).*

We give the definition of Smarandache L-fuzzy ideal of a S-semiring.

**DEFINITION 5.1.14:** *Let S be S-semiring.* $\mu \in F(S)$. *Clearly F(S) is also a S-semiring.* $\mu$ *is called a Smarandache L-fuzzy left (resp. right) ideal (S-L-fuzzy left (resp. right) ideal) of S if for all* $x, y \in P \subset S$, *P a semifield we have* $\mu$ *is a L-fuzzy subsemigroup of* $(P,+)$ *that is*

  i. $\mu(x + y) \geq \min \{\mu(x), \mu(y)\}$.
  ii. $\mu(xy) \geq \mu(y)$ *(resp.* $\mu(xy) \geq \mu(x)$*).*

*Thus we assume* $\mu \in F(P) \subseteq F(S)$.

**DEFINITION 5.1.15:** *A S-L-fuzzy left (resp. right) ideal* $\mu$ *of S; S a S-semiring is Smarandache normal relative (S-normal relative) to P in S if* $\mu(0) = 1$.

**Note:** The 0 in P will be different for different P's in S.

For instance if we take $S = \{0, a, b, 1\}$ to be the distributive lattice which is a S-semiring.

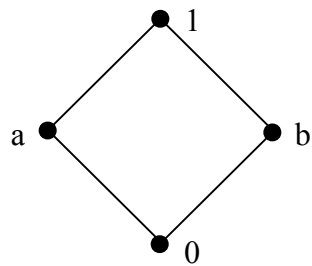

**Figure 5.1.2**

$P_1 = \{a, 1\}$, $P_2 = \{b, 1\}$ and $P_3 = \{0,1\}$ are semifields. Also one can take $P_4 = \{0, a\}$ and $P_5 = \{0, b\}$ to be semifields of S. Thus S has five distinct semifields.

Thus in case of S-L-fuzzy ideals $\mu(0)$ is the largest value of $\mu$ only relative to each $P_i$ for varying i we will have $\mu(0)$ to be varying seen from the above example. Recall, let µ be an S-L-fuzzy left (resp. right) ideal of the S-semiring S and let $\mu^+$ be an L-fuzzy set in S defined by $\mu^-(x) = \mu(x) + 1 - \mu(0)$ for all x ,∈ P ⊂ S then $\mu^+$ will be a Smarandache normal L-fuzzy left (resp. right) ideal (S-normal L-fuzzy left (resp. right) ideal) of S relative to P containing µ  Thus $\mu^+$ will vary from P to P i.e. depends on the subsemifield P in S.

Now we proceed on to develop an analogous results of isomorphism theorems and fuzzy k-ideals given by [52]. He calls a semiring S to be a k- semiring if for any a, b ∈ S there exists a unique element c in S such that either b = a + c or a = b + c but not both.

Now we define the Smarandache analogue of it.

**DEFINITION 5.1.16:** *Let S be a S-semiring. We call S a Smarandache k-semiring (S-k-semiring) if for any a, b ∈ P ⊂ S (P a semifield in S) there exists a unique element c in P such that either b = a + c or a = b + c but not both [52] calls a non-empty subset I of a semiring S to be subsemiring of I if I itself is a semiring with respect to the binary operations defined in S. An ideal I of a semiring S is called a k-ideal if r + a ∈ I implies r ∈ I for each r ∈ S and each a ∈ I.*

*We define Smarandache k-ideal of a S-semiring.*

**DEFINITION 5.1.17:** *Let S be a S-semiring. A S-ideal I of the S-semiring S is called a Smarandache k-ideal (S-k-ideal) if r + a ∈ I implies r ∈ I for each r ∈ P ⊂ S and each a ∈ I.*

The extension k-semiring is defined by [52] as follows.

*Let S be a k-semiring. Let S' be a set of the same continuality with S \ { 0 } such that S ∩ S' = φ and let us denote the image of a ∈ S \ {0} under a given bijection by a'. Let ⊕ and ⊙ denote addition and multiplication respectively on the set $\overline{S} = S \cup S'$ as follows:*

$$a \oplus b = \begin{cases} a+b & \text{if } a,b \in S \\ (x+y)' & \text{if } a = x' \text{ and } b = y' \in S' \\ c & \text{if } a \in S, b = y' \in S' \text{ and } a = y+c \\ c' & \text{if } a \in S, b = y' \in S', a+c = y \end{cases}$$

*where c is the unique element in S such that either a = y + c or a + c = y but not both and*

$$a \odot b = \begin{cases} ab & \text{if } a,b \in S \\ (xy)' & \text{if } a = x', b = y' \in S' \\ (ay)' & \text{if } a \in S, b = y' \in S' \\ (xb)' & \text{if } a = x' \in S' \text{ and } b \in S. \end{cases}$$

*(It is easily verified that these operations are well defined). [52] calls if S is K-semiring then ($\bar{S}, \oplus, \odot$) is a ring called the extension ring of S. For he denotes $\ominus a$ the additive inverse of any element $a \in S$ and write $a \oplus \ominus b$ simply as $a \ominus b$. then it is clear that $a' = \ominus a$ and $a = \ominus a'$ for all $a \in S$.*

Now we define the Smarandache extension ring of a Smarandache k-semiring.

**DEFINITION 5.1.18:** *Let S be a S-k-semiring ($\bar{S}, \oplus, \odot$) be a extension ring of S. If ($\bar{S}, \oplus, \odot$) is a S-ring then we call ($\bar{S}, \oplus, \odot$) a Smarandache extension ring (S-extension ring) of the S-k-semiring S. For a S-k-semiring S, I a S-ideal and I' = $\{a' \in S' \mid a \in I\}$. Then I is a S-k-ideal of S if and only if $\bar{I} = I \cup I'$ is an S-ideal of the S-extension ring, $\bar{R}$. We call $\bar{I}$ the Smarandache extension ideal (S-extension ideal) of I.*

**DEFINITION 5.1.19:** *A mapping f from a S-k-semiring S into a S-k-semiring R is called a S-k-semiring homomorphism if*

$$f(a + b) = f(a) + f(b) \text{ and}$$
$$f(a b) = f(a) f(b)$$

*for all a, b $\in P \subset S$ (P a semifield in S) and f(a), f(b) $\in P_1 \subset R$, $P_1$ a semifield in R.*

*Now for f : S $\to$ R a S-k-semiring homomorphism. Let $\bar{S}$ and $\bar{R}$ be S-extension rings of S and R respectively. Define a map $\bar{f} : \bar{S} \to \bar{R}$ by*

$$\bar{f}(x) = \begin{cases} f(x) & \text{if } x \in S \\ (f(x))' & \text{if } x \in S'. \end{cases}$$

*Then if $\bar{f}$ is a S-ring homomorphism, we call $\bar{f}$ the S-extension ring homomorphism of f.*

Several results in this direction can be obtained by an innovative reader.

**DEFINITION 5.1.20:** *Let A be a S-fuzzy ideal of a S-semiring S. Then A is called a Smarandache k-fuzzy ideal (S-k-fuzzy ideal) of S if $A(x + y) = A(0)$ and $A(y) = A(0)$ imply $A(x) = A(0)$ where A is defined relative to the subsemifield P in S and 0 corresponds to the 0 of P. For A the S-k-fuzzy ideal of S, the set $A_t = \{x \in P \subset S \mid A(x) \geq t\}$ ($t \in [0, 1]$) is called the level subset of S with respect to A.*

Now we proceed on to define the concept of S-fuzzy prime and primary ideal of a S-semiring R.

**DEFINITION 5.1.21:** *A S-fuzzy ideal A of a S-semiring R is called Smarandache fuzzy prime (S-fuzzy prime) if for all a, b $\in P \subset R$ either $A(ab) = A(a)$ or else $A(ab) = A(b)$, A defined relative to $P \subset R$.*

**DEFINITION 5.1.22:** *A S-fuzzy ideal A of a S-semiring R is called Smarandache fuzzy semiprime (S-fuzzy semiprime) if $A(a^m) = A(a)$ for all $a, \in P \subset R$ and for all $m \in Z_+$. (P is the semifield of R relative to which A is defined).*

Next we give the concept of S-fuzzy maximal ideals of the S-semiring.

**DEFINITION 5.1.23:** *Let A be a S-fuzzy ideal of a S-semiring R such that all of its subsets are S-k-ideals of R. A S-fuzzy ideal A of R is called S-fuzzy maximal if*

   i.  $A(0) = 1$,
   ii. $A(e) < A(0)$ and
   iii. *whenever $A(b) < A(0)$ where some $b \in P \subset R$ (P the semifield relative to which A is defined) then $\overline{A}(e_P \oplus (rb)') = A(0)$ for some $r \in P \subset R$ where $e_P$ is the identity of P.*

**DEFINITION 5.1.24:** *A semiring S is said to be semiregular if for each $a \in S$ there are $x_1, x_2 \in S$ such that $a + ax_1a = ax_2a$. We call a S-semiring, S to be Smarandache semiregular (S-semiregular) if for each $a \in P \subset S$, P a proper subset which is a semifield we have $x_1, x_2 \in P$ such that $a + ax_1 a = ax_2 a$.*

**DEFINITION 5.1.25:** *A semiring S is said to have k-closure $\overline{A}$ of a subset $A \subset S$ if $\overline{A} = \{x \in S \mid x + a_1 = a_2, a_1, a_2 \in A\}$. If A is a left ideal of S, then $\overline{A}$ is the smallest left k-ideal of S containing A. We also have $\overline{\overline{A}} = \overline{A}$ for each $A \subset S$ and $A \subset B \subset S$ implies $\overline{A} \subseteq \overline{B}$.*

*For a S-semiring S, $A \subset S$ (A a proper subset of S), the Smarandache k-closure (S-k-closure) $\overline{A}$ of A is defined by $\overline{A} = \{x \in P \subset S \mid P$ a subsemifield of $S, x + a_1 = a_2, a_1, a_2 \in A\}$.*

Several interesting results in this direction can be had and the reader is expected to define and study them.

### 5.2. Smarandache Fuzzy semivector spaces

The notion of semivector spaces over semifields is a recent one. Further the study of Smarandache semivector spaces is a very recent one [134]. Here we just recall the definition of Smarandache semivector spaces and some of its basic properties. The main motivation of this section is to introduce the notion of Smarandache fuzzy semivector spaces.

**DEFINITION [134]:** *Let G be a semigroup under the operation '+'. S any semifield G, be a semivector space over S. G is said to be a Smarandache semivector space (S-semivector space) over S if G is a Smarandache semigroup.*

***Example 5.2.1:*** Let $S = Q \times Z^o$ be a semigroup under componentwise addition. S is a semivector space over $Z^o$, the semifield. Now we see S is a S-semivector space over $Z^o$. It is important to note $S = Q \times Z^o$ is not a semivector space over the semifield $Q^o$.

***Example 5.2.2:*** Let $Q^o \times Q^o \times Q^o = S$ be a semigroup under component-wise addition. Clearly S is a semivector space over $Q^o$ but S is not a S-semivector space as $S = Q^o \times Q^o \times Q^o$ is not a S-semigroup.

**THEOREM 5.2.1**: *All S-semivector spaces over a semifield S are semivector spaces but all semivector spaces need not be S-semivector spaces.*

*Proof:* By the very definition of S-semivector spaces we see all S-semivector spaces are semivector spaces. We note that all semivector spaces need not in general be S-semivector spaces as seen from the above Example 5.2.2.

**DEFINITION 5.2.1**: *Let V be a S-semigroup which is a S-semivector space over a semifield S. A proper subset W of V is said to be Smarandache subsemivector space (S-subsemivector space) of V if W is a S-subsemigroup or W itself is a S-semigroup.*

***Example 5.2.3:*** Let $V = Q^o \times Z^o \times Z$, V is a S-semivector space over $Z^o$. $W = Q^o \times Z^o \times 2Z$ is a S-subsemivector space of V. In fact $W_1 = Q^o \times \{0\} \times Z \subseteq V$ is also a S-subsemivector space of V. But $W_2 = Q^o \times Z^o \times Z^o \subset V$ is not a S-subsemivector space of V over $Z^o$. But $W_2$ is a subsemivector space of V over $Z^o$.

**THEOREM 5.2.2**: *Let V be a S-semivector space over the semifield F. Every S-subsemivector space over S is a subsemivector space over S. But all subsemivector spaces of a S- semivector space need not be S-subsemivector space over S.*

*Proof:* By the very definition of S-subsemivector spaces $W \subset V$ we see W is a subsemivector space of V. But every subsemivector space W of V in general is not a S-subsemivector space as is evidenced from example 5.2.3 the subsemivector space $W_2 = Q^o \times Z^o \times Z^o \subset V$ is not a S-subsemivector space of V. Hence the claim.

***Example 5.2.4:*** Consider $V = Z \times Z^o$, V is a S-semigroup. V is a S-semivector space over $Z^o$. We see the set $\{(-1, 1), (1, 1)\}$ will not span V completely $\{(-1, 0) (1, 0), (0, 1)\}$ will span V. It is left for the reader to find out sets, which can span V completely. Can we find a smaller set, which can span V than the set, $\{(-1, 0), (1, 0), (0, 1)\}$?

Let V be any S-semivector space over the semifield S. Suppose $v_1, \ldots, v_n$ be n set of elements in V then we say

$$\alpha = \sum_{i=1}^{n} \alpha_i v_i$$

in V to be a linear combination of the $v_i$'s. We see when V is just a semivector space given in chapter I we could find semivector spaces using finite lattices but when we have made the definition of S-semivector spaces we see the class of those semivector spaces built using lattices can never be S-semivector spaces as we cannot make even semilattices into S-semigroups as x . x = x for all x in a semilattice. So we are left only with those semivector spaces built using $Q^o$, $Z^o$ and $R^o$ as semifields.

**DEFINITION 5.2.2**: *Let V be a S-semigroup which is a S-semivector space over a semifield S. Let P = {$v_1$, ..., $v_n$} be a finite set which spans V and the $v_i$ in the set P are such that no $v_i$'s in P can be expressed as the linear combination of the other elements in P \ {$v_i$}. In this case we say P is a linearly independent set, which span V.*

**DEFINITION 5.2.3**: *Let V be a S-semigroup which is a S-semivector space over a semifield S. If only one finite set P spans V and no other set can span V and if the elements of that set is linearly independent, that is no one of them can be expressed in terms of others then we say V is a finite dimensional S-semivector space and the cardinality of P is the dimension of V.*

*We see in the case of S-semivector spaces V the number of elements which spans V are always larger than the number of elements which spans the semivector spaces, which are not S-semivector spaces.*

**DEFINITION 5.2.4**: *Let V be a semigroup which is a S-semivector space over a semifield S. A Smarandache basis for V (S-basis for V) is a set of linearly independent elements, which span a S-subsemivector space P of V, that is P, is a S-subsemivector space of V, so P is also a S-semigroup. Thus for any semivector space V we can have several S-basis for V.*

***Example 5.2.5***: Let V = $Z^o \times Z$ be a S-semivector space over Z. Let P = {0} × {pZ} be a S-subsemivector space of V. Now the S-basis for P is {(0, p), (0, -p)}. We see for each prime p we can have S-subsemivector space which have different S-basis.

**DEFINITION 5.2.5**: *Let V and W be any two S-semigroups. We assume P ⊂ V and C ⊂ W are two proper subsets which are groups in V and W respectively. V and W be S-semivector spaces over the same semifield F. A map T: V → W is said to be a Smarandache linear transformation (S-linear transformation) if T(c$p_1$ + $p_2$) = cT$p_1$ + T$p_2$ for all $p_1$, $p_2$ ∈ P and T$p_1$, T$p_2$ ∈ C i.e. T restricted to the subsets which are subgroups acts as linear transformation.*

***Example 5.2.6***: Let V = $Z^o \times Z^o \times Z^o$ be a semigroup under addition. Clearly V is a semivector space over $Z^o$ but V is never a S-semivector space.

In view of this we have got a distinct behaviour of S-semivector space. We know if F is a field V = F × F × ... × F (n times) is a vector space over F. If S is a semifield then W = S × S × ... S = (n times) is a semivector over S. But for a S- semivector space we cannot have this for we see none of the standard semifields defined using $Z^o$, $Q^o$ and $R^o$ are S-semigroups. They are only semigroups under addition and they do not contain a proper subset which is a group under addition.

**THEOREM 5.2.3**: *Let V be a S-semivector space over $Q^o$ or $Z^o$ or $R^o$, then we can always find a subspace in V which is not a S-semivector space.*

*Proof*: If V is to be a S-semivector space the only possibility is that we should take care to see that V is only a semigroup having a subset which is a group i.e. our basic assumption is V is not a group but V is a S-semigroup. Keeping this in view, if V is to be a S-semivector space over $Z^o$ (or $Q^o$ or $R^o$) we can have $V = Z^o \times Z^o \times Z^o \times Q \times R \times \ldots \times Z^o$ i.e. V has at least once $Z^o$ occurring or $Q^o$ occurring or $R^o$ occurring and equally it is a must that in the product V, either Z or Q or R must occur for V to be a S-semigroup. Thus we see if V is a S-semivector space over $Z^o$. Certainly $W = Z^o \times \ldots \times Z^o \subset V$ is a semivector space over $Z^o$ and is not a S-semivector space over $Z^o$. Hence the claim.

**THEOREM 5.2.4**: *Let $V = S_1 \times \ldots \times S_n$ is a S-semivector space over $Z^o$ or $R^o$ or $Q^o$ where $S_i \in \{Z^o, Z, Q^o, Q, R^o, R\}$.*

    i.    *If one of the $S_i$ is Z or $Z^o$ then V can be a S-semivector space only over $Z^o$.*

    ii.    *If none of the $S_i$ is Z or $Z^o$ and one of the $S_i$ is Q or $Q^o$, V can be a S-semivector space only over $Z^o$ or $Q^o$.*

    iii.    *If none of the $S_i$ is Z or $Z^o$ or Q or $Q^o$ only R or $R^o$ then V can be a S-semivector space over $Z^o$ or $Q^o$ or $R^o$.*

*Proof*: It is left for the reader to verify all the three possibilities.

**THEOREM 5.2.5**: *Let $V = S_1 \times \ldots \times S_n$ where $S_i \in \{Z^o, Z, Q^o, Q, R$ or $R^o\}$ be a S-semigroup.*

    i.    *If V is a S-semivector space over $Z^o$ then $W = Z^o \times \ldots \times Z^o$ (n times) is a subsemivector space of V which is not a S-subsemivector space of V.*

    ii.    *If V is a S-semivector space over $Q^o$ then $W = Q^o \times \ldots \times Q^o$ (n times) is a subsemivector space of V which is not a S-subsemivector space of V.*

    iii.    *If V is a S-semivector space over $R^o$ then $W = R^o \times \ldots \times R^o$ (n times) is a subsemivector space of V and is not a S-subsemivector space of V.*

*Proof*: Left for the reader to do the proof as an exercise.

**THEOREM 5.2.6**: *Let $V = S_1 \times \ldots \times S_n$ where $S_i \in \{Z^o, Z, R^o, R, Q^o, Q\}$ if V is a S-semivector space over $Q^o$. Then $W = Z^o \times \ldots \times Z^o$ (n times) $\subset V$ is only a subset of V but never a subspace of V.*

*Proof*: Use the fact V is defined over $Q^o$ and not over $Z^o$.

We define a new concept called Smarandache pseudo subsemivector space.

**DEFINITION 5.2.6**: *Let V be a vector space over S. Let W be a proper subset of V. If W is not a subsemivector space over S but W is a subsemivector space over a proper subset $P \subset S$, then we say W is a Smarandache pseudo semivector space (S- pseudo semivector space) over $P \subset S$.*

*Example 5.2.7*: Let $V = Q^o \times R^o \times Q$ be a S-semivector space over $Q^o$. Now $W = Z^o \times Z^o \times Z^o$ and $W_1 = Q^o \times Q^o \times Q^o$ are S-pseudo semivector spaces over $Z^o \subset Q^o$.

Once the notion of S-semivector spaces are given the concept of Smarandache anti-semivector spaces becomes vital. So we now proceed on to define Smarandache anti-semivector spaces.

**DEFINITION 5.2.7:** *Let V be a vector space over the field F. We say V is a Smarandache anti semivector space (S-anti semivector space) over F if there exists a subspace $W \subset V$ such that W is a semivector space over the semifield $S \subset F$. Here W is just a semigroup under '+' and S is a semifield in F.*

*Example 5.2.8*: Let $V = Q \times R \times Q$ be a vector space over Q. We see $W = Q^o \times R^o \times Q$ is a S-semivector space over $Q^o$. $W_1 = Z \times Z^o \times Z^o$ is not a S-semivector space over $Q^o$. But V is a S-anti semivector space over Q as $P = Z^o \times Z^o \times Z^o$ is a semivector space over $Z^o$.

*Example 5.2.9:* Let $V = Q \times Q \times Q \times Q \times Q$, (5 times) is a vector space over Q. Now $W = Z \times Z^o \times Z^o \times Z^o \times Z^o$ is a S-semivector space over $Z^o$. So V is a S-anti semivector space. The basis for $V = Q \times Q \times Q \times Q \times Q$ is {(1, 0, 0, 0, 0) (0, 1, 0, 0, 0), (0, 0, 1, 0, 0), (0, 0, 0, 0, 1), (0, 0, 0, 1, 0)} as a vector space over Q.

Now what is the basis or the set which spans $W = Z \times Z^o \times Z^o \times Z^o \times Z^o$ over $Z^o$. Clearly the set of 5 elements cannot span W. So we see in case of S-anti semivector spaces the basis of V cannot generate W. If we take $W_1 = Q^o \times Q^o \times Q^o \times Q^o \times Z$ as a S-semivector space over $Z^o$. Clearly $W_1$ cannot be finitely generated by a set. Thus a vector space, which has dimension 5, becomes infinite dimensional when it is a S-anti semivector space.

**DEFINITION 5.2.8**: *Let V and W be any two vector spaces over the field F. Suppose $U \subset V$ and $X \subset W$ be two subsets of V and W respectively which are S-semigroups and so are S-semivector spaces over $S \subset F$ that is V and W are S-anti semivector spaces. A map $T: V \to W$ is said to be a Smarandache T-linear transformation (S-T-linear transformation) of the S-anti semivector spaces if $T: U \to X$ is a S-linear transformation.*

*Example 5.2.10:* Let $V = Q \times Q \times Q$ and $W = R \times R \times R \times R$ be two vector spaces over Q. Clearly, $U = Z \times Z^o \times Z^o \subset V$ and $X = Q \times Z \times Z^o \times Z^o \subset W$ are S-semigroups and U and X are S-semivector spaces so V and W are S-anti semivector spaces.

$T : V \to W$ be defined by $T(x, y, z) = (x, x, z, z)$ for $(x, y, z) \in Z \times Z^o \times Z^o$ and $(x, x, z, z) \in X$ is a S-T-linear operator.

Now we proceed onto define the notion of Smarandache fuzzy semivector spaces and Smarandache fuzzy anti-semivector spaces and obtain some interesting results about them.

**DEFINITION 5.2.9:** *A Smarandache fuzzy semivector space (S-fuzzy semivector space) $(G, \eta)$ is or $\eta_G$ is an ordinary S-semivector space G with a map $\eta : G \to [0, 1]$ satisfying the following conditions:*

   i.   $\eta (a + b) \geq \min (\eta (a), \eta(b))$.
   ii.  $\eta(-a)\ \eta(a)$.
   iii. $\eta(0) = 1$.
   iv.  $\eta (r a ) \geq \eta (a)$

*for all $a, b \in P \subset G$ where P is a subgroup under the operations of G and $r \in S$ where S is the semifield over which G is defined.*

*Thus it is important to note that in case of S-semivector spaces $\eta$ we see that $\eta$ is dependent solely on a subgroup P of the S-semigroup G that for varying P we may see that $\eta : G \to [0, 1]$ may or may not be a S-fuzzy semivector space of V. Thus the very concept of S-fuzzy semivector space is a relative concept.*

**DEFINITION 5.2.10:** *A S-fuzzy semivector space $(G, \eta)$ or $\eta_G$ is an ordinary semivector space G with a map $\eta : G \to [0, 1]$ satisfying the conditions of the Definition 5.2.9. If $\eta : G \to [0, 1]$ is a S-fuzzy semivector space for every subgroup $P_i$ of G then we call $\eta$ the Smarandache strong fuzzy semivector space (S-strong fuzzy semivector space) of G.*

The following theorem is immediate from the two definitions.

**THEOREM 5.2.7:** *Let $\eta : G \to [0, 1]$ be a S-strong fuzzy semivector space of G over the semifield S, then $\eta$ is a S-fuzzy semivector space of G.*

*Proof:* Straightforward by the very definitions hence left as an exercise for the reader to prove.

Now we proceed on to define S-fuzzy subsemivector space.

**DEFINITION 5.2.11**: *Let $(G, \eta)$ be a S-fuzzy semivector space related a subgroup $P \subset G$ over the semifield S. We call $\sigma : H \subset G \to [0, 1]$ a S-fuzzy subsemivector space of $\eta$ relative to $P \subset H \subset G$ where H is a S-subsemigroup G; and $\sigma \subset \eta$ that is $\sigma : G \to [0, 1]$ is a S-fuzzy semivector space relative to the same $P \subset H$ which we denote by $\eta_H$ i.e. $\sigma = \eta_H \subset \eta_G$.*

Now we define Smarandache fuzzy quotient semivector space.

**DEFINITION 5.2.12:** *For an arbitrary S-fuzzy semivector space $\eta_G$ and its S-fuzzy subsemivector space $\eta_H$ the Smarandache fuzzy semivector space $(G/H, \check{\eta})$ or $\eta_{G/H}$ determined by*

$$\overset{\vee}{\eta}(u+H) = \begin{cases} 1 & u \in H \\ \underset{\omega \in H}{\sup} \eta(u+\omega) & \text{otherwise} \end{cases}$$

is called the Smarandache fuzzy quotient semivector space (S-fuzzy quotient semivector space) of $\eta_G$ by $\eta_H$. or equivalently we can say $\overset{\vee}{\eta}$ i.e. the S-fuzzy quotient semivector space of $\eta_G$ by $\eta_H$ is determined by

$$\overset{\vee}{\eta}(v+H) = \begin{cases} 1 & v \in H \\ \underset{\omega \in H}{\inf}(v+\omega) & v \notin H \end{cases}$$

it will be denoted by $\overline{\eta}_{G/H}$. Let $A_{S_1}$ denote the collection of all S-semivector spaces of G; G a S-semigroup, relative to the subgroup $P \subset G$ with respect to the fuzzy subset $S_i$ of the semifield S.

**DEFINITION 5.2.13:** *Let A, $A_1$,..., $A_n$ be fuzzy subsets of G and let K be any fuzzy subset of S*

i. *Define the fuzzy subset $A_1+...+A_n$ of G by the following for all $x \in H \subset G$ (H a subgroup of G) $(A_1+...+A_n)(x) = \sup \{\min \{A_1(x_1),..., A_n(x_n) \mid x = x_1 +...+ x_n, x_i \in H \subset G\}$.*

ii. *Define the fuzzy subset K o A of G by, for all $x \in H \subset G$ (K o A)(x) = $\sup\{\min \{K(c), A(y)\} \mid c \in S, y \in H \subset G, x = cy\}$.*

*Fuzzy singletons are defined as in case of any other fuzzy subset.*

Further almost all results related to S-fuzzy vector spaces can be developed in case of S-fuzzy semivector spaces.

**DEFINITION 5.2.14:** *Let $\{A_i \mid i \in I\}$ be a non-empty collection of fuzzy subsets of S. Then the fuzzy subset $\underset{i \in I}{\bigcap} A_i$ of G is defined by the following for all $x \in H \subset G$ (H a subgroup of G)*

$$\left(\underset{i \in I}{\bigcap} A_i\right)(x) = \inf \{A_i(x) \mid i \in I\}.$$

*Let $A \in \mathcal{A}_{S_1}$ where $S_1$ is a fuzzy subset of the semifield S. Let X be a fuzzy subset of G such that $X \subset A$. (relative to some fixed subgroup, H of G) Let $\langle X \rangle$ denote the intersection of all fuzzy subsemispaces of G (over $S_1$) that contain X and are contained in A. Then $\langle X \rangle$ is called the Smarandache fuzzy subsemispaces (S-fuzzy subsemispaces) of A fuzzily spanned by X. We define the notion of fuzzy freeness in case of Smarandache fuzzy semivector spaces.*

Let $\xi$ denote a set of fuzzy singletons of H in G such that $x_\lambda$, $x_\nu \in \xi$ then $\lambda = \nu > 0$. Define the fuzzy subset $X(\xi)$ of H in G by the following for all $x \in H \subset G$, $X(\xi)(x) = \lambda$ if $x_\lambda \in \xi$ and $X(\xi)(x) = 0$ otherwise. Define $\langle \xi \rangle = \langle X(\xi) \rangle$. Let X be a fuzzy subset of H in G. Define $\xi(X) = \{x_\lambda / x \in H \subset G, \lambda = X(x) > 0\}$. Then $X(\xi(X)) = X$ and $\xi(X(\xi)) = \xi$. If there are only a finite number of $x_\lambda \in \xi$ with $\lambda > 0$, we call $\xi$ finite. If $X(x) > 0$ for only a finite number of $x \in X$, we call X finite. Clearly $\xi$ is finite if and only if $X(\xi)$ is finite and X is finite if and only if $\xi(X)$ is finite. For $x \in H \subset G$, let $X\backslash\{x\}$ denote the fuzzy subset of H in G defined by the following, for all $y \in H \subset G$, $(X \backslash x)(y) = X(y)$ if $y \neq x$ and $(X \backslash x)(y) = 0$ if $y = x$. Let $S_1$ be a fuzzy subset of the semifield S. Let $A \in \mathcal{A}_{S_1}$ and let X be a fuzzy subset of $H \subset G$ (H a subgroup of the S-semigroup G) such that $X \subset A$. Then X is called a Smarandache fuzzy system of generator (S-fuzzy system of generator) of A over $S_1$ if $\langle X \rangle = A$.

X is said to be Smarandache fuzzy free (S-fuzzy free) over $S_1$ if for all $x_\lambda \subseteq X$, where $\lambda = X(x)$, $x_\lambda \not\subset \langle X \backslash x \rangle$. X is said to be a Smarandache fuzzy basis (S-fuzzy basis) for A if X is a S-fuzzy system of generators of A and X is S-fuzzy free. Let $\xi$ denote a set of fuzzy singletons of $H \subset G$ such that $x_\lambda\, x_k \in \xi$ then $\lambda = k$, and $x_\lambda \subseteq A$, then $\xi$ is called a Smarandache fuzzy singletons system of generators (S-fuzzy singletons system of generators) of A over $S_1$ if, $\langle \xi \rangle = A$. $\xi$ is said to be S-fuzzy free over $S_1$ if for all $x_\lambda \in \xi$, $x_\lambda \not\subseteq \langle \xi \backslash \{x_\lambda\} \rangle$, $\xi$ is said to be fuzzy free over $S_1$ if for all $x_\lambda \in \xi$, $x_\lambda \not\subseteq \langle \xi \backslash \{x_\lambda\} \rangle$, $\xi$ is said to be a fuzzy basis of singletons for A if $\xi$ is a S-fuzzy singleton system of generators of A and $\xi$ is S-fuzzy free.

For $\lambda = \langle \xi \rangle (0)$, $0_\lambda \subseteq \langle \xi \rangle$ for every set $\xi$ of fuzzy singletons of H in G. Also $x_0 \subseteq \langle \xi \rangle$ for every such $\xi$ where $x \in H \subset G$. Thus if $\xi$ is a set of fuzzy singletons of $H \subset G$ such that either $x_0$ or $0_\lambda \in \xi$ then $\xi$ is not S-fuzzy free over $S_1$.

Let $A \in \mathcal{A}_{S_1}$. Set $A^* = \{x \in H \subset G\, /A(x) > 0\}$ and $S_1^* = \{c \in S / S_1(c) > 0\}$

It is easy to prove the following theorem hence left for the reader as an exercise.

**THEOREM 5.2.8:** *Suppose $A \in \mathcal{A}_{S_1}$. Then*

  i. $S_1^*$ is a subsemifield of S.
  ii. $A^*$ is a S subsemispace of the S-semivector space $H \subset G$ over $S_1^*$.

Now we proceed on to define the notion of Smarandache fuzzy linearly independent set over a fuzzy subset $S_1$ of a semifield S.

**DEFINITION 5.2.15:** *Let $A \in \mathcal{A}_{S_1}$, and let $\xi \subseteq \{x_\lambda / x \in A^*, \lambda \leq A(x)\}$ be such that if $x_\lambda$, $x_k \in \xi$ then $\lambda = k$. Then $\xi$ is said to be a Smarandache fuzzy linearly independent (S-fuzzy linearly independent) over $S_1$ if and only if for every finite subset $\{x_{1_{\lambda_1}}, x_{2_{\lambda_2}}, \cdots, x_{n_{\lambda_n}}\}$ of $\xi$, whenever*

$$\left(\sum_{i=1}^{n} c_{i\mu_i}\, o\, x_{i_{\lambda_i}}\right)(x) = 0$$

*for all $x \in H\setminus\{0\} \subset G$ (0 is the additive identity of the subgroup H contained in the S-semigroup G) where $c_i \in S$, $0 < \mu_1 \leq S_1(c_i)$ for $i = 1, 2, ..., n$ then $c_1 = c_2 = ... = c_n = 0$. Several analogous results can be obtained.*

Following the definitions of [79] we give some definitions of S-fuzzy semivector spaces. It is left for the reader to obtain a necessary and sufficient condition for these concepts to be equivalent or counter examples to show the non-equivalence of these definitions.

From now onwards S will be a semifield and G a S-semigroup and G a S-semivector space over S.

**DEFINITION 5.2.16:** *A fuzzy subset $\mu$ of the S semivector space G is a Smarandache subsemispace (S-subsemispace) of G if for any $a, b \in S$ and $x, y \in H \subset G$ (H a subgroup relative to which $\mu$ is defined ) the following conditions holds good. $\mu$ (ax + by) $\geq \mu(x) \wedge \mu(y)$. If $\mu$ is a S-fuzzy subsemispace of the S-semivector space G and $\alpha \in [0, 1]$ then define $G_H = \{x \in H \subset G \mid \mu(x) \geq \alpha\}$.*

*The subspaces $G_H^\alpha, \alpha \in \text{Im } \mu$ are S-level subspaces of $\mu$ relative to $H \subset G$. A set of vectors B is S-fuzzy linearly independent if*

  i.   *B is S-linear independent.*
  ii.  $\mu\left(\sum_{i=1}^{n} a_i x_i\right) = \bigwedge_{i=1}^{n} \mu(a_i x_i)$ *for finitely many distinct element $x_1, ..., x_n$ of B.*

*A S-fuzzy basis for the S-fuzzy subsemispace $\mu$ is a S-fuzzy linearly independent basis for $H \subset G$.*

Now we define Smarandache linear maps of S-semivector spaces.

**DEFINITION 5.2.17:** *Let G and L be S-semivector spaces over the same semifield S and let $\mu : G \to [0, 1]$ and $\lambda : L \to [0, 1]$ be S-fuzzy subsemispaces.*

*The S-linear map relative to subgroup H in G, $\phi : G \to L$, from the fuzzy subsemispaces $\mu$ to S-fuzzy subsemispaces $\lambda$ if $\lambda (\phi(x)) \geq \mu (x)$ for all $x \in H \subset G$. The space of S-fuzzy linear maps from $\mu$ to $\lambda$ is denoted by S F Hom ($\mu, \lambda$).*

*Suppose SF Hom ($\mu, \lambda$) be a collection of S-fuzzy linear maps from $\mu$ to $\lambda$ defined relative to the subgroup H in G. We define a function $\nu$ from SF (Hom ($\mu, \lambda$)) into unit interval [0, 1] where $\mu : G \to [0, 1]$ and $\lambda : G \to [0, 1]$ are S-fuzzy subsemispaces (relative to the subgroup H in G) G and L respectively such that $\nu$ determines the S-fuzzy subsemispace of SF Hom ($\mu, \lambda$).*

*The fuzzy subset $\nu$ of SF Hom ($\mu, \lambda$) is defined by*

$$\nu(\phi) = \begin{cases} \inf\{\lambda(\phi(x)) - \mu(x)\} \mid x \in H \subset G, \ \phi(x) \neq 0 \text{ if } \phi \neq 0 \\ \sup\{\lambda\phi(x)) - \mu(x) \mid x \in H \subset G \text{ if } \phi = 0. \end{cases}$$

*Thus if $\phi \neq 0$ then $\nu(\phi)$ is the greatest real number $\alpha$ such that for every vector $x$ in $H \subset G$ we have either $\lambda(\phi(x)) - \mu(x) \geq \alpha$ or $\phi(x) = 0$.*

*We assume from now onwards $G$ and $L$ are finite dimensional S-semivector spaces defined over the same semifield S.*

$$\mu : G \to [0, 1]$$
$$\lambda : G \to [0, 1]$$

*are S-fuzzy subsemispaces with S-fuzzy bases $\{e_1, ..., e_n\}$ and $\{f_1, ..., f_n\}$ respectively.*

*Then the Smarandache dual functionals (S-dual functionals) $\{e^1, e^2, ..., e^n\}$ and $\{f^1, f^2, ..., f^m\}$ form S-fuzzy basis of S-fuzzy subsemivector spaces.*

$$\mu^* : G^* \to [0, 1] \text{ and}$$
$$\lambda^* : G^* \to [0, 1].$$

*If $\phi \in S \text{ Hom }(G, L)$ then the dual map $\phi^* \in S \text{ Hom }(G^*, L^*)$ defined by $\phi'(g)(x) = g(\phi(x))$ for every $g \in P^* \subset L^*$ and $x \in H \subset G$ where $P^*$ is the related subgroup of $H^*$ in $G^*$ in the S-semigroup $L^*$. It is left for the reader to find out whether*

$$\phi'_{ij}(f^t)(e_s) = f^t(\phi_{ij}(e_s))$$
$$= f^t(\delta_{js} f_i)$$
$$= \delta_{it} \delta_{js}.$$

*and*

$$\phi'_{ij}(f^t) = \delta_{it} e^j.$$

*Now we will just test whether we can define Smarandache fuzzy continuity of linear maps of S-semivector spaces.*

*We consider $X$ to be a S-semivector space over the semifiled $Q^0$ or $R^0$. Fuzzy subsets of $X$ are denoted by greek letters; in general $\chi_A$ denotes the characteristic function of a set A.*

*By a fuzzy point (singleton) $\mu$ we mean a fuzzy subset $\mu : X \to [0, 1]$ such that*

$$[z] = \begin{cases} \alpha & \text{if } z = x \\ 0 & \text{otherwise} \end{cases}$$

where $\alpha \in (0,1)$, $I^X = \{\mu \mid \mu : X \to I = [0,1]\}$ Here $I$ denotes the closed interval $[0, 1]$. For any $\mu, \nu \in I^X$ $\mu + \nu \in I^X$ is defined by $(\mu + \nu)(x) = \sup \{\mu(\nu) \wedge \nu(\upsilon) / \mu + \nu = x, H \subset X, H$ a subgroup in $X\}$.

If $\mu \in I^*$, $t \in Q^0$, $R^0$, $t \neq 0$, then $(t\mu)(x) = \mu(X_H / t) = \mu(H/t)$, ($H \subset X$ is an additive subgroup of $X$ relative to which $\mu$ is defined)

$$(0 \cdot \mu)(x) = \begin{cases} 0 & \text{if } x \neq 0 \\ \bigvee_{y \in H \subset X} \mu(y) & \text{if } x = 0. \end{cases}$$

For any two sets $X$ and $Y$, $f : X \to Y$ then the fuzzification of $f$ denoted by $f$ itself is defined by

i. $f(\mu)(y) = \begin{cases} \bigvee_{x \in f^{-1}(y)} \mu(x) & \text{if } f^{-1}(y) \neq 0 \\ 0 & \text{otherwise for all } y \in Y \text{ and for all } \mu \in I^X \end{cases}$

ii. $f^{-1}(\mu)(x) = \mu(f(x))$ for all $x \in X$, for all $\mu \in I^X$. $\mu \in I^X$ is said to be a Smarandache fuzzy subsemispace (S-fuzzy subsemispace) if

    i. $\mu + \mu \leq \mu$ and
    ii. $t\mu < \mu$ for all $t \in Q^0$ or $R^0$ ($\mu : H \subset X \to [0, 1]$ is a S-fuzzy subsemivector space).

S-convex if $t\mu + \overline{(1 - t\mu)} \leq \mu$ for each $t \in [0, 1]$
S-balanced if $t\mu \leq \mu$ for $t \in Q^o$ or $R^o$, $|t| \leq 1$
S-absorbing if $\bigvee_{t>0} t\mu(x) = 1$ for all $x \in H \subset X$ ($H$ a subgroup of $X$).

Recall $(X, \tau)$ be a topological space and $\omega(\tau) = \{f \mid f : (X, \tau) \to [0, 1]$ is lower semicontinuous$\}$, $\omega(\tau)$ is a fuzzy topology on $X$.

This will called as fuzzy topology induced by $\tau$ on $X$. Let $(X, \delta)$ be a fuzzy topology $i(\delta)$ the ordinary topology on $X$. A base for $i(\delta)$ is provided by the finite intersection $\bigcap_{i=1}^{n} \overline{\nu}_i^1 (\in_i, 1]$, $\nu_i \in \delta$, $\in_i \in I$. A map $f : (X, \tau) \to (Y, \tau')$ is continuous if and only if for every $\mu \in \tau'$ in $f^{-1}(\mu) \in \tau$ in $X$ where $(X, \tau)$ and $(Y, \tau)$ are fuzzy topological spaces.

**DEFINITION 5.2.18:** *A fuzzy topology $\tau_x$ on a S-semivector space $X$ is said to be S-fuzzy semivector topology if addition and scalar multiplication are continuous from $H \times H$ and $Q^o \times H$ respectively to $H \subset X$ ($H$ an additive subgroup of the S-semigroup $X$) with the fuzzy topology induced by the usual topology on $Q^o$. A S-semivector space together with a S-fuzzy semivector topology is called a S-fuzzy topological semivector space.*

*A fuzzy seminorm on X is an absolutely convex absorbing fuzzy subset of X. A fuzzy seminorm on X is a fuzzy norm if $\bigwedge_{t>0} tp(x) = 0$ for $x \neq 0$.*

*If $\rho$ is a fuzzy seminorm on X we define $P_\epsilon : X \to R_+$ by $P_\epsilon(x) = \inf \{t > 0 \mid tp(x) > \epsilon\}$. Clearly $P_\epsilon$ is a seminorm on X for each $\epsilon \in (0,1)$.*

**DEFINITION 5.2.19:** *A Smarandache fuzzy seminorm (S-fuzzy seminorm) on a S-seminorm on a S-semivector space X is an S-absolutely, S-convex absorbing fuzzy subset of X.*

Obtain some interesting results in this direction using S-semivector spaces in the place of vector spaces.

We just define the notion of Smarandache fuzzy anti-semivector spaces.

**DEFINITION 5.2.20:** *Let V be a vector space over the field F. We say a fuzzy subset $\mu$ on V is a Smarandache fuzzy anti-semivector space (S- fuzzy anti-semivector space) over F if there exists a subspace W of V such that W is a semivector space over the semifield S contained in F. Here $\mu_W : W \to [0, 1]$ is a S-fuzzy semivector space then we call $\mu$ the S-fuzzy anti-semivector space.*

All results can be derived in an analogous way for S-fuzzy anti-semivector spaces.

## 5.3 Smarandache Fuzzy non-associative semirings

In this section for the first time we introduce both the notion of non-associative semirings, Smarandache non-associative semirings and above all the notion of Smarandache fuzzy non-associative semirings and derive several interesting properties about them.

**DEFINITION 5.3.1:** *Let $(S, +, \bullet)$ be a semiring. We call S a non-associative semiring if*

  i.   *$(S, +)$ is a commutative monoid and*
  ii.  *$(S, \bullet)$ is a groupoid i.e. a non-associative semigroup and*
  iii. *$a \bullet (b + c) = a \bullet b + a \bullet c$ and $(a + b) \bullet c = a \bullet c + b \bullet c$ for all $a, b, c \in S$.*

**Example 5.3.1:** Let $Z^0 = \{$set of all positive integers with zero$\}$, $Z^0$ is a semiring. Let $(L, \bullet)$ be a loop given by the following table:

| * | 1 | $g_1$ | $g_2$ | $g_3$ | $g_4$ | $g_5$ |
|---|---|---|---|---|---|---|
| 1 | 1 | $g_1$ | $g_2$ | $g_3$ | $g_4$ | $g_5$ |
| $g_1$ | $g_1$ | 1 | $g_3$ | $g_5$ | $g_2$ | $g_4$ |
| $g_2$ | $g_2$ | $g_5$ | 1 | $g_4$ | $g_1$ | $g_3$ |
| $g_3$ | $g_3$ | $g_4$ | $g_1$ | 1 | $g_5$ | $g_2$ |
| $g_4$ | $g_4$ | $g_3$ | $g_5$ | $g_2$ | 1 | $g_1$ |
| $g_5$ | $g_5$ | $g_2$ | $g_4$ | $g_1$ | $g_3$ | 1 |

Clearly the loop semiring $Z^0L$ is a non-associative semiring. A non-associative semiring can also be constructed using groupoids; i.e. groupoid semirings are non-associative semiring.

Now we proceed on to define Smarandache non-associative semirings.

**DEFINITION 5.3.2:** *Let (S, +, •) be a non-associative semiring. S is said to be a Smarandache non-associative semiring (S-non-associative semiring) if S contains a proper subset P such that (P, +, •) is an associative semiring.*

The class of S-na-semirings is non-empty. For consider class of loop semirings SL for varying loops and varying semirings. SL is a S-na-semiring.

Now as our main motivation is to define Smarandache fuzzy-na-semirings, we having defined S-na-semiring just give the definitions of Smarandache non-associative subsemirings, Smarandache non-associative ideals in semirings and finally the notion of Smarandache non-associative k-semirings.

**DEFINITION 5.3.3:** *Let (S, +, •) be a non-associative semiring. A subsemiring (P, +, •) is said to be a Smarandache non-associative subsemiring (S-non-associative subsemiring) if (P, +, •) itself is a S-na semiring.*

We have the following nice result which is of course direct.

**THEOREM 5.3.1:** *Let (S, +, •) be a na semiring if (S, +, •) has a S-na-subsemiring then (S, +, •) itself is a S-na-semiring.*

**DEFINITION 5.3.4:** *Let (S, +, •) be a S-na semiring. An non-empty subset I of S is said to be a Smarandache na-ideal (S-na-ideal) if*

  i. *I is a S-subsemiring.*
  ii. *$pI \subset I$ and $Ip \subset I$ for all $p \in P \subset S$ where P is an associative subsemiring of S relative to which I is defined.*

*Thus it has become noteworthy to mention that I may not be a S-ideal of S relative to all associative subsemirings of S.*

*If in case I happens to be a S-ideal relative to every associative subsemiring of S then we call I a Smarandache strong non-associative ideal (S-strong non-associative ideal) of S.*

**THEOREM 5.3.2:** *Every S-strong na-ideal of S is a S-na ideal of S.*

*Proof*: Straightforward by the very definition hence left as an exercise for the reader to prove.

Now we proceed on to define the notion of Smarandache semiregular na-semirings.

**DEFINITION 5.3.5:** *Let (S, +, •) be a non-associative semiring which is a S- na semiring. Let $P \subset S$ be such that (P, +, •) be an associative semiring.*

*If for each $a \in P$ there exists $x_1, x_2 \in S$ such that $a + ax_1a = ax_2a$ we call S a Smarandache non-associative semiregular semiring (S-non-associative semiregular semiring).*

*Let S be a semiring and $A \subseteq S$. The k-closure $\overline{A}$ of A is defined by $\overline{A} = \{ x \in S / x + a_1 = a_2, a_1, a_2 \in A\}$ We say, $A \subseteq P \subseteq S$ where P is an associative semiring of S to be Smarandache k-closure (S-k-closure) of A if $\overline{A} = \{ x \in P / x + a_1 = a_2, a_1, a_2 \in A\}$. If A is a left ideal of S then $\overline{A}$ is the smallest left k-ideal of S containing A. If A is a S-na-left ideal of S then $\overline{A}$ which is a S-k-closure of a is smallest S-left k-ideal of S containing A.*

**DEFINITION 5.3.6:** *Let S be a S-na-semiring. S is called a Smarandache na-k-semiring (S-na-k-semiring) if for any $a, b \in P \subset S$ (P an associative semiring in S) there exists a unique c in P such that either $b = a + c$ or $a = b + c$ but not both. Thus we see in case of S-na-k-semiring unlike semirings we can have several S-na-k-semirings depending on the number of associative semirings in the S-na-semiring which happen to be k- semirings.*

Now we proceed on to define Smarandache fuzzy non-associative ideals (S-fuzzy na-ideals) and their properties.

**DEFINITION 5.3.7:** *Let $(S, +, \bullet)$ be a non-associative semiring. A fuzzy subset µ of the semiring S is said to be a Smarandache fuzzy left (resp. right) ideal of S if $\mu(x + y) \geq \min \{\mu(x), \mu(y)\}$ and $\mu(xy) > \mu(y)$ (resp. $\mu(xy) \geq \mu(x)$) for all $x, y \in P \subset S$ where P is an associative semiring relative to which µ is defined. µ is a Smarandache fuzzy ideal of S if it is both a S-fuzzy left and a S-fuzzy right ideal of S.*

**DEFINITION 5.3.8:** *A S-fuzzy ideal µ of a S-na-semiring S is said to be a Smarandache fuzzy k-ideal (S-fuzzy k-ideal) of S if $\mu(x) \geq \min \{\max \{\mu(x + y), \mu(y + x)\}, \mu(y)\}$ for all $x, y \in P \subset S$.*

All other notions and results related to S-na semirings can be developed in case of S-semirings. So we leave rest of the result to be proved by the reader.

## 5.4 Smarandache fuzzy bisemirings and its properties

In this section we introduce the concept of Smarandache fuzzy bisemirings. The study of Smarandache bisemirings and the notion of bisemirings is very recent. Here we recall some basic properties about Smarandache bisemirings and then proceed on to define Smarandache fuzzy bisemirings.

The notion of bisemirings is introduced in Chapter I. Here we just recall the definition of Smarandache bisemiring. For more about S-bisemirings please refer [135].

**DEFINITION 5.4.1:** *Let $(S, +, \bullet)$ be a bisemiring. We call $(S, +, \bullet)$ a Smarandache bisemiring (S-bisemiring) if S has a proper subset P such that P under the operations of S is a bisemifield.*

**DEFINITION 5.4.2:** *Let S be a bisemiring. A non-empty proper subset A of S is said to be a Smarandache sub-bisemiring (S-sub-bisemiring) if A is a S-bisemiring i.e. A has a proper subset P such that P is a bisemifield under the operations of S.*

**THEOREM 5.4.1:** *If (S, +, •) is a bisemiring having a S-sub-bisemiring then S is a S-bisemiring.*

*Proof:* Follows from the very definitions.

**DEFINITION 5.4.3:** *Let S be a S-bisemiring. We say S is a Smarandache commutative bisemiring (S-commutative bisemiring) if S has a S-sub-bisemiring, which is commutative. If S has no commutative S-sub-bisemirings then we call S to be a Smarandache non-commutative bisemiring (S-non-commutative bisemiring). If every S-sub-bisemiring of S is commutative then we call S a Smarandache strongly commutative bisemiring (S-strongly commutative bisemiring).*

**DEFINITION 5.4.4:** *Let (S, +, •) be a bisemiring. A non-empty subset P of S is said to be a Smarandache right (left) bi-ideal (S-right (left) bi-ideal) of S if the following conditions are satisfied.*

  i.  *P is a S-sub-bisemiring.*
  ii. *For every p ∈ P and A ⊂ P where A is a bisemifield of P we have for all a ∈ A and p ∈ P, a • p (p • a) is in A. If P is simultaneously both a S-right bi-ideal and S-left bi-ideal then we say P is a Smarandache bi-ideal (S-bi-ideal) of S.*

**DEFINITION 5.4.5:** *Let (S, +, •) be a bisemiring. A proper subset A of S is said to be Smarandache pseudo sub-bisemiring (S-pseudo sub-bisemiring) if the following conditions are true:*

*If there exists a subset of P of S such that A ⊂ P, where P is a S-sub-bisemiring, i.e. P has a subset B such that B is a bisemifield under the operations of S or P itself is a bisemifield under the operations of S.*

**DEFINITION 5.4.6:** *Let S be a bisemiring. A non-empty subset P of S is said to be a Smarandache pseudo right (left) bi-ideal (S-pseudo right (left) bi-ideal) of the bisemiring S if the following conditions are true:*

  i.  *P is a S-pseudo sub-bisemiring i.e. P ⊂ A, A is a bisemifield in S.*
  ii. *For every p ∈ P and every a ∈ A, a • p ∈ P (p • a ∈ P). If P is simultaneously both a S-pseudo right bi-ideal and S-pseudo left bi-ideal then we say P is a Smarandache pseudo bi-ideal (S-pseudo bi-ideal).*

**DEFINITION 5.4.7:** *Let S be a bisemiring. A nonempty subset P of S is said to be a Smarandache dual bi-ideal (S-dual bi-ideal) of S if the following conditions hold good:*

  i.  *P is a S-sub-bisemiring.*
  ii. *For every p ∈ P and a ∈ A \ {0}; a + p is in A, where A ⊂ P.*

**DEFINITION 5.4.8:** *Let S be a bisemiring. A nonempty subset P of S is said to be a Smarandache pseudo dual bi-ideal (S-pseudo dual bi-ideal) of S if the following conditions are true:*

   i.    *P is a S-pseudo sub-bisemiring i.e. P ⊂ A, A is a bisemifield in S or A contains a bisemifield.*
   ii.   *For every p ∈ P and a ∈ A, p + a ∈ P. Clearly P is simultaneously left and right S-pseudo dual bi-ideal of S as S is additively commutative.*

**DEFINITION 5.4.9:** *Let S be a S-bisemiring. S is said to be a Smarandache bisemidivision ring (S-bisemidivision ring) if the proper subset A ⊂ S is such that*

   i.    *A is a S-sub-bisemiring.*
   ii.   *A contains a subset P such that P is a bisemidivision ring, that is P has no zero divisors and P is a non-commutative bisemiring.*

**DEFINITION 5.4.10:** *Let S be a bisemiring we say S is a Smarandache right chain bisemiring (S-right chain bisemiring) if the S-right bi-ideals of S are totally ordered by inclusion.*

*Similarly we define Smarandache left chain bisemirings (S-left chain bisemirings).*

*If all the S-bi-ideals of the bisemiring are ordered by inclusion we say S is a Smarandache chain bisemiring (S-chain bisemiring).*

**DEFINITION 5.4.11:** *Let S be a bisemiring. If $S_1 \subset S_2 \subset \ldots$ is a monotonic ascending chain of S-bi-ideals $S_i$ and there exists a positive integer r such that $S_r = S_n$ for all $r \geq n$, then we say the bisemiring S satisfies the Smarandache ascending chain conditions (S-acc) for S-bi-ideals in the bisemiring S.*

*We say S satisfies Smarandache descending chain conditions (S-dcc) on S-bi-ideals $S_i$, if every strictly decreasing sequence of S-ideals $N_1 \supset N_2 \supset \ldots$ in S is of finite length. The Smarandache min conditions (S-min conditions) for S-bi-ideals holds in S if given any set P of S-bi-ideals in S, there is a bi-ideal of P that does not properly contain any other bi-ideal in the set P.*

*Similarly one can define Smarandache maximum condition (S-maximum condition) for S-bi-ideals in case of bisemirings.*

**DEFINITION 5.4.12:** *Let S be a bisemiring. S is said to be a Smarandache idempotent bisemiring (S-idempotent bisemiring) if a proper subset P of S that is a sub-bisemiring of S satisfies the following conditions:*

   i.    *P is a S-sub-bisemiring.*
   ii.   *P is an idempotent bisemiring.*

**DEFINITION 5.4.13:** *Let S be a bisemiring. S is said to be a Smarandache e-bisemiring (S-e-bisemiring) if S contains a proper subset A satisfying the following conditions:*

i. *A is a sub-bisemiring.*
   ii. *A is a S-sub-bisemiring.*
   iii. *A is a e-bisemiring.*

**DEFINITION 5.4.14:** *A bisemiring S is said to be Smarandache bisemiring of level II (S-bisemiring of level II) if S contains a proper subset P that is a bifield.*

**DEFINITION 5.4.15:** *Let R be a biring. R is said to be a Smarandache anti-bisemiring (S-anti-bisemiring) if R contains a subset S such that S is just a bisemiring.*

**DEFINITION 5.4.16:** *Let S be a bisemifield. S = $S_1 \cup S_2$ is said to be a Smarandache bisemifield (S-bisemifield) if a proper subset P = $P_1 \cup P_2$ of S is a S-bisemialgebra with respect to the same induced operations and an external operator (i.e. $P_1$ is a S-semialgebra and $P_2$ is a k-semialgebra, P = $P_1 \cup P_2$ is a S-bisemialgebra).*

**DEFINITION 5.4.17:** *Let S be a bisemifield. S is said to be a Smarandache bisemifield of level II (S-bisemifield of level II) if S has a proper subset P where P is a bifield.*

**DEFINITION 5.4.18:** *Let S be a bisemifield. A proper subset P of S is said to be a Smarandache sub-bisemifield I (II) (S-sub-bisemifield I (II)) if P is a S-bisemifield of level I (or level II).*

**DEFINITION 5.4.19:** *Let S be a bifield or a biring. S is said to be a Smarandache anti-bisemifield (S-anti-bisemifield) if S has a proper subset, A which is a bisemifield.*

**DEFINITION 5.4.20:** *Let S be a biring or a bifield. A proper subset P in S is said to be a Smarandache anti sub-bisemifield (S-anti sub-bisemifield) of S if P is itself a S-anti bisemifield.*

**DEFINITION 5.4.21:** *Let S be a bifield or a biring, which is a S-anti bisemifield. If we can find a subset P in the sub-bisemifield T in S such that*

   i. *P is a bisemiring.*
   ii. *for all p ∈ P and t ∈ T, pt ∈ P, then P is called the Smarandache anti bi-ideal (S-anti bi-ideal) of the S-anti-bisemifield.*

Now we just define na bisemiring before we proceed onto define Smarandache na bisemirings.

**DEFINITION 5.4.22:** *Let (S, +, •) be a bisemiring. We say S is a non-associative bisemiring if S = $S_1 \cup S_2$ where $S_1$ and $S_2$ are semirings where at least one of $S_1$ or $S_2$ is a non-associative semiring.*

**DEFINITION 5.4.23:** *Let (S, +, •) be an associative bisemiring. Let L be a loop, the loop bisemiring SL = $S_1L \cup S_2L$ (where S = $S_1 \cup S_2$; $S_1$ and $S_2$ are commutative semirings with unit) and $S_1L$ and $S_2L$ are the loop semirings of the loop L over the semirings $S_1$ and $S_2$ respectively.*

*Clearly the loop bisemiring is a non-associative bisemiring.*

**DEFINITION 5.4.24:** *Let $(S, +, \bullet)$ be a commutative bisemiring with 1 and G any groupoid. The groupoid bisemiring $SG = S_1G \cup S_2G$ is a non-associative bisemiring. ($S_1G$ and $S_2G$ are the groupoid semirings of the groupoid G over the semirings $S_1$ and $S_2$ respectively).*

**DEFINITION 5.4.25:** *Let $(S, +, \bullet)$ be a non-associative bisemiring. We say S is a strict bisemiring if $a + b = 0$ then $a = 0$ and $b = 0$. We call $(S, +, \bullet)$ a na-bisemiring, to be a commutative bisemiring if both $S_1$ and $S_2$ are commutative semirings. $(S, +, \bullet)$ is said to be a semiring with unit if there exists $1 \in S$ such that $1 \bullet s = s \bullet 1 = s$ for all $s \in S$.*

*An element $0 \neq x \in S$ is said to be a zero divisor if there exists $0 \neq y \in S$ such that $x \bullet y = 0$. A na-bisemiring which has no zero-divisors but which is commutative with unit is called a non-associative bisemifield. If the operation in S is non-commutative we call S a na-bidivision semiring.*

**DEFINITION 5.4.26:** *Let $(S, +, \bullet)$ be a na-bisemiring. We say S is a Moufang bisemiring if all elements of S satisfy the Moufang identity i.e. $(xy)(zx) = (x(yz))x$ for all $x, y, z \in S$. A na-bisemiring S is said to be a Bruck bisemiring if $(x(yz))z = x(y(xz))$ and $(xy)^{-1} = x^{-1} y^{-1}$ for all $x, y, z \in S$. A na-bisemiring S is called a Bol bisemiring if $((xy)z)y = x((yz)y)$ for all $x, y, z \in S$. We call a na-bisemiring S to be right alternative if $(xy)y = x(yy)$ for all $x, y \in S$; left alternative if $(xx)y = x(xy)$ and alternative if it is both right and left alternative.*

**DEFINITION 5.4.27:** *Let $(S, +, \bullet)$ be a na-bisemiring. $x \in S$ is said to be right quasi regular (r.q.r) if there exists a $y \in R$ such that $x \circ y = 0$ and x is said to be left quasi regular (l.q.r) if there exists a $y' \in R$ such that $y' \circ x = 0$.*

*An element is quasi regular (q.r) if it is both right and left quasi-regular.*

*y is known as the right quasi inverse (r.q.i) of x and y' as the left quasi inverse (l.q.i) of x. A right ideal or left ideal in R is said to be right quasi regular (l.q.r or q.r respectively) if each of its elements is right quasi regular (l.q.r or q.r respectively).*

**DEFINITION 5.4.28:** *Let S be a na-bisemiring. An element $x \in S$ is said to be right regular (left regular) if there exists a $y \in S$ ($y' \in S$) such that $x(yx) = x$ (($xy')x = x$).*

**DEFINITION 5.4.29:** *Let S be a na-bisemiring. The Jacobson radical J(S) of a bisemiring S is defined as follows: $J(S) = \{x \in S \mid xS$ is right quasi-regular ideal of $S\}$. A bisemiring is said to be semisimple if $J(S) = \{0\}$ where J(S) is the Jacobson radical of S.*

**DEFINITION 5.4.30:** *Let S be a na-bisemiring. We say S is prime if for any two ideals A, B in S, $AB = 0$ implies $A = 0$ or $B = 0$.*

***Example 5.4.1:*** Let $S = S_1G \cup S_2G$ be a bisemiring where $S_1$ and $S_2$ are semirings given by the following figures:

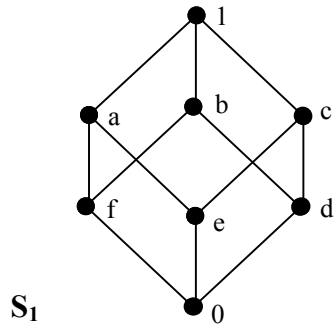

$S_1$

**Figure 5.4.1**

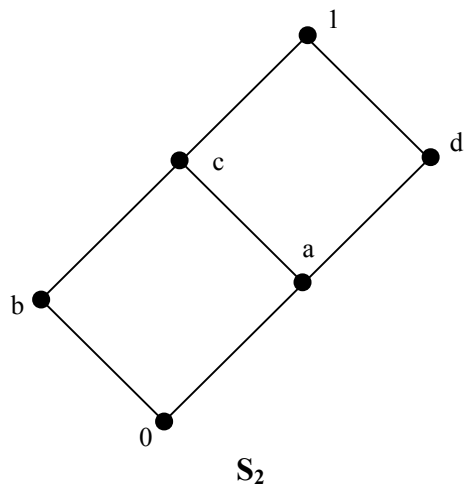

$S_2$

**Figure 5.4.2**

and G is the groupoid given by the following table:

| • | e | $a_0$ | $a_1$ | $a_2$ |
|---|---|---|---|---|
| e | e | $a_0$ | $a_1$ | $a_2$ |
| $a_0$ | $a_0$ | e | $a_2$ | $a_1$ |
| $a_1$ | $a_1$ | $a_2$ | e | $a_0$ |
| $a_2$ | $a_2$ | $a_1$ | $a_0$ | e |

$SG = S_1G \cup S_2G$ is a na-bisemiring having non-trivial zero divisors and idempotents.

Having described the concept non-associative bisemirings we now proceed onto recall the definition of Smarandache non-associative bisemirings.

**DEFINITION 5.4.31:** *Let $(S, +, \bullet)$ be a na-bisemiring. We say S is a Smarandache na bisemiring (S-na bisemiring) if S has a proper subset P, $(P \subset S)$ such that P is an associative bisemiring.*

**DEFINITION 5.4.32:** *Let (S, +, •) be a na- bisemiring. A proper subset P of S is said to be a Smarandache na sub-bisemiring (S-na sub-bisemiring) if P itself under the operations of S is a S-na-bisemiring.*

**DEFINITION 5.4.33:** *Let (S, +, •) be a na bisemiring. We say S is a Smarandache commutative na bisemiring (S-commutative na bisemiring) if S has a S-sub-bisemiring P that is commutative. If every S-sub-bisemiring P of S is commutative we call S a Smarandache strongly commutative bisemiring (S-strongly commutative bisemiring) even if S is commutative and S has no S-sub-bisemiring then also S is not Smarandache commutative (S-commutative).*

**DEFINITION 5.4.34:** *Let (S, +, •) be a na bisemiring. We call S a Smarandache Moufang bisemiring (S-Moufang bisemiring) if S has a proper subset P where P is a S-sub-bisemiring of S and P satisfies the Moufang identities i.e. P is a Moufang bisemiring.*

**DEFINITION 5.4.35:** *Let (S, +, •) be a non-associative bisemiring.*

*If every S-sub-bisemiring P of S satisfies the Moufang identity i.e. every P is a Moufang bisemiring then we call S a Smarandache strong Moufang bisemiring (S-strong Moufang bisemiring).*

**DEFINITION 5.4.36:** *Let (S, +, •) be a S-bisemiring. A fuzzy subset $\mu : S \to [0, 1]$ is said to be a Smarandache fuzzy bisemiring (S-fuzzy bisemiring) if $\mu : P \to [0, 1]$ where P is a bisemifield of S satisfies the following:*

  i. $\mu(x + y) \geq \min\{\mu(x), \mu(y)\}$.
  ii. $\mu(xy) \geq \mu(y)$ and $\mu(xy) \geq \mu(x)$, $x, y \in P \subset S$.

*Thus while defining S-fuzzy bisemirings $\mu$ we have some semifield $P \subset S$ associated with it.*

*Now if $\mu : S \to [0,1]$ is a fuzzy subset of S such that $\mu$ related to every subsemifield $P_i$ contained in S satisfies the conditions*

$$\mu(x + y) \geq \min\{\mu(x), \mu(y)\} \text{ and}$$
$$\mu(xy) \geq \mu(y) \text{ and}$$
$$\mu(xy) > \mu(x)$$

*for all $x, y \in P_i \subset S$ then we call the fuzzy subset $\mu$ to be the Smarandache strong fuzzy bisemiring (S-strong fuzzy bisemiring).*

**THEOREM 5.4.2:** *Every S-strong fuzzy bisemiring is a S-fuzzy bisemiring.*

*Proof*: Straightforward by the very definitions.

We can also define Smarandache fuzzy bisemirings in a different way.

**DEFINITION 5.4.37:** *Let $(S, +, \bullet)$ be a S-bisemiring. $\mu : S \to [0, 1]$ be a fuzzy subset of S suppose $S = S_1 \cup S_2$ is the S - bisemiring having $P = P_1 \cup P_2$ where $P_1 \subset S_1$ and $P_2 \subset S_2$ are semifields i.e. P is a bisemifield of S. We say $\mu$ is a Smarandache fuzzy bisemiring (S-fuzzy bisemiring) if $\mu = \mu_1 \cup \mu_2$ and $\mu_1 : S_1 \to [0, 1]$ is a S-fuzzy semiring relative to $P_1$ and $\mu_2 : S_2 \to [0,1]$ is a S-fuzzy semiring relative to $S_2$.*

The reader is requested to prove or disprove or obtain conditions for the equivalence of the two definitions.

Now we proceed on to define the notion of Smarandache L-fuzzy bi-ideals Smarandache fuzzy k-bi-ideals of k-bisemirings. From now onwards we will let L denote the complete distributive lattice.

$L = (L, \leq, \wedge, \vee)$ will denote a complete distributive lattice. Let $(S = S_1 \cup S_2, +, \bullet)$ be a S-bisemiring. Let F(S) denote the set of all L-fuzzy sets of S i.e. $F(S) = \{\mu : S \to L\}$. For $\mu, \nu \in F(S)$ $\mu \subseteq \nu$ if and only if $\mu(x) \leq \nu(x)$ for all $x \in P = P_1 \cup P_2 \subset S = S_1 \cup S_2$ $\overline{0}(x) = 0$, $T(x) = 1$ for all $x \in P$. It is easily verified $(F(S), \subseteq, \wedge, \vee)$ is again a S- bisemiring.

Now we take S a S-bisemiring F(S) the S-bisemiring associated with S and define S-fuzzy bi-ideals.

**DEFINITION 5.4.38:** *Let $\mu, \in F(S)$. Then $\mu$ is called an Smarandache L-fuzzy left (resp. right) bi-ideal of S if for all $x, y \in P = P_1 \cup P_2 \subset S$ we have*

> i. *$\mu$ is an S-L- fuzzy sub-bisemigroup of $(S, +)$ i.e.*
> ii. *$\mu (x + y) \geq \min \{\mu (x), \mu (y)\}$.*
> iii. *$\mu (x y) \geq \mu (y)$ (resp. $\mu (x y) \geq \mu (x)$).*

*that is we can say if $\mu = \mu_1 \cup \mu_2$ and $\mu_1 : P_1 \to L$ and $\mu_2 : P_2 \to L$ then $\mu_1$ and $\mu_2$ are S-L-fuzzy left (resp. right) ideals of $P_1$ and $P_2$ respectively.*

All results pertaining to S-L-fuzzy bisemirings can be obtained using the definitions of S-L-fuzzy ideals in S-semirings as a matter of routine.

**DEFINITION 5.4.39:** *Let $\mu$ be a S-L-fuzzy bi-ideal of the S-bisemiring R. We say $\mu$ is S-L-fuzzy normal if $\mu_1 (0_1) = 1$ and $\mu_2(0_2) = 1$ where '$0_1$' '$0_2$' are zero elements i.e. additive identity of $P_1$ and $P_2$ respectively.*

**DEFINITION 5.4.40:** *A S-fuzzy left bi-ideal A of a S-bisemiring S is called a S-fuzzy left k-bi-ideal of S if for any $x, y, z \in P \subset S$, $x + y = z$ implies $A(x) \geq \min \{A(y), A(z)\}$ where A is defined relative to P, the bisemifield of S. Let A and B be two fuzzy subsets of a S-bisemiring $S = S_1 \cup S_2$. The k-product $A \circ_k B$ is defined by for $A = A^1 \cup A^2$ and $B = B^1 \cup B^2$ where $A^1 : S_1 \to [0, 1]$ $A^2 : S_2 \to [0, 1]$, $B^1 : S_1 \to [0, 1]$ and $B^2 : S_2 \to [0, 1]$*

$$A^{j^o} o_k B^j(x) = \begin{cases} \sup_{xa_1b_1=a_2b_2} [\min\{A^j(a_i), B^j(b_i)\}, i = 1,2] \\ 0 \text{ if } x \text{ cannot be written as } x + a_1b_1 = a_2b_2 \end{cases}$$

*for j = 1, 2.*

One can derive several interesting results in this case like $A^{j^o} o_k B^j(x) \subseteq A^j \cap B^j$, j = 1, 2. But we proceed on to define semiregular bisemirings and S-semiregular bisemirings.

**DEFINITION 5.4.41:** *A bisemiring S is said to be semiregular if for each $a_i \in S_i$ ( i = 1, 2 where $S = S_1 \cup S_2$) there are $x_{1i}, x_{2i} \in S_i$; i = 1, 2; such that $a_1 + a_1x_{11}a_1 = a_1x_{21}a_1$ and $a_2 + a_2x_{12}a_2 = a_2x_{22}a_2$.*

We now proceed on to define Smarandache semiregular bisemiring.

**DEFINITION 5.4.42:** *A bisemiring ($S = S_1 \cup S_2$, +, •) is said to be Smarandache semiregular bisemiring (S-semiregular bisemiring) if S has a proper subset $P = P_1 \cup P_2$ where P is a semiregular sub-bisemiring of S.*

*For the bisemiring S we define k-closure of $A \subset S$ by*

$$\overline{A} = \begin{cases} x_1 \in S_1 \\ x_2 \in S_2 \end{cases}$$

*such that*

$$x_1 + a_{11} = a_{21}$$
$$x_2 + a_{12} = a_{22};$$

*where $A = A_1 \cup A_2$ with $A_1 \subset S_1$, $A_2 \subset S_2$ and $a_{11}, a_{21} \in A_1$ and $a_{12}, a_{22}, \in A_2$} where have $\overline{A} = \overline{A}_1 \cup \overline{A}_2$. Thus we can equivalently formulate the k-closure $\overline{A}$ of A of a bisemiring $S = S_1 \cup S_2$ as $A_1 \subset S_1$, has $\overline{A}_1$ and $A_2 \subset S_2$ has $\overline{A}_2$ where $S_1$ and $S_2$ are semirings and $\overline{A} = \overline{A}_1 \cup \overline{A}_2$ is a k-closure of the bisemiring.*

*It is left for the reader to prove in case of bisemiring $S = S_1 \cup S_2$, A, B $\subset$ S*

$$\overline{AB} = \overline{\overline{A}\overline{B}}.$$

*Further if $\overline{A}$ is the smallest left k-bi-ideal of S containing A and if in a bisemiring S, if A and B are right and left bi-ideals of the bisemiring S then it can be proved $\overline{AB} = A \cap B$.*

*Let S be a bisemiring. A k-bi-ideal I of S is called bi prime if $I \neq S$ and for any two bi-ideals A and B of S, $AB \subset I$ implies either $A \subset I$ or $B \subset I$. We say in the S-bisemiring S*

*a k-bi-ideal I of S to be S- bi prime if I is a bi prime bi-ideal relative to a bisemifield P in S.*

*We say a fuzzy k-bi-ideal P of a bisemiring S is said to be bi prime if P is not a constant function and for any two fuzzy bi-ideals A and B of S A $o_k$ B $\subseteq$ P implies either A $\subseteq$ P of B $\subset$ P.*

*We say the fuzzy k bi-ideal P of S-bisemiring S is a Smarandache fuzzy k-bi-ideal (S-fuzzy k-bi-ideal) of S if P is a fuzzy k-bi-ideal relative a proper subset T $\subset$ S where T is a bisemifield under the operations of S.*

Several interesting results in this direction can be obtained by any innovative reader.

Now we know in case of Smarandache bisemirings of level II we can have all Smarandache fuzzy concepts related to fuzzy fields. All the more several types of these concepts can also developed in case of Smarandache anti bisemifield.

We know a Smarandache na-bisemiring has a subset which is an associative bisemiring. Thus we have in case of S-na-bisemiring all properties true of bisemirings can be easily adopted and studied.

For we can define the very fuzzy set µ from S to [0, 1] by µ = $\mu_1 \cup \mu_2$ where $\mu_1$ : $S_1$ → [0, 1] and $\mu_2$ : $S_2$ → [0, 1]. So for defining the notions of Smarandache L-fuzzy bi-ideals µ in case of S-na-bisemirings, we can define µ = $\mu_1 \cup \mu_2$ : $S_1 \cup S_2$ → [0, 1] and by $\mu_1$ : $S_1$ → [0, 1] and $\mu_2$ : $S_2$ → [0, 1] to be S-fuzzy k ideals and similarly by defining $\mu_i$ : $S_i$ → [0, 1] for i = 1,2 as Smarandache normal L-fuzzy ideals in semirings; so that µ = $\mu_1 \cup \mu_2$ be comes Smarandache normal L-fuzzy bi-ideals in S-bisemirings.

Thus in case of Smarandache non-associative bisemirings, we can define and derive the notions of S-fuzzy k-bi-ideals and S-fuzzy k-bi-ideals in S-k-bisemirings.

## 5.5 Problems

This section gives fifty four problems for the reader to solve. By solving these problems the author hopes that the reader would have mastered both Smarandache semirings and Smarandache fuzzy semirings. The problems will certainly throw light on the subject.

**Problem 5.5.1:** Give an example of a S-strong fuzzy semiring .

**Problem 5.5.2:** Construct an example of a S-fuzzy semiring which is not a S-strong fuzzy semiring.

**Problem 5.5.3:** Give an example of a S-fuzzy subsemiring which is not a S-fuzzy ideal.

**Problem 5.5.4:** What is relation between S-fuzzy subsemiring and S-fuzzy ideals of a semiring?

**Problem 5.5.5:** If $\mu$ is a S-fuzzy ideal of a S-semiring relative to $P \subseteq S$. Will $\mu(x + y) = \mu(0) \Rightarrow \mu(x) = \mu(y)$ for all $x, y \in P \subset R$?

**Problem 5.5.6:** Let $\mu$ be a fuzzy subset of a S-semiring S. Is it true $\mu$ of the S-semiring S is a S-fuzzy-subsemiring (S-fuzzy ideal) of S if and only if the level subsets $\mu_t$, $t \in \text{Im } \mu$ are S-subsemirings (S-ideals) of S relative to a fixed proper subset P in S, P a semifield in S?

**Problem 5.5.7:** Let S and S' be S-semirings and f a S-homomorphism of the S-semirings from S onto S' then prove for each S-fuzzy subsemiring (S-fuzzy ideal) related to a subsemifield $P \subset S$. $f(\mu)$ is a S-fuzzy subsemiring (S-fuzzy ideal) of S' and for each S-fuzzy subsemiring (S-fuzzy ideal) $\mu'$ of S', $f^{-1}(\mu')$ is a S-fuzzy subsemiring. (S-fuzzy ideal) of S.

**Problem 5.5.8:** Is it true if S is a S-regular semiring then $\sigma\theta = \sigma \cap \theta$ where $\sigma$ and $\theta$ are S-fuzzy ideals of S relative to a fixed semifield P in S. If $\sigma\theta = \sigma \cap \theta$ does it imply S is a S-regular semiring?

**Problem 5.5.9:** If $\mu$ is any S-fuzzy prime ideal of a S-semiring S then prove the ideal $\mu_t$, $t = \mu(0)$ is a S-prime ideal of S.

**Problem 5.5.10:** Prove for a S-semiring S the two S-level subsemirings (S-level ideals) $\mu_s$ and $\mu_t$ (with $s < t$) of a S-fuzzy subsemiirng (S-fuzzy ideal) of $\mu$ of the S-semiring S are equal if and only if there is no x in $P \subset R$ such that $s \leq \mu(x) < t$.

**Problem 5.5.11:** Prove / disprove the intersection of any family of S-fuzzy subsemirings (S-fuzzy ideals) of a S-semiring S is a S-fuzzy subsemiring (S-fuzzy ideal) of S.

**Problem 5.5.12:** Let $\mu$ be any S-fuzzy ideal of a S-subsemiring such that Im $\mu = \{t\}$ or $\{0, s\}$ where $t = [0, 1]$ and $s \in (0, 1]$. If $\mu = \sigma \cup \theta$ where $\sigma$ and $\theta$ are S-fuzzy ideals of a S-semiring S defined relative to the same P then prove either $\sigma \subseteq \theta$ or $\theta \subseteq \sigma$.

**Problem 5.5.13:** Let f be a S-homomorphism from the S-semiring S onto a S-semiring S'. If $\mu$ and $\sigma$ are S-fuzzy ideals of S then the prove the following are true:

  i. $f(\mu + \sigma) = f(\mu) + f(\sigma)$.
  ii. $f(\mu\sigma) = f(\mu) f(\sigma)$.
  iii. $f(\mu \cap \sigma) \subseteq f(\mu) \cap f(\sigma)$.

with equality if atleast one of $\mu$ are $\sigma$ is f-invariant ($\mu$ and $\sigma$ defined relative to the same semifield P in S).

**Problem 5.5.14:** Suppose $\mu$ is a S-L-fuzzy left (respectively right) ideal of S relative to $P \subset S$ (P a semifield), will $S_\mu = \{x \in P \mid \mu(x) = \mu(0)\}$ is a left (resp. right) S- ideal of P? Justify your claim.

**Problem 5.5.15:** Prove if $\mu$ is a S-L-fuzzy left (resp. right) ideal of a S-semiring S satisfying $\mu^+(x) = 0$ for some $x \in P \subset S$ (P is the subsemifield relative to which $\mu$ is defined). Prove $\mu(x) = 0$.

**Problem 5.5.16:** Let $\mu$ and $\nu$ be S-L-fuzzy left (resp. right) ideals of S relative to P a semifield contained in S. If $\mu \subset \nu$ and $\mu(0) = \nu(0)$ then prove $P_\mu \subset P_\nu$.

**Problem 5.5.17:** If $\mu$ and $\nu$ are S-normal L-fuzzy left (resp. right) ideals of the S-semiring S containing the semifield P relative to which $\mu$ and $\nu$ are defined, satisfying $\mu \subset \nu$ then prove $P_\mu \subset P_\nu$.

**Problem 5.5.18:** An S-L-fuzzy left (right) ideal $\mu$ of a S-semiring S, prove is S-normal if and only if $\mu^- = \mu$.

**Problem 5.5.19:** $\mu$ a S-L-fuzzy left (right) ideal of S.

    i.                         Will $(\mu^-)^- = \mu^-$?
    ii.                        Will $(\mu^{-1})^+ = \mu$?

**Problem 5.5.20:** Let $\mu$ be a S-L-fuzzy left (resp. right) ideal of a S-semigroup S relative to $P \subset S$. If there exists a S-L-fuzzy left (resp. right) ideal $\nu$ of S relative to same P satisfying $\nu^- \subset \mu$ then prove $\mu$ is S-normal.

**Problem 5.5.21:** Let $\mu$ be a S-L-fuzzy left (resp. right) ideal of S relative to P, $P \subset S$. If there exists a S-L-fuzzy left (respectively right) ideal $\nu$ of S relative to $P \subset S$ satisfying $\nu^{-1} \subset \mu$ then prove $\mu^+ = \mu$.

**Problem 5.5.22:** Let S be a k-semiring, I be an ideal and $I' = \{a' \in S' \mid a \in I\}$. Then prove I is a k-ideal of S if and only if $\bar{I} = I \cup I'$ is an ideal of the extension ring $\bar{S}$ called the extension ideal of I.

**Problem 5.5.23:** Let I be a S-k-ideal of a S-k-semiring S, $\bar{I}$ the S-extension ideal of the S-ring $\bar{S}$, Will $a \oplus I = (a + \bar{\bar{I}}) \cap S$ where $a \in S$?

**Problem 5.5.24:** Let I be a S-k-ideal of the S-k-semiring S. Then prove $S / I = \{a \oplus I \mid a \in P \subset S\}$ is a S-k-semiring under the operations $(a \oplus I) \oplus (b \oplus I) = a \oplus b \oplus I$ and $(a \oplus I) \odot (b \oplus I) = (ab) \oplus I$.

**Problem 5.5.25:** Let A be a S-fuzzy ideal of a S-semiring S. Then the level set $A_t$ ($t \leq A(0)$) is the S-ideal of S.

**Problem 5.5.26:** Let A be a S-fuzzy ideal of a S-semiring S. IF $A_t$ is a S-k-ideal of S for each $t \leq A(0)$. Then prove A is a S-k-fuzzy ideal of S.

**Problem 5.5.27:** Let $f : R \to S$ be an S-epimorphism of semirings (R and S are S-semirings) and A an f-invariant S-fuzzy ideal of R. Then prove f(A) is a S-fuzzy ideal of S.

**Problem 5.5.28:** A S-fuzzy ideal A of a S-semiring S is called S-fuzzy semiprimary if for all a, b $\in$ S either $A(ab) \leq A(a^n)$ for some $n \in Z_+$ or else $A(ab) \leq A(b^m)$ for some $m \in Z_+$ where $Z_+$ is the set of all non-negative integers.

**Problem 5.5.29:** Prove an S-ideal I of a S-semiring R with identity is S-semiprimary if and only if $\chi_I$, the characteristic function of I is a S-fuzzy semiprimary ideal of R.

**Problem 5.5.30:** If R is a S-commutative semiring with identity, then prove A is any S-semiprimary fuzzy ideal of R if and only if $A_t$, where $t \in$ Im (A) Im is a semiprimary ideal of R.

**Problem 5.5.31:** Prove if A is a S-fuzzy semiprime ideal of a S-k-semiring R such that all of its level subsets are S-k-ideals relative to the same $P \subset R$. Prove R/A has no non-zero nilpotent elements.

**Problem 5.5.32:** Let A be a S-fuzzy primary ideal of the S-semiring R such that all of its level subsets are S-k-ideals of R. Then prove every zero divisor of R/A is nil potent.

**Problem 5.5.33:** Let $f : R \to S$ be an S-epimorphism of S-k-semirings and B a S-fuzzy ideal of S. Then prove B is a S-fuzzy maximal k-ideal of S if and only if $f^{-1}B$ is a S-fuzzy maximal k-ideal of R relative to P.

**Problem 5.5.34:** Prove for a S-semiring S and A, B $\subset$ P $\subset$ S; we have $\overline{AB} = \overline{A}\,\overline{B}$.

**Problem 5.5.35:** Suppose A and B are respectively right and left S-k-ideals of a S-semiring S then $\overline{AB} \subseteq A \cap B$.

**Problem 5.5.36:** Prove a S-semiring S is S-semigular if and only if for any right S-k-ideal A and for any left S-k-ideal B, $\overline{AB} = A \cap B$.

**Problem 5.5.37:** Prove any S-semiring S, is S-semiregular if and only if for any S-fuzzy right k-ideal A and any S-fuzzy left k-ideal B, $A o_k B = A \cap B$.

$$(A \circ_k B)(x) = \begin{cases} \sup_{x+a_1b_1=a_2b_2} [\min\{A(a_i), B(b_i) \mid i=1, 2\}] \\ 0 \text{ if } x \text{ cannot be written as } x + a_1b_1 = a_2b_2. \end{cases}$$

where A and B are two fuzzy subsets of a S-semiring.

**Problem 5.5.38:** Prove in case of S-fuzzy semivector spaces, for $c, d \in S$, $x, y \in H \subset G$, $0 \le k, \lambda, \nu, \mu \le 1$. Then $d_\mu \circ x_\lambda = (dx)_{\min(\mu,\lambda)}$, $x_\lambda + y_\nu = (x+y)_{\min(\lambda,\nu)}$, $d\mu \circ x_\lambda + c_R \circ y_\nu = (dx + cy)_{\min(\mu, \lambda, k, \nu)}$.

**Problem 5.5.39:** Prove the fuzzy subset $\nu$ of SF (Hom $(\mu, \lambda)$) is a S-fuzzy subsemispace of SF Hom $(\mu, \lambda)$.

**Problem 5.5.40:** Give an example of a S-strong fuzzy bisemigroup.

**Problem 5.5.41:** Illustrate by an example a S-fuzzy bisemigroup which is not a S-strong fuzzy bisemigroup.

**Problem 5.5.42:** Prove if $\mu$ is a S-L-fuzzy left (resp. right) bi-ideal of the S-bisemiring S and let $\mu^+$ be a L-fuzzy set in S defined by $\mu^+(x) = \mu(x) + 1 - \mu(0)$ for all $x \in P \subset S$. Prove $\mu^+$ is S-L-fuzzy normal left (resp. right) bi-ideal of S containing $\mu$.

**Problem 5.5.43:** Let $\mu$ be a S-L-fuzzy left (resp. right) bi-ideal of S satisfying $\mu^+(x) = 0$ for some $x \in P \subset S$ then prove $\mu(x) = 0$.

**Problem 5.5.44:** Let $\mu$ and $\upsilon$ be S-L-fuzzy left (resp. right) bi-ideals of S (S a S-bisemiring. If $\mu \subset \nu$ and $\mu(0) = \nu(0)$ prove $S_\mu \subset S_\nu$.

**Problem 5.5.45:** Prove in case of S-L-fuzzy left (or right) bi-ideals of $S(\mu^+)^+ = \mu^+$. Also prove $(\mu^-)^+$ if $\mu$ is S-L-fuzzy normal left (resp. right) bi-ideal of S.

**Problem 5.5.46:** Prove a bisemiring S is semiregular if and only if for any right k-bi-ideal A and for any left k-bi-ideal B, $\overline{AB} \subseteq A \cap B$.

**Problem 5.5.47:** A S- bisemiring S is S-semiregular if and only if for any S-fuzzy k-bi-ideal A and any fuzzy k-bi-ideal B, $A \circ_k B = A \cap B$.

**Problem 5.5.48:** Define for a S-bisemiring the S-fuzzy pseudo ideal. Illustrate this by an example.

**Problem 5.5.49:** Does there exist an example of a non-associative semiring other than the semirings got using loops and groupoids i.e. loop semirings and groupoid semirings?

**Problem 5.5.50:** If S is a S-na-semiring and A, B ⊆ P ⊆ S have the S-k-closure $\overline{A}$ and $\overline{B}$, will $\overline{AB} = \overline{\overline{AB}}$?

**Problem 5.5.51:** Let A and B be S-right and S-left k-ideals of the S-na-semiring S, will $\overline{AB} \subseteq A \cap B$? Justify your claim.

**Problem 5.5.52:** Prove or disprove. A S-na-semiring S is S-semiregular if and only if for any S-right k-ideal A and for any S-left k-ideal B, $\overline{AB} = A \cap B$.

**Problem 5.5.53:** Prove or disprove if μ is a S-fuzzy k-ideal of S, S a S-na-semiring then $\mu_t = \{x \in S \mid \mu(x) \geq t\}$ is a S-fuzzy k-ideal.

**Problem 5.5.54:** Let A be a S-fuzzy k-ideal of S, S a S-fuzzy na semiring S. Then $A(x) \leq A(0)$ for all $x \in S$.

**Chapter Six**

# SMARANDACHE FUZZY NEAR-RING AND ITS PROPERTIES

This chapter has four sections. In the first section we recall the definition of Smarandache near-rings and define the new notion of Smarandache fuzzy near-rings. Properties like Smarandache normal fuzzy-R-subgroups and S-fuzzy congruence of a near-ring module are defined in section 1. In section two we define these notions in the context of non-associative near-rings and study them. Study of bistructures is very important and inevitable as we do not see any proper structure given by union of two algebraic structures and union of two distinct algebraic structures. So section three is devoted to introduction and study of Smarandache fuzzy bistructures. The final section is completely devoted to giving problem about S-fuzzy near-rings and their properties. These problems will make the researcher to get more ideas and they can construct more and more Smarandache notions and find also suitable applications about them.

## 6.1 Smarandache Fuzzy Near-rings

In this section we just recall the definition of Smarandache near-rings and introduce the notions of Smarandache fuzzy near-rings and their properties. The concept of Smarandache fuzzy congruence of a near-ring module and Smarandache normal fuzzy R-subgroups in near-rings are mainly introduced and studied in this section.

As in this book we have only recalled the concept of fuzzy near-rings and its properties, in this section we recall the notions of Smarandache near-rings and Smarandache seminear-rings. We also define for the first time both notions of Smarandache fuzzy seminear-rings and fuzzy seminear-rings.

**DEFINITION 6.1.1:** *N is said to be a Smarandache near-ring (S-near-ring) if (N, +, •) is a near-ring and has a proper subset A such that (A, +, •) is a near-field.*

**Example 6.1.1:** Let $Z_2 = \{0, 1\}$ be the near-field. Take any group G such that $Z_2G$ is the group near-ring of the group G over the near-field $Z_2$. $Z_2G$ is a S-near-ring as $Z_2 \subset Z_2G$ and $Z_2$ is a near-field.

It is important to note that in case of S-near-rings we can get several or a new class of S-near-ring by using group near-rings and semigroup near-rings using basically the near-field $Z_2$ or any $Z_p$.

**Example 6.1.2:** Let $Z_2 = \{0, 1\}$ be a near-field and S(n) any symmetric semigroup. The semigroup near-ring $Z_2S(n)$ is a S-near-ring. This S-near-ring can be of any order finite or infinite depending on the order of the group or the semigroup which is used. Similarly they can be commutative or non-commutative which depends basically on the groups and semigroups.

**DEFINITION 6.1.2:** *Let (N, +, •) be a S-near-ring. A non-empty proper subset T of N is said to be a Smarandache subnear-ring (S-subnear ring) if (T, +, •) is a S-near-ring; i.e. T has a proper subset which is a nearfield.*

**DEFINITION 6.1.3:** *Let (P, +) be a S-semigroup with zero 0 and let N be a S-near-ring. Let $\mu : N \times Y \to Y$ where Y is a proper subset of P which is a group under the operations of P. (P, $\mu$) is called the Smarandache N-group (S-N-group) if for all $y \in Y$ and for all n, $n_1 \in N$ we have $(n + n_1) y = ny + n_1 y$ and $(nn_1) y = n(n, y)$. $S(N^P)$ stands for S-N-groups.*

**DEFINITION 6.1.4:** *A S-semigroup M of a near-ring N is called a Smarandache-quasi subnear-ring (S-quasi subnear-ring) of N if $X \subset M$ where X is a subgroup of M which is such that $X.X \subseteq X$.*

**DEFINITION 6.1.5:** *A S-subsemigroup Y of $S(N^P)$ with $NY \subset Y$ is said to be a Smarandache N-subgroup (S-N-subgroup) of P.*

**DEFINITION 6.1.6:** *Let N and $N_1$ be two S-near-rings. P and $P_1$ be S-N-subgroups.*

i. $h : N \to N_1$ is called a Smarandache near-ring homomorphism (S-near-ring homomorphism) if for all m, $n \in M$ (M is a proper subset of N which is a near-field) we have

$$h (m + n) = h(m) + h(n),$$
$$h (mn) = h(m) h(n) \text{ where } h(m) \text{ and}$$

$h(n) \in M_1$ ($M_1$ is a proper subset of $N_1$ which is a near-field). It is to be noted that h need not even be defined on whole of N.

ii. $h : P \to P_1$ is called the Smarandache N-subgroup-homomorphism (S-N-subgroup-homomorphism) if for all p, q in S (S the proper subset of P which is a S-N-subgroup of the S-semigroup P) and for all $m \in M \subset N$ (M a subfield of N) $h(p + q) = h(p) + h(q)$ and $h(mp) = mh(p)$; $h(p)$, $h(q)$ and $mh(p) \in S_1$ ($S_1$ a proper subset of $P_1$ which is a S-$N_1$ subgroup of S-semigroup $P_1$).

*Here also we do not demand h to be defined on whole of P.*

**DEFINITION 6.1.7:** *Let N be a S-near-ring. A normal subgroup I of (N, +) is called a Smarandache ideal (S-ideal) of N related to X if*

i. $IX \subseteq I$.
ii. $\forall x, y \in X$ and for all $i \in I$; $x(y + i) - xy \in I$,

*where X is the near-field contained in N.*

A subgroup may or may not be a S-ideal related to all near-fields. Thus while defining S-ideal it is important to mention the related near-field.

Normal subgroup T of (N, +) with (a) is called S-right ideals of N related to X while normal subgroups L of (N, +) with (b) are called S-left ideals of N related to X.

**DEFINITION 6.1.8:** *A proper subset S of P is called a Smarandache ideal of $S(N^P)$ (S-ideal of $S(N)^P$) related to Y if*

  i.   *S is a S-normal subgroup of the S-semigroup P.*
  ii.  *For all $s_1 \in S$ and $s \in Y$ (Y is the subgroup of P) and for all $m \in M$ (M the near-field of N)*
$$n ( s + s_1) - ns \in S.$$

*A S-near ring is Smarandache simple (S-simple) if it has no S-ideals. $S(N^P)$ is called Smarandache N-simple (S-N-simple) if it has no S-normal subgroups except 0 and P.*

**DEFINITION 6.1.9:** *A S-subnear-ring M of a near-ring N is called Smarandache invariant (S-invariant) related to the near-field X in N if $M X \subset M$ and $X M \subset M$ where X is a S-near-field of N. Thus in case of S-invariance it is only a relative concept as a S-subnear-ring M may not be invariant related to every near-fields in the near-ring N.*

**DEFINITION 6.1.10:** *Let X and Y be S-semigroups of $S(N^P)$. $(X{:}Y) = \{n \in M \,/\, nY \subset X\}$ where M is a near-field contained in N. (0, x) is called the Smarandache annihilator (S-annihilator) of X.*

**DEFINITION 6.1.11:** *Let N be a near-ring and S a S-subsemigroup of (N, +). The near-ring $N_s$ is called the Smarandache-near-ring of left (right) quotients (S-near-ring of left (right) quotients) of N with respect to S if*

  i.   *$N_s$ has identity.*
  ii.  *N is embeddable in $N_s$, by a homomorphism h.*
  iii. *For all $s \in S$; h(s) is invertible in $(N_s, \bullet)$.*
  iv.  *For all $q \in N_s$, there exists $s \in S$ and there exist $n \in N$ such that $q = h(n) h(s)^{-1}$ $(q = h(s)^{-1} h(n))$.*

The problem whether $N_s$ is a S-near-ring is left as an open problem.

**DEFINITION 6.1.12:** *The near-ring N is said to fulfil the Smarandache left (right) ore conditions (S-left(right) ore condition) (ore (1)) with respect to a given S-subsemigroup P of (N, $\bullet$) if for $(s, n) \in S \times N$ there exists $ns_1 = sn_1$ $(s_1 n = n_1 s)$.*

**DEFINITION 6.1.13:** *If $S = \{s \in N \,/\, s$ is cancellable$\}$, the $N_S$ if it exists and if $N_S$ is a S-near-ring then $N_s$ is called the Smarandache left (right) quotient (S-left(right) quotient) near-ring of N.*

**THEOREM 6.1.1:** *If $N_S$ is a S-quotient near-ring then $N_S$ is a quotient near-ring.*

*Proof*: Obvious by the very definition.

It is once again an interesting problem to obtain a necessary and sufficient condition for a left (right) quotient near-ring N to be a S-quotient near-ring.

Let S(V) denote the collection of all S-near-rings and X be any nonempty subset which is S-semigroup under '+' or '•'.

We define Smarandache free near-ring as follows.

**DEFINITION 6.1.14:** *A S-near-ring $F_x \in S(V)$ is called a Smarandache free near-ring (S-free near-ring) in V over X if there exists $f : X \to F$ for all $N \in V$ and for all $g : X \to N$ there exists a S-near-ring homomorphism $h \in S (Hom\ F_x , N)$ [Here S (Hom ($F_x$ , N)) denotes the collection of S-homomorphism from $F_x$ to N] such that $h \circ f = g$.*

**THEOREM 6.1.2:** *Let $F_x$ be a S-free near-ring then $F_x$ is a free near-ring.*

*Proof*: Obvious by the very definition.

**DEFINITION 6.1.15:** *$z \in N$ is called Smarandache quasi regular (S-quasi regular) if $z \in L_z$. An S-ideal $P \subset N$ is called S-quasi regular if and only if for all $s \in P$, s is S-quasi regular.*

**DEFINITION 6.1.16:** *A ring R is said to be Smarandache biregular (S-biregular) if each S-principal ideal is generated by an idempotent.*

**DEFINITION 6.1.17:** *Let N be a near-ring. N is said to by Smarandache biregular (S-biregular) if there exists some set E of central S-idempotents with*

    i.        For all $e \in N$; Ne is an S-ideal of N.
    ii.       For all $n \in N$ there exists $e \in E$; Ne = (n) (principal ideal generated by n).
    iii.      For all $e, f \in E$   $e + f = f + e$.
    iv.      For all $e, f \in E$, $ef \in E$ and $e + f - ef \in E$.

**DEFINITION 6.1.18:** *Let I be an S-ideal. The intersection of all S-prime ideal P, such that $I \subseteq P$ is called the Smarandache prime radical (S-prime radical) of I, i.e. $S(I) = \bigcap_{I \subseteq P} P$.*

**THEOREM 6.1.3**: *Let I be a S-ideal of N. S(I) be the S-prime radical of I then S(I) is a prime radical of I.*

Now we proceed on to define the notions of Smarandache seminear-rings.

**DEFINITION 6.1.19**: *Let N be a seminear-ring. N is said to be a Smarandache seminear-ring of level II (S-seminear-ring of level II) if N contains a proper subset P which is a semiring. Clearly the S-mixed direct product gives S-seminear-ring II.*

**DEFINITION 6.1.20**: *Let N be a S-near-ring. N is said to fulfil the Smarandache insertion of factors property (S-IFP for short) if for all a, b ∈ N we have a . b = 0 implies anb = 0 for all n ∈ P, P ⊂ N, where P is a near-field.*

We say N has S-strong IFP property if every homomorphic image of N has the IFP property or we redefine this as:

**DEFINITION 6.1.21**: *Let N be a S-near-ring and I a S-ideal of N. N is said to fulfil Smarandache strong IFP property (S-strong IFP property) if and only if for all a, b ∈ N, ab ∈ I implies anb ∈ I where n ∈ P, P ⊂ N and P is a near-field.*

**DEFINITION 6.1.22**: *Let N be a near-ring. A S-right ideal I of N is called Smarandache right quasi reflexive (S-right quasi reflexive) if whenever A and B are S-ideals of N with AB ⊂ I, then b(b' + a) – bb' ∈ I for all a ∈ A and for all b, b' ∈ B.*

**DEFINITION 6.1.23**: *Let N be a near-ring. N is said to be Smarandache strongly subcommutative (S-strongly subcommutative) if every S-right ideal of it is S-right quasi reflexive.*

**DEFINITION 6.1.24**: *Let N be a near-ring. S a S-subnormal subgroup of (N, +). S is called a Smarandache quasi-ideal (S-quasi ideal) of N if SN ⊂ S and NS ⊂ S where by NS we mean elements of the form {n(n' + s) – nn'/ for all s ∈ S and for all n, n' ∈ N} = NS.*

**DEFINITION 6.1.25**: *A left-near-ring N is said to be Smarandache left-self-distributive (S-left self distributive) if the identity abc = abac is satisfied for a, b, c ∈ A, A ⊂ N and A is a S-subnear-ring of N.*

**DEFINITION 6.1.26**: *Let R be a S-near-ring. R is said to be Smarandache equiprime (S-equiprime) if for all 0 ≠ a ∈ P ⊂ R. P the near-field in R and for x, y ∈ R, arx = ary for all r ∈ R implies x = y. If B is a S-ideal of R, B is called a Smarandache equiprime ideal (S-equiprime ideal) if R/B is an S-equiprime near-ring.*

**DEFINITION 6.1.27**: *Let (N, +, •) be a triple. N is said to be a Smarandache infra-near-ring (S-INR) where*

  i. *(S, +) is a S-semigroup.*
  ii. *(N, •) is a semigroup.*
  iii. *(x + y)z = xz – 0z + yz for all x, y, z ∈ N.*

**DEFINITION 6.1.28**: *Let I be a S-left ideal of N. Suppose I satisfies the following conditions:*

  i. *a, x, y ∈ N, anx – any ∈ I for all n ∈ N implies x – y ∈ I.*
  ii. *I is left invariant.*
  iii. *0N ⊂ I.*

*Then I is called a Smarandache equiprime left-ideal (S-equiprime left ideal) of N.*

**DEFINITION 6.1.29**: *Let N be a S-near-ring. N is said to be Smarandache partially ordered (S-partially ordered) by $\leq$ if*

  i.  *$\leq$ makes (N, +) into a partially ordered S-semigroup.*
  ii. *for all n, n' $\in$ N, n $\geq$ 0 and n' > 0 implies nn' $\geq$ 0.*

**DEFINITION 6.1.30**: *Let M be a S-$\Gamma$-near-ring, then a S-normal subsemigroup I of (M, +) is called*

  i.   *a S-left ideal if $a\alpha(b + i) - a\alpha b \in I$ for all a, b $\in$ M, $\alpha \in \Gamma$ and i $\in$ I.*
  ii.  *a S-right ideal if $i \alpha a \in I$ for all a $\in$ M, $\alpha \in \Gamma$ and i $\in$ I and*
  iii. *a S-ideal if it is both S-left and S-right ideal.*

**DEFINITION [130]:** *A non-empty set N is said to be a Smarandache seminear-ring (S-seminear-ring) if (N, +, •) is a seminear-ring having a proper subset A, (A $\subset$ N) such that A under the same operations of N is a near-ring, that is (A, +, •) is a near-ring.*

**DEFINITION 6.1.31:** *Let (N, +, •) be a seminear-ring. If in the S-semigroup (N, +) every proper subset (A, +) which is a group is commutative then we say the S-seminear-ring N is Smarandache commutative (S-commutative). Thus if (N, +, •) is commutative and if (N, +) is S-semigroup then trivially (N, +, •) is S-commutative. Secondly if (N, +, •) is S-commutative seminear-ring then N need not in general be commutative.*

**DEFINITION 6.1.32**: *Let (N, +, •) be a seminear-ring. If (N, +) is a S-semigroup and if N has at least one proper subset which is a subgroup that is commutative then we say the seminear-ring (N, +, •) is a Smarandache weakly commutative seminear-ring (S-weakly commutative seminear-ring).*

**THEOREM 6.1.4**: *Let (N, +, •) be a seminear-ring which is S-commutative then N is S-weakly commutative.*

*Proof*: Proof is direct and the reader is expected to prove.

**DEFINITION 6.1.33**: *Let (N, +, •) be a seminear-ring. (N, +) be a S-semigroup such that every proper subset A of N which is a group is a cyclic subgroup then we say the seminear-ring (N, +, •) is a Smarandache cyclic seminear-ring (S-cyclic seminear-ring). In particular (N, +, •) has at least one proper subset which is a cyclic group, we call (N, +, •) a Smarandache weakly cyclic seminear-ring (S-weakly cyclic seminear-ring).*

**DEFINITION 6.1.34:** *Let N and $N_1$ be any two S-seminear-rings. We say a map $\phi$ from N to $N_1$ is a Smarandache seminear-ring homomorphism (S-seminear-ring homomorphism) from A to $A_1$ where A $\subset$ N is a near-ring and $A_1 \subset N_1$ is a near-ring and $\phi(x + y) = \phi(x) + \phi(y)$, $\phi(xy) = \phi(x)\phi(y)$ where x, y $\in$ A and $\phi(x), \phi(y) \in A_1$ and it is true for all x, y $\in$ A. We need not even have the map $\phi$ to be well-defined or even defined on the whole of N.*

*The concept of Smarandache isomorphism (S-isomorphism) etc. are defined in a similar way.*

**DEFINITION 6.1.35:** *Let (N, +, •) be a seminear-ring. We say N is a Smarandache seminear-ring of level II (S-seminear-ring II) if N has a proper subset P where P is a semiring.*

**DEFINITION [130]:** *N is said to be a Smarandache pseudo seminear-ring (S-pseudo seminear-ring) if N is a near-ring and has a proper subset A of N which is a seminear-ring under the operations of N.*

**DEFINITION 6.1.36:** *Let N and $N_1$ be the S-pseudo seminear-ring. A mapping $\phi : N \to N_1$ is said to be a Smarandache pseudo seminear-ring homomorphism (S-pseudo seminear-ring homomorphism) if $\phi$ restricted from A to $A_1$ is a seminear-ring homomorphism where A and $A_1$ are proper subsets of N and $N_1$ which are seminear-rings. Thus $\phi$ need not even be defined on the whole of N.*

We define the concept of S-pseudo subseminear-ring and Smarandache pseudo ideals of a near-ring.

**DEFINITION 6.1.37**: *Let N be a near-ring if N has a proper subset A which is a subnear-ring A and if A itself a S-pseudo seminear-ring then we say A is a Smarandache pseudo subseminear-ring (S-pseudo subseminear-ring).*

**DEFINITION 6.1.38**: *Let N be a near-ring. A proper subset M of N is said to be a Smarandache pseudo ideal (S-pseudo ideal) if M is a S-ideal of the near-ring N.*

**THEOREM 6.1.5:** *Let N be a near-ring. If N has a S-pseudo subseminear-ring then N is a S-pseudo seminear-ring.*

*Proof*: Straightforward by the very definitions.

It is left for the reader to construct a S-pseudo seminear-rings which has no S-pseudo subseminear-rings.

**DEFINITION 6.1.39**: *Let N be a S-pseudo seminear-ring, if N has no proper S-pseudo subseminear-ring then we say N is a Smarandache pseudo simple seminear-ring (S-pseudo simple seminear-ring).*

**DEFINITION 6.1.40**: *N is said to be Smarandache pseudo seminear-ring (S-pseudo seminear-ring) if N is a near-ring and has a proper subset A of N such that A is a seminear-ring under the operations of N.*

**DEFINITION 6.1.41**: *Let N and $N_1$ be two S-pseudo seminear-rings. h : $N \to N_1$ is a Smarandache pseudo seminear ring homomorphism (S-pseudo seminear-ring homomorphism) if h restricted from A to $A_1$ is a seminear-ring homomorphism.*

**DEFINITION 6.1.42**: *Let N be a S-seminear-ring II. An additive subgroup A of N is called a Smarandache N-subgroup II (S-N-subgroup II) if NA ⊂ N (AN ⊂ A) where NA = {na / n ∈ N and a ∈ A}.*

Thus we see in case of S-seminear-rings we have the concept of S-N-subgroup. We define S-left ideal for S-seminear-ring.

**DEFINITION 6.1.43:** *Let (N, +, •) be a S-seminear-ring. A proper subset I of N is called Smarandache left ideal (S-left ideal II) in N if*

  i.  *(I, +) is a normal subgroup of A ⊂ N where (A, +) is a group.*
  ii. *n ($n_1$ + i) + $n_r$ $n_1$ ∈ I for each i ∈ I, n, $n_1$ ∈ A where $n_r$ denotes the unique inverse of n.*

**DEFINITION 6.1.44:** *A nonempty subset I of (N, +, •); N a S-seminear-ring is called a Smarandache ideal II (S-ideal II) in N if*

  i.  *I is a S-left ideal.*
  ii. *IA ⊂ I (A ⊂ N; A is a near-ring).*

**DEFINITION 6.1.45:** *Let N be a S-seminear ring. N is said to be a Smarandache s-seminear-ring II (S-s-seminear-ring II) if a ∈ Na for each a ∈ A ⊂ N; where A is a near ring.*

Now we define Smarandache fuzzy near-rings.

**DEFINITION 6.1.46:** *Let R be a S-near-ring. A fuzzy set μ : R → [0, 1] is called a Smarandache fuzzy near-ring (S-fuzzy near-ring) related to P ⊂ R where P is a near-field satisfying the following conditions:*

  i.   *μ(x + y) ≥ min{(μ(x), μ(y))}.*
  ii.  *μ(–x) ≥ μ(x).*
  iii. *μ(xy) ≥ min{μ(x), μ(y)}.*
  iv.  *μ(x) = μ($x^{-1}$) for all x, y ∈ P ⊂ R.*

*Thus while defining S-fuzzy near-rings it is also a relative concept related to which near-field it is defined.*

In view of this we have the following definition.

**DEFINITION 6.1.47:** *Let R be a S-near-ring. μ a fuzzy subset of R. We call μ a Smarandache strong fuzzy near-ring (S-strong fuzzy near-ring) if μ is a S-fuzzy near-ring relative to every near-field in R.*

**THEOREM 6.1.6:** *If R is a S-near-ring and μ : R → [0, 1] is a S-strong fuzzy near-rings of R then μ a S-fuzzy near-ring.*

*Proof*: Straightforward by the very definitions, hence left for the reader to prove.

Now we proceed on to define the concept of Smarandache fuzzy near-ring module.

**DEFINITION 6.1.48:** *Let R be a S-near-ring and $\mu$ a S-fuzzy near-ring in R. Let Y be a S-near-ring module over the S-near-ring R and $\sigma$ a fuzzy set in Y. Then $\sigma$ is called the Smarandache fuzzy near-ring module (S-fuzzy near-ring module) in Y if*

    i.   $\sigma(x + y) \geq \min\{\sigma(x), \sigma(y)\}$.
    ii.  $\sigma(\lambda x) \geq \min\{\mu(\lambda), \sigma(x)\}$ *for all* $x, y \in Y$ *and* $\lambda \in P \subset R$.
    iii. $\sigma(0) = 1$.

*It is worthwhile to mention in S-near-rings the condition $\sigma(\lambda x) \geq \min\{\mu(\lambda), \sigma(x)\}$ is replaced by $\sigma(\lambda x) > \sigma(x)$ for all $\lambda \in P \subset R$ and $x \in Y$.*

**DEFINITION 6.1.49:** *Let R be a S-near-ring. A fuzzy subset $\mu$ of R is called the Smarandache fuzzy right (resp. left) R-subgroup of R if*

    i.   $\mu$ *is a S-fuzzy subgroup of* $(P, +)$; $P \subset R$ *a subnear-field of R*.
    ii.  $\mu(xr) > \mu(x)$ *(respectively $\mu(rx) > \mu(x)$) for all* $r, x \in P \subset R$.

*The S-fuzzy right R-subgroup $\mu$ is said to be fuzzy normal whenever $\mu(a) = 1$.*

**DEFINITION 6.1.50:** *Let $\mu$ be a S-fuzzy right (resp. left) R subgroup of a S-near-ring R and let $\mu^+$ be a fuzzy set in R defined by $\mu^+(x) = \mu(x) + 1 - \mu(0)$ for all $x \in P \subset R$ (P a near-field of the S-near-ring R), $\mu^+$ is called the Smarandache fuzzy right (resp. left) R subgroup (S-fuzzy right (resp. left) R-subgroup) of R containing $\mu$.*

The following is in easy consequence of the definition.

**THEOREM 6.1.7:** *Let $\mu$ be a S-fuzzy right (resp. left) S-fuzzy R-subgroup of R satisfying $\mu^+(x) = 0$ for some $x \in P \subset R$ then $\mu(x) = 0$ also (0 is the additive identity of P).*

*Proof*: Easy and straightforward, hence left as an exercise for the reader to prove.

**THEOREM 6.1.8:** *Let $\mu$ and $\nu$ be S-fuzzy right (resp. left) R-subgroups of a S-near-ring R. If $\mu \subset \nu$ and $\mu(0) = \nu(0)$ then $P_\mu \subset P_\nu$ (where $P \subset R$ is a near-field relative to which both $\mu$ and $\nu$ are defined).*

*Proof*: Assume $\mu \subset \nu$ and $\mu(0) = \nu(0)$ (0 is the zero of the near-field P in R). If $x \in P_\mu$ then $\nu(x) \geq \mu(x) = \mu(0) = \nu(0)$, noting $\nu(x) \leq \nu(0)$ for all $x \in P \subset R$ we have $\nu(x) = \nu(0)$ hence $x \in P_\nu$. Hence the claim.

**THEOREM 6.1.9:** *A S-fuzzy right (resp. left) R-subgroup $\mu$ of a S-near-ring R is S-normal if and only if $\mu^+ = \mu$.*

*Proof*: Using the definition the proof is straightforward, hence left for the reader as an exercise.

**THEOREM 6.1.10:** *If $\mu$ is a S-fuzzy right (resp. left) R-subgroup of a S-near-ring R then $(\mu^+)^+ = \mu^+$.*

*Proof*: For any $x \in P \subset R$ relative to which $\mu$ is defined we have $(\mu^+)^+(x) = (\mu^+(x)) + 1 - \mu^+(0) = \mu^+(x)$. Hence the claim.

**Note**: If we change the near-field P from P to say $P_1$ and $\mu$ is not defined relative to $P_1$ then the above result is not true unless otherwise $\mu$ is a S-strong fuzzy subgroup. Further even if $x \in R \setminus P$ the result may not in general be true as all the axioms by $\mu$ are satisfied only in P and not on the whole of R. Several other interesting results in this direction can be defined and proved by any reader. Since when we define $\mu$ a S-fuzzy R-subgroup or $\sigma$ a S-fuzzy right (resp. left) near-ring module we do not demand in any way $\mu$ or $\sigma$ to be completely well defined on the whole of the S-near-ring R what we only expect is that $\mu$ or $\sigma$ is defined on the proper subsets P in R where P is a near-field. In this way all results will be distinct and different from fuzzy near-ring modules and fuzzy R-subgroups.

Only when $\mu : R \to [0, 1]$ happen to be well defined on R and on every proper subset P in R which is a near field we will have the coincidence of fuzzy and S-fuzzy concept.

Now using the results of [28] we give the definition of S-fuzzy congruence of a near-ring module which is slightly different from the definitions of [71].

**DEFINITION 6.1.51:** *Let $\mu$-be a non-empty subset of a S-R-module M of a S-near ring R. Then $\mu$ is said to be Smarandache fuzzy submodule (S-fuzzy module) of M if*

  i. $\mu(x + y) \geq M$ in $\{\mu(x), \mu(y)\}$ for all $x, y \in M$.
  ii. $\mu(-x) = \mu(x) \; \forall x \in M$.
  iii. $\mu(y + x - y) = \mu(y)$ for all $x, y \in M$.
  iv. $\mu((x + y)r - xr) \geq \mu(y)$ for all $x, y \in M$ and $r \in P \subset R$.

*P is the subnear field of the S-near ring relative to which the S-R-module M is defined.*

The following results are direct by the very definition.

**THEOREM 6.1.11:** *Let $\mu$ be a S-fuzzy submodule of a S-R-module M. Then the level subset $\mu_t = \{x \in M \mid \mu(x) \geq t\}$; $t \in \text{Im } \mu$ is a S-sub module of M.*

*(Hint: $\mu$ is a S-fuzzy submodule of a S-R-module M. Then the S-submodule $\mu_t$'s are called the S-level submodule of M).*

**THEOREM 6.1.12:** *Let $\mu$ be a S-fuzzy normal subgroup of M. Then $x + \mu = \mu + x$ if and only if $\mu(x - y) = \mu(0)$ for all x, y in M.*

*Proof*: Here M is a S-R-module of the S-near ring R. The proof is direct by the definition.

**DEFINITION 6.1.52:** *Let M be a S-R-module μ be a S-fuzzy submodule of M. M / μ = {all fuzzy cosets of μ with operations (x + μ) + (y + μ) = (x + y) + μ and ( x + μ) r = xr + μ for all x, y ∈ M and for all r ∈ P; P a proper subset of R which is a near field under the operations of R}, M / μ is defined as the Smarandache fuzzy quotient R-module (S-fuzzy quotient R-module) of M over the S-fuzzy submodule μ.*

**DEFINITION 6.1.53:** *Let M be a S-R-module. A non-empty fuzzy relation on M i.e. a mapping α from M × M → [0, 1] is called the Smarandache fuzzy equivalence relation (S-fuzzy equivalence relation) that is a fuzzy equivalence relation α defined on a S-R-module will be called as a S-fuzzy equivalence relation. A S-fuzzy equivalence relation α on an S-R-module M is called a Smarandache fuzzy congruence (S-fuzzy congruence) if α (a + c, b + d) ≥ Min [α (a, b), α (c, d)] and α (ar, br) ≥ α (a, b) for all a, b, c, d in M and for all r ∈ P ⊂ R (P is the proper subset in the S-near ring R relative to which the S-R-module M is defined). Let α be a S-fuzzy relation on a S-R-module M. For each t ∈ [0, 1] the set $\alpha_t$ = {(a, b) ∈ M × R | α (a, b) ≥ t} is called the Smarandache level relation (S-level relation) of α.*

Now we proceed on to define the concept of Smarandache fuzzy seminear ring for the first time. Before we define Smarandache fuzzy seminear ring we will first define the concept of fuzzy seminear rings and its properties.

**DEFINITION 6.1.54:** *Let T be a seminear ring. N a fuzzy set in T. Then N is called a fuzzy seminear ring in N if*

  i. *N(x + y) ≥ min {N(x), N(y)}.*
  ii. *N(xy) ≥ min {N (x), N(y)}.*

We now define the concept of fuzzy seminear-ring module.

**DEFINITION 6.1.55:** *Let T be a seminear ring. N a fuzzy seminear ring in T. Let Y be a seminear ring module over T and M, a fuzzy set in Y. Then M is called a fuzzy seminear ring module in Y if*

  i. *M (x + y) ≥ min {M(x), M(y)}.*
  ii. *M (λ x) ≥ min {N(λ), M(x)} for all x, y ∈ Y and for all λ ∈ t.*

**DEFINITION 6.1.56:** *Let (T, +, •) be a seminear ring. A fuzzy subset μ in T is called a fuzzy right (resp. left) T-subsemigroup of T if*

  i. *μ is a fuzzy subsemigroup of (T, +).*
  ii. *μ(xr) ≥ μ(x) (resp. μ(rx) ≥ μ(x)) for all x, r ∈ T.*

**THEOREM 6.1.13:** *Let μ be a fuzzy right (resp. left) T-subsemigroup of a seminear ring T then the set $T_μ$ = {x ∈ T / μ(x) = μ (0)} is a right (resp. left) R subsemigroup of T.*

*Proof***:** Straightforward by the definition.

**DEFINITION 6.1.57:** *Let $\mu$ be a nonempty fuzzy subset of a T-semimodule M. Then $\mu$ is said to be a fuzzy semisubmodule of M if*

    i.    $\mu(x + y) \geq M$ *in* $[\mu(x), \mu(y)]$.
    ii.   $\mu((x + y) r + xr)) \geq \mu(y)$

*for all x, y in M and r in T.*

On similar lines we can define fuzzy equivalence relation and fuzzy congruence on G an T-semimodule over a seminear ring T.

*If G is a T-semimodule. A non empty fuzzy relation $\alpha$ on G is a mapping $\alpha : G \times G \to [0, 1]$ is called a fuzzy equivalence relation if*

$$\alpha(x, x) = \sup_{y,z \in G} \alpha(y, z) \text{ for all } x, y, z \text{ in } G.$$

$$\alpha(x, y) = \alpha(y, x) \text{ for all } x, y \text{ in } G$$

$$\alpha(x, y) \geq \sup_{y,z \in G} [Min(\alpha(x\ z), \alpha(z, y)]\ x, y \in G \text{ for all } x, y \text{ in } G.$$

*A fuzzy equivalence relation $\alpha$ on an T-semimodule G is called fuzzy congruence if $\alpha(a + c, b + d) \geq M$ in $[\alpha(a, b), \alpha(c, d)]$ and $\alpha(ar, br) \geq \alpha(a, b)$ for all a, b, c, d in G and for all r in T.*

**THEOREM 6.1.14:** *Let $\alpha$ be a fuzzy congruence on a T-semimodule G and $\mu_\alpha$ be the fuzzy subsemimodule induced by $\alpha$. Let $t \in Im\ \alpha$. Then $(\mu_\alpha)_t = \{ x \in G\ /\ x \equiv 0\ (\alpha_t)\}$ is the subsemimodule induced by the congruence $\alpha_t$.*

*Proof*: Let $a \in G$. Now $a \in (\mu_a)_t \Leftrightarrow (\mu_\alpha)(a) \geq t \Leftrightarrow \alpha(a, 0) \geq t \Leftrightarrow (a, 0)\ t, \alpha_t \Leftrightarrow a \equiv 0(\alpha_t) \Leftrightarrow a \in \{x \in G\ |\ x \equiv 0\ (\alpha_t)\}$. Hence the result.

**DEFINITION 6.1.58:** *Let G be a T-semimodule and $\alpha$ be a fuzzy congruence on G. A fuzzy congruence $\mu$ on G. A fuzzy congruence $\beta$ on G is said to be $\alpha$ - invariant if $\alpha(x, y) \equiv \alpha(\upsilon, \nu)$ implies $\beta(x, y) = \beta(\upsilon, \nu)$ for all $(x, y), (\upsilon, \nu) \in G \times G$.*

Now we proceed on to define Smarandache fuzzy seminear rings.

**DEFINITION 6.1.59:** *Let R be a S-seminear ring. A fuzzy subset $\mu$ on R is said to be a Smarandache fuzzy seminear ring (S-fuzzy seminear-ring) if $\mu : P \to [0, 1]$ is a S-fuzzy semiring where $P \subset R$ and P is a semiring.*

We can define S-fuzzy seminear ring yet in another way.

**DEFINITION 6.1.60:** *Let R be a S-seminear ring. A fuzzy subset $\mu$ of R is said to be a Smarandache fuzzy seminear ring (S-fuzzy seminear ring) if*

$$\mu(x + y) \geq \mu(x) + \mu(y)$$
$$\mu(-x) \geq \mu(x)$$
$$\mu(xy) \geq \min\{\mu(x), \mu(y)\}$$

*for all $x, y \in N \subset R$ where $N$ a proper subset of $R$ is a near ring in $R$.*

The S-fuzzy seminear ring II definition will be distinguished from the other S-fuzzy seminear ring definition by denoting it as S-fuzzy seminear ring II.

In case of S-fuzzy seminear ring II all definitions and results proved in case of fuzzy near rings can be easily transformed into S-fuzzy seminear rings II.

Now we define Smarandache fuzzy seminear ring semimodule II.

**DEFINITION 6.1.61:** *Let $R$ be a S-seminear ring with $P$ a proper subset of $R$ which is a near ring and $N$ a fuzzy near ring in $P$. Let $Y$ be a near ring module over $P$ and $M$ a fuzzy set in $Y$. Then $M$ is called the Smarandache fuzzy seminear ring semimodule (S-fuzzy seminear ring semimodule) in $Y$ if*

  i. *$M(x + y) \geq \min \{M(x), M(y)\}$.*
  ii. *$M(\lambda x) \geq \min \{N(\lambda), M(x)\}$ for all $x, y \in Y$ for all $\lambda \in P$.*
  iii. *$M(0) = 1$.*

The following theorem is direct.

**THEOREM 6.1.15:** *Let $Y$ be a S-seminear ring semimodule over a S-seminear ring $R$ with identity. If $M$ is a S-fuzzy seminear ring semimodule in $Y$ and if $\lambda \in P \subset R$ ($Y$ defined as module over the near ring $P$ contained in $R$) is invertible then $M(\lambda x) = M(x)$ for all $x \in y$.*

*Proof*: As in case of S-near rings using the definition the result is got as a matter of routine.

Several other results can be defined for S-seminear rings and extended for S-fuzzy seminear ring of many other levels which we defined in [130]. The reader is advised to go through the results in [130] and obtain some more interesting results about these structures.

## 6.2 Smarandache Non-associative Fuzzy near-ring

In this section we just recall the definition of S-non-associative seminear-ring and their Smarandache fuzzy analogues. Also we proceed on to define Smarandache non-associative near-rings and their Smarandache fuzzy analogue. We derive some interesting results about them.

**DEFINITION 6.2.1**: *Let $(N, \text{`}+\text{'}, \text{`}\bullet\text{'})$ be a non-empty set endowed with two binary operation '$+$' and '$\bullet$' satisfying the following:*

     i.     *(N, +) is a semigroup.*
    ii.    *(N, •) is a groupoid.*
   iii.   *(a + b) c = a • c + b • c for all a, b, c ∈ N; (N, '+', '•') is called the right seminear-ring which is non-associative.*

*If we replace (c) by a • (b + c) = a • b + a • c for all a, b, c ∈ N. Then (N, '+', '•') is a non-associative left seminear-ring.*

**DEFINITION 6.2.2**: *Let (N, +, •) be a seminear-ring which is not associative. A subset P of N is said to be a subseminear-ring if (P, +, •) is a seminear-ring.*

**DEFINITION 6.2.3**: *Let N be a non-associative seminear-ring. An additive subsemigroup A of N is called the N-subsemigroup (right N-subsemigroup) if NA ⊆ A (AN ⊂ A) where NA = {na / n ∈ N, a ∈ A}.*

**DEFINITION 6.2.4**: *A non-empty subset I of N is called left ideal in N if*

     i.          *(I, +) is a normal subsemigroup of (N, +).*
    ii.         *$n (n_1 + i) + nn_1 \in I$ for each $i \in I$, $n, n_1 \in N$.*

**DEFINITION 6.2.5**: *Let N be a non-associative seminear-ring. A nonempty subset I of N is called an ideal in N if*

   i. *I is a left ideal.*
  ii. *IN ⊂ I.*

**DEFINITION 6.2.6**: *A non-associative seminear-ring N is called left bipotent if Na = $Na^2$ for a in N.*

**DEFINITION 6.2.7**: *A non-associative seminear-ring N is said to be a s-seminear-ring if a ∈ Na for each a in N.*

The following definitions about strictly prime ideals would be of interest when we develop Smarandache notions.

**DEFINITION 6.2.8**: *An ideal P (≠ N) is called strictly prime if for any two N-subsemigroups A and B of N such that AB ⊂ P then A ⊂ P or B ⊂ P.*

**DEFINITION 6.2.9**: *An ideal B of a non-associative (NA for short) seminear-ring N is called strictly essential if B ∩ K ≠ {0} for every non-zero N-subsemigroup K of N.*

**DEFINITION 6.2.10**: *An element x in N is said to be singular if there exists a non-zero strictly essential left ideal A in N such that Ax = {0}.*

Several other analogous results existing in near-rings and seminear-rings can also be defined for non-associative seminear-rings. Now we proceed on to define Smarandache notions.

**DEFINITION 6.2.11**: *Let (N, +, •) be a non associative seminear-ring. N is said be a Smarandache non associative seminear-ring of level I (S-NA seminear-ring I) if*

  i. *(N, +) is a S-semigroup.*
  ii. *(N, •) is a S-groupoid.*
  iii. *(a + b) c = a • c + b • c for all a, b, c ∈ N.*

**DEFINITION 6.2.12**: *Let (N, +, •) be a non-associative seminear-ring. N is said to have a Smarandache subseminear-ring (S-subseminear-ring) P ⊂ N if P is itself a S-seminear-ring.*

**DEFINITION 6.2.13**: *Let (N, +, •) be a non-associative seminear-ring. An additive S-subsemigroup A of N is called the Smarandache left N-subsemigroup (S-N-left subsemigroup) (right N-subsemigroup) if NA ⊂ A (AN ⊂ A) where Na = {na / n ∈ N, a ∈ A}.*

**DEFINITION 6.2.14**: *A non-empty subset I of N is called Smarandache left ideal (S-left ideal) in N if*

  i. *(I, +) is a normal S-subsemigroup of (N, +).*
  ii. *n ($n_1$ + i) + $nn_1$ ∈ I for each i ∈ I, n, $n_1$ ∈ N.*

**DEFINITION 6.2.15**: *Let N be a NA-seminear-ring. A nonempty subset I of N is called a Smarandache ideal (S-ideal) in N if*

  i. *I is a S-left ideal.*
  ii. *IN ⊂ I.*

**DEFINITION 6.2.16**: *A S-NA-seminear-ring N is called Smarandache left bipotent (S-left bipotent) if Na = $Na^2$ for every a in N.*

**DEFINITION 6.2.17**: *Let (N, +, •) be a NA seminear-ring. N is said to be a Smarandache seminear-ring I of type A (S-seminear-ring I of type A) if N has a proper subset P such that (P, +, •) is an associative seminear-ring.*

**DEFINITION 6.2.18**: *Let N be a NA-seminear-ring. N is said to be a Smarandache NA seminear-ring I of type B (S-NA seminear-ring I of type B) if N has a proper subset P where P is a near-ring.*

**DEFINITION 6.2.19**: *Let N and $N_1$ be two S-seminear-rings; we say a map φ from N to $N_1$ is a Smarandache non-associative seminear ring homomorphism (S-NA seminear-ring homomorphism) if*

$$\phi (x + y) = \phi (x) + \phi (y)$$
$$\phi (xy) = \phi (x) \phi (y)$$

*for all x, y ∈ N.*

**DEFINITION 6.2.20**: *Let (N, +, •) be a NA seminear-ring, we say N is a Smarandache NA seminear-ring II (S-NA-seminear-ring II) if N has a proper subset P which is a associative seminear-ring.*

**THEOREM 6.2.1**: *Let N be a S-NA seminear-ring II then N is a S-NA-seminear-ring I of type A.*

*Proof*: Obvious by the very definition. Now using the S-mixed direct product definition of seminear-rings we can extend it to the case of non-associative seminear-rings. This method will help to build a class of S-NA seminear-rings of type II.

**DEFINITION 6.2.21**: *Let N be a S-NA seminear-ring II. An additive S-semigroup A of N is said to be a S-left N subsemigroup (right-N subsemigroup) if $PA \subset A$ and ($AP \subset A$). where P is a proper subset of N and P is a seminear-ring which is associative. $PA = \{pa / p \in P$ and $a \in A\}$.*

**DEFINITION 6.2.22**: *Let N and $N_1$ be any two S-NA-seminear-rings A mapping $\phi$ from N to $N_1$ is called a Smarandache-NA seminear-ring homomorphism (S-NA seminear-ring homomorphism) if $\phi$ maps every $p \in P \subset N$ (p a associative seminear-ring-associative) into a unique element $\phi(p) \in P_1 \subset N_1$ where $P_1$ is an associative seminear-ring such that $\phi(p + p^1) = \phi(p) + \phi(p^1)$ and $\phi(p_1 p_2) = \phi(p_1) \phi(p_2)$ for every $p_1, p_2 \in P \subseteq N$.*

It is important to note that $\phi$ need not be defined on the whole of N it is sufficient if $\phi$ is defined on a subset P of N where P is an associative seminear-ring. Thus if $P \subset N$ and $P_1 \subset N_1$ and $\phi : N \to N_1$ is such that $\phi$ is one to one and on to from P to $P_1$ the two S NA seminear-rings would become isomorphic even if they are not having same number of elements in them.

**DEFINITION 6.2.23**: *Let N be a S-NA-seminear-ring. An additive S-semigroup A of N is said to be a Smarandache left ideal (S-left ideal) of N if it is an ideal of the semigroup (N, +) with the conditions.*

$p_1 (p_2 + a) – p_1 p_2 \in A$ *for each $a \in A$; $p_1, p_2 \in P \subset N$; P a seminear-ring.*

**DEFINITION 6.2.24**: *A subset I of N is called a S-ideal if it is a S-left ideal and $IP \subset I$ where $P \subset N$ is the associative seminear-ring.*

**DEFINITION 6.2.25**: *A S-NA seminear-ring N has Smarandache IFP (S-IFP) (insertion of factor property) if for a, b $\in$ P, ab = 0 implies a.p.B = 0 for all $p \in P \subset N$ since N is non-associative we have to restrict our selves only to the associative substructure to define IFP property as a (nb) $\neq$ (an) b in general for all a, n, b $\in$ N.*

Now we proceed onto define Smarandache fuzzy seminear-rings and its properties.

**DEFINITION 6.2.26**: *Let N be a S-non-associative seminear-ring $\mu : N \to [0, 1]$ be a fuzzy subset of N such that if $\mu : P \to [0, 1]$ is a Smarandache fuzzy seminear-ring (S-fuzzy seminear-ring) where P is a proper subset of N which is a S-seminear-ring* *under the operations of N. Then we call μ a Smarandache fuzzy non-associative seminear-ring (S-fuzzy non-associative seminear-ring) of N.*

**DEFINITION 6.2.27**: *Let R be a S-non-associative seminear-ring. N a fuzzy seminear-ring of R. Let Y be a S-seminear-ring semimodule over R. M a fuzzy set in Y. Then M is called a Smarandache fuzzy non-associative seminear-ring semimodule (S-fuzzy non-associative seminear-ring semimodule) (or just module) in Y if*

    i.    $M(x + y) \geq \min\{M(x), M(y)\}$.
    ii.    $M(\lambda x) \geq \min\{N(\lambda), M(x)\}$.

*for all $x, y \in Y$, $\lambda \in P \subset R$. P is the near-ring over which the module Y is defined.*

All properties enjoyed by S-fuzzy seminear-rings and S-fuzzy near-rings can be transferred in to the definitions of S-fuzzy non-associative seminear-ring with appropriate modifications. As we have given all Smarandache properties about non-associative seminear-rings it would be easy for any reader to define the Smarandache fuzzy analogue for them. Now we proceed on to define some properties about Smarandache non-associative near-ring thereby paving way for the definition and study of Smarandache fuzzy non-associative near-rings.

**DEFINITION 6.2.28**: *The system $N = (N, '+', '\bullet', 0)$ is called a Smarandache right loop half groupoid near-ring (S-right loop half groupoid near-ring) provided.*

    i.    *(N, +, 0) is a Smarandache loop.*

    ii.    *(N, '•') is a half groupoid.*

    iii.    *$(n_1 \bullet n_2) \bullet n_3 = n_1 \bullet (n_2 \bullet n_3)$ for all $n_1, n_2, n_3 \in N$ for which $n_1 \bullet n_2, n_2 \bullet n_3, (n_1 \bullet n_2) \bullet n_3$ and $n_1 \bullet (n_2 \bullet n_3) \in N$.*

    iv.    *$(n_1 + n_2) n_3 = n_1 \bullet n_3 + n_2 \bullet n_3$ for all $n_1, n_2, n_3 \in N$ for which $(n_1 + n_2) \bullet n_3, n_1 \bullet n_3, n_2 \bullet n_3 \in N$. If instead of (iv) in N the identity $n_1 \bullet (n_2 + n_3) = n_1 \bullet n_2 + n_1 \bullet n_3$ is satisfied then we say N is a Smarandache left half groupoid near-ring (S-left half groupoid near-ring).*

*We just say (L, +) is a S-loop if L has a proper subset P such that (P, +) is an additive group.*

**DEFINITION 6.2.29**: *A Smarandache right loop near-ring (S-right loop near-ring) N is a system (N, +, •) of double composition '+' and '•' such that*

    i.    *(N, +) is a S-loop.*
    ii.    *(N, •) is a S-semigroup.*
    iii.    *The multiplication '•' is right distributive over addition that is for all $n_1, n_2, n_3 \in N$ $(n_1 + n_2) \bullet n_3 = n_1 \bullet n_3 + n_2 \bullet n_3$.*

**THEOREM 6.2.2**: *Let N be a S-right loop near-ring then N is a S-right loop half groupoid near-ring.*

*Proof*: Obvious by the very definitions.

**THEOREM 6.2.3**: *Let N be a S-right loop half groupoid near-ring, then N is not in general a S-right loop near-ring.*

*Proof*: Obvious. (N, •) is only a half groupoid so it can never be a S-semigroup.

**DEFINITION 6.2.30**: *A nonempty subset M of a S-loop (N, +, '•', 0) is said to be a Smarandache subloop near-ring (S-subloop near-ring) of N if and only if (M, '+', '•', 0) is a S-loop near-ring.*

**DEFINITION 6.2.31**: *Let N be a non-associative near-ring we say N is a Smarandache quasi non-associative near-ring (S-quasi non-associative near-ring) if N has a proper subset which is a ring under the operations of N.*

**DEFINITION 6.2.32**: *Let N be a non-associative seminear-ring we say N is a Smarandache quasi non-associative seminear-ring (S-quasi non-associative seminear-ring) if N has a proper subset P such that P is a semiring.*

**DEFINITION 6.2.33**: *Let (N, +, •) be a S-quasi seminear-ring. We call a non-empty subset I to be a Smarandache quasi left ideal (S-quasi left ideals) in N if*

  i.    *(I, +) is a S-semigroup.*
  ii.   $n(n^1 + i) + nn^1 \in I$ *for each* $i \in I$ *and* $n, n^1 \in P$; *P a semiring in N.*

*We say I is a S-quasi ideal if* $IP \subset I$.

**DEFINITION 6.2.34**: *Let (N, +, •) be a S-quasi near-ring. We say a non-empty subset I of N to be a Smarandache quasi left ideal (S-quasi left ideal) in N if*

  i.    *(I, +) is a subgroup.*
  ii.   $n(n^1 + i) + nn^1 \in I$ *for each* $i \in I$ *and* $n, n^1 \in R$, $R \subset N$ *and R a ring.*

*We say I is a S-quasi ideal if I is a S-quasi left ideal of N and* $IR \subset I$.

**DEFINITION 6.2.35**: *Let N be a S-quasi near-ring (S-quasi seminear-ring). We say N is Smarandache quasi bipotent (S-quasi bipotent) if* $Pa = Pa^2$ *where* $P \subset N$ *and P is a ring* ($P \subset N$ *and P is a semiring) for every a in N.*

**DEFINITION 6.2.36:** *Let N be a S-quasi near-ring (S-quasi seminear-ring) N is said to be a Smarandache quasi s-near-ring (S-quasi s-near-ring) (S-quasi s-seminear-ring, in short, S-quasi s-seminear-ring) if* $a \in Pa$ *for each a in N where P is a proper subset N which is a ring* ($P \subset N$ *and P is a semiring).*

**DEFINITION 6.2.37:** *Let N be a Smarandache quasi near-ring (S-quasi near-ring) N is said to be Smarandache quasi regular (S-quasi regular) if for each a in N there exists x in P; $P \subset N$, P a ring ($P \subset N$ and P a semiring) such that a = a(xa) = (ax)a.*

**DEFINITION 6.2.38:** *A Smarandache quasi seminear-ring (S-quasi near-ring) N is called Smarandache quasi irreducible (S-quasi irreducible) (Smarandache quasi simple) if it contains only the trivial S-quasi N-subgroups (S quasi N-subsemigroups).*

**DEFINITION 6.2.39:** *Let N be a S-quasi near-ring (S-quasi seminear-ring) an element x is said to be quasi central if xy = yx for all $y \in R$; $R \subset N$ is a ring (or $R \subset N$ and R is a semiring).*

**DEFINITION 6.2.40:** *Let N be a S-quasi near-ring (or S-quasi seminear-ring) N is said to be Smarandache quasi subdirectly irreducible (S-quasi subdirectly irreducible) if the intersection of all nonzero S-quasi ideals of N is nonzero.*

**DEFINITION 6.2.41**: *Let N be a S-quasi near-ring (S-quasi seminear-ring) N is said to have Smarandache quasi insertion of factors property (S-quasi IFP) if a, b $\in$ N, ab = 0 implies arb = 0 where $r \in R$, $R \subset N$ and R is a ring (or $r \in R$, $R \subset N$, R is a semiring).*

**DEFINITION 6.2.42**: *Let N be a non-associative S-near-ring we say N is Smarandache weakly divisible (S-weakly divisble) if for all x, y $\in$ N there exists a $z \in P$; $P \subset N$ where P is an associative ring or P is a near-field such that xz = y or zx = y.*

**DEFINITION 6.2.43**: *Let N be a non-associative S-seminear-ring we say N is Smarandache weakly divisible (S-weakly divisible) if for all x, y $\in$ N there exists $z \in P$, $P \subset N$ where P is an associative seminear-ring such that xz = y or zx = y.*

**DEFINITION 6.2.44**: *Let N be a S-quasi near-ring (S-quasi seminear-ring). We say N is Smarandache quasi weakly divisible (S-quasi weakly divisible) if for all x, y $\in$ N there exists $z \in R$; R a ring $R \subset N$ (R a semiring $R \subset N$) such that xz = y or yz = x.*

**DEFINITION 6.2.45:** *Let P be a seminear pseudo ring (SNP-ring) we say P is a Smarandache SNP-ring I (S-SNP-ring I) if P has a proper subset $T \subset P$ such that T is a seminear-ring. Smarandache SNP-ring II (S-SNP-ring II) if P has a proper subset $R \subset P$ such that R is a near-ring. Smarandache SNP-ring III (S-SNP-ring III) if P has a proper $W \subset P$ such (W, $\oplus$, $\odot$) is a semiring. Thus we have 3 levels of S-SNP rings. A Smarandache SNP subring (S-SNP subring) is defined as a proper subset U of P such that (U, $\oplus$, $\odot$) is a S-SNP-ring.*

**DEFINITION 6.2.46**: *Let (P, $\oplus$, $\odot$) be a SNP-ring. A proper subset I of P is called a Smarandache SNP- ideal (S-SNP-ideal) if*

  i. *for all p, q $\in$ I, p $\oplus$ q $\in$ I.*
  ii. *$0 \in I$.*
  iii. *For all p $\in$ I and r $\in$ P we have p $\odot$ r or r $\odot$ p $\in$ I.*
  iv. *I is a S-SNP-ring.*

**DEFINITION 6.2.47**: *Let (R, ⊕, ⊙) be a quasi SNP-ring. R is said to be a Smarandache quasi SNP-rings (S-quasi SNP-ring) if and only if R is a S-SNP-ring.*

**DEFINITION 6.2.48**: *Let (R, ⊕, ⊙) and ($R_1$, ⊕, ⊙) be any two S-SNP-ring, we say a map $\phi$ is a Smarandache SNP-homomorphism I (II or III) (S-SNP homomorphism I, II or III) if $\phi$: S to $S_1$ where $S \subset R$ and $S_1 \subset R_1$ are seminear-ring (or near-ring or semiring) respectively and $\phi$ is a seminear-ring homomorphism from S to $S_1$ (or near-ring homomorphism from S to $S_1$ or a semiring homomorphism from S to $S_1$). $\phi$ need not be defined on the entire set R or $R^1$ it is sufficient if it is well defined on S to $S_1$.*

*Now we proceed on to define Smarandache right quasi regular element. We just recall that an element $x \in R$, R a ring is said to be right quasi regular if there exist $y \in R$ such that $x \circ y = x + y - xy = 0$ and left quasi regular if there exist $y^1 \in R$ such that $y^1 \circ x = 0 = y^1 + x - y^1 x$.*

*Quasi regular if it is right and left quasi regular simultaneously. We say an element $x \in R$ is Smarandache right quasi regular (S-right quasi regular) if there exist y and $z \in R$ such that $x \circ y = x + y - xy = 0$, $x \circ z = x + z - xz = 0$ but $y \circ z = y + z - yz \neq 0$ and $z \circ y = y + z - zx \neq 0$.*

*Similarly we define Smarandache left quasi regular (S-left quasi regular) and x will be Smarandache quasi regular (S-quasi regular) if it is simultaneously S-right quasi regular and S-left quasi regular.*

Using these definitions the reader can define the Smarandache fuzzy analogue of them. Just for the sake to interest the reader we give a few definitions here.

**DEFINITION 6.2.49:** *Let N be a S-non-associative near-ring. A fuzzy subset $\mu : N \to [0, 1]$ is said to be a Smarandache fuzzy non-associative near-ring (S-fuzzy non-associative near-ring) if $\mu (x + y) \geq \min \{\mu (x), \mu (y)\}$, $\mu (-x) \geq \mu (x)$, $\mu (xy) \geq \min\{\mu (x), \mu (y)\}$ for all $x, y \in P$, P a proper subset of N and P is an associative near-ring or a near-field.*

Thus we in case of S-fuzzy non-associative near-rings do not demand that $\mu$ be defined on the whole of the S-non-associative near-ring N, it is enough if $\mu$ is defined on a proper subset of N which is an associative near-ring or a near-field. Thus having defined the notion of S-na-fuzzy near-ring our interest would be for a given S-na-near-ring N; how many S-na-fuzzy near-rings can be defined.

**DEFINITION 6.2.50**: *Let N be a S-na-near-ring. Let $\mu$ be a fuzzy subset of N i.e. $\mu : N \to [0, 1]$. $\mu$ is said to be a Smarandache strongly fuzzy non-associative near-ring (S-strongly fuzzy non-associative near-ring) if $\mu$ is a S-fuzzy non-associative near-ring if for every proper subset P in N where P is a near-field.*

**THEOREM 6.2.4:** *Let N be a S-non-associative near-ring. If $\mu : N \to [0, 1]$ is a S-strongly fuzzy non-associative near-ring then $\mu$ is a S-fuzzy non-associative near-ring.*

*Proof*: Straightforward by the very definitions and hence left for the reader to prove.

**DEFINITION 6.2.51:** *Let R be a S-non-associative near-ring and N a S-fuzzy-na-near-ring in R. Let Y be a S-non-associative near-ring module over R and M a fuzzy set in Y. Then M is called a S-fuzzy na-near-ring module in Y if*

  i.   $M(x + y) \geq Min\{M(x), M(y)\}$.
  ii.  $M(\lambda x) \geq Min\{N(\lambda), M(x)\}$.
  iii. $M(0) = 1$

*for all $x, y \in Y$ and $\lambda \in P \subset R$ where P is a proper subset of R and P is an associative near-ring or a near-field.*

**THEOREM 6.2.5:** *Let $\{M_i / i \in I\}$ be a family of S-fuzzy non-associative near-ring modules relative to a fixed associative subnear-ring P in R, all the $M_i$'s are defined in Y. Then $\bigcap_{i \in I} M_i$ is a S-fuzzy non-associative near-ring module in Y.*

*Proof*: Using the fact if $M = \bigcap_{i \in I} M_i$, then we have $\lambda \in P \subset R$ (P relative to which all $M_i$'s are defined) and for all $x, y \in Y$

$$
\begin{aligned}
M(x + y) &= \inf_{i \in I} M_i(x+y) \\
&\geq \inf_{i \in I} \{min\{M_i(x), M_i(y)\}\} \\
&= \min\{\inf M_i(x), \inf M_i(y)\} \\
&= \min\{M(x), M(y)\}
\end{aligned}
$$

On similar lines we have $M(\lambda x) = Min\{N(\lambda), M(x)\}$. Hence the theorem.

Several other properties can be defined and studied in case of S-fuzzy na-near-rings.

## 6.3 Smarandache Fuzzy Binear-rings

In this section we introduce just the notions of Smarandache fuzzy binear-rings, Smarandache fuzzy biseminear-rings, Smarandache non-associative binear-rings and Smarandache non-associative biseminear-rings and their Smarandache fuzzy analogues. Several results in this direction are given.

**DEFINITION 6.3.1:** *Let $(N, +, \bullet)$ be a binear-ring. $(N = N_1 \cup N_2)$. We say N is a Smarandache binear-ring (S-binear-ring) if N contains a proper subset P such that P under the operations '+' and '•' is a binear field; i.e., $(P = P_1 \cup P_2, +, \bullet)$ is a binear field.*

**DEFINITION 6.3.2:** *Let $(N, +, \bullet)$ be a binear-ring. A proper subset P of N is said to be a Smarandache subbinear-ring (S-subbinear-ring) if P itself is a S-binear-ring.*

**DEFINITION 6.3.3:** *Let P be a S-bisemigroup with 0 and let N be a S-binear-ring. A map $\mu: N \times Y \to Y$ where Y is a proper subset of P which is a bigroup under the operations of P, (P, $\mu$) is called the Smarandache N-bigroup (S-N-bigroup) if for all $y \in Y$ and for all $n, n_l \in N$ we have $(n + n_l)y = ny + n_l y$ and $(nn_l) y = n(n_l y)$. $S(N^P)$ stands for the S-N-bigroups.*

**DEFINITION 6.3.4:** *A S-bisemigroup M of the binear-ring N is called a Smarandache quasi sub-binear-ring (S-quasi sub-binear-ring) of N if $X \subset M$ where X is a bisubgroup of M which is such that $XX \subset X$.*

**DEFINITION 6.3.5:** *A S-sub-bisemigroup Y of $S(N^P)$ with $NY \subset Y$ is said to be a Smarandache N-sub-bigroup (S-N-sub-bigroup) of P.*

**DEFINITION 6.3.6:** *Let N and N' be two S-binear-rings P and P' be S-N-sub-bigroups*

   i.                   *$h : N \to N'$ is called a Smarandache binear-ring homomorphism (S-binear-ring homomorphism) if for all $m, n \in M$ we have $h(m + n) = h(m) + h(n)$, $h(mn) = h(m)h(n)$ where $h(m), h(n) \in M'$ (M' is a proper subset of N' which is a binearfield). It is to be noted that h need not even be defined on whole of N.*

   ii.                  *$h : P \to P'$ is called the Smarandache N-sub-bigroup homomorphism (S-N-sub-bigroup homomorphism) if for all p, q in S (S the proper subset of P which is a S-N-sub-bigroup of the S-bisemigroup P) and for all $m \in M \subset N$ (M a nearfield of N); $h(p + q) = h(p) + h(q)$ and $h(mp) = mh(p)$, $h(p), h(q)$ and $mh(p) \in S'$ (S' is a proper subset of P' which is S-N'-sub-bigroup of S-bisemigroup P').*

*Here also we do not demand h to be defined on whole of P.*

**DEFINITION 6.3.7:** *Let $(N, +, \bullet)$ be a S-binear-ring. A normal sub-bigroup I of (N, +) is called a Smarandache bi-ideal (S-bi-ideal) of N related to X, where X is a binearfield contained in N if ($X = X_1 \cup X_2$, $X_1$ and $X_2$ are near fields).*

   i.              $I_1 X_1 \subset I_1$; $I_2 X_2 \subset I_2$ *where* $I = I_1 \cup I_2$.
   ii.             $\forall x_i, y_i \in X_i$ *and for all* $k_i \in I_i$, $x_i (y_i + k_i) - x_i y_i \in I_i$; $i = 1, 2$.

**DEFINITION 6.3.8:** *A proper subset S of P is called a S-bi-ideal of $S(N^P)$ related to Y if*

   i.      *S is a S-normal sub-bigroup of the S-bisemigroup.*
   ii.     *For all $s_1 \in S$ and $s \in Y$ and for all $m \in M$ (M the binear field of N), $n(s + s_1) - ns \in S$.*

*A S-binear-ring is S-bi-ideal if it has no S-bi-ideals. $S(N^P)$ is called Smarandache N-bisimple (S-N-bisimple) if it has no S-normal sub-bigroups expect 0 and P.*

**DEFINITION 6.3.9:** *A S-sub-binear-ring M of a binear-ring N is called Smarandache bi-invariant (S-bi-invariant) related to the binearfield X in N if $MX \subset M$ and $XM \subset M$ where X is a S-binear field of N. Thus in case of Smarandache bi-invariance (S-bi-* *invariance) it is only a relative concept as a S-sub-binear-ring M may not be invariant related to every binear field in the binear-ring N.*

**DEFINITION 6.3.10:** *The binear-ring N is said to fulfill the Smarandache left ore condition (S-left ore condition) with respect to a given S-sub-bisemigroup P of (N, •) if for (s, n) ∈ S × N there exists n • $s_1$ = s • $n_1$ ($s_1$ • n = $n_1$ • s).*

**DEFINITION 6.3.11:** *Let S(BV) denote the set of all Smarandache binear-rings. A S-binear-ring $F_X$ ∈ S(BV) is called a Smarandache free binear-ring (S-free binear-ring) in BV over X if there exists f : X → $F_X$ and for all N ∈ BV and for all g : X → N there exists a S-binear-ring homomorphism h ∈ S(Hom $F_X$ , N) such that h o f = g.*

**DEFINITION 6.3.12:** *A finite sequence N = $N_0$ ⊃ $N_1$ ⊃ $N_2$ ⊃ ... ⊃ $N_t$ = {0} of S-sub-binear-rings $N_i$ of N is called a Smarandache normal sequence (S-normal sequence) of N if and only if for all i ∈ {1, 2, ... , n}, $N_i$ is an S-bi-ideal of $N_{i-1}$. In the special case that all $N_i$ is an S-bi-ideal of N then we call the normal sequence a Smarandache bi-invariant sequence (S-bi-invariant sequence) and t is called the Smarandache length (S-length) of the sequence $N_{i-1}$ /$N_i$ are called the Smarandache bi-factors (S-bi-factors) of the sequence as*

$$N_{i-1} = N_{i-1}^1 \cup N_{i-1}^2, N_i = N_i^1 \cup N_i^2.$$

*So*

$$N_{i-1}/N_i = \left(N_{i-1}^1/N_i^1\right) \cup \left(N_{i-1}^2/N_i^2\right).$$

**DEFINITION 6.3.13:** *Let P be a S-bi-ideal of the binear-ring N. P is called Smarandache prime bi-ideal (S-prime bi-ideal) if for all S-ideals I and J of N, IJ ⊂ P implies I ⊂ P or J ⊂ P.*

**DEFINITION 6.3.14:** *Let P be a S-left bi-ideal of a binear-ring N. P is called Smarandache bimodular (S-bimodular) if and only if there exists a S-idempotent e ∈ N and for all n ∈ N; n – ne ∈ P.*

**DEFINITION 6.3.15:** *Let N be a binear-ring (N = $N_1$ ∪ $N_2$). z ∈ N is called Smarandache quasi regular (S-quasi regular) if z ∈ S($L_Z$) where S($L_Z$) = {$n_1$ – $n_1$z | $n_1$ ∈ $N_1$} ∪ {$n_2$ – $n_2$z | $n_2$ ∈ $N_2$}. An S-bi-ideal P ⊂ N is called S-quasi regular if and only if for all s ∈ Ps is S-quasi regular. An S-bi-ideal P is a Smarandache bi-principal bi-ideal (S-bi-principal bi-ideal) if P = $P_1$ ∪ $P_2$ and each of $P_1$ and $P_2$ are S-principal ideals of $N_1$ and $N_2$ where N = $N_1$ ∪ $N_2$ is a binear-ring. A binear-ring N is Smarandache biregular (S-biregular) if there exists some set E of central S-idempotents with*

   i.   *For all e ∈ N, Ne is a S-bi-ideal of N.*
   ii.  *For all n ∈ N there exists e ∈ E; Ne = (n).*
   iii. *For all e, f ∈ E, e + f = f + e.*
   iv.  *For all e, f ∈ E, ef and e + f – ef ∈ E.*

Now we define Smarandache binear-ring of level II.

**DEFINITION 6.3.16:** *Let N be a binear-ring. We say N is a Smarandache binear-ring of level II (S-binear-ring of level II) if N contains a proper subset P such that P is a biring under the operations of N.*

**DEFINITION 6.3.17:** *Let (N, +, •) be a biseminear-ring. We say N is a Smarandache biseminear-ring (S-biseminear-ring) if N contains a proper subset P, such that P is a binear-ring under the operations of N. We call this S-biseminear-ring a S-biseminear-ring of level I.*

Now we define Smarandache biseminear-ring of level II.

**DEFINITION 6.3.18:** *Let (N, +, •) be a seminear-ring. We call N a Smarandache biseminear-ring of level II (S-biseminear-ring of level II) if N contains a proper subset P, such that P under the operations of N is a bisemiring.*

We can obtain classes of S-biseminear-ring II by defining a Smarandache mixed direct product of a bisemirings, i.e., $W = V \times Z$ where $V = V_1 \cup V_2$ is a biseminear-ring and $Z = Z_1 \cup Z_2$ is a bisemiring.

**DEFINITION 6.3.19:** *Let N be a S-binear-ring. N is said to fulfill the Smarandache bi-insertions factors property (S-bi-insertions factors property) if for all a, b ∈ N = $N_1 \cup N_2$ we have a • b = 0 implies anb = 0 for all n ∈ P, P ⊂ N is a P-binear field. (if a, b ∈ $N_1$ then for all n ∈ $P_2$ we have ab = 0 implies anb = 0 and if a, b ∈ $N_2$ with ab = 0, then ∀n ∈ $P_2$ we have anb = 0).*

*We say the S-binear-ring satisfies the Smarandache strong bi-insertion factors property (S-strong bi-insertion factors property) if for all a, b ∈ $N_i$ (N = $N_1 \cup N_2$), ab ∈ $I_i$ (I = $I_1 \cup I_2$ a S-bi-ideal of N) implies anb ∈ $I_i$ where n ∈ $P_i$ (P ⊂ N, P a binear-ring with P = $P_1 \cup P_2$); i = 1, 2.*

**DEFINITION 6.3.20:** *Let p be a prime. A S-binear-ring N is called a Smarandache p-binear-ring (S-p-binear-ring) provided for all x ∈ $P_i$, $x^p = x$ and px = 0, where P = $P_1 \cup P_2$, P ⊂ N; P is a bisemifield.*

**DEFINITION 6.3.21:** *Let N be a binear-ring. A S-right bi-ideal I of N is called Smarandache right biquasi reflexive (S-right biquasi reflexive) if whenever A and B are S-bi-ideals of N with AB ⊂ I, then b(b' + a ) – bb' ∈ I for all a ∈ A and for all b, b' ∈ B. A binear-ring N is Smarandache strongly bi-subcommutative (S-strongly bi-commutative) if every S-right bi-ideal of it is S-right biquasi reflexive. Let N be a binear-ring P a S-subnormal bi-subgroup of (N, +). P is called Smarandache quasi bi-ideal (S-quasi bi-ideal) of N if PN ⊂ P and NP ⊂ P where by NP we mean elements of the from {n(n' + s) – nn' for all s ∈ P and for all n, n' ∈ N} = NP.*

**DEFINITION 6.3.22:** *Let (R, +, •) be a S-biring and M a S-right R-bimodule. Let W = R × M and define (α, s) ⊙ (β, t) = (αβ, sβ + t). Then (W, +, ⊙) is a S-left binear-ring the abstract Smarandache affine binear-ring (S-affine binear-ring) inducted by R and M. All other notions introduced for bi-ideals can also be defined and extended in case of bi-R-modules.*

**DEFINITION 6.3.23:** *Let N and N' be any two S-biseminear-rings. We say a map $\phi$ from N to N' is a Smarandache biseminear-ring homomorphism (S-biseminear-ring homomorphism) from A to A' where $A \subset N$ and $A' \subset N'$ is a binear-ring if $\phi(x + y) = \phi(x)\phi(y)$, $\phi(xy) = \phi(x)\phi(y)$ where $x, y \in A_i$; $i = 1, 2$; and $\phi(x), \phi(y) \in A'_1$; $i = 1, 2$.*

*Thus for a Smarandache biseminear-ring homomorphism we do not require $\phi$ to be defined on whole of N, it is sufficient if it is defined on a proper near-ring A which is a subset of N. We say the biseminear-ring N is a Smarandache strict biseminear-ring if the subset $A \subset N$ where A is a strict bisemiring.*

*So we can define Smarandache biseminear-ring homomorphism II between two S-biseminear-rings N and N' as $\phi : N \to N'$ is a S-biseminear-ring homomorphism, $\phi : A \to A'$ where A and A' are bisemirings contained in N and N' respectively is a bisemiring homomorphism. Clearly as in case of S-homomorphism I. $\phi$ need not in general be defined on the whole of N.*

**DEFINITION 6.3.24:** *Let N be a Smarandache pseudo biseminear-ring (S-pseudo biseminear-ring) if N is a binear-ring and has proper subset A of N, which is a biseminear-ring under the operations of N.*

**DEFINITION 6.3.25:** *Let N be a S-binear-ring an element $a \in P = P_1 \cup P_2$ (P a binear-field contained in N) is called Smarandache binormal element (S-binormal element) of N, if $aN_1 = N_1a$, $aN_2 = N_2a$, where $N = N_1 \cup N_2$; if $aN_1 = N_1a$ and $aN_2 = N_2a$ for every $a \in P = P_1 \cup P_2$ then N is called the Smarandache binormal near-ring (S-binormal near-ring). Let B(N) denote the set of all S-binormal bielements of N. N is called a S-binormal near-ring if and only if B(N) = P where $P \subset N$ and P is a binear-field.*

**DEFINITION 6.3.26:** *Let $(N, +, \bullet)$ be a na-binear-ring. We say N is a Smarandache non-associative binear-ring (S-non-associative binear-ring) if*

  i. *(N, +) is a S-bigroup.*
  ii. *$(N, \bullet)$ is a S-bigroupoid.*
  iii. *$(a + b) \bullet c = a \bullet c + b \bullet c$ for all $a, b, c \in N$.*

**DEFINITION 6.3.27:** *Let $(N, +, \bullet)$ be a non-empty set, we say N is a Smarandache non-associative biseminear-ring (S-non-associative biseminear-ring) if*

  i. *(N, +) is a S-bisemigroup.*
  ii. *$(N, \bullet)$ is a S-bigroupoid.*
  iii. *$(a + b) \bullet c = a \bullet c + b \bullet c$ for all $a, b, c \in N$.*

**DEFINITION 6.3.28:** *Let $(N, +, \bullet)$ be a na-binear-ring if N has a proper subset P, where P is an associative binear-ring then we call N a Smarandache non-associative binear-ring of level II (S-non-associative binear-ring of level II).*

**DEFINITION 6.3.29:** *Let (N, +, •) be a na-biseminear-ring, if N has a proper subset P such that P is an associative biseminear-ring then we call N a Smarandache non-associative biseminear-ring of level II (S-non-associative biseminear-ring of level II).*

**DEFINITION 6.3.30:** *Let (N, +, •) be a non-empty set we say N is a biquasi ring if the following conditions are true.*

   i.   *$N = N_1 \cup N_2$ ($N_1$ and $N_2$ are proper subsets of N).*
   ii.  *$(N_1, +, •)$ is a ring.*
   iii. *$(N_2, +, •)$ is a near-ring.*

**DEFINITION 6.3.31:** *Let (P, +, •) be a non-empty set with two binary operations '+' and '•' P is said to be a biquasi semiring if*

   i.   *$P = P_1 \cup P_2$ are proper subsets.*
   ii.  *$(P_1, +, •)$ is a ring.*
   iii. *$(P_2, +, •)$ is a semiring.*

**DEFINITION 6.3.32:** *Let (T, +, •) be a non-empty set with two binary operation. We call T a biquasi near-ring if*

   i.   *$T = T_1 \cup T_2$ where $T_1$ and $T_2$ are proper subsets of T.*
   ii.  *$(T_1, +, •)$ is a semiring.*
   iii. *$(T_2, +, •)$ is a near-ring.*

**DEFINITION 6.3.33:** *Let ( N, +, •) be a biquasi ring. A proper subset $P \subset N$, $P = P_1 \cup P_2$ is said to be a sub-biquasi ring if (P, +, •) is itself a biquasi ring.*

**DEFINITION 6.3.34:** *Let (N, +, •) be a biquasi semiring. A proper subset $P \subset N$ where $P = P_1 \cup P_2$ is a sub-bi-quasi semiring if P itself is a biquasi semiring.*

**DEFINITION 6.3.35:** *Let (N, +, •) be a biquasi ring we call N a non-associative biquasi ring if $(N_1, +, •)$ which is a ring is a non-associative ring and $(N_2, +, •)$ is a near-ring.*

**DEFINITION 6.3.36:** *Let (N, +, •) be a biquasi semiring we say N a non-associative biquasi semiring if $(N_2, +, •)$ which is a semiring is a non-associative semiring.*

**DEFINITION 6.3.37:** *Let (N, +, •) be a biquasi near-ring. We say N is a non-associative biquasi near-ring where $N = Z^0 \cup ZL_{11}(3)$ where $Z^0$ is a semiring and $ZL_{11}(3)$ is the loop near-ring which is non-associative, hence N is na biquasi near-ring.*

Now we proceed for their Smarandache analogue.

**DEFINITION 6.3.38:** *Let (N, +, •) where $N = N_1 \cup N_2$ be a biquasi ring. We call N a Smarandache biquasi ring (S-biquasi ring) if $(N_1, +, •)$ is a S-ring.*

**DEFINITION 6.3.39:** *Let (N, +, •) be a biquasi semiring. N is said to be a Smarandache biquasi semiring (S-biquasi semiring) if the semiring ($N_2$, +, •) is a S-semiring.*

**DEFINITION 6.3.40:** *Let (T, +, •) be a non-empty set which is a biquasi near-ring. We call T a Smarandache biquasi near-ring (S-biquasi near-ring) if the near-ring ($T_2$, +, •) is a S-near-ring.*

We can define the concept of S-sub-biquasi ring, S-sub-biquasi semiring and S-sub-biquasi near-rings and the corresponding ideas as follows.

**DEFINITION 6.3.41:** *Let (N, +, •) be a biquasi ring a proper subset, (P, +, •) of N is said to be a Smarandache sub-biquasi ring (S-sub-biquasi ring) if $P_1 \cup P_2$ is the union of proper subsets; ($P_1$, +, •) is a S-subring of $N_1$ and $P_2$ is a subnear-ring of $P_2$.*

**DEFINITION 6.3.42:** *Let (P, +, •) be a biquasi semiring. A proper subset (X, +, •) of P is said to be a Smarandache sub-biquasi semiring (S-biquasi near-ring) if ($X_2$, +, •) is a S-semiring where $X = X_1 \cup X_2$.*

**DEFINITION 6.3.43:** *Let (S, +, •) be a biquasi near-ring. A proper subset (L, +, •) of S is said to be a Smarandache sub-bi quasi near-ring (S-sub-bi quasi near-ring) if ($L_2$, +, •) is a S-near-ring where $L = L_1 \cup L_2$ and ($L_1$, +, •) is a semiring.*

Now we define Smarandache bipseudo structures.

**DEFINITION 6.3.44:** *(Y, +, •) is a Smarandache bi pseudo ring (S-bi pseudo ring) if the following conditions are true.*

    i.    *$Y = Y_1 \cup Y_2$ are proper subsets of Y.*
    ii.    *($Y_1$, +, •) is a S-ring.*
    iii.    *($Y_2$, +, •) is a S-near-ring.*

**DEFINITION 6.3.45:** *Let (P, +, •) be a non-empty set. P is said to be a Smarandache bi pseudo semiring (S-bi pseudo semiring) if the following conditions are true.*

    i.    *$P = P_1 \cup P_2$ where $P_1$ and $P_2$ are proper subsets of P.*
    ii.    *($P_1$, +, •) is a S-ring.*
    iii.    *($P_2$, +, •) is a S-semiring.*

**DEFINITION 6.3.46:** *Let (X, +, •) be a non-empty set, X is said to be a Smarandache bi pseudo near-ring (S-bi pseudo near-ring) if the following conditions are true.*

    i.    *$X = X_1 \cup X_2$ is the union of two proper subsets.*
    ii.    *($X_1$, +, •) is a S-semiring.*
    iii.    *($X_2$, +, •) is a S-near-ring.*

Now we proceed onto define Smarandache fuzzy binear-rings.

**DEFINITION 6.3.47:** *Let $(N = N_1 \cup N_2, +, \bullet)$ be a S-binear-ring A fuzzy subset $\mu : N \to [0, 1]$ is said to be a Smarandache fuzzy binear-ring (S-fuzzy binear-ring) if $\mu = \mu_1 \cup \mu_2 : N_1 \cup N_2 \to [0, 1]$ satisfies the following conditions.*

    i.                     *$\mu_1 : N_1 \to [0, 1]$ is a S-fuzzy near-ring and*
    ii.                   *$\mu_2 : N_2 \to [0, 1]$ is a S-fuzzy ring or S-fuzzy near-ring.*

**DEFINITION 6.3.48:** *Let $(N, +, \bullet)$ be a S-binear-ring. A fuzzy subset $\mu : N \to [0, 1]$ is a S-fuzzy binear-ring. A fuzzy subset $\sigma \subseteq \mu$ is said to be a S-fuzzy sub-binear-ring if $\sigma(x) \leq \mu(x)$ for all $x \in P = P_1 \cup P_2$, $P \subset N$ is a binear-field.*

**DEFINITION 6.3.49:** *Let $(N = N_1 \cup N_2, +, \bullet)$ be a S-binear-ring. A fuzzy subset $\mu = \mu_1 \cup \mu_2$ in N is called a S-fuzzy right (resp. left) N-sub-bigroup of N if*

    i.    *$\mu$ is a S-fuzzy sub-bigroup of $(N, +)$.*
    ii.   *$\mu(xr) \geq \mu(x)$ [resp. $\mu(rx) > \mu(x)$] for all $r, x \in P = P_1 \cup P_2 \subset N_1 \cup N_2$, the proper subset of N which is a binear-field.*

**THEOREM 6.3.1:** *If $\mu$ is a S-fuzzy right (respectively left) N-sub-bigroup of a S-binear-ring N then the set $N_\mu = \{x \in P \subset N / \mu(x) = \mu(0)\}$ is a S-right (resp. left) N-sub-bigroup of N.*

*Proof:* Proof is a matter of routine and left for the reader as an exercise.

Now we proceed on to define the notion of S-fuzzy non-associative binear-ring. We just give the definitions and expect the reader to develop more properties about them as a matter of routine.

**DEFINITION 6.3.50:** *Let $(N = N_1 \cup N_2, +, \bullet)$ be a S-non-associative binear-ring. $\mu$ a fuzzy subset of N. $\mu = \mu_1 \cup \mu_2$ is said to be a Smarandache fuzzy non-associative binear-ring (S-fuzzy non-associative binear-ring) if $\mu_1 : N_1 \to [0, 1]$ is a S-fuzzy near-ring relative to a near-field $P_1 \subset N_1$ and $\mu_2 : N_2 \to [0, 1]$ is a S-fuzzy near-ring relative to a near-field $P_2 \subset N_2$. Then $\mu = \mu_1 \cup \mu_2 : N_1 \cup N_2 \to [0, 1]$ is a S-fuzzy non-associative binear-ring relative to the binear-field $P_1 \cup P_2$.*

**DEFINITION 6.3.51:** *Let $(N = N_1 \cup N_2, +, \bullet)$ be a S-non-associative binear-ring. A fuzzy set $\mu$ in N is called a Smarandache fuzzy non-associative right or left N sub-bigroup (S-fuzzy non-associative right or left N-sub-bigroup) of N relative to a bifield $P = P_1 \cup P_2 \subset N$ if*

    i.   *$\mu$ is a S-fuzzy sub-bigroup of $(N, +)$.*
    ii.   *$\mu(xr) \geq \mu(x)$ (resp. $\mu(rx) > \mu(x)$) for all $x, r \in P$.*

*Thus depending on every bifield we will have different S-non-associative fuzzy sub-bigroups.*

Several interesting results in this direction can be derived. Now we proceed onto define Smarandache fuzzy biseminear-rings.

**DEFINITION 6.3.52:** *Let* $(S = S_1 \cup S_2, +, \bullet)$ *be a Smarandache associative biseminear-ring (S-associative biseminear-ring). A fuzzy subset* $\mu : S \to [0, 1]$ *is said to be a Smarandache fuzzy associative biseminear-ring relative to a proper subset (S-fuzzy associative biseminear-ring relative to a proper subset) P in S where P is a biseminear-field satisfying the following conditions:*

  i.  $\mu(x + y) \geq \min\{\mu(x), \mu(y)\}$.
  ii. $\mu(xy) \geq \min\{\mu(x), \mu(y)\}$

*for all x, y in P.*

*If* $\mu$ *is a fuzzy subset such that* $\mu$ *is a S-fuzzy associative biseminear-ring relative to every proper subset P which is a biseminear-field then we call* $\mu$ *a Smarandache strong fuzzy associative biseminear-ring (S-strong fuzzy associative biseminear-ring).*

**THEOREM 6.3.2:** *Let* $(S = S_1 \cup S_2, +, \bullet)$ *be an S-biseminear-ring.* $\mu : S \to [0, 1]$ *be a fuzzy subset of S which is a S-strong fuzzy biseminear-ring of S then S is a S-fuzzy biseminear-ring of S.*

*Proof*: Direct by the very definitions, hence left as an exercise for the reader to prove.

**DEFINITION 6.3.53:** *Let* $(S = S_1 \cup S_2, +, \bullet)$ *be a S-non-associative biseminear-ring. A fuzzy subset* $\mu$ *of S is said to be a Smarandache fuzzy non-associative biseminear-ring if the following conditions are true:*

  i. $\mu(x + y) \geq \min\{\mu(x), \mu(y)\}$.
  ii. $\mu(xy) \geq \min\{\mu(x), \mu(y)\}$

*for all* $x, y \in P$.

*P a proper subset of S which is a biseminear-field. If* $\mu$ *is a fuzzy subset of S-non-associative biseminear-ring S and if* $\mu$ *is a S-fuzzy non-associative biseminear-ring for every proper subset P of S which is a biseminear-field then we call* $\mu$ *a Smarandache strong fuzzy non-associative biseminear-ring (S-strong fuzzy non-associative biseminear-ring).*

**DEFINITION 6.3.54:** *Let* $(S, +, \bullet)$ *be a S-non-associative biseminear-ring. A fuzzy set* $\mu$ *of S is called a Smarandache fuzzy right (resp. left) S-sub-bisemigroup (S-fuzzy right (resp. left) S-sub-bisemigroup) of S if*

  i.  $\mu$ *is a S-fuzzy sub-bisemigroup of* $(S, +)$.
  ii. $\mu(xr) > \mu(x)$ *(resp.* $\mu(rx) > \mu(x)$*) for all* $r, x \in P \subset S$

*where P is the proper subset which is a biseminear-field relative to which* $\mu$ *is defined.*

Define S-fuzzy normal S-sub-bisemigroup of S and define $\mu^+$ and obtain the results about $\mu^+$ when S is S-biseminear-ring. Also introduce the concept of S-fuzzy congruence of a S-biseminear-ring bimodule.

## 6.4 Problems

This section which is solely devoted to problems gives 25 problems and proposes the reader to define some more new notions about Smarandache fuzzy structures. The reader is expected to solve these problems to get an insight into the subject.

**Problem 6.4.1:** If it true if Y is a S-near ring module over a S-fuzzy near ring $\mu$ in R. Then $\sigma$ is a S-fuzzy near ring module in Y of and only if $\sigma(sx + ry) \geq \min\{\min \mu(s), \sigma(x)\}, \min\{u(r), \sigma(y)\}\}$ for all s, r $\in \mu$ and all x, y $\in$ Y.

**Problem 6.4.2:** Let Y be a S-near-ring module over a S-near ring R with identity. Prove if $\sigma$ is a S-fuzzy near ring module in Y and if $\lambda \in P \subset R$ is invertible then $\sigma(\lambda x) = \sigma(x)$ for all x $\in$ Y.

**Problem 6.4.3:** Let $\{\sigma_i \mid i \in I\}$ be a family of S-fuzzy near-ring modules in Y. Then prove $\bigcap_{i \in I} \sigma_i$ is a S-fuzzy near ring module in Y.

**Problem 6.4.4:** Let Y and W be S-near ring modules over a S-fuzzy near ring, $\mu$ in a S-near ring R and $\theta$ a S-homomorphism of Y into W. Let $\sigma$ be a S-fuzzy near ring module in W. Will the inverse image $\theta^{-1}$ (M) of M a S-fuzzy near ring module in Y?

**Problem 6.4.5.** Let Y and W be S-near ring modules over a S-fuzzy near ring N in a S-near ring R and $\theta$ be a S-near ring homomorphism of Y into W. Let W be a S-fuzzy near ring module in Y that has the sup property. Then prove $\theta(m)$ of m is a S-fuzzy near ring module in W.

**Problem 6.4.6:** Let $\mu$ be a S-fuzzy right (resp. left) R-subgroup of a S-near ring R. ($\mu$ defined relative to near field P in R). If there exists a S-fuzzy right (resp. left) R-subgroup $\nu$ of R ($\nu$ defined only relative this P) satisfying $\nu^+ \subset \mu$ then prove $\mu$ is S-normal.

**Problem 6.4.7:** Let $\alpha$ be a S-fuzzy relation on an S-R-module M. Then prove $\alpha$ is a S-fuzzy congruence on M if and only if $\alpha_i$ is a S-congruence on m for each t $\in$ Im $\alpha$.

**Problem 6.4.8:** Let $\alpha$ be a S-fuzzy congruence on an S-R-module M and $\mu_a$ be the S-fuzzy submodule induced by $\alpha$. Let t $\in$ Im$\alpha$. Prove $(\mu_\alpha)_t = \{x \in M \mid x \equiv 0\ (\alpha_t)\}$ is the S-submodule induced by the congruence $\alpha_t$.

**Problem 6.4.9:** Define S-fuzzy biregular near rings illustrate it with examples. Obtain some nice results on them.

**Problem 6.4.10:** Define Smarandache fuzzy S-prime radical of a near ring. Give an example.

**Problem 6.4.11:** Define Smarandache fuzzy free near ring can there exists any relation between S-fuzzy free near ring and S-free near ring?

**Problem 6.4.12:** Define Smarandache fuzzy quasi regular near ring. Find some interesting results about them.

**Problem 6.4.13:** In case of S-non-associative semi near rings define S-fuzzy quasi-regular semi near ring.

**Problem 6.4.14:** Give an example of a S-strongly fuzzy non associative near ring.

**Problem 6.4.15:** Give an example of a S-fuzzy non associative near ring which is not a S-strongly fuzzy non associative near ring.

**Problem 6.4.16:** If $\mu$ is a S-fuzzy right (resp. left) N-subbigroup of N satisfying $\mu^+(x) = 0$ for some $x \in N$ then $\mu(x) = 0$ Prove.

**Problem 6.4.17:** Let $\mu$ and $\nu$ be a S-fuzzy right (resp. left) N-subbigroups of a S-binear ring N. If $\mu \subset \nu$ and $\mu(0) = \nu(0)$ then prove $P_\mu \subset P_\nu$ (Here P is the binear field contained in R relative to which $\mu$ and $\nu$ are defined).

**Problem 6.4.18:** A S-fuzzy right (resp. left) N-subbigroup $\mu$ of a S-binear ring N is S-normal of and only if $\mu^+ = \mu$. Prove.

**Problem 6.4.19:** Under the condition of the above problem prove

 i. $(\mu^+)^+ = \mu$ ($\mu$ S-normal).
 ii. $(\mu^+)^+ = \mu^+$ (if $\mu$ is not S-normal).

**Problem 6.4.20:** Give an example of a S-fuzzy non-associative biseminear ring, which is not a S-strong, fuzzy non-associative biseminear ring.

**Problem 6.4.21:** Define Smarandache fuzzy congruence of a binear ring module and illustrate it with example.

**Problem 6.4.22:** Define the notion of Smarandache fuzzy quotient R-module of a S-binear ring.

**Problem 6.4.23:** Define a Smarandache fuzzy equivalence relation on

   i. S-biseminear ring.
   ii. S-non associative binear ring.
   iii. S-non-associative biseminear ring.

**Problem 6.4.24:** Define Smarandache fuzzy congruence relation on S-binear ring, N both when N is associative and N is non-associative.

**Problem 6.4.25:** Obtain some nice results on S-non-associative binear rings.

**CHAPTER SEVEN**

# APPLICATIONS OF SMARANDACHE FUZZY ALGEBRAIC STRUCTURES

This chapter has three sections. In the first section we just recall the applications of Smarandache algebraic structures. In section two we will recall the application of fuzzy algebraic structures and Smarandache fuzzy algebraic structures in the application of automaton theory. Section three proposes problems.

## 7.1 Applications of Smarandache Algebraic Structures

In this section we give some of the applications of Smarandache algebraic structures viz. Smarandache groupoids, Smarandache near-rings and Smarandache semirings in automaton theory, in error correcting codes and in the construction of S-sub-biautomaton which are given in about forty definitions. For more about these concepts refer [128, 130, 135].

The semi-automaton and automaton are built using the fundamental algebraic structure semigroups. In this chapter we use the generalized concept of semigroups viz. Smarandache groupoids which always contains a semigroup in it are used to construct Smarandache semi automaton and Smarandache automaton. Thus, the Smarandache groupoids find its application in the construction of finite machines. Here we introduce the concept of Smarandache semi automaton and Smarandache automaton using Smarandache free groupoids. This chapter starts with the definition of Smarandache free groupoids.

**DEFINITION 7.1.1:** *Let S be non empty set. Generate a free groupoid using S and denote it by ⟨S⟩. Clearly the free semigroup generated by the set S is properly contained in ⟨S⟩; as in ⟨S⟩ we may or may not have the associative law to be true.*

*Remark:* Even $(ab)c \neq a(bc)$ in general for all $a, b, c \in \langle S \rangle$. Thus unlike a free semigroup where the operation is associative, in case of free groupoid we do not assume the associativity while placing them in juxtaposition.

**THEOREM 7.1.1:** *Every free groupoid is a S-free groupoid.*

*Proof:* Clearly if S is the set which generates the free groupoid then it will certainly contain the free semigroup generated by S, so every free groupoid is a S-free groupoid.

We just recall the definition of semi automaton and automaton from [95, 128].

**DEFINITION 7.1.2:** *A semi-automaton is a triple $Y = (Z, A, \delta)$ consisting of two non - empty sets Z and A and a function $\delta: Z \times A \to Z$, Z is called the set of states, A the input alphabet and $\delta$ the "next state function" of Y.*

Let A = {$a_1$, ..., $a_n$} and Z = {$z_1$, ..., $z_k$}. The semi automaton Y = (Z, A, $\delta$). The semi-automaton can also be described by the transition table.

Description by Table:

| $\delta$ | $a_1$ | ... | $a_n$ |
|---|---|---|---|
| $z_1$ | $\delta(z_1, a_1)$ | ... | $\delta(z_1, a_n)$ |
| $\vdots$ | $\vdots$ | | $\vdots$ |
| $z_k$ | $\delta(z_k, a_1)$ | ... | $\delta(z_k, a_n)$ |

as $\delta : Z \times A \to Z$, $\delta(z_i, a_j) \in Z$. The semi automaton can also be described by graphs.

Description by graphs:

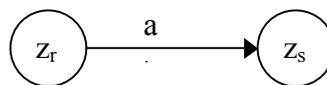

**Figure 7.1.1**

We depict $z_1$, ..., $z_k$ as 'discs' in the plane and draw an arrow labeled $a_i$ from $z_r$ to $z_s$ if $\delta(z_r, a_i) = z_s$. This graph is called the state graph.

*Example 7.1.1:* Z – set of states and A – input alphabet. Let Z = {0, 1, 2} and A = {0, 1}. The function $\delta: Z \times A \to Z$ defined by $\delta(0, 1) = 1 = \delta(2, 1) = \delta(1, 1)$, $\delta(0, 0) = 0$, $\delta(2, 0) = 1$, $\delta(1, 0) = 0$.

This is a semi automaton, having the following graph.

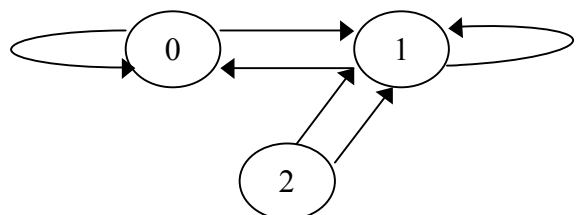

Figure 7.1.2

The description by table is:

| $\delta$ | 0 | 1 |
|---|---|---|
| 0 | 0 | 1 |
| 1 | 0 | 1 |
| 2 | 1 | 1 |

**DEFINITION 7.1.3:** *An Automaton is a quintuple K = (Z, A, B, $\delta$, $\Lambda$) where (Z, A, $\delta$) is a semi automaton, B is a non-empty set called the output alphabet and $\lambda: Z \times A \to B$ is the output function.*

If $z \in Z$ and $a \in A$, then we interpret $\delta(z, a) \in Z$ as the next state into which z is transformed by the input a, $\lambda(z, a) \in B$ is the output resulting from the input a.

Thus if the automaton is in state z and receives input a, then it changes to state $\delta(z, a)$ with output $\lambda(z, a)$.

If $A = \{a_1, ..., a_n\}$, $B = \{b_1, b_2, ..., b_m\}$ and $Z = \{z_1, ..., z_k\}$, $\delta : Z \times A \to Z$ and $\lambda : Z \times A \to B$ given by the description by tables where $\delta(z_k, a_i) \in Z$ and $\lambda(z_k, a_j) \in B$.

Description by Tables:

| $\delta$ | $a_1$ | ... | $a_n$ |
|---|---|---|---|
| $z_1$ | $\delta(z_1, a_1)$ | ... | $\delta(z_1, a_n)$ |
| $\vdots$ | $\vdots$ | | $\vdots$ |
| $z_k$ | $\delta(z_k, a_1)$ | ... | $\delta(z_k, a_n)$ |

In case of automaton, we also need an output table.

| $\lambda$ | $a_1$ | ... | $a_n$ |
|---|---|---|---|
| $z_1$ | $\lambda(z_1, a_1)$ | ... | $\lambda(z_1, a_n)$ |
| $\vdots$ | $\vdots$ | | $\vdots$ |
| $z_k$ | $\lambda(z_k, a_1)$ | ... | $\lambda(z_k, a_n)$ |

Description by Graph:

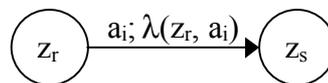

**Figure 7.1.3**

*Example 7.1.2:* Let $Z = \{z_1, z_2\}$ $A = B = \{0, 1\}$ we have the following table and graph:

| $\delta$ | 0 | 1 |
|---|---|---|
| $z_0$ | $z_0$ | $z_1$ |
| $z_1$ | $z_1$ | $z_0$ |

| $\lambda$ | 0 | 1 |
|---|---|---|
| $z_0$ | 0 | 1 |
| $z_1$ | 0 | 1 |

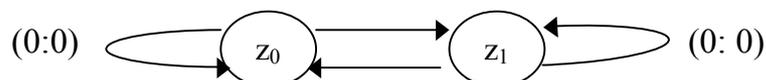

403
**Figure 7.1.4**

This automaton is known as the parity check automaton.

Now it is important and interesting to note that Z, A and B are only non-empty sets. They have no algebraic operation defined on them. The automatons and semi automatons defined in this manner do not help one to perform sequential operations. Thus, it is reasonable to consider the set of all finite sequences of elements of the set A including the empty sequence $\Lambda$.

In other words, in our study of automaton we extend the input set A to the free monoid $\overline{A}$ and $\overline{\lambda}: Z \times \overline{A} \to \overline{B}$ where $\overline{B}$ is the free monoid generated by B. We also extend functions $\delta$ and $\lambda$ from $Z \times A$ to $Z \times \overline{A}$ by defining $z \in Z$ and $a_1, ..., a_n \in A$ by $\overline{\delta}: Z \times \overline{A} \to Z$.

$$\overline{\delta}(z, \Lambda) = z$$
$$\overline{\delta}(z, a_1) = \delta(z, a_1)$$
$$\overline{\delta}(z, a_1 a_2) = \delta(\overline{\delta}(z, a_1), a_2)$$
$$\vdots$$
$$\overline{\delta}(z, a_1 a_2 ... a_n) = \delta(\overline{\delta}(z, a_1 a_2 ... a_{n-1}), a_n)$$

and $\lambda: Z \times A \to B$ by $\overline{\lambda}: Z \times \overline{A} \to \overline{B}$ by

$$\overline{\lambda}(z, \Lambda) = \Lambda$$
$$\overline{\lambda}(z, a_1) = \lambda(z, a_1)$$
$$\overline{\lambda}(z, a_1 a_2) = \lambda(z, a_1)\overline{\lambda}(\delta(z, a_1), a_2)$$
$$\vdots$$
$$\overline{\lambda}(z, a_1 a_2 ... a_n) = \lambda(z, a_1)\overline{\lambda}(\delta(z, a_1), a_2 ... a_r).$$

The semi-automaton $Y = (Z, A, \delta)$ and automaton $K = (Z, A, B, \delta, \lambda)$ is thus generalized to the new semi-automaton $Y = (Z, \overline{A}, \overline{\delta})$ and automaton $K = (Z, \overline{A}, \overline{B}, \overline{\delta}, \overline{\lambda})$. We throughout in this book mean by new semi automaton $Y = (Z, \overline{A}, \overline{\delta})$ and new automaton $K = (Z, \overline{A}, \overline{B}, \overline{\delta}, \overline{\lambda})$.

The concept of Smarandache semi automaton and Smarandache automaton was first defined by the author in the year 2002. The study of Smarandache semi automaton and uses the concept of groupoids and free groupoids.

Now we define Smarandache semi automaton and Smarandache automaton:

**DEFINITION 7.1.4:** $Y_s = (Z, \overline{A}_s, \overline{\delta}_s)$ is said to be a Smarandache semi automaton (S-semi automaton) if $\overline{A}_s = \langle A \rangle$ is the free groupoid generated by A with $\Lambda$ the empty element adjoined with it and $\overline{\delta}_s$ is the function from $Z \times \overline{A}_s \to Z$. Thus the Smarandache semi automaton contains $Y = (Z, \overline{A}, \overline{\delta})$ as a new semi automaton which is a proper sub-structure of $Y_s$.

*Or equivalently, we define a S-semi automaton as one, which has a new semi automaton as a sub-structure.*

The advantages of the S-semi automaton are if for the triple $Y = (Z, A, \delta)$ is a semi-automaton with Z, the set of states, A the input alphabet and $\delta : Z \times A \to Z$ is the next state function. When we generate the S-free groupoid by A and adjoin with it the empty alphabet $\Lambda$ then we are sure that $\overline{A}$ has all free semigroups. Thus, each free semigroup will give a new semi automaton. Thus by choosing a suitable A we can get several new semi automaton using a single S-semi automaton.

We now give some examples of S-semi-automaton using finite groupoids. When examples of semi automaton are given usually the books use either the set of modulo integers $Z_n$ under addition or multiplication we use groupoid built using $Z_n$.

**DEFINITION 7.1.5:** $\overline{Y}'_S = (Z_1, \overline{A}_s, \overline{\delta}'_s)$ is called the Smarandache subsemi automaton (S-subsemi-automaton) of $\overline{Y}_S = (Z_2, \overline{A}_s, \overline{\delta}'_s)$ denoted by $\overline{Y}'_S \leq \overline{Y}_S$ if $Z_1 \subset Z_2$ and $\overline{\delta}'_s$ is the restriction of $\overline{\delta}_s$ on $Z_1 \times \overline{A}_s$ and $\overline{Y}'_S$ has a proper subset $\overline{H} \subset \overline{Y}'_S$ such that $\overline{H}$ is a new semi automaton.

*Example 7.1.3:* Let $Z = Z_4(2, 1)$ and $A = Z_6(2, 1)$. The S-semi automaton $(Z, A, \delta)$ where $\delta : Z \times A \to Z$ is given by $\delta(z, a) = z \bullet a \pmod 4$.

We get the following table:

| $\delta$ | 0 | 1 | 2 | 3 |
|---|---|---|---|---|
| 0 | 0 | 1 | 2 | 3 |
| 1 | 2 | 3 | 0 | 1 |
| 2 | 0 | 1 | 2 | 3 |
| 3 | 2 | 3 | 0 | 1 |

We get the following graph for this S-semi automaton:

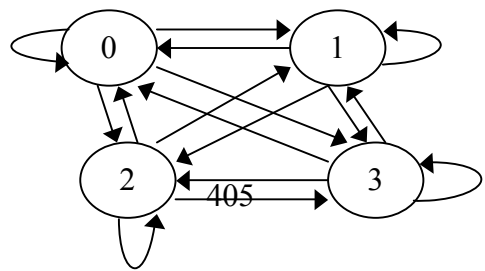

This has a nice Smarandache subsemi-automaton given by the table.

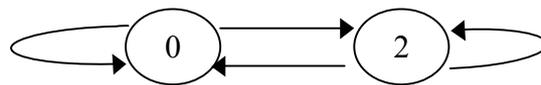

**Figure 7.1.6**

This Smarandache semi automaton is nothing but given by the states {0, 2}.

*Example 7.1.4:* Let $Z = Z_3(1, 2)$ and $A = Z_4(2, 2)$ the triple $(Z, A, \delta)$ is a S-semi automaton with $\delta(z, a) = (z * a) \pmod 3$ ' $*$ ' the operation in A given by the following table:

| δ | 0 | 1 | 2 | 3 |
|---|---|---|---|---|
| 0 | 0 | 2 | 1 | 0 |
| 1 | 2 | 1 | 0 | 2 |
| 2 | 1 | 0 | 2 | 1 |

The graphical representation of the S-semi automaton is

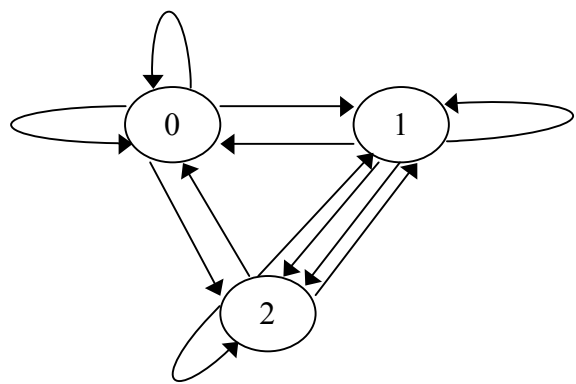

Figure 7.1.7

*Example 7.1.5:* Now let $Z = Z_4(2, 2)$ and $A = Z_3(1, 2)$ we define $\delta(z, a) = (z * a) \pmod 4$ '*' as in Z. The table for the semi automaton is given by

| δ | 0 | 1 | 2 |
|---|---|---|---|
| 0 | 0 | 2 | 0 |
| 1 | 2 | 0 | 2 |
| 2 | 0 | 2 | 0 |
| 3 | 2 | 0 | 2 |

The graph for it is

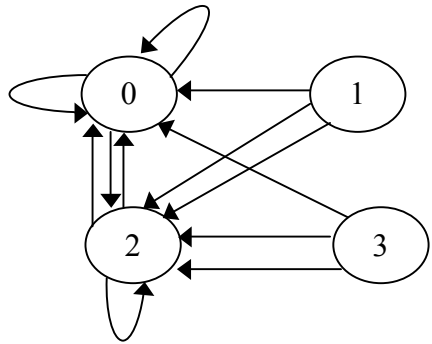

Figure 7.1.8

Thus this has a S-subsemi automaton $Z_1$ given by $Z_1 = \{0, 2\}$ states.

**DEFINITION 7.1.6:** $\overline{K}_S = (Z, \overline{A}_s, \overline{B}_s, \overline{\delta}_s, \overline{\lambda}_s)$ is defined to be a S-automaton if $\overline{K} = (Z, \overline{A}_s, \overline{B}_s, \overline{\delta}_s, \overline{\lambda}_s)$ is the new automaton and $\overline{A}_s$ and $\overline{B}_s$, the Smarandache free groupoids so that $\overline{K} = (Z, \overline{A}_s, \overline{B}_s, \overline{\delta}_s, \overline{\lambda}_s)$ is the new automaton got from K and $\overline{K}$ is strictly contained in $\overline{K}_S$.

Thus S-automaton enables us to adjoin some more elements which is present in A and freely generated by A, as a free groupoid; that will be the case when the compositions may not be associative.

Secondly, by using Smarandache automaton we can couple several automaton as

$$Z = Z_1 \cup Z_2 \cup \ldots \cup Z_n$$
$$A = A_1 \cup A_2 \cup \ldots \cup A_n$$
$$B = B_1 \cup B_2 \cup \ldots \cup B_n$$
$$\lambda = \lambda_1 \cup \lambda_2 \cup \ldots \cup \lambda_n$$
$$\delta = \delta_1 \cup \delta_2 \cup \ldots \cup \delta_n.$$

where the union of $\lambda_i \cup \lambda_j$ and $\delta_i \cup \delta_j$ denote only extension maps as '$\cup$' has no meaning in the composition of maps, where $K_i = (Z_i, A_i, B_i, \delta_i, \lambda_i)$ for i = 1, 2, 3, ..., n and $\overline{K} = \overline{K}_1 \cup \overline{K}_2 \cup \ldots \cup \overline{K}_n$. Now

$$\overline{K}_S = (\overline{Z}_s, \overline{A}_s, \overline{B}_s, \overline{\lambda}_s, \overline{\delta}_s)$$

is the S-automaton.

A machine equipped with this Smarandache automaton can use any new automaton as per need.

We give some examples of S-automaton using Smarandache groupoids.

***Example 7.1.6:*** Let $Z = Z_4 (3, 2)$, $A = B = Z_5 (2, 3)$. $K = (Z, A, B, \delta, \lambda)$ is a S-automaton defined by the following tables where

$\delta (z, a) = z * a \pmod 4$ and
$\lambda (z, a) = z * a \pmod 5$.

| λ | 0 | 1 | 2 | 3 | 4 |
|---|---|---|---|---|---|
| 0 | 0 | 3 | 1 | 4 | 2 |
| 1 | 2 | 0 | 3 | 1 | 4 |
| 2 | 4 | 2 | 0 | 3 | 1 |
| 3 | 1 | 4 | 2 | 0 | 3 |

| δ | 0 | 1 | 2 | 3 | 4 |
|---|---|---|---|---|---|
| 0 | 0 | 2 | 0 | 2 | 0 |
| 1 | 3 | 1 | 3 | 1 | 3 |
| 2 | 2 | 0 | 2 | 0 | 2 |
| 3 | 1 | 3 | 1 | 3 | 1 |

We obtain the following graph:

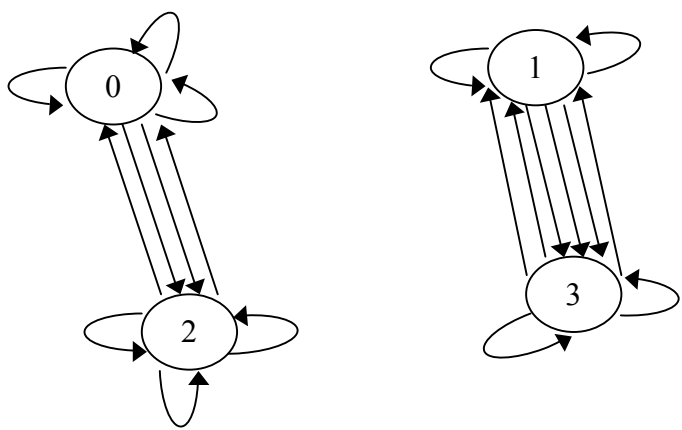

Figure 7.1.9

Thus we see this automaton has 2 S-sub automatons given by the states {0, 2} and {1, 3}.

***Example 7.1.7:*** Consider $Z_5 = Z_5(3, 2)$, $A = Z_3(0, 2)$, $B = Z_4(2, 3)$, $\delta (z, a) = z * a \pmod 5$, $\lambda (z, a) = z * a \pmod 4$. The S-automaton $(Z, A, B, \delta, \lambda)$ is given by the following tables:

| δ | 0 | 1 | 2 |
|---|---|---|---|
| 0 | 0 | 3 | 1 |
| 1 | 3 | 1 | 4 |
| 2 | 1 | 4 | 2 |
| 3 | 4 | 2 | 0 |
| 4 | 2 | 0 | 3 |

| λ | 0 | 1 | 2 | 3 | 4 |
|---|---|---|---|---|---|
| 0 | 0 | 3 | 2 | 1 | 0 |
| 1 | 2 | 1 | 0 | 3 | 2 |
| 2 | 0 | 3 | 2 | 1 | 0 |
| 3 | 2 | 1 | 0 | 3 | 2 |
| 4 | 0 | 3 | 2 | 1 | 0 |

This has no proper S-sub automaton.

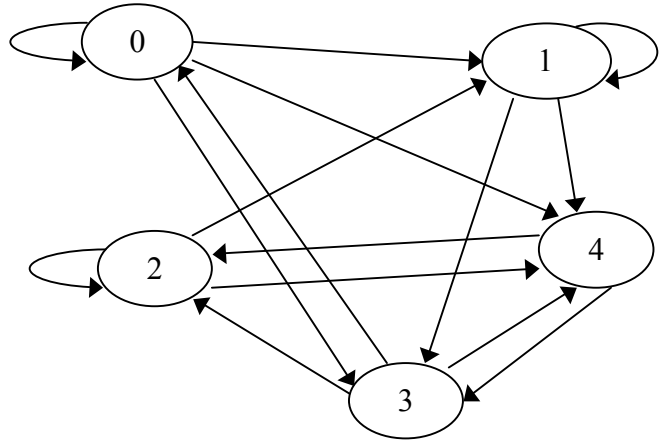

**Figure 7.1.10**

**DEFINITION 7.1.7:** $\overline{K}_s' = (Z_1, \overline{A}_s, \overline{B}_s, \overline{\delta}_s, \overline{\lambda}_s)$ is called S-sub-automaton of $\overline{K}_s = (Z_2, \overline{A}_s, \overline{B}_s, \overline{\delta}_s, \overline{\lambda}_s)$ denoted by $\overline{K}_s' \leq \overline{K}_s$ if $Z_1 \subseteq Z_2$ and $\overline{\delta}_s'$ and $\overline{\lambda}_s'$ are the restriction of $\overline{\delta}_s$ and $\overline{\lambda}_s$ respectively on $Z_1 \times \overline{A}_s$ and has a proper subset $\overline{H} \subset \overline{K}_s'$ such that $\overline{H}$ is a new automaton.

**DEFINITION 7.1.8:** Let $\overline{K}_1$ and $\overline{K}_2$ be any two S-automaton where $\overline{K}_1 = Z_1, \overline{A}_s, \overline{B}_s, \overline{\delta}_s, \overline{\lambda}_s)$ and $\overline{K}_2 = (Z_2, \overline{A}_s, \overline{B}_s, \overline{\delta}_s, \overline{\lambda}_s)$. A map $\phi : \overline{K}_1$ to $\overline{K}_2$ is a Smarandache automaton homomorphism (S-automaton homomorphism) if $\phi$ restricted from $K_1 = (Z_1, A_1, B_1, \delta_1, \lambda_1)$ and $K_2 = (Z_2, A_2, B_2, \delta_2, \lambda_2)$ denoted by $\phi_r$ is a automaton homomorphism from $K_1$ to $K_2$. $\phi$ is called a monomorphism (epimorphism or isomorphism) if there is an isomorphism $\phi_r$ from $K_1$ to $K_2$.

**DEFINITION 7.1.9:** Let $\overline{K}_1$ and $\overline{K}_2$ be two S-automatons, where $\overline{K}_1 = (Z_1, \overline{A}_s, \overline{B}_s, \overline{\delta}_s, \overline{\lambda}_s)$ and $\overline{K}_2 = (Z_2, \overline{A}_s, \overline{B}_s, \overline{\delta}_s, \overline{\lambda}_s)$.

*The Smarandache automaton direct product (S-automaton direct product) of $\overline{K}_1$ and $\overline{K}_2$ denoted by $\overline{K}_1 \times \overline{K}_2$ is defined as the direct product of the automaton $K_1 = (Z_1, A_1, B_1, \delta_1, \lambda_1)$ and $K_2 = (Z_2, A_2, B_2, \delta_2, \lambda_2)$ where $K_1 \times K_2 = (Z_1 \times Z_2, A_1 \times A_2, B_1 \times B_2, \delta, \lambda)$ with $\delta$ $((z_1, z_2), (a_1, a_2)) = (\delta_1 (z_1, a_1), \delta_2 (z_2, a_2)), \lambda ((z_1, z_2), (a_1, a_2)) = (\lambda_1 (z_1, a_1), \lambda_2 (z_2, a_2))$ for all $(z_1, z_2) \in Z_1 \times Z_2$ and $(a_1, a_2) \in A_1 \times A_2$.*

*Remark:* Here in $\overline{K}_1 \times \overline{K}_2$ we do not take the free groupoid to be generated by $A_1 \times A_2$ but only free groupoid generated by $\overline{A}_1 \times \overline{A}_2$.

Thus the S-automaton direct product exists wherever a automaton direct product exists. We have made this in order to make the Smarandache parallel composition and Smarandache series composition of automaton extendable in a simple way.

**DEFINITION 7.1.10:** *A S-groupoid $G_1$ divides a S-groupoid $G_2$ if the corresponding semigroups $S_1$ and $S_2$ of $G_1$ and $G_2$ respectively divides, that is, if $S_1$ is a homomorphic image of a sub-semigroup of $S_2$.*

*In symbols $G_1 | G_2$. The relation divides is denoted by '|'.*

**DEFINITION 7.1.11:** *Let $\overline{K}_1 = (Z_1, \overline{A}_s, \overline{B}_s, \overline{\delta}_s, \overline{\lambda}_s)$ and $\overline{K}_2 = (Z_2, \overline{A}_s, \overline{B}_s, \overline{\delta}_s, \overline{\lambda}_s)$ be two Smarandache automaton. We say the S-automaton $\overline{K}_1$ divides the S-automaton $\overline{K}_2$ if in the automatons $K_1 = (Z_1, A, B, \delta_1, \lambda_1)$ and $K_2 = (Z_2, A, B, \delta_2, \lambda_2)$, if $K_1$ is the homomorphic image of a sub-automaton of $K_2$. Notationally $K_1 | K_2$.*

**DEFINITION 7.1.12:** *Two S-automaton $\overline{K}_1$ and $\overline{K}_2$ are said to be equivalent if they divide each other. In symbols $\overline{K}_1 \sim \overline{K}_2$.*

We proceed on to define direct product Smarandache automaton and study about them. We can extend the direct product of semi automaton to more than two Smarandache automatons.

Using the definition of direct product of two automaton $K_1$ and $K_2$ with an additional assumption we define Smarandache series composition of automaton.

**DEFINITION [128]:** *Let $K_1$ and $K_2$ be any two Smarandache automatons where $\overline{K}_1 = (Z_1, \overline{A}_s, \overline{B}_s, \overline{\delta}_s, \overline{\lambda}_s)$ and $\overline{K}_2 = (Z_2, \overline{A}_s, \overline{B}_s, \overline{\delta}_s, \overline{\lambda}_s)$ with an additional assumption $A_2 = B_1$.*

*The Smarandache automaton composition series denoted by $\overline{K}_1 \Vdash \overline{K}_2$ of $\overline{K}_1$ and $\overline{K}_2$ is defined as the series composition of the automaton $K_1 = (Z_1, A_1, B_1, \delta_1, \lambda_1)$ and $K_2 = (Z_2, A_2, B_2, \delta_2, \lambda_2)$ with $\overline{K}_1 \Vdash \overline{K}_2 = (Z_1 \times Z_2, A_1, B_2, \delta, \lambda)$ where $\delta((z_1, z_2), a_1) = (\delta_1(z_1, a_1), \delta_2(z_2, \lambda_1(z_1, a_1)))$ and $\lambda((z_1, z_2), a_1) = (\lambda_2(z_2, \lambda_1(z_1, a_1))) ((z_1, z_2) \in Z_1 \times Z_2, a_1 \in A_1)$.*

*This automaton operates as follows: An input $a_1 \in A_1$ operates on $z_1$ and gives a state transition into $z_1' = \delta_1(z_1, a_1)$ and an output $b_1 = \lambda_1(z_1, a_2) \in B_1 = A_2$. This output $b_1$ operates on $Z_2$ transforms a $z_2 \in Z_2$ into $z_2' = \delta_2(a_2, b_1)$ and produces the output $\lambda_2(z_2, b_1)$.*

*Then $\overline{K}_1 \Vdash \overline{K}_2$ is in the next state $(z_1', z_2')$ which is clear from the following circuit:*

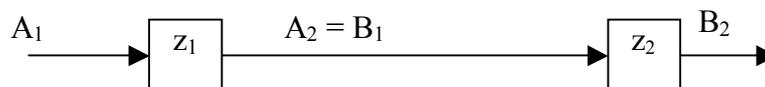

**Figure 7.1.11**

Now a natural question would be do we have a direct product, which corresponds to parallel composition, of the 2 Smarandache automatons $\overline{K}_1$ and $\overline{K}_2$. Clearly the Smarandache direct product of automatons $\overline{K}_1 \times \overline{K}_2$ since $Z_1$, and $Z_2$ can be interpreted as two parallel blocks. $A_i$ operates on $Z_i$ with output $B_i$ ($i \in \{1,2\}$), $A_1 \times A_2$ operates on $Z_1 \times Z_2$, the outputs are in $B_1 \times B_2$.

The circuit is given by the following diagram:

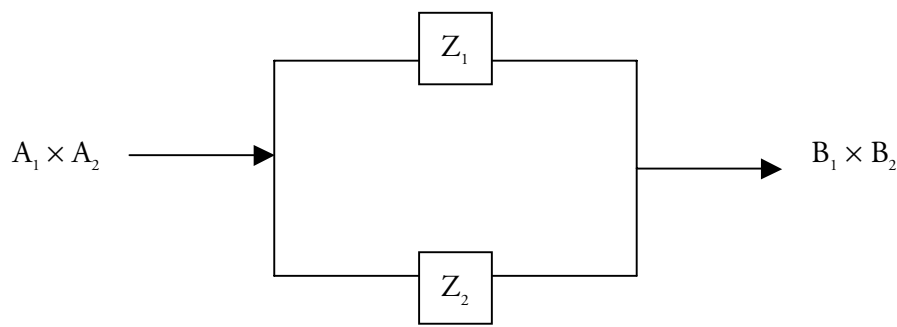

Figure 7.1.12

The concept of group semi-automatons is associated with the S-syntactic near-ring to get a finite machine which enjoys a varied structural property.

**DEFINITION [31]**: *A group semi-automaton is an ordered triple $S = (G, I, \delta)$ where $G$ is a group, $I$ is a set of inputs and $\delta : G \times I \to G$ is a transition function.*

***Example 7.1.8:*** Let $G = S_3$, where $S = \{p_0, p_1, \ldots, p_5\}$ $I = \{0, 1, 2,\}$ be the set of inputs. The transition function $\delta : G \times I \to G$ defined by

$$\begin{aligned}
\delta(p_0, 0) &= p_0 \\
\delta(p_i, 0) &= p_i \text{ for } i = 1, 2, \ldots, 5 \\
\delta(p_i, 1) &= p_{i+1} \text{ and } \delta(p_5, 1) = 0, \text{ for } i = 1, 2, \ldots, 4 \\
\delta(p_i, 2) &= p_{i+2}, \text{ for } i = 1, 2, 3 \text{ and } \delta(p_4, 2) = p_0 \; \delta = (p_5, 2) = p_1
\end{aligned}$$

Clearly $(G, I, \delta)$ is a group semi-automaton.

Now using the notions of [130] we define the following. [130] will give the papers used in this book. By using cross reference one can find it. As the main motivation of the book is to study Smarandache notions we just give the Smarandache books as main references.

**DEFINITION [31]**: *A semi-automaton $S = (Q, X, \delta)$ with state set $Q$ input monoid $X$ and state transition function $\delta$; $\delta : Q \times X \to Q$ is called a group semi-automaton if $Q$ is an additive group. $\delta$ is called additive if there is some $x_0 \in X$ with $\delta(q, x) = \delta(q, x_0) + \delta(0, x)$ and $\delta(q - q^1, x_0) = \delta(q, x_0) - \delta(q^1, x_0)$ for all $q, q^1 \in Q$ and $x \in X$. Then there is some homomorphism*

$$\psi : Q \to Q \text{ and some map}$$

$$\alpha : X \to Q \text{ with } (x_0) = 0 \text{ and } \delta(q, x) = \psi(q) + \alpha(x).$$

*Let $\delta_x$ for a fixed, $x \in X$ be the map $Q \to Q$; $q \to \delta(q, x)$ then $\{\delta_x / x \in X\}$ generates a subnear-ring $N(S)$ of near-ring $(M(Q), +, \bullet)$ of all mappings on Q. This near-ring $N(S)$ is called the syntactic near-ring.*

In this section we introduce the concept of Smarandache semigroup semi-automaton analogous to the group semi-automaton. Further for each group semi-automaton an associated syntactic near-ring is also given. Here we define the Smarandache analogous of it.

**DEFINITION 7.1.13**: *Let G be S-semigroup. A Smarandache S-semigroup semiautomaton (S-S-semigroup semiautomaton) is an ordered triple $S(S) = (G, I, \delta)$ where G is a S-semigroup; I is a set of inputs and $\delta : G \times I \to G$ is a transition function.*

Clearly all group semi-automatons are trivially S-S-semigroup semiautomaton. Thus the class of S-S semigroup semiautomaton propertly contains the class of group semi-automaton.

***Example 7.1.9***: Let $G = \{0, 1, 2, \ldots, 5\} = Z_6$ and $I = \{0, 1, 2\}$ be the set of states; define the transition function $\delta : G \times I \to G$ by $\delta(g, i) = g \cdot i \pmod 6$ thus we get the S-semigroup semi-automaton. We give the following state graph representation for this S-S semigroup semi-automaton.

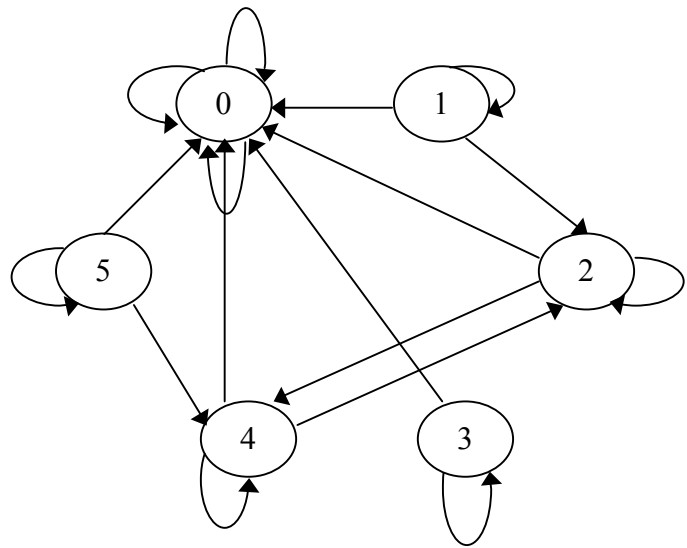

**FIGURE 7.1.13**

Now we see when we change the function $\delta_1$ for the same S-semigroup, G and I = {0, p1, 2} we get a different S-S-semigroup semiautomaton which is given by the following figure; here the transition fuctions. $\delta : G \times I \to G$ is defined by $\delta_1(g, i) = (g + i) \mod 6$. We get a nice symmetric S-S-semigroup semiautomaton.

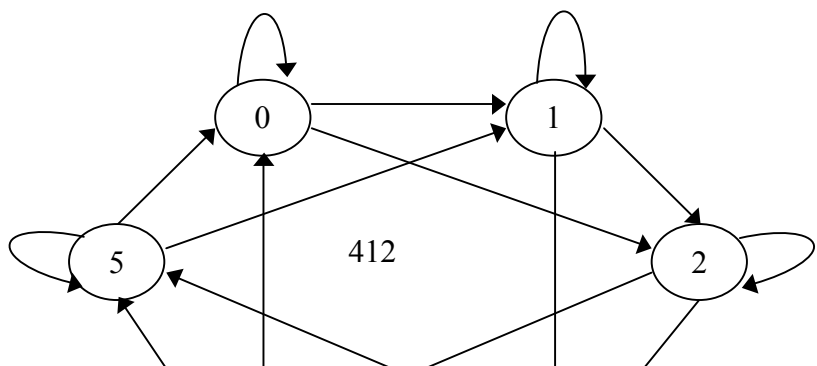

**DEFINITION 7.1.14**: *Let S (S) = (G, I, δ) where G is a S-semigroup I the set of input alphabets and δ the transition fuction as in case of S-S semigroup semi-automaton. We call a proper subset B ⊂ G with B a S-subsemigroup of G and $I_1$ ⊂ I denoted (B, $I_1$, δ) = S(SB) is called the Smarandache subsemigroup semi-automaton (S-subsemigroup semiautomaton) if δ : B × $I_1$ → B. A S-S-semigroup semi-automaton may fail to have sometimes S-S-subsemigroup subsemi-automaton.*

*Example 7.1.10*: In example 7.1.9 with S = {$Z_6$, I, δ} when we take B = {0, 2, 4} we get a S-S-semigroup subsemi-automaton

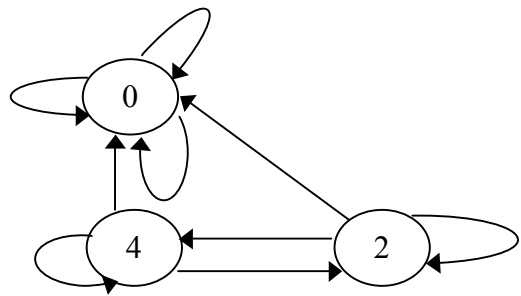

**FIGURE 7.1.15**

as B = {0, 2, 4} is a S-subsemigroup of $Z_6$.

Now to the best of our knowledge we have not seen the definition of group automaton we define it and then proceed on to get corresponding near-rings.

**DEFINITION 7.1.15**: *A group automaton is an ordered quintuple A = {G, A, B, δ, λ} where G is a group, A is a set of inputs, B is a set of outputs δ : G × A → G and λ : G × A → B is the transition and output functions respectively.*

The Smarandache analogous of this would be as follows.

**DEFINITION 7.1.16**: *A Smarandache S-semigroup automaton (S-S semigroup automaton) is an ordered quintuple S(A) = (G, A, B, δ, λ) where G is a S-semigroup and A set of input alphabets and B set of output alphabets, δ and λ are the transition and output functions respectively.*

*If δ restricted form G × A → G maps a subgroup T ⊂ G into itself δ : T × A → T then we can say the S-S-semigroup automaton has (T ⊂ G, A, B, δ, λ) to be a group automaton.*

We call this level I S-S-semigroup automaton. Now we proceed on to define S-S-semigroup semiautomaton and S-S-semigroup automaton of level II.

**DEFINITION 7.1.17**: *Let S = (G, I, δ) be an ordered triple where G is only a S-semigroup. If B ⊂ G and B is a group under the operations of G, and if $S^1$ = (B, I, δ) is a group semi-automaton with B a group, I is a set of inputs and δ: B × I → B then S = (G, I, δ) is called the Smarandache S-semigroup semi-automaton of level II (S-S-semigroup semiautomaton II). Similarly we define Smarandache S-semigroup automaton of level II (S-S-semigroup automaton II).*

It is left for the reader to find any relation existing between the two levels of automaton. Now we proceed on to see how these concepts will produce for us the S-near-rings.

**DEFINITION 7.1.18**: *Let S = (Q, X δ) be a S-S-semigroup semi-automaton where Q is a S-semigroup under addition '+' i.e. G ⊂ Q is such that G is a S-semigroup under '+' and G a proper subset of Q Let δ : G × X → G is well defined i.e. restriction of δ to G is well defined so (G, X, δ) = C becomes a semi-automaton C is called additive if there is some $x_0$ ∈ X with*

$$\delta(q, x) = \delta(q, x_0) + \delta(0, x) \text{ and}$$
$$\delta(q - q^1, x_0) = \delta(q_1 x_0) - \delta(q^1, x_0) \text{ for all } q, q^1 \in G \text{ and } x \in X.$$

*Then there is some homomorphism*

$$\psi : G \to G \text{ and some map}$$
$$\alpha : X \to G \text{ with } \psi(x_0) = 0 \text{ and } \delta(q, x) = \psi(q) + \alpha(x).$$

*Let $\delta_x$ for a fixed x ∈ X be the map from G → G; q → δ (q, x). Then {$\delta_x$ / x ∈ X} generates a subnear-ring N (P) of the near-ring (M (G), +, •) of all mappings on G; this near-ring N(P) is called as the Smarandache syntactic near-ring (S-syntactic near-ring).*

The varied structural properties are enjoyed by this S-syntactic near-ring for varying S-semigroups. The most interesting feature about these S-syntactic near-ring is that for one S-S-semigroup semi-automaton we have several group semiautomaton depending on the number of valid groups in the S-semigroup. This is the vital benefit in defining S-S-semigroup automaton and S-syntactic near-ring. So for a given S-S-semigroup semi-automaton we can have several S-syntactic near-ring.

Thus the Smarandache notions in this direction has evolved several group automaton for a given S-semigroup.

Another application is in the error-correcting codes.

The study of how experiments can be organized systematically so that statistical analysis can be applied in an interesting problem which is carried out by several researchers. In the planning of experiments it often occurs that results are influenced by phenomena outside the control of the experimenter. The introduction of balanced incomplete block design (BIBD) helps in avoiding undesirable influences in the experiment. In general, if we have to test the effect of r different conditions with m possibilities for each conditions this leads to a set of r orthogonal latin squares.

A planar near-ring can be used to construct balanced incomplete block designs (BIBD) of high efficiency; we just state how they are used in developing error correcting codes as codes which can correct errors is always desirable than the ones which can only detect errors. In view of this we give the following definition. Just we recall the definition of BIBD for the sake of completeness.

**DEFINITION 7.1.19**: *A balanced incomplete block design (BIBD) with parameters ($v$, $b$, $r$, $k$, $\lambda$) is a pair (P, B) with the following properties.*

 a. *P is a set with $v$ elements.*
 b. *B = ($B_1$, ..., $B_b$) is a subset of p (P) with b elements.*
 c. *Each $B_i$ has exactly k elements where $k < v$ each unordered pair (p, q) with p, q $\in$ P, p $\neq$ q occurs in exactly $\lambda$ elements in B.*

*The set $B_1$, ..., $B_b$ are called the blocks of BIBD. Each a $\in$ P occurs in exactly r sets of B. Such a BIBD is also called a ($v$, $b$, $r$, $k$, $\lambda$) configuration or 2 – ($v$, $k$, $\lambda$) tactical configuration or design. The term balance indicates that each pair of elements occurs in exactly the same number of block, the term incomplete means that each block contains less than $v$ - elements. A BIBD is symmetric if $v = b$.*

*The incidence matrix of a ($v$, $b$, $r$, $k$, $\lambda$) configuration is the $v \times b$ matrix $A = (a_{ij})$ where*

$$a_{ij} = \begin{cases} 1 \text{ if } i \in B_j \\ 0 \text{ otherwise} \end{cases},$$

*here i denotes the $i^{th}$ element of the configuration. The following conditions are necessary for the existence of a BIBD with parameters $v$, $b$, $r$, $k$, $\lambda$.*

 i. $bk = rv$.
 ii. $r(k – 1) = \lambda (v – 1)$.
 iii. $b \geq v$.

*Recall a near-ring N is called planar (or Clay near-ring) if for all equation x o a = x o b + c. (a, b, c $\in$ N, a $\neq$ b) have exactly one solution x $\in$ N.*

**Example 7.1.11**: Let ($Z_5$, '+', *) where '+' is the usual '+' and '*' is n * 0 = 0, n * 1 = n * 2 = n, n * 3 = n * 4 = 4n for all n $\in$ N. Then 1 $\equiv$ 2 and 3 $\equiv$ 4. N is planar near-ring; the equation x * 2 = x * 3 + 1 with 2 $\neq$ 3 has unique solution x = 3.

A planar near-ring can be used to construct BIBD of high efficiency where by high efficiency 'E' we mean E = λv / rk; this E is a number between 0 and 1 and it estimates the quality of any statistical analysis if E ≥ 0.75 the quality is good.

Now, how does one construct a BIBD from a planar ring? This is done in the following way.

Let N be a planar ring. Let a ∈ N. Define $g_a$ : N → N by n → n o a and form G = {$g_a$ / a ∈ N}. Call a ∈ N "group forming" if a o N is a subgroup of (N, +). Let us call sets a o N + b (a ∈ N*, b ∈ N) blocks. Then these blocks together with N as the set of "points" form a tactical configuration with parameters

$$(\upsilon, b, r, k, \lambda) = (\upsilon, \frac{\alpha_1 \upsilon}{|G|} + \alpha_2 \upsilon, \alpha_1 + \alpha_2 |G|, |G|, \lambda)$$

where $\upsilon$ = |N| and $\alpha_1$ ($\alpha_2$) denote the number of orbits of F under the group G \ {0} which consists of entirely of group forming elements. The tactical configuration is a BIBD if and only if either all elements are group forming or just 0 is group forming.

Now how are they used in obtaining error correcting codes. Now using the planar near-ring one can construct a BIBD. By taking either the rows or the columns of the incidence matrix of such a BIBD one can obtain error correcting codes with several nice features. Now instead of using a planar near-ring we can use Smarandache planar near-ring. The main advantage in using the Smarandache planar near-ring is that for one S-planar near-ring we can define more than one BIBD. If there are m-near-field in the S-planar near-ring we can build m BIBD. Thus this may help in even comparison of one BIBD with the another and also give the analysis of the common features.

Thus the S-planar near-ring has more advantages than the usual planar near-rings and hence BIBD's constructed using S-planar near-rings will even prove to be an efficient error correcting codes.

Finally we give the applications of Smarandache algebraic bistructures. Using the notion of Smarandache bisemigroup we build bisemiautomaton and bi-automatons.

**DEFINITION 7.1.20:** *A bisemi-automaton is a triple X = (Z, A, δ) consisting of two non-empty sets Z and A = $A_1 \cup A_2$ and a function δ : Z × ($A_1 \cup A_2$) → Z, Z is called the set of states and A = $A_1 \cup A_2$ the input alphabet and δ the next state function of X. We insist that if $X_1$ = (Z, $A_1$, δ) and $X_2$ = (Z, $A_2$, δ) then $X_1$ is not a subsemi-automaton of $X_2$ or $X_2$ is not a subsemi-automaton of $X_1$ but both $X_1$ and $X_2$ are semi-automatons. Thus in case of bisemi-automaton we consider the set A i.e. input alphabets as a union of two sets which may enjoy varied properties.*

*Hence we see in A = $A_1 \cup A_2$, $A_1$ may take states in a varied way and $A_2$ in some other distinct way. So when we take a set A and not the set $A_1 \cup A_2$ = A there is difference for*

$$Z_r \xrightarrow{a_1} Z_t \quad a_1 \in A_1$$

$$Z_r \xrightarrow{a_1} Z_p \quad a_1 \in A_2$$

*Hence the same $a_1 \in A$ may act in two different ways depending on, to which $A_i$ it belongs to, or it may happen in some cases in $A = A_1 \cup A_2$ with $A_1 \cap A_2 = \phi$, so that they will have a unique work assigned by the function $\delta$.*

**DEFINITION 7.1.21:** *A biautomaton is a quintuple $Y = (Z, A, B, \delta, \lambda)$ where $(Z, A, \delta)$ is a bisemi-automaton and $A = A_1 \cup A_2$, $B = B_1 \cup B_2$ and $\lambda : Z \times A \to B$ is the output function. Here also $Y_1 = (Z, A_1, B_1, \delta, \lambda)$ and $Y_2 = (Z, A_2, B_2, \delta, \lambda)$ then $Y_1$ is not a subautomaton of $Y_2$ and vice versa but both $Y_1$ and $Y_2$ are automatons. The bisemi-automaton and biautomaton will still find its application in a very interesting manner when we replace the sets by free bisemigroups. We denote by $\overline{A} = \overline{A}_1 \cup \overline{A}_2$, $\overline{A}_1$ and $\overline{A}_2$ are free monoids generated by $A_1$ and $A_2$ so $\overline{A}$ denotes the free bimonoid as both have the common identity, $\Lambda$ the empty sequence.*

$$\overline{\delta}(z, \Lambda) = z.$$
$$\overline{\delta}(z, a_1) = \delta(z, a_1)$$
$$\overline{\delta}(z, a_2) = \delta(z, a_2)$$
$$\vdots$$
$$\overline{\delta}(z, a_n) = \delta(z, a_n)$$

*where $a_1, a_2, ..., a_n \in A$. Here $\overline{\delta}(z, a_i)$ may have more than one $z_i$ associated with it.*

$$\overline{\delta}(z, a_1 a_2) = \delta(\overline{\delta}(z, a_1), a_2)$$
$$\vdots$$
$$\overline{\delta}(z, a_1 a_2 \cdots a_r) = \delta(\overline{\delta}(z, a_1 a_2 \cdots a_{r-1}), a_r)$$

*Thus the bisemi-automaton may do several work than the semi-automaton and can also do more intricate work.*

*Now we define biautomaton as follows:*

$$\overline{Y} = (Z, \overline{A}, \overline{B}, \overline{\delta}, \overline{\lambda})$$

*where $(Z, \overline{A}, \overline{\delta})$ is a bisemi-automaton. Here $B = B_1 \cup B_2$ and $\overline{B}$ is a free bimonoid.*

$$\overline{\lambda} : Z \times \overline{A} \to \overline{B}$$
$$\overline{\lambda}(z, \Lambda) = \Lambda$$
$$\overline{\lambda}(z, a_1) = \lambda(z, a_1)$$
$$\overline{\lambda}(z, a_1 a_2) = \lambda(z, a_1)\overline{\lambda}(\delta(z, a_1), a_2)$$
$$\vdots$$
$$\overline{\lambda}(z, a_1 a_2 \cdots a_n) = \lambda(z, a_1)\overline{\lambda}(\delta(z, a_1), a_2 a_3 \cdots a_n).$$

Thus we see the new biautomaton $\overline{Y}$ can perform multifold and multichannel operations than the usual automaton.

**DEFINITION 7.1.22:** $Y_1 = (Z_1, A, B, \delta_1, \lambda_1)$ is called a sub-biautomaton of $Y_2 = (Z, A, B, \delta, \lambda)$ (in symbols $Y_1 \subset Y_2$ if $Z_1 \subset Z$ and $\delta_1$ and $\lambda_1$ are the restrictions of $\delta$ and $\lambda$ respectively on $Z_1 \times A$ (Here $A = A_1 \cup A_2$ and $B = B_1 \cup B_2$ are subsets).

**DEFINITION 7.1.23:** Let $Y_1 = (Z_1, A_1, B_1, \delta_1, \lambda_1)$ and $Y_2 = (Z_2, A_2, B_2, \delta_2, \lambda_2)$. A biautomaton homomorphism $\phi : Y_1 \to Y_2$ is a triple $\phi = (\xi, \alpha, \beta) \in Z_2^{Z_1} \times A_2^{A_1} \times B_2^{B_1}$ with the property (where $A_1 = A_{11} \cup A_{12}$, $A_2 = A_{21} \cup A_{22}$, $B_1 = B_{11} \cup B_{12}$ and $B_2 = B_{21} \cup B_{22}$).

$$\xi(\delta_1(z,a)) = \delta_2(\xi(z), \alpha(a))$$
$$\beta(\lambda_1(z,a)) = \lambda_2(\xi(z), \alpha(a))$$
$$(z \in Z, a \in A_1 = A_{11} \cup A_{12}).$$

$\phi$ is called a monomorphism (epimorphism, isomorphism) if all functions $\xi$, $\alpha$, $\beta$ are injective (surjective or bijective) where we make relevant adjustments when for a single $a \in A_1$ there are two state changes.

Similarly for a single $a \in A_1$ we may have outputs depending on the section where we want to send the message. Thus we see the biautomaton can simultaneously do the work of two automatons, which is very rare in case of automatons. So in a single machine we have double operations to be performed by it.

*Example 7.1.12:* Let $Z = \{0, 1, 2, 3\}$, $A = A_1 \cup A_2$, where $A_1 = Z_2 = \{0, 1\}$ and $A_2 = \{0, 1, 2\} = Z_3$. $\delta : (z, a) = z.a$. Is $Y = (Z, A, \delta)$ a bisemi-automaton?

The table for Y is

|  | δ | 0 | 1 | 2 | 3 |
|---|---|---|---|---|---|
| $A_1$ | 0 | 0 | 0 | 0 | 0 |
|  | 1 | 0 | 1 | 2 | 3 |
| $A_2$ | 0 | 0 | 0 | 0 | 0 |
|  | 1 | 0 | 1 | 2 | 3 |
|  | 2 | 0 | 2 | 0 | 2 |

The revised table is

| δ | 0 | 1 | 2 | 3 |
|---|---|---|---|---|
| 0 | 0 | 0 | 0 | 0 |
| 1 | 0 | 1 | 2 | 3 |
| 2 | 0 | 2 | 0 | 2 |

But if we consider the semi-automaton $Y_1 = \{Z, A_1, \delta\}$ and $(Z, A_2, \delta) = Y_2$ we get the following graphs:

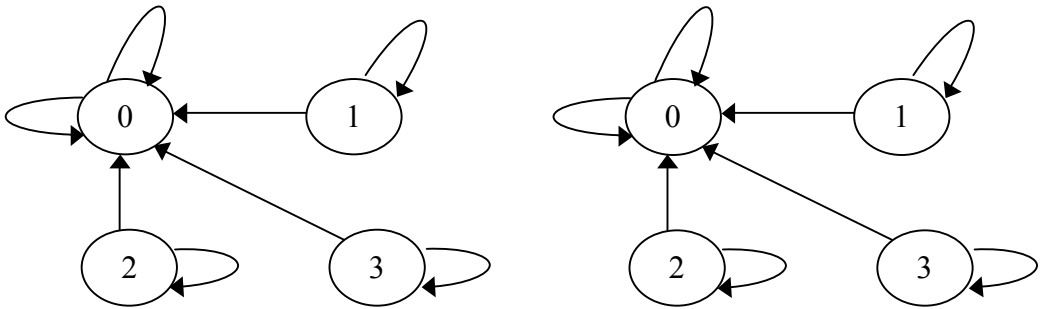

**Figures 7.1.16**

Suppose we take A = {0, 1, 2} and divide it as $A'_1$ = {0, 1} and $A'_2$ = {2} then we get the graphs as

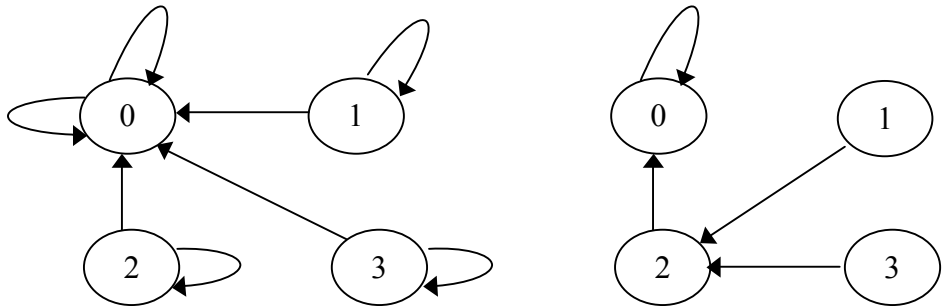

**Figures 7.1.17**

But $X'_1 = \{Z, A'_1, \delta\}$ and $X'_2 = \{Z, A'_2, \delta\}$ is a bisemi-automaton.

**Remark:** It is very important to note that how we choose the subsets $A_1$ and $A_2$ from A if $A_1 \subset A_2$ or $A_2 \subset A_1$. The semi-automaton will cease to be a bisemi-automaton.

Thus the very notion of bisets should have been defined even it may look trivial. So we at this juncture define bisets.

**DEFINITION 7.1.24:** *Let A be a set we say $A = A_1 \cup A_2$ is a biset if $A_1$ and $A_2$ are proper subsets of A such that $A_1 \not\subset A_2$ and $A_2 \not\subset A_1$.*

Unless this practical problem was not illustrated, the very notion of bisets would have been taken up researchers as a trivial concept. Only when bisets plays a vital role in determining whether a semi-automaton or an automaton is a bisemi-automaton or biautomaton.

Thus we give another formulation of the bisemi-automaton and biautomaton. Thus one of the necessary conditions for a semi-automaton to be a bisemi-automaton is that $A = A_1 \cup A_2$ is a biset. If $A = A_1 \cup A_2$ is not a biset then we do not have the semi-automaton to be a bisemi-automaton. But when we construct a bisemi-automaton A =

$A_1 \cup A_2$ where $\delta$ is a special mapping then the map gives different positions to the same input symbol depending on which subset the input symbol is.

To this end we give the following example:

***Example 7.1.13:*** Let $X = \{Z, A = A_1 \cup A_2, \delta\}$ be the bisemi-automaton, here $Z = \{0, 1, 2\}$, $A_1 = \{0, 1\}$ and $A_2 = \{0, 1, 2, 3\}$. $\delta(z, 0) = z + 0$ if $0 \in A_1$, $\delta(z, 0) = z.0$ if $0 \in A_2$, $\delta(z, 1) = z + 1$ if $1 \in A_1$, $\delta(z, 1) = z.1$ if $1 \in A_2$, $\delta(z, 2) = z.2 \pmod 3$, $\delta(z, 3) = z.3 \pmod 3$.

Now the diagram for the bisemi-automaton is as follows:

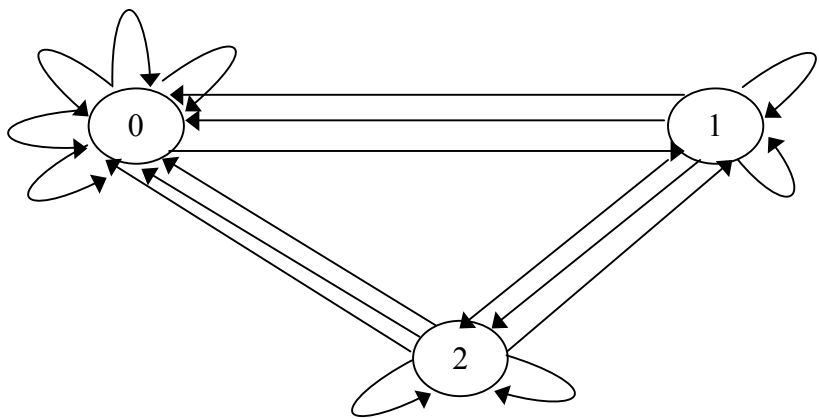

**Figure 7.1.18**

Thus this is a bisemi-automaton but not a semi-automaton as the function is not well defined. Thus we see all semiautomaton cannot be made into bisemiautomaton and in general a bisemiautomaton is not a semiautomaton. What best can be observed is that one can see that by building in a finite machine a bisemiautomaton structure one has the opportunity to make the machine perform more actions at a least cost and according to demand the bisets can be used with minute changes in the function. This bisemi-automaton will serve a greater purpose than that of the semi-automaton, secondly when we generate the free semigroups using the bisets or the sets in a bisemi-automaton we get nice results. Several applications can be found in this case.

Now at this juncture we wish to state that for the automaton also to make it a biautomatons both the input alphabets A and the output alphabets B must be divided into bisets to be the basic criteria to be a bi-automaton. The theory for this is similar to that of bisemi-automaton where in case of bisemi-automaton we used only the input alphabet A but in this case we will use both the input alphabets A and the output alphabets B.

Thus while converting a automaton to biautomaton we basically demand the sets A and B be divided into bisets such that the bisets give way to two distinct subautomaton. But if we are constructing biautomaton we define on $A = A_1 \cup A_2$ and $B = B_1 \cup B_2$ such that $\delta$ and $\lambda$ are not identical on $A_1 \cap A_2$ and $B_1 \cap B_2$ if they are non-empty. The reader is given the task of finding results in this direction for this will have a lot of application in constructing finite machines. Similarly when they are made into new biautomaton we take A to the free bisemigroup generated by $A_1$ and $A_2$ and B also to be the free bisemigroup generated by $B_1$ and $B_2$. We give two examples, one a automaton made into a biautomaton another is a constructed biautomaton.

***Example 7.1.14:*** Let $Y = (Z, A, B, \delta, \lambda)$ where $Z = \{0, 1, 2, 3\}$, $A = \{0, 1, 2, 3, 4, 5\}$ $= \{0, 2, 4\} \cup \{1, 3, 5\} = A_1 \cup A_2$, $B = \{0, 1, 2, 3\} = \{0, 2\} \cup \{1, 3\} = B_1 \cup B_2$.

$$\delta(z, a) = z + a \pmod 4 \text{ if } a \in A_1,$$
$$\delta(z, a) = z.a \pmod 4 \text{ if } a \in A_2.$$
$$\lambda(z, a) = z.a \pmod 4 \text{ if } a \in B_1,$$
$$\lambda(z, a) = z + a \pmod 4 \text{ if } a \in B_2.$$

Now we find the graphical representation of them.

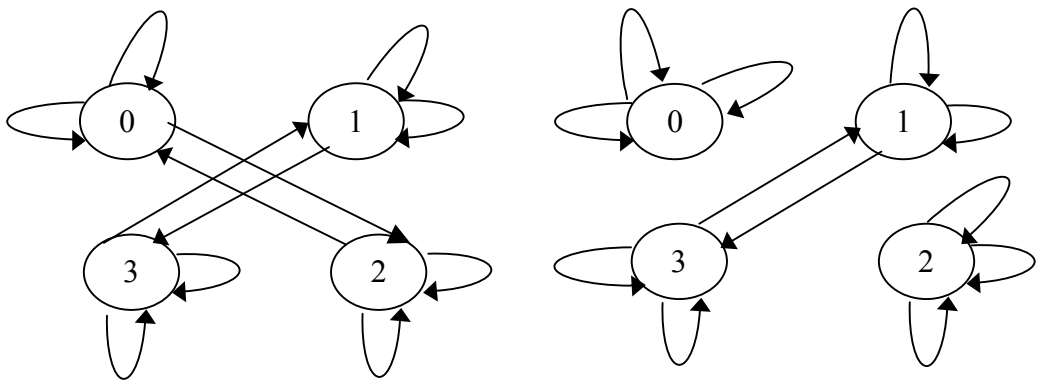

**Figures 7.1.19**

Output table:

| λ | 0 | 1 | 2 | 3 |
|---|---|---|---|---|
| 0 | 0 | 0 | 0 | 0 |
| 2 | 0 | 2 | 0 | 2 |

| λ | 0 | 1 | 2 | 3 |
|---|---|---|---|---|
| 1 | 0 | 1 | 2 | 3 |
| 3 | 0 | 3 | 2 | 1 |

We claim Y is the biautomaton as the output functions do not tally in both the tables. We give one more example of an automaton, which is a biautomaton.

***Example 7.1.15:*** Let $Y = (Z, A, B, \delta, \lambda)$ where $A = \{0, 1, 2\}$ and $B = \{0, 1, 2, 3\}$. $Z = \{0, 1, 2, 3, 4\}$. $\delta(z, a) = z.a \pmod 5$, $\lambda(z, a) = 2za \pmod 4$.

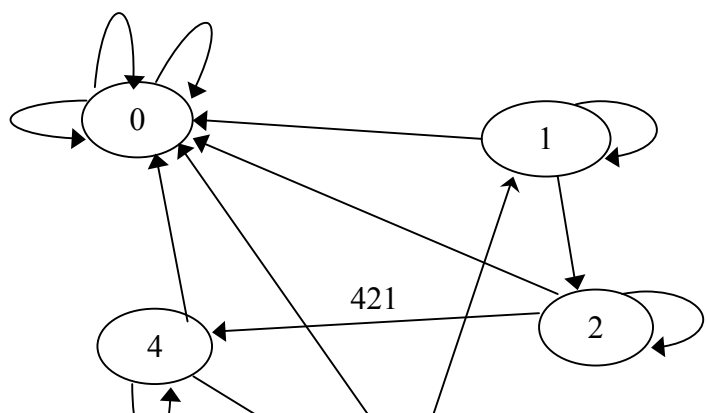

Next state function table:

| δ | 0 | 1 | 2 | 3 | 4 |
|---|---|---|---|---|---|
| 0 | 0 | 0 | 0 | 0 | 0 |
| 1 | 0 | 1 | 2 | 3 | 4 |
| 2 | 0 | 2 | 4 | 1 | 3 |

*Output table:*

| δ | 0 | 1 | 2 | 3 | 4 |
|---|---|---|---|---|---|
| 0 | 0 | 0 | 0 | 0 | 0 |
| 1 | 0 | 2 | 0 | 2 | 0 |
| 2 | 0 | 0 | 0 | 0 | 0 |
| 3 | 0 | 2 | 0 | 2 | 0 |

If we take $B = B_1 \cup B_2$ where $B_1 = \{0, 1\}$ and $B_2 = \{2, 3\}$ and $A = A_1 \cup A_2 = \{0, 1\} \cup \{2\}$ we get a biautomaton.

Thus using bisemigroups we obtain the notion of S-bisemi-automaton.

**DEFINITION 7.1.25:** $\overline{Y} = (Z, \overline{A}_s, \overline{\delta}_s)$ *is said to be a Smarandache bisemi-automaton (S-bisemi-automaton) if* $\overline{A}_s = \langle A \rangle$ *is the free groupoid generated by A with $\Lambda$ the empty element adjoined with it and $\overline{\delta}_s$ is the function from $Z \times \overline{A}_s \to Z$. Thus the S-bisemi-automaton contains $Y = (Z, \overline{A}, \overline{\delta})$ as a new bi semi-automaton which is a proper substructure of $\overline{Y}$.*

*Or equivalently, we define a S-semi-automaton as one, which has a new bisemi-automaton as its substructure.*

**DEFINITION 7.1.26:** $X_t = (Z_t, \overline{A}_s, \overline{\delta}_s)$ *is called the Smarandache subsemi-automaton (S-subsemi-automaton) of $X = (Z_p, \overline{A}_s, \overline{\delta}_s)$ denoted by $X_t \leq X$ if $Z_t \subset Z_p$ and $\overline{\delta}_s'$ is the restriction of $\overline{\delta}_s$ on $Z_t \times A_s$ and $X_t$ has a proper subset $\overline{H} \subset X_t$ such that H is a new semi-automaton.*

**DEFINITION 7.1.27:** $\overline{X} = (Z, \overline{A}_s, \overline{B}_s, \overline{\delta}_s, \overline{\lambda}_s)$ is said to be Smarandache biautomaton (S-biautomaton) if $\overline{A}_s = \langle A \rangle$ is the free bigroupoid generated by A with $\Lambda$ the empty element adjoined with it and $\overline{\delta}_s$ is the function from $Z \times \overline{A}_s \to Z$ and $\overline{B}_s$. Thus the S-biautomaton contains $X = (Z, \overline{A}, \overline{B}, \overline{\delta}, \overline{\lambda})$ as the new biautomaton which is a proper substructure of $\overline{X}$.

The notion of Smarandache sub-bi-automaton (S-sub-bi-automaton) and Smarandache sub-bi-semi-automaton (S-sub-bisemi-automaton) are defined in a similar way as that of S-subautomaton and S-subsemi-automaton.

We give the use of binear-rings in case of automatons in case of automatons.

Now when we proceed on to use the bisemigroup and binear-ring structure to define these concepts we get many multipurpose S-S-semi automatons i.e. a single S-S-semigroup automaton can serve several purposes as a finite machine.

**DEFINITION 7.1.28:** *Let $(G, +, \bullet)$ be a S-bisemigroup. A Smarandache S-bisemigroup semi bi-automaton (S-S-bisemigroup semi bi-automaton) is an ordered triple SB (S) = $(G = G_1 \cup G_2, I = I_1 \cup I_2, \delta = \delta_1 \cup \delta_2)$ where G is a S-bisemigroup, $I = I_1 \cup I_2$ is a set of inputs $I_1$ and $I_2$ are proper subsets of I $\delta_i$: $G_i \times I_i \to G_i$ is a transition function i = 1, 2.*

**Example 7.1.16:** $G = Z_3 \cup Z_4$, $I = \{0, 1, 2\} \cup \{a_1, a_2\}$; $\delta = (\delta_1, \delta_2)$ $\delta_1 : Z_3 \times \{2\} \to Z_3$, $\delta_2 : Z_4 \times \{a_1, a_2\} \to Z_4$ given by the following tables:

| $\delta_1$ | 0 | 1 | 2 |
|---|---|---|---|
| 0 | 0 | 0 | 0 |
| 1 | 0 | 1 | 2 |
| 2 | 0 | 2 | 1 |

| $\delta_2$ | 0 | 1 | 2 | 3 |
|---|---|---|---|---|
| $a_1$ | 0 | 2 | 0 | 2 |
| $a_2$ | 1 | 3 | 1 | 3 |

**(G, I, δ) is a S-bisemigroup bisemi-automaton. G is a S-bisemigroup under multiplication.**

Thus these two graphs together give the S-bisemigroup bisemi-automaton.

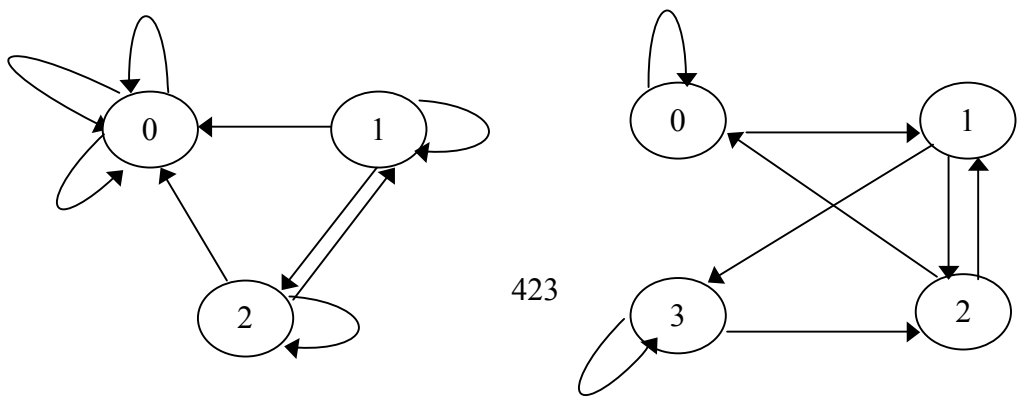

**DEFINITION 7.1.29:** *Let $SB(S) = (G = G_1 \cup G_2, I = I_1 \cup I_2, \delta = (\delta_1\ \delta_2))$ be a S-S-semigroup bisemi-automaton. A proper subset $B \subset G$ where $B = B_1 \cup B_2$, $(B_1 \subset G_1$ and $B_2 \subset G_2)$ is called the Smarandache S-sub-bisemigroup bisemi-automaton if $\delta/B_i : B_i \times I_i \to B_i$, $i = 1, 2$.*

The above example has no S-S-sub-bisemigroup bisemi-automaton.

**DEFINITION 7.1.30:** *A bigroup biautomaton is an ordered quintuple, $A = \{G, A, B, \delta = (\delta_1\ \delta_2), \lambda = (\lambda_1, \lambda_2)\}$ where $G$ is a bigroup i.e. $G = G_1 \cup G_2$, $A$ set of inputs, $B$ is set of outputs $\delta_i : G_i \times A \to G_i$, $\lambda_i : G_i \times A \to B$, $i = 1, 2$; is the transition and output functions respectively.*

The Smarandache analogue would be as follows:

**DEFINITION 7.1.31:** *A Smarandache S-bisemigroup biautomaton (S-S-bisemigroup biautomaton) is an ordered quintuple $SS(A) = (G = G_1 \cup G_2, A = A_1 \cup A_2, B = B_1 \cup B_2, \delta = (\delta_1, \delta_2), \lambda = (\lambda_1, \lambda_2))$ where $G$ is a S-bisemigroup and $A$ and $B$ are bisets which are input and output alphabets $\delta$ and $\lambda$ are the transition and output functions respectively.*

*If $\delta_i : G_i \times A_i \to G_i$ maps a sub-bigroup $T_i \subset G_i$ into itself $\delta_i : T_i \times A_i \to T_i$ then we can say the SS(A) bi-automaton has $(T \subset G, A, B, \delta, \lambda)$ to be a group bi-automaton.*

We can define level II bi-automaton in a similar way. The task of doing this is left as an exercise to the reader.

**DEFINITION 7.1.32:** *Let $S = (Q, X, \delta)$ be a S-S bisemigroup bisemi-automaton where $Q$ is a S-bisemigroup under '+' i.e. $G \subset Q$ is such that $G$ is a S-bisemigroup under '+' and $G$ is a proper subset of $Q$. Let $\delta_i : G_i \times X \to G_i$ is well defined that is restriction of $\delta$ to $G_i$ is well defined so that $(G, X, \delta) = C$ becomes a bisemi-automaton. $C$ is called additive if there is some $x_0 \in X$ with $\delta(q, x) = \delta(q, x_0) + \delta(0, x)$ and $\delta(q - q', x) = \delta(q, x) - \delta(q', x)$ for all $q, q' \in G = G_1 \cup G_2$ and $x \in X$.*

*Then there is some bi-homomorphism*

$$\psi : G \to G \text{ and some map}$$
$$\alpha : X \to G \text{ with}$$
$$\psi(x_0) = 0 \text{ and } \delta(q, x) = \psi(q) + \alpha(x).$$

*Let $\delta_x$ for a fixed $x \in X$ be the map from $G \to G$; $q \to \delta (q, x)$; then $\{\delta_x \mid x \in X\}$ generates a sub-binear-ring, $N(P)$ of the binear-ring. $\{M(G) = M(G_1) \cup M(G_2), +, \bullet\}$ of all mappings on $G = G_1 \cup G_2$, this binear-ring $N(P)$ is called the Smarandache bisyntactic binear-ring (S-bisyntactic binear-ring).*

Just as in case of S-semigroups we see in case of S-bisemigroups also we have several bigroup, bisemi-automaton depending on the number of valid bigroups in the S-semigroup. So for a given S-S-bisemigroup bisemi-automaton we can have several S-bisyntactic binear-ring.

Now we see the probable applications of binear-rings to error correcting codes and their Smarandache analogue.

**DEFINITION 7.1.33:** *A binear-ring $N = N_1 \cup N_2$ is called biplanar if for all equations $x_i \, o a = x_i o b + c$ $(a, b, c \in N_i \, a \neq b)$ have exactly one solution $x_i \in N_i$, $i = 1, 2$.*

Now we show how a biplanar binear-ring can be used to construct BIBD of high efficiency where by high efficiency E we mean $E = \lambda v / rk$ this E is a number between 0 and 1 and it estimates the quality of any statistical analysis if $E > 0.75$ the quality is good.

Here we give the construction of BIBD from a biplanar binear-ring.

Let N be a biplanar binear-ring. Let $a \in N = N_1 \cup N_2$. Define $g_a : N_1 \to N_1$ if $a \in N_1$ or $g_a : N_2 \to N_2$ if $a \in N_2$, by

$$n_1 \to n_1 \, o \, a \text{ if } a \in N_1$$
$$n_2 \to n_2 \, o \, a \text{ if } a \in N_2$$

and from $G = G_1 \cup G_2 = \{g_a \mid a \in N_1\} \cup \{g_a \mid a \in N_2\}$. Call $a \in N$, bigroup forming if $a \, o \, N_1$ is a subgroup or $a \, o \, N_2$ is a sub group of $(N_1, +)$ and $(N_2, +)$ respectively. Let us call sets $a.N_i + b$ $(a \in N_i; b \in N_i)$ blocks $i = 1, 2$. Then these blocks together with N as the set of points form a tactical configuration with parameters $(\upsilon, b, r, k, \lambda) =$

$$\left( (\upsilon_1, \frac{\alpha_1 \upsilon_1}{|G_1|} + \alpha_2 \upsilon_1, \alpha_2 + \alpha_2 |G_1|, |G_1|, \lambda), (\upsilon_2, \frac{\alpha_1 \upsilon_2}{|G_2|} + \alpha_2 \upsilon_2, \alpha_2 + \alpha_2 |G_2|, |G_2|, \lambda) \right)$$

where $\upsilon_i = |N_i|$ and $\alpha_1$ denotes the number of orbits of F under the group $G_1 \setminus \{0\}$. $\alpha_2$ denotes the number of orbits of F under the group $G_2 \setminus \{0\}$ which consists of entirely of bigroup forming elements. The tactical configuration is a BIBD if and only if either all elements are bigroup forming or just 0 is bigroup forming.

The main advantage in using the Smarandache biplanar binear-ring is that for one S-biplanar binear-ring we define more than one BIBD. If there are m-binear field in the S-biplanar binear-ring we can build m BIBD.

Thus the S-biplanar binear-ring has more advantages than the usual planar binear-rings and hence BIBD's construed using S-biplanar binear-rings will even prove to be an efficient error correcting codes.

Now we proceed on to recall some of the applications of Fuzzy algebra in section two.

## 7.2 Some applications of fuzzy algebraic structures and Smarandache fuzzy algebraic structures

The application of fuzzy algebra is mainly found in finite machines or which we choose to call as Smarandache fuzzy automaton. Apart from this we see that there is no direct application of Fuzzy algebra. When we say Fuzzy algebraic structures we mean only fuzzy groupoids, fuzzy semigroups, fuzzy groups, fuzzy rings, fuzzy near-rings, fuzzy seminear-rings, fuzzy semirings and fuzzy vector spaces. We don't mix up the concepts of fuzzy logic with fuzzy algebraic structures. The book by [95] entitled *Fuzzy Automaton and Languages, Theory and Applications*, gives an elaborate account of both the theory of automaton and algebraic fuzzy automaton. We would recall some notions of algebraic fuzzy automaton theory which will be directly helpful to us in constructing Smarandache algebraic fuzzy automatons.

We also give the probable Smarandache fuzzy algebraic structures then and there so that the reader can compare them and study them and make further research on these concepts.

*A fuzzy finite state machine (FFSM) is a triple $M = (Q, X, \mu)$ where $Q$ and $X$ are finite non-empty sets and $\mu$ is a fuzzy subset of $Q \times X \times Q$ i.e $\mu : Q \times X \times Q \to [0, 1]$. As usual $X^*$ denotes the set of all words of elements of $X$ of finite length $Q$ is called the set of states and $X$ is called the set of input symbols. Let $\wedge$ denote the empty words in $X^*$ and $|x|$ denote the length of $x$ for all $x \in X^*$, $X^*$ is a free semigroup with identity $\wedge$ with respect to the binary concatenation of two words.*

**DEFINITION [95]**: *Let $M = (Q, X, \mu)$ be a FFSM. Define $\mu^* : Q \times X^* \times Q \to [0, 1]$ by*

$$\mu^*(q, \wedge, p) = \begin{cases} 1 & \text{if } q = p \\ 0 & \text{if } q \neq p \end{cases}$$

*and $\mu^*(g, xa, p) = \vee \{\mu^*(q, x, r) \wedge \mu(r, a, p) \mid r \in Q\}$ for all $x \in X^*$ and $a \in X$. Let $X^+ = X^* \setminus \{\wedge\}$. Then $X^+$ is a semigroup. For $\mu^*$ given in here we let $\mu^+ = \mu^*$ restricted to $Q \times X^+ \times Q$.*

Several results about FFSM can be had from [95].

Analogous to the above definition of [95] we define Smarandache fuzzy finite state machine (SFFSM). Here we replace $X^*$ the free semigroup by $S(X^*)$ the free groupoid generated by X so that it is immediate $X^* \subseteq S(X^*)$, further $S(X^*)$ is a S-groupoid and has several free semigroups.

**DEFINITION 7.2.1:** *Let $M = (Q, X, \mu)$ be a FFSM. Define $S(\mu^*) : Q \times S(X^*) \times Q \to [0, 1]$ by $S(\mu^*)(q \wedge p) = \mu^*(q \wedge p)$ on $X^*$ and on elements in $S(X^*) \setminus X^*$. $S(\mu^*)$ need not be even defined or it can take any value between 0 and 1 i.e. (0, 1).*

*$S(\mu^*(q, xa, p)) = \mu^*(q, xa, p)$ for all $x \in X^*$, $a \in X$. We call $S(\mu^*(X))$ or $S(M)$ the SFFSM.*

*Now we can take the set X is FFSM to be say n sets of input symbols (n-finite). Say $X = X_1 \cup X_2 \cup ... \cup X_n$. '$\cup$' only for notation. Now we saw from definition 7.2.1, the machine FFSM is completely contained in SFFSM. But if we are interested to couple several machines in a single piece then SFFSM becomes handy. For if we have say $X_1, ..., X_n$; $n$ ($n < \infty$), set of distinct output symbols then the Smarandache fuzzy finite state machines (SFFSM) comes handy. For let $S(X^*) = \langle X_1 \cup ... \cup X_n \rangle$ denote the free groupoid generated by the set X. Clearly $S(X^*)$ is a S-groupoid and we see $S(X^*)$ has free semigroups $X_1^*, ..., X_n^*$ so that $S(\mu^*) : Q \times S(X^*) \times Q \to [0, 1]$ we give all the n-FFSM as*

$$\mu_1^* : Q \times X_1^* \times Q \to [0, 1]$$
$$\mu_2^* : Q \times X_2^* \times Q \to [0, 1]$$
$$\vdots$$
$$\mu_n^* : Q \times X_n^* \times Q \to [0, 1]$$

*thus we see all the n-FFSM are contained in the SFFSM. Further if we are interested in compartmenting a FFSM into submachines that also can be easily done by SFFSM. For if $X = \{x_1 ... x_n\}$ we want to divide this into say t machines with $t < n$ such that the t sets $T_1, ..., T_t$ are subsets of X with no conditions on $T_i$'s the only condition being*

$$X \subset \bigcup_{i=1}^{t} T_i$$

*$T_i \cap T_j \neq \phi$ even if $i \neq j$ in general. Then now we see $S(X^*)$ is a S-free groupoid having $T_i^*$ to be a free semigroups for $i = 1, 2, ..., n$. $S(\mu^*) : Q \times (X^*) \times Q \to [0, 1]$ gives restriction maps as*

$$\mu_1^* : Q \times T_1^* \times Q \to [0, 1]$$
$$\mu_2^* : Q \times T_2^* \times Q \to [0, 1]$$
$$\vdots$$
$$\mu_t^* : Q \times T_t^* \times Q \to [0, 1]$$

*We see these t-machines acts as sub-FFSM to SFFSM and one can easily use this method in dividing the networking and still use the FFSM as a single machines which can perform compartmental operations.*

Thus in the author's opinion this will be of a great help when applied to real world problems. We have already defined the concept of Smarandache semi-automaton we now proceed on to define the notion of Smarandache fuzzy bistate finite semi-automaton using the concept of S-bigroupoids.

Suppose we are interested in constructing finite machines where X is the set of input symbols should operate with two distinct binary operations that too not on the whole set X but on some part of X then we adopt Smarandache free bigroupoids to construct such machines.

*Given Q set of state, X the collection of input symbols $X = X_1 \cup X_2$ where $X_1$ and $X_2$ are proper subsets of X may be with some overlap but still they have to perform actions in a different way. Now we make the free bigroupoid to be generated by X is a very special way. If $S(X^*)$ denotes the free bigroupoid then $S(X^*) = S(X_1^*) \cup S(X_2^*)$ i.e. $S(X_1^*)$ and $S(X_2^*)$ are free groupoids on $X_1$ and $X_2$. Using this we build the Smarandache fuzzy finite bistate machines we say bistate because it can give two states for a same input symbols depending on from which free groupoids we pick them.*

*Let $SB(M) = (Q, S(X^*), S(\mu^*))$ where $S(X^*) = S(X_1^*) \cup S(X_2^*)$ with $S(X_1^*)$ and $S(X_2^*)$ free groupoids so that $S(X^*)$ is a Smarandache free bigroupoid (S-free bigroupoid).*

*Define $\mu^* = (\mu_1^* \cup \mu_2^*) : Q \times S(X^*) \times Q \to [0, 1]$ by*

$$\mu_1^*(q, \wedge, p) = \begin{cases} 1 & \text{if } q = p \\ 0 & \text{if } q \neq p \end{cases}$$

$$\mu_2^*(q, \wedge, p) = \begin{cases} 1 & \text{if } q = p \\ 0 & \text{if } q \neq p \end{cases}$$

*and $S(\mu_1^*)(q, xa, p) = \vee \{\mu_1^*(q, x, r) \wedge \mu_1(r, a, p) \mid r \in Q\}$ for all $x \in X_1^*$, $a \in X_1$, $\mu_1 : Q \times X_1 \times Q \to [0, 1]$ and $(\mu_2^*)(q, xa, p) = \vee \{\mu_2^*(q, x, r) \wedge \mu_2(r, a, p) \mid r \in Q\}$ for all $x \in X_2^*$, $a \in X_2$, $\mu_2 : Q \times X_2 \times Q \to [0, 1]$. Thus we call SB(M) to be a Smarandache fuzzy finite bistate machine (SFFBSM). Thus the notion of Smarandache fuzzy bigroupoids helps in the construction of SFFBSM.*

*All results pertaining to FFSM can be easily extended to SFFSM in an analogous way.*

Now we will just recall the definition of T-generalized state machine.

*A triple $(Q, X, \tau)$ is called a generalized state machine if Q and X are finite sets and $\tau : Q \times X \times Q \to [0, 1]$ is such that*

$$\sum_{q \in Q} r(p, a, q) < 1$$

*for all $p \in Q$ and $a \in X$. Let M be a generalized machine. Then the condition*

$$\sum_{q \in Q} \tau(p,a,p) \leq 1$$

*contains the crisp case in the sense that $\tau$ can be considered a partial function i.e. if Im ($\tau$) $\subseteq$ {0, 1} then $\tau$ is a partial function if there exists $p \in Q$ such that $\tau(p, a, q) = 0$ for all $q \in Q$ and $a \in X$. M is called complete if*

$$\sum_{q \in Q} \tau(p,a,p) = 1$$

*for all $p \in Q$ and $a \in X$. If M is not complete, it can be completed in the following measure. Let $Q' = Q \cup \{x\}$ where $x \notin Q$. Let $M^C = (Q', \times, \tau')$ where*

$$\tau'(p',a,q') = \begin{cases} \tau(p',a,q') & \text{if } p', q' \in Q \\ 1 - \sum_{q \in Q} \tau(p',a,q') & \text{if } p' \in Q \text{ and } q' = z \\ 0 & \text{if } p' = z \text{ and } q' \in Q \\ 1 & \text{if } p' = z \text{ and } q' = z \end{cases}$$

*for all $a \in X$. Let $X^+ = X^* \setminus \{\Lambda\}$. Let T be a t-norm on [0, 1]. Define $\tau^+ : Q \times X^+ \times Q \to [0, 1]$ by $\tau^+(p, a_1, ..., a_n, q) = \vee \{\tau(p, a, T_1) / \tau(r_1, a_2, r_2)T ... T\tau(r_{n-1} a_n, q) / r_i \in Q,, i = 1, 2, ..., n – 1\}$ where p and $q \in Q$ and $a_1 ... a_n \in X$. M is called a T-generalized state machine when $\tau^+$ is defined in terms of T.*

*Now we to get the Smarandache T-generalized state machine (S-T-generalized state machine T-generalized state machine) we replace $X^*$ the free semigroup by the free groupoid $S(X^*)$. Since free groupoids are S-free groupoids containing free semigroup we replace $X^+ = X^* \setminus \{\Lambda\}$ by $S(X^+) = S(X^*) \setminus (\Lambda)$. By using $S(X^*)$ in the $X^*$ we will get the Smarandache generalized state machine (S-generalized state machine) as SFFSM.*

Clearly we can define the Smarandache generalized state machine as follows:

**DEFINITION 7.2.2:** *A Smarandache generalized state machine $S(X^*)$ is one that contains a proper subset which is a generalized state machine. Thus we see all our Smarandache generalized state machines (S-generalized state machines) contain the generalized state machines properly.*

All results pertaining to these can be worked out as exercises.

Let $M = (Q, X, \mu)$ be a FFSM. Let E(M) be equal to $\{[x] \mid x \in X^*\}$ where $[x] = \{y \in X^* \mid x \equiv y\}$. (E(M), $*$) is a finite semigroup with the identity.

(Q, E(M), $\rho$) is a fuzzy transformation semigroup that we denote by FTS(M) where E(M) is a finite semigroup with identity [$\Lambda$]. Now using the FTS(M) we construct Smarandache fuzzy transformation S-groupoid (S-fuzzy transformation S-groupoid) associated with M. This is carried out in the following steps:

    i.      $X^*$ will be replaced by $S(X^*)$ the S-fuzzy free groupoid.

Subsequent changes are made.

**DEFINITION 7.2.3:** *Let $M = (Q, X, \mu)$ be a FFSM. Let $p, q \in Q$, $p$ is called an immediate successor of $q$ if there exists $a \in X$ such that $\mu(q, a, p) > 0$, $p$ is called a successor of $q$ if there exists $x \in X^+$ such that $\mu^*(q, x, p) > 0$. If $S(M) = (Q, S(X^*), \mu^*)$ be a SFFSM. Let $p, q \in Q$, $p$ is called an Smarandache immediate (S-immediate) of $q$ if there exists $a \in X$ such that $\mu(q, a, p) > 0$, $p$ is called a successor of $q$ if there exists $x \in S(X^*)$ such that $S(\mu^*(qxp)) > 0$.*

**DEFINITION 7.2.4:** *Let $M = (Q, X, \mu)$ be a FFSM. $M$ is called strongly connected if for all $p, q \in Q$, $p \in S(q)$. For more properties in this direction please refer [95]. Let $M = (Q, X, \mu)$ be a FFSM. Let $\delta$ be a fuzzy subset of $Q$. Then $(Q, \delta, X, \mu)$ is called a subsystem of $M$ if for all $p, q \in Q$ and for all $a \in X$. $\delta(q) \geq \delta(p) \wedge \mu(p, a, q)$. Let $S(M) = (Q, S(X^*), \mu^*)$ be a SFFSM. Let $\delta$ be a fuzzy subset of $Q$. Then $(Q, \delta, S(X^*), \mu^*)$ is called a Smarandache subsystem (S-subsystem) of $S(M)$ if for all $p, q \in Q$ and for all $a \in X$, $\delta(q) \geq \delta(p) \wedge \mu^*(p, a, q)$.*

Properties on S-subsystems can be obtained analogous to subsystems.

**DEFINITION [95]:** *A fuzzy semi-automaton (FSA) over a finite group $(Q, *)$ is a triple $(Q, X, \mu)$ where $X$ is a finite set and $\mu$ is a fuzzy subset of $Q \times X \times Q$.*

*We call $(Q, X, \mu)$ a Smarandache fuzzy semi-automaton (S-FSA) over a S-semigroup $(Q, *)$ where $X$ is a finite set and $\mu$ is a fuzzy subset of $Q \times X \times Q$. Thus if $(Q, X, \mu)$ is a S-fuzzy semi-automaton then it contains proper subsets which are fuzzy semi-automatons. Let $S = (Q, X, \mu)$ and $T = (Q_1, X_1, \mu_1)$ be a S-FSA over a S-semigroup $Q$. A pair of functions $(f, g)$ where $f : Q \to Q_1$, $g : X \to X_1$ are Smarandache homomorphism (S-homomorphism) from $S$ into $T$ written $(f, g) : S \to T$ if the following condition hold good*

    i.   *$f$ is a S-semigroup homomorphism.*
    ii.  *$\mu(p, x, q) \leq \mu_1(f(p), g(x), f(q))$ for all $p, q \in Q$, $x \in X$.*

*The pair $(f, g)$ is a Smarandache strong homomorphism (S-strong homomorphism) from $S$ into $T$ if it satisfies (i) of the above definition and the added condition $\mu_1(f(p), g(x), f(q)) = \vee \{\mu(p, x, r) \mid r \in Q, f(r) = f(q)\}$ for all $p, q \in Q' \subset Q$, $x \in X$, $Q'$ a subgroup of $Q$.*

**DEFINITION 7.2.5:** *A fuzzy subset $\lambda$ of $Q$ is called a Smarandache fuzzy subsemi-automaton (S-subsemi-automaton) of $S$ if*

    i.   *$\lambda$ is a S-fuzzy subsemigroup of $Q$.*
    ii.  *$\lambda(p) \geq \mu(q, x, p) \wedge \lambda(q)$*

*for all $p, q \in Q' \subset Q$ ($Q'$ a subgroup of $Q$) and $x \in X$.*

Having defined Smarandache analogue on similar lines an innovative reader can develop the Smarandache fuzzy analogue of all the Smarandache applications dealt in the section 7.1. This work is left as an exercise for the reader.

## 7.3 Problems

In this section we introduce only 16 problems. Each can be viewed or visualized as a research problem on the applications of Smarandache fuzzy algebraic structures. Not only these 16 problems still we can have several related ones and new ones by the reader.

**Problem 7.3.1:** Define and construct Smarandache cascades using S-automaton.

**Problem 7.3.2:** Give an example of a S-syntactic near-ring. Use it to define S-semigroup semi-automaton.

**Problem 7.3.3:** For a S-planar near-ring, $N = (Z_{144}, +, \bullet)$ using the BIBD associated with N construct error-correcting codes.

**Problem 7.3.4:** Give a general method of constructing BIBD for planar near-ring $(Z_p, +, \bullet)$, p a prime.

**Problem 7.3.5:** Give an example of a bi-automaton with 7 states.

**Problem 7.3.6:** Define Smarandache bicascades. Illustrate with examples.

**Problem 7.3.7:** Can two finite S-automaton be combined to get a S-bi-automaton?

**Problem 7.3.8:** Construct a Smarandache fuzzy finite state machine.

**Problem 7.3.9:** Construct a S-FFBSM.

**Problem 7.3.10:** Construct the Smarandache generalized state machine using $(Q, S(X^*), \tau)$.

**Problem 7.3.11:** Suppose $(Q, S(X^*), \tau)$ be a S-generalized state machine define congruence relation on $S(X^+)^* = S(X^* \setminus \Lambda)$.

**Problem 7.3.12:** Define Smarandache cascade using SFFSM.

**Problem 7.3.13:** Define the concept of Smarandache strongly connected.

**Problem 7.3.14:** Define the notion of Smarandache submachine. Illustrate it with examples.

**Problem 7.3.15:** Define Smarandache strong subsystems by giving explicit examples.

**Problem 7.3.16:** Can we define Smarandache fuzzy bi-automaton? Justify your claim.

# INDEX

# ABOUT THE AUTHOR

Dr. W. B. Vasantha is an Associate Professor in the Department of Mathematics, Indian Institute of Technology Madras, Chennai, where she lives with her husband Dr. K. Kandasamy and daughters Meena and Kama. Her current interests include Smarandache algebraic structures, fuzzy theory, coding/ communication theory. In the past decade she has completed guidance of seven Ph.D. scholars in the different fields of non-associative algebras, algebraic coding theory, transportation theory, fuzzy groups, and applications of fuzzy theory to the problems faced in chemical industries and cement industries. Currently, six Ph.D. scholars are working under her guidance. She has to her credit 241 research papers of which 200 are individually authored. Apart from this she and her students have presented around 262 papers in national and international conferences. She teaches both undergraduate and postgraduate students at IIT and has guided 41 M.Sc. and M.Tech projects. She has worked in collaboration projects with the Indian Space Research Organization and with the Tamil Nadu State AIDS Control Society. She is currently authoring a ten book series on Smarandache Algebraic Structures in collaboration with the American Research Press.

She can be contacted at vasantha@iitm.ac.in
You can visit her on the web at: http://mat.iitm.ac.in/~wbv